\documentclass{amsbook}
\usepackage[applemac]{inputenc}
\usepackage[french]{babel}
\usepackage{bbold,euscript}
\usepackage{stmaryrd,scalefnt,xfrac}
\usepackage{MnSymbol}
\usepackage{graphicx}
\usepackage[all,pdftex,2cell]{xy}
\UseAllTwocells
\xyoption{import}
\entrymodifiers={+!!<0pt,\fontdimen22\textfont2>}

\evensidemargin=\oddsidemargin
\renewcommand{\labelenumi}{(\roman{enumi})}
\newtheorem{theoreme}{Théorème}[subsection] 
\newtheorem{proposition}[theoreme]{Proposition}
\newtheorem{corollaire}[theoreme]{Corollaire}
\newtheorem{lemme}[theoreme]{Lemme}

\newtheorem{innercustomthm}{Proposition}
\newenvironment{customthm}[1]
  {\renewcommand\theinnercustomthm{#1}\innercustomthm}
  {\endinnercustomthm}

\theoremstyle{definition}

\usepackage{fancyhdr}
\headwidth=\textwidth

\newcommand{\cat}{\mathbf{cat}}
\newcommand{\Cat}{\mathbf{Cat}}
\newcommand{\CAT}{\mathbf{CAT}}
\newcommand{\ens}{\mathbf{Ens}}
\newcommand{\simp}{\widehat{\mathbf{\,{\bf\Delta}\,}}}
\newcommand{\ssimp}{\widehat{\mathbf{{\bf\Delta}}\times\mathbf{{\bf\Delta}}}}

\newcommand{\M}{\mathcal{M}}
\newcommand{\C}{\mathcal{C}}

\newcommand{\F}{\mathcal{F}}
\newcommand{\D}{\mathcal{D}}

\newcommand{\A}{\mathcal{A}}
\newcommand{\B}{\mathcal{B}}
\newcommand{\G}{\mathcal{G}}
\renewcommand{\H}{\mathcal{H}}

\entrymodifiers={+!!<0pt,\fontdimen22\textfont2>}
\makeindex

\begin{document}
\title{Le déterminant de Deligne \\et le $2$-groupe d'homotopie \\ des espaces connexes}
\author{Elhoim Sumano}
\maketitle


\tableofcontents
\phantom{a}\thispagestyle{empty}\newpage
\pagenumbering{arabic}
\chapter*{Introduction}  

Soit $\C = \big(\C,{\bf cof},{\bf W}\big)$ une catégorie de Waldhausen. Si $\G$ est un $2$-groupe (voir \ref{2groupessection}) un \emph{determinant} (unitaire et non-commutatif) {de $\C$ à valeurs dans $\G$} est un couple $(D,T)$ formé d'un foncteur et une fonction:
$$
\xymatrix{\C_{\bf W} \ar[r]^-{D} & \G}
\qquad\text{et}\qquad
\xymatrix{\mathrm{s}_2(\C) \ar[r]^-{T} & \{\text{Morphismes de $\G$}\}}\,,
$$
où $\C_{\bf W}$ note la sous-catégorie de $\C$ dont les morphisme sont les élément de ${\bf W}$ et $\mathrm{s}_2(\C)$ est l'ensemble des suites cofibres de $\C$:
$$
\xymatrix@C+3pt{A\; \ar@{>->}[r]& B\ar@{->>}[r]&C}\;,
$$ 
tels que:
\begin{enumerate}
\item (Compatibilité) \, Si $\xi = \big(\vcenter{\xymatrix@C-10pt{A\;\ar@{>->}[r]&B\ar@{->>}[r]&C}}\big)$ est une suite cofibre de $\C$ le morphisme $T(\xi)$ de $\G$ est un morphisme de la forme:
$$
\xymatrix@C+8pt{D_0(A) \otimes D_0(C)  \ar[r]^-{T(\xi)} & D_0(B) \,.}
$$
\item (Unitaire) \, On a que $D_0\big( 0 \big)\,=\,\mathbb{1}$ et $T\big(\xymatrix@C-10 pt{0\;\ar@{>->}[r]&0\ar@{->>}[r]&0}\big)\,=\,l_{\mathbb{1}}^{-1}\,=\,r_{\mathbb{1}}^{-1}$.
\item (Naturalité)\, Si on se donne un morphisme de suites cofibre de $\C$:
$$
\vcenter{\xymatrix@C+10pt@R-5pt{\xi \ar[d]\\ \xi'}}
\; = \; \vcenter{\xymatrix@C+10pt@R-5pt{
A\;\ar[d]_-{f_A}\ar@{>->}[r]&B\ar[d]|-{f_B}\ar@{->>}[r]&C\ar[d]^-{f_C}\\
A'\;\ar@{>->}[r]&B'\ar@{->>}[r]&C'}}
$$ 
où $f_A$, $f_B$ et $f_C$ appartient à ${\bf W}$, on a un carré commutatif:
$$
\xymatrix@C+10pt{
D_0(A) \otimes D_0(C)\ar[d]_-{D_1(f_A)\otimes D_1(f_C)} \ar[r]^-{T(\xi)} & D_0(B)   \ar[d]^-{D_1(f_B)}\\
D_0(A') \otimes D_0(C') \ar[r]_-{T(\xi')} & D_0(B') 
}
$$
\item (Associativité)\, Étant un diagramme de $\C$:
$$
\xymatrix@R-4pt@C+1pt{
&& A_{3,2}\\ & A_{2,1} \; \ar[r] & A_{3,1}\ar[u] \\
A_{1,0} \; \ar[r] & A_{2,0} \; \ar[r] \ar[u]  & A_{3,0} \ar[u] 
}
$$
tel que: 
$$
\xi_0 = \big(\vcenter{\xymatrix@C-10pt{A_{2,1}\;\ar@{>->}[r]&A_{3,1}\ar@{->>}[r]&A_{3,2}}}\big)\;,\qquad \xi_1 = \big(\vcenter{\xymatrix@C-10pt{A_{2,0}\;\ar@{>->}[r]&A_{3,0}\ar@{->>}[r]&A_{3,2}}}\big)\;,
$$ 
$$
\xi_2 = \big(\vcenter{\xymatrix@C-10pt{A_{1,0}\;\ar@{>->}[r]&A_{3,0}\ar@{->>}[r]&A_{3,1}}}\big)\qquad \text{et} \qquad 
\xi_3 = \big(\vcenter{\xymatrix@C-10pt{A_{1,0}\;\ar@{>->}[r]&A_{2,0}\ar@{->>}[r]&A_{2,1}}}\big)
$$ 
sont de suites cofibre de $\C$, on a un diagramme commutatif de $\G$:
\begin{equation}\label{pentafinipliintro}
\vcenter{\xymatrix@R-10pt@C+35pt{
D_0(A_{3,0}) &D_0(A_{1,0})\otimes D_0(A_{3,1}) \ar[l]_-{T(\xi_2)}\\&\\
&D_0(A_{1,0}) \otimes \big(D_0(A_{2,1})\otimes D_0(A_{3,2})\big) \ar[uu]_-{D_0(A_{1,0})\otimes T(\xi_0) }\ar@{}[d]_-{a}|-{\text{\rotatebox[origin=c]{90}{$\cong$}}}\\
D_0(A_{2,0})\otimes D_0(A_{3,2}) \ar[uuu]^-{T(\xi_1)} &\big(D_0(A_{1,0})\otimes D_0(A_{2,1}))\otimes D_0(A_{3,2}\big)\ar[l]^-{T(\xi_3)\otimes D_0(A_{3,2})}\;.}}
\end{equation}
\end{enumerate}

Si $(D,T)$ et $(D',T')$ sont deux déterminants de $\C$ á valeurs dans $\G$ un \emph{morphisme de déterminants} de $\C$ à valeurs dans $\G$ est une transformation naturelle $\xymatrix@C+8pt{\C_{\bf W}\rtwocell^{D}_{D'}{\alpha}& \G}$ vérifiant les propriétés:
\begin{enumerate}
\item[(v)] (Unitaire) \, $\alpha_0 \, = \, \mathrm{id}_\mathbb{1}$ 
\item[(vi)] (Compatibilité) \, Si $\xi = \big(\vcenter{\xymatrix@C-10pt{A\;\ar@{>->}[r]&B\ar@{->>}[r]&C}}\big)$ est une suite cofibre de $\C$, on a un carré commutatif de $\G$:
$$
\xymatrix@C+15pt@R-3pt{
D_0(A) \otimes D_0(C)\ar[d]_-{\alpha_A  \otimes \alpha_C} \ar[r]^-{T(\xi)} & D_0(B)   \ar[d]^-{\alpha_B}\\
D'_0(A) \otimes D'_0(C) \ar[r]_-{T'(\xi)} & D'_0(B) 
}
$$
\end{enumerate}

On vérifie sans peine que les déterminants de $\C$ à valeurs dans $\G$ et les morphismes entre eux forment un groupoïde qu'on note $\underline{{\bf det}}_{\C}(\G)$. La composition étant la composition de transformation naturelles.

Remarquons qu'on a un foncteur:
\begin{equation}\label{tututumaintro}
\xymatrix@R=5pt@C+35pt{
\text{$2$-${\bf Grp}$} \ar[r]^-{\underline{\bf det}_{\C}}  &  {\bf Grpd}\,,\\
\G \;\, \ar@{}[r]|-{\longmapsto} & \underline{\bf det}_{\C}(\G)}
\end{equation}
défini comme suit: Si $(\varphi,m^{\varphi})\colon\xymatrix@C-10pt{\G\ar[r]&\H}$ est un morphisme de $2$-groupes le foncteur:
$$
\xymatrix@C+35pt@R=5pt{
\underline{\bf det}_\C(\G) \ar[r]^-{\underline{\bf det}_\C(\varphi,m^\varphi)} & \underline{\bf det}_\C(\H)
}
$$
est définie dans les objets par la formule $\underline{\bf det}_\C\big(\varphi,m^\varphi\big)(D,T)\,=\,(\overline{D},\overline{T})$ où $\overline{D}$ est le foncteur composé:
$$
\xymatrix@C+5pt{
\C_{\bf W}  \ar[r]^-{D} & \G \ar[r]^-{\varphi} & \H}
$$
et $\overline{T}$ est la fonction:
$$
\xymatrix@C+15pt@R=2pt{
\mathrm{s}_2(\C)  \ar[r]^-{\overline{T}} & \{\text{Morphismes de $\H$}\}  \;. \\
\xi \ar@{}[r]|-{\longmapsto} &  \varphi\big(T(\xi)\big) \circ m^{\varphi} 
}
$$

On note $h{\bf Grpd}$ la \emph{catégorie homotopique des groupoïdes} \emph{i.e.} la catégorie qu'on obtient de la catégorie des groupoïdes ${\bf Grpd}$ en identifient deux foncteurs s'il existe une transformation naturelle entre eux. 

Si $\pi\colon\xymatrix@C-8pt{{\bf Grpd}  \ar[r]& h{\bf Grpd}}$ est le foncteur canonique, l'énoncé suivant est bien connu (voir \cite{deligne}, \cite{knudsen} et \cite{muro}):

\begin{customthm}{A}\label{intropropW}
Si $\C$ est une catégorie de Waldhausen le foncteur composé:
$$
\xymatrix@R=5pt@C+3pt{
\text{$2$-${\bf Grp}$} \ar[rr]^-{\underline{{\bf det}}_\C(\,\bullet\,)} & &  {\bf Grpd}  \ar[r]^-{\pi}& h{\bf Grpd}
}
$$
est représentable par le $2$-groupe d'homotopie (voir la fin de \ref{lenerfSegalss}) de l'ensemble simplicial:
\begin{equation}\label{WaldhKK}
\vcenter{\xymatrix@C+12pt@R=1pt{
\Delta^{op} \ar[r] & \ens\\ [p]\ar@{}[r]|-{\longmapsto} & \mathrm{Hom}_{\bf cat}\big([p],{\bf W}\mathrm{S}_{p}(\C)\big) }}
\end{equation}
dont les groupes d'homotopie sont les $K$-groupes de $\C$.

En particulier si $(\varphi,m^{\varphi})\colon \xymatrix@C-10pt{\G\ar[r]& \H}$ est une $2$-équivalence faible de $2$-groupes le foncteur:
$$
\xymatrix@C+35pt{
 \underline{{\bf det}}_X(\G) \ar[r]^-{ \underline{{\bf det}}_X(\varphi,m^{\varphi})} &  \underline{{\bf det}}_X(\H)
}
$$
est une équivalence faible de groupoïdes et le foncteur induit:
$$
\xymatrix@R=5pt@C+40pt{
\text{$2$-$h{\bf Grp}$} \ar[r]^-{\pi_0\big(\underline{{\bf det}}_X(\,\bullet\,)\big)}   &  {\bf Ens}}
$$
est représentable aussi par le $2$-groupe d'homotopie de l'ensemble simplicial \eqref{WaldhKK}.
\end{customthm} 

On trouve dans l'article \cite{muro2} une preuve de la Proposition \ref{intropropW} pour tous les objets $\A_\bullet$ de la catégorie des foncteurs ${\bf cat}^{\Delta^{op}}$ tels que $\A_0=\star$. 

\vspace{.5cm}

\renewcommand{\thesubsection}{\S\thesection.\arabic{subsection}}
\subsubsection*{{\bf I}}
\renewcommand{\thesubsection}{\thesection.\arabic{subsection}}

Soit ${\bf MS}$ la sous-catégorie pleine de la catégorie des ensembles bisimpliciaux dont les objets sont les $X$ tels que $X_{\bullet,0}=\star$ (voir \ref{premonsegal}). La catégorie $\mathbf{MS}$ admet l'enrichissement simpliciaux suivant: 
\begin{equation}\label{hom1intro}
\begin{split}
\underline{\mathrm{Hom}}_{\mathbf{MS}}^{(1)} \big( X, Y \big)_n
\,=&\,\mathrm{Hom}_{\mathbf{MS}}\Big( \big(X\times \mathrm{p}_1^*(\Delta^n)\big)\big/ \big(\star\times \mathrm{p}_1^*(\Delta^n) \big)  ,Y  \Big) \\
\,=&\,\mathrm{Hom}_{\ssimp}\Big(\big(X\times \mathrm{p}_1^*(\Delta^n)\big)\big/ \big(\star\times \mathrm{p}_1^*(\Delta^n) \big) ,Y\Big) \\
\,=&\,\mathrm{Hom}_{\ssimp}\big(X\times \mathrm{p}_1^*(\Delta^n) ,Y\big) \;.
\end{split}
\end{equation}

Si $0\leq n\leq \infty$ dans la Proposition \ref{moduno} on note que la catégorie ${\bf MS}$ munie de l'enrichissement $\underline{\mathrm{Hom}}^{(1)}_{{\bf MS}}$ admet une structure de catégorie de modèles simplicial pointée lorsque:
\begin{align*}
\big\{\,\text{équivalences faibles}\,\big\} \quad &=\quad 
\big\{\, f:\xymatrix@-14pt{X\ar[r]&Y}\,\big|\; \text{$\mathrm{diag}(f)$ est une $n$-équivalence}\\
 &\phantom{=}\;\;\quad \text{\phantom{$\{\,$}faible d'ensembles simpliciaux}\,\big\}\;\,=\;\, {\bf W}_{n}^{diag}\,,\\
\big\{\,\text{cofibrations}\,\big\} \quad &=\quad  \big\{\,\text{monomorphismes}\,\big\} \;\;\,= \;\,  {\bf mono}\,,\\
\big\{\,\text{fibrations}\,\big\} \quad &=\quad  \big\{\,\text{morphismes avec la propriété de relèvement }\\
 &\phantom{=}\;\;\quad \text{\phantom{$\{\,$}à droite par rapport à ${\bf mono}\,\cap\,{\bf W}_{n}^{diag}$}\,\big\}\;\,= \;\, {\bf fib}_{n}^{diag}\,.
\end{align*}

En fait la catégorie de modèles $({\bf MS},{\bf W}^{diag}_{n}, {\bf mono},{\bf fib}_{n}^{diag})$ est le remplacement simplicial de \cite{dugger} de la catégorie de modèles $\big(\simp_0,{\bf W}^{red}_{n}, {\bf mono},{\bf fib}^{red}_{n}\big)$ des ensembles simpliciaux réduits (voir \ref{ensreduitlacat}). En particulier la catégorie homotopique ${\bf MS}\big[({\bf W}^{diag}_{n})^{-1}\big]$ est équivalente à la catégorie homotopique des $n$-types d'homotopie connexes.

\vspace{.5cm}

\renewcommand{\thesubsection}{\S\thesection.\arabic{subsection}}
\subsubsection*{{\bf II}}
\renewcommand{\thesubsection}{\thesection.\arabic{subsection}}

Dans le Chapitre 3 de cette thèse on se rappelle que le foncteur $2$-nerf des bicatégories défini dans \cite{paoli} induit par restriction aux $2$-groupes un foncteur pleinement fidèle: 
\begin{equation}\label{nerfSintro}
\xymatrix@C+20pt{\text{$2$-${\bf Grp}$} \ar[r]^-{\mathcal{N}_{\mathcal{S}}} & {\bf MS}\;.}
\end{equation}

Si $\G$ est un $2$-groupe on montre dans la Proposition \ref{fibrachii} que l'ensemble bisimplicial $\mathcal{N}_{\mathcal{S}}(\G)$ est un objet fibrant de la catégorie de modèles simplicial $\big({\bf MS},{\bf W}^{diag}_{2}, {\bf mono},{\bf fib}_{2}^{diag},\underline{\mathrm{Hom}}^{(1)}_{{\bf MS}}\big)$ et dans le Lemme \ref{igrphom1MS} on montre que pour tout objet $X$ de ${\bf MS}$ l'ensemble simplicial des morphismes $\underline{\mathrm{Hom}}^{(1)}_{{\bf MS}}(X,\mathcal{N}_{\mathcal{S}}(\G))$ est isomorphe au nerf d'un groupoïde (un groupoïde de déterminants de $X$ à valeurs dans $\G$ d'après le Corollaire \ref{fininininifer}).

Enfin on vérifie que si $\C$ est une catégorie de Waldhausen et $X$ note l'ensemble bisimplicial:
$$
\vcenter{\xymatrix@C+12pt@R=1pt{
(\Delta\times \Delta)^{op} \ar[r] & \ens\\ \big([p],[q]\big)\ar@{}[r]|-{\longmapsto} & \mathrm{Hom}_{\bf cat}\big([p],{\bf W}\mathrm{S}_{q}(\C)\big) }}
$$
alors l'ensemble simplicial $\underline{\mathrm{Hom}}^{(1)}_{{\bf MS}}(X,\mathcal{N}_{\mathcal{S}}(\G))$ est isomorphe au nerf du groupoïde $\underline{{\bf det}}_\C(\G)$ ci-dessus. 

Vu que le foncteur \eqref{nerfSintro} induit une équivalence de catégories (voir le Corollaire \ref{equiv2grpMS2}):
$$
\xymatrix@C+15pt{
\text{$2$-$h{\bf Grp}$} \ar[r]^-{h\mathcal{N}_{\mathcal{S}}} & {\bf MS}\big[({\bf W}^{diag}_{n})^{-1}\big]
}
$$
on en déduit une autre preuve de la Proposition \ref{intropropW} (valable pour tout objet de la catégorie ${\bf MS}$).

\vspace{.5cm}

\renewcommand{\thesubsection}{\S\thesection.\arabic{subsection}}
\subsubsection*{{\bf III}}
\renewcommand{\thesubsection}{\thesection.\arabic{subsection}}

Rappelons qu'un catégorie de modèles simpliciale pointée est une catégorie pointée $\C$ munie d'une structure de catégorie de modèles $(\C,{\bf W}, {\bf cof},{\bf fib})$ et d'un enrichissement $\underline{\mathrm{Hom}}_{\C}$ sur la catégorie de modèles monoïdale $\big(\simp,{\bf W}_{\infty}, {\bf mono},{\bf fib}_{\infty}\big)$ lequel vérifie certains conditions (voir \ref{catmodsimpnn}). 

Si $n\geq 1$ considérons la catégorie de modèles monoïdale $\big(\simp,{\bf W}_{n-1}, {\bf mono},{\bf fib}_{n-1}\big)$ dont la catégorie sous-jacent est la catégorie cartésienne fermée des ensembles simpliciaux, les cofibrations sont les monomorphismes et les équivalences faibles sont les $(n-1)$-équivalences faibles. 

Dans le Chapitre 2 de cette thèse on montre (Théorème \ref{cassi1}) l'énoncé suivant: 

\begin{customthm}{B}
Si $n\geq 1$ et $(\C,{\bf W}, {\bf cof},{\bf fib},\underline{\mathrm{Hom}}_{\C})$ est une catégorie de modèles simpliciale pointée telle que tous ses objets sont cofibrants, les énoncés suivants sont équi\-va\-lents:
\begin{enumerate}
\item L'enrichissement $\underline{\mathrm{Hom}}_{\C}$ de la catégorie de modèles $(\C,{\bf W}, {\bf cof},{\bf fib})$ est en fait un enrichissement sur la catégorie de modèles monoïdale $\big(\simp,{\bf W}_{n-1}, {\bf mono},{\bf fib}_{n-1}\big)$.
\item Si $m\geq n$ la $m$-ième itération d'un foncteur suspension $\Sigma^m\colon \xymatrix@C-10pt{\mathrm{Ho}(\C)\ar[r]&\mathrm{Ho}(\C)}$ est null.
\item Si $m\geq n$ la $m$-ième itération d'un foncteur espace de lacets $\Omega^m\colon \xymatrix@C-10pt{\mathrm{Ho}(\C)\ar[r]&\mathrm{Ho}(\C)}$ est null.
\item Si $X$ et $Y$ sont des objets de $\C$ quelconques, $\pi_{m} \big(\mathbf{R}\mathrm{Hom}_{\C}(X,Y),\star_{X}^Y\big)=0$ pour $m\geq n$.
\item Si $Y$ est un objet fibrant de $(\C,{\bf W}, {\bf cof},{\bf fib})$ l'ensemble simplicial $\underline{\mathrm{Hom}}_{\C}(X,Y)$ est un objet fibrant de la catégorie de modèles $(\simp,{\bf W}_{n-1}, {\bf mono},{\bf fib}_{n-1})$ pour tout objet $X$ de $\C$.
\end{enumerate} 
\end{customthm}

On en déduit (voir la Proposition \ref{moduno}) que l'enrichissement $\underline{\mathrm{Hom}}_{\mathbf{MS}}^{(1)}$ de la catégorie de modèles $\big({\bf MS},{\bf W}^{diag}_{n}, {\bf mono},{\bf fib}_{n}^{diag}\big)$ se restreint à $\big(\simp,{\bf W}_{n-1}, {\bf mono},{\bf fib}_{n-1}\big)$. En particulier si $Y$ est un objet fibrant de la catégorie de modèles $\big({\bf MS},{\bf W}^{diag}_{n}, {\bf mono},{\bf fib}_{n}^{diag}\big)$, pour tout objet $X$ de ${\bf MS}$ l'ensemble simplicial des morphismes $\underline{\mathrm{Hom}}_{\mathbf{MS}}^{(1)} \big( X, Y \big)$ est un objet fibrant de la catégorie de modèles $\big(\simp,{\bf W}_{n-1}, {\bf mono},{\bf fib}_{n-1}\big)$ \emph{i.e.} un complexe de Kan dont les groupes d'homotopie $\pi_i$ sont nuls pour tout $i\geq n$. 

Autrement dit l'espace des morphisme dérivé ${\bf R}_n\mathrm{Hom}$ de $\big({\bf MS},{\bf W}^{diag}_{n}, {\bf mono},{\bf fib}_{n}^{diag}\big)$ est un enrichissement de la catégorie homotopique ${\bf MS}\big[({\bf W}^{diag}_{n})^{-1}\big]$ sur la catégorie des $(n-1)$-types d'homotopie $\simp\big[({\bf W}_{n-1})^{-1}\big]$. 

Remarquons que pour $n=2$ on constat par ce moyen que si $X$ est un objet quelconque de la catégorie ${\bf MS}$ le foncteur des morphismes dérivé ${\bf R}_2\mathrm{Hom}(X, \bullet)$ est le foncteur des "déterminants de $X$" de la catégorie homotopique des $1$-groupes $\big($qui est équivalente à la catégorie ${\bf MS}\big[({\bf W}^{diag}_{2})^{-1}\big]\big)$ vers la catégorie homotopique des groupoïdes $\big($qui est équivalente à la catégorie des $1$-types d'ho\-mo\-to\-pie$\big)$, et si $n=1$ le foncteur ${\bf R}_1\mathrm{Hom}(X, \bullet)$ est le foncteur des "fonctions additives de $X$" de la catégorie des groupes $\big($qui est équivalente à la catégorie ${\bf MS}\big[({\bf W}^{diag}_{1})^{-1}\big]\big)$ vers la catégorie des ensembles $\big($qui est équivalente à la catégorie des $0$-types d'homotopie$\big)$ (voir \ref{ffonctionadditetdet}).

Donc si $n\geq 3$ le foncteur ${\bf R}_n\mathrm{Hom}(X, \bullet)$ est un foncteur des "déterminants supérieurs de $X$" de  la catégorie homotopique des $n$-types d'homotopie connexes (la dite catégorie homotopique des $n$-groupes) vers la catégorie des $(n-1)$-types d'homotopie.

\chapter{Les $K$-groupes des catégories de modèles et des dérivateurs}

\section{Catégories de modèles}


\renewcommand{\thesubsection}{\S\thesection.\arabic{subsection}}
\subsection{}\label{catmode}\; 
\renewcommand{\thesubsection}{\thesection.\arabic{subsection}}

Rappelons qu'une \emph{catégorie de modèles} $\C = (\C,{\bf W}, {\bf cof},{\bf fib})$ (voir \cite{hirschhorn} ou \cite{hovey}), est la donnée d'une catégorie $\C$ avec trois familles distinguées de morphismes: ${\bf W}$, ${\bf cof}$ et ${\bf fib}$ dont les éléments sont dits respectivement équivalences faibles, cofibrations et fibrations; et qui satisfont aux propriétés suivantes:
\begin{itemize}\label{momode}
\item[{\bf CM1}.] La catégorie $\C$ admet de petites limites et colimites.
\item[{\bf CM2}.] La famille ${\bf W}$ satisfait la \emph{propriété de deux-sur-trois}; c'est-à-dire, si deux des trois morphismes dans un triangle commutatif:
$$
\xymatrix@R-8pt@C-6pt{X\ar[rr]^-f\ar[dr]_-{h}&&Z\,,\ar[ld]^-g\\&Y&}
$$
appartiennent à ${\bf W}$; il en est de même pour le troisième. 
\item[{\bf CM3}.] Les familles ${\bf W}$, ${\bf cof}$ et ${\bf fib}$ sont stables par rétractes; plus explicitement, étant donné un diagramme commutatif:
$$
\xymatrix{
A \ar[r]\ar@/^12pt/[rr]^{\mathrm{id}}\ar[d]_-{g} & X \ar[d]|-f\ar[r] & A\ar[d]^-{g}\,\\
B \ar[r]\ar@/_12pt/[rr]_{\mathrm{id}} & Y \ar[r] & B\,;}
$$
si $f$ est un élément de ${\bf W}$, ${\bf cof}$ ou ${\bf fib}$, il en est de même pour $g$. 
\item[{\bf CM4}.] Les éléments de ${\bf fib}$ (resp. de ${\bf cof}$) vérifient la \emph{propriété de relèvement à droite} (resp. \emph{à gauche}) \emph{par rapport aux éléments} de ${\bf cof}\cap{\bf W}$ (resp. ${\bf fib}\cap{\bf W}$); c'est-à-dire, étant donné un carré commutatif:
$$
\xymatrix{A\ar[r]\ar[d]_-{i}&\,X\,\ar[d]^-p\\B\ar@{-->}[ru]|-{\varphi}\ar[r]&\,Y\,,}
$$
avec $i\in{\bf cof}$ et $p\in{\bf fib}$; si $i\in{\bf W}$ (resp. si $p\in{\bf W}$), il existe le morphisme $\varphi$ comme dans le diagramme, rendant commutatifs les triangles qui en résultent. 
\item[{\bf CM5}.]  Il existe de factorisations fonctorielles de tout morphisme $f$ de $\C$:
\begin{enumerate}
\item $f=p\, i$, où $p\in{\bf fib}$ et $i\in{\bf cof}\cap{\bf W}$. 
\item $f=q\, j$, où $q\in{\bf fib}\cap {\bf W}$ et $j\in{\bf cof}$.
\end{enumerate} 
\end{itemize}

On note respectivement $\xymatrix@-10pt{\phantom{a}\ar@{->}[r]^-{\text{\rotatebox[origin=c]{180}{\Large $\widetilde{\phantom{w}}$}}}&\phantom{a}}$, $\xymatrix@-10pt{\phantom{a}\ar@{>->}[r]&\phantom{a}}$ et $\xymatrix@-10pt{\phantom{a}\ar@{->>}[r]&\phantom{a}}$, les flèches des familles ${\bf W}$, ${\bf cof}$ et ${\bf fib}$.


\renewcommand{\thesubsection}{\S\thesection.\arabic{subsection}}
\subsection{}\;\label{remplaS}
\renewcommand{\thesubsection}{\thesection.\arabic{subsection}}

Si $\C$ est une catégorie de modèles, grâce à la propriété {\bf CM1} on peut choisir dans $\C$ un objet initial $\emptyset$ et un objet final $\star$. Un objet $X$ de $\C$ est dit \emph{cofibrant} (resp. \emph{fibrant}), si l'unique flèche $\xymatrix@C-8pt{\emptyset\ar[r]&X}$ est une cofibration $\big($resp. $\xymatrix@C-8pt{X\ar[r]&\star}$ est une fibration$\big)$. 

Il résulte de la propriété {\bf CM5} que pour tout objet $X$ de $\C$, on peut trouver un diagramme commutatif:
\begin{equation}\label{rempla}
\vcenter{\xymatrix@R-9pt@C+5pt{\underset{\phantom{a}}{\emptyset\,}\ar[r]\ar@{>->}[dr]_{j_{X}}&X\\&QX\ar[u]_{q_{X}}^{\text{\rotatebox[origin=c]{270}{\Large $\widetilde{\phantom{w}}$}}}}}
\qquad\quad\left(\;\text{resp.}\qquad
\vcenter{\xymatrix@R-9pt@C+5pt{X\ar[r]\ar[d]_{i_{X}}^{\text{\rotatebox[origin=c]{90}{\Large $\widetilde{\phantom{w}}$}}}&\star\\RX\ar@{->>}[ru]_{p_{X}}&}}
\right),
\end{equation}
où $QX$ est un objet cofibrant (resp. $RX$ fibrant) de $\C$. On appelle une telle donnée, formée d'un objet cofibrant $QX$ (resp. fibrant $RX$) et d'une équivalence faible $q_X$ (resp. $i_X$), un \emph{remplacement cofibrant de $X$} (resp. un \emph{remplacement fibrant}).

Soulignons que par hypothèse, on peut supposer qu'on s'est donné dans $\C$ un foncteur et une transformation naturelle:
\begin{equation}\label{laQ}
Q: \xymatrix@-6pt{\C\ar[r]&\C}\qquad \text{et} \qquad q: \xymatrix@-12pt{Q\ar@{=>}[r]&\mathrm{id}_{\C}}
\end{equation}
$$
\big( \text{resp.} \quad R: \xymatrix@-6pt{\C\ar[r]&\C} \qquad \text{et} \qquad i: \xymatrix@-12pt{\mathrm{id}_{\C}\ar@{=>}[r]&R} \big)\,,
$$
avec la propriété que pour tout objet $X$ de $\C$, le morphisme $q_{X}$ (resp. $i_{X}$) soit une fibration (resp. cofibration) et le couple $(Q(X),\,q_{X})$ $\big($resp. $(R(X),\,i_{X})$ $\big)$ un remplacement cofibrant (resp. fibrant) de $X$; on vérifie en particulier que l'objet $QX$ (resp. $RX$) est aussi fibrant (resp. cofibrant) si $X$ l'est.

D'autre part, rappelons qu'un \emph{objet cylindre} d'un objet $X$ de $\C$ (resp. \emph{objet cocylindre} de $Y$), est la donnée d'un objet ${\bf Cyl}(X)$ $\big($resp. ${\bf coCyl}(Y)$ $\big)$ et d'une factorisation:
\begin{equation}\label{cyl}
\vcenter{\xymatrix@-8pt{\underset{\phantom{.}}{X \sqcup X}\ar@<-4pt>@{>->}@/_5pt/[rd]_-{i_{0}+i_{1}}\ar[rr]^-{\mathrm{id}+\mathrm{id}} && X \\ 
&{\bf Cyl}(X) \ar@<-4pt>@/_5pt/[ru]_-{p}\ar@{}[ru]|-{\text{\rotatebox[origin=c]{225}{\Large $\widetilde{\phantom{w}}$}}}&}}
\qquad\qquad\left( \, \text{resp.}\qquad
\vcenter{\xymatrix@-8pt{Y \ar@<-4pt>@/_5pt/[rd]_-{i}\ar@{}[rd]|-{\text{\rotatebox[origin=c]{135}{\Large $\widetilde{\phantom{w}}$}}}\ar[rr]^-{(\mathrm{id},\mathrm{id})} && Y\times Y \\ &{\bf coCyl}(Y)\ar@<-4pt>@{->>}@/_5pt/[ru]_-{(p_{0},p_{1})}&}}
\, \right)\,.
\end{equation}

Étant donnés deux morphismes $\xymatrix@-7pt{X\ar[r]^-{f,\,g}&Y}$ de $\C$, on dit que $f$ est \emph{homotope à gauche} à $g$ (resp. \emph{homotope à droite}), s'il existe une \emph{homotopie à gauche} (resp. \emph{homotopie à droite}) de $f$ vers $g$; c'est-à-dire, s'il existe un objet cylindre de $X$ (resp. cocylindre de $Y$) et un morphisme $H$ qui s'insère dans un triangle commutatif: 
\begin{equation}\label{homo}
\vcenter{\xymatrix{X\sqcup X\ar[r]^-{f+g}\ar[d]_-{i_{0}+i_{1}} & Y \\ {\bf Cyl}(X)\ar@<-4pt>@/_5pt/@{-->}[ru]_-{H}&}}
\qquad\qquad\left(\text{resp.}\qquad
\vcenter{\xymatrix{X\ar[r]^-{(f,g)}\ar@<-5pt>@/_5pt/@{-->}[rd]_-{H} & Y\times Y \\ & {\bf coCyl}(Y)\ar[u]_-{(i_{0},i_{1})}}}
\right).
\end{equation}

On vérifie que dans la sous-catégorie $\C_{cf}$ de $\C$, dont les objets sont les objets à la fois fibrants et cofibrants, la relation d'homotopie à gauche et d'homotopie à droite coïncident et sont une relation d'équivalence. Vu qu'en plus cette relation est compatible à la composition, on en déduit une catégorie quotient $\pi\C_{cf}$ munie d'un foncteur $\xymatrix@C-6pt{\C_{cf}\ar[r]&\pi\C_{cf}}$. 

Enfin, si $\C[{\bf W}^{-1}]$ désigne la catégorie qu'on obtient de $\C$ après avoir inversé les équivalences faibles (voir \cite{GZ}), on montre l'existence de carrés commutatifs:
$$
\vcenter{\xymatrix@C+6pt{\, \C_{cf} \,\ \ar[d] \ar@{^(->}[r] &\C\ar[d]\\ \pi\C_{cf} \ar@{-->}[r]&\C[{\bf W}^{-1}] }}\qquad\text{et}\qquad
\vcenter{\xymatrix@C+6pt{\C\ar[d] \ar[r]^{RQ} &\C_{cf}\ar[d]\\ \C[{\bf W}^{-1}]\ar@{-->}[r]& \pi\C_{cf}}}\,,
$$
qui définissent une équivalence de catégories $\pi\C_{cf} \simeq \C[{\bf W}^{-1}]$.

Il se suit en particulier qu'un morphisme $\xymatrix@-5pt{X\ar[r]^-{f}&Y}$ dans une catégorie de modèles est une équivalence faible, si et seulement si son image dans $\C[{\bf W}^{-1}]$ est un isomorphisme. 

Dans le cadre des catégories de modèles, on note $\mathrm{Ho}_{\bf W}(\C)=\mathrm{Ho}(\C)$ la catégorie des fractions $\C[{\bf W}^{-1}]$ et on l'appelle la \emph{catégorie homotopique} de $\C$. On note aussi $[X,Y]_{\C}$, l'ensemble des morphismes dans $\mathrm{Ho}_{\bf W}(\C)$ entre deux objets $X$ et $Y$ de $\C$. 

\renewcommand{\thesubsection}{\S\thesection.\arabic{subsection}}
\subsection{}\;
\renewcommand{\thesubsection}{\thesection.\arabic{subsection}}

Soit $\C$ une catégorie quelconque munie d'une famille distinguée de morphismes ${\bf W}$ (par exemple une catégorie de modèles avec ses équivalences faibles). Notons $\C[{\bf W}^{-1}]$ la catégorie qu'on obtient de $\C$ après avoir inversé les éléments de ${\bf W}$, et $\gamma:\xymatrix@-10pt{\C\ar[r]&\C[{\bf W}^{-1}]}$ le foncteur canonique.

Si $\xymatrix@C-8pt{\C\ar[r]^F&\A}$ est un foncteur de but une catégorie arbitraire, rappelons qu'un \emph{foncteur dérivé à gauche} (resp. \emph{à droite}) de $F$ est la donnée  d'un foncteur $\xymatrix@C-6pt{\C[{\bf W}^{-1}]\ar[r]^-{{L}F}&\A}$ $\Big($resp. $\xymatrix@C-6pt{\C[{\bf W}^{-1}]\ar[r]^-{{R}F}&\A}$ $\Big)$, et d'une transformation naturelle $\xymatrix@-10pt{{\bf L}F\circ \gamma \ar@{=>}[r]^-{\alpha}&F}$ $\Big($resp. $\xymatrix@-10pt{F \ar@{=>}[r]^-{\alpha}&{\bf R}F\circ \gamma }\Big)$, telle que pour tout foncteur $\xymatrix@C-10pt{\C[{\bf W}^{-1}]\ar[r]^-{G}&\A}$, l'application suivante soit bijective:
$$
\xymatrix@C+4pt@R=1pt{
\mathrm{Hom}_{\text{\scriptsize{$\A^{\C[{\bf W}^{-1}]}$}}}\big(G,{L}F\big)\ar[r]&
\mathrm{Hom}_{\text{\scriptsize{$\A^{\C}$}}}\big(G\circ\gamma,F\big)\\
\qquad\beta\qquad\ar@{}[r]|{\longmapsto}& \alpha\circ(\beta\star\gamma)}
$$
$$
\left(\vcenter{\xymatrix@C+4pt@R=1pt{
\text{resp.}\;\;&
\mathrm{Hom}_{\text{\scriptsize{$\A^{\C[{\bf W}^{-1}]}$}}}\big({R}F,G\big)\ar[r]&
\mathrm{Hom}_{\text{\scriptsize{$\A^{\C}$}}}\big(F,G\circ\gamma\big)\\&
\qquad\beta\qquad\ar@{}[r]|{\longmapsto}& (\beta\star\gamma)\circ \alpha}}\right).
$$

Il en résulte que si $F$ envoie les éléments de ${\bf W}$  vers les isomorphismes de $\A$, c'est-à-dire, si on a une factorisation:
$$
\xymatrix@C+10pt{ \C\ar[d]_-{\gamma}\ar[r]^F & \A \,,\\   \C[{\bf W}^{-1}]\ar@<-4pt>@/_4pt/[ru]_-{\widehat{F}}&}
$$
alors $(\widehat{F},\mathrm{id})$ est à la fois un foncteur dérivé à gauche et à droite de $F$.

Si $\D$ est une autre catégorie munie d'une famille distinguée de morphismes ${\bf W}$, un \emph{foncteur dérivé total à gauche} (resp. \emph{à droite}) d'un foncteur  $\xymatrix@C-8pt{\C\ar[r]^F&\D}$, est par définition un foncteur dérivé à gauche (resp. à droite) $\xymatrix@C-6pt{\C[{\bf W}^{-1}]\ar[r]^-{{\bf L}F}&\D[{\bf W}^{-1}]}$ (resp. $\xymatrix@C-6pt{\C[{\bf W}^{-1}]\ar[r]^-{{\bf R}F}&\D[{\bf W}^{-1}]}$) du composé:
$$
\xymatrix@C-2pt{\C\ar[r]^-{F}&\D\ar[r]&\D[{\bf W}^{-1}]}.
$$ 

Supposons maintenant que $\C$ et $\D$ sont des catégories de modèles. Rappelons qu'une adjonction: 
\begin{equation}\label{adjQ}
\xymatrix@C+15pt{\C\ar@/^12pt/[r]^{F}\ar@{}[r]|-{\perp}&\D \,,\ar@/^12pt/[l]^{G}}
\end{equation}
est dite une \emph{adjonction de Quillen}, si le foncteur adjoint à gauche $F$ préserve les éléments des familles ${\bf cof}$ et ${\bf cof}\cap{\bf W}$, ou de fa\c con équivalente, si le foncteur adjoint à droite $G$ préserve les éléments des familles ${\bf fib}$ et ${\bf fib}\cap{\bf W}$.

Étant donnée une adjonction de Quillen \eqref{adjQ}, on peut construit un foncteur dérivé total à gauche de $F$ (resp. à droite de $G$):
$$
\vcenter{\xymatrix{\mathrm{Ho}(\C)\ar[r]^-{{\bf L}F}&\mathrm{Ho}(\D)}} \qquad\;
\left(resp. \quad \vcenter{\xymatrix{\mathrm{Ho}(\D)\ar[r]^-{{\bf R}G}&\mathrm{Ho}(\C)\,}}\right),
$$
avec la propriété que le morphisme $\xymatrix{{\bf L}F(X)\ar[r]^-{\alpha_{X}}&F(X)}$ $\big($resp. $\xymatrix{G(Y)\ar[r]^{\alpha_{Y}}&{\bf R}G(Y)}$$\big)$ soit un isomorphisme, pour tout objet $X$ de $\C$ cofibrant (resp. objet $Y$ de $\D$ fibrant). En effet, pour cela on pose ${\bf L}F(X)=F(QX)$ (resp. ${\bf R}G(Y)=G(RY)$) et $\alpha_{X}=F(q_{X})$ (resp. $\alpha_{Y}=G(i_{Y})$).

On en déduit en particulier une adjonction (voir par exemple \cite{quillenhh}):
\begin{equation}
\xymatrix@C+15pt{\mathrm{Ho}(\C)\;\ar@/^15pt/[r]^{{\bf L}F}\ar@{}[r]|-{\perp}&\;\mathrm{Ho}(\D) \ar@/^15pt/[l]^{{\bf R}G}\,. }
\end{equation}

Enfin, rappelons qu'une adjonction de Quillen \eqref{adjQ} est dite une \emph{équivalence de Quillen}, si le foncteur ${\bf L}F$ (ou le foncteur ${\bf R}F$) est une équivalence de catégories. De fa\c con équivalente, l'adjonction \eqref{adjQ} est une équivalence de Quillen si pour tout objet cofibrant $X$ de $\C$ et tout objet fibrant $Y$ de $\D$, un morphisme $\xymatrix@-5pt{FX\ar[r]&Y}$ est une équivalence faible de $\D$ si et seulement si, son morphisme adjoint $\xymatrix@-5pt{X\ar[r]&GY}$ est une équivalence faible de $\C$. 

\section{Extensions de Kan homotopiques}

\renewcommand{\thesubsection}{\S\thesection.\arabic{subsection}}
\subsection{}\;
\renewcommand{\thesubsection}{\thesection.\arabic{subsection}}

Soit $\C$ une catégorie arbitraire. Si $I$ est une petite catégorie, on note $\C^{I}$ la catégorie des \emph{$I$-diagrammes de $\C$}, c'est-à-dire, des foncteurs $\mathfrak{X}$ de $I$ vers $\C$. Étant donné un objet $a$ de $I$, si $\mathfrak{X}$ est un $I$-diagramme de $\C$, on désigne par $\mathfrak{X}_{a}$ la valeur du foncteur $\mathfrak{X}$ dans $a$.

Supposons maintenant que $\C$ est une catégorie munie d'une famille distinguée de morphismes ${\bf W}$. La catégorie $\C^{I}$ a alors une famille distinguée de morphismes ${\bf W}_I$, dont les éléments sont les \emph{équivalences faibles argument par argument}; \emph{i.e.} les morphismes de $I$-diagrammes $F:\xymatrix@C-12pt{\mathfrak{X}\ar[r]&\mathfrak{X}'}$ tels que pour tout objet $a$ de $I$, la flèche $F_{a}:\xymatrix@C-12pt{\mathfrak{X}_{a}\ar[r]&\mathfrak{X}'_{a}}$ soit une équivalence faible.

On construit en particulier un diagramme commutatif:
\begin{equation}\label{WW}
\xymatrix{\C^{I}\ar[r]^-{(\gamma_{\C})^I}\ar[d]_-{\gamma_{\C^I}}&\big(\C[{\bf W}^{-1}]\big)^I\,.\\\C^I[{\bf W}_{I}^{-1}]\ar@/_10pt/@{-->}[ru]_{\Phi_{\C,I}}&}
\end{equation}

\begin{lemme}\label{reflee}   
Soit $\C$ une catégorie munie d'une famille distinguée de morphismes ${\bf W}$ fortement saturée\footnote{On dit que la famille ${\bf W}$ est \emph{fortement saturée} si on a la propriété suivante: Un morphisme de $\C$ appartient à ${\bf W}$ si et seulement si son image dans $\C[{\bf W}^{-1}]$ est un isomorphisme.}, alors pour toute petite catégorie $I$ le foncteur canonique:
\begin{equation}\label{cons} 
\xymatrix{\C^I[{\bf W}_{I}^{-1}] \ar[r]^-{\Phi_{\C,I}}& \big(\C[{\bf W}^{-1}]\big)^I\,,}
\end{equation}
est conservatif; \emph{i.e.} si $F:\xymatrix@C-12pt{\mathfrak{X}\ar[r]&\mathfrak{X}'}$ est un morphisme de $\C^I[{\bf W}_{I}^{-1}]$ tel que pour tout objet $a$ de $I$ le morphisme $F_{a}:\xymatrix@C-12pt{\mathfrak{X}_{a}\ar[r]&\mathfrak{X}'_{a}}$ soit un isomorphisme de $\C[{\bf W}]$, alors $F$ est un isomorphisme de $\C^I[{\bf W}_{I}^{-1}] $.
\end{lemme}
\begin{proof}
Rappelons que se donner un morphisme $F:\xymatrix@C-12pt{\mathfrak{X}\ar[r]&\mathfrak{X}'}$ de $\C^I[{\bf W}_{I}^{-1}]$ signifie se donner une chaîne de morphismes de $\C^I$:
$$
\xymatrix{\mathfrak{X}=\mathfrak{X}^0\ar@{-}[r]^-{F^1}&\mathfrak{X}^1 \ar@{-}[r]^-{F^2}&\cdots \ar@{-}[r]^{F^{n-1}}& \mathfrak{X}^{n-1}\ar@{-}[r]^-{F^n}&\mathfrak{X}^n=\mathfrak{X}'}
$$ 
où pour chaque $1\leq k\leq n$, soit $F^k$ est un morphisme arbitraire de $\C^I$ de la forme $F^k:\xymatrix@C-12pt{\mathfrak{X}^{k-1}\ar[r]&\mathfrak{X}^k}$, soit $F^k$ est un morphisme de ${\bf W}_I$ de la forme $F^k:\xymatrix@C-12pt{\mathfrak{X}^{k}\ar[r]&\mathfrak{X}^{k-1}}$. 

Supposons que $\Phi_{\C,I}(F)$ est un isomorphisme de $\big(\C[{\bf W}^{-1}]\big)^I$ c'est-à-dire supposons que pour tout objet $a$ de $I$ la chaîne:
$$
\xymatrix{\mathfrak{X}_a=\mathfrak{X}_a^0\ar@{-}[r]^-{F_a^1}&\mathfrak{X}^1 _a\ar@{-}[r]^-{F^2_a}&\cdots \ar@{-}[r]^{F^{n-1}_a}& \mathfrak{X}^{n-1}_a\ar@{-}[r]^-{F^n_a}&\mathfrak{X}^n_a=\mathfrak{X}'_a}
$$  
est un isomorphisme de $\C[{\bf W}^{-1}]$. 

Vu que ${\bf W}$ est une famille de morphismes fortement saturée, on déduit que pour tout objet $a$ de $I$ et tout $0\leq k \leq n$, le morphisme $F^k_a$ appartient à ${\bf W}$. Donc, $F^k$ appartient à ${\bf W}_I$ pour tout $0\leq k \leq n$, c'est-à-dire $F$ est un isomorphisme de $\C^I[{\bf W}_{I}^{-1}] $.
\end{proof}

D'autre part, étant donné un foncteur $u:\xymatrix@C-10pt{I\ar[r]&J}$ entre petites catégories, vu que le foncteur induit $u^*:\xymatrix@C-10pt{\C^J\ar[r]&\C^I}$ respecte les équivalences faibles argument par argument, on déduit un foncteur qu'on note par le même symbole:
\begin{equation}\label{eto}
\xymatrix@C+15pt{\C^J[{\bf W}_{J}^{-1}] \ar[r]^{u^*}&\C^I[{\bf W}_{I}^{-1}]\,.}
\end{equation}

On dit que $\C$ \emph{admet des extensions de Kan homotopiques à gauche} (resp. \emph{droite}) \emph{le long de $u$}  (\emph{par rapport aux morphismes de ${\bf W}$}) si le foncteur $u^*$ admet un adjoint à gauche (resp. à droite):
$$
\xymatrix@C+15pt{
\C^J[{\bf W}_{J}^{-1}]\ar@/_14pt/[r]_{u^*}\ar@{}[r]|-{\perp}&\C^I[{\bf W}_{I}^{-1}]\ar@/_14pt/[l]_{{u_{!}}}}
\qquad\left(\;\text{resp.}\quad
\vcenter{
\xymatrix@C+15pt{
\C^J[{\bf W}_{J}^{-1}]\ar@/^14pt/[r]^{u^*}\ar@{}[r]|-{\perp}&\C^I[{\bf W}_{I}^{-1}]\ar@/^14pt/[l]^{{u_{*}}}}}\right)\,.
$$
On appelle un tel foncteur $u_{!}$ $($resp. $u_{*}$ $)$ un (ou le) \emph{foncteur extension de Kan homotopique à gauche} (resp. \emph{droite}) \emph{le long de $u$} (\emph{par rapport aux morphismes de ${\bf W}$}).

Étant donnée une petite catégorie $I$, on désigne par $p_{I}:\xymatrix@C-10pt{I\ar[r]&e}$ le seul foncteur de $I$ vers la catégorie ponctuelle $e$. Si $\C$ admet des extensions de Kan homotopiques à gauche (resp. droite) le long d'un foncteur de la forme $p_{I}:\xymatrix@C-10pt{I\ar[r]&e}$:
$$
\xymatrix@C+18pt{
\C[{\bf W}^{-1}]\ar@/_14pt/[r]_{p_{I}^*}\ar@{}[r]|-{\perp}&\C^I[{\bf W}_{I}^{-1}]\ar@/_14pt/[l]_{{(p_{I})_{!}}\,=\,\mathrm{hocolim}_{I}}}
\qquad\left(\;\text{resp.}\quad
\vcenter{
\xymatrix@C+15pt{
\C[{\bf W}^{-1}]\ar@/^14pt/[r]^{p_{I}^*}\ar@{}[r]|-{\perp}&\C^I[{\bf W}_{I}^{-1}]\ar@/^14pt/[l]^{{(p_{I})_{*}}\,=\,\mathrm{holim}_{I}}}}\right)\,,
$$
on dit que $\C$ \emph{admet des colimites homotopiques} (resp. \emph{des limites homotopiques}) \emph{de type $I$}, et on appelle l'objet $\mathrm{hocolim}_{I}(\mathfrak{X})\big($resp. $\mathrm{holim}_{I}(\mathfrak{X})\big)$ une (ou la) \emph{colimite homotopique} (resp. la \emph{limite homotopique}) du $I$-diagramme $\mathfrak{X}$ de $\C$. 

Rappelons:

\begin{lemme}\label{prosomm}
Soit $\C$ une catégorie munie d'une famille distinguée de morphismes ${\bf W}$ laquelle contient les identités. Alors pour toute petite catégorie $I$ isomorphe à une réunion disjoint finie $I\cong I_1\bigsqcup \cdots \bigsqcup I_n$ où $n\geq 0$, le foncteur canonique:
\begin{equation}\label{fonctcanoniqq}
\xymatrix@C+25pt{\C^I[{\bf W}_{I}^{-1}] \ar[r]^-{(u_1^*,\dots,u_n^*)}& \C^{I_1}[{\bf W}_{I_1}^{-1}] \times \cdots \times \C^{I_n}[{\bf W}_{I_n}^{-1}]\,,}
\end{equation}
est un isomorphisme de catégories. 
\end{lemme}
\begin{proof}
Si $n=0$ c'est-à-dire si $I$ est la catégorie vide on a que $\C^I$ est la catégorie ponctuelle, alors \eqref{fonctcanoniqq} est un isomorphisme de catégories parce que le produit vide de catégories est aussi la catégorie ponctuelle. D'un autre si $n=1$ le foncteur \eqref{fonctcanoniqq} est l'identité. 

Finalement rappelons que si $\A$ et $\B$ sont de catégories munies des familles distinguées de morphismes ${\bf W}_{\A}$ et ${\bf W}_{\B}$ respectivement lesquelles contient les morphismes identités, on montre que le foncteur canonique:
$$
\xymatrix@C+25pt{
(\A\times \B) \big[({\bf W}_{\A}\times {\bf W}_{\B})^{-1}\big]\ar[r] & 
\A[{\bf W}_{\A}^{-1}] \times \B[{\bf W}_{\B}^{-1}]}
$$
est un isomorphisme de catégories, en vérifiant que pour tout foncteur $\xymatrix@C-5pt{\A\times \B \ar[r]^-{F} & \C}$ avec la propriété:
$$
F\big({\bf W}_{\A}\times {\bf W}_{\B}\big) \; \text{\Large  $\subseteq$} \; \big\{ \; \text{isomorphismes  de  $\C$} \; \big\}
$$ 
il existe un seul foncteur:
$$
\xymatrix@C+25pt{
\A[{\bf W}_{\A}^{-1}] \times \B[{\bf W}_{\B}^{-1}] \ar[r]^-{\widetilde{F}} & \C}
$$
tel que $\widetilde{F}\circ \big(\gamma_{\A}\times \gamma_{\B}\big) = F$.

Il se suit en particulier que le foncteur \eqref{fonctcanoniqq} est un isomorphisme si $n\geq 2$.
\end{proof}

D'après ce dernier affirmation, si $I$ est une petite catégorie discrète finie\footnote{$I$ est une petite catégorie discrète finie si elle n'admet pas de morphismes autres que les identités et l'ensemble de ses objets est fini (possiblement vide).} et $\C$ est une catégorie de modèles; alors $\C$ admet des colimites homotopiques (resp. limites homotopiques) de type $I$ si et seulement si la catégorie $\C[{\bf W}^{-1}]$ admet des sommes (resp. des produits) indexées par l'ensemble d'objet $I_{0}$ de $I$. 

Dans ce cas, une colimite homotopique (resp. une limite homotopique) d'un $I$-diagramme $\mathfrak{X}=\{X_{i}\}_{i\in I_{0}}$ de $\C$ est canoniquement isomorphe à une somme (resp. un produit) de la famille $\{X_{i}\}_{i\in I_{0}}$ dans $\C[{\bf W}^{-1}]$.

Montrons maintenant:

\begin{lemme}\label{cwgen}
Soit $\C$ une catégorie munie d'une famille distinguée de morphismes ${\bf W}$. On considère un foncteur $u:\xymatrix@C-10pt{I\ar[r]&J}$ entre petites catégories, admettant un adjoint à droite (resp. à gauche) $v:\xymatrix@C-10pt{J\ar[r]&I}$. Alors, l'adjonction $v^* \dashv u^* : \xymatrix@-10pt{\C^J\ar[r]&\C^I}$ $\big($resp. $u^* \dashv v^* : \xymatrix@-10pt{\C^I\ar[r]&\C^J}\big)$ induit une adjonction:
\begin{equation} 
\xymatrix@C+15pt{
\C^J[{\bf W}_{J}^{-1}]\ar@/_14pt/[r]_{u^*}\ar@{}[r]|-{\perp}&\C^I[{\bf W}_{I}^{-1}]\ar@/_14pt/[l]_{v^*}}
\qquad\left(\;\text{resp.}\quad
\vcenter{
\xymatrix@C+15pt{
\C^J[{\bf W}_{J}^{-1}]\ar@/^14pt/[r]^{u^*}\ar@{}[r]|-{\perp}&\C^I[{\bf W}_{I}^{-1}]\ar@/^14pt/[l]^{v^*}}}\right)\,,
\end{equation}
dont l'unité (resp. counité) est un isomorphisme si $u$ est pleinement fidèle.

En particulier, le foncteur: 
$$
\vcenter{\xymatrix@+20pt{\C^I[{\bf W}_{I}^{-1}]\ar[r]^-{{u_{!}} \, = \, v^*}&\C^J[{\bf W}_{J}^{-1}]}}\qquad 
\left(\,\text{resp.}\qquad \vcenter{\xymatrix@+20pt{\C^I[{\bf W}_{I}^{-1}]\ar[r]^-{{u_{*}} \, = \, v^*}&\C^J[{\bf W}_{J}^{-1}]}}\right),
$$
est un foncteur des extensions de Kan homotopiques à gauche (resp. à droite) le long de $u$, lequel est pleinement fidèle si $u$ l'est.
\end{lemme}
\begin{proof}
Notons $\underline{\Cat}_{W}$ la $2$-catégorie dont les objets sont les couples $(\A,{\bf W})$, formés d'une catégorie $\A$ et d'une famille distinguée de morphismes ${\bf W}$ de $\A$. Une $1$-flèche de $(\A,{\bf W})$ vers $(\A',{\bf W}')$, est un foncteur $\xymatrix@-12pt{\A\ar[r]^F&\A'}$ vérifiant $F({\bf W})\subset {\bf W}'$; et si $F$ est $G$ sont deux $1$-flèches, une $2$-flèches de $F$ vers $G$ est une transformation naturelle entre les foncteurs. 

D'après les propriétés des catégories de fractions \cite{GZ}, on vérifie aussitôt qu'il y a un $2$-foncteur de $\underline{\Cat}_{W}$ vers la $2$-catégorie des grandes catégories $\underline{\CAT}$, qui associe au couple $(\A,{\bf W})$ la catégorie $\A[{\bf W}^{-1}]$. En composant ce $2$-foncteur avec le $2$-foncteur (contravariant en $1$-flèches) de diagrammes de $\C$, qui est défini de la $2$-catégorie des petites catégories $\underline{\cat}$ vers $\underline{\Cat}_{W}$, par la règle $I\mapsto (\C^I,{\bf W}_{I})$; on obtient un $2$-foncteur: 
\begin{equation}\label{2foncter}
\vcenter{\xymatrix@+10pt@R=1pt{
\underline{\cat}\ar[r]^-{\C^{-}[{\bf W}^{-1}_{-}]}&\underline{\CAT}\,.\\
\phantom{aaaa}I\phantom{aaaa}\ar@{|->}[r] & \phantom{aa}\C^I[{\bf W}_{I}^{-1}]\phantom{aa}}}
\end{equation}

Le résultat désiré est une conséquence formelle, des propriétés des $2$-foncteurs par rapport aux adjonctions.
\end{proof}

On déduit sans difficulté:

\begin{corollaire}
Soit $\C$ une catégorie munie d'une famille distinguée de morphismes ${\bf W}$. On considère un foncteur $u:\xymatrix@C-10pt{I\ar[r]&J}$ entre petites catégories, et une adjonction: 
$$
\vcenter{\xymatrix@+7pt{I\ar@/_10pt/[r]_{v}\ar@{}[r]|-{\perp}&I'\ar@/_10pt/[l]_{w}}}\qquad\quad
\left(\text{resp.}\quad \vcenter{\xymatrix@+7pt{I\ar@/^10pt/[r]^{v}\ar@{}[r]|-{\perp}&I'\ar@/^10pt/[l]^{w}}} \right)\,.
$$

Alors, si $\C$ admet des extensions de Kan homotopiques à gauche (resp. à droite) le long de $u$, $\C$ admet des extensions de Kan homotopiques à gauche (resp. à droite) le long du composé $uw$, et 
$$
(uw)_{!}\,\cong\, u_{!}\,\circ\,v^*\qquad\quad\big(\text{resp.}\quad 
(uw)_{*}\,\cong\, u_{*}\,\circ\,v^*
\big)
$$
\end{corollaire}

\renewcommand{\thesubsection}{\S\thesection.\arabic{subsection}}
\subsection{}\;
\renewcommand{\thesubsection}{\thesection.\arabic{subsection}}

Dans le cadre de catégories de modèles on peut montrer l'existence des extensions de Kan homotopiques à gauche (resp. droite) le long d'un foncteur quelconque, et que le système de foncteurs qu'il en découle satisfait des propriétés similaires à celles des extensions de Kan habituelles (voir \cite{cide} ou le paragraphe \S 20.2 de \cite{hirkan}). Dans le présent chapitre cependant, on va s'intéresser aux extensions de Kan homotopiques le long de foncteurs entre catégories appartenant à une famille assez restreinte où la théorie est beaucoup plus simple. 

Étant donnée une petite catégorie $I$, rappelons qu'une \emph{structure de catégorie de Reedy} sur $I$ est la donnée d'une fonction $\xymatrix@+5pt{\mathrm{Obj}(I)\ar[r]^-\lambda&\mathbb{N}}$ et de deux sous-catégories distinguées $I_{+}$ et $I_{-}$ de $I$ vérifiant les propriétés:
\begin{enumerate}
\item  Si $\xymatrix@-10pt{a\ar[r]&b}$ est un morphisme de $I_{+}$ (resp. de $I_{-}$) différent du morphisme identité, alors $\lambda(b) > \lambda(a)$ $\big($resp. $\lambda(a) >\lambda(b)$ $\big)$. 
\item Tout morphisme $u$ de $I$ s'écrit de fa\c con unique comme un composé $u=u_{+}u_{-}$, où $u_{+}$ (resp. $u_{-}$) est un morphisme de $I_{+}$ (resp. $I_{-}$).
\end{enumerate}

On dit alors que $I=(I,I_{+},I_{-},\lambda)$ est une \emph{catégorie de Reedy}, et on appelle $I_{+}$ (resp. $I_{-}$) la \emph{sous-catégorie directe} (resp. \emph{sous-catégorie inverse}) de $I$. 

Il résulte de la propriété (ii) que si $I$ est une catégorie de Reedy vérifiant $I=I_{+}$ (resp. $I=I_{-}$), alors $I_{-}$ (resp. $I_{+}$) est la sous-catégorie de $I$ de même ensemble d'objets, dont les seuls morphismes sont les identités. On appelle dans ce cas la donnée $I=(I_{+},\lambda)$ $\big($resp. $I=(I_{-},\lambda)$ $\big)$ une \emph{catégorie directe} (resp. \emph{inverse}).

\begin{proposition}\label{hocomp}
Soit $\C$ une catégorie de modèles; si $I$ est une catégorie de Reedy, alors la catégorie de diagrammes $\C^{I_{+}}$ (resp. $\C^{I_{-}}$) admet une structure de catégorie de modèles, dont les équivalences faibles et les fibration (resp. les équivalences faibles et les cofibrations) sont définies argument par argument. 

De même, la catégorie de diagrammes $\C^I$ admet une structure de catégorie de modèles dont les équivalences faibles sont définies argument par argument, et les fibration (resp. cofibrations) sont les morphismes dont les restrictions à $I_{+}$ et  $I_{-}$ sont des fibrations (resp. cofibrations) de $\C^{I_{+}}$ et $\C^{I_{-}}$ respectivement.
\end{proposition}
\begin{proof}
Voir le Théorème 16.3.4 de \cite{hirschhorn} ou le Théorème 5.2.5 de \cite{hovey}.
\end{proof}

L'énoncé suivant est facilement démontré à partir de la Proposition \ref{hocomp}.

\begin{corollaire}\label{hocompco}
Soit $\C$ une catégorie de modèles. Si $I$ est une catégorie directe (resp. inverse), et si on muni $\C^{I}$ de la structure de catégorie de modèles de la Proposition \ref{hocomp}, l'adjonction:
\begin{equation}\label{ann}
\xymatrix@C+15pt{
\C\ar@/_12pt/[r]_{p_{I}^*}\ar@{}[r]|-{\perp}&\C^{I}\ar@/_12pt/[l]_{\mathrm{colim}_{I}}}
\qquad\left(\;\text{resp.}\quad
\vcenter{
\xymatrix@C+15pt{
\C\ar@/^12pt/[r]^{p_{I}^*}\ar@{}[r]|-{\perp}&\C^{I}\ar@/^12pt/[l]^{\mathrm{lim}_{I}}}}\right)\,,
\end{equation}
est une adjonction de Quillen. En particulier, $\C$ admet des colimites homotopiques (resp. limites homotopiques) de type $I$.

Plus encore, si $u:\xymatrix@-10pt{I\ar[r]&J}$ est un foncteur entre catégories directes (resp. inverses) alors l'adjonction:
\begin{equation}\label{ann2}
\xymatrix@C+15pt{
\C^{J}\ar@/_12pt/[r]_{u^*}\ar@{}[r]|-{\perp}&\C^{I}\ar@/_12pt/[l]_{u_{!}}}
\qquad\left(\;\text{resp.}\quad
\vcenter{
\xymatrix@C+15pt{
\C^{J}\ar@/^12pt/[r]^{u^*}\ar@{}[r]|-{\perp}&\C^{I}\ar@/^12pt/[l]^{u_{*}}}}\right)\,;
\end{equation}
où $u_{!}$ (resp. $u_{*}$) est le foncteur des extensions de Kan à gauche (resp. à droite) le long de $u$ au sens habituel, est une adjonction de Quillen. En particulier, $\C$ admet des extensions de Kan homotopiques à gauche (resp. droite) le long de $u$.
\end{corollaire}
\begin{proof}
On vérifie facilement que les foncteurs $p_{I}^*$ et $u^*$ sont de foncteurs Quillen à droite (resp. à gauche), lorsque $I$ et $J$ sont de catégories directes (resp. inverses). 
\end{proof}

\section{Transformations de (co)changement de base}

\renewcommand{\thesubsection}{\S\thesection.\arabic{subsection}}
\subsection{}\;\label{conjug}    
\renewcommand{\thesubsection}{\thesection.\arabic{subsection}}

Soit $\C$ une catégorie munie d'une famille distinguée de morphismes ${\bf W}$ et supposons qu'on a une transformation naturelle  de la forme $\Gamma \colon\xymatrix@C-10pt{u\circ a \ar@{=>}[r] &b\circ v}$ c'est-à-dire:
\begin{equation}\label{Aca1}
\vcenter{\xymatrix@+8pt{
K\drtwocell<\omit>{\phantom{aa}\Gamma} \ar[d]_-{v}\ar[r]^-{a}&I\ar[d]^-{u}\\ L\ar[r]_-{b}&J
}}\,.
\end{equation}

Vu que la couple $(\C,{\bf W})$ induit un $2$-foncteur \eqref{2foncter} on en déduit un diagramme entre les catégories de fractions, des catégories de diagrammes de $\C$ correspondantes:
\begin{equation}\label{Aca2}
\vcenter{\xymatrix@+10pt{
\C^{K}[{\bf W}_{K}^{-1}]\drtwocell<\omit>{\phantom{aa}{\Gamma^*}}&\C^{I}[{\bf W}_{I}^{-1}]\ar[l]_-{{a^*}}\\ \C^L[{\bf W}_L^{-1}]\ar[u]^-{{v^*}}&\C^J[{\bf W}_{J}^{-1}]\ar[l]^-{{b^*}}\ar[u]_-{{u^*}}
}}\,.
\end{equation}

Si on suppose l'existence d'extensions de Kan homotopiques à gauche de $\C$ le long de $u$ et $v$ (resp. d'extensions de Kan homotopiques à droite de $\C$ le long de $a$ et $b$); plus précisément  si on choisit  des foncteurs adjoints comme ci-dessous:
\begin{equation}\label{Achois1}
\C^L[{\bf W}_L^{-1}]\vcenter{\xymatrix@C+10pt{
\phantom{\cdot}\ar@{}[r]|-{\perp}\ar@<-4pt>@/_10pt/[r]_-{{v^*}}&
\phantom{\cdot}\ar@<-4pt>@/_10pt/[l]_-{v_{!}}}}\C^{K}[{\bf W}_{K}^{-1}]
\qquad\text{et}\qquad
\C^J[{\bf W}_{J}^{-1}]
\vcenter{\xymatrix@+10pt{
\phantom{\cdot}\ar@{}[r]|-{\perp}\ar@<-4pt>@/_10pt/[r]_-{{u^*}}&
\phantom{\cdot}\ar@<-4pt>@/_10pt/[l]_-{u_{!}}}}
\C^{I}[{\bf W}_{I}^{-1}]
\end{equation}
$$
\left(\,\text{resp.}\qquad\;\,
\C^I[{\bf W}_I^{-1}]
\vcenter{\xymatrix@C+10pt{
\phantom{\cdot}\ar@{}[r]|-{\perp}\ar@<+4pt>@/^10pt/[r]^-{{a^*}}&
\phantom{\cdot}\ar@<+4pt>@/^10pt/[l]^-{a_{*}}}}
\C^{K}[{\bf W}_{K}^{-1}]
\qquad\text{et}\qquad
\C^J[{\bf W}_{J}^{-1}]
\vcenter{\xymatrix@+10pt{
\phantom{\cdot}\ar@{}[r]|-{\perp}\ar@<+4pt>@/^10pt/[r]^-{{b^*}}&
\phantom{\cdot}\ar@<+4pt>@/^10pt/[l]^-{b_{*}}}}
\C^{L}[{\bf W}_{L}^{-1}]
\right),
$$
et respectivement une unité et une counité (resp. une counité et une unité):
\begin{equation}\label{Achois2}
\xymatrix@-10pt{\mathrm{id}\ar@{=>}[r]^-{\eta^u}&{u^*}\circ u_{!}}\qquad\text{et}\qquad
\xymatrix@-10pt{ v_{!} \circ {v^*}\ar@{=>}[r]^-{\varepsilon^v}&\mathrm{id}}
\end{equation}
$$ 
\Big(\text{resp.}  \quad \xymatrix@-10pt{{b^{*}}\circ b_*\ar@{=>}[r]^-{\varepsilon^b}&\mathrm{id}}
\qquad\text{et}\qquad
\xymatrix@-10pt{\mathrm{id}\ar@{=>}[r]^-{\eta^a}&a_*\circ{a^{*}}}\Big)\,;
$$ 
on définit une transformation naturelle :
\begin{equation}\label{Aca3}
\vcenter{\xymatrix@+10pt{
\C^{K}[{\bf W}_{K}^{-1}]\ar[d]_-{v_{!}}&\C^{I}[{\bf W}_{I}^{-1}]\ar[d]^-{u_{!}}\ar[l]_-{{a^*}}\\ 
\C^L[{\bf W}_L^{-1}]\urtwocell<\omit>{\phantom{aaaa}{\Gamma_!}}&\C^J[{\bf W}_{J}^{-1}]\ar[l]^-{{b^*}}
}}\qquad\;\qquad
\left(\,\text{resp.}\qquad
\vcenter{\xymatrix@+10pt{
\C^{K}[{\bf W}_{K}^{-1}]   \ar[r]^-{a_{*}}  &
\C^{I}[{\bf W}_{I}^{-1}]    \dltwocell<\omit>{{\Gamma_*}\phantom{aaaa}}\\ 
\C^L[{\bf W}_L^{-1}]   \ar[r]_-{b_{*}}  \ar[u]^-{{v^*}} &\C^J[{\bf W}_{J}^{-1}]   \ar[u]_-{{u^*}}   
}}\right)\,,
\end{equation}
comme le composé:
$$
\vcenter{\xymatrix@C+40pt@R-5pt{
 v_{!}\circ {a^*}  \ar@{==>}[d]_-{{\Gamma_!}}\ar@{=>}[r]^-{(v_{!}\circ{a^*})\star \eta^u}& v_{!}\circ {a^*} \circ {u^*} \circ u_{!} \ar@{=>}[d]^-{v_{!}\star{\Gamma^*} \star u_{!}}\\
{b^*} \circ u_{!} &
v_{!}\circ{v^*} \circ {b^*} \circ u_{!} \ar@{=>}[l]^-{\epsilon^v\star( {b^*} \circ u_{!})}}}
$$
$$
\left(\,\text{resp.}\qquad
\vcenter{\xymatrix@C+40pt@R-5pt{
 {u^*} \circ b_{*} \ar@{==>}[d]_-{{\Gamma_*}}\ar@{=>}[r]^-{\eta^a\star({u^*}\circ b_{*})}&
 a_{*}\circ{a^*} \circ {u^*} \circ  b_{*} \ar@{=>}[d]^{ a_{*}\star{\Gamma^*} \star  b_{*}}\\
 a_{*}\circ{v^*}&
 a_{*}\circ{v^*} \circ {b^*} \circ  b_{*}
\ar@{=>}[l]^-{( v_{*}\circ{a^*})\star \epsilon^b}}}\right)\,.
$$

La transformation naturelle \eqref{Aca3} ainsi définie est appelée le \emph{conjuguée à gauche} (resp. \emph{à droite}) \emph{de la transformation naturelle} \eqref{Aca2} \emph{par rapport aux données} \eqref{Achois1} et \eqref{Achois2}. Cependant, dans le cas où on veut mettre l'accent sur le $2$-foncteur \eqref{2foncter}, on appelle la transformation naturelle \eqref{Aca3} \emph{un morphisme de cochangement de base} (resp, \emph{changement de base}) \emph{de} $(\C,{\bf W})$ \emph{le long de} \eqref{Aca1}. 

Remarquons que si on prend le conjuguée à gauche $\Gamma_!$ (resp. à droite $\Gamma_*$) de la transformation naturelle \eqref{Aca2} par rapport aux donnés \eqref{Achois1} et \eqref{Achois2} et après on considère la conjuguée à droite (resp. à gauche) de $\Gamma_!$ (resp. de $\Gamma_*$) par rapport aux foncteurs adjonctions \eqref{Achois1} et aux transformations naturelles:
\begin{equation}\label{Achois5}
\xymatrix@-10pt{{u^{*}}\circ u_*\ar@{=>}[r]^-{\varepsilon^u}&\mathrm{id}}
\;\,\text{et}\;\,
\xymatrix@-10pt{\mathrm{id}\ar@{=>}[r]^-{\eta^v}&v_*\circ{v^{*}}}\qquad
\Big(\text{resp.}  \quad 
\xymatrix@-10pt{\mathrm{id}\ar@{=>}[r]^-{\eta^b}&{b^*}\circ b_{!}}
\;\,\text{et}\;\,
\xymatrix@-10pt{  a_{!} \circ {a^*}\ar@{=>}[r]^-{\varepsilon^a}&\mathrm{id}}
\Big)
\end{equation}
inverses\footnote{Une unité et une counité d'une même adjonction sont dites \emph{inverses} si elles vérifient les identités triangulaires.} des transformations naturelles \eqref{Achois2},  alors on obtient à nouveau la transformation naturelle d'origine \eqref{Aca2}.

Montrons maintenant:

\begin{lemme}\label{isolemme}
Deux conjuguée à gauche (resp. à droite) ${\Gamma_!}$ et $\overline{{\Gamma_!}}$ (resp. ${\Gamma_*}$ et $\overline{{\Gamma_*}}$) de la même transformation naturelle \eqref{Aca2} sont liées par des isomorphismes de la forme:
\begin{equation}\label{isoconju}
\vcenter{\xymatrix@C+15pt@R-3pt{
 \overline{v}_{!} \circ {a^*} \ar@{=>}[r]^-{\beta\star {a^*}} _-{\cong}\ar@{=>}[d]_{\overline{{\Gamma_!}}}&
 v_{!} \circ {a^*}  \ar@{=>}[d]^-{{\Gamma_!}} \\
{b^*} \circ \overline{u}_{!}\ar@{=>}[r]_-{{b^*} \star \alpha}^-{\cong} &{b^*} \circ u_{!}}}
\qquad \quad
\left(\,\text{resp.} \qquad \vcenter{\xymatrix@C+15pt@R-3pt{
{b^*} \circ \overline{ u}_{*}\ar@{=>}[r]^-{{b^*} \star \alpha}_-{\cong} \ar@{=>}[d]_-{\overline{{\Gamma_*}}}&{b^*} \circ u_{*}\ar@{=>}[d]^-{{\Gamma_*}}\\
 \overline{ v}_{*} \circ {a^*} \ar@{=>}[r]_-{\beta\star {a^*}}^-{\cong} &
 v_{!} \circ {a^*}}}\right)\,,
\end{equation}
\end{lemme}
\begin{proof}
En effet, on sait que si on considère d'autres foncteurs adjoints à gauche (resp. droite):
\begin{equation}\label{Achois3}
\C^L[{\bf W}_L^{-1}]\vcenter{\xymatrix@C+10pt{
\phantom{\cdot}\ar@{}[r]|-{\perp}\ar@<-4pt>@/_10pt/[r]_-{{v^*}}&
\phantom{\cdot}\ar@<-4pt>@/_10pt/[l]_-{ \overline{v}_{!} } }}\C^{K}[{\bf W}_{K}^{-1}]
\qquad\text{et}\qquad
\C^J[{\bf W}_{J}^{-1}]
\vcenter{\xymatrix@+10pt{
\phantom{\cdot}\ar@{}[r]|-{\perp}\ar@<-4pt>@/_10pt/[r]_-{{u^*}}&
\phantom{\cdot}\ar@<-4pt>@/_10pt/[l]_-{\overline{u}_{!}}}}
\C^{I}[{\bf W}_{I}^{-1}]
\end{equation}
$$
\left(\,\text{resp.}\qquad\;\,
\C^L[{\bf W}_L^{-1}]
\vcenter{\xymatrix@C+10pt{
\phantom{\cdot}\ar@{}[r]|-{\perp}\ar@<+4pt>@/^10pt/[r]^-{{v^*}}&
\phantom{\cdot}\ar@<+4pt>@/^10pt/[l]^-{\overline{v}_{*}}}}
\C^{K}[{\bf W}_{K}^{-1}]
\qquad\text{et}\qquad
\C^J[{\bf W}_{J}^{-1}]
\vcenter{\xymatrix@+10pt{
\phantom{\cdot}\ar@{}[r]|-{\perp}\ar@<+4pt>@/^10pt/[r]^-{{u^*}}&
\phantom{\cdot}\ar@<+4pt>@/^10pt/[l]^-{\overline{u}_{*}}}}
\C^{I}[{\bf W}_{I}^{-1}]
\right),
$$
avec de respectives transformations naturels:
\begin{equation}\label{chois4}
\xymatrix@-10pt{\mathrm{id}\ar@{=>}[r]^-{\overline{\eta}}&{u^*}\circ\overline{u}_{!}}\qquad\text{et}\qquad
\xymatrix@-10pt{\overline{v}_!\circ{v^{*}}\ar@{=>}[r]^-{\overline{\varepsilon}}&\mathrm{id}}
\end{equation}
$$ 
\Big(\text{resp.}  \quad \xymatrix@-10pt{{u{*}}\circ\overline{u}_*\ar@{=>}[r]^-{\overline{\varepsilon}}&\mathrm{id}}
\qquad\text{et}\qquad
\xymatrix@-10pt{\mathrm{id}\ar@{=>}[r]^-{\overline{\eta}}&\overline{v}_{!}\circ{v^*}}\Big)\,;
$$ 
ils existent des isomorphismes naturels:
$$
\xymatrix@C-10pt{\overline{u}_{!}\ar@{=>}[r]^{\alpha}&u_{!} }  \quad \text{et} \quad 
\xymatrix@C-10pt{\overline{v}_{!}\ar@{=>}[r]^{\beta}&v_{!} } \quad\quad
\Big(\text{resp.}\quad
\xymatrix@C-10pt{\overline{u}_{*}\ar@{=>}[r]^{\alpha}&u_{*} } \quad \text{et} \quad 
\xymatrix@C-10pt{\overline{v}_{*}\ar@{=>}[r]^{\beta}&v_{*} }\Big)\,,
$$
tels que les diagrammes suivants commutent:
$$
\vcenter{\xymatrix@R-20pt@C+5pt{  & {u^*} \circ u_{!} \\
\mathrm{id} \ar@{=>}[ru]^-{\eta} \ar@{=>}[rd]_-{\overline{\eta}} & \\
&{u^*} \circ \overline{u}_{!} \ar@{=>}[uu]_{{u^*}\star \alpha} }}\qquad\text{et} \qquad 
\vcenter{\xymatrix@R-20pt@C+5pt{v_!\circ{v_{*}} \ar@{=>}[rd]^-{\epsilon}& \\& \mathrm{id}\\
\overline{v}_!\circ {v_{*}} \ar@{=>}[ru]_-{\overline{\epsilon}} \ar@{=>}[uu]^-{ \beta\star {v^*}} }}
$$
$$
\left(\,\text{resp.}\qquad \quad \vcenter{\xymatrix@R-20pt@C+5pt{  {u^*} \circ u_{*} \ar@{=>}[rd]^-{\epsilon}&\\
& \mathrm{id}  \\
{u^*} \circ \overline{u}_{*} \ar@{=>}[uu]^{{u^*}\star \alpha} \ar@{=>}[ru]_-{\overline{\epsilon}}}} \qquad \text{et} \qquad 
\vcenter{\xymatrix@R-20pt@C+5pt{ & v_{!}\circ{v^*}\\
  \mathrm{id}   \ar@{=>}[ru]^-{\eta}\ar@{=>}[rd]_-{\overline{\eta}}&\\
&\overline{v}_{!}\circ{v^*} \ar@{=>}[uu]_-{ \beta\star {v^*}} }}\right)\,.
$$

On montre alors que le carré \eqref{isoconju} est bien commutatif vu qu'il s'insère dans le cube suivant:
$$
\xymatrix@C+15pt@R-3pt{
\overline{v}_{!} \circ {a^*}  \ar@{=>}[dd]\ar@{=>}[rr]|-{\phantom{aa}( \overline{v}_{!} \circ {a^*})\star \overline{\eta}\phantom{aa}} 
\ar@{=>}[rd]|-{\beta\star {a^*}} & &   
\overline{v}_{!}\circ {a^*} \circ {u^*} \circ \overline{u}_{!} \ar@{=>}[rd]|-{\beta\star({a^*} \circ {u^*})\star \alpha} 
\ar@{=>}'[d]|(.80){ \underset{\phantom{aa}}{\underset{\phantom{aa}}{\overset{\phantom{aa}}{\overline{v}_{!}\star {\Gamma^*} \star \overline{u}_{!}}} }}[dd]    & \\
&  v_{!} \circ {a^*} \ar@{=>}[rr]|(.30){\phantom{aa}(v_{!} \circ {a^*})\star \eta\phantom{aa}} \ar@{=>}[dd] && 
v_{!}\circ {a^*} \circ {u^*} \circ u_{!} 
\ar@{=>}[dd]|-{\underset{\phantom{aa}}{\overset{\phantom{aa}}{v_{!}\star {\Gamma^*} \star u_{!}}} }  
\\{b^*} \circ \overline{u}_{!}\ar@{=>}[rd]|-{{b^*} \star \alpha}&& 
\overline{v}_{!}\circ{v^*} \circ {b^*} \circ \overline{u}_{!} \ar@{=>}[rd]|-{\beta\star({v^*} \circ {b^*})\star \alpha} 
\ar@{=>}'[l]|(.70){\phantom{.}\overline{\epsilon}\star ( {b^*} \circ \overline{u}_{!})\phantom{.}}[ll]&\\
&{b^*} \circ u_{!}&& v_{!}\circ{v^*} \circ {b^*} \circ u_{!} 
\ar@{=>}[ll]|-{\epsilon\star ({b^*} \circ u_{!})} 
}
$$
dont toutes les autres faces sont effectivement des carrés commutatifs.
\end{proof}

Voyons aussi:

\begin{lemme}\label{matesconjureal}
Si la transformation naturelle \eqref{Aca2} admet une transformation naturelle conjuguée à gauche $\xymatrix{ v_{!} \circ{a^*}\ar@{=>}[r]^-{{\Gamma_!}} &{b^*} \circ u_{!} }$  et une conjuguée à droite $\xymatrix{{u^*} \circ b_{*}  \ar@{=>}[r]^-{{\Gamma_*}} & a_{*}\circ{v^*} }$, ce deux transformations naturelles sont conjuguée l'une de l'autre dans le sens IV.7 de \cite{lane} ou \cite{mate}. 

Explicitement on a un carré commutatif:
\begin{equation} \label{conjuMAC}
\xymatrix@C=2pt{
\mathrm{Hom}_{\C^{L}[{\bf W}_{L}^{-1}]} \big( v_{!} \circ{a^*} (\mathfrak{X}), \mathfrak{Y} \big)& \cong & 
\mathrm{Hom}_{\C^{I}[{\bf W}_{I}^{-1}]} \big(\mathfrak{X},  a_{*} \circ{v^*} (\mathfrak{Y}) \big) \\
\mathrm{Hom}_{\C^{L}[{\bf W}_{L}^{-1}]} \big( {b^*} \circ u_{!} (\mathfrak{X}),\mathfrak{Y}\big) \ar[u]^-{-\,\circ({\Gamma_!})_\mathfrak{X}}& \cong & 
\mathrm{Hom}_{\C^{I}[{\bf W}_{I}^{-1}]} \big({X},{u^*} \circ b_{*}( \mathfrak{Y}) \big)\ar[u]_-{({\Gamma_*})_\mathfrak{Y}\circ\,-} }
\end{equation}
pour tous les $I$-diagrammes $\mathfrak{X}$ et les $L$-diagrammes $\mathfrak{Y}$, où les isomorphismes horizontaux sont définis à partir des transformations naturelles \eqref{Achois2} qui définissent $\Gamma_!$ et $\Gamma_*$ plus ces inverses \eqref{Achois5}.

En particulier la transformation $\xymatrix{ v_{!} \circ{a^*}\ar@{=>}[r]^-{{\Gamma_!}} &{b^*} \circ u_{!} }$ est un isomorphisme si et seulement si la transformation $\xymatrix{{u^*} \circ b_{*}  \ar@{=>}[r]^-{{\Gamma_*}} & a_{*}\circ{v^*} }$ est un isomorphisme.
\end{lemme}
\begin{proof}
D'après le Théorème 2 de IV.7 dans \cite{lane}, il suffit de voir que le carré ci-dessous est commutatif:
\begin{equation}\label{carrre}
\xymatrix@C+30pt{
\mathrm{id}\ar@{=>}[r]^-{\eta^2}\ar@{=>}[d]_-{\eta^1} & 
(a_{*}\circ {v^*})\circ(v_{!}\circ {a^*}) \ar@{=>}[d]^-{(a_{*}\circ {v^*})\star{\Gamma_!}}\\
({u^*}\circ b_{*})\circ({b^*}\circ u_{!}) \ar@{=>}[r]_-{{\Gamma_*}\star({b^*}\circ u_{!}) }&   
( a_{*}\circ {v^*})\circ({b^*}\circ u_{!})}
\end{equation}
où $\eta^1$ et $\eta^2$ sont respectivement les transformations naturelles suivantes:
$$
\xymatrix@C+5pt{\mathrm{id} \ar@{=>}[r]^-{\eta^u}&{u^*}\circ u_{!}\ar@{=>}[rr]^-{{u^*}\star\eta^b\star u_{!}}&&
{u^*}\circ (  b_{*}\circ{b^*} ) \circ u_{!} = ({u^*}\circ  b_{*})\circ({b^*}  \circ u_{!} )} 
$$
$$\text{et}$$

$$
\xymatrix@C+5pt{\mathrm{id} \ar@{=>}[r]^-{\eta^a}&  a_{*}\circ {a^*}\ar@{=>}[rr]^-{ u_{*} \star \eta^v \star {a^*}}&& 
a_{*}\circ ({u^*}\circ u_{!})\circ {a^*} = (  a_{*}\circ {u^*})\circ ( u_{!}\circ {a^*})  }
$$

Pour montrer cela écrivons le deux chemin possibles $\big[{\Gamma_*}\star({b^*}\circ u_{!})\big] \circ \eta^1$  et $\big[( a_{*}\circ {v^*})\star{\Gamma_!}\big]\circ\eta^2$  du diagramme \eqref{carrre} comme suit:
$$
\vcenter{\xymatrix@+10pt{
&\C^{I}[{\bf W}_{I}^{-1}] \dtwocell<\omit>{<8>\phantom{aa}\eta^a} &\\
\C^{K}[{\bf W}_{K}^{-1}]\ar@/^15pt/[ru]^-{a_*}\drtwocell<\omit>{\phantom{aa}{\Gamma}}
&\C^{I}[{\bf W}_{I}^{-1}]\ar[l]|-{{a^*}}\ar[u]_-{\mathrm{id}}  \rtwocell<\omit>{<5>\phantom{aa}\eta^u} &
\C^{I}[{\bf W}_{I}^{-1}]\ar[l]_-{\mathrm{id}}\ar[d]^-{u_{!}} \\ 
\C^L[{\bf W}_L^{-1}]\ar[u]^-{{v^*}}\drtwocell<\omit>{<-2>\phantom{a}\varepsilon^b}&
\C^J[{\bf W}_{J}^{-1}]\ar[l]|-{{b^*}}\ar[u]|-{{u^*}}\rtwocell<\omit>{<5>\phantom{aa}\eta^b}  &
\C^{J}[{\bf W}_{J}^{-1}]\ar[l]|-{\mathrm{id}}   \ar@/^15pt/[dl]^-{{b*}}\\
& \C^L[{\bf W}_L^{-1}] \ar@/^15pt/[lu]^-{\mathrm{id}} \ar[u]|-{b_*}&
}}
$$
$$\text{et} \qquad \quad
\vcenter{\xymatrix@R+23pt{
&\C^{K}[{\bf W}_{K}^{-1}]  \ar[r]^-{a_*}  \dtwocell<\omit>{<6>\phantom{aa}\eta^v}  \dtwocell<\omit>{<-8>\phantom{aa}\eta^a}   & 
\C^{I}[{\bf W}_{I}^{-1}]\\
\C^L[{\bf W}_L^{-1}]\ar@/^15pt/[ru]^-{{v^*}}\drtwocell<\omit>{\phantom{a}\varepsilon^v}&
\C^{K}[{\bf W}_{K}^{-1}]\ar[l]|-{ v_!}\ar[u]|-{\mathrm{id}}\drtwocell<\omit>{\phantom{aa}{\Gamma}}&
\C^{I}[{\bf W}_{I}^{-1}]\ar[l]|-{{a^*}}\ar[u]_-{\mathrm{id}}\rtwocell<\omit>{<5>\phantom{a}\eta^u}&
    \C^{I}[{\bf W}_{I}^{-1}]    \ar[l]_-{\mathrm{id}} \\ 
&\C^L[{\bf W}_L^{-1}]\ar[u]|-{{v^*}}\ar@/^15pt/[lu]^-{\mathrm{id}}&
\C^J[{\bf W}_{J}^{-1}]\ar[l]^-{{b^*}}\ar[u]|-{{u^*}}\ar@{<-}@/_15pt/[ru]_-{ u_!}&
}}\,,
$$
respectivement. 

On déduit alors des dites identités triangulaires que ceux deux chemins sont égaux au composé:
$$
\xymatrix@C+20pt{\mathrm{id} \ar@{=>}[r]^-{\eta^a\star\eta^u}  & a_{*}\circ{a^*} \circ{u^*}\circ u_{!}
\ar@{=>}[r]^-{ a_{*}\star\Gamma\star  u_{!}}&
 a_{*}\circ{v^*} \circ{b^*}\circ u_{!}}\,,
$$
donc \eqref{carrre} est bien un carré commutatif.
\end{proof}

\renewcommand{\thesubsection}{\S\thesection.\arabic{subsection}}
\subsubsection{}\;\label{evalus}     
\renewcommand{\thesubsection}{\thesection.\arabic{subsection}}

Si $\xymatrix@-10pt{I\ar[r]^-u&J}$ est un foncteur entre petites catégories, rappelons que la catégorie  $I\,|\,b$ (resp. $b\,|\,I$) des \emph{objets de $I$ au-dessus} (resp. \emph{au-dessous}) \emph{d'un objet $b$ de $J$}, est par définition la catégorie dont les objets sont les couples $(a,\alpha)$, où $a$ est un objet de $I$ et $\alpha:\xymatrix@-10pt{ua\ar[r]&b}$ $\big($resp. $\alpha:\xymatrix@-10pt{b\ar[r]&ua}$ $\big)$ est une flèche de $J$. Un morphisme de $(a,\alpha)$ vers $(a',\alpha')$ est une flèche $f:\xymatrix@-10pt{a\ar[r]&a'}$ de $I$ telle que $\alpha' u(f)=\alpha$ (resp. $u(f)\alpha=\alpha'$).

On vérifie en particulier qu'il y a une transformation naturelle canonique: 
\begin{equation}\label{ca1}
\vcenter{\xymatrix@+8pt{
I\,|\,b\drtwocell<\omit>{\phantom{aa}\Gamma\,|\,b} \ar[d]_-{p_{I\,|\,b}}\ar[r]^-{\pi\,|\,b}&I\ar[d]^-{u}\\ e\ar[r]_-{b}&J
}}\qquad\;\qquad
\left(\,\text{resp.}\qquad
\vcenter{\xymatrix@+8pt{
b\,|\,I\drtwocell<\omit>{\phantom{aa}b\,|\,\Gamma}\ar[r]^-{p_{b\,|\,I}}\ar[d]_-{b\,|\,\pi}&e\ar[d]^-{b}\\ I\ar[r]_-{u}&J
}}\right),
\end{equation}
définie pour tout objet $(a,\alpha)$ de $I\,|\,b$ (resp. $b\,|\,I$) par la formule:
$$
(\Gamma\,|\,b)_{(a,\alpha)}=\alpha \qquad \quad \big(\text{resp.} \quad (b\,|\,\Gamma)_{(a,\alpha)}=\alpha \big).
$$

Si $\C$ est une catégorie munie d'une famille distinguée de morphismes ${\bf W}$, on déduit de \eqref{ca1} une transformation naturelle:
\begin{equation}\label{ca2}
\vcenter{\xymatrix@+10pt{
\C^{I\,|\,b}[{\bf W}_{I\,|\,b}^{-1}]\drtwocell<\omit>{\phantom{aa}{\Gamma\,|\,b}}&\C^{I}[{\bf W}_{I}^{-1}]\ar[l]_-{{(\pi\,|\,b)^*}}\\ \C[{\bf W}^{-1}]\ar[u]^-{{(p_{I\,|\,b})^*}}&\C^J[{\bf W}_{J}^{-1}]\ar[l]^-{{b^*}}\ar[u]_-{{u^*}}
}}\qquad\;\qquad
\left(\,\text{resp.}\qquad
\vcenter{\xymatrix@+10pt{
\C^{b\,|\,I}[{\bf W}_{b\,|\,I}^{-1}]\drtwocell<\omit>{\phantom{aa}{b\,|\,\Gamma}}&\C[{\bf W}^{-1}]\ar[l]_-{{(p_{b\,|\,I})^*}}\\ 
\C^I[{\bf W}^{-1}_I]\ar[u]^-{{(b\,|\,\pi)^*}}&\C^J[{\bf W}_{J}^{-1}]\ar[l]^-{{u^*}}\ar[u]_-{{b^*}}
}}\right)\,.
\end{equation}

Une conjuguée à gauche de la transformation naturelle ${\Gamma\,|\,b}$ (resp. une conjuguée á droite de ${b\,|\,\Gamma}$):
\begin{equation}\label{ca3}
\vcenter{\xymatrix@+10pt{
\C^{I\,|\,b}[{\bf W}_{I\,|\,b}^{-1}]\ar[d]_-{\mathrm{hocolim}_{I\,|\,b}}&\C^{I}[{\bf W}_{I}^{-1}]\ar[d]^-{u_{!}}\ar[l]_-{{(\pi\,|\,b)^*}}\\ 
\C[{\bf W}^{-1}]\urtwocell<\omit>{\phantom{aaaa}{(\Gamma\,|\,b)_!}}&\C^J[{\bf W}_{J}^{-1}]\ar[l]^-{{b^*}}
}}\qquad\;\qquad
\left(\,\text{resp.}\qquad
\vcenter{\xymatrix@+10pt{
\C^{b\,|\,I}[{\bf W}_{I\,|\,b}^{-1}]\ar[r]^-{\mathrm{holim}_{b\,|\,I}}& \C[{\bf W}^{-1}]\dltwocell<\omit>{{(b\,|\,\Gamma)_*}\phantom{aaaa}} \\ 
\C^I[{\bf W}^{-1}_I]\ar[u]^-{{(b\,|\,\pi)^*}}\ar[r]_-{u_*}& \C^J[{\bf W}_{J}^{-1}] \ar[u]_-{b^{*}}
}}\right)\,;
\end{equation}
est appelé un \emph{morphisme d'évaluation} de $u_{!}$ $\big($resp. $u_{*} \big)$ \emph{en $b$}. 

Si $\mathfrak{X}$ est un $I$-diagramme de $\C$, on appelle le morphisme de la catégorie $ \C[{\bf W}^{-1}]$:
\begin{equation}\label{evalu}
\vcenter{\xymatrix@+15pt{\mathrm{hocolim}_{I\,|\,b} \big(\mathfrak{X}\circ(\pi\,|\,b) \big)\ar[r]^-{\big( ({b\,|\,\Gamma})_!\big)_{\mathfrak{X}}}&u_{!} (\mathfrak{X})_{b}}}
\end{equation}
$$
\left(\,\text{resp.}\qquad\;\vcenter{\xymatrix@+15pt{
u_{*} (\mathfrak{X})_{b}\ar[r]^-{\big( ({b\,|\,\Gamma})_*\big)_{\mathfrak{X}}}&
\mathrm{holim}_{b\,|\,I} \big(\mathfrak{X}\circ(b\,|\,\pi) \big)}} \right),
$$
un \emph{morphisme d'évaluation de} $u_{!}(\mathfrak{X})$ $\big($resp. $u_{*}(\mathfrak{X})$ $\big)$ \emph{en $b$}.

\begin{proposition}\label{local} 
Si $\C$ est une catégorie de modèles et $u:\xymatrix@C-10pt{I\ar[r]&J}$ est un foncteur entre petites catégories directes (resp. inverses), alors pour tout objet $b$ de $J$ et tout $I$-diagramme $\mathfrak{X}$ de $\C$ le morphisme d'évaluation \eqref{evalu} de $u_{!}(\mathfrak{X})$ $\big($resp. $u_{*}(\mathfrak{X})\big)$ en $b$ est un isomorphisme (voir \ref{evalus}). Autrement dit la transformation naturelle \eqref{ca3} est un isomorphisme (voir le Lemme \ref{reflee}).
\end{proposition} 
\begin{proof}
Pour commencer rappelons que si $I$ est une catégorie directe (resp. indirecte); alors la catégorie $I\,|\,b$ (resp. $b\,|\,I$) est munie d'une structure canonique de catégorie directe (resp. indirecte).

Ensuite, remarquons qu'il suffit de montrer le résultat pour $\mathfrak{X}$ fibrant et cofibrant dans $\C^I$, lorsque $\C^{I}$ est munie de la structure de catégorie de modèles de la Proposition \ref{hocomp} (voir aussi le Corollaire \ref{hocompco}). 

Finalement, on vérifie que dans ce cas le morphisme d'évaluation \eqref{evalu} est l'isomorphisme:
$$
\vcenter{\xymatrix@+15pt{\mathrm{colim}_{I\,|\,b} \big((\pi\,|\,b)^* \, \mathfrak{X} \big)\ar[r]^-{(\Gamma\,|\,b)_{\mathfrak{X}}}&u_{!} (\mathfrak{X})_{b}}}
$$
$$
\left(\,\text{resp.}\qquad\;\vcenter{\xymatrix@+15pt{
u_{*} (\mathfrak{X})_{b}\ar[r]^-{(b\,|\,\Gamma)_{\mathfrak{X}}}&
\mathrm{lim}_{b\,|\,I} \big((b\,|\,\pi)^* \, \mathfrak{X} \big)}} \right),
$$
qui donne les formules des extensions de Kan ordinaires de \cite{lane}.
\end{proof}

\renewcommand{\thesubsection}{\S\thesection.\arabic{subsection}}
\subsection{}\;    
\renewcommand{\thesubsection}{\thesection.\arabic{subsection}}

Supposons maintenant qu'on a deux transformations naturelles avec un côté en commun:
$$
\vcenter{\xymatrix@+8pt{
K\drtwocell<\omit>{\phantom{aa}\Gamma^1} \ar[d]_-{v}\ar[r]^-{a}&I\ar[d]^-{u}\\ L\ar[r]_-{b}&J
}}\qquad\text{et}\qquad
\vcenter{\xymatrix@+8pt{
P\drtwocell<\omit>{\phantom{aa}\Gamma^2} \ar[d]_-{w}\ar[r]^-{c}&K\ar[d]^-{v}\\ Q\ar[r]_-{d}&L}}
\qquad\qquad\left(\,\text{resp.}\qquad \quad 
\vcenter{\xymatrix@+8pt{ 
P\drtwocell<\omit>{\phantom{aa}\Gamma^2} 
\ar[d]_-{k}\ar[r]^-{c}&
Q\ar[d]^-{w}\\ K\ar[r]_-{a}&I}}\right)\,
$$
et considérons la transformation naturelle qu'on obtient en collant les carrés comme suit:
$$
\vcenter{\xymatrix@C+18pt@R+6pt{
P\drtwocell<\omit>{\phantom{aa}\Gamma^2} \ar[d]_-{w}\ar[r]^-{c}& 
K\drtwocell<\omit>{\phantom{aa}\Gamma^1} \ar[d]|-{\overset{\phantom{a}}{\underset{\phantom{a}}{v}}} \ar[r]^-{a} & 
I\ar[d]^-{u}\\
Q\ar[r]_-{d}&L \ar[r]_-{b}& J}}\quad = \quad
\vcenter{\xymatrix@C+15pt@R+2pt{
P\drtwocell<\omit>{\phantom{aaa}\Gamma^3} \ar[d]_-{w}\ar[r]^-{a\circ c}&I\ar[d]^-{u}\\ Q\ar[r]_-{b\circ d}&J
}}
$$
$$
\left(\,\text{resp.}\quad \qquad
\vcenter{\xymatrix@C+8pt@R+5pt{ 
P \drtwocell<\omit>{\phantom{aa}\Gamma^2}\ar[d]_-{k}\ar[r]^-{c}&Q\ar[d]^-{w}\\
K\drtwocell<\omit>{\phantom{aa}\Gamma^1} \ar[d]_-{v}\ar[r]|-{a}&I\ar[d]^-{u}\\
L\ar[r]_-{b}&J}}\quad = \quad
\vcenter{\xymatrix@C+15pt@R+5pt{ 
P \drtwocell<\omit>{\phantom{aaa}\Gamma^3}\ar[d]_-{v\circ k}\ar[r]^-{c}&Q\ar[d]^-{u\circ w}\\
L\ar[r]_-{b}&J
}}\right)\,
$$
c'est -à-dire, on définit $\Gamma^3= (b\star \Gamma^2)\circ (\Gamma^1\star c)\quad\left(\text{resp.}\quad \Gamma^3 = (\Gamma^1\star k) \circ (u\star \Gamma^2) \right)$.

Si $\C$ est une catégorie munie d'une famille distinguée de morphismes ${\bf W}$, en appliquant le $2$-foncteur $\C^{-}[{\bf W}_{-}^{-1}]$ on obtient les diagrammes suivants entre les catégories de fractions:
\begin{equation}\label{gamagama}
\vcenter{\xymatrix@C+8pt@R+6pt{
\C^{P}[{\bf W}_{P}^{-1}]\drtwocell<\omit>{\phantom{aa}\Gamma^2} & 
\C^{K}[{\bf W}_{K}^{-1}]\drtwocell<\omit>{\phantom{aa}\Gamma^1} \ar[l]_-{c^*}& 
\C^{I}[{\bf W}_{I}^{-1}] \ar[l]_-{a^*}\\
\C^{Q}[{\bf W}_{Q}^{-1}]\ar[u]^-{w^*}&\C^{L}[{\bf W}_{L}^{-1}] \ar[u]|-{v^*} \ar[l]^-{d^*}& \C^{J}[{\bf W}_{J}^{-1}]  \ar[u]_-{u^*} \ar[l]^-{b^*} }}\; = \;
\vcenter{\xymatrix@C+15pt@R+2pt{
\C^{P}[{\bf W}_{P}^{-1}]\drtwocell<\omit>{\phantom{aa}\Gamma^3} & \C^{I}[{\bf W}_{I}^{-1}]\ar[l]_-{c^*\circ a^*}\\
\C^{Q}[{\bf W}_{Q}^{-1}]\ar[u]^-{w^*}&\C^{J}[{\bf W}_{J}^{-1}]\ar[u]_-{u^*}\ar[l]^-{d^*\circ b^*}
}}
\end{equation}
$$
\left(\,\text{resp.}\qquad 
\vcenter{\xymatrix@C+8pt@R+5pt{ 
\C^{P}[{\bf W}_{P}^{-1}]\drtwocell<\omit>{\phantom{aa}\Gamma^2} &\C^{Q}[{\bf W}_{Q}^{-1}] \ar[l]_-{c^*} \\ 
\C^{K}[{\bf W}_{K}^{-1}]\ar[u]^-{k^*}\drtwocell<\omit>{\phantom{aa}\Gamma^1} &\C^{I}[{\bf W}_{I}^{-1}] \ar[u]_-{w^*}\ar[l]|-{a^*}\\  
\C^{L}[{\bf W}_{L}^{-1}]\ar[u]^-{v^*}&\C^{J}[{\bf W}_{J}^{-1}] \ar[u]_-{u^*}\ar[l]^-{b^*}\\}}\;=\;
\vcenter{\xymatrix@C+18pt@R+5pt{ 
\C^{P}[{\bf W}_{P}^{-1}]  \drtwocell<\omit>{\phantom{aa}\Gamma^3}&\C^{Q}[{\bf W}_{Q}^{-1}]\ar[l]_-{c^*}\\
\C^{L}[{\bf W}_{L}^{-1}]\ar[u]^-{k^*\circ v^*}&\C^{J}[{\bf W}_{J}^{-1}]\ar[l]^-{b^*}\ar[u]_-{w^*\circ u^*}
}}\right)\,
$$

Étant donnés des foncteurs adjoints à gauche des foncteurs $u^*$, $v^*$ et $w^*$ (resp. des adjoints à droite de $a^*$, $b^*$ et $c^*$):
\begin{equation}\label{elecc}
u_!   \dashv  u^* \, , \quad v_!   \dashv  v^*\quad \text{et}\quad w_!  \dashv  w^* \qquad \quad
\left( \text{resp.}\qquad  a^*   \dashv  a_*\,,\quad b^*  \dashv  b_*\quad \text{et} \quad c^*  \dashv  c_* \right)\,,
\end{equation}
accompagnés des unités et counités suivantes:
\begin{equation} \label{eleccc}
\xymatrix@-10pt{\mathrm{id}\ar@{=>}[r]^-{\eta^u}&{u^*}\circ u_{!}}\,,\quad 
\xymatrix@-10pt{\mathrm{id}\ar@{=>}[r]^-{\eta^v}&{v^*}\circ v_{!}}\,,\quad
\xymatrix@-10pt{ v_{!} \circ {v^*}\ar@{=>}[r]^-{\varepsilon^v}&\mathrm{id}}\qquad\text{et}\qquad
\xymatrix@-10pt{ w_{!} \circ {w^*}\ar@{=>}[r]^-{\varepsilon^w}&\mathrm{id}}
\end{equation}
$$ 
\Big(\text{resp.}  \quad 
\xymatrix@-10pt{{b^{*}}\circ b_*\ar@{=>}[r]^-{\varepsilon^b}&\mathrm{id}}\,,\quad 
\xymatrix@-10pt{{a^{*}}\circ a_*\ar@{=>}[r]^-{\varepsilon^a}&\mathrm{id}}\,,\quad 
\xymatrix@-10pt{\mathrm{id}\ar@{=>}[r]^-{\eta^a}&{a_*}\circ a^*}
\qquad\text{et}\qquad
\xymatrix@-10pt{\mathrm{id}\ar@{=>}[r]^-{\eta^c}&c^*\circ{c_{*}}}\Big)\,;
$$ 
si on pose $\Gamma^1_!$, $\Gamma^2_!$ et $\Gamma^3_!$ $\left(\text{resp.}\quad\Gamma^1_*,\;\Gamma^2_*\;\text{et}\;\Gamma^3_*\right)$ pour noter les transformations naturelles conjuguées à gauche (resp. à droite) des transformations \eqref{gamagama} par rapport aux donnés \eqref{elecc} et \eqref{eleccc}, on vérifie aussi-tôt que:

\begin{lemme}\label{paste}
On a une égalité:
$$
\vcenter{\xymatrix@C+8pt@R+6pt{
\C^{P}[{\bf W}_{P}^{-1}] \ar[d]_-{w_!}& 
\C^{K}[{\bf W}_{K}^{-1}]\ar[l]_-{c^*}\ar[d]|-{v_!}& 
\C^{I}[{\bf W}_{I}^{-1}] \ar[l]_-{a^*}\ar[d]^-{u_!}\\
\C^{Q}[{\bf W}_{Q}^{-1}]\urtwocell<\omit>{\phantom{aa}\Gamma^2_!}&
\C^{L}[{\bf W}_{L}^{-1}] \urtwocell<\omit>{\phantom{aa}\Gamma^1_!}  \ar[l]^-{d^*}&
\C^{J}[{\bf W}_{J}^{-1}]   \ar[l]^-{b^*} }}\; = \;
\vcenter{\xymatrix@C+15pt@R+2pt{
\C^{P}[{\bf W}_{P}^{-1}]\ar[d]_-{w_!}& \C^{I}[{\bf W}_{I}^{-1}]\ar[l]_-{c^*\circ a^*}\ar[d]^-{u_!}\\
\C^{Q}[{\bf W}_{Q}^{-1}]\urtwocell<\omit>{\phantom{aaa}\Gamma^3_!} &\C^{J}[{\bf W}_{J}^{-1}]\ar[l]^-{d^*\circ b^*}
}}
$$
$$
\left(\,\text{resp.}\qquad 
\vcenter{\xymatrix@C+8pt@R+5pt{ 
\C^{P}[{\bf W}_{P}^{-1}]\ar[r]^-{c_*}&\C^{Q}[{\bf W}_{Q}^{-1}] \dltwocell<\omit>{\Gamma^2_*\phantom{aa}} \\
\C^{K}[{\bf W}_{K}^{-1}]  \ar[u]^-{k^*}\ar[r]|-{a_*} &\C^{I}[{\bf W}_{I}^{-1}]\ar[u]_-{w^*} \dltwocell<\omit>{\Gamma^1_*\phantom{aa}} \\ 
\C^{L}[{\bf W}_{L}^{-1}] \ar[u]^-{v^*}\ar[r]_-{b_*} &\C^{J}[{\bf W}_{J}^{-1}] \ar[u]_-{u^*} }}\;=\;
\vcenter{\xymatrix@C+18pt@R+5pt{ 
\C^{P}[{\bf W}_{P}^{-1}]\ar[r]^-{c_*}&\C^{Q}[{\bf W}_{Q}^{-1}]\dltwocell<\omit>{\Gamma^3_*\phantom{aa}}\\
\C^{L}[{\bf W}_{L}^{-1}]\ar[r]_-{b_*}\ar[u]^-{k^*\circ v^*}&\C^{J}[{\bf W}_{J}^{-1}]\ar[u]_-{w^*\circ u^*}
}}\right)\,
$$
c'est-à-dire
$$
\Gamma_!^3 = (d^*\star\Gamma^1_!)\circ (\Gamma^2_!\star a^*) \qquad \quad \left(\text{resp.} \quad \Gamma_*^3 = (\Gamma^2_!\star v^*)\circ (w^*\star\Gamma^1_*)\right)\,,
$$
toujours que les transformations naturelles $\eta^v$ et $\varepsilon^v$ (resp. $\eta^a$ et $\varepsilon^a$) soient inverses l'une de l'autre, c'est-à-dire s'ils vérifient les identités triangulaires. 
\end{lemme}

\renewcommand{\thesubsection}{\S\thesection.\arabic{subsection}}
\subsubsection{}\;     
\renewcommand{\thesubsection}{\thesection.\arabic{subsection}}

Considérons un foncteur pleinement fidèle $u:\xymatrix@C-10pt{I\ar[r]&J}$ et remarquons que pour tout objet $a$ de $I$ on a une égalité de transformations naturelles:
\begin{equation}\label{upf}
\vcenter{\xymatrix@C+18pt@R+6pt{
I\,|\,ua\drtwocell<\omit>{\phantom{aa}\Gamma} \ar[d]_-{p_{I\,|\,ua}}\ar[r]^-{\pi\,|\,ua}&
I\drtwocell<\omit>{\phantom{aa}\mathrm{id}} \ar[d]|-{\mathrm{id}} \ar[r]^-{\mathrm{id}} & I\ar[d]^-{u}\\
 e\ar[r]_-{a}&I \ar[r]_-{u}& J}}\quad = \quad
\vcenter{\xymatrix@C+15pt@R+2pt{
I\,|\,ua\drtwocell<\omit>{\phantom{aaaa}\Gamma\,|\,ua} \ar[d]_-{p_{I\,|\,ua}}\ar[r]^-{\pi\,|\,ua}&I\ar[d]^-{u}\\ e\ar[r]_-{ua}&J
}}
\end{equation}
$$
\left(\,\text{resp.}\quad \qquad
\vcenter{\xymatrix@C+8pt@R+5pt{ 
ua\,|\,I \drtwocell<\omit>{\phantom{aa}\Gamma}\ar[d]_-{ua\,|\,\pi}\ar[r]^-{p_{ua\,|\,I}}&e\ar[d]^-{a}\\
I\drtwocell<\omit>{\phantom{aa}\mathrm{id}} \ar[d]_-{\mathrm{id}}\ar[r]|-{\mathrm{id}}&I\ar[d]^-{u}\\
I\ar[r]_-{u}&J}}\quad = \quad
\vcenter{\xymatrix@C+8pt@R+5pt{ 
ua\,|\,I \drtwocell<\omit>{\phantom{aaaa}ua\,|\,\Gamma}\ar[d]_-{ua\,|\,\pi}\ar[r]^-{p_{ua\,|\,I}}&e\ar[d]^-{ua}\\
I\ar[r]_-{u}&J
}}\right)\,
$$
où la transformation naturelle:
\begin{equation}\label{pleinef}
\vcenter{\xymatrix@C+15pt@R+2pt{
I\,|\,ua\drtwocell<\omit>{\phantom{aa}\Gamma} \ar[d]_-{p_{I\,|\,ua}}\ar[r]^-{\pi\,|\,ua}&I\ar[d]^-{\mathrm{id}}\\ e\ar[r]_-{a}&I
}}
\qquad\;\qquad
\left(\,\text{resp.}\qquad
\vcenter{\xymatrix@+8pt{
ua\,|\,I\drtwocell<\omit>{\phantom{aa}\Gamma}\ar[r]^-{p_{ua\,|\,I}}\ar[d]_-{ua\,|\,\pi}&e\ar[d]^-{a}\\ I\ar[r]_-{\mathrm{id}}&I
}}\right),
\end{equation}
est définie dans un objet $(c,\alpha)$ de la catégorie $I\,|\,ua$,  comme le seul morphisme de $I$ dont son image par $u$ est égale à $\alpha$. 

Si $\C$ est une catégorie munie d'une famille distinguée de morphismes ${\bf W}$, on vérifie que toute conjuguée à gauche (resp. à droite):
\begin{equation}\label{estamera}
\vcenter{\xymatrix@+10pt{
\C^{I\,|\,ua}[{\bf W}_{I\,|\,ua}^{-1}]\ar[d]_-{\mathrm{hocolim}_{I\,|\,ua}}&\C^{I}[{\bf W}_{I}^{-1}]\ar[d]^-{\mathrm{id}}\ar[l]_-{{(\pi\,|\,ua)^*}}\\ 
\C[{\bf W}^{-1}]\urtwocell<\omit>{\phantom{aaaa}{\Gamma_!}}&\C^I[{\bf W}_{I}^{-1}]\ar[l]^-{{ua^*}}
}}\qquad\;\qquad
\left(\,\text{resp.}\qquad
\vcenter{\xymatrix@+10pt{
\C^{ua\,|\,I}[{\bf W}_{I\,|\,ua}^{-1}]\ar[r]^-{\mathrm{holim}_{ua\,|\,I}}& \C[{\bf W}^{-1}]\dltwocell<\omit>{{\Gamma_*}\phantom{aaaa}} \\ 
\C^I[{\bf W}^{-1}_I]\ar[u]^-{{(b\,|\,\pi)^*}}\ar[r]_-{\mathrm{id}}& \C^I[{\bf W}_{I}^{-1}] \ar[u]_-{ua^{*}}
}}\right)
\end{equation}
de la transformation naturelle qu'on déduit de \eqref{pleinef}:
\begin{equation}\label{fpleno}
\vcenter{\xymatrix@+10pt{
\C^{I\,|\,ua}[{\bf W}_{I\,|\,ua}^{-1}]\drtwocell<\omit>{\phantom{aaaa}{\Gamma}}&\C^{I}[{\bf W}_{I}^{-1}]\ar[l]_-{{(\pi\,|\,ua)^*}}\\ 
\C[{\bf W}^{-1}]\ar[u]^-{{(p_{I\,|\,ua})^*}}&\C^I[{\bf W}_{I}^{-1}]\ar[l]^-{{ua^*}}\ar[u]_-{{\mathrm{id}}}
}}\qquad\;\qquad
\left(\,\text{resp.}\qquad
\vcenter{\xymatrix@+10pt{
\C^{ua\,|\,I}[{\bf W}_{ua\,|\,I}^{-1}]\drtwocell<\omit>{\phantom{aa}{\Gamma}}&\C[{\bf W}^{-1}]\ar[l]_-{{(p_{ua\,|\,I})^*}}\\ 
\C^I[{\bf W}^{-1}_I]\ar[u]^-{{(ua\,|\,\pi)^*}}&\C^I[{\bf W}_{I}^{-1}]\ar[l]^-{{\mathrm{id}}}\ar[u]_-{{ua^*}}
}}\right)\,,
\end{equation}
est un isomorphismes.

En effet, cela se suit alors du Lemme \ref{isolemme} parce que le conjuguée à gauche de \eqref{fpleno} par rapport aux adjonctions:
$$
\C[{\bf W}^{-1}]\vcenter{\xymatrix@C+20pt{
\phantom{\cdot}\ar@{}[r]|-{\perp}\ar@<-4pt>@/_10pt/[r]_-{{(p_{I\,|\,ua})^*}}&
\phantom{\cdot}\ar@<-4pt>@/_10pt/[l]_-{{(a,\mathrm{id}_{ua})^*}}}}\C^{I\,|\,ua}[{\bf W}_{I\,|\,ua}^{-1}]
\qquad \text{et} \qquad
\C^I[{\bf W}_I^{-1}]\vcenter{\xymatrix@C+20pt{
\phantom{\cdot}\ar@{}[r]|-{\perp}\ar@<-4pt>@/_10pt/[r]_-{\mathrm{id}}&
\phantom{\cdot}\ar@<-4pt>@/_10pt/[l]_-{\mathrm{id}}}}\C^{I}[{\bf W}_{I}^{-1}]
$$
$$
\left( \text{resp.} \qquad \C[{\bf W}^{-1}]\vcenter{\xymatrix@C+20pt{
\phantom{\cdot}\ar@{}[r]|-{\perp}\ar@<+4pt>@/^10pt/[r]^-{{(p_{I\,|\,ua})^*}}&
\phantom{\cdot}\ar@<+4pt>@/^10pt/[l]^-{{(a,\mathrm{id}_{ua})^*}}}}\C^{I\,|\,ua}[{\bf W}_{I\,|\,ua}^{-1}]
\qquad \text{et} \qquad
\C^I[{\bf W}_I^{-1}]\vcenter{\xymatrix@C+20pt{
\phantom{\cdot}\ar@{}[r]|-{\perp}\ar@<+4pt>@/^10pt/[r]^-{\mathrm{id}}&
\phantom{\cdot}\ar@<+4pt>@/^10pt/[l]^-{\mathrm{id}}}}\C^{I}[{\bf W}_{I}^{-1}]\right)
$$
et aux transformations naturelles:
$$
\xymatrix@-10pt{ {(a,\mathrm{id}_{ua})^*} \circ {(p_{I\,|\,ua})^*}\ar@{=>}[r]^-{\mathrm{id}}&\mathrm{id}}\qquad \text{et} \qquad 
\xymatrix@-10pt{ \mathrm{id}\ar@{=>}[r]^-{\mathrm{id} }&\mathrm{id}\circ\mathrm{id}}
$$
$$
\left( \text{resp.} \qquad \xymatrix@-10pt{ {(a,\mathrm{id}_{ua})^*} \circ {(p_{I\,|\,ua})^*}\ar@{=>}[r]^-{\mathrm{id}}&\mathrm{id}}\qquad \text{et} \qquad 
\xymatrix@-10pt{ \mathrm{id}\ar@{=>}[r]^-{\mathrm{id} }&\mathrm{id}\circ\mathrm{id}}\right)
$$
est la transformation naturelle identité du foncteur $u:\xymatrix@C+5pt{\C^{I}[{\bf W}_{I}^{-1}]\ar[r]^{a^*}&\C[{\bf W}^{-1}]}$. 

Donc d'après l'égalité des transformations naturelles \eqref{upf} et selon les énoncés Lemme \ref{paste}, Lemme \ref{reflee} et Proposition \ref{local}, on en conclut que toutes les transformations conjuguées à gauche (resp. à droite) de la transformation naturelle:
$$
\vcenter{\xymatrix@+10pt{
\C^{I}[{\bf W}_{I}^{-1}]\drtwocell<\omit>{\phantom{aaaa}{\mathrm{id}}}&\C^{I}[{\bf W}_{I}^{-1}]\ar[l]_-{{\mathrm{id}}}\\ 
\C^I[{\bf W}_I^{-1}]\ar[u]^-{{\mathrm{id}}}&\C^I[{\bf W}_{I}^{-1}]\ar[l]^-{{u^*}}\ar[u]_-{{u^*}}
}}
$$
sont d'isomorphismes. 

En particulier:

\begin{corollaire}\label{localco}  
Soit $\C$ une catégorie de modèles et $u:\xymatrix@C-10pt{I\ar[r]&J}$ un foncteur pleinement fidèle entre petites catégories directes (resp. inverses), si $u_{!}$ $($resp. $u_{*})$ est un foncteur extension de Kan homotopique à gauche (resp. à droite) le long de $u$  (par rapport aux équivalences faibles ${\bf W}$ de $\C$):
$$
\xymatrix@C+15pt{
\C^J[{\bf W}_{J}^{-1}]\ar@/_14pt/[r]_{u^*}\ar@{}[r]|-{\perp}&\C^I[{\bf W}_{I}^{-1}]\ar@/_14pt/[l]_{{u_{!}}}}
\qquad\left(\;\text{resp.}\quad
\vcenter{
\xymatrix@C+15pt{
\C^J[{\bf W}_{J}^{-1}]\ar@/^14pt/[r]^{u^*}\ar@{}[r]|-{\perp}&\C^I[{\bf W}_{I}^{-1}]\ar@/^14pt/[l]^{{u_{*}}}}}\right)\,,
$$
alors $u_{!}$ $($resp. $u_{*})$ est aussi un foncteur pleinement fidèle. 
\end{corollaire}

Montrons aussi:

\begin{corollaire}\label{necimescimagg}
Soit $\C$ une catégorie de modèles, $u:\xymatrix@C-10pt{I\ar[r]&J}$ un foncteur pleinement fidèle entre petites catégories directes (resp. inverses) et $u_{!}$ $($resp. $u_{*})$ un foncteur extension de Kan homotopique à gauche (resp. à droite) le long de $u$  (par rapport aux équivalences faibles ${\bf W}$ de $\C$):
$$
\xymatrix@C+15pt{
\C^J[{\bf W}_{J}^{-1}]\ar@/_14pt/[r]_{u^*}\ar@{}[r]|-{\perp}&\C^I[{\bf W}_{I}^{-1}]\ar@/_14pt/[l]_{{u_{!}}}}
\qquad\left(\;\text{resp.}\quad
\vcenter{
\xymatrix@C+15pt{
\C^J[{\bf W}_{J}^{-1}]\ar@/^14pt/[r]^{u^*}\ar@{}[r]|-{\perp}&\C^I[{\bf W}_{I}^{-1}]\ar@/^14pt/[l]^{{u_{*}}}}}\right)\,.
$$

Si $\mathfrak{Y}$ est un $J$-diagramme de $\C$ alors $\mathfrak{Y}$ appartient à l'image essentielle du foncteur $u_{!}$ $($resp. $u_{*})$ si et seulement si, pour toute counité $\epsilon$ (resp. toute unité $\eta$) de l'adjonction $u_{!} \dashv {u^*}$ $\Big($resp. ${u^*}\dashv u_{*}\Big)$ et tout objet $x$ de $J$ qui n'appartient pas à l'image de $u$ le morphisme:
\begin{equation}\label{now}
\vcenter{\xymatrix@+5pt{
\big(u_{!}\circ{u^*} (\mathfrak{Y})\big)_{x} \ar[r]^-{(\epsilon_\mathfrak{Y})_x}  & \mathfrak{Y}_x}}\qquad \quad
\Bigg(\text{resp.} \quad \vcenter{\xymatrix@+5pt{
\mathfrak{Y}_x  \ar[r]^-{(\eta_\mathfrak{Y})_x} & \big(u_{*}\circ{u^*} (\mathfrak{Y})\big)_{x}}}\Bigg)\,,
\end{equation}
est un isomorphisme de $\C[{\bf W}^{-1}]$.
\end{corollaire}
\begin{proof}
Si $\xymatrix@-10pt{ \mathrm{id}\ar@{=>}[r]^-{\eta}&u^*\circ u_!}$ et $\xymatrix@-10pt{ u_!\circ u^* \ar@{=>}[r]^-{\varepsilon}&\mathrm{id}}$ $\Big(\text{resp.} \; \xymatrix@-10pt{ \mathrm{id}\ar@{=>}[r]^-{\eta}&u_*\circ u^*}$ et $\xymatrix@-10pt{ u^*\circ u_* \ar@{=>}[r]^-{\varepsilon}&\mathrm{id}}\Big)$ sont une unité et une counité respectivement de l'adjonction $u_{!} \dashv {u^*}$ $\Big($resp. ${u^*}\dashv u_{*}\Big)$ vérifiant les identités triangulaires, pour tout $I$-diagramme $\mathfrak{X}$ de $\C$ et tout objet $a$ de $I$ on a que le morphisme:
\begin{equation}\label{now2}
\vcenter{\xymatrix@+5pt{
\mathfrak{X}_a \ar[r]^-{(\eta_\mathfrak{X})_a} & u^*\circ u_! (\mathfrak{X})_a 
}}\qquad \; \Bigg(\text{resp.} \quad 
\vcenter{\xymatrix@+5pt{
 u_!\circ u^* ( \mathfrak{X})_a \ar[r]^-{(\varepsilon_\mathfrak{X})_a} &\mathfrak{X}_a 
}}\Bigg)
\end{equation}
est un isomorphisme de $\C[{\bf W}^{-1}]$ parce que d'après le Corollaire \ref{localco} les foncteurs $u_!$ et $u_*$ sont pleinement fidèle.

D'un autre côté, on note que si $x$ est un objet de $J$ tel que $x=u(a)$ pour certain objet $a$ de $I$, le morphisme \eqref{now} prend  la forme: 
$$
\vcenter{\xymatrix@+5pt{
{u^*}\big(u_{!}\circ {u^*} (\mathfrak{Y})\big)_{a} \ar[r]^-{{u^*}(\epsilon_\mathfrak{Y})_a}  & {u^*}(\mathfrak{Y})_a}}\qquad \quad
\Bigg(\text{resp.} \quad \vcenter{\xymatrix@+5pt{
u^*(\mathfrak{Y})_a  \ar[r]^-{u^*(\eta_\mathfrak{Y})_a} & u^*\big( u_{*}\circ {u^*} (\mathfrak{Y})\big)_{a}}}\Bigg)\,;
$$
et il est donc un isomorphisme car il est l'inverse à gauche (resp. à droite) de l'isomorphisme \eqref{now2} où $\mathcal{X}= {u^*}\mathcal{Y}$.
\end{proof}

\renewcommand{\thesubsection}{\S\thesection.\arabic{subsection}}
\subsubsection{}\;     
\renewcommand{\thesubsection}{\thesection.\arabic{subsection}}

Montrons l'énoncé suivant:

\begin{corollaire}\label{eqcisi}
Si $\C$ est une catégorie de modèles et $\xymatrix@C-8pt{I\ar[r]^u&J}$ est un foncteur entre petites catégories lesquelles sont directes et inverses, alors les conditions qui suit sont équivalentes:
\begin{enumerate}
\item Toute conjuguée à gauche $\xymatrix@C-10pt{u_!\circ p_I^*\ar@{=>}[r]& p_J^*}$ $\big($resp. à droite $\xymatrix@C-10pt{p_J^*\ar@{=>}[r]& u_*\circ p_I^*}\big)$ de la transformation naturelle identité $\xymatrix@C-10pt{p_I^*\circ\mathrm{id}^*\ar@{=>}[r]& u^*\circ p_J^*}$  est un isomorphismes.
\item Toute conjuguée à droite $\xymatrix@C-10pt{\underset{J}{\mathrm{holim}}\ar@{=>}[r]& \underset{I}{\mathrm{holim}}\circ u^*}$ $\big($resp. à gauche $\xymatrix@C-10pt{\underset{I}{\mathrm{hocolim}}\circ u^*\ar@{=>}[r]& \underset{J}{\mathrm{hocolim}}}\big)$ de la transformation naturelle identité $\xymatrix@C-10pt{p_I^*\circ\mathrm{id}^*\ar@{=>}[r]& u^*\circ p_J^*}$ $\big($resp. $\xymatrix@C-10pt{u^*\circ p_J^*\ar@{=>}[r]&p_I^*\circ\mathrm{id}^*}\big)$ est un isomorphismes.
\item Si $b$ est un objet quelconque de $J$ et $\underset{I\,|\,b}{\mathrm{hocolim}}$ $\Big(\text{resp.}$ $\underset{I\,|\,b}{\mathrm{holim}}\Big)$ est un adjoint à gauche de $p_{I\,|\,b}^*$ $\Big(\text{resp.}$ à droite de $p_{b\,|\,I}^*\Big)$ toute counité (resp. toute unité): 
$$
\xymatrix@C-10pt{\underset{I\,|\,b}{\mathrm{hocolim}}\circ p_{I\,|\,b}^*\ar@{=>}[r]&\mathrm{id}_{\C[{\bf W}^{-1}]}} \qquad 
\Big(\text{resp.} \qquad 
\xymatrix@C-10pt{\mathrm{id}_{\C[{\bf W}^{-1}]}\ar@{=>}[r]& \underset{b\,|\,I}{\mathrm{holim}}\circ p_{b\,|\,I}^*}
\Big)
$$
est un isomorphisme.
\item Si $b$ est un objet quelconque de $J$ le foncteur: 
$$
\xymatrix{\C[{\bf W}^{-1}]\ar[r]^-{p_{I\,|\,b}^*}&\C^{I\,|\,b}[{\bf W}_{I\,|\,b}^{-1}]}\qquad\;
\bigg(\text{resp.}\;\xymatrix{\C[{\bf W}^{-1}]\ar[r]^-{p_{b\,|\,I}^*}&\C^{b\,|\,I}[{\bf W}_{b\,|\,I}^{-1}]}\bigg)
$$ 
est pleinement fidèle pour tout objet $b$ de $J$.
\end{enumerate}

Plus encore la condition (iv) se vérifie si on suppose un des deux énoncés:
\renewcommand{\labelenumi}{\arabic{enumi}.}
\begin{enumerate}
\item Le foncteur $\xymatrix@C-8pt{I\ar[r]^u&J}$ admet un adjoint à droite $u\dashv v$ (resp. un adjoint à gauche $v\dashv u$).
\item Les catégories $I$ et $J$ admettent un objet initial (resp. final) préservé par le foncteur $\xymatrix@C-8pt{I\ar[r]^u&J}$.
\end{enumerate}
\renewcommand{\labelenumi}{(\roman{enumi})}
\end{corollaire}
\begin{proof}
Les conditions (i) et (ii) sont équivalentes d'après le Lemme \ref{matesconjureal} et les conditions (iii) et (iv) sont équivalentes par des propriétés des adjonctions.

D'un autre côté si $b$ est un objet quelconque de $J$ considérons l'égalité suivante entre transformations naturelles:
$$
\vcenter{\xymatrix@C+18pt@R+6pt{
I\,|\,b\drtwocell<\omit>{\phantom{aa}\Gamma\,|\,b} \ar[d]_-{p_{I\,|\,b}}\ar[r]^-{\pi\,|\,b}&
I\drtwocell<\omit>{\phantom{aa}\mathrm{id}} \ar[d]|-{u} \ar[r]^-{p_I} & e\ar[d]^-{\mathrm{id}}\\
 e\ar[r]_-{b}&J \ar[r]_-{p_J}& e}}\quad = \quad
\vcenter{\xymatrix@C+15pt@R+2pt{
I\,|\,b\drtwocell<\omit>{\phantom{aa}\mathrm{id}} \ar[d]_-{p_{I\,|\,b}}\ar[r]^-{p_{I\,|\,b}}&e\ar[d]^-{\mathrm{id}}\\ e\ar[r]_-{\mathrm{id}}&e
}}
$$
$$
\left(\,\text{resp.}\quad \qquad
\vcenter{\xymatrix@C+8pt@R+5pt{ 
b\,|\,I \drtwocell<\omit>{\phantom{aa}b\,|\,\Gamma}\ar[d]_-{b\,|\,\pi}\ar[r]^-{p_{b\,|\,I}}&e\ar[d]^-{b}\\
I\drtwocell<\omit>{\phantom{aa}\mathrm{id}} \ar[d]_-{p_I}\ar[r]|-{u}&J\ar[d]^-{p_J}\\
e\ar[r]_-{\mathrm{id}}&e}}\quad = \quad
\vcenter{\xymatrix@C+8pt@R+5pt{ 
b\,|\,I \drtwocell<\omit>{\phantom{aa}\mathrm{id}}\ar[d]_-{p_{b\,|\,I}}\ar[r]^-{p_{b\,|\,I}}&e\ar[d]^-{\mathrm{id}}\\
e\ar[r]_-{\mathrm{id}}&e
}}\right)\,
$$

On déduit que (i)$\Leftrightarrow$(iii) de la Proposition \ref{local} et les Lemmes \ref{paste} et \ref{reflee}. 

Enfin remarquons que si on suppose un des deux conditions 1. ou 2. de l'énoncé du Corollaire, alors le foncteur $p_{I\,|\,b}$ $\big($resp. $p_{b\,|\,I} \big)$ admet un adjoint à droite (resp. à gauche) pleinement fidèle. Donc le foncteur $p_{I\,|\,b}^*$ $\big($resp. $p_{b\,|\,I}^* \big)$ est pleinement fidèle d'après le Lemme \ref{cwgen}.
\end{proof}

\section{Carrés homotopiquement (co)cartésiens}

\renewcommand{\thesubsection}{\S\thesection.\arabic{subsection}}
\subsection{}\;\label{homcartcarrff}
\renewcommand{\thesubsection}{\thesection.\arabic{subsection}}

Soit $\C$ une catégorie de modèles. On considère la petite catégorie engendrée par le diagramme:
$$
\boxempty \; = \;
\def\objectstyle{\scriptstyle}
\def\labelstyle{\scriptstyle}
\vcenter{\xymatrix@-10pt{a\ar[r]\ar[d] & b\ar[d] \\ d\ar[r] & c}}\,,
$$
soumis à la relation qu'identifie les deux flèches de $a$ vers $c$; autrement dit, la donnée d'un $\boxempty$-diagramme $\mathfrak{X}$ dans $\C$ équivaut à la donnée d'un carré commutatif dans $\C$:
\begin{equation}\label{cX}
\mathfrak{X} \; = \; \vcenter{\xymatrix{X_{a}\ar[r]\ar[d]&X_{b}\,\ar[d]\\X_{d}\ar[r]&X_{c}\,.}}
\end{equation}

Si on considère la sous-catégorie suivante de $\boxempty$:
$$
\text{\Large $\lefthalfcap$} \; = \;
\def\objectstyle{\scriptstyle}
\def\labelstyle{\scriptstyle}
\vcenter{\xymatrix@-20pt{
a\ar[rr]\ar[dd] && b &&&&& a\ar[rr]\ar[dd]&&b\ar[dd]\\
&& &\phantom{.}\ar@{^(->}[rrr]^-{q_{\lefthalfcap}}&&&&&&\\
d&& & & &&& d\ar[rr] & & c}} 
$$
$$
\left(\text{resp.} \qquad 
\text{\Large $\righthalfcup$}  \; = \;
\def\objectstyle{\scriptstyle}
\def\labelstyle{\scriptstyle}
\vcenter{\xymatrix@-20pt{
&& b\ar[dd]&&&&& a\ar[rr]\ar[dd]&&b \ar[dd]\\
&& &\phantom{.}\ar@{^(->}[rrr]^-{q_{\righthalfcup}}&&&&&&\\
d\ar[rr]&& c& & &&& d\ar[rr] & & c}} \right);
$$
on déduit facilement du Corollaire \ref{hocompco}, l'existence d'une adjonction:
$$
\C^{\lefthalfcap}[{\bf W}_{\lefthalfcap}^{-1}]
\vcenter{\xymatrix@C+10pt{
\phantom{\cdot}\ar@{}[r]|-{\perp}\ar@<+4pt>@/^10pt/[r]^-{q_{\lefthalfcap\,!}}&
\phantom{\cdot}\ar@<+4pt>@/^10pt/[l]^-{{q_{\lefthalfcap}^*}}}}
\C^{\boxempty}[{\bf W}_{\boxempty}^{-1}]
\qquad\qquad
\left(\;\text{resp.}\quad 
\C^{\righthalfcup}[{\bf W}_{\righthalfcup}^{-1}]
\vcenter{\xymatrix@C+10pt{
\phantom{\cdot}\ar@{}[r]|-{\perp}\ar@<-4pt>@/_10pt/[r]_-{{q_{\righthalfcup}^*}}&
\phantom{\cdot}\ar@<-4pt>@/_10pt/[l]_-{q_{\righthalfcup\,*}}}}
\C^{\boxempty}[{\bf W}_{\boxempty}^{-1}]
\right).
$$

Étant donné un carré commutatif \eqref{cX} comme ci-dessus, on dit que $\mathfrak{X}$ est un carré \emph{homotopiquement cocartésien} (resp. {homotopiquement cartésien}) si la valeur en $\mathfrak{X}$ d'une counité $\xymatrix@-10pt{q_{\lefthalfcap\, !}\; {q_{\lefthalfcap}^*}\ar@{=>}[r]^-{\varepsilon} & \mathrm{id}}$ $\big($resp. de l'unité  \xymatrix@-10pt{\mathrm{id}\ar@{=>}[r]^-{\eta}&q_{\righthalfcup}^*\,{q_{\righthalfcup\; *}}}$\big)$ de l'adjonction $q_{\lefthalfcap\; !} \dashv {q_{\lefthalfcap}^*}$ $\big($resp. ${q_{\righthalfcup}^*} \dashv q_{\righthalfcup\; *}$ $\big)$:
$$
\vcenter{\xymatrix@+2pt{
q_{\lefthalfcap\,!} \; {q_{\lefthalfcap}^*} (\mathfrak{X})\ar[r]^-{\varepsilon_{\mathfrak{X}}}& \mathfrak{X}}}
\qquad\quad
\left(\,\text{resp.}\quad
\vcenter{\xymatrix@+2pt{
\mathfrak{X} \ar[r]^-{\eta_{\mathfrak{X}}} & q_{\righthalfcup\,*} \; {q_{\righthalfcup}^*} (\mathfrak{X})}}\right)\,.
$$
est un isomorphisme dans la catégorie $\C^{\boxempty}[{\bf W}_{\boxempty}^{-1}]$. 

Pour donner une définition équivalente considérons la transformation naturelle:
\begin{equation}\label{Varva}
\def\objectstyle{\scriptstyle}
\def\labelstyle{\scriptstyle}
\vcenter{\xymatrix@C-13pt@R-18pt{
\C^{\lefthalfcap}[{\bf W}_{\lefthalfcap}^{-1}]\;\rrtwocell<\omit>{<3> {\Phi}}&&\C^{\boxempty}[{\bf W}_{\boxempty}^{-1}]\ar[ll]_-{{q_{\lefthalfcap}^*}}\ar[ldd]^{{c^*}}\\
&&\\&\C[{\bf W}^{-1}] \ar[luu]^{{p_{\lefthalfcap}^*}} &}}
\qquad\qquad
\left(\;\text{resp.}\quad 
\def\objectstyle{\scriptstyle}
\def\labelstyle{\scriptstyle}
\vcenter{\xymatrix@C-13pt@R-18pt{
\C^{\righthalfcup}[{\bf W}_{\righthalfcup}^{-1}]\;&&\C^{\boxempty}[{\bf W}_{\boxempty}^{-1}]\ar[ll]_-{{q_{\righthalfcup}^*}}\ar[ldd]^{{a^*}}\lltwocell<\omit>{<-3> {\Psi}}\\
&&\\&\C[{\bf W}^{-1}] \ar[luu]^{{p_{\lefthalfcup}^*}} &}}\right),
\end{equation}
induite par la transformation naturelle canonique:
\begin{equation}\label{varva}
\def\objectstyle{\scriptstyle}
\def\labelstyle{\scriptstyle}
\vcenter{\xymatrix@C-5pt@R-18pt{
\text{\Large $\lefthalfcap$}\;\ar[rdd]_{p_{\lefthalfcap}}\ar@{^(->}[rr]^-{q_{\lefthalfcap}}\rrtwocell<\omit>{<3> \Phi}&&\text{\large $\boxempty$}\\&&\\&e\ar[ruu]_{c}&}}
\qquad\qquad
\left(\;\text{resp.}\quad 
\def\objectstyle{\scriptstyle}
\def\labelstyle{\scriptstyle}
\vcenter{\xymatrix@C-5pt@R-18pt{
\text{\Large $\righthalfcup$}\;\ar[rdd]_{p_{\righthalfcup}}\ar@{^(->}[rr]^-{q_{\righthalfcup}}&&\text{\large $\boxempty$}\lltwocell<\omit>{<-3> \Psi}\\
&&\\&e\ar[ruu]_{a}&}}\right);
\end{equation}
et prenons un foncteur adjoint à gauche (resp. à droite) de ${p_{\lefthalfcap}^*}$ $\big($resp. ${p_{\righthalfcup}^*}$ $\big)$:
$$
\vcenter{\xymatrix@+10pt{\C^{\lefthalfcap}[{\bf W}_{\lefthalfcap}^{-1}]\ar[r]^-{\mathrm{hocolim}_{\lefthalfcap}}&\C[{\bf W}^{-1}]}}
\qquad\;\left(\,\text{resp.}\quad 
\vcenter{\xymatrix@+10pt{\C^{\righthalfcup}[{\bf W}_{\righthalfcup}^{-1}]\ar[r]^-{\mathrm{holim}_{\righthalfcup}}&\C[{\bf W}^{-1}]}}
\right),
$$
lequel existe grâce au Corollaire \ref{hocompco}, et qu'on appelle un (ou le) foncteur \emph{somme amalgamée homotopique} (resp. \emph{produit fibré homotopique}) de $\C$. 

\begin{lemme}\label{carrlele}
Soit $\C$ une catégorie de modèles. Si $\mathfrak{X}$ est un carré commutatif de $\C$:
$$
\mathfrak{X} \; = \; \vcenter{\xymatrix{X_{a}\ar[r]\ar[d]&X_{b}\,\ar[d]\\X_{d}\ar[r]&X_{c}\,,}}
$$ 
les énoncés suivants sont équivalents:
\begin{enumerate}
\item $\mathfrak{X}$ est un carré homotopiquement cocartésien (resp. homotopiquement cartésien), \emph{i.e.} le morphisme canonique:
$$
\vcenter{\xymatrix@+2pt{
q_{\lefthalfcap\,!} \; {q_{\lefthalfcap}^*} (\mathfrak{X})\ar[r]^-{\varepsilon_{\mathfrak{X}}}& \mathfrak{X}}}
\qquad\quad
\left(\,\text{resp.}\quad
\vcenter{\xymatrix@+2pt{
\mathfrak{X} \ar[r]^-{\eta_{\mathfrak{X}}} & q_{\righthalfcup\,*} \; {q_{\righthalfcup}^*} (\mathfrak{X})}}\right)\,,
$$
est un isomorphisme de la catégorie $\C^{\boxempty}[{\bf W}_{\boxempty}^{-1}]$.
\item $\mathfrak{X}$  est isomorphe à un objet dans l'image du foncteur $q_{\lefthalfcap\,!}$ (resp. $q_{\righthalfcup\,*}$). 
\item Le morphisme canonique:
$$
\vcenter{\xymatrix@+2pt{
q_{\lefthalfcap\,!} \; {q_{\lefthalfcap}^*} (\mathfrak{X})_{c}\ar[r]^-{(\varepsilon_{\mathfrak{X}})_{c}}& X_{c}}}
\qquad\quad
\left(\,\text{resp.}\quad
\vcenter{\xymatrix@+2pt{
X_{a} \ar[r]^-{(\eta_{\mathfrak{X}})_{a}} & q_{\righthalfcup\,*} \; {q_{\righthalfcup}^*} (\mathfrak{X})_{a}}}\right)\,,
$$
est un isomorphisme de la catégorie $\C[{\bf W}^{-1}]$. 
\item Le morphisme suivant:
$$
\vcenter{\xymatrix@+55pt{\mathrm{hocolim}_{\lefthalfcap}\big( {q_{\lefthalfcap}^*} (\mathfrak{X})\big)\ar[r]^-{\mathrm{hocolim}_{\lefthalfcap}({\Phi}_{\mathfrak{X}})} & \mathrm{hocolim}_{\lefthalfcap}\big(  {p_{\lefthalfcap}^*} \; ( X_{c} )\big)}}
$$
$$
\left(\,\text{resp.}\quad
\vcenter{\xymatrix@+40pt{\mathrm{holim}_{\righthalfcup} \big( {p_{\righthalfcup}^*} \;  ( X_{a} )\big)
\ar[r]^-{\mathrm{holim}_{\righthalfcup}({\Phi}_{\mathfrak{X}})} &
\mathrm{holim}_{\righthalfcup}\big({q_{\righthalfcup}^*}(\mathfrak{X})\big)}}\right)\,,
$$
est un isomorphisme de la catégorie homotopique de $\C$. 
\item Le morphisme suivant:
\begin{equation}\label{defca}
\vcenter{\xymatrix{\mathrm{hocolim}_{\lefthalfcap}\big( {q_{\lefthalfcap}^*} (\mathfrak{X}) \big)\ar[r]&  X_{c}}}\qquad\,
\left(\text{resp.}\quad
\vcenter{\xymatrix{X_{a} \ar[r]&\mathrm{holim}_{\righthalfcup}\big( {q_{\righthalfcup}^*} (\mathfrak{X}) \big)}}
\right)\,,
\end{equation}
obtenu du morphisme $\xymatrix@+5pt{{q_{\lefthalfcap}^*} \, (\mathfrak{X}) \ar[r]^-{{\Phi}_{\mathfrak{X}}}&  {p_{\lefthalfcap}^*} \; ( X_{c} )}$ $\big($resp. $\xymatrix@+5pt{{p_{\righthalfcup}^*} \; (X_{a}) \ar[r]^-{{\Psi}_{\mathfrak{X}}}& {q_{\righthalfcup}^*} \, (\mathfrak{X})}$ $\big)$ par adjonction, est un isomorphisme de la catégorie homotopique de $\C$. 
\end{enumerate}
\end{lemme}
\begin{proof}
D'après le Corollaire \ref{localco} le foncteur $q_{\lefthalfcap\,!}$ (resp. $q_{\righthalfcup\,*}$) est pleinement fidèle, car $q_{\lefthalfcap}$ (resp. $q_{\righthalfcup}$) l'est. Alors les énoncés (i) et (ii) sont équivalents grâce à un argument général concernant les adjonctions. 

On vérifie que les énoncés (ii) et (iii) sont équivalents d'après le Corollaire \ref{necimescimagg}. 

Remarquons maintenant que d'après le Lemme \ref{cwgen}, vu que $\text{\Large $\lefthalfcap$}$ (resp. $\text{\Large $\righthalfcup$}$) admet un objet initial (resp. final), pour tout objet $A$ de $\C$ la counité de l'adjonction $\mathrm{hocolim}_{\lefthalfcap} \dashv{p_{\lefthalfcap}^*}$ (resp. l'unité de ${p_{\righthalfcup}^*}\dashv\mathrm{holim}_{\righthalfcup}$):
$$
\vcenter{\xymatrix{\mathrm{hocolim}_{\lefthalfcap}\big(  {p_{\lefthalfcap}^*} A \big)\ar[r]& A}}\qquad
\left(\text{resp.}\quad
\vcenter{\xymatrix{
A\ar[r]& \mathrm{holim}_{\righthalfcup} \big( {p_{\righthalfcup}^*} A\big)}}
\right)\,,   
$$
est un isomorphisme.  Donc, les énoncés (iv) et (v) sont équivalents.

Montrons finalement (iii)$\iff$(v) dans le cas qui concerne la catégorie $\lefthalfcap$. Pour cela remarquons que les carrés extérieurs des diagrammes suivants sont égaux: 
\begin{equation}\label{cc}  
\vcenter{\xymatrix@+8pt{
\lefthalfcap\,|\,c\ar@{}[rd]|-{=} \ar[d]_-{p_{\lefthalfcap\,|\,c}}\ar[r]^-{\pi\,|\,c}&\lefthalfcap\drtwocell<\omit>{\phantom{aa}\Phi}\ar[d]|-{p_{\lefthalfcap}}\ar[r]^-{q_{\lefthalfcap}}&\boxempty\ar[d]^-{\mathrm{id}}\\ 
e\ar[r]_-{\mathrm{id}}&e\ar[r]_{c}& \boxempty } }
\;\qquad\text{et}\qquad\;
\vcenter{\xymatrix@+8pt{
\lefthalfcap\,|\,c\drtwocell<\omit>{\phantom{aa}\Gamma\,|\,c} \ar[d]_-{p_{\lefthalfcap\,|\,c}}\ar[r]^-{\pi\,|\,c}&\lefthalfcap\ar[d]|-{q_{\lefthalfcap}}\ar[r]^-{q_{\lefthalfcap}}\ar@{}[rd]|-{=}&\boxempty\ar[d]^-{\mathrm{id}}\\ 
e\ar[r]_-{c}&\boxempty\ar[r]_{\mathrm{id}}& \boxempty }}\,;
\end{equation}
ce qui permet de voir que le morphisme $\xymatrix{\mathrm{hocolim}_{\lefthalfcap}\big( {q_{\lefthalfcap}^*}(\mathfrak{X})\big)\ar[r]&\mathfrak{X}_{c}}$ de l'énoncé (v) qu'on construit à partir de la transformation naturelle $\Phi$, s'insère dans un carré commutatif de la catégorie $\C[{\bf W}^{-1}]$:
\begin{equation}\label{carrebi}
\vcenter{\xymatrix@R+6pt@C+15pt{
\mathrm{hocolim}_{\lefthalfcap\,|\,c}\big( {\pi\,|\,c^*}\,{q_{\lefthalfcap}^*}\, \mathfrak{X}\big)\ar[r]\ar[d]_-{({\Gamma\,|\,c})_{\mathfrak{X}}}&\mathrm{hocolim}_{\lefthalfcap}\big( {q_{\lefthalfcap}^*} \mathfrak{X} \big)\ar[d]\\
q_{\lefthalfcap\,!} \; {q_{\lefthalfcap}^*} (\mathfrak{X})_{c}\ar[r]_-{(\varepsilon_{\mathfrak{X}})_{c}}& \mathfrak{X}_{c}\,.}}
\end{equation} 

Enfin, il résulte que le morphisme de la ligne supérieure de \eqref{carrebi} est un isomorphisme parce que $\pi\,|\,c$ est un isomorphisme de catégories, et que $({\Gamma\,|\,c})_{\mathfrak{X}}$ l'est aussi d'après la Proposition \ref{local}. Donc, (iii)$\iff$(v).
\end{proof}

\renewcommand{\thesubsection}{\S\thesection.\arabic{subsection}}
\subsection{}\;
\renewcommand{\thesubsection}{\thesection.\arabic{subsection}}

L'énoncé suivant est immédiat avec la définition du paragraphe précédent.

\begin{lemme}\label{cubelele}
Soit $\C$ une catégorie de modèles. On considère un morphisme de $\boxempty$-diagrammes $\xymatrix{\mathfrak{X}\ar[r]^F&\mathfrak{X}'}$  de $\C$, tel que $\xymatrix{{q_{\lefthalfcap}^*} (\mathfrak{X}) \ar[r]^-{{q_{\lefthalfcap}^*} F}&{q_{\lefthalfcap}^*} (\mathfrak{X}')}$ $\big($resp. $\xymatrix{{q_{\righthalfcup}^*} (\mathfrak{X}) \ar[r]^-{{q_{\righthalfcup}^*} F}&{q_{\righthalfcup}^*} (\mathfrak{X}')}$$\big)$ soit un isomorphisme dans $\C^{\lefthalfcap}[{\bf W}_{\lefthalfcap}^{-1}]$ $\big($resp. $\C^{\righthalfcup}[{\bf W}_{\righthalfcup}^{-1}]$$\big)$; autrement dit, on a un cube:
\begin{equation}\label{cubebe} 
\xymatrix@-15pt{
&&X_{a}'\ar[rrr]\ar[ddd]&&&X_{b}'\,\ar[ddd]\\
X_{a}\ar[rru]|-{\,F_{a}\,}\ar[rrr]\ar[ddd]&&&X_{b}\ar[ddd]\ar[rru]|-{\,F_{b}\,}&&\\
&&&&&\\
&&X_{d}'\ar[rrr]&&&X_{c}'\,,\\
X_{d}\ar[rrr]\ar[urr]|-{\,F_{d}\,}&&&X_{c}\ar[rru]|-{\,F_{c}\,}}
\end{equation}
dont les faces sont des carrés commutatifs de $\C$, et les morphismes $F_{a}$, $F_{b}$ et $F_{d}$ $($resp. $F_{b}$, $F_{c}$ et $F_{d}$$)$ d'équivalences faibles de $\C$.

Alors, si $\mathfrak{X}$ est un carré homotopiquement cocartésien $($resp. homotopiquement cartésien$)$, il en est de même pour $\mathfrak{X}'$, si et seulement si $F$ est un isomorphisme dans $\C^{\boxempty}[{\bf W}_{\boxempty}^{-1}]$, \emph{i.e.} si et seulement si le morphisme $F_{c}$ $($resp. $F_{a}$$)$ est une équivalence faible de $\C$. 
\end{lemme}

Le Lemme qui suit donne des conditions suffisants pour qu'un carré cocartésien (resp. cartésien) soit homotopiquement cocartésien (resp. homotopiquement cartésien).

\begin{lemme}\label{fibrabra}
Soit $\C$ une catégorie de modèles. On considère un carré commutatif de $\C$:
\begin{equation}\label{sisiX}
\mathfrak{X} \; = \; \vcenter{\xymatrix{X_{a}\ar[r]^f\ar[d]&X_{b}\,\ar[d]^g\\X_{d}\ar[r]&X_{c}\,.}}
\end{equation}

Alors, si $\mathfrak{X}$ est un carré cocartésien (resp. cartésien), pour que $\mathfrak{X}$ soit homotopiquement cocartésien (cartésien) il suffit que les sommets de $\mathfrak{X}$ soient des objets cofibrants (resp. fibrants), et que la flèche $f$ soit une cofibration (resp. $g$ soit une fibration)\footnote{Voir aussi \S10 de \cite{DS}.}.
\end{lemme}
\begin{proof}
Esquissons une preuve dans le cas de carrés cocartésiens. Pour commencer on considère la catégorie de Reedy $(\lefthalfcap,\lefthalfcap_{+},\lefthalfcap_{-},\lambda)$, où $\lambda:\xymatrix@-6pt{\{a,b,d\}\ar[r]&\mathbb{N}}$ est la fonction définie par $\lambda(d)=0$, $\lambda(a)=1$ et $\lambda(b)=2$, et $\xymatrix@-5pt{\lefthalfcap_{+}\ar[r]^{\nu_{+}}&\lefthalfcap}$ (resp. $\xymatrix@-5pt{\lefthalfcap_{-}\ar[r]^{\nu_{-}}&\lefthalfcap}$) est la sous-catégorie de $\lefthalfcap$ engendrée par le morphisme $f$ (resp. $g$).

Si on munie $\C^{\lefthalfcap}$ de la structure de catégorie de modèles de la Proposition \ref{hocomp}, on montre qu'un $\lefthalfcap$-diagramme $\widetilde{\mathfrak{X}} \; = \; \vcenter{\xymatrix{X_{a}\ar[r]^f\ar[d]&X_{b}\\X_{d}&}}$ de $\C$ est cofibrant, si et seulement si les sommets $X_{a}$, $X_{b}$, et $X_{d}$ sont des objets cofibrants et le morphisme $f$ est une cofibration de $\C$.

L'énoncé qu'on veut montrer est alors une conséquence de la propriété (iii) du Corollaire \ref{carrlele}.


\end{proof}

Montrons maintenant:

\begin{lemme} \label{equivasinffmmf}
Soit $\C$ une catégorie de modèles. On considère un carré commutatif de $\C$:
$$
\mathfrak{X} \; = \; \vcenter{\xymatrix{X_{a}\ar[r]^-f\ar[d]_{g'}&X_{b}\,\ar[d]^-{g}\\X_{d}\ar[r]_{f'}&X_{c}\,}}.
$$

\begin{enumerate} 
\item Si $\mathfrak{X}$ est homotopiquement cocartésien (resp. homotopiquement cartésien), alors:
$$\mathfrak{X}' \; = \; \vcenter{\xymatrix{X_{a}\ar[r]^-{g'}\ar[d]_-{f}&X_{d}\,\ar[d]^{f'}\\X_{b}\ar[r]_{g}&X_{c}\,}}$$
l'est aussi.
\item Si $f$ et $f'$ sont des équivalences faibles, alors $\mathfrak{X}$ est homotopiquement cartésien et homotopiquement cocartésien.
\end{enumerate}
\end{lemme}
\begin{proof}
On va montrer la partie des énoncés qui correspond aux carrés homotopiquement cocartésiens. 

Pour vérifier (i) considérons l'isomorphisme de catégories $\xymatrix@-5pt{\boxempty\ar[r]^{\theta_{\boxempty}}&\boxempty}$, qui échange les objets $c$ et $b$, et laisse invariants $a$ et $d$. Soit $\theta_{\lefthalfcap}$ la restriction de $\theta_{\boxempty}$ à la sous-catégorie $\text{\Large $\lefthalfcap$}$ de $\boxempty$.

Dans le carré commutatif:
$$
\vcenter{\xymatrix@C+10pt{
\C^{\lefthalfcap}[{\bf W}_{\lefthalfcap}^{-1}]\ar[d]_{{\theta_{\lefthalfcap}^*}}^{\cong}&\ar[l]_-{{q_{\lefthalfcap}^*}}
\C^{\boxempty}[{\bf W}_{\boxempty}^{-1}]\ar[d]^-{{\theta_{\boxempty}^*}}_{\cong}\\
\C^{\lefthalfcap}[{\bf W}_{\lefthalfcap}^{-1}]&\ar[l]^-{{q_{\lefthalfcap}^*}}
\C^{\boxempty}[{\bf W}_{\boxempty}^{-1}]\,,}}
$$
les foncteur ${\theta_{\lefthalfcap}^*}$ et ${\theta_{\boxempty}^*}$ sont des isomorphismes; donc, on a un carré commutatif à isomorphisme près:
$$
\vcenter{\xymatrix@C+10pt{
\C^{\lefthalfcap}[{\bf W}_{\lefthalfcap}^{-1}]\ar[d]_{{\theta_{\lefthalfcap}^*}}^{\cong}\ar[r]^-{{\bf L}{q_{\lefthalfcap\,!}}}&
\C^{\boxempty}[{\bf W}_{\boxempty}^{-1}]\ar[d]^-{{\theta_{\boxempty}^*}}_{\cong}\\
\C^{\lefthalfcap}[{\bf W}_{\lefthalfcap}^{-1}]\ar[r]_-{{\bf L}{q_{\lefthalfcap\,!}}}&
\C^{\boxempty}[{\bf W}_{\boxempty}^{-1}]\,.}}
$$ 

Puisque $\mathfrak{X}'= {\theta_{\boxempty}^*} (\mathfrak{X})$, on déduit que $\mathfrak{X}$ est isomorphe à un objet dans l'image du foncteur ${q_{\lefthalfcap\,!}}$, si et seulement si il en est de même pour $\mathfrak{X}'$. Le résultat désiré est alors une conséquence du Lemme \ref{carrlele}.

Pour montrer l'énoncé (ii), on va construire une équivalence faible argument par argument de $\boxempty$-diagrammes $F:\xymatrix@-8pt{\widetilde{\mathfrak{X}}\ar[r]&\mathfrak{X}}$, tel que le carré source $\widetilde{\mathfrak{X}}$ soit homotopiquement cocartésien. Le morphisme sera donné par un cube:
\begin{equation}\label{cubeproof}
\xymatrix@-15pt{
&&X_{a}\ar[rrr]^f\ar'[d][ddd]&&&X_{b}\,\ar[ddd]\\
QX_{a}\ar[rru]|-{\,F_{a}\,}\ar@<-1.5pt>[rrr]^(.65){\varphi}\ar[ddd]_-{\psi}&&&X_{b}'\ar@<+1.5pt>[ddd]\ar[rru]|-{\,F_{b}\,}&&\\
&&&&&\\
&&X_{d}\ar'[r][rrr]^(.28){f'}&&&X_{c}\,,\\
X_{d}'\ar[rrr]_{\varphi'}\ar[urr]|-{\,F_{d}\,}&&&X_{d}'\underset{QX_{a}}{\sqcup}X_{b}'\ar[rru]|-{\,F_{c}\,}}
\end{equation}
construit de la fa\c con suivante:  Pour commencer on prend $F_{a}:\xymatrix@-5pt{QX_{a}\ar[r]^-{\text{\rotatebox[origin=c]{180}{\Large $\widetilde{\phantom{w}}$}}}&X_{a}}$, un remplacement cofibrant de $X_{a}$; puis, on considère des factorisations:
$$
\vcenter{\xymatrix@R-19pt{\underset{\phantom{a}}{QX_{a}}\ar@/_5pt/@<-4pt>@{>->}[rd]_-{\varphi}\ar[r]^-{\text{\rotatebox[origin=c]{180}{\Large $\widetilde{\phantom{w}}$}}}&X_{a}\ar[r]&X_{b}\\&X_{b}'\ar@<-4pt>@/_3pt/[ru]_{F_{b}}^-{\text{\rotatebox[origin=c]{215}{\Large $\widetilde{\phantom{w}}$}}}&}}
\qquad\;\text{et}\;\qquad
\vcenter{\xymatrix@R-19pt{\underset{\phantom{a}}{QX_{a}}\ar@/_5pt/@<-4pt>@{>->}[rd]_-{\psi}\ar[r]^-{\text{\rotatebox[origin=c]{180}{\Large $\widetilde{\phantom{w}}$}}}&X_{a}\ar[r]&X_{d}\,.\\&X_{d}'\ar@<-4pt>@/_3pt/[ru]_{F_{d}}^-{\text{\rotatebox[origin=c]{215}{\Large $\widetilde{\phantom{w}}$}}}&}}
$$

On choisit ensuite un carré cocartésien:
\begin{equation}\label{cak}
\widetilde{\mathfrak{X}} \; = \; 
\vcenter{\xymatrix{
\underset{\phantom{.}}{QX_{a}}\,\,\ar@{>->}[r]^{\varphi}\ar@{>->}[d]_{\psi}&X_{b}'\,\ar[d]^-{\psi'}\\X_{d}'\ar[r]_-{\varphi'}&X_{d}'\underset{QX_{a}}{\sqcup}X_{b}'\,,}}
\end{equation}
qui par le Lemme \ref{fibrabra} est forcement homotopiquement cocartésien. La propriété universelle de la somme amalgamée, fournit un morphisme $F_{c}:\xymatrix@-5pt{X_{d}'\underset{QX_{a}}{\sqcup}X_{b}'\ar[r]&X_{c}}$ qui fait du cube \eqref{cubeproof} un morphisme de $\boxempty$-diagrammes $F:\xymatrix@-8pt{\widetilde{\mathfrak{X}}\ar[r]&\mathfrak{X}}$; \emph{i.e.} dont toutes les faces sont de carrés commutatifs.

Puisque $\widetilde{\mathfrak{X}}$ est un carré homotopiquement cocartésien, et par construction les flèches $F_{a}$, $F_{b}$ et $F_{d}$ sont d'équivalences faibles, il nous reste à vérifier que $F_{c}$ est une équivalence faible.

Effectivement, remarquons par ailleurs que le morphisme $\varphi$ de \eqref{cubeproof} est une équivalence faible, puisque par hypothèse $f$ l'est aussi. D'autre part, du fait que dans une catégorie de modèles, les morphismes qui sont à la fois des cofibrations et des équivalences faibles sont stables par cochangement de base, on déduit que la flèche $\varphi'$ du carré \eqref{cak} est en particulier une équivalence faible. Enfin, vu qu'on a supposé que dans le cube \eqref{cubeproof} la flèche $f'$ est une équivalence faible; on en déduit que $F_{c}$ l'est aussi.
\end{proof}

Considérons finalement la petite catégorie engendrée par le diagramme:
$$
\boxempty\!\!\boxempty \; = \;
\def\objectstyle{\scriptstyle}
\def\labelstyle{\scriptstyle}
\vcenter{\xymatrix@-10pt{a\ar[r]\ar[d] & b\ar[d]\ar[r] & c\ar[d]\\ f\ar[r] & e\ar[r] &d}}\,,
$$
soumis à la relation qu'identifie les deux flèches de $a$ vers $e$, les deux flèches de $b$ vers $d$ et toutes les flèches de $a$ vers $d$.

Il se suit que se donner un $\boxempty\!\!\boxempty$-diagramme $\mathfrak{X}$ de $\C$, équivaut à se donner un diagramme commutatif de $\C$ de la forme:
\begin{equation}\label{cX1}
\mathfrak{X} \; = \; \vcenter{\xymatrix@-10pt{X_{a}\ar[r]\ar[d] & X_{b}\ar[d]\ar[r] & X_{c}\ar[d]\\ X_{f}\ar[r] & X_{e}\ar[r] &X_{d}\,.}}
\end{equation}

Pour $1\leq i\leq 3$, notons $\alpha_{i}\xymatrix@C-5pt{\boxempty\ar[r]&\boxempty\!\!\boxempty}$ respectivement les foncteurs d'inclusion canoniques des sous-catégories:
$$
\def\objectstyle{\scriptstyle}
\def\labelstyle{\scriptstyle}
\vcenter{\xymatrix@-10pt{a\ar[r]\ar[d] &b \ar[d]\\f\ar[r]&e}}\,,\qquad
\def\objectstyle{\scriptstyle}
\def\labelstyle{\scriptstyle}
\vcenter{\xymatrix@-10pt{b\ar[r]\ar[d] &c\ar[d]\\e\ar[r]&d}}\qquad
\text{et}\qquad
\def\objectstyle{\scriptstyle}
\def\labelstyle{\scriptstyle}
\vcenter{\xymatrix@-10pt{a\ar[r]\ar[d] &c \ar[d]\\f\ar[r]&d}}.
$$

\begin{lemme}\label{exterio}
Soit $\mathfrak{X}$ un $\boxempty\!\!\boxempty$-diagramme de $\C$ et supposons que le carré $\alpha^*_{1}(\mathfrak{X})$ $\big($resp. $\alpha^*_{2}(\mathfrak{X})\Big)$ est homotopiquement cocartésien (resp. cartésien), alors le carré $\alpha^*_{2}(\mathfrak{X})$ $\big($resp. $\alpha^*_{1}(\mathfrak{X})\Big)$ est homotopiquement cocartésien (resp. cartésien) si et seulement si $\alpha^*_{3}(\mathfrak{X})$ est un carré homotopiquement cocartésien (resp. cartésien).
\end{lemme}
\begin{proof}
Montrons l'énoncé concernant les carrés homotopiquement cocartésiens. 

Pour commencer considérons les foncteurs $\tau_1\colon\xymatrix@C-10pt{I\ar[r]&J}$ et $\tau_2\colon\xymatrix@C-10pt{J\ar[r]&\boxempty\!\!\boxempty}$ définis par:
$$
\def\objectstyle{\scriptstyle}
\def\labelstyle{\scriptstyle}\vcenter{\xymatrix@-10pt{a\ar[r]\ar[d] &b \ar[r] & c \\f&&}}
\;
\def\objectstyle{\scriptstyle}
\def\labelstyle{\scriptstyle}
\vcenter{\xymatrix{\phantom{A}\ar@{^(->}[r]^-{\tau_{1}}&}}\;
\def\objectstyle{\scriptstyle}
\def\labelstyle{\scriptstyle}
\vcenter{\xymatrix@-10pt{a\ar[r]\ar[d] &b \ar[d]\ar[r] & c \\f\ar[r]&e&}}
\;
\def\objectstyle{\scriptstyle}
\def\labelstyle{\scriptstyle}
\vcenter{\xymatrix{\phantom{A}\ar@{^(->}[r]^-{\tau_{2}}&}}\;
\vcenter{\xymatrix@-10pt{a\ar[r]\ar[d] &b \ar[d]\ar[r] &c\ar[d]\\f\ar[r]&e\ar[r]&d}}\,,
$$
et notons $\tau_3\colon\xymatrix@C-10pt{I\ar[r]&\boxempty\!\!\boxempty}$ le composé $\tau_{2}\circ \tau_{1}$.

Si $1\leq i\leq 3$ supposons qu'on s'est donné une adjonction $(\tau_i)_{!} \dashv {\tau_i^*}$ munie d'une counité $\varepsilon^i\colon\xymatrix@-10pt{ (\tau_i)_!\circ \tau_i^* \ar@{=>}[r]&\mathrm{id}}$. Vu que $\tau_2\circ \tau_1=\tau_3$ on sait qu'il existe un isomorphisme de foncteurs:
$$
\xymatrix@+10pt{ \tau_{2\,!} \circ \tau_{1\,!} \ar@{=>}[r]^-{\Psi}_-{\cong}&\tau_{3\,!} }
$$ 
tel que pour tout $\boxempty\!\!\boxempty$-diagramme $\mathfrak{X}$ de $\C$: 
\begin{equation}\label{triang}
\xymatrix@R=12pt@C+25pt{
\tau_{2\,!} \circ \tau_{1\,!} \circ \tau_{1}^* \circ \tau_{2}^* (\mathfrak{X})
\ar@{}[d]|-{\text{\rotatebox[origin=c]{90}{\Large$\cong$}}}_-{\Psi_{\tau_{3}^* \, \mathfrak{X}}\phantom{a.}} 
\ar[r]^-{\tau_{2\,!} \big(\varepsilon^{1}_{ \tau_{2}^* (\mathfrak{X})}\big)} &
\tau_{2\,!} \circ \tau_{2}^* (\mathfrak{X})
\ar[r]^-{\varepsilon^{2}_{\mathfrak{X}}} & 
\mathfrak{X}\\
\tau_{3\,!} \circ \tau_{3}^* (\mathfrak{X})\ar@/_10pt/[rru]_{\varepsilon^3_{\mathfrak{X}}}&&}
\end{equation}
est un diagramme commutatif de $\C^{\boxempty\!\!\boxempty}[{\bf W}_{\boxempty\!\!\boxempty}^{-1}]$.

Donc pour montrer l'énoncé désiré il suffit de vérifier que pour tout $\boxempty\!\!\boxempty$-diagramme $\mathfrak{X}$ de $\C$ le carré $\alpha^*_{1}(\mathfrak{X})$ est homotopiquement cocartésien si et seulement si le morphisme:
$$
\xymatrix@C+15pt{
(\tau_1)_!\circ \tau_1^*\big( \tau_2^*\mathfrak{X}\big)_e\ar[r]^-{(\varepsilon_{\tau_2^*\mathfrak{X}}^{1})_e} &
(\tau_2^*\mathfrak{X})_e
}
$$
est un isomorphisme de $\C[{\bf W}^{-1}]$ (voir les Corollaires \ref{localco} et \ref{necimescimagg}), et si $2\leq i\leq 3$ que le carré $\alpha^*_{i}(\mathfrak{X})$ est homotopiquement cocartésien si et seulement si le morphisme:
$$
\xymatrix@C+15pt{
(\tau_i)_!\circ \tau_i^*\big( \mathfrak{X}\big)_d\ar[r]^-{(\varepsilon_{\mathfrak{X}}^{i})_d} & \mathfrak{X}_d
}
$$
est un isomorphisme de $\C[{\bf W}^{-1}]$.

Montrons que si $\mathfrak{X}$ est un $\boxempty\!\!\boxempty$-diagramme de $\C$, alors le carré $\alpha^*_{1}(\mathfrak{X})$ est homotopiquement cocartésien si et seulement si le morphisme:
$$
\xymatrix@C+15pt{
(\tau_1)_!\circ \tau_1^*\big( \tau_2^*\mathfrak{X}\big)_e\ar[r]^-{(\varepsilon_{\tau_2^*\mathfrak{X}}^{1})_e} &
(\tau_2^*\mathfrak{X})_e
}
$$
est un isomorphisme de $\C[{\bf W}^{-1}]$.

En effet, prenons le carré commutatif:
\begin{equation}\label{Fdoubcu1}
\xymatrix@+8pt{
\lefthalfcap \ar@{}[rd]|-{=}\ar[d]_-{q_{\lefthalfcap}}\ar[r]^-{\gamma_1}&I\ar[d]^-{\tau_1}\\ 
\boxempty \ar[r]_{\beta_1}& J }
\end{equation}
où $\beta_1$ et $\gamma_1$ sont les seuls foncteurs tels que $\tau_2\circ \beta_1=\alpha_1$ et $\tau_2\circ\tau_1\circ\gamma_1=\alpha_1\circ q_\lefthalfcap$ respectivement, et considérons une adjonction $(q_{\lefthalfcap})_{!} \dashv {q_{\lefthalfcap}^*}$ munie d'une counité $\varepsilon^{q_\lefthalfcap}\colon \xymatrix@-10pt{ (q_{\lefthalfcap})_!\circ q_{\lefthalfcap}^* \ar@{=>}[r]&\mathrm{id}}$.

Si on considère la conjuguée à gauche de \eqref{Fdoubcu1} par rapport aux transformations naturelles $\eta^{1}$ et $\varepsilon^{q_\lefthalfcap}$:
\begin{equation}\label{Fdoubcu2}
\vcenter{\xymatrix@+10pt{
\C^{\lefthalfcap}[{\bf W}_{\lefthalfcap}^{-1}]\ar[d]_-{(q_{\lefthalfcap})_!}&\C^{I}[{\bf W}_{I}^{-1}]\ar[d]^-{(\tau_1)_{!}}\ar[l]_-{{\gamma_1^*}}\\ 
\C^\boxempty[{\bf W}_{\boxempty}^{-1}]\urtwocell<\omit>{\phantom{aaaa}{(\mathrm{id})_!}}&\C^J[{\bf W}_{J}^{-1}]\ar[l]^-{{\beta_1^*}}
}}
\end{equation}
où $\eta^{1}$ est une unité de l'adjonction $(\tau_1)_{!} \dashv {\tau_1^*}$ telle que $\eta^1$ et $\varepsilon^1$ vérifient les identités triangulaires, on vérifie sans peine que pour tout $J$-diagramme $\mathfrak{Y}$ de $\C$ on a un diagramme commutatif de $\boxempty$-diagrammes de $\C$:
\begin{equation}\label{Fdoubcu3}
\vcenter{\xymatrix@R-7pt@C+10pt{
(q_{\lefthalfcap})_!\circ\gamma_1^*\circ\tau_1^* \, \mathfrak{Y} \ar[r]^-{(\mathrm{id}_!)_{\tau_1^*\mathfrak{Y}}} \ar@{=}[d] &
\beta_1^*\circ (\tau_1)_!\circ \tau_1^* \, \mathfrak{Y}\ar[d]^-{\beta_1^* (\varepsilon^{1}_{\mathfrak{Y}})}\\
(q_{\lefthalfcap})_!\circ q_{\lefthalfcap}^*\circ\beta_1^* \, \mathfrak{Y}\ar[r]_-{\varepsilon^{q_{\lefthalfcap}}_{\beta_1^*\mathfrak{Y}}} &\beta_1^*\mathfrak{Y}\,.
}}
\end{equation} 

Donc pour montrer que pour tout $\boxempty\!\!\boxempty$-diagramme $\mathfrak{X}$ de $\C$ le carré $\alpha^*_{1}(\mathfrak{X})$ est homotopiquement cocartésien si et seulement si le morphisme:
$$
\xymatrix@C+15pt{
(\tau_1)_!\circ \tau_1^*\big( \tau_2^*\mathfrak{X}\big)_e\ar[r]^-{(\varepsilon_{\tau_2^*\mathfrak{X}}^{1})_e} &
(\tau_2^*\mathfrak{X})_e
}
$$
est un isomorphisme de $\C[{\bf W}^{-1}]$, il suffit de montrer que le composé de la transformations \eqref{Fdoubcu2} et le foncteur d'évaluation $e^*\colon\xymatrix{\C^\boxempty[{\bf W}_{\boxempty}^{-1}]\ar[r] & \C[{\bf W}^{-1}]}$ est un isomorphisme.

Ceci est une conséquence du Lemme \ref{paste} et la Proposition \ref{local} vu qu'on a une égalité de transformations naturelles:
$$
\vcenter{\xymatrix@+8pt{
\lefthalfcap\,|\,e\drtwocell<\omit>{\phantom{aa}\Gamma\,|\,e}\ar[d]_-{p_{\lefthalfcap\,|\,e}}\ar[r]^-{\pi\,|\,e}&\lefthalfcap \ar@{}[rd]|-{=} \ar[d]|-{q_{\lefthalfcap}}\ar[r]^-{\gamma_1}&I\ar[d]^-{\tau_1}\\ 
e\ar[r]_-{e}&\boxempty\ar[r]_{\beta_1}& J} }
\;\qquad = \qquad\;
\vcenter{\xymatrix@+8pt{
\lefthalfcap\,|\,e\ar@{}[rd]|-{=}\ar[d]_-{p_{\lefthalfcap\,|\,e}}\ar[r]^-{\gamma_1\,|\,e}_{\cong}&I\,|\,e\ar[d]|-{p_{I\,|\,e}}\ar[r]^-{\pi\,|\,e}
\drtwocell<\omit>{\phantom{aa}\Gamma'\,|\,e} &I\ar[d]^-{\tau_1}\\ 
e\ar[r]_-{\mathrm{id}}&e\ar[r]_{e}& J }}\,,
$$
où $\gamma_1\,|\,e\,\colon \;\lefthalfcap\,|\,e \cong I\,|\,e$ est un isomorphisme de catégories.

Si $\mathfrak{X}$ est un $\boxempty\!\!\boxempty$-diagramme de $\C$, on montre de façon analogue que le carré $\alpha^*_{2}(\mathfrak{X})$ est homotopiquement cocartésien si et seulement si le morphisme:
$$
\xymatrix@C+15pt{
(\tau_2)_!\circ \tau_2^*\big( \mathfrak{X}\big)_d\ar[r]^-{(\varepsilon_{\mathfrak{X}}^{2})_d} & \mathfrak{X}_d
}
$$
est un isomorphisme de $\C[{\bf W}^{-1}]$, en considérant l'égalité des transformations naturelles:
$$
\vcenter{\xymatrix@+8pt{
\lefthalfcap\,|\,d\drtwocell<\omit>{\phantom{aa}\Gamma\,|\,d}\ar[d]_-{p_{\lefthalfcap\,|\,d}}\ar[r]^-{\pi\,|\,d}&\lefthalfcap \ar@{}[rd]|-{=} \ar[d]|-{q_{\lefthalfcap}}\ar[r]^-{\gamma_2}&J\ar[d]^-{\tau_2}\\ 
e\ar[r]_-{d}&\boxempty\ar[r]_{\alpha_2}&   \boxempty\!\!\boxempty   } }
\;\qquad = \qquad\;
\vcenter{\xymatrix@+8pt{
\lefthalfcap\,|\,d\ar@{}[rd]|-{=}\ar[d]_-{p_{\lefthalfcap\,|\,d}}\ar[r]^-{\gamma_2\,|\,d}&J\,|\,d\ar[d]|-{p_{J\,|\,d}}\ar[r]^-{\pi\,|\,d}
\drtwocell<\omit>{\phantom{aa}\Gamma'\,|\,d} &J\ar[d]^-{\tau_2}\\ 
e\ar[r]_-{\mathrm{id}}&e\ar[r]_{d}&     \boxempty\!\!\boxempty    }}\,,
$$
où $\gamma_2$ est le seul foncteur tel que $\tau_2\circ\gamma_2=\alpha_2\circ q_{\lefthalfcap}$. 

D'un ce cas il faut juste noter que le foncteur $\gamma_2\,|\,d$ n'est pas un isomorphisme mais il admet un adjoint à gauche (voir (ii) et 1 du Corollaire \ref{eqcisi}). 

Enfin si $\mathfrak{X}$ est un $\boxempty\!\!\boxempty$-diagramme de $\C$ pour montrer que le carré $\alpha^*_{3}(\mathfrak{X})$ est homotopiquement cocartésien si et seulement si le morphisme:
$$
\xymatrix@C+15pt{
(\tau_3)_!\circ \tau_3^*\big( \mathfrak{X}\big)_d\ar[r]^-{(\varepsilon_{\mathfrak{X}}^{3})_d} & \mathfrak{X}_d
}
$$
est un isomorphisme de $\C[{\bf W}^{-1}]$ on considère l'égalité des transformations naturelles:
$$
\vcenter{\xymatrix@+8pt{
\lefthalfcap\,|\,d\drtwocell<\omit>{\phantom{aa}\Gamma\,|\,d}\ar[d]_-{p_{\lefthalfcap\,|\,d}}\ar[r]^-{\pi\,|\,d}&\lefthalfcap \ar@{}[rd]|-{=} \ar[d]|-{q_{\lefthalfcap}}\ar[r]^-{\gamma_3}&I\ar[d]^-{\tau_3}\\ 
e\ar[r]_-{d}&\boxempty\ar[r]_{\alpha_3}&   \boxempty\!\!\boxempty   } }
\;\qquad = \qquad\;
\vcenter{\xymatrix@+8pt{
\lefthalfcap\,|\,d\ar@{}[rd]|-{=}\ar[d]_-{p_{\lefthalfcap\,|\,d}}\ar[r]^-{\gamma_3\,|\,d}&I\,|\,d\ar[d]|-{p_{I\,|\,d}}\ar[r]^-{\pi\,|\,d}
\drtwocell<\omit>{\phantom{aa}\Gamma'\,|\,d} &I\ar[d]^-{\tau_3}\\ 
e\ar[r]_-{\mathrm{id}}&e\ar[r]_{d}&     \boxempty\!\!\boxempty    }}\,,
$$
où $\gamma_3$ est le seul foncteur tel que $\tau_3\circ\gamma_3=\alpha_3\circ q_{\lefthalfcap}$, en particulier on vérifie que $\gamma_3\,|\,d$ admet un adjoint à gauche (voir (ii) et 1 du Corollaire \ref{eqcisi}). 
\end{proof}





\renewcommand{\thesubsection}{\S\thesection.\arabic{subsection}}
\subsection{}\;\label{propreg}
\renewcommand{\thesubsection}{\thesection.\arabic{subsection}}

Soit $\C$ une catégorie de modèles. Rappelons qu'on dit que $\C$ est \emph{propre à gauche} (resp. \emph{propre à droite}), si les équivalences faibles de $\C$ sont stables par cochangement de base (resp. par changement de base) le long des cofibrations (resp. fibrations); autrement dit, si l'on considère un carré cocartésien (resp. cartésien) de $\C$:
$$
\vcenter{\xymatrix{X_{a}\ar[r]^-{f}\ar[d]_{g'}&X_{b}\,\ar[d]^g\\X_{d}\ar[r]_{f'}&X_{c}\,,}}
$$
tel que $f$ soit une cofibration (resp. $g$ soit une fibration) et $g'$ une équivalence faible (resp. $f'$ une équivalence faible), alors la flèche $g$ (resp. $f$) est aussi une équivalence faible.

\begin{corollaire}
Si $\C$ est une catégorie de modèles propre à gauche (resp. droite), alors un carré commutatif de $\C$:
\begin{equation}\label{carreX}
\mathfrak{X} \; = \; \vcenter{\xymatrix{X_{a}\ar[r]^-{f}\ar[d]_{g'}&X_{b}\,\ar[d]^g\\X_{d}\ar[r]_{f'}&X_{c}\,,}}
\end{equation}
est homotopiquement cocartésien (resp. cartésien) si et seulement s'il satisfait la propriété suivante:
\begin{enumerate}
\item[{\bf CP}]
Pour toute factorisation:
\begin{equation}\label{fafa} 
\vcenter{\xymatrix@R-19pt{\underset{\phantom{a}}{X_{a}}\ar@/_5pt/@<-4pt>@{>->}[rd]_{i}\ar[rr]^f&&X_{b}\\&A\ar@<-4pt>@{->>}@/_3pt/[ru]_{p}\ar@{}[ru]|-{\text{\rotatebox[origin=c]{212}{\Large $\widetilde{\phantom{w}}$}}}&}}
\qquad\quad\left(\text{resp.}\quad
\vcenter{\xymatrix@R-19pt{\underset{\phantom{a}}{X_{b}}\ar@/_5pt/@<-4pt>@{>->}[rd]_{j}\ar@{}[rd]|-{\text{\rotatebox[origin=c]{135}{\Large $\widetilde{\phantom{w}}$}}}\ar[rr]^{g}&&X_{c}\\&A\ar@<-4pt>@{->>}@/_3pt/[ru]_{q}&}}
\right)\,
\end{equation}
le morphisme $\varphi$ dans le diagramme commutatif suivant: 
\begin{equation}\label{globo}  
\vcenter{\xymatrix@C+6pt@R-1pt{
X_{a}\,\,\ar[d]_-{g'}\ar@{>->}[r]^-{i}\ar@/^20pt/[rr]^-{f} &
A\ar[d]\ar@{->>}[r]^p & X_{b}\ar[d]^g\\
X_{d} \ar[r]\ar@/_20pt/[rr]_-{f'}&X_{d}\underset{X_{a}}{\sqcup} A \ar@{-->}[r]_{\varphi} & X_{c}}}
\qquad\;\left(\text{resp.}\quad
\vcenter{\xymatrix@C+6pt@R-1pt{
X_{a}\ar[r]^-{f}\ar@{-->}[d]^-{\varphi}\ar@/_32pt/[dd]_-{g'} & X_{\underset{\phantom{.}}{b}} \ar@{>->}[d]_-{j}\ar@/^23pt/[dd]^-{g} \\
X_{d}\underset{X_{c}}{\times} A\ar[r]\ar[d] & A\ar@{->>}[d]_-{q}\\
X_{d}\ar[r]_-{f'}& X_{c}}}
\right)
\end{equation}
est une équivalence faible.
\end{enumerate}
\end{corollaire}
\begin{proof}
On va montrer la partie de l'énoncé qui correspond aux carrés homotopiquement cocartésiens; on vérifie le cas dual de fa\c con analogue. 

Pour cela considérons un carré commutatif \eqref{carreX} de $\C$, et une factorisation \eqref{fafa} du morphisme $f$.
On va construire à partir de cette donnée un diagramme commutatif:
$$
\vcenter{\xymatrix@C-5pt@R-10pt{
&&X_{a}\ar'[d]'[dd][ddd]|-{g'}\ar[rrr]^-f&&& X_{b}\ar@{->>}[ddd]|-g\\
&  X_{a}\,\, \ar@{=}[ru] \ar@<-2.5pt>@{>->}[rrr]|-{i}\ar' [d] [ddd]|-{g'} &&&A_{\phantom{.}}\ar@<+2pt>[ddd]\ar[ru]^(.4){\text{\rotatebox[origin=c]{218}{\Large $\widetilde{\phantom{w}}$}}}_{p}&\\
\underset{\phantom{.}}{QX_{a}}\,\,\ar@{>->}[ddd]\ar@<+1pt>[ru]^{\text{\rotatebox[origin=c]{218}{\Large $\widetilde{\phantom{w}}$}}}\ar@{>->}@<-2pt>[rrr] &&& C_{\phantom{.}} \ar@<-2pt>[ru]^{\text{\rotatebox[origin=c]{218}{\Large $\widetilde{\phantom{w}}$}}}\ar@<+2pt>[ddd] &&\\ 
&& X_{d}\ar'[rr]'[rr][rrr]|-{f'}&&&X_{c} \,,\\
&X_{d}\ar'[rr][rrr]\ar@{=}[ur] &&&Z_{\phantom{.}}\ar@{-->}[ru]_{\varphi}\\
B_{\phantom{.}}\ar[ru]^-{\text{\rotatebox[origin=c]{218}{\Large $\widetilde{\phantom{w}}$}}}\ar[rrr]&&& W_{\phantom{.}}\ar@{-->}[ru]_-{\psi}&&}}
$$
par le procédé suivant: Le cube à l'arrière plan est simplement le diagramme \eqref{globo}, où on note $X_{d}\underset{X_{a}}{\sqcup}A \, = \,  Z$. Ensuite, on prend un remplacement cofibrant $\xymatrix@-10pt{QX_{a}\ar[r]^-{\text{\rotatebox[origin=c]{180}{\Large $\widetilde{\phantom{w}}$}}}&X_{a}}$ de l'objet $X_{a}$, et on considère des factorisations:
$$
\vcenter{\xymatrix@R-19pt{\underset{\phantom{a}}{QX_{a}}\ar@/_5pt/@<-4pt>@{>->}[rd]\ar[r]^-{\text{\rotatebox[origin=c]{180}{\Large $\widetilde{\phantom{w}}$}}}&X_{a}\ar[r]^-{g'}&X_{d}\\&B\ar@{}[ru]|-{\text{\rotatebox[origin=c]{212}{\Large $\widetilde{\phantom{w}}$}}}\ar@<-4pt>@{->>}@/_3pt/[ru]&}}
\qquad\;\text{et}\;\qquad
\vcenter{\xymatrix@R-19pt{\underset{\phantom{a}}{QX_{a}}\ar[r]^-{\text{\rotatebox[origin=c]{180}{\Large $\widetilde{\phantom{w}}$}}}\ar@/_5pt/@<-4pt>@{>->}[rd]&X_{a}\ar[r]&A\,.\\&C\ar@{}[ru]|-{\text{\rotatebox[origin=c]{212}{\Large $\widetilde{\phantom{w}}$}}}\ar@<-4pt>@{->>}@/_3pt/[ru]&}}
$$

Enfin, on définit l'objet $W$ comme étant la somme amalgamée des flèches $\vcenter{\xymatrix@-10pt{QX_{a}\ar[r]\ar[d]&C\\B&}}$; en particulier, on déduit le morphisme $\psi$, qui est forcement une équivalence faible car $\C$ est propre à gauche (voir par exemple la Proposition 13.2.4 de \cite{hirschhorn}).

On vérifie finalement que le carré  $\vcenter{\xymatrix@-10pt{QX_{a}\ar[r]\ar[d]&C\ar[d]\\B\ar[r]&W}}$ est un carré homotopiquement cocartésien grâce au lemme \ref{fibrabra}. Donc, il résulte du lemme \ref{cubelele} que $\mathfrak{X}$ est homotopiquement cocartésien, si et seulement si le composé $\psi\,\varphi$ est une équivalence faible, \emph{i.e.} si et seulement si $\varphi$ est une équivalence faible.
\end{proof}

\section{Foncteurs espace de lacets et suspension}

\renewcommand{\thesubsection}{\S\thesection.\arabic{subsection}}
\subsection{}\;\label{stapre}  
\renewcommand{\thesubsection}{\thesection.\arabic{subsection}}

Soit $\C$ une catégorie de modèles et considérons les foncteur suivants:
\begin{equation}\label{somega}
\def\objectstyle{\scriptstyle}
\def\labelstyle{\scriptstyle}
\vcenter{\xymatrix@-21pt{
&&&&&&& b\ar[dd]&&&&& a\ar[rr]\ar[dd]&&b \ar[dd]\\
e&\phantom{.}\ar@{^(->}[rrr]^-{t}&&& &&& &\phantom{.}\ar@{^(->}[rrr]^-{q_{\righthalfcup}}&&&&&&\\
&&&&&d\ar[rr]&& c& & &&& d\ar[rr] & & c}}
\qquad\text{et}\qquad
\def\objectstyle{\scriptstyle}
\def\labelstyle{\scriptstyle}
\vcenter{\xymatrix@-21pt{
&&&&&a\ar[rr]\ar[dd] && b &&&&& a\ar[rr]\ar[dd]&&b\ar[dd]\\
e&\phantom{.}\ar@{^(->}[rrr]^-{s}&&& &&& &\phantom{.}\ar@{^(->}[rrr]^-{q_{\lefthalfcap}}&&&&&&\\
&&&&&d&& & & &&& d\ar[rr] & & c}};
\end{equation}
où $\xymatrix@-10pt{e\ar[r]^-t&\righthalfcup}$ et  $\xymatrix@-10pt{e\ar[r]^-s&\lefthalfcap}$ désignent respectivement les foncteurs associés aux objets $c$ et $a$.

Supposons aussi qu'on s'est donné des foncteurs adjoints (ce qu'on peut faire toujours d'après le Corollaire \ref{hocompco}):
\begin{equation}\label{adj1}
\C[{\bf W}^{-1}]
\vcenter{\xymatrix@C+6pt{
\phantom{\cdot}\ar@{}[r]|-{\perp}\ar@<+4pt>@/^8pt/[r]^-{t_{!}}&
\phantom{\cdot}\ar@<+4pt>@/^8pt/[l]^-{{t^*}}}}
\C^{\righthalfcup}[{\bf W}_{\righthalfcup}^{-1}]
\vcenter{\xymatrix@C+6pt{
\phantom{\cdot}\ar@{}[r]|-{\perp}\ar@<-4pt>@/_8pt/[r]_-{q_{\righthalfcup\,*}}&
\phantom{\cdot}\ar@<-4pt>@/_8pt/[l]_-{{q_{\righthalfcup}^*}}}}
\C^{\boxempty}[{\bf W}_{\boxempty}^{-1}]
\end{equation}$$\text{et}$$\begin{equation}\label{adj2}
\C[{\bf W}^{-1}]
\vcenter{\xymatrix@C+6pt{
\phantom{\cdot}\ar@{}[r]|-{\perp}\ar@<-4pt>@/_8pt/[r]_-{s_{*}}&
\phantom{\cdot}\ar@<-4pt>@/_8pt/[l]_-{{s^*}}}}
\C^{\lefthalfcap}[{\bf W}_{\lefthalfcap}^{-1}]
\vcenter{\xymatrix@C+6pt{
\phantom{\cdot}\ar@{}[r]|-{\perp}\ar@<+4pt>@/^8pt/[r]^-{q_{\lefthalfcap\,!}}&
\phantom{\cdot}\ar@<+4pt>@/^8pt/[l]^-{{q_{\lefthalfcap}^*}}}}
\C^{\boxempty}[{\bf W}_{\boxempty}^{-1}].
\end{equation}

On définit \emph{le foncteur espace de lacets} de $\C$ (resp. le foncteur \emph{suspension} de $\C$) \emph{par rapport aux adjonctions} \eqref{adj1} (resp. \eqref{adj2}) par la formule suivante:
\begin{equation}\label{omesigma}
\vcenter{\xymatrix@R=1pt@C+10pt{\mathrm{Ho}(\C)\ar[r]^-{\Omega(\,\cdot\,)}&\mathrm{Ho}(\C)\\X\ar@{}[r]|-{\longmapsto} & (q_{\righthalfcup\,*}\,t_{!} \; X)_{a}}}
\qquad\qquad\;\left(\text{resp.}\;\qquad
\vcenter{\xymatrix@R=1pt@C+10pt{\mathrm{Ho}(\C)\ar[r]^-{\Sigma(\,\cdot\,)}&\mathrm{Ho}(\C)\,.\\ X\ar@{}[r]|-{\longmapsto} & (q_{\lefthalfcap\,!}\,s_{*} \; X)_{c}}}
\right)\,.
\end{equation}

Remarquons que de choix différentes des adjonctions \eqref{adj1} (resp. \eqref{adj2}) dérivent dans de foncteurs espace des lacets (resp. foncteurs suspensions) canoniquement isomorphes. \emph{Un foncteur espace de lacets} (resp. \emph{un foncteur suspension}) de la catégorie de modèles $\C$ est un foncteur $\xymatrix@C-10pt{\mathrm{Ho}(\C)\ar[r]&\mathrm{Ho}(\C)}$ isomorphe à un foncteur de la forme \eqref{omesigma}. 

\begin{lemme} \label{isoomega}
Soit $\C$ une catégorie de modèles et $\mathcal{A}$ une sous-catégorie pleine de $\C$. Supposons que on a un foncteur:
\begin{align*}
&\;\quad \xymatrix@C+80pt{ \mathcal{A} \ar[r]^-{\Phi} & \C^{\boxempty}  } \\
&\vcenter{\def\objectstyle{\scriptstyle}\def\labelstyle{\scriptstyle}\xymatrix@-20pt{&&Y\\X\ar[rru]^-{\,f\,}&&}}\quad \longmapsto \quad
\vcenter{\def\objectstyle{\scriptstyle}\def\labelstyle{\scriptstyle}\xymatrix@-20pt{
&&Y_{a}\ar[rrr]\ar'[d][ddd]&&&Y_{b}\,\ar[ddd]\\
X_{a}\ar[rru]^-{\,f_{a}\,}\ar@<-2pt>[rrr]\ar[ddd]&&&X_{b}\ar[ddd]\ar[rru]|-{\,f_{b}\,}&&\\
&&&&&\\
&&Y_{d}\ar'[r][rrr]&&&Y_{c}\\
X_{d}\ar[rrr]\ar[urr]^-{\,f_{d}\,}&&&X_{c}\ar[rru]_-{\,f_{c}\,}}}
\end{align*}
vérifiant les propriétés:
\begin{enumerate}
\item Pour tout objet $X$ de $\mathcal{A}$ on a que $X_c=X$ (resp. $X_a=X$), les morphismes $\xymatrix@-12pt{X_b\ar[r]&\star&X_d\ar[l]}$ sont des équivalences faibles et le carré:
$$\vcenter{\def\objectstyle{\scriptstyle}\def\labelstyle{\scriptstyle}\xymatrix@-8pt{X_{a}\ar[r]\ar[d]&X_{b}\,\ar[d]\\X_{d}\ar[r]&X_{c}\,.}}$$
est homotopiquement cartésien (resp. homotopiquement cocartésien).
\item Pour tout morphisme $f$ de $\mathcal{A}$ on a que $f_c=f$ (resp. $f_a=f$).
\end{enumerate}

Alors le foncteur: 
$$
\vcenter{\xymatrix@C+45pt@R=1pt{
\A\ar[r]^-{\Phi_a}&\C\\
X\ar@{}[r]|-{\longmapsto} & X_a
}}
\qquad\qquad\left(\text{resp.}\qquad
\vcenter{\xymatrix@C+45pt@R=1pt{
\A\ar[r]^-{\Phi_c}&\C\\
X\ar@{}[r]|-{\longmapsto} & X_c
}}\;\right)
$$
respect les équivalences faibles et le foncteur induit:
$$
\vcenter{\xymatrix@C+10pt{\A\big[\big({\bf W}\cap \mathcal{A}\big)^{-1}\big] \ar[r]^-{\widetilde{\Phi_a}} &\C\big[{\bf W}^{-1}\big]}}
\qquad\qquad\bigg(\text{resp.}\qquad 
\vcenter{\xymatrix@C+10pt{\A\big[\big({\bf W}\cap \mathcal{A}\big)^{-1}\big] \ar[r]^-{\widetilde{\Phi_c}} &\C\big[{\bf W}^{-1}\big]}}
\bigg)
$$
est isomorphe au foncteur composé:
$$
\vcenter{\xymatrix@+3pt{\A\big[\big({\bf W}\cap \mathcal{A}\big)^{-1}\big]  \ar[r] & \C\big[{\bf W}^{-1}\big] \ar[r]^-{\Omega} & \C\big[{\bf W}^{-1}\big]\,,}}
$$
$$
\bigg(\text{resp.}\qquad 
\vcenter{\xymatrix@+3pt{\A\big[\big({\bf W}\cap \mathcal{A}\big)^{-1}\big]  \ar[r] & \C\big[{\bf W}^{-1}\big] \ar[r]^-{\Sigma} & \C\big[{\bf W}^{-1}\big]\,,}}
\bigg)
$$
où $\Omega$ (resp. $\Sigma$) est un foncteur espace de lacets (resp. un foncteur suspension) de $\C$.
\end{lemme}
\begin{proof}
Montrons le résultat concernant le foncteur espace de lacets; l'autre cas est analogue. Pour commencer remarquons que le foncteur: 
$$
\vcenter{\xymatrix@C+45pt@R=1pt{
\A\ar[r]^-{\Phi_a}&\C\\
X\ar@{}[r]|-{\longmapsto} & X_a
}}
$$
respect les équivalences faibles d'après le Lemme \ref{cubelele}. 

D'un autre, pour montrer que le foncteur induit:
$$
\vcenter{\xymatrix@C+10pt{\A\big[\big({\bf W}\cap \mathcal{A}\big)^{-1}\big] \ar[r]^-{\widetilde{\Phi_a}} &\C\big[{\bf W}^{-1}\big]}}
$$
est isomorphe au foncteur composé:
$$
\vcenter{\xymatrix@+3pt{\A\big[\big({\bf W}\cap \mathcal{A}\big)^{-1}\big]  \ar[r] & \C\big[{\bf W}^{-1}\big] \ar[r]^-{\Omega} & \C\big[{\bf W}^{-1}\big]}}
$$
où $\Omega$ est un foncteur espace de lacets de $\C$, il se suit des propriétés des catégories de fractions qu'il suffit de montrer qu'on a un isomorphisme de foncteurs:
$$
\xymatrix@C+11pt@R-8pt{ 
\underset{\phantom{a}}{\A} \ar@{^(->}[d]\ar[r]^-{\Phi_a} & \C \ar[dd]^-{\gamma} \\ 
\C\ar[d]_-{\gamma} \ar[d]_-{\gamma}& \\
\C[{\bf W}^{-1}] \ar[r]_{\Omega} & \C[{\bf W}^{-1}] \uultwocell<\omit>{\alpha}}
$$

Pour définir $\alpha$ remarquons que d'après la définition des foncteurs \eqref{somega} on a pour tout objet $X$ de $\A$ que ${t^*}\,\circ\,{q_{\righthalfcup}^*} \, \big(\Phi(X)\big)=X_{c}=X$. On obtient par adjonction un morphisme de la catégorie $\C^{\righthalfcup}[{\bf W}_{\righthalfcup}^{-1}]$:
\begin{equation}\label{1natutu}
\xymatrix{ t_{!} (X) \ar[r]^-{\ell_X} & {q_{\righthalfcup}^*}\big(\Phi(X)\big) }
\end{equation}
lequel selon la Proposition \ref{local} et le Lemme \ref{cwgen} est un isomorphisme, parce que par hypothèse les morphismes $\xymatrix@-12pt{X_b\ar[r]&\star&X_d\ar[l]}$ sont des équivalences faibles de $\C$ et $X=X_c$.

D'un autre côté on a que la flèche canonique est un isomorphisme:
\begin{equation} \label{2natutu}
\xymatrix@+5pt{
\Phi(X)  \ar[r]^-{m_X}&q_{\righthalfcup\, *} \,\circ\, {q_{\righthalfcup}^*} \,\big(\Phi(X)\big)\,,} 
\end{equation}
parce que $\Phi(X)$ est un carré homotopiquement cartésien pour tout objet $X$ de $\A$ (voir le Lemme \ref{carrlele}).

On définit l'isomorphisme $\alpha_X$ de la catégorie $\C[{\bf W}^{-1}]$ par le composé: 
$$
\xymatrix@+5pt{
\Omega(X) = q_{\righthalfcup\, *}\,\circ\,  t_{!} (X)_a
\ar[rr]^-{ q_{\righthalfcup\, *}(\ell_X)_{a}}&&
\big(q_{\righthalfcup\, !} \,\circ\, {q_{\righthalfcup}^*} \,\big(\Phi(X)\big)\,\big)_{a} & \big(\Phi(X)\big)_a = X_{a}\ar[l]_-{(m_X)_a}\,.}
$$

Vu que les morphismes \eqref{1natutu} et \eqref{2natutu} sont naturelles en $X$ on obtient bien une transformation naturelle $\alpha$.
\end{proof}

Comme une application du Lemme \ref{isoomega} qu'on vient de montrer, remarquons que si $\C$ est une catégorie de modèles vérifiant la propriété que le morphisme canonique $\xymatrix@-8pt{\emptyset\ar[r]&\star}$ soit une fibration (resp. une cofibration), par exemple si $\C$ est une catégorie pointée \emph{i.e.} $\emptyset\cong\star$, alors la définition qu'on a donnée des foncteurs $\Omega$ et $\Sigma$ coïncide avec la définition plus habituelle. 

En effet, si $X$ est un objet fibrant (resp. cofibrant) de $\C$ considérons une factorisation: 
$$
\vcenter{\xymatrix@R-19pt{\underset{\phantom{a}}{\emptyset}\ar@/_5pt/@<-4pt>[rd]^-{\text{\rotatebox[origin=c]{145}{\Large $\widetilde{\phantom{w}}$}}}\ar[rr]&&X\\&PX\ar@<-4pt>@{->>}@/_3pt/[ru]&}}
\qquad\quad
\left(\text{resp.}\quad
\vcenter{\xymatrix@R-19pt{\underset{\phantom{a}}{X}\ar@/_5pt/@<-4pt>@{>->}[rd]\ar[rr]&&\star\\&CX\ar@<-4pt>@/_3pt/[ru]^-{\text{\rotatebox[origin=c]{215}{\Large $\widetilde{\phantom{w}}$}}}&}}
\right)\,.
$$

On déduit du Lemme \ref{fibrabra} qu'un carré cartésien (resp. cocartésien) de $\C$:
$$
\vcenter{\xymatrix{\emptyset\underset{X}{\times}PX\ar[r]\ar[d]&PX\ar@{->>}[d]\\\emptyset\ar[r]&X}}
\qquad\quad
\left(\text{resp.}\qquad
\vcenter{\xymatrix{X\,\,\ar@{>->}[r]\ar[d]&CX\ar[d]\\\star\ar[r]& \star\underset{X}{\sqcup}CX}}\right)\,
$$
est un carré homotopiquement cartésien (resp. cocartésien) vu que par hypothèse $\emptyset$ est un objet fibrant (resp. $\star$ est un objet cofibrant) de $\C$.

Si $\C_{fib}$ (resp. $\C_{cofib}$) note la sous-catégorie pleine de $\C$ dont les objets sont les objets fibrants (resp. cofibrants) on peut définir par cette procédure un foncteur:
\begin{equation}\label{lacccclass}
\vcenter{\xymatrix@C+10pt@R=1pt{
\C_{fib} \ar[r] & \C_{fib}\\
X\ar@{}[r]|-{\longmapsto} & \emptyset\underset{X}{\times} PX
}}\qquad\qquad\Bigg(\text{resp.} \qquad 
\vcenter{\xymatrix@C+10pt@R=1pt{
\C_{cofib} \ar[r] & \C_{cofib}\\
X\ar@{}[r]|-{\longmapsto} & \star\underset{X}{\sqcup} CX
}}\Bigg)
\end{equation}
qui respect les équivalences faibles. D'après le Lemme \ref{isoomega} le foncteur induit:
$$
\vcenter{\xymatrix@C+10pt@R=1pt{
Ho \big(\C\big) \cong \C_{fib}\big[{\bf W}^{-1}\big] \ar[r] & \C_{fib}\big[{\bf W}^{-1}\big] \cong Ho \big(\C\big)
}}
$$
$$
\left(\text{resp.}\qquad\quad
\vcenter{\xymatrix@C+10pt@R=1pt{
Ho \big(\C\big) \cong \C_{cofib}\big[{\bf W}^{-1}\big]  \ar[r] & \C_{cofib}\big[{\bf W}^{-1}\big]  \cong Ho \big(\C\big)
}}\right)
$$
est un foncteur espace de lacets (resp. suspension) de $\C$.

Montrons:

\begin{lemme}\label{quillenomegasus}
Supposons que on a une adjonction de Quillen:
$$
\C\,
\vcenter{\xymatrix@C+15pt{
\phantom{a}\ar@{}[r]|{\perp}\ar@<-4pt>@/_10pt/[r]_{G}&\phantom{a}
\ar@<-4pt>@/_10pt/[l]_-{F}}}
\D\,.
$$
entre catégorie de modèles pointées \emph{i.e.} telles que le morphisme canonique $\xymatrix@-10pt{\emptyset\ar[r]&\star}$ soit un isomorphisme. Alors si $\Omega^\C$ et $\Omega^\D$ (resp. $\Sigma^\C$ et $\Sigma^\D$) sont de foncteurs espaces de lacets (resp. de foncteurs suspension) de $\C$ et $\D$ respectivement, le carré suivant commute á isomorphisme près:
$$
\vcenter{\xymatrix@C+5pt@R+2pt{
\mathrm{Ho}(\C)\ar[r]^-{{\bf R}G} \ar@{}[rd]|-\cong\ar[d]_-{\Omega^\C}&\mathrm{Ho}(\D)\ar[d]^-{\Omega^D}\\ 
\mathrm{Ho}(\C)\ar[r]_-{{\bf R}G}&\mathrm{Ho}(\D)
}}
\qquad\quad\qquad\left(\,\text{resp.}\qquad
\vcenter{\xymatrix@C+5pt@R+2pt{
\mathrm{Ho}(\C)\ar@{}[rd]|-\cong\ar[d]_-{\Sigma^\C}&\mathrm{Ho}(\D)\ar[d]^-{\Sigma^D}\ar[l]_-{{\bf L}F}\\ 
\mathrm{Ho}(\C)&\mathrm{Ho}(\D)\ar[l]^-{{\bf L}F}
}}
\right)
$$
où ${\bf R}G$ (resp. ${\bf L}F$) est un foncteur dérivé total à droite de $G$ (resp. à gauche de $F$).
\end{lemme}
\begin{proof}
Montrons le résultat concernant les foncteurs espaces de lacets; le cas des foncteurs suspension est analogue.  Pour cela on va construire une transformation naturelle :
\begin{equation}\label{lecacarreOS}
\vcenter{\xymatrix@C+5pt@R+2pt{
\C_{fib}\ar[r]^-{G} \drtwocell<\omit>{\alpha}\ar[d]_-{\eqref{lacccclass}}&\D_{fib}\ar[d]^-{\eqref{lacccclass}}\\ 
\C_{fib}\ar[r]_-{G}&\D_{fib}
}}
\end{equation}
telle que $\alpha_X$ soit une équivalence faible de $\D$ pour tout objet fibrant $X$ de $\C$.

Remarquons que si on se donne de factorisations fonctorielles dans $\C$ et $\D$:
$$
\vcenter{\xymatrix{X}}\qquad\longmapsto\qquad
\left(\vcenter{\xymatrix@R-25pt@C-10pt{\underset{\phantom{a}}{\emptyset^\C}\ar@/_5pt/@<-4pt>[rd]_-{\text{\rotatebox[origin=c]{45}{\LARGE $\wr$}}}\ar[rr]&&X\\&P^{\C}X\ar@<-4pt>@{->>}@/_3pt/[ru]&}}\right)\qquad \text{et}\qquad
\vcenter{\xymatrix{A}}\quad\longmapsto\quad
\left(\vcenter{\xymatrix@R-25pt@C-10pt{\underset{\phantom{a}}{\emptyset^\D}\ar@/_5pt/@<-4pt>[rd]_-{\text{\rotatebox[origin=c]{45}{\LARGE $\wr$}}}\ar[rr]&&X\\&P^{\D}A\ar@<-4pt>@{->>}@/_3pt/[ru]&}}\right)
$$
respectivement, ainsi que des foncteurs produits fibrés:
$$
\left(\vcenter{\def\objectstyle{\scriptstyle}\def\labelstyle{\scriptstyle}\xymatrix@-10pt{
& C\ar[d]\\ A\ar[r] & B }}\right) \qquad \longmapsto \qquad 
\left(\vcenter{\def\objectstyle{\scriptstyle}\def\labelstyle{\scriptstyle}\xymatrix@-10pt{
A\underset{B}{\times}C\ar[r]\ar[d]& C\ar[d] \\ A\ar[r] & B }}\right)\,;
$$
alors l'image d'un objet fibrant $X$ de $\C$ par le foncteur composé $G\circ \eqref{lacccclass}$ est l'objet $G\big(\emptyset^\C\underset{X}{\times}P^\C X \big)$ dans le carré cartésien:
\begin{equation}\label{ccrree1}
\xymatrix@+3pt{
G\big(\emptyset^\C\underset{X}{\times}P^\C X \big)\ar[r]\ar[d]& G\big(P^\C X\big)     \ar[d]\\
G(\emptyset^\C) \ar[r]&GX}
\end{equation}
et l'image de $X$ par le composé $\eqref{lacccclass}\circ G$ est l'objet $\emptyset^\D \underset{GX}{\times} P^\D (GX)$ dans le carré cartésien:
\begin{equation}\label{ccrree2}
\xymatrix@+3pt{
\emptyset^\D \underset{GX}{\times} P^\D (GX)\ar[r]\ar[d]&P^\D(GX)\ar[d]\\
\emptyset^\D\ar[r]&GX\,.}
\end{equation}

On définit le morphisme $\alpha_X$ en complétant le morphisme de carrés cartésiens:
$$
\xymatrix@-20pt{
&&G\big(\emptyset^\C\underset{X}{\times}P^\C X \big)\ar[rrr]\ar'[d][ddd]&&& G\big(P^\C X\big)     \ar[ddd]\\
\emptyset^\D \underset{GX}{\times} P^\D (GX)\ar@{-->}[rru]^-{\,\alpha_X\,}\ar@<-2pt>[rrr]\ar[ddd]&&&P^\D(GX)\ar[ddd]\ar[rru]|-{}&&\\
&&&&&\\
&&G(\emptyset^\C) \ar'[r][rrr]&&&GX\ar@{=}[lld]\\
\emptyset^\D\ar[rrr]\ar[urr]^-{}&&&GX}
$$
où $\xymatrix@-10pt{\emptyset^{\D} \ar[r] & G(\emptyset^{\C})}$ est le seul morphisme possible (lequel est un isomorphisme d'après les hypothèses) et $\xymatrix@-5pt{P^\D(GX)\ar[r] & G\big(P^\C X\big)}$ est le composé:
$$
\xymatrix@+5pt{P^\D(GX)\ar[r] & \star^\D \cong \emptyset^\D \ar[r] & G\big(P^\C X\big)}
$$
(lequel est une équivalence faible de $\D$).

Il se suit que $\alpha_X$ ainsi défini est une équivalence faible de $\D$ parce que les carrés \eqref{ccrree1} et \eqref{ccrree2} sont en fait des carrés homotopiquement cartésiens.
\end{proof}

Remarquons aussi:

\begin{proposition}\label{adjsigma}
Soit $\C$ est une catégorie de modèles pointée, \emph{i.e.} telle que le morphisme canonique $\xymatrix@-10pt{\emptyset\ar[r]&\star}$ soit un isomorphisme. Alors si $\Omega$ et $\Sigma$ sont un foncteur espace de lacets et un foncteur suspension de $\C$, on a un adjonction:
$$
\mathrm{Ho}(\C)\,
\vcenter{\xymatrix@C+15pt{
\phantom{a}\ar@{}[r]|{\perp}\ar@<-4pt>@/_10pt/[r]_{\Omega}&\phantom{a}
\ar@<-4pt>@/_10pt/[l]_-{\Sigma}}}
\mathrm{Ho}(\C)\,.
$$ 

\end{proposition}
\begin{proof}
Considérons de foncteurs $\Omega$ et $\Sigma$ définis à partir des adjonctions  \eqref{adj1} et \eqref{adj2}. Si $X$ et $Y$ sont des objets de $\C$ remarquons que par adjonction on a d'isomorphismes naturels:
\begin{align*}
\mathrm{Hom}_{\C^{\righthalfcup}[{\bf W}_{\righthalfcup}^{-1}]}\Big(
{q_{\righthalfcup}^*} \, q_{\lefthalfcap\,!}\,s_{*} \; X \;,\, t_{!} \; Y\Big)
&\;\cong\;
\mathrm{Hom}_{\C^{\boxempty}[{\bf W}_{\boxempty}^{-1}]}\Big(
q_{\lefthalfcap\,!}\,s_{*} \; X \;,\, q_{\righthalfcup\,*}\,t_{!} \; Y\Big)\\
&\;\cong\;
\mathrm{Hom}_{\C^{\lefthalfcap}[{\bf W}_{\lefthalfcap}^{-1}]}\Big(
s_{*} \; X \;,\, {q_{\lefthalfcap}^*}\, q_{\righthalfcup\,*}\,t_{!} \; Y\Big)\,.
\end{align*}

D'un autre part, les égalités:
$$
\Sigma\; X \, = \, {t^*} \; {q_{\righthalfcup}^*} \; q_{\lefthalfcap\,!}\,\,s_{*} \; X  \qquad\text{et}\qquad 
{s^*} \; {q_{\lefthalfcap}^*}\; q_{\righthalfcup\,*}\,\,t_{!} \; Y \, = \, \Omega\; Y;
$$ 
induisent par adjonction des morphismes:
$$
\xymatrix{ t_{!} \; \Sigma\; X\ar[r] & {q_{\righthalfcup}^*} \; q_{\lefthalfcap\,!}\,\,s_{*} \; X}\qquad\text{et}\qquad 
\xymatrix{{q_{\lefthalfcap}^*}\; q_{\righthalfcup\,*}\,\,t_{!} \; Y\ar[r]&  s_{*} \; \Omega\; Y}\,,
$$ 
lesquels sont des isomorphismes d'après la Proposition \ref{local} et le Lemme \ref{cwgen}, vu que par hypothèse $\star\cong\emptyset$.

On obtient ainsi que:
$$
\mathrm{Hom}_{\C^{\righthalfcup}[{\bf W}_{\righthalfcup}^{-1}]}\Big(t_{!} \; \Sigma\; X \;,\, t_{!} \; Y\Big)
\;\cong\; \mathrm{Hom}_{\C^{\lefthalfcap}[{\bf W}_{\lefthalfcap}^{-1}]}\Big(s_{*} \; X \;,\,  s_{*} \; \Omega\; Y\Big)\,.
$$

Donc, $\Sigma$ est adjoint à gauche de $\Omega$, parce que les foncteurs $t_{!}$ et $s_{*}$ sont pleinement fidèles d'après le Lemme \ref{localco}. 
\end{proof}


\section{Les $K$-groupes des catégories modèles}


\renewcommand{\thesubsection}{\S\thesection.\arabic{subsection}}
\subsection{}\;
\renewcommand{\thesubsection}{\thesection.\arabic{subsection}}


Si $p\geq 0$ posons $\text{\Large $\righthalfcup$}_{p}$ pour noter la sous-catégorie pleine de la catégorie produit $\mathbb{Z}\times\mathbb{Z}$ engendrée par les flèches du diagramme:
\begin{equation}\label{tringrosse}
\def\objectstyle{\scriptstyle}
\def\labelstyle{\scriptstyle}
\vcenter{\xymatrix@C-1.2pc@R-1.8pc{
   &&&&&&\\
   &&&&&(p,p)\\
   &&&&&\\
   &&&&\ar@{.}[r]&\ar[uu]\phantom{.}\\&&&&&&&\\
   &&&\ar@{.}[r]&\ar@{.}[r]\ar@{.}[uu]&\ar@{.}[uu]\\
   &&&&&\\
   &&(1,1)\ar[r]&\ar@{.}[r]\ar@{.}[uu]&\ar@{.}[uu]\ar[r]&(np,1)\ar[uu]\\
   &&&&&\\
   &(0,0)\ar[r]&(1,0)\ar[uu]\ar[r]&\ar@{.}[r]\ar@{.}[uu]&\ar@{.}[uu]\ar[r]&(p,0)\ar[uu]\\
   &&&&&&}}.
\end{equation}

Remarquons qu'on a un foncteur:
\begin{equation}\label{lefoncttripe}
\xymatrix@R=1pt@C+5pt{ \Delta \ar[r]^-{\text{\Large $\righthalfcup$}_{?}} & {\bf cat} \\ 
\text{\small $[p]$} \ar@{}[r]|-{\longmapsto}  & \text{\Large $\righthalfcup$}_{p}}
\end{equation}
défini dans une flèche $\xymatrix@C-14pt{[p]\ar[r]^{f}& [p']}$ de $\Delta$ par le foncteur:
$$
\xymatrix@R=1pt@C+5pt{
\text{\Large $\righthalfcup$}_{p}\ar[r]^-{\text{\Large $\righthalfcup$}_{f}} & \text{\Large $\righthalfcup$}_{p'} \\ 
(i,j) \ar@{}[r]|-{\longmapsto} & \big(f(i),f(j)\big)\,;
}
$$
donc on obtient un foncteur:
\begin{equation}\label{triangleburdes}
\xymatrix@R=1pt@C+5pt{ \Delta^{op} \ar[r] & {\bf cat} \\ \text{\small $[p]$} \ar@{}[r]|-{\longmapsto}  & \C^{\righthalfcup_{p}}}
\end{equation}
pour toute catégorie $\C$. 

Si $p\geq 0$ et $\C$ est une catégorie de modèles munie d'un objet nul $0$ fixe et d'une famille distinguée $\mathcal{A}$ d'objets de $\C$ contenant l'objet $0$, un $(p-1)$-\emph{triangle à gauche} (resp. \emph{à droite}) $\mathfrak{X}$ \emph{de $\C$} (\emph{à valeurs dans $\A$}) est par définition un $\text{\Large $\righthalfcup$}_{p}$-diagramme $\mathfrak{X}$ de $\C$:
$$
\vcenter{\xymatrix@C+15pt@R=3pt{
\text{\Large $\righthalfcup$}_{q} \ar[r]^-{\mathfrak{X}} & \C_\mathcal{A}\\
(i,j) \ar@{}[r]|-{\longmapsto} & X_{i,j}
}}
$$
vérifiant les propriétés suivantes:
\begin{enumerate}
\item $X_{i,j}$ appartient à $\A$ pour tout objet $(i,j)$ de $\text{\Large $\righthalfcup$}_{p}$ \emph{i.e.} toujours que $0\leq j \leq  i \leq p$.
\item $X_{k,k}= 0$ pour tout $0\leq k\leq p$.
\item Si $0\leq m < k\leq p-1$ le carré  de $\C$:
$$
 \vcenter{\xymatrix@-8pt{X_{k,m+1}\ar[r]&X_{k+1,m+1}\\X_{k,m}\ar[u]\ar[r]&X_{k+1,m}\ar[u]}}
$$ 
 induit par le foncteur $\mathfrak{X}$ est homotopiquement cocartésien (resp. homotopiquement cartésien). 
\end{enumerate}

Il se suit de la définition qu'il y a un seul $(-1)$-triangles à gauche (resp. à droite) de $\C$ qu'on identifie avec l'objet $0$. D'un autre, se donner un $0$-triangle à gauche (resp. à droite) de $\C$ $\def\objectstyle{\scriptstyle}\def\labelstyle{\scriptstyle}\vcenter{\xymatrix@-5pt{&0\\0\ar[r]&\ar[u]X}}$ équivaut à se donner un objet $X$ de $\C$ appartenant à la famille $\A$; et se donner un $1$-triangle à gauche (resp. à droite):
$$
\def\objectstyle{\scriptstyle}
\def\labelstyle{\scriptstyle}
\vcenter{\xymatrix@-5pt{&&0\\&0\ar[r]&Z\ar[u] \\0\ar[r]&X\ar[r]\ar[u]&Y\ar[u]}}\,,
$$
équivaut à se donner une suite de morphismes de $\C$:
$$
\xymatrix@C-6pt{0\ar[r] & X \ar[r]^-f & Y \ar[r]^-g & Z\ar[r] & 0}
$$
telle que les objets $X$, $Y$ et $Z$ appartient à $\A$ et:
$$
\def\objectstyle{\scriptstyle}
\def\labelstyle{\scriptstyle}
\vcenter{\xymatrix@-5pt{0\ar[r]&Z \\ X\ar[r]_-f\ar[u]&Y\ar[u]_-g}}\,.
$$
soit un carré homotopiquement cocartésien (resp. homotopiquement cartésien) de $\C$.

Si $p\geq 0$ on note $\mathrm{s}^g_{p}(\C,\A)$ $\big($resp. $\mathrm{s}^d_{p}(\C,\A)\big)$ l'ensemble des $(p-1)$-triangles à gauche (resp. à droite) de $\C$ à valeurs dans $\A$, et $\mathrm{S}^g_{p}(\C,\A)$ $\big($resp. $\mathrm{S}^d_{p}(\C,\A)\big)$ la sous-catégorie pleine de la catégorie des diagrammes $\C^{{\righthalfcup}_{p}}$ dont l'ensemble d'objets est $\mathrm{s}^g_{p}(\C,\A)$ $\big($resp. $\mathrm{s}^d_{p}(\C,\A)\big)$.

Posons aussi ${\bf W}\mathrm{S}^{g}_{p}(\C,\mathcal{A})$ $\big($resp. ${\bf W}\mathrm{S}^d_{p}(\C,\A)\big)$ pour noter la sous-catégorie de $\mathrm{S}^{g}_{p}(\C,\A)$ $\big($resp. $\mathrm{S}^d_{p}(\C,\A)\big)$ du même ensemble d'objets dont les morphismes sont les morphismes de $\text{\Large $\righthalfcup$}_{p}$-diagrammes qui sont des équivalences faibles argument par argument. 

Remarquons que d'après le Lemme \ref{exterio} le foncteur \eqref{triangleburdes} induit des foncteurs:
\begin{equation}\label{sss}
\vcenter{\xymatrix@C+8pt@R=1pt{
\Delta^{op} \ar[r]^-{\mathrm{s}^g_{\bullet}(\C,\mathcal{A})} & \ens\\ [p]\ar@{}[r]|-{\longmapsto} & \mathrm{s}^g_{p}(\C,\mathcal{A})}}
\qquad\quad
\left(\text{resp.}\quad
\vcenter{\xymatrix@C+8pt@R=1pt{
\Delta^{op} \ar[r]^-{\mathrm{s}^d_{\bullet}(\C,\mathcal{A})} & \ens\\ [p]\ar@{}[r]|-{\longmapsto} & \mathrm{s}^d_{p}(\C,\mathcal{A})}}\right)\,,
\end{equation}
\begin{equation}\label{ssS}
\vcenter{\xymatrix@C+8pt@R=1pt{
\Delta^{op} \ar[r]^-{\mathrm{S}^g_{\bullet}(\C,\mathcal{A})} & \cat\\ [p]\ar@{}[r]|-{\longmapsto} & \mathrm{S}^g_{p}(\C,\mathcal{A})}}
\qquad\quad
\left(\text{resp.}\quad
\vcenter{\xymatrix@C+8pt@R=1pt{
\Delta^{op} \ar[r]^-{\mathrm{S}^d_{\bullet}(\C,\mathcal{A})} & \cat\\ [p]\ar@{}[r]|-{\longmapsto} & \mathrm{S}^d_{p}(\C,\mathcal{A})}}\right)
\end{equation}
$$\text{et}$$
\begin{equation}\label{ssW}
\vcenter{\xymatrix@C+8pt@R=1pt{
\Delta^{op} \ar[r]^-{{\bf W}\mathrm{S}^g_{\bullet}(\C,\mathcal{A})} & \cat\\ [p]\ar@{}[r]|-{\longmapsto} & {\bf W}\mathrm{S}^g_{p}(\C,\mathcal{A})}}
\qquad\quad
\left(\text{resp.}\quad
\vcenter{\xymatrix@C+8pt@R=1pt{
\Delta^{op} \ar[r]^-{{\bf W}\mathrm{S}^d_{\bullet}(\C,\mathcal{A})} & \cat\\ [p]\ar@{}[r]|-{\longmapsto} & {\bf W}\mathrm{S}^d_{p}(\C,\mathcal{A})}}\right)\,.
\end{equation}

Si $n\geq 0$ et $\C$ est une catégorie de modèles munie d'un objet nul $0$ fixe et d'une famille distinguée $\mathcal{A}$ d'objets de $\C$ contenant l'objet $0$, le $n$-\emph{ème} $K$-\emph{groupe de Kan} (à gauche) de $(\C,\A,0)$ est par définition la classe de l'ensemble simplicial réduit:
$$
\vcenter{\xymatrix@C+8pt@R=1pt{
\Delta^{op} \ar[r]^-{\mathrm{s}^g_{\bullet}(\C,\mathcal{A})} & \ens\\ [p]\ar@{}[r]|-{\longmapsto} & \mathrm{s}^g_{p}(\C,\mathcal{A})}}
$$
dans l'ensemble des classes à isomorphismes près des objets de la catégorie homotopique $\mathrm{Ho}_{n+1}(\simp_0)$ des $(n+1)$-groupes de Kan (voir après le Corollaire \ref{dimminus}).

Le $n$-\emph{ème} $K$-\emph{groupe de Segal} (à gauche) de $(\C,\A,0)$ est la classe de l'ensemble bisimplicial:
$$
\vcenter{\xymatrix@C+12pt@R=1pt{
\big(\Delta\times \Delta\big)^{op} \ar[r] & \ens\\ \big([p],[q]\big)\ar@{}[r]|-{\longmapsto} & \mathrm{Hom}_{\simp}\big([p],{\bf W}\mathrm{S}^g_{q}(\C,\mathcal{A}) \big)  }}
$$
dans l'ensemble des classes à isomorphismes près des objets de la catégorie homotopique $\mathrm{Ho}_{n+1}({\bf MS})$ des $(n+1)$-groupes de Segal (voir autour des équivalences \eqref{duggerQeqHO} et \eqref{diagoadjHO}).

\begin{proposition}\label{waldproppremier}
Si $n\geq 0$ et $\C$ est une catégorie de modèles munie d'un objet nul $0$ fixe et d'une famille distinguée $\mathcal{A}$ d'objets de $\C$ contenant l'objet $0$, le $n$-ième $K$-groupe de Kan de $(\C,\A,0)$ est égale à la classe de l'ensemble simplicial réduit:
$$
\vcenter{\xymatrix@C+12pt@R=1pt{
\Delta^{op} \ar[r] & \ens\\ [p]\ar@{}[r]|-{\longmapsto} & \mathrm{Hom}_{\simp}\big([p],{\bf W}\mathrm{S}^g_{p}(\C,\mathcal{A}) \big)  }}
$$
dans l'ensemble des classes à isomorphismes près des objets de la catégorie homotopique $\mathrm{Ho}_{n+1}(\simp_0)$ des $(n+1)$-groupes de Kan.

En particulier les équivalences de catégories:
\begin{equation*}
\mathrm{Ho}_{n+1}({\bf MS}) \xymatrix@C+15pt{\phantom{a}\ar@/_12pt/[r]_-{{\bf R}(\,\cdot\,)_{0,\bullet}}\ar@{}[r]|-{\simeq}&\ar@/_12pt/[l]_-{{\bf L}p^*}\phantom{a}}
\mathrm{Ho}_{n+1}(\simp_0)
\qquad\text{et}\qquad
\mathrm{Ho}_{n+1}({\bf MS}) \xymatrix@C+15pt{\phantom{a}\ar@/_12pt/[r]_-{{\bf L}diag}\ar@{}[r]|-{\simeq}&\ar@/_12pt/[l]_-{{\bf R}r}\phantom{a}}
\mathrm{Ho}_{n+1}(\simp_0)
\end{equation*}
(voir \eqref{duggerQeqHO} et \eqref{diagoadjHO}) identifient les $n$-ièmes $K$-groupes de Kan et de Segal de $(\C,\A,0)$.
\end{proposition}
\begin{proof}
On va montrer en suivant Waldhausen (voir par exemple \S 2.3 de \cite{goodw}) que les ensembles simpliciaux diagonaux des ensembles bisimpliciaux:
$$
\vcenter{\xymatrix@C+12pt@R=1pt{
\big(\Delta\times \Delta\big)^{op} \ar[r] & \ens\\ \big([p],[q]\big)\ar@{}[r]|-{\longmapsto} & \mathrm{s}^g_{q}(\C,\mathcal{A})}}
\qquad \text{et} \qquad
\vcenter{\xymatrix@C+12pt@R=1pt{
\big(\Delta\times \Delta\big)^{op} \ar[r] & \ens\\ \big([p],[q]\big)\ar@{}[r]|-{\longmapsto} & \mathrm{Hom}_{\simp}\big([p],{\bf W}\mathrm{S}^g_{q}(\C,\mathcal{A}) \big)  }}
$$
c'est-à-dire les ensembles simpliciaux réduits:
$$
\vcenter{\xymatrix@C+12pt@R=1pt{
\Delta^{op} \ar[r] & \ens\\ [p]\ar@{}[r]|-{\longmapsto} & \mathrm{s}^g_{p}(\C,\mathcal{A}) \; ,}}
\qquad \text{et} \qquad
\vcenter{\xymatrix@C+12pt@R=1pt{
\Delta^{op} \ar[r] & \ens\\ [p]\ar@{}[r]|-{\longmapsto} & \mathrm{Hom}_{\simp}\big([p],{\bf W}\mathrm{S}^g_{p}(\C,\mathcal{A}) \big)  }}
$$
induisent la même classe dans l'ensemble des classes à isomorphismes près des objets de la catégorie homotopique $\mathrm{Ho}_{\infty}(\simp_0)$ des $\infty$-groupes de Kan,  

Remarquons pour commencer que si $\C_\A$ note la sous-catégorie pleine de $\C$ dont les objets sont les éléments de $\A$, alors l'ensemble simplicial $\mathrm{s}^g_{\bullet}(\C,\mathcal{A})$ est isomorphe au sous-ensemble simplicial ${\bf K}$ de l'ensemble simplicial:
\begin{equation}\label{triangleburdes1}
\xymatrix@R=1pt@C+5pt{ \Delta^{op} \ar[r] & {\bf Ens} \\ \text{\small $[q]$} \ar@{}[r]|-{\longmapsto}  & \big\{ \, \text{Foncteurs \, de \, $\righthalfcup_{q}$ \, vers \, $\C_\A$} \, \big\} }
\end{equation}
défini dans chaque $q\geq 0$ par l'ensemble ${\bf K}_q$ des foncteurs:
$$
\vcenter{\xymatrix@C+15pt@R=3pt{
\text{\Large $\righthalfcup$}_{q} \ar[r]^-{\mathfrak{X}} & \C_\mathcal{A}\\
(i,j) \ar@{}[r]|-{\longmapsto} & X_{i,j}
}}
$$
tels que:
\begin{enumerate}
\item $X_{i,i}= 0$ pour tout $0\leq i\leq q$.
\item Si $0\leq m <k\leq q-1$ le carré:
$$
 \vcenter{\xymatrix@-8pt{X_{k,m+1}\ar[r]&X_{k+1,m+1}\\X_{k,m}\ar[u]\ar[r]&X_{k+1,m}\ar[u]}}
$$ 
induit par le foncteur $\mathfrak{X}$ est homotopiquement cocartésien. 
\end{enumerate}

D'un autre côté si $p\geq 0$ posons $ p\C_\A$ pour noter la sous-catégorie pleine de la catégorie des foncteurs $\C^{[p]}$ dont les objets sont les éléments de l'ensemble:
$$
\Big\{ \; A^0\overset{f^1}{\to} \cdots \overset{f^p}{\to} A^p \; \Big| \;  A^\ell \, \in \, \A \; \text{et} \; f^\ell \,\in \, {\bf W} \; \Big\}\,,
$$
alors l'ensemble simplicial $\mathrm{Hom}_{\simp}\big([p],{\bf W}\mathrm{S}^g_{\bullet}(\C,\mathcal{A})\big)$ est isomorphe au sous-ensemble simplicial ${\bf L}$ de l'ensemble simplicial:
\begin{equation}\label{triangleburdes2}
\xymatrix@R=1pt@C+5pt{ \Delta^{op} \ar[r] & {\bf Ens} \\ \text{\small $[q]$} \ar@{}[r]|-{\longmapsto}  & \big\{ \, \text{Foncteurs \, de \, $\righthalfcup_{q}$ \, vers \, $p\C_\A$} \, \big\} }
\end{equation}
défini dans chaque $q\geq 0$ par l'ensemble ${\bf L}_q$ des foncteurs:
\begin{equation}\label{unsigmasinfff}
\vcenter{\xymatrix@C+15pt@R=3pt{
\text{\Large $\righthalfcup$}_{q} \ar[r]^-{\sigma} & p\C_A\\
(i,j) \ar@{}[r]|-{\longmapsto} & \sigma_{i,j}=\big(A_{i,j}^0\overset{f^1_{i,j}}{\to} \cdots \overset{f^p_{i,j}}{\to} A^p_{i,j}\big)
}}
\end{equation}
tels que:
\begin{enumerate}
\item[(iii)] $\sigma_{i,i} \, = \, \big(0\overset{\mathrm{id}}{\to} \cdots \overset{\mathrm{id}}{\to} 0\big)$ pour tout $0\leq i\leq q$.
\item[(iv)] Si $0\leq m < k\leq q-1$ et $0\leq \ell\leq p$  le carré induit par le foncteur $\sigma$: 
$$
 \vcenter{\xymatrix@R-6pt@C+6pt{
 A_{k,m+1}^{\ell}\ar[r]^-{h^\ell_{k,m+1}}&A_{k+1,m+1}^\ell\\A_{k,m}^\ell\ar[u]^-{v^\ell_{k,m}}\ar[r]_-{h^\ell_{k,m}}&A_{k+1,m}^\ell\ar[u]_-{v^\ell_{k+1,m}}}}
$$ 
est un carré homotopiquement cocartésien de $\C$.
\end{enumerate}

Considérons maintenant les foncteurs:
\begin{equation}\label{lesfonctnnns}
\vcenter{\xymatrix@R=3pt@C+15pt{
\C_\A \ar[r]^-{\Phi} & p\C_\A \\
A \ar@{}[r]|-{\longmapsto} & \big(A\overset{\mathrm{id}}{\to} \cdots \overset{\mathrm{id}}{\to} A\big)}}
\qquad\text{et}\qquad
\vcenter{\xymatrix@R=3pt@C+15pt{
p\C_\A \ar[r]^-{\Psi} & \C_\A \\
\big(A^0\overset{f^1}{\to} \cdots \overset{f^p}{\to} A^p\big) \ar@{}[r]|-{\longmapsto} & A^p}}
\end{equation}
et remarquons que $\Psi\circ\Phi \, = \, \mathrm{id}_{\C_\A}$ et qu'il existe une transformation naturelle:
$$
\xymatrix@C+12pt{p\C_\A\rtwocell^{\mathrm{id}}_{\Phi(p)\circ\Psi(p)}{\alpha} & p\C_\A}
$$
définie dans un objet $A^{\bullet}$ de $p\C_\A$ par le morphisme:
$$
\vcenter{\xymatrix@C-5pt{A^{\bullet} \ar[d]^-{\alpha_{A^\bullet}} \\ \Phi\circ\Psi (A^\bullet)}} \; = \; 
\left(
\vcenter{\xymatrix@C+10pt@R-5pt{
A^0\ar[r]^-{f^1} \ar[d]_-{f^p\circ \cdots \circ f^1} & A^1 \ar[r]^-{f^2} \ar[d]|-{f^p\circ\cdots\circ f^2} & \cdots \ar[r]^-{f^p} & A^p \ar[d]^-{\mathrm{id}}     \\
A^p \ar[r]_-{\mathrm{id}} & A^p \ar[r]_-{\mathrm{id}} & \cdots \ar[r]_-{\mathrm{id}} & A^p
}}
\right)
$$

On vérifie sans difficulté que les foncteurs \eqref{lesfonctnnns} induisent des morphismes d'ensembles simpliciaux:
$$
\vcenter{\xymatrix@R=3pt@C+15pt{
{\bf K} \ar[r]^-{\widetilde{\Phi}} & {\bf L} }}
\qquad\text{et}\qquad
\vcenter{\xymatrix@R=3pt@C+15pt{
{\bf L} \ar[r]^-{\widetilde{\Psi}} & {\bf K}}}
$$
tels que $\widetilde{\Psi}\circ\widetilde{\Phi} = \mathrm{id}_{\bf K}$.  

Plus encore si on note $\Gamma\colon\xymatrix@C-8pt{p\C_\A\times [1] \ar[r] & p\C_\A}$ le foncteur définit par la transformation naturelle $\alpha\colon\xymatrix@C-10pt{\mathrm{id}_{p\C^\A}\ar@{=>}[r]&\Phi\circ\Psi}$, on définit un morphisme d'ensembles simpliciaux ${H}\colon\xymatrix@C-8pt{{\bf L} \times \Delta^1 \ar[r] & {\bf L}}$ de la façon suivante: Si $q\geq 0$ et $(\sigma,\varphi)$ est un $q$-simplexe de $L\times \Delta^1$ le foncteur $H_q(\sigma,\varphi)$ est le composé:
$$
\xymatrix@C+10pt@R=2pt{
\text{\Large $\righthalfcup$}_{q}   \ar[r] & \text{\Large $\righthalfcup$}_{q}  \times [q] \ar[r]^-{\sigma\times \varphi} & p\C_\A \times [1] \ar[r]^-{\Gamma} & p\C_\A\\
(i,j) \ar@{}[r]|-{\longmapsto}& \big((i,j),i\big)&&}
$$

Montrons que le foncteur $H_q(\sigma,\varphi)\colon\xymatrix@C-10pt{ \text{\Large $\righthalfcup$}_{q}\ar[r]&p\C_\A}$ ainsi défini vérifie les propriétés (iii) et (iv) ci-dessus. En effet on note sans peine que $H_q(\sigma,\varphi)_{i,i} \, = \, \big(0\overset{\mathrm{id}}{\to} \cdots \overset{\mathrm{id}}{\to} 0\big)$ pour tout $0\leq i\leq q$. D'un autre côté si $0\leq m < k\leq q-1$ et $0\leq \ell\leq p$ on vérifie que le carré de (iv) défini par le foncteur $H_q(\sigma,\varphi)$ est le carré homotopiquement cocartésien: 
$$
\vcenter{\xymatrix@R-6pt@C+6pt{
A_{k,m+1}^{\ell}\ar[r]^-{h^\ell_{k,m+1}}&A_{k+1,m+1}^\ell\\A_{k,m}^\ell\ar[u]^-{v^\ell_{k,m}}\ar[r]_-{h^\ell_{k,m}}&A_{k+1,m}^\ell\ar[u]_-{v^\ell_{k+1,m}}}}
\qquad\text{si \; $\varphi(m+1) = \varphi(k+1) = 0$\;,}
$$ 
le carré homotopiquement cocartésien (voir les Lemmes \ref{exterio} et \ref{equivasinffmmf}):
$$
\vcenter{\xymatrix@R-6pt@C+65pt{
A_{k,m+1}^{\ell}\ar[r]^-{f^{p}_{k+1,m+1} \circ \cdots \circ f^{\ell+1}_{k+1,m+1} \circ h^\ell_{k,m+1}}&A_{k+1,m+1}^p\\ 
A_{k,m}^\ell\ar[u]^-{v^\ell_{k,m}}\ar[r]_-{f^{p}_{k+1,m} \circ \cdots \circ f^{\ell+1}_{k+1,m} \circ h^\ell_{k,m}}&A_{k+1,m}^p\ar[u]_-{v^p_{k+1,m}}}}
\qquad\text{si \; $0=\varphi(m+1) < \varphi(k+1) = 1$}
$$ 
et le carré homotopiquement cocartésien:
$$
\vcenter{\xymatrix@R-6pt@C+6pt{
A_{k,m+1}^{p}\ar[r]^-{h^p_{k,m+1}}&A_{k+1,m+1}^p\\A_{k,m}^p\ar[u]^-{v^p_{k,m}}\ar[r]_-{h^p_{k,m}}&A_{k+1,m}^p\ar[u]_-{v^p_{k+1,m}}}}
\qquad\text{si \; $\varphi(m+1) = \varphi(k+1) = 1$\,. }
$$ 

Enfin remarquons que si $\nu_i$ note le morphisme composé:
$$
\xymatrix@C+12pt{{\bf L}\,\cong \, {\bf L}\times \Delta^0 \ar[r]^-{{\bf L}\times \delta_i} & {\bf L}\times \Delta^1\; ,}
$$
alors on a un diagramme commutatif d'ensembles simpliciaux:
$$
\xymatrix@R=2pt@C+10pt{
{\bf L} \ar[rd]|-{\nu_1} \ar@/^12pt/[rrd]^-{\mathrm{id}_{\bf L}}& &\\
& {\bf L}\times \Delta^1 \ar[r]|-{H} & {\bf L}\;. \\
{\bf L} \ar[ru]|-{\nu_0} \ar@/_12pt/[rru]_-{\widetilde{\Phi}\circ \widetilde{\Psi}} & &
}
$$

Autrement dit, si on note par $F$ le morphisme de l'ensemble bisimplicial:
$$
\vcenter{\xymatrix@C+12pt@R=1pt{
\big(\Delta\times \Delta\big)^{op} \ar[r] & \ens\\ \big([p],[q]\big)\ar@{}[r]|-{\longmapsto} & \mathrm{s}^g_{q}(\C,\mathcal{A}) \, = \, \mathrm{Hom}_{\simp}\big([0],{\bf W}\mathrm{S}^g_{q}(\C,\mathcal{A}) \big)}}
$$
vers l'ensemble bisimplicial:
$$
\vcenter{\xymatrix@C+12pt@R=1pt{
\big(\Delta\times \Delta\big)^{op} \ar[r] & \ens\\ \big([p],[q]\big)\ar@{}[r]|-{\longmapsto} & \mathrm{Hom}_{\simp}\big([p],{\bf W}\mathrm{S}^g_{q}(\C,\mathcal{A}) \big)  }}
$$
induit par les morphismes $\xymatrix@C-12pt{[p]\ar[r]&[0]}$, on a montré que $F_{p,\bullet}$ est une $\infty$-équivalence faible d'ensembles simpliciaux réduits pour tout $p\geq 0$.
\end{proof}

Si $n\geq 0$ d'après la Proposition \ref{waldproppremier} le $(n+1)$-ième groupe d'homotopie de la réalisation géométrique des ensembles simpliciaux réduits:
$$
\vcenter{\xymatrix@C+12pt@R=1pt{
\Delta^{op} \ar[r] & \ens\\ [p]\ar@{}[r]|-{\longmapsto} & \mathrm{Hom}_{\simp}\big([p],{\bf W}\mathrm{S}^g_{p}(\C,\mathcal{A}) \big)  }}
\qquad\;\text{et}\qquad\;
\vcenter{\xymatrix@C+12pt@R=1pt{
\Delta^{op} \ar[r] & \ens\\ [p]\ar@{}[r]|-{\longmapsto} & \mathrm{s}^g_{p}(\C,\mathcal{A}) }}
$$
sont le même élément dans l'ensemble des groupes à isomorphisme près. On appelle la classe de ce groupe le \emph{$n$-ième $K$-groupe classique de $(\C,\A,0)$}.

\renewcommand{\thesubsection}{\S\thesection.\arabic{subsection}}
\subsection{}\;
\renewcommand{\thesubsection}{\thesection.\arabic{subsection}}

Si $p\geq 0$ et $\C$ est une catégorie de modèles munie d'un objet nul $0$ fixe et d'une famille distinguée $\mathcal{A}$ d'objets de $\C$ contenant l'objet $0$, un $(p-1)$-\emph{triangle strict à gauche} (resp. \emph{à droite}) $\mathfrak{X}$ \emph{de $\C$} (\emph{à valeurs dans $\A$}) est par définition un $(p-1)$-triangle à gauche (resp. à droite) de $\C$ (à valeurs dans $\A$) vérifiant les propriétés:
\begin{enumerate}
\item[(iv)]  Pour tout $1\leq k\leq p-1$ et $0\leq m\leq p-2$, si $\xymatrix{\boxempty\ar[r]^-w&\text{\Large $\righthalfcup$}_{p}}$ désigne le foncteur:
$$
w\,\def\objectstyle{\scriptstyle}
\def\labelstyle{\scriptstyle}
\left(\vcenter{\xymatrix@-10pt{c\ar[r]&d\\a\ar[u]\ar[r]&b\ar[u]}}\right) \; = \; 
\def\objectstyle{\scriptstyle}
\def\labelstyle{\scriptstyle}
\vcenter{\xymatrix@-8pt{(k,m+1)\ar[r]&(k+1,m+1)\\(k,m)\ar[u]\ar[r]&(k+1,m)\ar[u]}}\,,
$$
alors $w^*(\mathfrak{X})$ est un carré cocartésien (resp. homotopiquement cartésien). 
\item[(v)] Pour tout $0 \leq k\leq p-1$ le morphisme $\xymatrix@-5pt{\mathfrak{X}_{k,0}\ar[r]&\mathfrak{X}_{k+1,0}}$ (resp.  $\xymatrix@-5pt{\mathfrak{X}_{p,k}\ar[r]&\mathfrak{X}_{p,k+1}}$) est une cofibration (resp. une fibration). En particulier, les morphismes $\xymatrix@-5pt{\mathfrak{X}_{i,j}\ar[r]&\mathfrak{X}_{i',j}}$ (resp.  $\xymatrix@-5pt{\mathfrak{X}_{i,j}\ar[r]&\mathfrak{X}_{i,j'}}$) sont de cofibrations (resp. fibrations) si $0\leq i< i'\leq p$ et $0\leq j\leq p-1$ (resp. si $1\leq i\leq p$ et $0\leq j<j'\leq p$) et les objets $\mathfrak{X}_{i,j}$ sont cofibrants (resp. fibrants) si $0\leq i\leq j\leq p$.
\end{enumerate}

On pose ${\bf W}\mathrm{S}^{g}_{p}(\C,\mathcal{A})^{cart}$ $\big($resp. ${\bf W}\mathrm{S}^d_{p}(\C,\A)^{cart}\big)$ pour noter la sous-catégorie pleine de ${\bf W}\mathrm{S}^{g}_{p}(\C,\mathcal{A})$ $\big($resp. ${\bf W}\mathrm{S}^d_{p}(\C,\A)\big)$ dont les objets sont les $(p-1)$-triangles stricts à gauche (resp. à droite) de $\C$ à valeurs dans $\A$. 

Il se suit aussi-tôt que le foncteur \eqref{ssW} induit par restriction un foncteur:
\begin{equation}\label{ssWcar}
\vcenter{\xymatrix@C+8pt@R=1pt{
\Delta^{op} \ar[r]^-{{\bf W}\mathrm{S}^g_{\bullet}(\C,\mathcal{A})^{cart}} & \cat\\ [p]\ar@{}[r]|-{\longmapsto} & {\bf W}\mathrm{S}^g_{p}(\C,\mathcal{A})^{cart}}}
\qquad\quad
\left(\text{resp.}\quad
\vcenter{\xymatrix@C+8pt@R=1pt{
\Delta^{op} \ar[r]^-{{\bf W}\mathrm{S}^d_{\bullet}(\C,\mathcal{A})^{cart}} & \cat\\ [p]\ar@{}[r]|-{\longmapsto} & {\bf W}\mathrm{S}^d_{p}(\C,\mathcal{A})^{cart}}}\right)\,.
\end{equation}

\begin{proposition}\label{propwaldkth}
Soit $n\geq 0$ et $\C$ une catégorie de modèles munie d'un objet nul $0$ fixe et d'une famille distinguée $\mathcal{A}$ d'objets de $\C$ contenant l'objet $0$. Si la famille d'objets $\A$ vérifie la propriété:
\begin{enumerate}
\item[(a)] $\A$ est fermée par des équivalences faibles \emph{i.e.} si $\xymatrix@C-5pt{A\ar@{->}[r]^-{\text{\rotatebox[origin=c]{180}{\Large $\widetilde{\phantom{w}}$}}}&B}$ est une équivalence faible de $\C$ alors $A$ appartient à $\mathcal{A}$ si et seulement s'il en est de même pour $B$.
\end{enumerate}
le $n$-ième $K$-groupe de Segal (à gauche) de $(\C,\A,0)$ est la classe de l'ensemble bisimplicial:
$$
\vcenter{\xymatrix@C+12pt@R=1pt{
\big(\Delta\times \Delta\big)^{op} \ar[r] & \ens\\ \big([p],[q]\big)\ar@{}[r]|-{\longmapsto} & \mathrm{Hom}_{\simp}\big([p],{\bf W}\mathrm{S}^g_{q}(\C,\mathcal{A})^{cart} \big)  }}
$$
dans l'ensemble des classes à isomorphismes près des objets de la catégorie homotopique $\mathrm{Ho}_{n+1}({\bf MS})$ des $(n+1)$-groupes de Segal, et le $n$-ième $K$-groupe de Kan de $(\C,\A,0)$ est égale à la classe de l'ensemble simplicial:
$$
\vcenter{\xymatrix@C+12pt@R=1pt{
\Delta^{op} \ar[r] & \ens\\ [p]\ar@{}[r]|-{\longmapsto} & \mathrm{Hom}_{\simp}\big([p],{\bf W}\mathrm{S}^g_{p}(\C,\mathcal{A})^{cart} \big)  }}
$$
dans l'ensemble des classes à isomorphismes près des objets de la catégorie homotopique $\mathrm{Ho}_{n+1}(\simp_0)$ des $(n+1)$-groupes de Kan.
\end{proposition}
\begin{proof}
Montrons par ailleurs que pour tout $q\geq 0$ il existe un foncteur:
\begin{equation}\label{phiynpmFF}
\xymatrix@+10pt{{\bf W}\mathrm{S}^g_{q}(\C,\mathcal{A})\ar[r]^-{\Phi_{q}}&{\bf W}\mathrm{S}^g_{q}(\C,\mathcal{A})}\qquad\;
\end{equation}
et une transformation naturelle $\vartheta_{q}:\xymatrix@-12pt{\Phi_{q}\ar@{=>}[r]&\mathrm{id}}$ tel que $\Phi_q(\mathfrak{X})$ est un $(q-1)$-triangle strict pour tout $(q-1)$-triangle à gauche $\mathfrak{X}$ de $\C$. La preuve est par induction sur $q$.

On définit $\Phi_{0}$ comme étant le foncteur identité de la catégorie ponctuelle. D'autre part $\Phi_{1}$ et $\vartheta_{1}$ sont définis à partir  d'un remplacement cofibrant fonctoriel \eqref{laQ}.

Supposons qu'on a déjà défini $\Phi_{q-1}$ et $\vartheta_{q-1}$ vérifiant la propriété désiré et considérons le foncteur d'inclusion:
$$
\xymatrix@+15pt@R=3pt{
\text{\Large $\righthalfcup$}_{q-1}\ar[r]^-{\nu_q}&\text{\Large $\righthalfcup$}_{q}\;.\\
(i,j)\ar@{}[r]|-{\longmapsto} &(i,j) 
}
$$

Si $\mathfrak{X}$ est un $q$-triangle à gauche de $\C$:
$$
\def\objectstyle{\scriptstyle}
\def\labelstyle{\scriptstyle}
\mathfrak{X}\;=\, \vcenter{\xymatrix@C-1.2pc@R-1.8pc{
   &&&&&&\\
   &&&&&0\\
   &&&&&\\
   &&&&\ar@{.}[r]&\ar[uu]\phantom{.}\\&&&&&&&\\
   &&&\ar@{.}[r]&\ar@{.}[r]\ar@{.}[uu]&\ar@{.}[uu]\\
   &&&&&\\
   &&0\ar[r]&\ar@{.}[r]\ar@{.}[uu]&\ar@{.}[uu]\ar[r]&X_{q+1,1}\ar[uu]\\
   &&&&&\\
   &0\ar[r]&X_{1,0}\ar[uu]\ar[r]&\ar@{.}[r]\ar@{.}[uu]&\ar@{.}[uu]\ar[r]&X_{q+1,0}\ar[uu]\\
   &&&&&&}},
$$   
on définit la restriction $\nu_{q}^*\big(\Phi_{q}\mathfrak{X}\big)$ comme étant le diagramme:
$$
\def\objectstyle{\scriptstyle}
\def\labelstyle{\scriptstyle}
\Phi_{q-1}(\nu_{q}^*\mathfrak{X}) \;=\, \vcenter{\xymatrix@C-1.2pc@R-1.8pc{
   &&&&&&\\
   &&&&&0\\
   &&&&&\\
   &&&&\ar@{.}[r]&\ar[uu]\phantom{.}\\&&&&&&&\\
   &&&\ar@{.}[r]&\ar@{.}[r]\ar@{.}[uu]&\ar@{.}[uu]\\
   &&&&&\\
   &&0\ar[r]&\ar@{.}[r]\ar@{.}[uu]&\ar@{.}[uu]\ar[r]&A_{q,1}\ar[uu]\\
   &&&&&\\
   &0\ar[r]&A_{1,0}\ar[uu]\ar[r]&\ar@{.}[r]\ar@{.}[uu]&\ar@{.}[uu]\ar[r]&A_{q,0}\ar[uu]\\
   &&&&&&}},
$$ 
et le morphisme $\nu_{q}^* \big((\vartheta_{q})_{\mathfrak{X}}\big)$ comme $(\vartheta_{q-1})_{\nu_{q}^*\mathfrak{X}}$.

La colonne de $\Phi_{q}\mathfrak{X}$ et la partie du morphisme $(\vartheta_{q})_{\mathfrak{X}}$ qui reste à décrire sont définies de façon inductive par de cubes:
\begin{equation}\label{cubekkfer}
\xymatrix@-15pt{
&&X_{q,k+1}\ar[rrr]\ar'[d][ddd]&&&X_{q+1,k+1}\,\ar[ddd]\\
A_{q,k+1}\phantom{.}\ar[rru]|-{\,F_{q,k+1}\,}\ar@{>->}@<-1pt>[rrr]\ar[ddd]&&&A_{q+1,k+1}\ar[ddd]\ar[rru]|-{\,F_{q+1,k+1}\,}&&\\
&&&&&\\
&&X_{q,k}\ar'[r][rrr]&&&X_{q+1,k}\,,\\
A_{q,k}\phantom{.}\ar@{>->}[rrr]\ar[urr]|-{\,F_{q,k}\,}&&&A_{q+1,k}\ar[rru]|-{\,F_{q+1,k}\,}}
\end{equation}
où $0\leq k\leq q$ comme suit: Premièrement on décompose de façon fonctorielle le morphisme composé:
$$
\xymatrix@-15pt{
&&X_{q,k}\ar[rrr]&&&X_{q+1,k}\,,\\
A_{q,k}\phantom{.}\ar[urr]|-{\,F_{q,k}\,}&&&}
$$
comme une cofibration suivie d'une fibration triviale de $\C$:
$$
\xymatrix@-15pt{
&&&&&X_{q+1,k}\,,\\
A_{q,k}\phantom{.}\ar@{>->}[rrr]&&&A_{q+1,k}\ar[rru]|-{\,F_{q+1,k}\,}}
$$

En suit on prend de façon fonctorielle un carré cocartésien de $\C$:
$$
\xymatrix@-15pt{
A_{q,k+1}\phantom{.}\ar@{>->}[rrr]\ar[ddd]&&&A_{q+1,k+1}\ar[ddd]&&\\
&&&&&\\
&&&&&\,,\\
A_{q,k}\phantom{.}\ar@{>->}[rrr]&&&A_{q+1,k}}
$$
et on définit le morphisme $F_{q+1,k+1}\colon\xymatrix@C-10pt{A_{q+1,k+1}\ar[r]&X_{q+1,k+1}}$ comme étant le seul morphisme tel que \eqref{cubekkfer} est un cube dont toutes les faces sont commutatives. En particulier $F_{q+1,k-1}$ est une équivalence faible d'après le Lemmes \ref{cubelele} et \ref{fibrabra}.

Remarquons d'un autre côté que si on pose:
$$
\vcenter{\xymatrix@+10pt{{\bf W}\mathrm{S}^{g}_{q}(\C,\mathcal{A})^{cart}\;\ar@{^(->}[r]^-{\Psi_{q}}&{\bf W}\mathrm{S}^g_{q}(\C,\mathcal{A})}},
$$
pour noter le foncteur d'inclusion canonique, le foncteur \eqref{phiynpmFF} qu'on vient de définir induit un foncteur:
$$
\xymatrix@+10pt{{\bf W}\mathrm{S}^g_{q}(\C,\mathcal{A})\ar[r]^-{\Phi_{q}}&{\bf W}\mathrm{S}^g_{q}(\C,\mathcal{A})^{cart}}\qquad\;
$$
et la transformation naturelle $\vartheta_{q}:\xymatrix@-12pt{\Phi_{q}\ar@{=>}[r]&\mathrm{id}}$ induite des transformations naturelles: 
$$
\xymatrix@-8pt{\Phi_{q}\Psi_{q}\ar@{=>}[r]&\mathrm{id}} \qquad \text{et} \qquad \xymatrix@-8pt{\Psi_{q}\Phi_{q}\ar@{=>}[r]&\mathrm{id}}\,.
$$

Autrement dit pour chaque $q\geq 0$ l'image de $\Psi_{q}$ par le foncteur nerf $\xymatrix@-8pt{\cat\ar[r]&\simp}$ est une équivalence homotopique d'ensembles simpliciaux. Donc si $\Psi$ note le morphisme d'inclusion canonique de l'ensemble bisimplicial:
$$
\vcenter{\xymatrix@C+12pt@R=1pt{
\big(\Delta\times \Delta\big)^{op} \ar[r] & \ens\\ \big([p],[q]\big)\ar@{}[r]|-{\longmapsto} & \mathrm{Hom}_{\simp}\big([p],{\bf W}\mathrm{S}^g_{q}(\C,\mathcal{A})^{cart} \big)  }}
$$
vers l'ensemble bisimplicial:
$$
\vcenter{\xymatrix@C+12pt@R=1pt{
\big(\Delta\times \Delta\big)^{op} \ar[r] & \ens\\ \big([p],[q]\big)\ar@{}[r]|-{\longmapsto} & \mathrm{Hom}_{\simp}\big([p],{\bf W}\mathrm{S}^g_{q}(\C,\mathcal{A})\big)  }}
$$
on a que $\Psi_{\bullet,q}$ est une $\infty$-équivalence faible d'ensembles simpliciaux pour tout $q\geq 0$. 

Donc l'image de $\Psi$ par le foncteur ensemble simplicial diagonal $\xymatrix@-8pt{\ssimp\ar[r]&\simp}$ est une $\infty$-équivalence faible d'ensembles simpliciaux.
\end{proof}

On vérifie sans peine:

\begin{lemme}\label{waldkth}
Soit $\C$ une catégorie de modèles munie d'un objet nul $0$ fixe et d'une famille distinguée $\mathcal{A}$ d'objets de $\C$ contenant l'objet $0$ et vérifiant la propriété:
\begin{enumerate}
\item[(b)] La famille $\A$ est stable par cochangement de base le long des cofibrations \emph{i.e.} si:
\begin{equation}\label{}
\vcenter{\xymatrix@-5pt{A\ar[r]^-f\ar[d]& B\ar[d]\\C\ar[r] & D}}\,,
\end{equation}
est un carré cocartésien de $\C$ tel que $f$ est une cofibration et les objets $A$, $B$ et $C$ appartient à $\mathcal{A}$ alors $D$ appartient à $\mathcal{A}$ aussi.
\end{enumerate}

Si $\mathcal{A}_{cof}$ note la sous-catégorie pleine de $\C$ dont les objets sont les objets cofibrants qui appartient à $\mathcal{A}$, alors $\mathcal{A}_{cof}$ admet une structure de catégorie de Waldhausen (voir \cite{waldhausen}) dont les cofibrations (resp. les équivalences faibles) sont les cofibrations (resp. les équivalences faibles) de $\C$ entre des objets de $\A_{cof}$.
\end{lemme}

On déduit du Corollaire \ref{waldkth} et des Propositions \ref{waldproppremier} et \ref{propwaldkth}:

\begin{corollaire}
Soit $\C$ une catégorie de modèles munie d'un objet nul $0$ fixe et d'une famille distinguée $\mathcal{A}$ d'objets de $\C$ contenant l'objet $0$ vérifiant les propriétés:
\begin{enumerate}
\item[(a)] $\A$ est fermée par des équivalences faibles \emph{i.e.} si $\xymatrix@C-5pt{A\ar@{->}[r]^-{\text{\rotatebox[origin=c]{180}{\Large $\widetilde{\phantom{w}}$}}}&B}$ est une équivalence faible de $\C$ alors $A$ appartient à $\mathcal{A}$ si et seulement s'il en est de même pour $B$.
\item[(b)] La famille $\A$ est stable par cochangement de base le long des cofibrations \emph{i.e.} si:
\begin{equation}\label{}
\vcenter{\xymatrix@-5pt{A\ar[r]^-f\ar[d]& B\ar[d]\\C\ar[r] & D}}\,,
\end{equation}
est un carré cocartésien de $\C$ tel que $f$ est une cofibration et les objets $A$, $B$ et $C$ appartient à $\mathcal{A}$ alors $D$ appartient à $\mathcal{A}$ aussi.
\end{enumerate}

Si $n\geq 0$ le $n$-ième $K$-groupe classique de $(\C,\A,0)$ est égal au $n$-ième $K$-groupe de Waldhausen de la catégorie de Waldhausen $\A_{cof}$ du Lemme \ref{waldkth} \emph{i.e.} à la classe du $n$-ième groupe d'homotopie de la réalisation géométrique de l'ensemble simplicial:
$$
\vcenter{\xymatrix@C+12pt@R=1pt{
\Delta^{op} \ar[r] & \ens\\ [p]\ar@{}[r]|-{\longmapsto} & \mathrm{Hom}_{\simp}\big([p],{\bf W}\mathrm{S}^g_{p}(\C,\mathcal{A})^{cart} \big) }}\,,
$$
dans l'ensemble des groupes à isomorphisme près.
\end{corollaire}

\renewcommand{\thesubsection}{\S\thesection.\arabic{subsection}}
\subsection{}\;\label{ffonctionadditetdet}
\renewcommand{\thesubsection}{\thesection.\arabic{subsection}}

Soit $\C$ une catégorie de modèles munie d'un objet nul $0$ fixe et d'une famille distinguée $\mathcal{A}$ d'objets de $\C$ contenant l'objet $0$. Si $G$ est un groupe une \emph{fonction additive} de $(\C,\A,0)$ \emph{à valeurs dans} $G$ est une fonction des ensembles:
$$
\xymatrix@C+10pt{
\A \ar[r]^-{F} & G
}
$$
telle que:
\begin{enumerate}
\item $F(0) \, = \, e_G$ où $e_G$ est le neutre de $G$.
\item Si $X$, $Y$ et $Z$ sont des objets de $\C$ appartenant à la famille $\A$ et on a un carré homotopiquement cocartésien de $\C$:
$$
\def\objectstyle{\scriptstyle}
\def\labelstyle{\scriptstyle}
\vcenter{\xymatrix@-5pt{0\ar[r]&Z \\ X\ar[r]\ar[u]&Y\ar[u]}}\,,
$$ 
alors $F(Y) \, = \, F(X) \, \cdot \, F(Z)$.
\end{enumerate}

Note que si $\varphi\colon\xymatrix@C-10pt{G\ar[r]&H}$ est un morphisme de groupes et $F\colon\xymatrix@C-10pt{\A\ar[r]&G}$ est une fonction additive de $(\C,\A,0)$ à valeurs dans $G$, alors le composé $\varphi\circ F\colon\xymatrix@C-10pt{\A\ar[r]&H}$ est un déterminant de $(\C,\A,0)$ à valeurs dans $H$.

\begin{proposition}
Si $\C$ est une catégorie de modèles munie d'un objet nul $0$ fixe et d'une famille distinguée $\mathcal{A}$ d'objets de $\C$ contenant l'objet $0$, le foncteur:
$$
\xymatrix@C+25pt@R=1pt{
{\bf Grp} \ar[r] & {\bf Ens} \\
G \ar@{}[r]|-{\longmapsto} & \Bigg\{ \overset{\text{\small $\text{Fonctions additives}$}}{\underset{\text{\small $\;\text{de}\;(\C,\A,0)\;\text{vers}\; G\;$}}{\phantom{AAA}}} \Bigg\}}
$$
est représentable par le $0$-ième $K$-groupe classique de $(\C,\A,0)$ \emph{i.e.} par le groupe fondamental de la réalisation géométrique de l'ensemble simplicial réduit:
$$
\vcenter{\xymatrix@C+12pt@R=1pt{
\Delta^{op} \ar[r] & \ens\\ [p]\ar@{}[r]|-{\longmapsto} & \mathrm{Hom}_{\simp}\big([p],{\bf W}\mathrm{S}^g_{p}(\C,\mathcal{A}) \big)  }}
$$
ou de l'ensemble simplicial réduit:
$$
\vcenter{\xymatrix@C+12pt@R=1pt{
\Delta^{op} \ar[r] & \ens\\ [p]\ar@{}[r]|-{\longmapsto} & \mathrm{s}^g_{p}(\C,\mathcal{A}) }}\,.
$$
\end{proposition}
\begin{proof}
Cette Proposition est un cas particulier du Corollaire \ref{additivffm}.
\end{proof}

\renewcommand{\thesubsection}{\S\thesection.\arabic{subsection}}
\subsubsection{}\;
\renewcommand{\thesubsection}{\thesection.\arabic{subsection}}

Soit $\C$ une catégorie de modèles munie d'un objet nul $0$ fixe et d'une famille distinguée $\mathcal{A}$ d'objets de $\C$ contenant l'objet $0$. Si $\G$ est un $2$-groupe (voir la section \ref{2groupeseection}) un \emph{déterminant de $(\C,\A,0)$ á valeurs dans $\G$} est une couple $D=(D,T)$ où:
$$
\xymatrix{{\bf W}\mathrm{S}^g_{1}(\C,\mathcal{A}) \ar[r]^-{D} & \G}
$$
est un foncteur de la sous-catégorie ${\bf W}\mathrm{S}^g_{1}(\C,\mathcal{A})$ de $\C$ dont les objets sont les éléments de $\A$ et les morphismes sont les équivalences faibles entre les objets de  $\A$ vers le groupoïde sous-jacent à $\G$, et:
$$
\xymatrix{\mathrm{s}^g_2(\C,\A) \ar[r]^-{T} & \{\text{Morphismes de $\G$}\}}
$$
est une fonction de l'ensemble des $1$-triangles (à gauche) de $\C$ (à valeurs dans $\A$) vers l'ensemble des morphismes de $\G$, vérifiant les propriétés:
\begin{enumerate}
\item (Compatibilité) \, Si $\mathfrak{X}$ est un $1$-triangle (à gauche) de $\C$ (à valeurs dans $\A$):
$$
\mathfrak{X} \, = \, \vcenter{\def\objectstyle{\scriptstyle}
\def\labelstyle{\scriptstyle}\xymatrix@-5pt{
0\ar[r] & X_0 \\
X_2 \ar[u]\ar[r] & X_1\ar[u]
}}\,,
$$
$T(\mathfrak{X})$ est un morphisme de $\G$ de la forme: 
$$
\xymatrix@C+15pt{D_0(X_2) \otimes D_0(X_0)\ar[r]^-{T(\mathfrak{X})} & D_0(X_1) }\,.
$$
\item (Unitaire) \, On a que $D_0\big(\,0\,\big)\,=\,\mathbb{1}_{\G}$ et 
$T\left(\,\vcenter{
\def\objectstyle{\scriptstyle}
\def\labelstyle{\scriptstyle}
\xymatrix@-5pt{
0\ar[r] & 0 \\
0 \ar[u]\ar[r] & 0\ar[u]
}}\,\right)\,=\,l_{\mathbb{1}}^{-1}\,=\,r_{\mathbb{1}}^{-1}$.
\item (Naturalité)\, Si $\zeta\colon\xymatrix@C-10pt{\mathfrak{X}\ar[r]&\mathfrak{Y}}$ est un morphisme de la catégorie ${\bf W}\mathrm{S}^g_{2}(\C,\mathcal{A})$:
$$
\zeta \, = \,\qquad \vcenter{
\def\objectstyle{\scriptstyle}
\def\labelstyle{\scriptstyle}
\xymatrix@R-20pt@C-17pt{
&&  X_0 \ar[rrr]^{\zeta_0}&&& Y_0 \\
&&&&&\\
0 \ar@<-2pt>@{=}[rrr] \ar[uurr]&&& 0 \ar[rruu] &&\\
&&  X_1\ar'[u][uuu] \ar'[r][rrr]|-{\zeta_1}&&& Y_1\ar[uuu] \\
&&&&&\\
X_2 \ar[uuu]\ar[rrr]_-{\zeta_2} \ar[uurr]&&& Y_2 \ar[uuu] \ar[rruu] &&
}}$$
on a un carré commutatif de $\G$:
$$
\vcenter{\xymatrix@C+10pt{
D_0(X_2) \otimes D_0(X_0)\ar[d]_-{D_1(\zeta_2)\otimes D_1(\zeta_0)} \ar[r]^-{T(\mathfrak{X})} & D_0(X_1) \ar[d]^-{D_1(\zeta_1)}\\
D_0(Y_0) \otimes D_0(Y_0) \ar[r]_-{T(\mathfrak{Y})} & D_0(Y_1)
}}\,.
$$

\item (Associativité)\, Si $\eta$ est un $2$-triangle (à gauche) de $\C$ (à valeurs dans $\A$):
$$
\eta \, = \,
\def\objectstyle{\scriptstyle}
\def\labelstyle{\scriptstyle}
\vcenter{\xymatrix@C-1.2pc@R-1.8pc{
   &0\ar[r]&A_{3,2}\\
   &&\\
   0\ar[r]&A_{2,1}\ar[r]\ar[uu]&A_{3,1}\ar[uu]\\
   &&\\
   A_{1,0}\ar[uu]\ar[r]&A_{2,0}\ar[r]\ar[uu]&A_{3,0}\ar[uu]}}
$$
on a un diagramme commutatif de $\G$:
$$
\def\objectstyle{\scriptstyle}
\def\labelstyle{\scriptstyle}
\vcenter{\xymatrix@R-10pt@C+45pt{
D_0(A_{30})   &D_0(A_{10})\otimes D_0(A_{31})  \ar[l]_-{T(d_2\,\eta)} \\&\\
&D_0(A_{10}) \otimes (D_0(A_{21})\otimes D_0(A_{32})) \ar[uu]_-{D_0(A_{10})\otimes T(d_0\, \eta) }  \ar@{}[d]_-{a}|-{\text{\rotatebox[origin=c]{90}{$\cong$}}}\\
D_0(A_{20})\otimes D_0(A_{32})  \ar[uuu]^-{T(d_1\,\eta)}&(D_0(A_{10})\otimes D_0(A_{21}))\otimes D_0(A_{32})\ar[l]^-{T(d_3\,\eta)\otimes D_0(A_{32})}\,,}}
$$
\begin{align*}
\text{où} \qquad 
A_{03} \, = \, & d_{1} d_{1} \eta \, = \, d_{1} d_{2} \eta \qquad
A_{01} \, = \,  d_{2} d_{2} \eta \, = \, d_{2} d_{3} \eta \qquad
A_{13} \, = \,  d_{1} d_{0} \eta \, = \, d_{0} d_{2} \eta \\
A_{02} \, = \, & d_{2} d_{1} \eta \, = \, d_{1} d_{3}  \eta \qquad
A_{23} \, = \,  d_{0} d_{0} \eta \, = \, d_{0} d_{1} \eta \qquad
A_{12} \, = \,  d_{2} d_{0} \eta \, = \, d_{0} d_{3} \eta \,.
\end{align*}
\end{enumerate}

Si $(D,T)$ et $(D',T')$ sont deux déterminants de $(\C,\A,0)$ à valeurs dans $\G$ un \emph{morphisme de déterminants} $\alpha\colon\xymatrix@C-10pt{(D,T)\ar[r] & (D',T')}$ est une transformation naturelle:
$$
\xymatrix@C+10pt{{\bf W}\mathrm{S}^g_{1}(\C,\mathcal{A})\phantom{aaa} \rtwocell<\omit>{\phantom{a}\alpha}  \ar@<+3pt>@/^13pt/[r]^-{D} \ar@<-3pt>@/_13pt/[r]_-{D'} & \G}
$$
telle que:
\begin{enumerate}
\item Le morphisme $\alpha_0$ est le morphisme identité de l'objet $\mathbb{1}$ de $\G$.
\item Si $\mathfrak{X}$ est un $1$-triangle (à gauche) de $\C$ (à valeurs dans $\A$):
$$
\mathfrak{X} \, = \, \vcenter{\def\objectstyle{\scriptstyle}
\def\labelstyle{\scriptstyle}\xymatrix@-5pt{
0\ar[r] & X_0 \\
X_2 \ar[u]\ar[r] & X_1\ar[u]
}}
$$
le carré:
$$
\vcenter{\xymatrix@C+10pt{
D_0(X_1) \ar[d]_-{\alpha_{X_1}} \ar[r]^-{T(\mathfrak{X})} & D_0(X_2) \otimes D_0(X_0)\ar[d]^-{\alpha_{X_2}\otimes \alpha_{X_0}}\\
D'_0(X_1) \ar[r]_-{T(\mathfrak{Y})} & D'_0(X_0) \otimes D'_0(X_0)
}}
$$
est un diagramme commutatif de $\G$.
\end{enumerate}

Remarquons que les déterminants de $(\C,\A,0)$ à valeurs dans $\G$ et les morphismes entre eux forment un groupoïde qu'on note $\underline{\bf det}_{(\C,\A,0)}\big(\G\big)$. La composition étant définie par la composition de transformations naturelles. 

On définit un $2$-foncteur:
$$
\vcenter{\xymatrix@C+40pt{
\text{$2$-$\underline{\bf Grp}$} \ar[r]^-{\underline{\bf det}_{(\C,\A,0)}}  & \underline{\bf Grpd}
}}
$$
comme suit: Si $(\varphi,m^{\varphi})\colon\xymatrix@C-10pt{\G\ar[r]&\H}$ est un morphisme de $2$-groupes le foncteur:
$$
\xymatrix@C+50pt@R=5pt{
\underline{\bf det}_{(\C,\A,0)}(\G) \ar[r]^-{\underline{\bf det}_{(\C,\A,0)}(\varphi,m^\varphi)} & \underline{\bf det}_{(\C,\A,0)}(\H)
}
$$
est définie dans un déterminant $(D,T)$ de $(\C,\A,0)$ à valeurs dans $\G$ par la formule: 
$$
\underline{\bf det}_{(\C,\A,0)}\big(\varphi,m^\varphi\big)(D,T)\,=\,(\overline{D},\overline{T})
$$ 
où $\overline{D}$ est le foncteur composé:
$$
\xymatrix@C-5pt{
{\bf W}\mathrm{S}^g_{1}(\C,\mathcal{A}) \ar[r]^-{D} & \G \ar[r]^-{\varphi} & \H}
$$
et $\overline{T}$ est la fonction composée:
$$
\xymatrix@C-5pt{
\mathrm{s}^g_{2}(\C,\mathcal{A}) \ar[r]^-{T} & \mathcal{N}_{\mathcal{S}}\big(\G\big)_{0,2} 
\ar[rr]^-{\mathcal{N}_\mathcal{S}(\varphi,m^\varphi)_{0,2}} && \mathcal{N}_{\mathcal{S}}\big(\H\big)_{0,2}\,\subset \, \big\{\,\text{Morphismes de $\H$}\,\big\}\,.}
$$

Si d'un autre $\xymatrix@C+20pt{\mathcal{G}\rtwocell^{(\varphi,m^{\varphi})}_{(\psi,m^{\psi})}{\eta}& \mathcal{H}}$ est un transformation entre morphismes de $2$-groupes, la transformation naturelle:
$$
\xymatrix@C+50pt{
\underline{\bf det}_{(\C,\A,0)}(\G)  \rtwocell^{{\underline{\bf det}_{(\C,\A,0)}(\varphi,m^\varphi)}}_{{\underline{\bf det}_{(\C,\A,0)}(\psi,m^\psi)}}{\overline{\eta}}& \underline{\bf det}_{(\C,\A,0)}(\H)
}
$$
est définie dans un déterminant $(D,T)$ de $(\C,\A,0)$ à valeurs dans $\G$ par la transformation naturelle:
$$
\xymatrix@C+10pt{{\bf W}\mathrm{S}^g_{1}(\C,\mathcal{A})\phantom{aaa} \rtwocell<\omit>{\phantom{a}D\star \eta}  \ar@<+3pt>@/^13pt/[r]^-{D\circ\varphi} \ar@<-3pt>@/_13pt/[r]_-{D\circ \psi} & \H\;.}.
$$

\begin{proposition}
Si $\C$ est une catégorie de modèles munie d'un objet nul $0$ fixe et d'une famille distinguée $\mathcal{A}$ d'objets de $\C$ contenant l'objet $0$, le foncteur:
\begin{equation}\label{lefonctdetfF}
\xymatrix@C+25pt@R=1pt{
\text{$2$-${\bf Grp}$} \ar[r]^-{\underline{\bf det}_{(\C,\A,0)}}   & {\bf Grpd} \ar[r]^-\pi & h{\bf Grpd}}
\end{equation}
est représentable par le $1$-ère $K$-groupe (de Kan ou de Segal) de $(\C,\A,0)$; plus précisément le foncteur \eqref{lefonctdetfF} est représentable par le $2$-groupe d'homotopie (voir le Chapitre 3) du pré-monoïde de Segal:
$$
\vcenter{\xymatrix@C+12pt@R=1pt{
\big(\Delta\times\Delta\big)^{op} \ar[r] & \ens\\ \big([p],[q]\big)\ar@{}[r]|-{\longmapsto} & \mathrm{Hom}_{\simp}\big([p],{\bf W}\mathrm{S}^g_{q}(\C,\mathcal{A}) \big)  }}
$$
ou de l'ensemble simplicial réduit:
$$
\vcenter{\xymatrix@C+12pt@R=1pt{
\Delta^{op} \ar[r] & \ens \,.\\ [p]\ar@{}[r]|-{\longmapsto} & \mathrm{s}^g_{p}(\C,\mathcal{A}) }}
$$

En particulier si $(\varphi,m^{\varphi})\colon \xymatrix@C-10pt{\G\ar[r]& \H}$ est une $2$-équivalence faible de $2$-groupes le foncteur:
$$
\xymatrix@C+45pt{
 \underline{{\bf det}}_{(\C,\A,0)}(\G) \ar[r]^-{ \underline{{\bf det}}_{(\C,\A,0)}(\varphi,m^{\varphi})} &  \underline{{\bf det}}_{(\C,\A,0)}(\H)
}
$$
est une équivalence faible de groupoïdes et le foncteur induit:
$$
\xymatrix@R=5pt@C+50pt{
\text{$2$-$h{\bf Grp}$} \ar[r]^-{\pi_0\big(\underline{{\bf det}}_{(\C,\A,0)}(\,\bullet\,)\big)}   &  {\bf Ens}}
$$
est représentable aussi par le $1$-ère $K$-groupe de $(\C,\A,0)$.
\end{proposition}
\begin{proof}
Voir le Corollaire \ref{repsteouni}.
\end{proof}

\renewcommand{\thesubsection}{\S\thesection.\arabic{subsection}}
\subsubsection{}\;
\renewcommand{\thesubsection}{\thesection.\arabic{subsection}}

Soit $\C$ une catégorie de modèles munie d'un objet nul $0$ fixe et d'une famille distinguée $\mathcal{A}$ d'objets de $\C$ contenant l'objet $0$. Si $\G$ est un $2$-groupe (voir la section \ref{2groupeseection}) un \emph{déterminant réduit de $(\C,\A,0)$ á valeurs dans $\G$} est une couple $D=(D,T)$ de fonctions:
$$
\xymatrix{\mathrm{s}^g_{1}(\C,\mathcal{A}) \ar[r]^-{D} & \{\text{Objets de $\G$}\}}\qquad\text{et}\qquad \xymatrix{\mathrm{s}^g_2(\C,\A) \ar[r]^-{T} & \{\text{Morphismes de $\G$}\}}
$$
vérifiant les propriétés:
\begin{enumerate}
\item (Compatibilité) \, Si $\mathfrak{X}$ est un $1$-triangle (à gauche) de $\C$ (à valeurs dans $\A$):
$$
\mathfrak{X} \, = \, \vcenter{\def\objectstyle{\scriptstyle}
\def\labelstyle{\scriptstyle}\xymatrix@-5pt{
0\ar[r] & X_0 \\
X_2 \ar[u]\ar[r] & X_1\ar[u]
}}\,,
$$
$T(\mathfrak{X})$ est un morphisme de $\G$ de la forme: 
$$
\xymatrix@C+15pt{D(X_2) \otimes D(X_0)\ar[r]^-{T(\mathfrak{X})} & D(X_1) }\,.
$$
\item (Unitaire) \, On a que $D\big(\,0\,\big)\,=\,\mathbb{1}_{\G}$ et 
$T\left(\,\vcenter{
\def\objectstyle{\scriptstyle}
\def\labelstyle{\scriptstyle}
\xymatrix@-5pt{
0\ar[r] & 0 \\
0 \ar[u]\ar[r] & 0\ar[u]
}}\,\right)\,=\,l_{\mathbb{1}}^{-1}\,=\,r_{\mathbb{1}}^{-1}$.
\item (Associativité)\, Si $\eta$ est un $2$-triangle (à gauche) de $\C$ (à valeurs dans $\A$):
$$
\eta \, = \,
\def\objectstyle{\scriptstyle}
\def\labelstyle{\scriptstyle}
\vcenter{\xymatrix@C-1.2pc@R-1.8pc{
   &0\ar[r]&A_{3,2}\\
   &&\\
   0\ar[r]&A_{2,1}\ar[r]\ar[uu]&A_{3,1}\ar[uu]\\
   &&\\
   A_{1,0}\ar[uu]\ar[r]&A_{2,0}\ar[r]\ar[uu]&A_{3,0}\ar[uu]}}
$$
on a un diagramme commutatif de $\G$:
$$
\def\objectstyle{\scriptstyle}
\def\labelstyle{\scriptstyle}
\vcenter{\xymatrix@R-10pt@C+45pt{
D(A_{30})   &D(A_{10})\otimes D(A_{31})  \ar[l]_-{T(d_2\,\eta)} \\&\\
&D(A_{10}) \otimes (D(A_{21})\otimes D(A_{32})) \ar[uu]_-{D(A_{10})\otimes T(d_0\, \eta) }  \ar@{}[d]_-{a}|-{\text{\rotatebox[origin=c]{90}{$\cong$}}}\\
D(A_{20})\otimes D(A_{32})  \ar[uuu]^-{T(d_1\,\eta)}&(D(A_{10})\otimes D(A_{21}))\otimes D(A_{32})\ar[l]^-{T(d_3\,\eta)\otimes D(A_{32})}\,,}}
$$
\begin{align*}
\text{où} \qquad 
A_{03} \, = \, & d_{1} d_{1} \eta \, = \, d_{1} d_{2} \eta \qquad
A_{01} \, = \,  d_{2} d_{2} \eta \, = \, d_{2} d_{3} \eta \qquad
A_{13} \, = \,  d_{1} d_{0} \eta \, = \, d_{0} d_{2} \eta \\
A_{02} \, = \, & d_{2} d_{1} \eta \, = \, d_{1} d_{3}  \eta \qquad
A_{23} \, = \,  d_{0} d_{0} \eta \, = \, d_{0} d_{1} \eta \qquad
A_{12} \, = \,  d_{2} d_{0} \eta \, = \, d_{0} d_{3} \eta \,.
\end{align*}
\end{enumerate}

Si $(D,T)$ et $(D',T')$ sont deux déterminants réduits de $(\C,\A,0)$ à valeurs dans $\G$ un \emph{morphisme de déterminants réduits} $H\colon\xymatrix@C-10pt{(D,T)\ar[r] & (D',T')}$ est une fonction:
$$
\xymatrix@C+10pt{\mathrm{s}^g_{1}(\C,\mathcal{A})\ar[r]^-{H} & \big\{\text{Morphismes de $\G$}\big\}}
$$
vérifiant les propriétés:
\begin{enumerate}
\item $H(A)$ est un morphisme en $\G$ de la forme:
$$
\xymatrix@C+10pt{D(A) \ar[r]^-{H(A)} & D'(A)} 
$$
pour tout objet $A$ de $\C$ appartenant à $\A$.
\item On a que $H(s_0 \star) = \mathrm{id}_{\mathbb{1}}$.
\item On a un diagramme commutatif:
$$
\xymatrix@+10pt{
D(d_2\mathfrak{X})\otimes D(d_0 \mathfrak{X})\ar[d]_-{H(d_2\mathfrak{X})\otimes H(d_0\mathfrak{X})}\ar[r]^-{T(\mathfrak{X})} & D(d_1\mathfrak{X})  \ar[d]^-{H(d_1\mathfrak{X})}  \\
D'(d_2\mathfrak{X})\otimes D'(d_0 \mathfrak{X})  \ar[r]_-{T'(\mathfrak{X})} &  D'(d_1\mathfrak{X})\,.}
$$
pour tout $1$-triangle (à gauche) $\mathfrak{X}$ de $\C$ (à valeurs dans $\A$). 
\end{enumerate}

Posons $\pi_0\big(\underline{\bf det}^{red}_{(\C,\A,0)}(\G)\big)$ pour noter l'ensemble des déterminants réduits de $(\C,\A,0)$ à valeurs dans $\G$ soumis à la relation d'équivalence qu'identifie deux déterminants réduits s'il existe un morphisme de déterminants réduits entre eux. 

On définit sans peine un foncteur (voir \ref{detrereduit}):
$$
\vcenter{\xymatrix@R=5pt@C+15pt{
\text{$2$-${\bf Grp}$} \ar[r]  &  {\bf Ens}\\
\G\;\, \ar@{}[r]|-{\longmapsto} & \pi_0\big(\underline{{\bf det}}^{red}_{(\C,\A,0)}(\G)\big)   } }\,.
$$

Il se suit des Corollaires \ref{repsteouni2} et \ref{repsteouni3}:

\begin{proposition}
Soit $\C$ une catégorie de modèles munie d'un objet nul $0$ fixe et d'une famille distinguée $\mathcal{A}$ d'objets de $\C$ contenant l'objet $0$. Si $(\varphi,m^{\varphi})\colon \xymatrix@C-10pt{\G\ar[r]& \H}$ est une $2$-équivalence faible de $2$-groupes la fonction:
$$
\xymatrix@C+55pt{
 \pi_0\big(\underline{{\bf det}}^{red}_{(\C,\A,0)}(\G)\big) \ar[r]^-{\pi_0\big(\underline{{\bf det}}^{red}_{(\C,\A,0)}(\varphi,m^{\varphi})\big)} &  \pi_0\big(\underline{{\bf det}}^{red}_{(\C,\A,0)}(\H)\big)
}
$$
est bijective et le foncteur induit:
$$
\xymatrix@R=5pt@C+60pt{
\text{$2$-$h{\bf Grp}$} \ar[r]^-{\pi_0\big(\underline{{\bf det}}^{red}_{(\C,\A,0)}(\,\bullet\,)\big)}   &  {\bf Ens}}
$$
est représentable par le $1$-ère $K$-groupe (de Kan ou de Segal) de $(\C,\A,0)$, c'est-à-dire par le $2$-groupe d'homotopie (voir le Chapitre 3) du pré-monoïde de Segal:
$$
\vcenter{\xymatrix@C+12pt@R=1pt{
\big(\Delta\times\Delta\big)^{op} \ar[r] & \ens\\ \big([p],[q]\big)\ar@{}[r]|-{\longmapsto} & \mathrm{Hom}_{\simp}\big([p],{\bf W}\mathrm{S}^g_{q}(\C,\mathcal{A}) \big)  }}
$$
ou de l'ensemble simplicial réduit:
$$
\vcenter{\xymatrix@C+12pt@R=1pt{
\Delta^{op} \ar[r] & \ens \,.\\ [p]\ar@{}[r]|-{\longmapsto} & \mathrm{s}^g_{p}(\C,\mathcal{A}) }}
$$
\end{proposition}

\section{Dérivateurs et ses $K$-groupes}

\renewcommand{\thesubsection}{\S\thesection.\arabic{subsection}}
\subsection{}\;
\renewcommand{\thesubsection}{\thesection.\arabic{subsection}}

Une \emph{catégorie de types de diagrammes} est une sous-catégorie pleine ${\bf Dia}$ de la catégorie des petites catégories $\cat$, satisfaisant les propriétés suivantes:
\begin{enumerate}
\item[{\bf Dia 0}] La catégorie augmentée des simplexes ${\bf \Delta}_{+}$ (voir le début de \ref{enssimp}) est une sous-catégorie de ${\bf Dia}$. 
\item[{\bf Dia 1}] ${\bf Dia}$ est stable par sommes finies et par produits fibrés.
\item[{\bf Dia 2}] Pour tout foncteur $u:\xymatrix@-7pt{I\ar[r]&J}$ de ${\bf Dia}$ et tout objet $b$ de $J$, la catégorie  $I\,|\,b$ (resp. $b\,|\,I$) des objets de $I$ au-dessus de $b$ (resp. au-dessous de $b$) appartient à ${\bf Dia}$.
\end{enumerate}


Si $\underline{\cat}$ note la $2$-catégorie des petites catégories, foncteurs et transformations naturelles, une \emph{$2$-catégorie de types diagrammes} est une sous-$2$-catégorie pleine $\underline{\bf Dia}$ de $\underline{\cat}$ dont la catégorie sous-jacente ${\bf Dia}$, \emph{i.e.} cela qu'on obtient de $\underline{\bf Dia}$ en oubliant ses $2$-flèches, est une catégorie de types de diagrammes.

Soit $\underline{\bf Dia}$ une $2$-catégorie de types de diagrammes. Un \emph{pré-dérivateur} de domaine $\underline{\bf Dia}$ est par définition un $2$-foncteur strict de la forme:
$$
\xymatrix@R=3pt@C+25pt{
\; \underline{\bf Dia}^{op} \; \ar[r]^-{\mathbb{D}}& \; \underline{{\bf CAT}} \; \\ 
I\ar@/^13pt/[dddd]^-u\ar@/_13pt/[dddd]_-{v}& \mathbb{D}(I) \\ &\\ 
\overset{\alpha}{\Rightarrow} \ar@{}[r]|-{\longmapsto} & \overset{\alpha_{*}}{\Rightarrow}\\ &\\ J&\mathbb{D}(J)\ar@/^13pt/[uuuu]^-{u^*} \ar@/_13pt/[uuuu]_{v^*}}.
$$ 

Si $\mathbb{D}$ est un pré-dérivateur de domaine $\underline{\bf Dia}$ et $I$ une petite catégorie de ${\bf Dia}$, on appelle les objets de la catégorie $\mathbb{D}(I)$ les $I$-\emph{diagrammes} de $\mathbb{D}$. La catégorie des $e$-diagrammes $\mathbb{D}(e)$ de $\mathbb{D}$, $e$ étant la catégorie ponctuelle, est appelée la \emph{catégorie des coefficients} de $\mathbb{D}$.

Le foncteur \emph{diagramme sous-jacent} de $\mathbb{D}$:
\begin{equation}\label{sousous}  
\xymatrix@+10pt{\mathbb{D}(I) \ar[r]^-{\mathrm{dgm}} & \mathbb{D}(e)^{I} \, = \, \underline{\mathrm{Hom}}_{\bf CAT}\big(I,\mathbb{D}(e)\big)},
\end{equation}
est par définition le foncteur qu'on obtient du composé:
$$
\phantom{a}\vcenter{\xymatrix@+23pt{\; I \; \ar[r] & \; \underline{\mathrm{Hom}}_{\bf Dia}\big(e,I\big) \ar[r]^-{\mathbb{D}} \; & \; \underline{\mathrm{Hom}}_{\bf CAT}\big(\mathbb{D}(I),\mathbb{D}(e)\big) \; }}
$$
$$
\vcenter{\xymatrix@R-8pt{a\ar[d]_{\varphi}\\b}} \qquad \longmapsto \qquad 
\left(\vcenter{\xymatrix{e\ar@/^10pt/[r]^-a \ar@/_10pt/[r]_-{b} \ar@{}[r]|-{\varphi\,\Downarrow } & I }}\right) \qquad \longmapsto \qquad
\left(\vcenter{\xymatrix@C+4pt{\mathbb{D}(I)\ar@/^10pt/[r]^-{a^*} \ar@/_10pt/[r]_-{b^*} \ar@{}[r]|-{\varphi_*\Downarrow } & \mathbb{D}(e) }}\right).
$$

Alors si $\mathfrak{X}$ est un $I$-diagramme de $\mathbb{D}$, le diagramme sous-jacente de $\mathfrak{X}$ est le foncteur de $I$ vers $\mathbb{D}(e)$:
$$
\xymatrix@R=3pt@C+25pt{
\; I \; \ar[r]^-{\mathrm{dgm}(\mathfrak{X})}& \; \mathbb{D}(e) \; \\ 
a\ar[dddd]_-{\varphi}& \mathfrak{X}_{a}=a^*(\mathfrak{X})\ar[dddd]_{\mathfrak{X}_{\varphi}}^-{\varphi^*(\mathfrak{X})}|-{=} \\ &\\ 
 \ar@{}[r]|-{\longmapsto} &\\ &\\ b&\mathfrak{X}_{b}=b^*(\mathfrak{X})}.
$$

\renewcommand{\thesubsection}{\S\thesection.\arabic{subsection}}
\subsubsection{}\; 
\renewcommand{\thesubsection}{\thesection.\arabic{subsection}}

Soit $\underline{\bf Dia}$ une $2$-catégorie de types de diagrammes. Étant donné un pré-dérivateur $\mathbb{D}$ de domaine $\underline{\bf Dia}$ et $u:\xymatrix@C-9pt{I\ar[r]&J}$ un foncteur de ${\bf Dia}$, on dit que $\mathbb{D}$ \emph{admet des extensions de Kan homotopiques à gauche} (resp. \emph{droite}) \emph{le long de $u$}, si le foncteur $u^*:\xymatrix@C-7pt{\mathbb{D}(J)\ar[r]&\mathbb{D}(I)}$ admet un foncteur adjoint à gauche (resp. à droite):
$$
\xymatrix@C+15pt{
\mathbb{D}(J)\ar@/_14pt/[r]_{u^*}\ar@{}[r]|-{\perp}&\mathbb{D}(I)\ar@/_14pt/[l]_{{u_{!}}}}
\qquad\left(\;\text{resp.}\quad
\vcenter{
\xymatrix@C+15pt{
\mathbb{D}(J)\ar@/^14pt/[r]^{u^*}\ar@{}[r]|-{\perp}&\mathbb{D}(J)\ar@/^14pt/[l]^{{u_{*}}}}}\right)\,,
$$
appelé un (ou le) \emph{foncteur extension de Kan homotopique à gauche} (resp. \emph{droite}) \emph{de $\mathbb{D}$ le long de $u$}.

Si $\mathbb{D}$ est un pré-dérivateur qui admet des extensions de Kan homotopiques à gauche (resp. droite) le long du foncteur $p_{I}:\xymatrix@C-10pt{I\ar[r]&e}$, pour $I$ un objet de ${\bf Dia}$:
$$
\xymatrix@C+18pt{
\mathbb{D}(e)\ar@/_14pt/[r]_{p_{I}^*}\ar@{}[r]|-{\perp}&\mathbb{D}(I)\ar@/_14pt/[l]_{{(p_{I})_{!}}\,=\,\mathrm{hocolim}_{I}}}
\qquad\left(\;\text{resp.}\quad
\vcenter{
\xymatrix@C+15pt{
\mathbb{D}(e)\ar@/^14pt/[r]^{p_{I}^*}\ar@{}[r]|-{\perp}&\mathbb{D}(I)\ar@/^14pt/[l]^{{(p_{I})_{*}}\,=\,\mathrm{holim}_{I}}}}\right)\,,
$$
on dit que $\mathbb{D}$ \emph{admet des colimites homotopiques} (resp. \emph{des limites homotopiques}) \emph{de type $I$}, et on appelle l'objet $\mathrm{hocolim}_{I}(\mathfrak{X})$ $\big($resp. $\mathrm{holim}_{I}(\mathfrak{X})\big)$ une (ou la) \emph{colimite homotopique} (resp. la \emph{limite homotopique}) du $I$-diagramme $\mathfrak{X}$ de $\mathbb{D}$. 

On vérifie sans peine:

\begin{lemme}\label{prosommD}
Soit $\mathbb{D}$ un pré-dérivateur de domaine ${\bf Dia}$. Si $I$ est une catégorie discrète finie $I$ de ${\bf Dia}$ telle que le foncteur:
$$
\xymatrix@+5pt{\mathbb{D}(I)\ar[r] & \underset{a\in \mathrm{Ob}(I)}{\prod} \mathbb{D}(e)},
$$ 
induit des foncteurs canoniques $\xymatrix@-7pt{e\ar[r]^-{a}&I}$ est une équivalence de catégories, alors $\mathbb{D}$ admet des colimites homotopiques (resp. limites homotopiques) de type $I$ si et seulement si la catégorie $\mathbb{D}(e)$ admet des sommes (resp. des produits) indexées par l'ensemble d'objet $I_{0}$ de $I$. Dans ce cas la colimite homotopique (resp. la limite homotopique) d'un $I$-diagramme $\mathfrak{X}$ de $\mathbb{D}$ est canoniquement isomorphe à la somme (resp. le produit) de la famille $\{\mathfrak{X}_{a}\}_{a\in I_{0}}$ dans $\mathbb{D}(e)$.
\end{lemme}


\renewcommand{\thesubsection}{\S\thesection.\arabic{subsection}}
\subsubsection{}\; 
\renewcommand{\thesubsection}{\thesection.\arabic{subsection}}

Soit $\underline{\bf Dia}$ une $2$-catégorie de types de diagrammes, et $\mathbb{D}$ un pré-dérivateur de domaine $\underline{\bf Dia}$. Étant donné un foncteur $u:\xymatrix@C-7pt{I\ar[r]&J}$ de ${\bf Dia}$ et un objet $b$ de $J$, on a obtient des transformations naturelles: 
$$
\vcenter{\xymatrix@+10pt{
\mathbb{D}(I\,|\,b)\drtwocell<\omit>{\phantom{aaaa}(\Gamma\,|\,b)^{*}}&\mathbb{D}(I)\ar[l]_-{(\pi\,|\,b)^*}\\ \mathbb{D}(e)\ar[u]^-{(p_{I\,|\,b})^*}&\mathbb{D}(J)\ar[l]^-{b^*}\ar[u]_-{u^*}
}}\qquad\text{et}\qquad
\vcenter{\xymatrix@+10pt{
\mathbb{D}(b\,|\,I)&\mathbb{D}(I)\ar[l]_-{(b\,|\,\pi)^*}\\ \mathbb{D}(e)\ar[u]^-{(p_{b\,|\,I})^*}&\mathbb{D}(J)\ar[l]^-{b^*}\ar[u]_-{u^*}\ultwocell<\omit>{(b\,|\,\Gamma)^{*}\phantom{aaaa}}
}},
$$
en appliquant le $2$-foncteur $\mathbb{D}$ aux transformations naturelles canoniques:
$$
\vcenter{\xymatrix@+8pt{
I\,|\,b\drtwocell<\omit>{\phantom{aa}\Gamma\,|\,b} \ar[d]_-{p_{I\,|\,b}}\ar[r]^-{\pi\,|\,b}&I\ar[d]^-{u}\\ e\ar[r]_-{b}&J
}}\qquad\;\text{et}\;\qquad
\vcenter{\xymatrix@+8pt{
b\,|\,I\ar[d]_-{p_{b\,|\,I}}\ar[r]^-{b\,|\,\pi}&I\ar[d]^-{u}\\ e\ar[r]_-{b}&J\ultwocell<\omit>{b\,|\,\Gamma\phantom{a}}
}}.
$$

En supposons l'existence de foncteurs d'extension de Kan homotopiques à gauche (resp. droite) de $\mathbb{D}$:
\begin{equation}\label{chois3}
\mathbb{D}(e)\vcenter{\xymatrix@C+10pt{
\phantom{\cdot}\ar@{}[r]|-{\perp}\ar@<-4pt>@/_10pt/[r]_-{(p_{I\,|\,b})^*}&
\phantom{\cdot}\ar@<-4pt>@/_10pt/[l]_-{\mathrm{hocolim}_{I\,|\,b}}}}\mathbb{D}(I\,|\,b)\,,
\qquad\quad
\mathbb{D}(J)
\vcenter{\xymatrix@+10pt{
\phantom{\cdot}\ar@{}[r]|-{\perp}\ar@<-4pt>@/_10pt/[r]_-{u^*}&
\phantom{\cdot}\ar@<-4pt>@/_10pt/[l]_-{u_{!}}}}
\mathbb{D}(I)\,,
\end{equation}
$$
\mathbb{D}(e)
\vcenter{\xymatrix@C+10pt{
\phantom{\cdot}\ar@{}[r]|-{\perp}\ar@<+4pt>@/^10pt/[r]^-{(p_{b\,|\,I})^*}&
\phantom{\cdot}\ar@<+4pt>@/^10pt/[l]^-{\mathrm{holim}_{b\,|\,I}}}}
\mathbb{D}(b\,|\,I)
\qquad\text{et}\qquad
\mathbb{D}(J)
\vcenter{\xymatrix@+10pt{
\phantom{\cdot}\ar@{}[r]|-{\perp}\ar@<+4pt>@/^10pt/[r]^-{u^*}&
\phantom{\cdot}\ar@<+4pt>@/^10pt/[l]^-{u_{*}}}}
\mathbb{D}(I)\,;
$$
on construit par le même procédé du paragraphe \ref{evalus}, de transformations naturelles:
$$
\vcenter{\xymatrix@+10pt{
\mathbb{D}(I\,|\,b)\ar[d]_-{\mathrm{hocolim}_{I\,|\,b}}&\mathbb{D}(I)\ar[d]^-{u_{!}}\ar[l]_-{(\pi\,|\,b)^*}\\ \mathbb{D}(e)\urtwocell<\omit>{\phantom{aaaa}(\Gamma\,|\,b)_{!}}&\mathbb{D}(J)\ar[l]^-{b^*}
}}\qquad\text{et}\qquad
\vcenter{\xymatrix@+10pt{
\mathbb{D}(b\,|\,I)\ar[d]_-{\mathrm{holim}_{b\,|\,I}}&\mathbb{D}(I)\ar[l]_-{(b\,|\,\pi)^*}\dltwocell<\omit>{(b\,|\,\Gamma)_{*}\phantom{aaaa}}\ar[d]^-{u_{*}}\\ \mathbb{D}(e)&\mathbb{D}(J)\ar[l]^-{b^*}
}}\,.
$$

Si $\mathfrak{X}$ est un $I$-diagramme de $\mathbb{D}(I)$, on appelle les morphismes:
$$
\vcenter{\xymatrix@+15pt{\mathrm{hocolim}_{I\,|\,b} \big((\pi\,|\,b)^* \, \mathfrak{X} \big)\ar[r]^-{(\Gamma\,|\,b)_{*\,\mathfrak{X}}}&u_{!} (\mathfrak{X})_{b}}}
$$
$$
\text{et}\qquad \vcenter{\xymatrix@+15pt{
u_{*} (\mathfrak{X})_{b}\ar[r]^-{(b\,|\,\Gamma)_{*\,\mathfrak{X}}}&
\mathrm{holim}_{b\,|\,I} \big((b\,|\,\pi)^* \, \mathfrak{X} \big)}},
$$
\emph{les morphisme d'évaluation} de $u_{!}(\mathfrak{X})$ et $u_{*}(\mathfrak{X})$ en $b$.

\renewcommand{\thesubsection}{\S\thesection.\arabic{subsection}}
\subsubsection{}\;
\renewcommand{\thesubsection}{\thesection.\arabic{subsection}}

Soit $\underline{\bf Dia}$ une $2$-catégorie de types de diagrammes. Un pré-dérivateur $\mathbb{D}$ de domaine $\underline{\bf Dia}$ est appelé un \emph{dérivateur} de domaine $\underline{\bf Dia}$ si $\mathbb{D}$ vérifie les propriétés suivants: 
\begin{enumerate}
\item[{\bf Der 1}] Si $I$ est une catégorie discrète finie de ${\bf Dia}$ le foncteur:
$$
\xymatrix@+5pt{\mathbb{D}(I)\ar[r] & \underset{a\in \mathrm{Ob}(I)}{\prod} \mathbb{D}(e)},
$$
induit des foncteurs canoniques $\xymatrix@-7pt{e\ar[r]^-{a}&I}$ est une équivalence de catégories (voir le Lemme \ref{prosommD}).
\item[{\bf Der 2}] Pour tout petite catégorie $I$ de ${\bf Dia}$ le foncteur diagramme sous-jacente:
$$
\xymatrix@+10pt{\mathbb{D}(I) \ar[r]^-{\mathrm{dgm}} & \mathbb{D}(e)^{I}},
$$ 
est conservatif; \emph{i.e.} si $F:\xymatrix@C-12pt{\mathfrak{X}\ar[r]&\mathfrak{X}'}$ est un morphisme de $\mathbb{D}(I)$, tel que pour tout objet $a$ de $I$, le morphisme $F_{a}:\xymatrix@C-12pt{\mathfrak{X}_{a}\ar[r]&\mathfrak{X}'_{a}}$ soit un isomorphisme de $\mathbb{D}(e)$, alors $F$ est un isomorphisme.
\item[{\bf Der 3}] $\mathbb{D}$ admet des extensions de Kan homotopiques à gauche et à droite le long de tout foncteur $u:\xymatrix@C-7pt{I\ar[r]&J}$ de ${\bf Dia}$:
$$
\xymatrix@C+15pt{
\mathbb{D}(J)\ar@/_14pt/[r]_{u^*}\ar@{}[r]|-{\perp}&\mathbb{D}(I)\ar@/_14pt/[l]_{{u_{!}}}}
\qquad\text{et}\qquad
\vcenter{
\xymatrix@C+15pt{
\mathbb{D}(J)\ar@/^14pt/[r]^{u^*}\ar@{}[r]|-{\perp}&\mathbb{D}(J)\ar@/^14pt/[l]^{{u_{*}}}}}.
$$
\item[{\bf Der 4}] Pour tout foncteur $u:\xymatrix@C-7pt{I\ar[r]&J}$ de ${\bf Dia}$, tout objet $b$ de $J$, et tout objet $\mathfrak{X}$ de $\mathbb{D}(I)$, les morphisme d'évaluation de $u_{!}(\mathfrak{X})$ et $u_{*}(\mathfrak{X})$ en $b$:
$$
\vcenter{\xymatrix@+15pt{\mathrm{hocolim}_{I\,|\,b} \big((\pi\,|\,b)^* \, \mathfrak{X} \big)\ar[r]&u_{!} (\mathfrak{X})_{b}}}
$$
$$
\text{et}\qquad \vcenter{\xymatrix@+15pt{
u_{*} (\mathfrak{X})_{b}\ar[r]&
\mathrm{holim}_{b\,|\,I} \big((b\,|\,\pi)^* \, \mathfrak{X} \big)}},
$$
respectivement, sont des isomorphismes. 
\end{enumerate}

On vérifie sans peine que si $\mathbb{D}$ est un dérivateur alors les Corollaires \ref{localco}, \ref{necimescimagg} et \ref{eqcisi} sont valables pour $\mathbb{D}$.

D'un autre remarquons que si $\underline{\bf Fins}$ note la $2$-catégorie de types de diagrammes dont les objets sont les catégories finies directes et inverses,  d'après la Proposition \ref{local}, le Corollaire \ref{hocompco} et les Lemmes \ref{reflee} et \ref{prosomm}, pour toute catégorie de modèles $\C$ le $2$-foncteur canonique:
$$
\xymatrix@R=3pt@C+25pt{
\; \underline{\bf Fins}^{op} \; \ar[r]^-{\mathbb{D}}& \; \underline{{\bf CAT}} \; \\ 
I   \ar@{}[r]|-{\longmapsto} & \C^I[{\bf W}_I^{-1}] }
$$ 
est bien un dérivateur de domaine ${\bf Fins}$.

\renewcommand{\thesubsection}{\S\thesection.\arabic{subsection}}
\subsection{}\;
\renewcommand{\thesubsection}{\thesection.\arabic{subsection}}

Considérons la petite catégorie engendrée par le diagramme:
\begin{equation}\label{carrinombrr}
\boxempty \; = \;
\def\objectstyle{\scriptstyle}
\def\labelstyle{\scriptstyle}
\vcenter{\xymatrix@-10pt{a\ar[r]\ar[d] & b\ar[d] \\ d\ar[r] & c}}\,,
\end{equation}
soumis à la relation qu'identifie les deux flèches de $a$ vers $c$ (voir \ref{homcartcarrff}) et prenons la sous-catégorie $\text{\Large $\lefthalfcap$}$ (resp. $\text{\Large $\righthalfcup$}$) de $\boxempty$:
$$
\text{\Large $\lefthalfcap$} \; = \;
\def\objectstyle{\scriptstyle}
\def\labelstyle{\scriptstyle}
\vcenter{\xymatrix@-20pt{
a\ar[rr]\ar[dd] && b &&&&& a\ar[rr]\ar[dd]&&b\ar[dd]\\
&& &\phantom{.}\ar@{^(->}[rrr]^-{q_{\lefthalfcap}}&&&&&&\\
d&& & & &&& d\ar[rr] & & c}} 
$$
$$
\left(\text{resp.} \qquad 
\text{\Large $\righthalfcup$}  \; = \;
\def\objectstyle{\scriptstyle}
\def\labelstyle{\scriptstyle}
\vcenter{\xymatrix@-20pt{
&& b\ar[dd]&&&&& a\ar[rr]\ar[dd]&&b \ar[dd]\\
&& &\phantom{.}\ar@{^(->}[rrr]^-{q_{\righthalfcup}}&&&&&&\\
d\ar[rr]&& c& & &&& d\ar[rr] & & c}} \right)\;.
$$

Si $\mathbb{D}$ est un dérivateur de domaine une catégorie de types de diagrammes ${\bf Dia}$, considérons des adjonctions:
\begin{equation}\label{laadjfmmf}
\mathbb{D}(\lefthalfcap)
\vcenter{\xymatrix@C+10pt{
\phantom{\cdot}\ar@{}[r]|-{\perp}\ar@<+4pt>@/^10pt/[r]^-{q_{\lefthalfcap\,!}}&
\phantom{\cdot}\ar@<+4pt>@/^10pt/[l]^-{{q_{\lefthalfcap}^*}}}}
\mathbb{D}({\boxempty})
\qquad\qquad
\left(\;\text{resp.}\quad 
\mathbb{D}({\righthalfcup})
\vcenter{\xymatrix@C+10pt{
\phantom{\cdot}\ar@{}[r]|-{\perp}\ar@<-4pt>@/_10pt/[r]_-{{q_{\righthalfcup}^*}}&
\phantom{\cdot}\ar@<-4pt>@/_10pt/[l]_-{q_{\righthalfcup\,*}}}}
\mathbb{D}({\boxempty})
\right).
\end{equation}

Un $\boxempty$-diagramme $\mathfrak{X}$ de $\mathbb{D}$ est dit \emph{un carré homotopiquement cocartésien de $\mathbb{D}$} (resp. {homotopiquement cartésien}) si la valeur en $\mathfrak{X}$ d'une counité $\xymatrix@-10pt{q_{\lefthalfcap\, !}\; {q_{\lefthalfcap}^*}\ar@{=>}[r]^-{\varepsilon} & \mathrm{id}}$ $\big($resp. de l'unité  \xymatrix@-10pt{\mathrm{id}\ar@{=>}[r]^-{\eta}&q_{\righthalfcup}^*\,{q_{\righthalfcup\; *}}}$\big)$ de l'adjonction \eqref{laadjfmmf}:
$$
\vcenter{\xymatrix@+2pt{
q_{\lefthalfcap\,!} \; {q_{\lefthalfcap}^*} (\mathfrak{X})\ar[r]^-{\varepsilon_{\mathfrak{X}}}& \mathfrak{X}}}
\qquad\quad
\left(\,\text{resp.}\quad
\vcenter{\xymatrix@+2pt{
\mathfrak{X} \ar[r]^-{\eta_{\mathfrak{X}}} & q_{\righthalfcup\,*} \; {q_{\righthalfcup}^*} (\mathfrak{X})}}\right)\,.
$$
est un isomorphisme dans la catégorie $\mathbb{D}({\boxempty})$ (voir aussi le Lemme \ref{carrlele}). 

On vérifie sans peine que la preuve du Lemme \ref{exterio} marche bien pour les dérivateurs.

\renewcommand{\thesubsection}{\S\thesection.\arabic{subsection}}
\subsubsection{}\;
\renewcommand{\thesubsection}{\thesection.\arabic{subsection}}

Soit $\mathbb{D}$ un dérivateur et supposons que $0$ est un objet nul fixe de $\mathbb{D}(e)$ et $\mathcal{A}$ est une famille distinguée d'objets de $\mathbb{D}(e)$ contenant l'objet $0$. Si $p\geq 0$ un $(p-1)$-\emph{triangle à gauche} (resp. \emph{à droite}) $\mathfrak{X}$ \emph{de $\mathbb{D}$} (\emph{à valeurs dans $\A$}) est par définition un $\text{\Large $\righthalfcup$}_{p}$-diagramme $\mathfrak{X}$ de $\mathbb{D}$ vérifiant les propriétés suivantes:
\begin{enumerate}
\item $\mathfrak{X}_{i,j}$ appartient à $\A$ pour tout objet $(i,j)$ de $\text{\Large $\righthalfcup$}_{p}$ \emph{i.e.} toujours que $0\leq j \leq  i \leq p$.
\item $\mathfrak{X}_{k,k}= 0$ pour tout $0\leq k\leq p$.
\item Si $0\leq m < k\leq p-1$ et $\alpha_{m,k}\colon\xymatrix@C-10pt{\boxempty \ar[r] & \text{\Large $\righthalfcup$}_p}$ note le foncteur d'inclusion de la sous-catégorie:
$$
\def\objectstyle{\scriptstyle}
\def\labelstyle{\scriptstyle}
\vcenter{\xymatrix@-8pt{(k,m+1)\ar[r]&(k+1,m+1)\\(k,m)\ar[u]\ar[r]&(k+1,m)\ar[u]}}
$$
le $\boxempty$-diagramme $\alpha_{m,k}^\star\mathfrak{X}$ est un carré homotopiquement cocartésien (resp. homotopiquement cartésien) de $\mathbb{D}$. 
\end{enumerate}

Par exemple il y a qu'un seul $(-1)$-triangles à gauche (resp. à droite) de $\mathbb{D}$ qu'on note $0$. D'un autre vu que $\text{\Large $\righthalfcup$}_{1} = \def\objectstyle{\scriptstyle}\def\labelstyle{\scriptstyle}\vcenter{\xymatrix@-5pt{&(1,1)\\(0,0)\ar[r]&\ar[u](1,0)}}$ le foncteur $\xymatrix@C+10pt{\mathbb{D}\big(\text{\Large $\righthalfcup$}_{1}\big) \ar[r]^-{(1,0)^*} & \mathbb{D}(e)}$ induit une équivalence entre la sous-catégorie pleine de $\mathbb{D}\big(\text{\Large $\righthalfcup$}_{1}\big) $ dont les objet sont les $0$-triangle à gauche (resp. à droite) de $\mathbb{D}$ et la sous-catégorie pleine de $\mathbb{D}(e)$ dont les objets sont les éléments de la famille $\A$. 

De façon analogue le foncteur $\xymatrix@C+10pt{\mathbb{D}\big(\text{\Large $\righthalfcup$}_{2}\big) \ar[r]^-{\alpha_{1,0}^*} & \mathbb{D}(\boxempty)}$ (voir (iii) ci-dessus) 
induit une équivalence entre la sous-catégorie pleine de $\mathbb{D}\big(\text{\Large $\righthalfcup$}_{2}\big)$ dont les objets sont les $1$-triangle à gauche (resp. à droite) de $\mathbb{D}$ et la sous-catégorie pleine de $\mathbb{D}(\boxempty)$ dont les objets sont les carrés $\mathfrak{X}$ homotopiquement cocartésien (resp. homotopiquement cartésien) de $\mathbb{D}$ tels que ( voir \eqref{carrinombrr}) $\mathfrak{X}_{a},\mathfrak{X}_{b}$ et $\mathfrak{X}_{c}$ sont des éléments de $\A$ et $\mathfrak{X}_{d}=0$.

On note $\mathrm{s}^g_{p}(\mathbb{D},\A)$ $\big($resp. $\mathrm{s}^d_{p}(\mathbb{D},\A)\big)$ l'ensemble des $(p-1)$-triangles à gauche (resp. à droite) de $\mathbb{D}$ à valeurs dans $\A$. Le $n$-\emph{ème} $K$-\emph{groupe} de $(\mathbb{D},\A,0)$ est par définition la classe de l'ensemble simplicial réduit (voir le foncteur):
\begin{equation}\label{sssddd}
\vcenter{\xymatrix@C+8pt@R=1pt{
\Delta^{op} \ar[r]^-{\mathrm{s}^g_{\bullet}(\C,\mathcal{A})} & \ens\\ [p]\ar@{}[r]|-{\longmapsto} & \mathrm{s}^g_{p}(\C,\mathcal{A})}}
\end{equation}
dans l'ensemble des classes à isomorphismes près des objets de la catégorie homotopique $\mathrm{Ho}_{n+1}(\simp_0)$ des $(n+1)$-groupes de Kan (voir après le Corollaire \ref{dimminus}).

\renewcommand{\thesubsection}{\S\thesection.\arabic{subsection}}
\subsubsection{}\;
\renewcommand{\thesubsection}{\thesection.\arabic{subsection}}

Soit $\mathbb{D}$ un dérivateur, $0$ un objet nul fixe de $\mathbb{D}(e)$ et $\mathcal{A}$ une famille distinguée d'objets de $\mathbb{D}(e)$ contenant l'objet $0$. Si $\G$ est un $2$-groupe (voir la section \ref{2groupeseection}) un \emph{déterminant réduit de $(\mathbb{D},\A,0)$ á valeurs dans $\G$} est une couple $(D,T)$ de fonctions:
$$
\xymatrix{\mathrm{s}^g_{1}(\mathbb{D},\mathcal{A}) \ar[r]^-{D} & \{\text{Objets de $\G$}\}}\qquad\text{et}\qquad \xymatrix{\mathrm{s}^g_2(\mathbb{D},\A) \ar[r]^-{T} & \{\text{Morphismes de $\G$}\}}
$$
vérifiant les propriétés:
\begin{enumerate}
\item (Compatibilité) \, Si $\mathfrak{X}$ est un $1$-triangle (à gauche) de $\mathbb{D}$ (à valeurs dans $\A$) le morphisme $T(\mathfrak{X})$ de $\G$ est une flèche de la forme: 
$$
\xymatrix@C+15pt{D\big(d_2\mathfrak{X}\big) \otimes D\big(d_0\mathfrak{X}\big)\ar[r]^-{T(\mathfrak{X})} & D(d_1\mathfrak{X}) }\,.
$$
\item (Unitaire) \, On a que $D\big(s_0 \, 0\big)\,=\,\mathbb{1}_{\G}$ et  $D\big(s_0s_0 \, 0\big)\,=\,l_{\mathbb{1}}^{-1}\,=\,r_{\mathbb{1}}^{-1}$.
\item (Associativité)\, Si $\eta$ est un $2$-triangle (à gauche) de $\mathbb{D}$ (à valeurs dans $\A$) on a un diagramme commutatif de $\G$:
$$
\def\objectstyle{\scriptstyle}
\def\labelstyle{\scriptstyle}
\vcenter{\xymatrix@R-10pt@C+45pt{
D(A_{30})   &D(A_{10})\otimes D(A_{31})  \ar[l]_-{T(d_2\,\eta)} \\&\\
&D(A_{10}) \otimes (D(A_{21})\otimes D(A_{32})) \ar[uu]_-{D(A_{10})\otimes T(d_0\, \eta) }  \ar@{}[d]_-{a}|-{\text{\rotatebox[origin=c]{90}{$\cong$}}}\\
D(A_{20})\otimes D(A_{32})  \ar[uuu]^-{T(d_1\,\eta)}&(D(A_{10})\otimes D(A_{21}))\otimes D(A_{32})\ar[l]^-{T(d_3\,\eta)\otimes D(A_{32})}\,,}}
$$
\begin{align*}
\text{où} \qquad 
A_{03} \, = \, & d_{1} d_{1} \eta \, = \, d_{1} d_{2} \eta \qquad
A_{01} \, = \,  d_{2} d_{2} \eta \, = \, d_{2} d_{3} \eta \qquad
A_{13} \, = \,  d_{1} d_{0} \eta \, = \, d_{0} d_{2} \eta \\
A_{02} \, = \, & d_{2} d_{1} \eta \, = \, d_{1} d_{3}  \eta \qquad
A_{23} \, = \,  d_{0} d_{0} \eta \, = \, d_{0} d_{1} \eta \qquad
A_{12} \, = \,  d_{2} d_{0} \eta \, = \, d_{0} d_{3} \eta \,.
\end{align*}
\end{enumerate}

Si $(D,T)$ et $(D',T')$ sont deux déterminants réduits de $(\mathbb{D},\A,0)$ à valeurs dans $\G$ un \emph{morphisme de déterminants réduits} $H\colon\xymatrix@C-10pt{(D,T)\ar[r] & (D',T')}$ est une fonction:
$$
\xymatrix@C+10pt{\mathrm{s}^g_{1}(\mathbb{D},\mathcal{A})\ar[r]^-{H} & \big\{\text{Morphismes de $\G$}\big\}}
$$
vérifiant les propriétés:
\begin{enumerate}
\item $H(A)$ est un morphisme en $\G$ de la forme:
$$
\xymatrix@C+10pt{D(A) \ar[r]^-{H(A)} & D'(A)} 
$$
pour tout $0$-triangle $A$ de $\mathbb{D}$ à valeurs dans $\A$.
\item On a que $H(s_0 \star) = \mathrm{id}_{\mathbb{1}}$.
\item On a un diagramme commutatif:
$$
\xymatrix@+10pt{
D(d_2\mathfrak{X})\otimes D(d_0 \mathfrak{X})\ar[d]_-{H(d_2\mathfrak{X})\otimes H(d_0\mathfrak{X})}\ar[r]^-{T(\mathfrak{X})} & D(d_1\mathfrak{X})  \ar[d]^-{H(d_1\mathfrak{X})}  \\
D'(d_2\mathfrak{X})\otimes D'(d_0 \mathfrak{X})  \ar[r]_-{T'(\mathfrak{X})} &  D'(d_1\mathfrak{X})\,.}
$$
pour tout $1$-triangle (à gauche) $\mathfrak{X}$ de $\mathbb{D}$ (à valeurs dans $\A$). 
\end{enumerate}

On pose  $\pi_0\big(\underline{\bf det}^{red}_{(\mathbb{D},\A,0)}(\G)\big)$ pour noter l'ensemble des déterminants réduits de $(\mathbb{D},\A,0)$ à valeurs dans $\G$ soumis à la relation d'équivalence qu'identifie deux déterminants réduits s'il existe un morphisme de déterminants réduits entre eux. Considérons le foncteur (voir \ref{detrereduit}):
$$
\vcenter{\xymatrix@R=5pt@C+15pt{
\text{$2$-${\bf Grp}$} \ar[r]  &  {\bf Ens}\\
\G\;\, \ar@{}[r]|-{\longmapsto} & \pi_0\big(\underline{{\bf det}}^{red}_{(\mathbb{D},\A,0)}(\G)\big)   } }\,.
$$

D'après les Corollaires \ref{repsteouni2} et \ref{repsteouni3}:

\begin{proposition}
Soit $\mathbb{D}$ un dérivateur, $0$ un objet nul fixe de $\mathbb{D}(e)$ et $\mathcal{A}$ une famille distinguée d'objets de $\mathbb{D}(e)$ contenant l'objet $0$. Si $(\varphi,m^{\varphi})\colon \xymatrix@C-10pt{\G\ar[r]& \H}$ est une $2$-équivalence faible de $2$-groupes, la fonction:
$$
\xymatrix@C+50pt{
 \pi_0\big(\underline{{\bf det}}^{red}_{(\mathbb{D},\A,0)}(\G)\big) \ar[r]^-{\pi_0\big(\underline{{\bf det}}^{red}_{(\mathbb{D},\A,0)}(\varphi,m^{\varphi})\big)} &  \pi_0\big(\underline{{\bf det}}^{red}_{(\mathbb{D},\A,0)}(\H)\big)}
$$
est bijective et le foncteur induit:
\begin{equation}\label{fonctrepfffin}
\xymatrix@R=5pt@C+50pt{
\text{$2$-$h{\bf Grp}$} \ar[r]^-{\pi_0\big(\underline{{\bf det}}^{red}_{(\mathbb{D},\A,0)}(\,\bullet\,)\big)}   &  {\bf Ens}}
\end{equation}
est représentable par le $1$-ère $K$-groupe de $(\mathbb{D},\A,0)$, plus précisément \eqref{fonctrepfffin} est représentable par le $2$-groupe d'homotopie (voir le Chapitre 3) de l'ensemble simplicial réduit:
$$
\vcenter{\xymatrix@C+12pt@R=1pt{
\Delta^{op} \ar[r] & \ens \,.\\ [p]\ar@{}[r]|-{\longmapsto} & \mathrm{s}^g_{p}(\mathbb{D},\mathcal{A}) \;.}}
$$
\end{proposition}

\chapter{Les $n$-types d'homotopie}  \label{chaptersimp}   
\setcounter{section}{7}

\section{Ensembles simpliciaux}\label{enssimp}

\renewcommand{\thesubsection}{\S\thesection.\arabic{subsection}}
\subsection{}\;
\renewcommand{\thesubsection}{\thesection.\arabic{subsection}} 

Soit $\cat$ la catégorie des petites catégories et foncteurs. On identifie les ensembles pré-ordonnés avec les catégories qui satisfont la propriété suivante: L'ensemble des morphismes entre deux quelconques objets, est soit vide soit un ensemble avec un seul objet. On trouve alors que les foncteurs entre telles catégories, sont les fonctions qui respectent l'ordre.

Notons ${\bf\Delta}$ la \emph{catégorie des simplexes}; c'est-à-dire, la sous-catégorie pleine de $\cat$ dont les objets sont les catégories $[n]=\{0<\dots<n\}$ pour $n\geq 0$. De fa\c con équivalente (voir par exemple II.2.2 de \cite{GZ}), ${\bf\Delta}$ est la catégorie libre engendrée par les morphismes faces et dégénérescences:
$$
\Big\{\xymatrix@C-4pt{[n]\ar[r]^-{\delta^n_{i}}&[n+1]}\,\Big|\,0\leq i \leq n+1 \Big\}_{n\geq 0}
\quad\,\text{et}\,\quad
\Big\{\xymatrix@C-4pt{[n+1]\ar[r]^-{\sigma^n_{i}}&[n]}\,\Big|\,0\leq i \leq n \Big\}_{n\geq 0},
$$
$$\text{soumis aux rélations définies pour $n\geq 0$:}$$
\begin{equation}\label{relsim}
\begin{split}
\delta^{n+1}_{l} \delta^{n}_{k}\,=\; & \delta^{n+1}_{k}\delta^{n}_{l-1} \;\quad \quad\quad  0\leq k< l \leq n+2,\\
\sigma^{n}_{l}\sigma^{n+1}_{k} \,=\; & \sigma^{n}_{k}\sigma^{n+1}_{l+1} \quad\quad\quad  0\leq k\leq l\leq n ,\\
\sigma^{n+1}_{k} \delta^{n+1}_{l}\,=\; & \begin{cases}
\delta^{n}_{l}\sigma^{n}_{k-1} &\quad\text{si\; $0\leq l<k \leq n+1$},\\
\mathrm{id}_{[n]}       &\quad\text{si \;$0\leq k \leq l \leq k+1\leq n+2$},\\
\delta^{n}_{l-1}\sigma^{n}_{k}    &\quad\text{si \;$0\leq k< l-1\leq n+1$}.
\end{cases}
\end{split}
\end{equation}

On note simplement: 
$${\bf\Delta}\;=\;\;
\xymatrix{
**[l]        [0] \ar@{<-}@/_1pc/[r]|-{\sigma^0_0}&
**[r]        [1] \,\ar@{<-}@/_2pc/[l]|-{\delta^0_1}
                \ar@{<-}@/_1pc/[l]|-{\delta^0_0} \ar@{}[r]|(.45){\dots\dots\dots\dots}&
**[l]        [n] \ar@{<-}@/_1pc/[r]|-{\sigma^n_{n}}_{\underset{\vdots}{\phantom{a}}}
                \ar@{<-}@/_2.8pc/[r]|-{\sigma^n_0}^{\overset{\vdots}{\phantom{a}}}&
             [n+1] \ar@{<-}@/_3pc/[l]|-{\delta^n_{n+1}}^{\underset{\vdots}{\phantom{a}}}
                \ar@{<-}@/_1pc/[l]|-{\delta^n_0}_{\overset{\vdots}{\phantom{a}}} \ar@{}[r]|-{\dots\dots\dots\dots}&}.$$

La \emph{catégorie augmentée des simplexes} ${\bf \Delta}_{+}$, est la sous-catégorie pleine de $\cat$ dont les objets sont les catégories $[n]_{+}=\{0<\dots<n-1\}$, pour $n\geq 0$; c'est-à-dire, $[n]_{+}=[n-1]$ si $n\geq 1$, et $[0]_{+}$ est la catégorie vide. On vérifie en particulier, que ${\bf \Delta}_{+}$ est obtenue de la catégorie des simplexes en ajoutant un objet initial. 

Donc, ${\bf \Delta}_{+}$ est la catégorie libre engendrée par l'unique fonction $\xymatrix@C-5pt{[0]_{+}\ar[r]^{\delta^{-1}_{0}}&[1]_{+}}$, plus les morphismes faces et dégénérescences:
$$
\Big\{\xymatrix@C-4pt{[n+1]_{+}=[n]\ar[r]^-{\delta^n_{i}}&[n+1]=[n+2]_{+}}\,\Big|\,0\leq i \leq n+1 \Big\}_{n\geq 0}
$$
$$
\text{et}\qquad\Big\{\xymatrix@C-4pt{[n+2]_{+}=[n+1]\ar[r]^-{\sigma^n_{i}}&[n]=[n+1]_{+}}\,\Big|\,0\leq i \leq n \Big\}_{n\geq 0};
$$
$$\text{soumis aux rélations définies pour $n\geq 0$:}$$
\begin{equation}\label{relsim+}
\begin{split}
\delta^n_{l} \delta^{n-1}_{k}\,=\; & \delta^n_{k}\delta^{n-1}_{l-1} \;\quad \quad\quad 0\leq k< l \leq n+1,\\
\sigma^n_{l}\sigma^{n+1}_{k} \,=\; & \sigma^n_{k}\sigma^{n+1}_{l+1} \quad\quad\quad 0\leq k\leq l\leq n,\\
\sigma^n_{k} \delta^{n}_{l}\,=\; & \begin{cases}
\delta^{n-1}_{l}\sigma^{n-1}_{k-1} &\,\text{si $0\leq l< k\leq n$},\\
\mathrm{id}_{[n]}       &\;\text{si $0\leq k \leq l \leq k+1\leq  n+1$},\\
\delta^{n-1}_{l-1}\sigma^{n-1}_{k}    &\;\text{si $0\leq k< l-1\leq n$}.
\end{cases}
\end{split}
\end{equation}

$${\bf\Delta}_{+}\;=\;\;
\xymatrix{
[0]_{+}\ar[r]^-{\delta^{-1}_{0}}&
**[l]        [1]_{+} \ar@{<-}@/_1pc/[r]|-{\sigma^0_0}&
**[r]        [2]_{+} \,\ar@{<-}@/_2pc/[l]|-{\delta^0_1}
                \ar@{<-}@/_1pc/[l]|-{\delta^0_0} \ar@{}[r]|(.45){\dots\dots\dots\dots}&
**[l]        [n+1]_{+} \ar@{<-}@/_1pc/[r]|-{\sigma^n_{n}}
                \ar@{<-}@/_2.8pc/[r]|-{\sigma^n_0}^\vdots&
             [n+2]_{+} \ar@{<-}@/_3pc/[l]|-{\delta^n_{n+1}}^{\vdots}
                \ar@{<-}@/_1pc/[l]|-{\delta^n_0} \ar@{}[r]|-{\dots\dots\dots\dots}&}.$$

Rappelons que la catégorie augmentée des simplexes ${\bf\Delta}_{+}$ a une structure de catégorie monoïdale stricte, où:
\begin{equation}\label{delta+}
\xymatrix@C-8pt@R=1pt{
\;{\bf\Delta}_{+}\times  {\bf\Delta}_{+} \;\ar[r]^-{\otimes}& \;{\bf\Delta}_{+}\;\\
\quad\big([n]_{+},[m]_{+}\big)\quad\ar@{|->}[r]&\quad[n+m]_{+}\quad}
\end{equation}
est défini pour deux morphismes $\xymatrix@C-8pt{[n]_{+}\ar[r]^f&[n']_{+}}$ et $\xymatrix@C-8pt{[m]_{+}\ar[r]^g&[m']_{+}}$, par la règle:
\begin{equation}\label{definidelta}
\big(f\otimes g\big) (k) =\begin{cases}
f(k)   &\; \text{si $0\leq k\leq n-1$}\\
g(k-n)+n'   &\; \text{si $n\leq k\leq n+m-1$\,.}
\end{cases}
\end{equation}

\renewcommand{\thesubsection}{\S\thesection.\arabic{subsection}}
\subsection{}\;
\renewcommand{\thesubsection}{\thesection.\arabic{subsection}} 

Notons $\simp$ la catégorie des \emph{ensembles simpliciaux}, \emph{i.e.} des foncteurs $\xymatrix@C-6pt{X:{\bf \Delta}^{op}\ar[r]&
\ens}$ et les transformations naturelles entre eux. 

Si $X$ est un ensemble simplicial, l'ensemble des \emph{$n$-simplexes de $X$} est noté $X_{n}=X\big([n]\big)$. Les morphismes face et dé\-géné\-re\-scen\-ce de la catégorie ${\bf \Delta}$, induisent les \emph{morphismes face et dé\-géné\-re\-scen\-ce de $X$}:
$$
\xymatrix@C-7pt{X_{n+1}\ar[r]^-{d^n_{i}}&X_{n}}\quad \text{et}\quad\xymatrix@C-7pt{X_{n}\ar[r]^-{s^n_{i}}&X_{n+1}}, \quad\text{respectivement.}
$$

Un $n$-simplexe $x$ d'un ensemble simplicial $X$ est dit \emph{dégénéré}, si $x$ est dans l'image d'un des morphismes $s^{n-1}_{0},\dots, s^{n-1}_{n-1}$; de façon équivalente, $x\in X_n$ es dégénéré s'il existe $\xymatrix@C-8pt{[n]\ar[r]^f&[m]}$ surjectif en objets et $y\in X_m$ tel que $f^*(y)=x$. Étant donné un $0$-simplexe $a$ de $X$, on écrit aussi:
$$
a \,= \, s_{0}^{n-1}\cdots\, s_{0}^0(a)\,\in\,X_{n}\qquad \text{où}\;\, n\geq 1,
$$
pour noter le $n$-simplexe dégénéré associé à $a$, et on l'appelle un $n$-simplexe \emph{totalement dégénéré} de $X$.

Si $m\geq 0$, on désigne par $\Delta^{m}$ l'ensemble simplicial représenté par $[m]$, \emph{i.e.} $\Delta^m_{n}$ est égal à l'ensemble des foncteurs de $[n]$ vers $[m]$. On note par abus:
$$
\xymatrix@C-7pt{\Delta^{m}\ar[r]^-{\delta^m_{i}}&\Delta^{m+1}}\quad \text{et}\quad\xymatrix@C-7pt{\Delta^{m+1}\ar[r]^-{\sigma^m_{i}}&\Delta^{m}},
$$
les morphismes d'ensembles simpliciaux définis pour $p\geq 0$ par les règles :
$$
\xymatrix@C+5pt@R=1pt{\Delta^m_{p}\ar[r]^-{(\delta^m_{i})_{p}}&\Delta^{m+1}_{p}\\
\varphi \ar@{}[r]|-{\longmapsto} & \delta_{i}^m\circ \varphi }\quad \text{et}\quad
\xymatrix@C+5pt@R=1pt{\Delta^{m+1}_{p}\ar[r]^-{(\sigma^m_{i})_{p}}&\Delta^{m}_{p}\\
\varphi \ar@{}[r]|-{\longmapsto} & \sigma_{i}^m\circ \varphi },\qquad \text{respectivement}.
$$

Enfin, rappelons que la catégorie des ensembles simpliciaux $\simp$ admet une structure canonique de catégorie cartésienne fermée (voir IV\S6 de \cite{lane}), où l'objet final est l'ensemble simplicial constant à valeur $\star$ (un ensemble ponctuel fixe), le produit d'ensembles simpliciaux est le produit cartésien argument par argument et l'adjoint du produit:
\begin{equation}\label{deltacarte}
\simp
\vcenter{\xymatrix@C+15pt{
\phantom{a}\ar@{}[r]|{\perp}\ar@<-4pt>@/_10pt/[r]_{\underline{\mathrm{Hom}}_{\simp}(X,\,\cdot\,)}&\phantom{a}
\ar@<-4pt>@/_10pt/[l]_-{X \times\,\cdot\, }}}
\simp 
\qquad\text{est définie par}\quad 
\underline{\mathrm{Hom}}_{\simp}(X,Y)_{n} = \mathrm{Hom}_{\simp}(X\times\Delta^n,Y)\,,
\end{equation}
pour tout ensemble simplicial $X$.

\renewcommand{\thesubsection}{\S\thesection.\arabic{subsection}}
\subsection{}\; \label{deltalambda}
\renewcommand{\thesubsection}{\thesection.\arabic{subsection}}

Si $m\geq -1$, on va définir \emph{le bord  $\partial\Delta^{m+1}$ de $\Delta^{m+1}$} comme le conoyau dans $\simp$ de la double flèche: 
$$
\xymatrix@C-1pt{
\underset{0\leq i<j\leq m+1}{\bigsqcup}\Delta^{m-1}
\ar@<+7pt>[rrr]^-{\underset{i<j}{\sqcup} \; \nu_{i}\,\circ\,\delta^{m-1}_{j-1}}
\ar@<-7pt>[rrr]_-{\underset{i<j}{\sqcup} \;\nu_{j}\,\circ\,\delta^{m-1}_{i}}&&&
\underset{0\leq l\leq m+1}{\bigsqcup}\Delta^{m}};
$$
où $\nu_{l}$ note l'inclusion de la $l$-ième composante $\xymatrix@-5pt{\Delta^m\;\ar@{^(->}[r]&\underset{0\leq l\leq m+1}{\bigsqcup}\Delta^{m}}$ $\big($ $\Delta^{-1}=\Delta^{-2}$ étant l'ensemble simplicial vide$\big)$.
 
De la même manière, pour $0\leq k\leq m+1$ on définit \emph{le $k$-horn $\Lambda^{m+1,k}$ de $\Delta^{m+1}$} comme le conoyau de la double flèche: 
$$
\xymatrix@C-1pt{
\underset{i,j\neq k}{\underset{0\leq i<j\leq m+1}{\bigsqcup}}\Delta^{m-1}
\ar@<+7pt>[rrr]^-{\underset{i<j}{\sqcup} \; \nu_{i}\,\circ\,\delta^{m-1}_{j-1}}
\ar@<-7pt>[rrr]_-{\underset{i<j}{\sqcup} \;\nu_{j}\,\circ\,\delta^{m-1}_{i}}&&&
\underset{l\neq k}{\underset{0\leq l\leq m+1}{\bigsqcup}}\Delta^{m}&}.
$$ 

On construit en particulier un triangle commutatif de monomorphismes d'ensembles simpliciaux:
\begin{equation}\label{sigma}
\xymatrix{\Lambda^{m+1,k} \ar@<-3pt>@/_10pt/[rr]_{\alpha^{m,k}}\ar[r]^{\widetilde{\alpha}^{m,k}}&\partial\Delta^{m+1} \ar[r]^{\alpha^m}& 
\Delta^{m+1}},
\end{equation}
où le morphisme $\alpha^m$ est induit à partir de: 
$$
\xymatrix@C+10pt{
\underset{0\leq l\leq m+1}{\bigsqcup}\Delta^{m}\ar[r]^-{\underset{l}{\sqcup} \;\delta^{m}_{i}}&
\Delta^{m+1}},
$$
et le morphisme $\widetilde{\alpha}^{m,k}$ découle de l'inclusion évidente:
$$
\vcenter{\xymatrix@C+10pt{
\underset{l\neq k}{\underset{0\leq l\leq m+1}{\bigsqcup}}\Delta^{m}\ar[r]&
\underset{0\leq l\leq m+1}{\bigsqcup}\Delta^{m}}}.
$$

Donc, pour tout ensemble simplicial $X$ on a des fonctions:
\begin{equation}\label{sigmaX}
\xymatrix@C+10pt{
\mathrm{Hom}_{\simp}\big(\Delta^{m+1},X\big)
\ar@<-3pt>@/_16pt/[rr]_{\alpha^{m,k}_{X}}\ar[r]^{\alpha^{m}_{X}}&
\mathrm{Hom}_{\simp}\big(\partial\Delta^{m+1},X\big)\ar[r]^{\widetilde{\alpha}^{m,k}_{X}}& 
\mathrm{Hom}_{\simp}\big(\Lambda^{m+1,k},X\big)}\,.
\end{equation}

De façon plus explicite, si on considère les isomorphismes canoniques qu'il découle du Lemme de Yoneda\footnote{Le cas $m=-1$ est compris si on définit $X_{-1}$ comme l'ensemble ponctuel, \emph{i.e.} $X_{-1}=\mathrm{Hom}_{\simp}(\Delta^{-1},X)$.}:
\begin{equation}
\mathrm{Hom}_{\simp}\big(\Delta^{m+1},X\big)\cong X_{m+1},
\end{equation}
\begin{equation}\label{bord}
\mathrm{Hom}_{\simp}\big(\partial\Delta^{m+1},X\big)\cong \Bigg\{
\text{\scriptsize{$\big(a_{0},\dots,a_{m+1}\big)$}}\;\Bigg|\;
\vcenter{\xymatrix@R=1pt{\text{\scriptsize{$a_{i}\in X_{m}$ \; et\; $d^{m-1}_{i}a_{j}=d^{m-1}_{j-1}a_{i}$}}\\ \text{\scriptsize{si $0\leq i<j\leq m+1
$.}}}}\Bigg\}\qquad\text{et}
\end{equation}
\begin{equation}\label{bord2}
\mathrm{Hom}_{\simp}\big(\Lambda^{m+1,k},X\big)\cong \Bigg\{
\text{\scriptsize{$\big(a_{0},\dots,a_{k-1},a_{k+1}\dots,a_{m+1}\big)$}}\;\Bigg|\;
\vcenter{\xymatrix@R=1pt{\text{\scriptsize{$a_{i}\in X_{m}$ \; et\; $d^{m-1}_{i}a_{j}=d^{m-1}_{j-1}a_{i}$}}\\ \text{\scriptsize{si $0\leq i<j\leq m+1$ et 
$i,j\neq k$.}}}} \Bigg\}\,,
\end{equation}
alors on vérifie que les fonctions du diagramme \eqref{sigmaX}, sont données par les règles suivantes:
\begin{equation}\label{sigm}
\alpha^m_{X}(a)\,=\,\big(d^m_{0}a,\dots,d^m_{m+1}a\big),
\end{equation}
\begin{equation}\label{sigm2}
\widetilde{\alpha}^{m,k}_{X}(a_{0},\dots,a_{m+1})\,=\,\big(a_{0},\dots,a_{k-1},a_{k+1},\dots,a_{m+1}\big)
\end{equation}
\begin{equation}\label{sigm3}
\text{et}\qquad\alpha^{m,k}_{X}(a)\,=\,\big(d^m_{0}a,\dots,d^m_{k-1}a,d^m_{k+1}a,\dots, d^m_{m}a\big).
\end{equation}

On aura besoin plus loin du Lemme suivant:

\begin{lemme}\label{glennlemme}
Soit $X$ un ensemble simplicial quelconque $m\geq 0$ et $0\leq k \leq m+1$, alors dans le diagramme \eqref{sigmaX} ci-dessus: 
\begin{enumerate}
\item 
Si la fonction $\alpha^{m,k}_{X}$ est injective, la fonction $\alpha^{m}_{X}$ est aussi injective.
\item
Si la fonction $\alpha^{m-1}_{X}$ est surjective (resp. injective), la fonction $\widetilde{\alpha}^{m,k}_{X}$ est aussi surjective (resp. injective). $($Voir le Lemme 1.7.1 de \cite{glenn}$)$.
\item
Si la fonction $\alpha^{m,k}_{X}$ est surjective et la fonction $\widetilde{\alpha}^{m,k}_{X}$ injective, alors la fonction ${\alpha}^{m}_{X}$ est surjective.
\end{enumerate} 
\end{lemme}
\begin{proof}
Soit $X$ un ensemble simplicial, $m\geq 0$ et $0\leq k \leq m+1$.

La preuve de (i) est immédiate. Pour montrer (ii) considérons un élément $a$ de l'ensemble \eqref{bord2}, c'est-à-dire:
\begin{align*}
a= \big(a_{0},\dots,a_{k-1},a_{k+1},\dots,a_{m+1}\big)\;\in \,  \overset{m+1}{\underset{\underset{i\neq k}{i=0}}{\prod}} X_{t}
\qquad\text{où}\;\;\\
d_{i}^{m-1}a_{j}=d^{m-1}_{j-1}a_{i}\quad \text{si}\quad 0\leq i<j \leq m+1\quad\text{et} \quad  i,j\neq k.
\end{align*}

Il n'est pas difficile de constater que si on pose $(b_{0},\dots,b_{m})$ pour noter:
$$
\Big(d_{k-1}^{m-1}a_{0},\dots,d_{k-1}^{m-1}a_{k-1},d_{k}^{m-1}a_{k+1},\dots, d_{k}^{m-1}a_{m+1}\Big)\;\in\, \overset{m}{\underset{\underset{i\neq k}{i=0}}{\prod}} X_{m-1};
$$
alors $d_{i}^{m-2}b_{j}=d^{m-2}_{j-1}b_{i}$ pour $0\leq i<j \leq m$.

Donc, vu que $\alpha^{m-1}_{X}$ est surjective par hypothèse, il existe $a_{k}\in X_{m}$ tel que:
$$
\text{$d^{m-1}_{i}a_{k}=d_{k-1}^{m-1}a_{i}$\quad si\;\, $0\leq i <k$ \quad et\quad  $d^{m-1}_{j-1}a_{k}=d_{k}^{m-1}a_{j}$ \quad si\;\, $0\leq k< j$};
$$
c'est-à-dire, $\widetilde{\alpha}^{m,k}_{X}$ est surjective.

Pour montrer maintenant (iii), considérons un élément $y$ de l'ensemble d'arrivé de la fonction $\alpha^m_{X}$. Si $x$ est un $(m+1)$-simplexe de $X$ tel que $\alpha^{m,k}_{X}(x)=\widetilde{\alpha}^{m,k}_{X}(y)$; vu que $\alpha^{m,k}_{X}(x)=\widetilde{\alpha}^{m,k}_{X} \big(\alpha^m_{X}(x)\big)$, on déduit que  $\alpha^m_{X}(x)$ et $y$ sont envoyés au même élément par la fonction injective $\widetilde{\alpha}^{m,k}_{X}$. Donc, $\alpha^m_{X}(x)=y$; autrement dit, $\alpha^m_{X}$ est une fonction surjective.
\end{proof}

Si $m\geq 0$ on dit qu'un ensemble simplicial \emph{$X$ satisfait la condition d'extension de Kan en dimension $m$} (resp. \emph{satisfait la condition d'extension de Kan de façon stricte en dimension $m$}), si la fonction $\alpha^{m,k}_{X}$ du diagramme \eqref{sigmaX} est surjective (resp. bijective) pour tout $0\leq k\leq m+1$. On vérifie sans peine que tout ensemble simplicial satisfait la condition d'extension de Kan en dimension $0$.

Un \emph{complexe de Kan} est un ensemble simplicial qui satisfait la condition d'extension de Kan en dimension $m$ pour tout $m\geq 0$.  Explicitement, un ensemble simplicial $X$ est un complexe de Kan si pour tout carré commutatif de $\simp$:
\begin{equation}\label{complexeKKAN}
\xymatrix@-3pt{
\Lambda^{m+1,k} \ar[d]_{\alpha^{m,k}} \ar[r] &  X\ar[d]\\
\Delta^{m+1} \ar[r] & \star}
\end{equation}
où $m\geq 0$ et $0\leq k \leq m+1$, il existe un morphisme diagonal:
\begin{equation}\label{complexeKKAN2}
\xymatrix@-3pt{
\Lambda^{m+1,k} \ar[d]_{\alpha^{m,k}} \ar[r] &  X\ar[d]\\
\Delta^{m+1}\ar@{-->}[ru] \ar[r] & \star}
\end{equation}
complétant \eqref{complexeKKAN} dans un diagramme commutatif.

Plus généralement, si $0\leq n\leq \infty$ un ensemble simplicial $X$ est appelé un \emph{$n$-groupoïde de Kan} (voir aussi le paragraphe \ref{groupoideDG}) si $X$ est un complexe de Kan qui satisfait la condition d'extension de Kan de fa\c con stricte en dimension $m\geq n$ \emph{i.e.} si le morphisme diagonal du diagramme \eqref{complexeKKAN2} est unique pour $m\geq n$. En particulier un $\infty$-groupoïde de Kan est simplement un complexe de Kan. 




\renewcommand{\thesubsection}{\S\thesection.\arabic{subsection}}
\subsection{}\;\label{nequival}
\renewcommand{\thesubsection}{\thesection.\arabic{subsection}} 

Soit $\mathbf{Top}$ la catégorie des espaces topologiques et fonctions continues. Rappelons que le foncteur:
\begin{equation}\label{stand}
\begin{split}
{\bf\Delta}\xymatrix@C+5pt{\ar[r]&} \;&  \mathbf{Top}\\
[n]\xymatrix@C-10pt{\;\;\;\ar@{|->}[r]&\;\;\;}&\,\Delta_{top}^n\,=\bigg\{ \big(t_{0},\dots,t_{n}\big)\in\mathbb{R}^{n+1} \;\bigg|\; \overset{n}{\underset
{i=0}{\sum}}\; t_{n}=1\,,\;\,t_{i}\geq 0 \bigg\},
\end{split}
\end{equation}
induit une adjonction:
\begin{equation}\label{singu}
\xymatrix@C+15pt{\mathbf{Top}\ar@/_12pt/[r]_{s(\;\cdot\;)}\ar@{}[r]|-{\perp}&\simp \ar@/_12pt/[l]_{|\;\cdot\;|} \,,}
\end{equation}
où $s(\,\cdot\,)$ est défini par la formule $s(\mathcal{X})_n=\mathrm{Hom}_{\bf Top}\big(\Delta_{top}^n,\mathcal{X}\big)$, et le foncteur $|\,\cdot\,|$ est une extension de Kan à gauche du foncteur  \eqref{stand} le long du plongement de Yoneda:
\begin{equation*}
\xymatrix@C+5pt@R=1pt{{\bf\Delta} \ar@{^(->}[r]&  \simp \\
[n] \ar@{}[r]|-{\longmapsto} & \Delta^n \,,}
\end{equation*}
laquelle on va choisir et fixer.

On appelle $s(\mathcal{X})$ l'\emph{ensemble simplicial singulier} de l'espace topologique $\mathcal{X}$, et $|X|$ la \emph{réalisation géométrique} de l'ensemble simplicial $X$.

Fixons aussi une unité $\eta$ de l'adjonction $|\,\cdot \,| \dashv s(\,\cdot\,)$ et remarquons que cela nous permet d'associer à chaque $0$-simplexes d'un ensemble simplicial $X$, un point de l'espace topologique $|X|$ lequel on va note par le même symbole:
$$
\xymatrix@R=1pt@C-3pt{ X_0 \ar[rrr]^-{(\eta_X)_0}& & & s\big(|\,X\,|\big)_0\\ a &\ar@{|->}[r] & &a  }= \mathrm{Hom}_{\bf Top}\big(\Delta_{top}^0,|\,X\,|\big) \,.
$$

Étant donné $a\in X_{0}$ et $m\geq 0$, on pose $\pi_{m}(X,a)$ pour noter l'ensemble des classes à homotopie pointée près, des fonctions continues pointées:
$$
\xymatrix{\big(\mathbb{S}^m_{top},\star\big)\ar[r]&\big(|X|,a\big)}\,,
$$
où $\mathbb{S}^m_{top}$ est la $m$-sphère topologique, \emph{i.e.} $\mathbb{S}^m_{top}$ est l'ensemble des $x\in\mathbb{R}^{m+1}$ tels que $||x||=1$. 

Pour $m=0$, on voit qu'en fait $\pi_{0}(X,a)$ est indépendant de $a$; il est noté simplement $\pi_{0}(X)$ et appelé l'ensemble de \emph{composantes connexes par arcs de $X$}. Pour $m\geq 1$,  l'ensemble $\pi_{m}(X,a)$ admet une structure de groupe, commutatif si $m\geq 2$. Il est appelé le \emph{$m$-ième groupe d'homotopie de $X$ basé en $a$}.

Un morphisme d'ensembles simpliciaux $\xymatrix@C-7pt{X\ar[r]^F&Y}$ est dit une \emph{$n$-équivalence faible}, s'il satisfait que pour tout $0$-simplexe $a$ de $X$, les fonctions:
$$
\xymatrix@C+16pt{\pi_{m}(X,a)\ar[r]^{\pi_{m}(F,a)}&\pi_{m}(Y,Fa)}, 
$$
induisent en composant avec $|F|\colon\xymatrix@C-5pt{|X|\ar[r]&|Y|}$, soient bijectives pour tout $0\leq m\leq n$. 

\begin{lemme}
La notion de $n$-équivalence faible est indépendante du foncteur réalisation géométrique qu'on a choisi.
\end{lemme} 
\begin{proof}
Si $|\,\cdot \,|'$ est une autre extension de Kan à gauche du foncteur \eqref{stand} (le long du plongement de Yoneda) et $\eta'$ une unité de l'adjonction $|\,\cdot\,|' \dashv s(\,\cdot\,)$, alors il existe un isomorphisme naturelle de foncteurs $\xymatrix@C-3pt{|\,\cdot\,| \ar@{=>}[r]^-{\alpha} & |\,\cdot \,|' }$ tel que le triangle qui suit commute:
$$
\vcenter{\xymatrix@R-20pt@C+10pt{  & s\big(|\,\cdot\,|\big)\ar@{=>}[dd]^{s(\alpha)} \\
\mathrm{id} \ar@/^8pt/@{=>}[ru]^-{\eta} \ar@/_8pt/@{=>}[rd]_-{{\eta'}} & \\
& s\big(|\,\cdot\,|'\big)  }}\,.
$$

Il se suit que pour tout ensemble simplicial $X$ on a un homéomorphisme d'espaces topologiques $\xymatrix@C-8pt{|\,X\,| \ar[r]^-{\alpha_X} & |\,X\,|' }$ faisant commutatif le diagramme suivant:
$$
\vcenter{\xymatrix@R-20pt@C+18pt{  
& \mathrm{Hom}_{{\bf Top}}\Big( \Delta_{top}^0,|\,X\,|\Big) \ar[dd]^{\alpha_X\circ -}\,\,  \\
X_0 \ar@/^8pt/[ru]^-{(\eta_X)_0} \ar@/_8pt/[rd]_-{{(\eta'_X)_0}} & \\
& \mathrm{Hom}_{{\bf Top}}\Big( \Delta_{top}^0,|\,X\,|'\Big)\,;  }}
$$
donc $\alpha_X$ induit un isomorphisme naturelle entre les groupes $\pi_{m}(X,a)$ définis à partir du foncteur $|\,\cdot \,|$ et la transformation naturelle $\eta$ ou de $|\,\cdot \,|'$ et $\eta'$. 
\end{proof}

Rappelons le résultat suivant:

\begin{theoreme}\label{ntypess}
Si $n\geq 0$, la catégorie des ensembles simpliciaux $\simp$ admet une structure de catégorie de modèles lorsque:
\begin{align*}
\big\{\,\text{équivalences faibles}\,\big\} \quad &=\quad \big\{\,\text{$n$-équivalences faibles}\,\big\}\;\,= \;\, {\bf W}_{n}\,,\\
\big\{\,\text{cofibrations}\,\big\} \quad &=\quad \big\{\,\text{monomorphismes}\,\big\}\;\,= \;\, {\bf mono}\,,\\
\text{et}\qquad\big\{\,\text{fibrations}\,\big\}\quad &=\quad \big\{\,\text{morphismes avec la propriété de relèvement }\\
 &\phantom{=}\;\;\quad \text{\phantom{$\{\,$}à droite par rapport à ${\bf mono}\,\cap\,{\bf W}_{n}$}\,\big\}\;\,= \;\, {\bf fib}_{n}\,.
 \end{align*}
 
En plus, un ensemble simplicial $X$ est un objet fibrant de $(\simp,{\bf W}_n, {\bf mono},{\bf fib}_n)$ si et seulement si $X$ est un complexe de Kan tel que $\pi_m(X,a)=0$ pour tout $m\geq n+1$ et tout $0$-simplexe $a$ de $X$.
\end{theoreme}
\begin{proof}
Dans le Chapitre 9 de \cite{ci} on montre que $(\simp,{\bf W}_n, {\bf mono},{\bf fib}_n)$ est une catégorie de modèles à engendrement cofibrant pour tout $0\leq n\leq \infty$. 

Si $0\leq n< \infty$, à partir de résultats plus connus de la catégorie de modèles $(\simp,{\bf W}_\infty, {\bf mono},{\bf fib}_\infty)$ et le Théorème 4.7 de \cite{barwick} on va montrer que $(\simp,{\bf W}_n, {\bf mono},{\bf fib}_n)$ est une localisation de Bousfield à gauche de $(\simp,{\bf W}_\infty, {\bf mono},{\bf fib}_\infty)$ par rapporte à l'ensemble de morphismes \eqref{elsnantes} ci-dessous. Plus explicitement:

On sait bien que $(\simp,{\bf W}_\infty, {\bf mono},{\bf fib}_\infty)$ est une catégorie de modèles dont les objets fibrants sont les complexes de Kan (voir par exemple le Théorème 3.6.5 de \cite{hovey} ou le Théorème 3 de II\S3 de \cite{quillenhh}). En fait $(\simp,{\bf W}_\infty, {\bf mono},{\bf fib}_\infty)$ est une catégorie de modèles propre à gauche et combinatoire\footnote{Une catégorie de modèles est dite \emph{combinatoire} si elle est localement présentable et à engendrement cofibrant. Voir la Définition 1.21 de \cite{barwick}}; donc d'après le Théorème 4.7 de \cite{barwick} on peut considérer la localisation de Bousfield à gauche de $(\simp,{\bf W}_\infty, {\bf mono},{\bf fib}_\infty)$ par rapport à l'ensemble des morphismes:
\begin{equation}\label{elsnantes}
S_{n}\,=\,\Big\{\xymatrix@C+6pt{\partial\Delta^{m+1}\ar@{^(->}[r]^-{\alpha^m}&\Delta^{m+1}}\Big\}_{m\geq n}\,.
\end{equation}

Dans le deux Lemmes suivants on va poser $[\;\cdot\;,\;\cdot\;]_\infty$ pour noter l'ensemble des morphismes dans la catégorie homotopique $\simp[{\bf W_{\infty}}^{-1}]$.

\begin{lemme}\label{lemmelocala}
Si $m\geq 0$ et $X$ est un ensemble simplicial alors $\pi_{m}(X,a)=0$ pour tout $0$-simplexe $a$ de $X$ si et seulement si la fonction:
$$
\xymatrix@C+10pt{[\Delta^{m+1},X]_\infty\ar[r]^-{(\alpha^m)^*}&[\partial\Delta^{m+1},X]_\infty}\,,
$$
induite par le morphisme $\xymatrix@C+6pt{\partial\Delta^{m+1}\ar@{^(->}[r]^-{\alpha^m}&\Delta^{m+1}}$ dans la catégorie homotopique $\simp[{\bf W_{\infty}}^{-1}]$ est bijective.  

En particulier, un ensemble simplicial $X$ vérifie que $\pi_{m}(X,a)=0$ pour tout $m\geq n$ et tout $0$-simplexe $a$ de $X$ si et seulement si, $X$ est un objet $S_{n}$-local de $(\simp,{\bf W}_\infty, {\bf mono},{\bf fib}_\infty)$.
\end{lemme}
\begin{proof}
Rappelons que si $A$ et $B$ sont deux ensembles simpliciaux, le foncteur réalisation géométrique $|\,\cdot\,|$ induit une bijection entre l'ensemble de morphismes $[A,B]_{\infty}$ de la catégorie homotopique $\simp[{\bf W}_{\infty}^{-1}]$, et l'ensemble $\big[|A|,|B|\big]_{{\bf Top}}$ des morphismes dans la catégorie de fractions  ${\bf Top}\big[({\bf W}_{\infty}^{top})^{-1}\big]$ (voir les Théorèmes 2.4.19, 2.4.23 et 3.6.7 de \cite{hovey}). D'un autre vu que $|A|$ est un complexe cellulaire, l'ensemble $\big[|A|,|B|\big]_{{\bf Top}}$ est simplement l'ensemble des classes d'homotopie habituelles des fonctions continues de $|A|$ vers $|B|$ (voir le Théorème 2.4.19 de \cite{hovey}).

En particulier, si $\mathbb{D}^{m+1}$ note l'ensemble des $x\in\mathbb{R}^{m+1}$ tels que $||x||\leq 1$ et on considère $\mathbb{S}_{top}^{m}$ comme étant le sous-ensemble des $x\in\mathbb{D}^{m+1}$ tels que $||x||=1$; la fonction:  
\begin{equation}\label{alphapra}
\xymatrix{[\Delta^{m+1},X]_{\infty}\ar[r]^-{{(\alpha^m)}^*}&[\partial\Delta^{m+1},X]_{\infty}}\,,
\end{equation}
s'identifie à la fonction:
\begin{equation}\label{alphapra2}
\xymatrix{[\mathbb{D}^{m+1},|X|]_{\bf Top}\ar[r]^-{}&[\mathbb{S}_{top}^{m},|X|]_{\bf Top}}\,,
\end{equation}
induite par l'inclusion $\xymatrix@C-8pt{\mathbb{S}^m_{top}\ar@{^(->}[r]&\mathbb{D}^{m+1}}$ où $\big[\;,\;\big]_{{\bf Top}}$ note l'ensemble des classes d'homotopie habituelles des fonctions continues.

D'un autre côté, on  se convaincre que \eqref{alphapra2} est toujours une fonction injective dont l'image sont les classes d'homotopie de fonctions $\xymatrix@C-8pt{\mathbb{S}_{top}^{m}\ar@{^(->}[r]&|X|}$ admettant une extension le long de l'inclusion $\xymatrix@C-8pt{\mathbb{S}^m_{top}\ar@{^(->}[r]&\mathbb{D}^{m+1}}$. En plus, on sait que si $f$ est une fonction continue de $\mathbb{S}_{top}^{m}$ vers un espace topologique quelconque $A$, pour $m\geq 1$ les énoncés suivants sont équivalents:
\begin{enumerate}
\item $f$ représente l'élément nul du groupe $\pi_{m}(A,a)$ où $a=f(1,0,\dots,0)\in A$.
\item $f$ admet une extension le long de l'inclusion $\xymatrix@C-8pt{\mathbb{S}^m_{top}\ar@{^(->}[r]&\mathbb{D}^{m+1}}$. 
\end{enumerate}

Tandis que pour $m=0$ on a les énoncés équivalents:
\begin{enumerate}
\item L'image de $f$ est contenue au même composante connexe par arcs de $A$. ($\mathbb{S}^0_{top}$ est constitué de deux points)
\item $f$ admet une extension le long de l'inclusion $\xymatrix@C-8pt{\mathbb{S}^0_{top}\ar@{^(->}[r]&\mathbb{D}^{1}}$. 
\end{enumerate}
\end{proof}

Montrons:

\begin{lemme}\label{lemmelocal2}
Si $f\colon\xymatrix@C-8pt{X\ar[r]&Y}$ est un morphisme d'ensembles simpliciaux et $n\geq 0$, alors $f$ est une $n$-équivalence faible d'ensembles simpliciaux si et seulement si, $f$ est une équivalence faible $S_{n}$-locale de $(\simp,{\bf W}_\infty, {\bf mono},{\bf fib}_\infty)$, \emph{i.e.} $f$ vérifie que pour tout ensemble simplicial $S_n$-local $Z$ (voir le Lemme \ref{lemmelocala} ci-dessus) la fonction:
$$
\xymatrix@C+15pt{[Y,Z]_{\infty}\ar[r]^-{f^*} & [X,Z]_{\infty}\,,}
$$
dans la catégorie homotopique $\simp[{\bf W}_\infty^{-1}]$ soit bijective.
\end{lemme}
\begin{proof}
Dans la preuve de cet énoncé on aura besoin des propriétés du cosquelette des ensembles simpliciaux qu'on verra dans \S\ref{tron} ci-dessous.

Rappelons pour commencer que la transformation naturelle de foncteurs:
$$
\xymatrix{ 
\mathrm{Hom}_{\simp}( \;\cdot\;,\;\cdot\;) \ar[r] & 
[\;\cdot\;,\;\cdot\;]_\infty}
$$
définie par le foncteur canonique $\xymatrix@C-8pt{\simp\ar[r]&\simp[{\bf W_{\infty}}^{-1}]}$, induit une bijection d'ensembles: 
$$
\pi_0\big(\underline{\mathrm{Hom}}_{\simp}(A,B)\big) \, \cong \, [A,B]_\infty
$$ 
pour tout ensemble simplicial $A$ et tout complexe de Kan $B$. En effet si $A$ est un ensemble simplicial quelconque alors:
$$
\vcenter{\xymatrix@-8pt{\underset{\phantom{.}}{A \sqcup A}\ar@<-4pt>@{>->}@/_5pt/[rd]_-{i_{0}+i_{1}}\ar[rr]^-{\mathrm{id}+\mathrm{id}} && A \\ 
&A\times \Delta^1 \ar@<-4pt>@/_5pt/[ru]_-{\mathrm{proj}}\ar@{}[ru]|-{\text{\rotatebox[origin=c]{225}{\Large $\widetilde{\phantom{w}}$}}}&}}
$$
est un objet cylindre de la catégorie de modèles $(\simp,{\bf W}_{\infty}, {\bf mono},{\bf fib}_{\infty})$, où les monomorphismes $i_0$ et $i_1$ sont par définition les morphismes composés:
$$
\vcenter{\xymatrix@C+8pt{A\,\cong\,A\times \Delta^0 \ar[r]^-{A\times \delta_0} & A\times \Delta^1}}
\qquad\text{et}\qquad
\vcenter{\xymatrix@C+8pt{A\,\cong\,A\times \Delta^1 \ar[r]^-{A\times \delta_1} & A\times \Delta^1}}
$$
respectivement.

Soit maintenant $f\colon\xymatrix@C-8pt{X\ar[r]&Y}$ un morphisme d'ensembles simpliciaux et $Z$ un ensemble simplicial quelconque. Remarquons que si on considère de remplacements fibrants de $f$ et $Z$ dans la catégorie de modèles $(\simp,{\bf W}_{\infty}, {\bf mono},{\bf fib}_{\infty})$, \emph{i.e.} un carré commutatif et une flèche:
$$
\vcenter{\xymatrix@R-6pt{X\ar[d]\ar[r]^f & Y\ar[d]\\ X' \ar[r]_{f'} & Y'}}\qquad\;\text{et}\qquad\;\vcenter{\xymatrix@R-6pt{Z\ar[d]\\Z'}}\,,
$$
où les morphismes verticaux sont des $\infty$-équivalences faibles de but un complexe de Kan; alors la fonction entre les ensembles de morphismes de la catégorie $\simp[{\bf W_{\infty}}^{-1}]$:
\begin{equation}\label{infty}
\xymatrix@C+15pt{
[Y,Z]_{\infty}\ar[r]^-{f^*} & [X,Z]_{\infty}\,,}
\end{equation}
s'identifie avec la fonction: 
\begin{equation}\label{infty1a}
\xymatrix@C+15pt{
\pi_{0}\Big(\underline{\mathrm{Hom}}_{\simp}\big(Y',Z'\big)\Big)\ar[r]^{(f')^*}& \pi_{0}\Big(\underline{\mathrm{Hom}}_{\simp}\big(X',Z'\big)\Big)\,.}
\end{equation}

Si on suppose que $Z$ soit un ensemble simplicial $S_{n}$-local, il se suit du Lemme \ref{lemmelocala} et du Corollaire \ref{lecoro} de \S\ref{tron} que le morphisme $\xymatrix@C-8pt{Z'\ar[r] & {\bf csq}_{n+1}(Z')}$ induit d'une unité quelconque de l'adjonction $\tau_{n\;*}  \dashv \tau_{n}^{\phantom{a} *}$ est une $\infty$-équivalence faible entre complexes de Kan. Donc, \eqref{infty1a} s'identifie avec la fonction:
$$
\xymatrix@C+15pt{
\pi_{0}\Big(\underline{\mathrm{Hom}}_{\simp}\big(Y',{\bf csq}_{n+1}(Z')\big)\Big)\ar[r]^-{(f')^*}& \pi_{0}\Big(\underline{\mathrm{Hom}}_{\simp}\big(X',{\bf csq}_{n+1}(Z')\big)\Big)\,,}
$$
qui d'après le Lemme \ref{lemehomo} s'identifie avec la fonction: 
\begin{equation}\label{infty2a}
\xymatrix@C+15pt{
\pi_{0}\Big(\underline{\mathrm{Hom}}_{\simp}\big({\bf csq}_{n+1}(Y'),{\bf csq}_{n+1}(Z')\big)\Big)\ar[r]^{{\bf csq}_{n+1}(f')^*}& \pi_{0}\Big(\underline{\mathrm{Hom}}_{\simp}\big({\bf csq}_{n+1}(X'),{\bf csq}_{n+1}(Z')\big)\Big)\,.}
\end{equation}

Si on suppose que $f$ soit une $n$-équivalence faible, il se suit du Corollaire \ref{lecoro} que ${\bf csq}_{n+1}(f')$ est une $\infty$-équivalence faible, \emph{i.e.} \eqref{infty2a} est une fonction bijective. Donc, si $f$ est une $n$-équivalence faible la fonction \eqref{infty} est bijective pour tout ensemble simplicial $S_{n}$-local $Z$, c'est-à-dire, $f$ est une équivalence faible $S_{n}$-locale.

Soit maintenant $f\colon\xymatrix@C-8pt{X\ar[r]&Y}$ une équivalence faible $S_{n}$-locale, et considérons un remplacement fibrant de $f$ dans la catégorie de modèles $(\simp,{\bf W}_{\infty}, {\bf mono},{\bf fib}_{\infty})$, \emph{i.e.} un carré commutatif:
$$
\vcenter{\xymatrix@R-6pt{X\ar[d]\ar[r]^f & Y\ar[d]\\ X' \ar[r]_{f'} & Y'\,,}}
$$
où les morphismes verticaux sont des $\infty$-équivalences faibles de but un complexe de Kan.

Remarquons par ailleurs que d'après le Corollaire \ref{lecoro} et le Lemme \ref{lemmelocala} les ensembles simpliciaux ${\bf csq}_{n+1}X'$ et ${\bf csq}_{n+1}Y'$ sont $S_{n}$-locaux; en particulier, la fonction:
\begin{equation}\label{deuxeqa}
\xymatrix@R=5pt@C+15pt{[Y,{\bf csq}_{n+1}X']_{\infty}\ar@{}[d]|-{\mathrel{\reflectbox{\rotatebox[origin=c]{90}{$\cong$}}}}
\ar[r]^-{f^*} & [X,{\bf csq}_{n+1}X']_{\infty}\ar@{}[d]|-{\mathrel{\reflectbox{\rotatebox[origin=c]{90}{$\cong$}}}}\\
\pi_{0}\Big(\underline{\mathrm{Hom}}_{\simp}\big(Y',{\bf csq}_{n+1}X'\big)\Big)
\ar[r]_-{(f')^*}& \pi_{0}\Big(\underline{\mathrm{Hom}}_{\simp}\big(X',{\bf csq}_{n+1}X'\big)\Big)\,,}
\end{equation}
est bijective.

D'un autre côté, grâce au Lemme \ref{lemehomo} la fonction bijective \eqref{deuxeqa} s'identifie avec la fonction:
\begin{equation}\label{deuxeq2}
\xymatrix@C+15pt@R=5pt{
\pi_{0}\Big(\underline{\mathrm{Hom}}_{\simp}\big({\bf csq}_{n+1}Y',{\bf csq}_{n+1}X'\big)\Big)
\ar@{}[d]|-{\mathrel{\reflectbox{\rotatebox[origin=c]{90}{$\cong$}}}}
\ar[r]^{{\bf csq}_{n+1}(f')^*}& \pi_{0}\Big(\underline{\mathrm{Hom}}_{\simp}\big({\bf csq}_{n+1}X',{\bf csq}_{n+1}X'\big)\Big)
\ar@{}[d]|-{\mathrel{\reflectbox{\rotatebox[origin=c]{90}{$\cong$}}}}\\
[{\bf csq}_{n+1}Y',{\bf csq}_{n+1}X'\big]_{\infty}\ar[r]_-{{\bf csq}_{n+1}(f')^*}& \big[{\bf csq}_{n+1}X',{\bf csq}_{n+1}X'\big]_{\infty}\,.}
\end{equation}

Soit $g:\xymatrix@C-8pt{{\bf csq}_{n+1}Y'\ar[r]&{\bf csq}_{n+1}X'}$ un morphisme dans la catégorie de fractions $\simp[{\bf W}_\infty^{-1}]$, tel que le composé $g\circ {\bf csq}_{n+1}(f')$ soit le morphisme identité de l'objet ${\bf csq}_{n+1}X'$ dans cette catégorie homotopique.

On vérifie alors que la fonction composée: 
$$
\xymatrix@C+6pt{
 [{\bf csq}_{n+1}X',{\bf csq}_{n+1}Y']_{\infty} \ar[r]^-{g^*}&
 [{\bf csq}_{n+1}Y',{\bf csq}_{n+1}Y']_{\infty}\ar[r]^-{{\bf csq}_{n+1}(f')^*} & 
 [{\bf csq}_{n+1}X',{\bf csq}_{n+1}Y']_{\infty}\,,}
$$
est égale à la fonction identité; donc:
$$
\xymatrix@C+6pt{
 [{\bf csq}_{n+1}Y',{\bf csq}_{n+1}Y']_{\infty} \ar[r]^-{{\bf csq}_{n+1}(f')^*}&
 [{\bf csq}_{n+1}X',{\bf csq}_{n+1}Y']_{\infty}\ar[r]^-{g^*} & 
 [{\bf csq}_{n+1}Y',{\bf csq}_{n+1}Y']_{\infty}\,,}
$$
est aussi égale à la fonction identité, parce que la flèche:
$$
\xymatrix@C+16pt{
[{\bf csq}_{n+1}Y',{\bf csq}_{n+1}Y']_{\infty} \ar[r]^-{{\bf csq}_{n+1}(f')^*} & [{\bf csq}_{n+1}X',{\bf csq}_{n+1}Y']_{\infty}\,,}
$$
est bijective.

Donc $g\circ{\bf csq}_{n+1}(f')$ et ${\bf csq}_{n+1}(f')\circ g$ sont des morphismes identité dans la catégorie $\simp[{\bf W}_\infty^{-1}]$. Autrement dit, ${\bf csq}_{n+1}(f')$ est une $\infty$-équivalence faible. Donc, $f$ est une $n$-équivalence faible d'après  le Corollaire \ref{lecoro}.
\end{proof}

D'après le Théorème 4.7 de \cite{barwick} et les Lemmes \ref{lemmelocala} et \ref{lemmelocal2} qu'on vient de montrer, on a que $(\simp,{\bf W}_n, {\bf mono},{\bf fib}_n)$ est une catégorie de modèles propre à gauche et combinatoire dont les objets fibrants sont les complexes de Kan $X$ tels que $\pi_m(X,a)=0$ pour tout $m\geq n+1$ et tout $0$-simplexe $a$ de $X$.
\end{proof}

Montrons aussi:

\begin{lemme}\label{ntypessp} 
La catégorie cartésienne fermée des ensembles simpliciaux $\simp$ munie de la structure de catégorie de modèles $(\simp,{\bf W}_n, {\bf mono},{\bf fib}_n)$ du Théorème \ref{ntypess} est une catégorie de modèles monoïdale symétrique dans le sens de \cite{hovey} $($voir Définition 4.2.6 de \cite{hovey} ou Définition 1.27 de \cite{barwick}$)$. 
\end{lemme}
\begin{proof}
D'après les Définitions  4.2.1 et 4.2.6 de \cite{hovey} il faut démontrer que si $\xymatrix@C-8pt{X\ar[r]^j&Y}$ et $\xymatrix@C-8pt{Z\ar[r]^q&W}$ sont de monomorphismes des ensembles simpliciaux, alors le morphisme $\varphi$ dans le diagramme somme amalgamée suivant:
\begin{equation}\label{mmcarre} 
\xymatrix@R-10pt@C-10pt{
Z\otimes X \ar[dd]_{Z\otimes \,j}\ar[rr]^{q\,\otimes X} && 
W\otimes X \ar[dd]|-{\widetilde{Z\otimes \,j}}\ar@/^15pt/[rddd]^{W\otimes \,j}&\\\\
Z\otimes Y\ar[rr]|-{\widetilde{q\,\otimes X}}\ar@/_15pt/[rrrd]_{q\,\otimes Y}&&
Z\otimes Y\underset{Z\otimes X}{\bigsqcup}W\otimes X\ar@{-->}[rd]|(.65)\varphi&\\
&&&W\otimes Y\,,}
\end{equation}
est un monomorphisme de $\simp$, laquelle est aussi une $n$-équivalence faible lorsque un de deux morphismes $j$ ou $q$ le soit.

Vu que dans ce cas le foncteur $\otimes$ est le produit cartésien argument par argument, il se suit que pour $\xymatrix@C-8pt{X\ar[r]^j&Y}$ et $\xymatrix@C-8pt{Z\ar[r]^q&W}$ de monomorphismes d'ensembles simpliciaux, les morphismes $Z\otimes j$, $W\otimes j$, $q\otimes X$ et $q\otimes Y$ dans le diagramme \eqref{mmcarre} sont aussi de monomorphismes. En particulier,  les morphismes $ \widetilde{Z\otimes \,j}$ et $\widetilde{q\,\otimes X}$ sont de monomorphismes; donc, $\varphi$ est aussi un monomorphisme (Il suffit de vérifier ces assertions pour la catégorie des ensembles).

Enfin, si on suppose que le morphisme $j$ (resp. $q$) soit en plus une $n$-équivalence faible; vu que les foncteurs $\pi_i$ commutent aux produit finis, les morphismes ${Z\otimes \,j}$ et ${W\otimes \,j}$ (resp. ${q\,\otimes X}$ et ${q\,\otimes Y}$) sont des monomorphismes et des $n$-équivalences faibles aussi. En particulier, $\widetilde{Z\otimes \,j}$ (resp. $\widetilde{q\,\otimes X}$) est un monomorphisme et une $n$-équivalence faible, car la famille de morphismes ${\bf mono}\cap{\bf W}_n$ est stable par cochangement de base (vrais pour toute catégorie de modèles). Donc, $\varphi$ est une $n$-équivalence faible parce que ${\bf W}_n$ satisfait à la propriété de deux-sur-trois.
\end{proof}

On pose $\mathrm{Ho}_{n}(\simp)$ pour noter la \emph{catégorie des $n$-types d'homotopie} (\emph{topologiques}) \emph{i.e.} la catégorie homotopique $\simp[{\bf W}_n^{-1}]$ de la catégorie de modèles $(\simp,{\bf W}_{n}, {\bf mono},{\bf fib}_{n})$ du Théorème \ref{ntypess}, et on désigne par $[\,\cdot\,,\,\cdot\,]_{n}$ l'ensemble des morphismes dans $\mathrm{Ho}_{n}(\simp)$. Si $X$ est un ensemble simplicial, la classe d'équivalence associée à $X$ dans l'ensemble (ou la classe) des objets à isomorphisme près de la catégorie $\mathrm{Ho}_{n}\big(\,\simp\,\big)$, est appelée le \emph{$n$-type d'homotopie de $X$}.

D'après le Lemme \ref{ntypessp} la catégorie des $n$-types d'homotopie $\mathrm{Ho}_{n}(\simp)$ est une catégorie carté\-sienne fermée (voir le Lemme \ref{prosomm} ci-dessus et \S 4.3 de \cite{hovey}), dont l'objet des morphismes (le "$\mathrm{\bf hom}$" interne) est le complexe de fonction dérivé ${\bf R}\mathrm{Hom}_{\simp}$ (défini par la localisation simpliciale de \cite{DK1,DK2,DK}, ou construit à partir de résolutions fibrants et cofibrants comme dans \S5.4 de \cite{hovey} ou le Chapitre 18 de \cite{hirschhorn}). 

La catégorie des $n$-types d'homotopie munie d'une telle structure de catégorie cartésienne fermée est dite quelquefois la \emph{catégorie homotopique des $n$-groupoïdes}. 

Vu qu'on est dans le cadre d'une catégorie de modèles dont tous les objets sont cofibrants $(\simp,{\bf W}_n, {\bf mono},{\bf fib}_n)$, il est possible de construire la catégorie homotopique des $n$-groupoïdes plus directement: Appelons un ensemble simplicial $X$ \emph{$n$-fibrant} si $X$ est un objet fibrant de la catégorie de modèles $(\simp,{\bf W}_{n}, {\bf mono},{\bf fib}_{n})$ du Théorème \ref{ntypess}, c'est-à-dire si $X$ est un complexe de Kan tel que $\pi_m(X,a)=0$ pour tout $m\geq n+1$ et tout $0$-simplexe $a$ de $X$. On pose ${\bf Fib}^n$ pour noter la sous-catégorie pleine de $\simp$ dont les objets sont les ensembles simpliciaux $n$-fibrants.

Remarquons:

\begin{lemme}\label{nfibcartesiene}
Soit $0\leq n\leq \infty$. Si $X$ et $Y$ sont des ensembles simpliciaux $n$-fibrants le produit cartésienne argument par argument $X\times Y$ et l'ensemble des morphismes $\underline{\mathrm{Hom}}_{\simp}(X,Y)$ sont des ensembles simpliciaux $n$-fibrants. Autrement dit, la catégorie ${\bf Fib}^n$ est une sous-catégorie cartésienne fermée de $\simp$.
\end{lemme}
\begin{proof}
Le produit de deux ensembles simpliciaux $n$-fibrants est un ensemble simplicial $n$-fibrant parce que la famille des objets fibrants d'une catégorie de modèles est toujours stable par produits. 

D'un autre si $X$ et $Y$ sont des objets fibrants de la catégorie de modèles $(\simp,{\bf W}_{n}, {\bf mono},{\bf fib}_{n})$, on déduit du Lemme \ref{ntypessp} (voir par exemple le Lemme 4.2.2 de \cite{hovey}) que l'ensemble simplicial $\underline{\mathrm{Hom}}_{\simp}\big(X,Y\big)$ est aussi un objet fibrant, c'est-à-dire un ensemble simplicial $n$-fibrant, vu que tous les objets de $(\simp,{\bf W}_{n}, {\bf mono},{\bf fib}_{n})$ sont des objets cofibrants.
\end{proof}

On définit la catégorie homotopique $h{\bf Fib}^n$ de la catégorie ${\bf Fib}^n$ des ensembles simpliciaux $n$-fibrants comme étant la catégorie dont les objets sont les ensembles simpliciaux $n$-fibrants et l'ensemble de morphismes est l'ensemble des composantes connexes par arcs $\pi_0\big(\underline{\mathrm{Hom}}_{\simp}\big)$. Il se suit aussi-tôt que le foncteur d'inclusion de ${\bf Fib}^n$ vers $\simp$ induit une équivalence de catégories entre $h{\bf Fib}^n$ et la catégorie des $n$-types d'homotopie $\mathrm{Ho}_{n}(\simp)$. 

Enfin remarquons que dans la catégorie cartésienne fermée $h{\bf Fib}^n$ le foncteur $\underline{\mathrm{Hom}}_{\simp}$ est l'objet des morphismes (le "$\mathrm{\bf hom}$" interne), tandis que dans la catégorie cartésienne fermée des $n$-types d'homotopie $\mathrm{Ho}_{n}(\simp)$ l'objet des morphismes est le complexe de fonction dérivé ${\bf R}\mathrm{Hom}_{\simp}$.

\renewcommand{\thesubsection}{\S\thesection.\arabic{subsection}}
\subsection{}\;\label{pointesec}
\renewcommand{\thesubsection}{\thesection.\arabic{subsection}} 

On désigne par $\simp_{\star}$ la catégorie des \emph{ensembles simpliciaux pointés}, \emph{i.e.} la catégorie des couples $X=(X,x)$, où $X$ est un ensemble simplicial et $x:\xymatrix@C-8pt{\star\ar[r]&X}$ est un morphisme de l'ensemble simplicial final $\star$ vers $X$. Un morphisme de $(X,x)$ vers $(Y,y)$ est une flèche $f:\xymatrix@-6pt{X\ar[r]&Y}$ de $\simp$ telle que $f x=y$. 

De fa\c con équivalente, si ${\bf Ens}_{\star}$ note la catégorie des \emph{ensembles pointés}, \emph{i.e.} la catégorie des couples $X=(X,x)$, où $X$ est un ensemble et $x:\xymatrix@C-8pt{\star\ar[r]&X}$ est une fonction de l'ensemble ponctuel $\star$ vers $X$; la catégorie $\simp_{\star}$ est isomorphe à la catégorie des préfaisceaux de la catégorie des simplexes $\Delta$ à valeurs dans la catégorie ${\bf Ens}_{\star}$.

En tout cas, on vérifie que le foncteur canonique $\pi :\xymatrix@-8pt{\simp_{\star}\ar[r]&\simp}$ défini par la règle $(X,x)\mapsto X$, admet un adjoint à  gauche $(\;\cdot\;)_{+}: \xymatrix@-5pt{\simp\ar[r]&\simp_{\star}}$, défini dans un objet $A$ de $\simp$ comme le coproduit $A_{+}=A\sqcup\star$, pointé par le morphisme canonique $\xymatrix@C-8pt{\star\ar[r]&A\sqcup\star}$. 

Rappelons que la catégorie $\simp_{\star}$ admet une structure de catégorie monoïdale symétrique fermée, lorsque le tenseur est le \emph{produit wedge}:
$$
\xymatrix@C+20pt{\simp_{\star} \times \simp_{\star} \ar[r]^-{\cdot\,\wedge\,\cdot} & \simp_{\star}}, 
$$
défini dans deux ensembles simpliciaux pointés $X=(X,x)$ et $Y=(Y,y)$, par un carré cocartésien dans la catégorie $\simp$: 
\begin{equation}\label{wedgedefi2}
\xymatrix{
\star \ar[r]& X\wedge Y\\ 
\big(X\times \star\big) \sqcup \big(\star\times \,Y\big) \ar[r] \ar[u]  & X \times Y\ar[u]\,.}
\end{equation}

L'unité du produit wedge est l'ensemble simplicial constant à valeurs $\star_{+}=\star\sqcup \star$, l'isomorphisme de symétrie est induit de l'isomorphisme canonique $X\times Y\cong Y\times X$ d'ensembles simpliciaux, et le cotenseur:
$$
\xymatrix@C+35pt{\simp_{\star}^{op} \times \simp_{\star} \ar[r]^-{\mathrm{hom}_{\simp_{\star}}^{\wedge}(\,\cdot\,,\,\cdot\,)} & \simp_{\star}},
$$
est le foncteur défini dans deux ensembles simpliciaux pointés $X=(X,x)$ et $Y=(Y,y)$ par la formule:
\begin{equation}\label{hompointe2}
\mathrm{hom}_{\simp_{\star}}^{\wedge}(X,Y)_{n} \, = \, \Big(\mathrm{Hom}_{\simp_{\star}}\big(X\wedge\Delta^n_{+},Y\big),y\Big)\qquad \text{si}\quad n\geq 0.
\end{equation}

\begin{corollaire}\label{pointeco2}
Soit $0\leq n\leq \infty$. La catégorie des ensembles simpliciaux pointés $\simp_{\star}$ munie de la structure de catégorie monoïdale symétrique donnée par le produit wedge, admet une structure de catégorie de modèles monoïdale symétrique lorsque:
\begin{align*}
\big\{\,\text{équivalences faibles}\,\big\} \quad &= \quad 
\big\{\, f:\xymatrix@-14pt{(X,x)\ar[r]&(Y,y)}\,\big|\; f:\xymatrix@-14pt{X\ar[r]&Y}\in {\bf W}_{n}\,\big\}\;\,= \;\,\pi^{-1}{\bf W}_{n}\,, \\
\big\{\,\text{cofibrations}\,\big\} \quad &=\quad \;\big\{\,\text{monomorphismes}\,\big\}\;\;\,= \;\,{\bf mono}\\
\text{et}\qquad\big\{\,\text{fibrations}\,\big\} \quad &=\quad  
\big\{\, f:\xymatrix@-14pt{(X,x)\ar[r]&(Y,y)}\,\big|\; f:\xymatrix@-14pt{X\ar[r]&Y}\in {\bf fib}_{n}\,\big\}\;\,= \;\,\pi^{-1}{\bf fib}_{n}\,.
\end{align*}

En particulier, 
\begin{equation}\label{olvidar2}
\simp_{\star}
\vcenter{\xymatrix@C+15pt{
\phantom{a}\ar@{}[r]|{\perp}\ar@<-4pt>@/_10pt/[r]_{\pi \,=\, \textit{Fonteur d'oubli}}&\phantom{a}
\ar@<-4pt>@/_10pt/[l]_-{(\,\cdot\,)_{+}}}}
\simp
\end{equation}
est une adjonction de Quillen.  
\end{corollaire}
\begin{proof}
D'après la Proposition 4.2.9 de \cite{hovey}, ceci est une conséquence du Théorème \ref{ntypess} et le Lemme \ref{ntypessp} ci-dessus.
\end{proof}

On désigne par $\mathrm{Ho}_{n}\big(\,\simp_{\star}\,\big)$ la \emph{catégorie des $n$-types d'homotopie pointés}, \emph{i.e.} la catégorie homotopique de la catégorie de modèles $(\simp_{\star},\pi^{-1}{\bf W}_{n}, {\bf mono},\pi^{-1}{\bf fib}_{n})$, et on note $[\,\cdot\,,\,\cdot\,]_{n}^{\text{\tiny{$\text{pointé}$}}}$ l'ensemble des morphismes dans $\mathrm{Ho}_{n}(\simp_{\star})$. Si $X$ est un ensemble simplicial pointé, la classe associée à $X$ dans l'ensemble des objets à isomorphisme près de la catégorie $\mathrm{Ho}_{n}\big(\,\simp_{\star}\,\big)$, est appelée le \emph{$n$-type d'homotopie pointé de $X$}. 

\renewcommand{\thesubsection}{\S\thesection.\arabic{subsection}}
\subsection{}\;\label{sphere}
\renewcommand{\thesubsection}{\thesection.\arabic{subsection}}\label{combina}

Si $m\geq 0$ posons $\mathbb{S}^{m}: = \Delta^{m}/  \partial\Delta^{m}$ pour noter la $m$-sphère simpliciale, \emph{i.e.} l'ensemble simplicial défini par un carré cocartésien de $\simp$:  
$$
\xymatrix{
\star \ar[r]^-{\star} &  \mathbb{S}^{m}\\ 
\partial\Delta^{m}\ar[r]_{\alpha^{m-1}} \ar[u]  & \Delta^{m}\ar[u]\,,}
$$
où $\star$ note l'ensemble simplicial constant à valeurs l'ensemble ponctuel. 

Étant donné un ensemble simplicial $X$ et $a\in X_{0}$; du Lemme de Yoneda et de la propriété universelle de tout carré cocartésien, on déduit une bijection entre l'ensemble $X_m$ des $m$-simplexes $x$ de $X$ avec la propriété:
$$
\alpha^{m-1}_{X}(x)=\big(d_{0}x,\dots,d_{m}x\big) = (a,\dots,a) \in \underbrace{X_{m-1}\times \dots \times X_{m-1}}_{m+1}\,,
$$
et l'ensemble des morphisme d'ensembles simpliciaux:
\begin{equation}\label{poin}
\xymatrix@C-3pt{\mathbb{S}^{m}\ar[r]^-{\sigma^{x}} & X} \qquad \text{vérifiant \, que} \qquad \sigma^x_{0}(\star)=a \in X_{0}\,.
\end{equation}

Si on se fixe une équivalence homotopique d'espaces topologiques, 
$$
\xymatrix@C+5pt{\mathbb{S}^{m}_{top} \ar[r]^-{\Phi}_-{\cong} &\big| \mathbb{S}^{m} \big| },
$$
on obtient une fonction naturelle en $X$:
\begin{equation}\label{kanpi}
\xymatrix@R=3pt@C+20pt{\Big\{ \, x\in X_{m} \;\Big|\;  d_i(x)=a\;\;\; 0\leq i \leq m \,\Big\} \ar[r]^-{}& \pi_{m}(X,a)\\
\qquad \qquad\qquad\qquad x\ar@{}[r]|-{\longmapsto} & \underset{\text{pointée de $|\sigma^{x}|\circ\Phi$.}}{\overset{\text{La classe}}{\text{\scriptsize{d'homotopie}}}}}.
\end{equation}

\begin{proposition}\label{kankan}
Si $X$ est un complexe de Kan, la fonction \eqref{kanpi} est surjective et identifie deux éléments $x$ et $y$ si et seulement si, il existe $w\in X_{m+1}$ tel que: 
$$
\text{ $d_{m+1}w=x$, \; $d_{m}w=y$ \; et \; $d_{i}w=a$ \; pour\;  $0\leq i \leq m-1$. }
$$
\end{proposition}
\begin{proof}
On a définit $\pi_{m}(X,a)$ comme étant l'ensemble des classes à homotopies pointées près des fonctions continues pointées:
$$
\xymatrix{\big(\mathbb{S}^m_{top},\star\big)\ar[r]&\big(|X|,a\big)}\,.
$$

Vu que l'espace topologique $\mathbb{S}^m_{top}$ admet une structure de complexe cellulaire, l'ensemble $\pi_{m}(X,a)$ s'identifie alors à l'ensemble des morphismes $[(\mathbb{S}^m_{top},\star),(|X|,a)]_{{\bf Top}_*}$ de la catégorie de fractions ${\bf Top_*}\big[({\bf W}_\infty^{top})^{-1}\big]$ (voir le Corollaire 2.4.20 et le Théorème 2.4.19 de \cite{hovey}). Donc, à l'aide de l'équivalence faible $\xymatrix@C-5pt{\mathbb{S}^{m}_{top} \ar[r]^-{\Phi}_-{\cong} &\big| \mathbb{S}^{m} \big| }$, on obtient une bijection $[(|\mathbb{S}^m|,\star),(|X|,a)]_{{\bf Top}_*} \cong \pi_{m}(X,a)$.

D'un autre rappelons que le foncteur $|\,\cdot\,|$ détermine une bijection entre  $[(|\mathbb{S}^m|,\star),(|X|,a)]_{{\bf Top_*}}$ et l'ensemble de morphismes $[(\mathbb{S}^m,\star),(X,a)]_{\infty}^{\text{pointé}}$ de la catégorie homotopique $\simp_\star[{\bf W}_\infty^{-1}]$ (Voir le Théorème 3.6.7 et le Corollaire 2.4.24 de \cite{hovey}). 

Si on suppose que $X$ soit un complexe de Kan, on déduit une bijection:
\begin{equation}\label{isosi}
\xymatrix@R=3pt@C-10pt{
\pi_0\Big(\underline{\mathrm{Hom}}_{\simp_\star}\big((\mathbb{S}^m,\star),(X,a)\big)\Big)   \ar[rrr]^-{\cong} &&& \pi_{m}(X,a)\,,\\
[\sigma] & \ar@{|->}[r] && [|\sigma| \circ \Phi]}
\end{equation}
définie dans la classe d'un morphisme d'ensembles simpliciaux:
$$
\xymatrix@C-3pt{\mathbb{S}^{m}\ar[r]^-{\sigma} & X} \qquad \text{vérifiant \, que} \qquad \sigma_{0}(\star)=a \,,
$$
comme la classe à homotopie pointée près du morphisme: 
$$
\xymatrix@C+5pt{\mathbb{S}^{m}_{top} \ar[r]^-{\Phi}_-{\cong} &\big| \mathbb{S}^{m+1} \big| \ar[r]^-{|\sigma|} & |X| },
$$

On obtient ainsi un carré commutatif:
\begin{equation} \label{arribaaaa}
\xymatrix@R=8pt@C-20pt{ 
\underline{\mathrm{Hom}}_{\simp_\star}\big((\mathbb{S}^m,\star),(X,a)\big)_0 \ar[r]^-{\cong} \ar[ddd] &
\Big\{ \, x\in X_{m} \;\Big|\;  d_i(x)=a\; \; \;0\leq i \leq m \,\Big\} \ar[ddd]^-{\eqref{kanpi}} \\
\\   \\ 
\pi_0\Big(\underline{\mathrm{Hom}}_{\simp_\star}\big((\mathbb{S}^m,\star),(X,a)\big)\Big) \ar[r]^-{\cong}_-{\eqref{isosi}} &  \pi_{m}(X,a) \,; }
\end{equation} 
en particulier la fonction \eqref{kanpi} est surjective.

Pour montrer la Proposition il suffit de voir que si $\xymatrix@C-5pt{\Delta^{m}\ar[r]^-{x, y} & X}$ sont deux morphismes tels que les composés $\xymatrix@C-4pt{\partial\Delta^m\ar[r]^{\alpha^{m-1}} & \Delta^{m}\ar[r]^-{x, y} & X}$ sont égaux au morphisme constant $\xymatrix@C-6pt{\partial\Delta^{m} \ar[r]&\star \ar[r]^-{a}& X}$, alors:

Il existe un morphisme $\xymatrix{\Delta^{m}\times \Delta^1 \ar[r]^-{H} & X}$ tel que $\xymatrix@C-6pt{\partial\Delta^m\times \Delta^1 \ar[rr]^{\alpha^{m-1}\times\Delta^1} && \Delta^{m}\times\Delta^1 \ar[r]^-{H} & X}$ est le morphisme constant $\xymatrix@C-6pt{\partial\Delta^m\times \Delta^1 \ar[r]&\star \ar[r]^-{a}& X}$ et:
\begin{equation}
\xymatrix@R=15pt@C+5pt{
\Delta^{m} \times \Delta^0 \ar[r]^-{\text{proj}}_-{\cong} \ar@/^2pt/[rd]_-{{\footnotesize \Delta^{m}\times\delta_1}} & \Delta^m   \ar@<+3pt>@/^8pt/[rrd]^-{x}&& \\
  & \Delta^{m}\times \Delta^1 \ar[rr]|-{\;H\;}& & X\\
\Delta^m \times \Delta^0    \ar[r]_-{\text{proj}}^-{\cong} \ar@/_2pt/[ru]^-{{\footnotesize \Delta^{m}\times\delta_0}} & \Delta^m  \ar@<-3pt>@/_8pt/[rru]_-{y} &&}
\end{equation}
sont de diagrammes commutatifs si et seulement si, il existe un morphisme $\xymatrix@C-3pt{\Delta^{m+1}\ar[r]^-{w} & X}$ tel que le composé $\xymatrix@C-3pt{\Delta^{m}\ar[r]^-{\delta_i} &\Delta^{m+1}\ar[r]^-{w} & X}$ est égal à $x$ si $i=m+1$, égal à $y$ si $i=m$ et égal au morphisme constant $\xymatrix@C-6pt{\Delta^{m} \ar[r]&\star \ar[r]^-{a}& X}$ si $0\leq i \leq m-1$.

Supposons pour commencer qu'on a un morphisme $\xymatrix{\Delta^{m}\times \Delta^1 \ar[r]^-{H} & X}$ comme ci-dessus. Pour montrer l'existence du morphisme $\xymatrix@C-3pt{\Delta^{m+1}\ar[r]^-{w} & X}$ avec les propriétés désirées, on va construire un carré commutatif:
\begin{equation} \label{carreKAN}
\xymatrix@R+6pt@C+5pt{
(\Delta^{m+1}\times \Delta^0) \underset{(\partial\Delta^{m+1}\times\Delta^0)}{\bigsqcup} (\partial\Delta^{m+1}\times\Delta^1) \ar[r]^-{\psi} \ar[d]_-{\phi} &X\ar[d]\\
\Delta^{m+1}\times \Delta^1 \ar[r] &\phantom{..}\star\,,}
\end{equation}
où l'ensemble simplicial {\scriptsize $(\Delta^{m+1}\times \Delta^0) \underset{(\partial\Delta^{m+1}\times\Delta^0)}{\bigsqcup} (\partial\Delta^{m+1}\times\Delta^1)$} est le coproduit fibré des morphismes:
$$
\xymatrix@C+25pt@R+6pt{
\partial\Delta^{m+1}\times\Delta^0\ar[r]^{\partial\Delta^{m+1}\times \delta_1}\ar[d]_{\alpha^m\times \Delta^0} & \partial\Delta^{m+1}\times \Delta^1\;,\\
\Delta^{m+1}\times\Delta^0 & }
$$
de la façon suivante:  

$\phi$ est déduit de la propriété universelle des carrés cocartésiens à partir du carré commutatif:
$$
\xymatrix@C+25pt@R+6pt{
\partial\Delta^{m+1}\times\Delta^0\ar[r]^{\partial\Delta^{m+1}\times \delta_1}\ar[d]_{\alpha^m\times \Delta^0} & \partial\Delta^{m+1}\times \Delta^1 \ar[d]^-{\alpha^m\times\Delta^1}\\
\Delta^{m+1}\times\Delta^0 \ar[r]_-{\Delta^{m+1}\times \delta_1}& \Delta^{m+1}\times\Delta^1\,.}
$$
et $\psi$ à partir du carré commutatif:
\begin{equation}\label{psipsi}
\xymatrix@C+25pt@R+6pt{
\partial\Delta^{m+1}\times\Delta^0\ar[r]^{\partial\Delta^{m+1}\times \delta_1}\ar[d]_{\alpha^m\times \Delta^0} & \partial\Delta^{m+1}\times \Delta^1 \ar[d]^-{\psi_1}\\
\Delta^{m+1}\times\Delta^0 \ar[r]_-{\psi_2}& X\,;}
\end{equation}
où $\psi_2$ est le composé $\xymatrix{\Delta^{m+1}\times\Delta^0 \ar[r]^-{\text{proj}} & \Delta^{m+1}  \ar[r]^-{\sigma_m} & \Delta^m \ar[r]^-{x} & X }$.

Pour définir $\psi_1$ on voit l'ensemble simplicial $\partial\Delta^{m+1}\times\Delta^1$ comme le conoyau de la double flèche:
$$
\xymatrix@C-1pt{
\underset{0\leq i<j\leq m+1}{\bigsqcup}\big(\Delta^{m-1}\times\Delta^1\big)
\ar@<+7pt>[rrr]^-{\underset{i<j}{\sqcup} \; (\nu_{i}\circ \delta^{m-1}_{j-1}\times \Delta^1) }
\ar@<-7pt>[rrr]_-{\underset{i<j}{\sqcup} \; (\nu_{j} \circ \delta^{m-1}_{i}\times \Delta^1)}&&&
\underset{0\leq l\leq m+1}{\bigsqcup}\big(\Delta^{m}\times\Delta^1\big)},
$$
où $\nu_{l}$ note l'inclusion de la $l$-ième composante $\xymatrix@-5pt{\Delta^m\;\ar@{^(->}[r]&\underset{0\leq l\leq m+1}{\bigsqcup}\Delta^{m}}$; alors  $\psi_1$ est définit comme le morphisme induit par $\xymatrix{\underset{0\leq l\leq m+1}{\bigsqcup}\big(\Delta^{m}\times\Delta^1\big) \ar[r]^-{\underset{l}{\sqcup} \rho_l} & X}$ où les $\rho_l$ sont donnés par les règles:
\begin{equation*}
\begin{split}
\rho_0\; = \;  \cdots \; = \; \rho_{m-1} \; = & \; \Big( \xymatrix@C+60pt{\Delta^m\times\Delta^1 \ar[r]^-{ \text{constant à valeurs $a$}} & X } \Big)\\
\rho_m\; = \;  \Big( \xymatrix{\Delta^m\times\Delta^1 \ar[r]^-{H} & X } \Big) & \qquad\text{et} \qquad
\rho_{m+1} = \Big( \xymatrix{\Delta^m\times\Delta^1 \ar[r]^-{\text{proj}} & \Delta^m  \ar[r]^-{x} & X } \Big) 
\end{split}
\end{equation*}

On vérifie que \eqref{psipsi} est bien un carré commutatif, en remarquant qu'on a de diagrammes commutatifs pour tout $0\leq i\leq m+1$:
$$
\xymatrix@C+5pt{ 
\Delta^m\times \Delta^0 \ar[r]^{\Delta^m\times\delta_1} \ar[d]_{\delta_i\times\Delta^0}&  \Delta^m\times \Delta^1 \ar[d]^-{\rho_i}\\
\Delta^{m+1}\times\Delta^0 \ar[r]_-{\psi_2} & X}
$$

Maintenant qu'on a construit les morphismes $\phi$ et $\psi$ du carré \eqref{carreKAN}, remarquons que $\phi$ est en fait un monomorphisme et une $\infty$-équivalence faible; en particulier, vu que par hypothèse $X$ est un complexe de Kan on déduit l'existence d'un morphisme $\xymatrix{\Delta^{m+1}\times\Delta^1\ar[r]^-{\xi}&X}$ dans un triangle commutatif: 
\begin{equation}\label{triangg}
\xymatrix@R+6pt@C+5pt{
(\Delta^{m+1}\times \Delta^0) \underset{(\partial\Delta^{m+1}\times\Delta^0)}{\bigsqcup} (\partial\Delta^{m+1}\times\Delta^1) \ar[r]^-{\psi} \ar[d]_-{\phi} &X\\
\Delta^{m+1}\times \Delta^1 \ar@/_6pt/@{-->}[ru]_-{\xi} &}
\end{equation}

Posons $w=\Big(\xymatrix@C-5pt{\Delta^{m+1} \ar[r]_-{\cong}^-{\mathrm{proj}^{-1}}& \Delta^{m+1}\times\Delta^0 \ar[rr]^-{\Delta^{m+1}\times \delta_0}  &&  \Delta^{m+1}\times \Delta^1 \ar[r]^-{\xi} & X}\Big)$. On vérifie alors que le composé $w\circ d_i$ est égal au morphisme $\xymatrix@C-5pt{\Delta^{m}\ar[r]_-{\cong}^-{\mathrm{proj}^{-1}}& \Delta^{m}\times\Delta^0 \ar[rr]^-{\delta_i\times \delta_0}  &&  \Delta^{m}\times \Delta^1 \ar[r]^-{\xi} & X}$, que d'après le triangle \eqref{triangg} est égale au morphisme  $\xymatrix@C-5pt{\Delta^{m}\ar[r]_-{\cong}^-{\mathrm{proj}^{-1}}& \Delta^{m}\times\Delta^0 \ar[rr]^-{\Delta^m\times \delta_0}  &&  \Delta^{m}\times \Delta^1 \ar[r]^-{\rho_i} & X}$ pour tout $0\leq i\leq m+1$. 

De la définition des morphismes $\rho_i$ on déduit que $w\circ d_i= a$ si $0\leq i\leq m-1$, $w\circ d_m = y$ et $w\circ d_{m+1} = x$ comme désiré.

Réciproquement supposons qu'on part d'un morphisme $\xymatrix@C-3pt{\Delta^{m+1}\ar[r]^-{w} & X}$ tel que $w\circ d_i= a$ si $0\leq i\leq m-1$, $w\circ d_m = y$ et $w\circ d_{m+1} = x$. Pour montrer l'existence du morphisme $\xymatrix{\Delta^{m}\times \Delta^1 \ar[r]^-{H} & X}$ cherché, on considère par ailleurs un carré commutatif:
\begin{equation} \label{carreKAN2}
\xymatrix@R+6pt@C+5pt{  
(\Delta^{m+1}\times \partial\Delta^1) \underset{(\Lambda^{m+1,m}\times\partial\Delta^1)}{\bigsqcup} (\Lambda^{m+1,m}\times\Delta^1) \ar[r]^-{\psi} \ar[d]_-{\phi} &X\ar[d]\\
\Delta^{m+1}\times \Delta^1 \ar[r] &\phantom{..}\star\,,}
\end{equation}
où l'ensemble simplicial {\scriptsize $(\Delta^{m+1}\times \partial\Delta^1) \underset{(\Lambda^{m+1,m}\times\partial\Delta^1)}{\bigsqcup} (\Lambda^{m+1,m}\times\Delta^1)$} est le coproduit fibré des morphismes:
$$
\xymatrix@C+25pt@R+6pt{
\Lambda^{m+1,m} \times\partial\Delta^1 \ar[r]^-{\Lambda^{m+1,m}\times\alpha^0 }\ar[d]_-{\alpha^{m,m}\times \partial\Delta^1}&
\Lambda^{m+1,m} \times \Delta^1  \\
\Delta^{m+1}\times \partial\Delta^1 
&}
$$

Le morphisme $\phi$ de \eqref{carreKAN2} est déduit du carré commutatif:
$$
\xymatrix@C+25pt@R+6pt{
\Lambda^{m+1,m} \times\partial\Delta^1 \ar[r]^-{\Lambda^{m+1,m}\times\alpha^0 }\ar[d]_-{\alpha^{m,m}\times \partial\Delta^1}&
\Lambda^{m+1,m} \times \Delta^1  \ar[d]^-{\alpha^{m,m}\times\Delta^1}\\
\Delta^{m+1}\times \partial\Delta^1 \ar[r]_-{\Delta^{m+1}\times \alpha^0}
& \Delta^{m+1}\times\Delta^1}
$$
et le morphisme $\psi$ du carré: 
\begin{equation}\label{carree}
\xymatrix@C+25pt@R+6pt{
\Lambda^{m+1,m} \times\partial\Delta^1 \ar[r]^-{\Lambda^{m+1,m}\times\alpha^0 }\ar[d]_-{\alpha^{m,m}\times \partial\Delta^1}&
\Lambda^{m+1,m} \times \Delta^1\ar[d]^-{\psi_1}\\  \Delta^{m+1}\times \partial\Delta^1
 \ar[r]_-{\psi_2}& X}
\end{equation}
qu'on définit comme suit: Pour définir $\psi_2$ remarquons qu'on a un carré cocartésien:
$$
\xymatrix@C+25pt@R+8pt{\emptyset \ar[r]\ar[d] &  \Delta^{m+1}\times\Delta^0\ar[d]^-{\Delta^{m+1}\times \overline{\delta_1}}\\  
\Delta^{m+1}\times \Delta^0 \ar[r]_-{\Delta^{m+1} \times \overline{\delta^0}}& \Delta^{m+1}\times\partial\Delta^1\\
}
$$
où $\overline{\delta_i}$ est la seul flèche tel que $\big(\xymatrix{\Delta^0\ar[r]^-{\delta_i} &\Delta^1}\big)=\big(\xymatrix{\Delta^0 \ar[r]^-{\overline{\delta_i}} & \partial\Delta^1 \ar[r]^-{\alpha^0} & \Delta^1}\big)$; alors $\psi_2$ est le morphisme induit par le carré commutatif:
$$
\xymatrix@R-10pt@C-20pt{\emptyset \ar[rrr]\ar[ddd] && & \Delta^{m+1}\times\Delta^0\ar[d]^-{\text{proj}}\\  
&&&\Delta^{m+1} \ar[d]^-{\sigma_m}\\ &&&\Delta^m\ar[d]^-{x}\\
\Delta^{m+1}\times \Delta^0 \ar[rr]_-{\text{proj}}&& \Delta^{m+1}\ar[r]_-{w}&X\\}
$$

Le morphisme $\psi_1$ est définit en voyant l'ensemble simplicial $\Lambda^{m+1,m}\times \Delta^1$ comme le conoyau de la double flèche:
$$
\xymatrix@C-1pt{
\underset{i,j\neq m}{\underset{0\leq i<j\leq m+1}{\bigsqcup}}\big(\Delta^{m-1}\times\Delta^1\big)
\ar@<+7pt>[rrr]^-{\underset{i<j}{\sqcup} \; (\nu_{i}\circ \delta^{m-1}_{j-1}\times \Delta^1) }
\ar@<-7pt>[rrr]_-{\underset{i<j}{\sqcup} \; (\nu_{j} \circ \delta^{m-1}_{i}\times \Delta^1)}&&&
\underset{l\neq m}{\underset{0\leq l\leq m+1}{\bigsqcup}}\big(\Delta^{m}\times\Delta^1\big)},
$$
où $\nu_{l}$ note l'inclusion de la $l$-ième composante $\xymatrix@-5pt{\Delta^m\;\ar@{^(->}[r]&\underset{l\neq m}{\underset{0\leq l\leq m+1}{\bigsqcup}}\Delta^{m}}$; alors $\psi_1$ est définit comme le morphisme induit par $\xymatrix{\underset{l\neq m}{\underset{0\leq l\leq m+1}{\bigsqcup}}\big(\Delta^{m}\times\Delta^1\big) \ar[r]^-{\underset{l}{\sqcup} \tau_l} & X}$ où les $\tau_l$ sont donnés par les règles:
$$
\tau_0\, = \,  \cdots \, = \, \tau_{m-1} \, = \, \Big( \xymatrix@C+55pt{\Delta^m\times\Delta^1 \ar[r]^-{ \text{constant à valeurs $a$}} & X } \Big)
$$
$$\text{et}$$
$$
\tau_{m+1} = \Big( \xymatrix{\Delta^m\times\Delta^1 \ar[r]^-{\text{proj}} & \Delta^m  \ar[r]^-{x} & X } \Big) 
$$

Pour montrer que le carré \eqref{carree} est commutatif, il suffit de remarquer qu'on a des diagrammes commutatifs:
$$
\vcenter{\xymatrix{\Delta^m\times\Delta^0  \ar[r]^-{\delta_i\times\Delta^0} \ar[d]_-{\Delta^m\times\delta_1} & \Delta^{m+1}\times\Delta^0 \ar[d]^-{x\circ\sigma_m\circ\text{proj}} \\
\Delta^m\times\Delta^1 \ar[r]_{\tau_i} &X  }}
\quad\qquad\text{et}\qquad\quad
\vcenter{\xymatrix{\Delta^m\times\Delta^0  \ar[r]^-{\delta_i\times\Delta^0} \ar[d]_-{\Delta^m\times\delta_0} & \Delta^{m+1}\times\Delta^0 \ar[d]^-{w\circ\text{proj}} \\
\Delta^m\times\Delta^1 \ar[r]_{\tau_i} &X  }}
$$
toujours que $0\leq i\leq m+1$ avec $i\neq m$. 

Après avoir défini les morphismes $\phi$ et $\psi$ du carré \eqref{carreKAN2}, on note que $\phi$ est un monomorphisme et une $\infty$-équivalence faible. Donc, vu que $X$ est un complexe de Kan il existe un morphisme $\xymatrix{\Delta^{m+1}\times\Delta^1\ar[r]^-{\xi}&X}$ comme dans le triangle commutatif: 
\begin{equation}\label{triangg2}
\xymatrix@R+6pt@C+5pt{
(\Delta^{m+1}\times \partial\Delta^1) \underset{(\Lambda^{m+1,m}\times\Delta^1)}{\bigsqcup} (\Lambda^{m+1,m}\times\Delta^1) \ar[r]^-{\psi} \ar[d]_-{\phi} &X\\
\Delta^{m+1}\times \Delta^1 \ar@/_6pt/@{-->}[ru]_-{\xi} &}
\end{equation}

On définit $H= \big(\xymatrix@C-10pt{\Delta^m\times\Delta^1\ar[rr]^-{\delta_m\times \Delta^1} && \Delta^{m+1}\times\Delta^1 \ar[r]^-{\xi} & X}\big)$. 

Pour montrer que $\xymatrix@C-7pt{\partial\Delta^m\times \Delta^1 \ar[rr]^{\alpha^{m-1}\times\Delta^1} && \Delta^{m}\times\Delta^1 \ar[r]^-{H} & X}$ est le morphisme constant à valeurs $a$, il suffit de montrer que $\xymatrix@C-3pt{\Delta^{m-1}\times\Delta^1\ar[rr]^-{(\delta_m\circ\delta_i)\times\Delta^1}&& \Delta^{m+1}\times \Delta^1 \ar[r]^-{\xi} & X}$ est le morphisme constant à valeurs $a$ quelque soit $0\leq i\leq m$. Mais du diagramme commutatif \eqref{triangg2} on déduit que le morphisme composé $\xi\circ\big((\delta_m\circ\delta_i)\times\Delta^1\big)$ est égale au morphisme $\xymatrix@C-10pt{\Delta^{m-1}\times\Delta^1\ar[rr]^-{\delta_{m,i}\times\Delta^1}&& \Lambda^{m+1,m}\times \Delta^1 \ar[r]^-{\psi_1} & X}$ où $\delta_{m,i}$  est le seul morphisme dans le triangle commutatif:
$$
\xymatrix{&&\Lambda^{m+1,m}\ar[d]^-{\alpha^{m,m}}\\
\Delta^{m-1}  \ar[r]_{\delta_i}\ar@/^5pt/@{-->}[rru]^-{\delta_{m,i}} & \Delta^m  \ar[r]_{\delta_m} & \Delta^{m+1} \,; }
$$
donc, il est constant à valeurs $a$ parce que $\xymatrix@C-8pt{\Delta^{m-1}\times \Delta^1 \ar[rr]^-{\delta_j\times\Delta^1} && \Delta^m\times \Delta^1 \ar[r]^-{\tau_i} & X}$ l'est quelques soient $0\leq i,j\leq m$.
 
Finalement pour montrer qu'on a un diagramme commutatif:
\begin{equation}
\xymatrix@R=15pt@C+5pt{
\Delta^{m} \times \Delta^0 \ar[r]^-{\text{proj}} \ar@/^2pt/[rd]_-{{\footnotesize \Delta^{m}\times\delta_1}} & \Delta^m   \ar@<+3pt>@/^8pt/[rrd]^-{x}&& \\
  & \Delta^{m}\times \Delta^1 \ar[rr]|-{\;H\;}& & X\;,\\
\Delta^m \times \Delta^0    \ar[r]|-{\text{proj}} \ar@/_2pt/[ru]^-{{\footnotesize \Delta^{m}\times\delta_0}} & \Delta^m  \ar@<-3pt>@/_8pt/[rru]_-{y} &&}
\end{equation}
remarquons que d'après le triangle commutatif \eqref{triangg2} on a des carrés commutatifs:
$$
\vcenter{\xymatrix{
\Delta^m\times\Delta^0 \ar[r]^-{\text{proj}} \ar[d]_-{\Delta^m\times\delta_0} & \Delta^m\ar[d]^-{w\circ\delta_m \, = \, y} \\
\Delta^m \times \Delta^1 \ar[r]_-{H} & X }}\qquad\text{et}\qquad 
\vcenter{\xymatrix{
\Delta^m\times\Delta^0 \ar[r]^-{\text{proj}} \ar[d]_-{\Delta^m\times\delta_1} & \Delta^m\ar[d]^-{x\circ \sigma_m\circ\delta_m \, = \, x} \\
\Delta^m \times \Delta^1 \ar[r]_-{H} & X \,.}}
$$
\end{proof}

Rappelons aussi:

\begin{proposition}\label{kankan2}
Soit $X$ un complexe de Kan. Si $x$, $y$ et $w$ sont de $m$-simplexes de $X$ tels que:  
$$
\alpha^m_{X}(x)=\alpha^m_{X}(y)=\alpha^m_{X}(w)=(a,\dots,a);
$$ 
alors, $\overline{\alpha}(x)\overline{\alpha}(y)=\overline{\alpha}(w)$ (où $\overline{\alpha}$ est la fonction \eqref{kanpi} ci-dessus) si et seulement si, il existe $z\in X_{m+1}$ tel que: 
$$
\text{ $d_{m+1}z=x$, \; $d_{m}z=w$, \; $d_{m-1}z=y$ \; et \; $d_{i}z=a$ \; pour \; $0\leq i \leq m-2$. }
$$
\end{proposition}

\section{Catégories de modèles simpliciaux d'ordre $n$}\label{catsimpn}



\renewcommand{\thesubsection}{\S\thesection.\arabic{subsection}}
\subsection{}\;  \label{catmodsimpnn}
\renewcommand{\thesubsection}{\thesection.\arabic{subsection}}

Si $0 \leq n \leq \infty$, une catégorie de modèles $(\C,{\bf W}, {\bf cof},{\bf fib},\underline{\mathrm{Hom}}_{\C})$ enrichie dans la catégorie de modèles monoïdale $(\simp,{\bf W}_{n}, {\bf mono},{\bf fib}_{n},\times)$ du Lemme \ref{ntypessp} est appelée une \emph{catégorie de modèles simpliciale d'ordre $n$}. 

De façon explicite (voir la Définition 4.2.18 de \cite{hovey}), une catégorie de modèles simpliciale d'ordre $n$ est une catégorie $\C$, munie d'une structure de catégorie de modèles $(\C,{\bf W}, {\bf cof},{\bf fib})$ et d'un enrichissement dans la catégorie cartésienne fermée des ensembles simpliciaux $\simp$ (voir \cite{kelly}):
$$
\xymatrix@C+35pt{\C^{op} \times \C \ar[r]^-{\underline{\mathrm{Hom}}_{\C}(\,\cdot\,,\,\cdot\,)} & \simp},
$$
vérifiant les deux propriétés suivantes:
\begin{itemize}
\item[{\bf CMS1.}] L'enrichissement est tensoré et cotensoré, \emph{i.e.} si $X$ et $Y$ sont des objets quelconques de $\C$, on a des adjonctions:
\begin{equation}\label{lesdeux}  
\C
\vcenter{\xymatrix@C+15pt{
\phantom{a}\ar@{}[r]|{\perp}\ar@<-4pt>@/_10pt/[r]_{\underline{\mathrm{Hom}}_{\C}(X,\,\cdot\,)}&\phantom{a}
\ar@<-4pt>@/_10pt/[l]_-{X\otimes\,\cdot\, }}}
\simp
\qquad\text{et}\qquad
\C^{op}
\vcenter{\xymatrix@C+15pt{
\phantom{a}\ar@{}[r]|{\perp}\ar@<-4pt>@/_10pt/[r]_{\underline{\mathrm{Hom}}_{\C}(\,\cdot\,,Y)}&\phantom{a}
\ar@<-4pt>@/_10pt/[l]_-{Y^{\,\cdot\,} }}}
\simp\,.
\end{equation}
\item[{\bf CMS2}.]  Si $\xymatrix@C-8pt{A\ar[r]^j&B}$ est un monomorphisme d'ensembles simpliciaux et $\xymatrix@C-8pt{X\ar[r]^q&Y}$ une cofibration de $\C$, alors le morphisme $\varphi$ dans le diagramme somme amalgamée suivant:
\begin{equation}\label{cardul}
\xymatrix@R-10pt@C-10pt{
X\otimes A \ar[dd]_{X\otimes \,j}\ar[rr]^{q\,\otimes A} && 
Y\otimes A \ar[dd]\ar@/^15pt/[rddd]^{Y\otimes \,j}&\\\\
X\otimes B\ar[rr]\ar@/_15pt/[rrrd]_{q\,\otimes B}&&
X\otimes B\underset{X\otimes A}{\bigsqcup}Y\otimes A\ar@{-->}[rd]|(.65)\varphi&\\
&&&Y\otimes B\,,}
\end{equation}
est une cofibration de $\C$, laquelle est aussi une équivalence faible si $j$ est une $n$-équivalence faible d'ensembles simpliciaux, ou si $q$ est une équivalence faible de $\C$.
\end{itemize}

\phantom{}

Remarquons que si $(\C,{\bf W}, {\bf cof},{\bf fib},\underline{\mathrm{Hom}}_{\C})$ est une catégorie de modèles simpliciale d'ordre $n$ et on note $[\;\cdot\;,\;\cdot\;]_{\bf W}$ l'ensemble des morphismes de la catégorie des fractions $\C\big[{\bf W}^{-1}\big]$; alors la transformation naturelle:
$$
\xymatrix{ 
\mathrm{Hom}_{\C}( \;\cdot\;,\;\cdot\;) \ar[r] & 
[\;\cdot\;,\;\cdot\;]_{\bf W}}
$$
induit une bijection:
$$
\pi_0\big( \underline{\mathrm{Hom}}_{\C}(X,Y) \big) \; \cong \; \big[ X, Y\big]_{\bf W} 
$$
toujours que $X$ soit un objet cofibrant de $\C$ et $Y$ un objet fibrant. 

En effet, de la propriété {\bf CMS2} il se suit que si $X$ est un objet cofibrant de $\C$ alors:
$$
\vcenter{\xymatrix@-8pt{\underset{\phantom{.}}{X \sqcup X}\ar@<-4pt>@{>->}@/_5pt/[rd]_-{i_{0}+i_{1}}\ar[rr]^-{\mathrm{id}+\mathrm{id}} && X \\ 
&X\otimes \Delta^1 \ar@<-4pt>@/_5pt/[ru]_-{p}\ar@{}[ru]|-{\text{\rotatebox[origin=c]{225}{\Large $\widetilde{\phantom{w}}$}}}&}}
$$
est un objet cylindre de la catégorie de modèles $(\C,{\bf W}, {\bf cof},{\bf fib})$, où $p$ est l'équivalence faible:
$$
\vcenter{\xymatrix@C+8pt{X\otimes \Delta^1 \ar[r]^-{X\otimes \sigma_0} & X\otimes\Delta^0 \, \cong \, X}}
$$
et les cofibrations $i_0$ et $i_1$ sont les morphismes composés:
$$
\vcenter{\xymatrix@C+8pt{X\,\cong\,X\otimes \Delta^0 \ar[r]^-{X\otimes \delta_0} & X\otimes \Delta^1}}
\qquad\text{et}\qquad
\vcenter{\xymatrix@C+8pt{X\,\cong\,X\otimes \Delta^0 \ar[r]^-{X\otimes \delta_1} & X\otimes \Delta^1}}
$$
respectivement.

Plus encore:

\begin{lemme}\label{RHOMloc}
Soit $0\leq n \leq \infty$ et $(\C,{\bf W}, {\bf cof},{\bf fib},\underline{\mathrm{Hom}}_{\C})$ une catégorie de modèles simpliciale d'ordre $n$. Si $X$ est un objet cofibrant de $\C$ et $Y$ un objet fibrant, le type d'homotopie de l'ensemble simplicial $\underline{\mathrm{Hom}}_{\C}(X,Y)$, est égale au \emph{complexe de fonction dérivé} ${\bf R}\mathrm{Hom}_{\C}(X,Y)$ défini par la localisation simpliciale $\mathcal{L}(\C,{\bf W})$ de \cite{DK1,DK2,DK}, ou construit à partir de résolutions fibrants et cofibrants comme dans \S5.4 de \cite{hovey} ou le Chapitre 18 de \cite{hirschhorn}. 
\end{lemme}
\begin{proof}
Si $(\C,{\bf W}, {\bf cof},{\bf fib},\underline{\mathrm{Hom}}_{\C})$ est une catégorie de modèles simpliciale d'ordre $n$ et on se donne un autre nombre $m$ tel que $n \leq m\leq \infty$, alors on a que $(\C,{\bf W}, {\bf cof},{\bf fib},\underline{\mathrm{Hom}}_{\C})$ est aussi une catégorie de modèles simpliciale d'ordre $m$ parce que dans ce cas ${\bf W}_{m}\subset {\bf W}_{n}$. Il se suit en particulier que $(\C,{\bf W}, {\bf cof},{\bf fib},\underline{\mathrm{Hom}}_{\C})$ est une catégorie de modèles simpliciale d'ordre $\infty$. 

Le résultat est une conséquence du Corollaire 4.7 de \S4.3 dans \cite{DK}.
\end{proof}

On définit une \emph{catégorie de modèles simpliciale pointée d'ordre $n$}, comme une catégorie de modèles enrichie dans la catégorie de modèles monoïdale $(\simp_{\star},\pi^{-1}{\bf W}_{n}, {\bf mono},\pi^{-1}{\bf fib}_{n},\wedge)$ du Corollaire \ref{pointeco2}. 

Explicitement, une catégorie de modèles simpliciale pointée d'ordre $n$ est une catégorie $\C$, munie d'une structure de catégorie de modèles $(\C,{\bf W}, {\bf cof},{\bf fib})$ et d'un enrichissement dans la catégorie monoïdale symétrique fermée $\big(\simp_{\star},\wedge\big)$:
$$
\xymatrix@C+45pt{\C^{op} \times \C \ar[r]^-{\underline{\mathrm{Hom}}^{\text{\tiny{$\text{pointé}$}}}_{\C}(\,\cdot\,,\,\cdot\,)} & \simp_{\star}},
$$
vérifiant les propriétés suivantes:
\begin{itemize}
\item[{\bf CMS$1^{\text{\tiny{$\text{pointé}$}}}$.}] L'enrichissement est tensoré et cotensoré, \emph{i.e.} si $X$ et $Y$ sont des objets quelconques de $\C$, on a des adjonctions:
\begin{equation}\label{lesdeuxp}  
\C
\vcenter{\xymatrix@C+15pt{
\phantom{a}\ar@{}[r]|{\perp}\ar@<-4pt>@/_10pt/[r]_{\underline{\mathrm{Hom}}^{\text{\tiny{$\text{pointé}$}}}_{\C}(X,\,\cdot\,)}&\phantom{a}
\ar@<-4pt>@/_10pt/[l]_-{X\wedge\,\cdot\, }}}
\simp_{\star}
\qquad\text{et}\qquad
\C^{op}
\vcenter{\xymatrix@C+15pt{
\phantom{a}\ar@{}[r]|{\perp}\ar@<-4pt>@/_10pt/[r]_{\underline{\mathrm{Hom}}^{\text{\tiny{$\text{pointé}$}}}_{\C}(\,\cdot\,,Y)}&\phantom{a}
\ar@<-4pt>@/_10pt/[l]_-{Y^{\wedge\,\cdot\,} }}}
\simp_{\star}\,.
\end{equation}
\item[{\bf CMS$2^{\text{\tiny{$\text{pointé}$}}}$.}]  Si $\xymatrix@C-4pt{A\ar[r]^j&B}$ est un monomorphisme d'ensembles simpliciaux pointés et $\xymatrix@C-8pt{X\ar[r]^q&Y}$ une cofibration de $\C$, alors le morphisme $\psi$ dans le diagramme somme amalgamée suivant:
\begin{equation}\label{cardulp}
\xymatrix@R-10pt@C-10pt{
X\wedge A \ar[dd]_{X\wedge \,j}\ar[rr]^{q\,\wedge A} && 
Y\wedge A \ar[dd]\ar@/^15pt/[rddd]^{Y\wedge \,j}&\\\\
X\wedge B\ar[rr]\ar@/_15pt/[rrrd]_{q\,\wedge B}&&
X\wedge B\underset{X\wedge A}{\bigvee}Y\wedge A\ar@{-->}[rd]|(.65)\psi&\\
&&&Y\wedge B\,,}
\end{equation}
est une cofibration de $\C$, laquelle est aussi une équivalence faible si $j$ est une $n$-équivalence faible d'ensembles simpliciaux pointés, ou si $q$ est une équivalence faible de $\C$.
\end{itemize}

\phantom{}

De fa\c con équivalente: 

\begin{lemme}\label{equipointe}
Si $(\C,{\bf W}, {\bf cof},{\bf fib})$ est une catégorie de modèles, munie d'un enrichissement dans la catégorie cartésienne fermée des ensembles simpliciaux:
$$
\xymatrix@C+35pt{\C^{op} \times \C \ar[r]^-{\underline{\mathrm{Hom}}_{\C}(\,\cdot\,,\,\cdot\,)} & \simp}\,;
$$
alors les énoncés suivants sont équivalents:
\begin{enumerate}
\item $(\C,{\bf W}, {\bf cof},{\bf fib},\underline{\mathrm{Hom}}_{\C})$ est une catégorie modèles simpliciale d'ordre $n$ et $\C$ est une catégorie pointée (\emph{i.e.} si $\emptyset$ est un objet initial de $\C$ et $\star$ est un objet final, le seul morphisme $\xymatrix@C-10pt{\emptyset\ar[r]&\star}$ est un isomorphisme).
\item Le foncteur $\underline{\mathrm{Hom}}_{\C}(\,\cdot\,,\,\cdot\,)$ admet une extension:
$$
\xymatrix@C+35pt{\C^{op} \times \C \ar[r]^-{\underline{\mathrm{Hom}}^{\text{\tiny{$\text{pointé}$}}}_{\C}(\,\cdot\,,\,\cdot\,)} & \simp_{\star}},
$$
le long du foncteur qu'oublie $\pi :\xymatrix@-8pt{\simp_{\star}\ar[r]&\simp}$ tel que $(\C,{\bf W}, {\bf cof},{\bf fib},\underline{\mathrm{Hom}}^{\text{\tiny{$\text{pointé}$}}}_{\C})$ soit une catégorie de modèles simpliciale pointée d'ordre $n$.
\end{enumerate}
\end{lemme}
\begin{proof}
Rappelons d'ailleurs la preuve de l'affirmation suivante: Étant donnée une catégorie $\C$ munie d'un objet initial $\emptyset$ et d'un objet final $\star$, si on se donne un enrichissement de $\C$ dans la catégorie cartésienne fermée $\simp$:
$$
\xymatrix@C+35pt{\C^{op} \times \C \ar[r]^-{\underline{\mathrm{Hom}}_{\C}(\,\cdot\,,\,\cdot\,)} & \simp\,,}
$$
les énoncés suivants sont équivalents:
\begin{enumerate}
\item La catégorie $\C$ est pointée, \emph{i.e.} le seul morphisme $\xymatrix@C-10pt{\emptyset\ar[r]&\star}$ est un isomorphisme.
\item Le foncteur $\underline{\mathrm{Hom}}_{\C}(\,\cdot\,,\,\cdot\,)$ admet une extension:
$$
\xymatrix@C+35pt{\C^{op} \times \C \ar[r]^-{\underline{\mathrm{Hom}}^{\text{\tiny{$\text{pointé}$}}}_{\C}(\,\cdot\,,\,\cdot\,)} & \simp_{\star}},
$$
le long du foncteur $\pi :\xymatrix@-8pt{\simp_{\star}\ar[r]&\simp}$. 
\end{enumerate}

En effet, si on suppose que $\C$ est une catégorie pointée, on définit pour $X$ et $Y$ des objets de $\C$ l'ensemble simplicial pointé :
\begin{equation}\label{HOMpointe}
\underline{\mathrm{Hom}}^{\text{\tiny{$\text{pointé}$}}}_{\C}(X,Y)=\big(\underline{\mathrm{Hom}}_{\C}(X,Y), \star_{X}^Y\big)\,,
\end{equation}
où $\star_{X}^Y$ note le morphisme null $\xymatrix@C-8pt{X\ar[r]&\star\cong\emptyset\ar[r]&Y}$. On montre facilement que $\underline{\mathrm{Hom}}^{\text{\tiny{$\text{pointé}$}}}_{\C}$ est bien un bifoncteur. 

Réciproquement, supposons que $\underline{\mathrm{Hom}}_{\C}$ admet une extension $\underline{\mathrm{Hom}}^{\text{\tiny{$\text{pointé}$}}}_{\C}$ le longe du foncteur d'oubli $\pi :\xymatrix@-8pt{\simp_{\star}\ar[r]&\simp}$; alors l'ensemble simplicial $\underline{\mathrm{Hom}}_{\C}(\star,\emptyset)$ est non vide, \emph{i.e.} le morphisme canonique $\xymatrix@C-10pt{\emptyset\ar[r]&\star}$ est un isomorphisme.

Deuxièmement, étant donnée une catégorie pointée $\C$ munie d'un enrichissement dans la catégorie cartésienne fermée $\simp$:
$$
\xymatrix@C+35pt{\C^{op} \times \C \ar[r]^-{\underline{\mathrm{Hom}}_{\C}(\,\cdot\,,\,\cdot\,)} & \simp\,,}
$$
montrons que $\underline{\mathrm{Hom}}_{\C}$ vérifie la propriété {\bf CMS1}, si et seulement si l'enrichissement $\underline{\mathrm{Hom}}_{\C}^{\text{\tiny{$\text{pointé}$}}}$ de \eqref{HOMpointe} vérifie la propriété {\bf CMS$1^{\text{\tiny{$\text{pointé}$}}}$}. 

En effet, soient $X$ et $Y$ des objets de $\C$ et considérons des adjonctions:
\begin{equation}\label{tensor}
\C
\vcenter{\xymatrix@C+15pt{
\phantom{a}\ar@{}[r]|{\perp}\ar@<-4pt>@/_10pt/[r]_{\underline{\mathrm{Hom}}_{\C}(X,\,\cdot\,)}&\phantom{a}
\ar@<-4pt>@/_10pt/[l]_-{X\otimes\,\cdot\, }}}
\simp
\qquad\text{et}\qquad
\C^{op}
\vcenter{\xymatrix@C+15pt{
\phantom{a}\ar@{}[r]|{\perp}\ar@<-4pt>@/_10pt/[r]_{\underline{\mathrm{Hom}}_{\C}(\,\cdot\,,Y)}&\phantom{a}
\ar@<-4pt>@/_10pt/[l]_-{Y^{\,\cdot\,} }}}
\simp\,;
\end{equation}
on définit des foncteurs adjoints:
\begin{equation}\label{wedge}
\C
\vcenter{\xymatrix@C+15pt{
\phantom{a}\ar@{}[r]|{\perp}\ar@<-4pt>@/_10pt/[r]_{\underline{\mathrm{Hom}}^{\text{\tiny{$\text{pointé}$}}}_{\C}(X,\,\cdot\,)}&\phantom{a}
\ar@<-4pt>@/_10pt/[l]_-{X\wedge\,\cdot\, }}}
\simp_{\star}
\qquad\text{et}\qquad
\C^{op}
\vcenter{\xymatrix@C+15pt{
\phantom{a}\ar@{}[r]|{\perp}\ar@<-4pt>@/_10pt/[r]_{\underline{\mathrm{Hom}}^{\text{\tiny{$\text{pointé}$}}}_{\C}(\,\cdot\,,Y)}&\phantom{a}
\ar@<-4pt>@/_10pt/[l]_-{Y^{\wedge\,\cdot\,} }}}
\simp_{\star}\,,
\end{equation}
en considérant pour tout ensemble simplicial pointé $A=(A,a_0)$ un carré cocartésien et un carré cartésien de $\C$:
\begin{equation}\label{lescwe}
\vcenter{\xymatrix{
\star\ar[r] & X\wedge A \\
X\otimes \star \ar[u]\ar[r]_-{X\otimes a_0} & X\otimes A\ar[u]}}
\,\;\qquad\text{et}\qquad\;\,
\vcenter{\xymatrix{
\star\ar[d] & Y^{\wedge\,A\,}\ar[l]\ar[d] \\
Y^{\,\star\,} & Y^{\, A\,}\ar[l]^-{Y^{\,a_0}} \,;}}
\end{equation}
respectivement.

Pour montrer que les foncteur \eqref{wedge} ainsi définissent sont effectivement des foncteurs adjoints, il suffit de appliquer les foncteurs $\mathrm{Hom}_{\C}(\,-\,,Y)$ et $\mathrm{Hom}_{\C}(X,\,-\,)$ aux carrés \eqref{lescwe} respectivement, ayant en tête que $\mathrm{Hom}_{\C}(\star,Y)\cong\star\cong\mathrm{Hom}_{\C}(X,\star)$ pour $\C$ pointée et qu'on a un carré cartésien d'ensembles: 
$$
\xymatrix{
\star \ar[d]_{\star_X^Y} & \mathrm{Hom}_{\simp_\star}\Big(A,\underline{\mathrm{Hom}}_{\C}^{\text{\tiny{$\text{pointé}$}}}(X,Y)\Big) \ar[l]\ar[d]^-{\pi} \\
\mathrm{Hom}_{\simp}\big(\star,\underline{\mathrm{Hom}}_{\C}(X,Y)\big) & \mathrm{Hom}_{\simp}\Big(A,\underline{\mathrm{Hom}}_{\C}(X,Y)\Big)  \ar[l]^-{-\circ a_0}}
$$

Remarquons aussi qu'on vérifie que les adjonctions composées:
\begin{equation}\label{tensorwedge}
\begin{split}
\C
\vcenter{\xymatrix@C+15pt{
\phantom{a}\ar@{}[r]|{\perp}\ar@<-4pt>@/_10pt/[r]_{\underline{\mathrm{Hom}}^{\text{\tiny{$\text{pointé}$}}}_{\C}(X,\,\cdot\,)}&\phantom{a}
\ar@<-4pt>@/_10pt/[l]_-{X\wedge\,\cdot\, }}}
\simp_{\star}&
\vcenter{\xymatrix@C+15pt{
\phantom{a}\ar@{}[r]|{\perp}\ar@<-4pt>@/_10pt/[r]_{\pi}&\phantom{a}
\ar@<-4pt>@/_10pt/[l]_-{(\,\cdot\,)_{+}}}}
\simp\\\text{et}\,&\\
\C^{op}
\vcenter{\xymatrix@C+15pt{
\phantom{a}\ar@{}[r]|{\perp}\ar@<-4pt>@/_10pt/[r]_{\underline{\mathrm{Hom}}^{\text{\tiny{$\text{pointé}$}}}_{\C}(\,\cdot\,,Y)}&\phantom{a}
\ar@<-4pt>@/_10pt/[l]_-{Y^{\wedge\,\cdot\,} }}}
\simp_{\star}&
\vcenter{\xymatrix@C+15pt{
\phantom{a}\ar@{}[r]|{\perp}\ar@<-4pt>@/_10pt/[r]_{\pi}&\phantom{a}
\ar@<-4pt>@/_10pt/[l]_-{(\,\cdot\,)_{+}}}}
\simp\,.
\end{split}
\end{equation}
sont isomorphes aux adjonctions \eqref{tensor} données. En particulier, si réciproquement on suppose que $\underline{\mathrm{Hom}}_{\C}^{\text{\tiny{$\text{pointé}$}}}$ vérifie la propriété {\bf CMS$1^{\text{\tiny{$\text{pointé}$}}}$}, on définit les adjonctions \eqref{tensor} par les composés \eqref{tensorwedge}.

Finalement, étant donnée une catégorie pointée $\C$ munie d'un enrichissement $\underline{\mathrm{Hom}}_{\C}(\,\cdot\,,\,\cdot\,)$ dans la catégorie cartésienne fermée $\simp$, vérifiant la propriété {\bf CMS1}; montrons que $\underline{\mathrm{Hom}}_{\C}$ vérifie la propriété {\bf CMS2}, si et seulement si l'enrichissement $\underline{\mathrm{Hom}}_{\C}^{\text{\tiny{$\text{pointé}$}}}$ vérifie la propriété {\bf CMS$2^{\text{\tiny{$\text{pointé}$}}}$}. 

En effet, supposons que $\underline{\mathrm{Hom}}_{\C}$ vérifie {\bf CMS2}. Si $\xymatrix@C-4pt{A\ar[r]^j&B}$ est un monomorphisme d'ensembles simpliciaux pointés et $\xymatrix@C-8pt{X\ar[r]^q&Y}$ est une cofibration de $\C$, on doit montrer que le morphisme $\psi$ dans le diagramme cocartésien:
\begin{equation}
\xymatrix@R-10pt@C-10pt{
X\wedge A \ar[dd]_{X\wedge \,j}\ar[rr]^{q\,\wedge A} && 
Y\wedge A \ar[dd]\ar@/^15pt/[rddd]^{Y\wedge \,j}&\\\\
X\wedge B\ar[rr]\ar@/_15pt/[rrrd]_{q\,\wedge B}&&
X\wedge B\underset{X\wedge A}{\bigsqcup}Y\wedge A\ar@{-->}[rd]|(.65)\psi&\\
&&&Y\wedge B\,,}
\end{equation}
est une cofibration, laquelle est aussi une équivalence faible si $j$ est une $n$-équivalence faible d'ensembles simpliciaux pointés, ou si $q$ est une équivalence faible de $\C$.

Pour cela considérons le diagramme cocartésien de $\C$: 
\begin{equation}\label{uneaf}
\xymatrix@R-10pt@C-10pt{
X\otimes \pi A \ar[dd]_{X\otimes \,\pi j}\ar[rr]^{q\,\otimes \pi A} && 
Y\otimes \pi A \ar[dd]\ar@/^15pt/[rddd]^{Y\otimes \,\pi j}&\\\\
X\otimes \pi B\ar[rr]\ar@/_15pt/[rrrd]_{q\,\otimes \pi B}&&
X\otimes \pi B\underset{X\otimes \pi A}{\bigsqcup}Y\otimes \pi A\ar@{-->}[rd]|(.65)\varphi&\\
&&&Y\otimes \pi B\,,}
\end{equation}
associé à la cofibration $\xymatrix@C-8pt{X\ar[r]^q&Y}$ de $\C$ et au monomorphisme d'ensembles simpliciaux $\xymatrix@C-4pt{\pi A\ar[r]^{\pi j}&\pi B}$ sous-jacent à $j$. 

Vu que par hypothèse $\underline{\mathrm{Hom}}_{\C}$ vérifie la propriété {\bf CMS2}, il se suit que le morphisme $\varphi$ de \eqref{uneaf} est une cofibration de $\C$, laquelle est aussi une équivalence faible si $\pi j$ est une $n$-équivalence faible d'ensembles simpliciaux, ou si $q$ est une équivalence faible de $\C$.

D'un autre côté, remarquons que dans le diagramme:
\begin{equation}\label{avantdernier}
\xymatrix{
\star\ar[r] &Y\wedge A\ar@{}[rd]|-{(i)}\ar[r] \ar@/^20pt/[rr]^-{Y\wedge j} &(X\wedge B)\underset{(X\wedge A)}{\bigsqcup}(Y\wedge A)\ar[r]^-{\psi}\ar@{}[rd]|-{(ii)}& Y\wedge B\\
Y\otimes \star\ar[r]\ar[u]&Y\otimes \pi A\ar[r]\ar[u]\ar@/_25pt/[rr]_-{Y\otimes \pi j}  & (X\otimes \pi B)\underset{(X\otimes \pi A)}{\bigsqcup}(Y\otimes \pi A)\ar[u]\ar[r]_-{\varphi}& Y\otimes \pi B\ar[u]\,,}
\end{equation}
le carré $(i)$ est cartésien vu qu'il se décompose dans deux carré cocartésiens de $\C$:
\begin{equation}
\xymatrix@R=6pt{
Y\wedge A \ar[r]^-{\cong} & \star\underset{\star}{\bigsqcup} (Y\wedge A)\ar[r] &(X\wedge B)\underset{(X\wedge A)}{\bigsqcup} (Y\wedge A)\\
& \\&\\
Y\otimes \pi A\ar[uuu]\ar[r]_-{\cong} & (X\otimes\star)\underset{(X\otimes\star)}{\bigsqcup}(Y\otimes\pi A)\ar[r]\ar[uuu]& (X\otimes \pi B)\underset{(X\otimes \pi A)}{\bigsqcup}(Y\otimes \pi A)\,;\ar[uuu]}
\end{equation}
où le carré a droite est construit comme une colimite des diagrammes cocartésiens:
$$
\vcenter{\xymatrix{\star\ar[r] & X\wedge B \\X\otimes \star \ar[r]\ar[u] & X\otimes \pi B\ar[u]\,,}}
\qquad\,\,
\vcenter{\xymatrix{Y\wedge A\ar[r] & Y\wedge A \\Y\otimes \pi A \ar[r]\ar[u] & Y\otimes \pi A\ar[u]}}
\quad\,\text{et}\quad\,
\vcenter{\xymatrix{\star\ar[r] & X\wedge B \\X\otimes \star \ar[r]\ar[u] & X\otimes \pi B\ar[u]\,.}}
$$

Il se suit que le carré $(ii)$ du diagramme \eqref{avantdernier}:
\begin{equation}\label{dernier}
\xymatrix{
(X\wedge B)\underset{(X\wedge B)}{\bigsqcup}(Y\wedge A)\ar[r]^-{\psi}& Y\wedge B\\
(X\otimes \pi B)\underset{(X\otimes \pi B)}{\bigsqcup}(Y\otimes \pi A)\ar[u]\ar[r]_-{\varphi}& Y\otimes \pi B\,,\ar[u]}
\end{equation}
est un carré cocartésien de $\C$.

Vu que les familles de morphismes ${\bf cof}$ et ${\bf cof} \cap {\bf W}$, sont stables par cochangement de base dans toute catégorie de modèles; on conclut que si le morphisme $\varphi$ est une cofibration de $\C$ (resp. une cofibration et une équivalence faible), alors le morphisme $\psi$ est aussi une cofibration (resp. une cofibration et une équivalence faibles). 

En particulier, le morphisme $\psi$ est une cofibration laquelle est aussi une équivalence faible si $j$ est une $n$-équivalence faible d'ensembles simpliciaux pointés, ou si $q$ est une équivalence faible de $\C$. Donc, si l'enrichissement $\underline{\mathrm{Hom}}_{\C}$ vérifie la propriété {\bf CMS2}, alors $\underline{\mathrm{Hom}}_{\C}^{\text{\tiny{$\text{pointé}$}}}$ la vérifie la propriété {\bf CMS$2^{\text{\tiny{$\text{pointé}$}}}$}.

Réciproquement, si $\underline{\mathrm{Hom}}_{\C}^{\text{\tiny{$\text{pointé}$}}}$ satisfait la propriété {\bf CMS$2^{\text{\tiny{$\text{pointé}$}}}$}, pour montrer que $\underline{\mathrm{Hom}}_{\C}$ vérifie la propriété {\bf CMS2} remarquons que $(X\otimes A)_{+} \cong X \wedge (A_{+})$ $\big($le foncteur $X\wedge \,-$ commute aux coproduits$\big)$ et qu'un morphisme $\xymatrix@C-5pt{W\ar[r]^-{\phi}&Z}$ d'un catégorie de modèles pointée est une cofibration (resp. une cofibration et une équivalence faible) si $\xymatrix@C-5pt{W_+\ar[r]^-{\phi_+}&Z_+}$ l'est $\big($parce que $\xymatrix@C-5pt{W\ar[r]^-{\phi}&Z}$ est un rétracte du morphisme $\xymatrix@C-5pt{W_+\ar[r]^-{\phi_+}&Z_+}$ dans toute catégorie de modèles pointée$\big)$.
\end{proof}

Par exemple vu que $(\simp_{\star},\pi^{-1}{\bf W}_{n}, {\bf mono},\pi^{-1}{\bf fib}_{n},\wedge)$ est une catégorie de modèles monoïdale, $(\simp_{\star},\pi^{-1}{\bf W}_{n}, {\bf mono},\pi^{-1}{\bf fib}_{n},\underline{\mathrm{Hom}}_{\simp_{*}})$ est une catégorie de modèles simplicial d'ordre $n$ lorsque:
\begin{equation}\label{hompointe23}
\xymatrix@C+40pt{\simp_{\star} \times \simp_{\star} \ar[r]^-{\underline{\mathrm{Hom}}_{\simp_{\star}}(\,\cdot\,,\,\cdot\,)} & \simp}\,,
\end{equation}
$$
\text{est défini par}\qquad
\underline{\mathrm{Hom}}_{\simp_{*}}(X,Y)_{n} \,  = \, \pi\Big(\mathrm{hom}^{\wedge}_{\simp_{\star}}(X,Y)\Big)_{n} \, = \, \mathrm{Hom}_{\simp_{*}}(X \wedge \Delta^n_{+},Y)\,,
$$
où $\pi$ est le foncteur d'oubli de $\simp_{\star}$ vers $\simp$, et $\mathrm{hom}^{\wedge}_{\simp_{\star}}(\,\cdot\,,\,\cdot\,)$ est le foncteur de \eqref{hompointe2}.

Les foncteurs tenseur et cotenseur étant donné par:
$$
\simp_{\star}
\vcenter{\xymatrix@C+15pt{
\phantom{a}\ar@{}[r]|{\perp}\ar@<-4pt>@/_10pt/[r]_{\underline{\mathrm{Hom}}_{\simp_{\star}}(X,\,\cdot\,)}&\phantom{a}
\ar@<-4pt>@/_10pt/[l]_-{X\wedge(\,\cdot\,)_{+} }}}
\simp
\qquad\text{et}\qquad
\simp_{\star}^{op}
\vcenter{\xymatrix@C+15pt{
\phantom{a}\ar@{}[r]|{\perp}\ar@<-4pt>@/_10pt/[r]_{\underline{\mathrm{Hom}}_{\simp_{\star}}(\,\cdot\,,Y)}&\phantom{a}
\ar@<-4pt>@/_10pt/[l]_-{Y^{\,\cdot\,} }}}
\simp\,,
$$
\begin{align*}
\text{où}\quad 
\big(Y^{K}\big)_{n} \, & = \, \mathrm{Hom}_{\simp_{\star}}\big(K_{+}\wedge\Delta^n_{+},Y\big)\, = \,  \mathrm{Hom}_{\simp}\big(K\times \Delta^n,\pi(Y)\big).
\end{align*}

\renewcommand{\thesubsection}{\S\thesection.\arabic{subsection}}
\subsection{}\;
\renewcommand{\thesubsection}{\thesection.\arabic{subsection}}

Rappelons la construction des foncteurs suspension et espace des lacets qu'on a dans le cadre des catégories de modèles simpliciales pointées (voir \ref{stapre} ci-dessus et le Chapitre 6 de \cite{hovey}):

\begin{lemme}\label{classique}
Soit $0 \leq n \leq \infty$ et $\C$ une catégorie de modèles simpliciale pointée d'ordre $n$. Si $\mathbb{S}^1$ note le cercle simplicial vu comme un ensemble simplicial pointé (voir \ref{sphere}), alors l'adjonction:
\begin{equation}\label{adjcercle}
\C
\vcenter{\xymatrix@C+15pt{
\phantom{a}\ar@{}[r]|{\perp}\ar@<-4pt>@/_10pt/[r]_{(\,\cdot\,)^{\wedge\mathbb{S}^1}}&\phantom{a}
\ar@<-4pt>@/_10pt/[l]_-{(\,\cdot\,)\wedge \mathbb{S}^1}}}
\C\,,
\end{equation}
induit des adjonctions \eqref{lesdeuxp} est une adjonction de Quillen, et:
$$
\mathrm{Ho}\big(\,\C\,\big) 
\vcenter{\xymatrix@C+15pt{
\phantom{a}\ar@{}[r]|{\perp}\ar@<-4pt>@/_10pt/[r]_{\mathbf{R}\big((\,\cdot\,)^{\wedge\mathbb{S}^1}\big)}&\phantom{a}
\ar@<-4pt>@/_10pt/[l]_-{\mathbf{L}\big( \,\cdot\,\wedge \mathbb{S}^1\big)}}}
\mathrm{Ho}\big(\,\C\,\big)\,,
$$ 
est isomorphe à l'adjonction de la Proposition \ref{adjsigma}:
\begin{equation}\label{adjlemme}
\mathrm{Ho}(\,\C\,)\,
\vcenter{\xymatrix@C+15pt{
\phantom{a}\ar@{}[r]|{\perp}\ar@<-4pt>@/_10pt/[r]_{\Omega(\,\cdot\,)}&\phantom{a}
\ar@<-4pt>@/_10pt/[l]_-{\Sigma(\,\cdot\,)}}}
\mathrm{Ho}(\,\C\,)\,.
\end{equation}

En particulier, si $X$ et $Y$ sont des objets de $\C$ quelconques il y a des isomorphismes naturels:
\begin{equation}\label{homsigma}
\big[\Sigma^m(X),Y\big]_{\C} \; \cong \; \pi_{m} \big(\mathbf{R}\mathrm{Hom}_{\C}(X,Y),\star_{X}^Y\big) \; \cong \; \big[X,\Omega^m(Y)\big]_{\C}\qquad\text{pour}\;\; m\geq 0\,,
\end{equation}
où $\star_{X}^Y$ note le morphisme null de $X$ vers $Y$ dans $\C$\footnote{Dans \cite{hovey} on trouve une preuve de la formule \eqref{homsigma} pour les catégories de modèles pas forcement simpliciales.}.
\end{lemme}
\begin{proof}
Montrons que \eqref{adjcercle} est une adjonction de Quillen. Pour cela remarquons que d'après la propriété {\bf CMS$2^{\text{\tiny{$\text{pointé}$}}}$}  si on se donne une cofibration $f:\xymatrix@C-8pt{X\ar[r]&Y}$ de $\C$ (resp. une cofibration et une équivalence faible), vu que $\mathbb{S}^1$ est un objet cofibrant de la catégorie de modèles $(\simp_{\star},\pi^{-1}{\bf W}_{n}, {\bf mono},\pi^{-1}{\bf fib}_{n})$, le morphisme $f\wedge \mathbb{S}^1$ est une cofibration de $\C$ (resp. une cofibration et une équivalence faible). Autrement dit, le foncteur $\,\cdot\,\wedge \mathbb{S}^1$ est Quillen à gauche.  

D'un autre côté, remarquons que par définition de $\mathbb{S}^1$ on a un carré cocartésien de $\simp$:
\begin{equation}\label{s1simp}
\xymatrix{
\star \ar[r] &  \mathbb{S}^{1}\\ 
\partial\Delta^{1}\ar[r] \ar[u]  & \Delta^{1}\ar[u]\,.}
\end{equation}

Si on note $\big(\partial\Delta^{1},[0]\big)$ $\big($resp. $\big(\Delta^{1},[0]\big)$ $\big)$ l'ensemble simplicial pointé d'ensemble simplicial sous-jacent $\partial\Delta^1$ (resp. $\Delta^1$), dont le point de base est défini par le $0$-simplexe $[0]$; on vérifie facilement que le carré \eqref{s1simp} induit un carré cocartésien de $\simp_{\star}$:
$$
\xymatrix{
\star \ar[r] &  \mathbb{S}^{1}\\ 
\big(\partial\Delta^{1},[0]\big)\ar[r] \ar[u]  & \big(\Delta^{1},[0]\big)\ar[u]\,.}
$$

Donc, on a un carré cocartésien de $\C$:
\begin{equation}\label{cocacoca}
\xymatrix{
X\wedge\star \ar[r] &  X\wedge\mathbb{S}^{1}\\ 
X\wedge\big(\partial\Delta^{1},[0]\big)\ar[r] \ar[u]  & X\wedge\big(\Delta^{1},[0]\big)\ar[u]\,,}
\end{equation}
pour tout objet $X$ de $\C$.

Remarquons que d'après {\bf CMS$2^{\text{\tiny{$\text{pointé}$}}}$} et le Lemme \ref{fibrabra}, le carré \eqref{cocacoca} est homotopiquement cocartésien si $X$ est un objet cofibrant. En particulier, il se suit du Lemme \ref{isoomega} que pour montrer le résultat désiré, il suffit de noter qu'on a des isomorphismes de $\C$: 
$$
\xymatrix@C-8pt{X\wedge\star\ar[r]&\star}\,\;\qquad \text{et} \,\;\qquad \xymatrix@C-8pt{X\wedge\big(\Delta^{1},[0]\big)\ar[r]&X\,,}
$$
et une équivalences faible:
$$
\xymatrix@C-8pt{X\wedge\big(\partial\Delta^{1},[0]\big)\ar[r]&X\,.}
$$

En effet pour commencer remarquons que si $X$ est un objet quelconque de $\C$, le produit $X\wedge \star$ est un objet null de $\C$ parce que le foncteur $X\wedge\,\cdot\,$ commute aux colimites et $\star$ est un objet null de $\simp_{\star}$. Deuxièmement, vu que l'ensemble simplicial pointé $\big(\partial\Delta^{1},[0]\big)$ est isomorphe à $\Delta^0_{+}$, l'unité de la catégorie monoïdale $\simp_{\star}$  avec le produit wedge (voir \ref{pointesec}), il se suit que l'objet $X\wedge\big(\partial\Delta^{1},[0]\big)$ est canoniquement isomorphe à $X$ dans $\C$. On obtient ainsi deux isomorphismes $\xymatrix@C-8pt{X\wedge\star\ar[r]&\star}$ et $\xymatrix@C-8pt{X\wedge\big(\partial\Delta^{1},[0]\big)\ar[r]&X}$ de $\C$.

D'un autre côté, remarquons que le morphisme $\xymatrix@C-8pt{\Delta^0\ar[r]^{d_{1}}&\Delta^1}$ est une cofibration et une $n$-équivalence faible d'ensembles simpliciaux pour tout $0\leq n\leq \infty$; donc le morphisme qu'on en déduit pour $X$ cofibrant $\xymatrix@C-3pt{X\wedge\big(\Delta^{0},[0]\big)\ar[r]^{X\wedge d_{1}}&X\wedge\big(\Delta^{1},[0]\big)}$, est une équivalence faible de $\C$ d'après la Propriété  {\bf CMS$2^{\text{\tiny{$\text{pointé}$}}}$}. Donc, le seul morphisme $\xymatrix@C-8pt{X\wedge\big(\Delta^{1},[0]\big)\ar[r]&\star}$ de $\C$, induisant un isomorphisme dans la catégorie homotopique $\mathrm{Ho}(\C)$.


Finalement, si $X'$ est un remplacement cofibrant de $X$ et $Y'$ est un remplacement fibrant de $Y$, on a des isomorphismes pour $m\geq 0$:
\begin{equation*}
\begin{split}
\big[\Sigma^m(X),Y\big]_{\C} \; &\cong \; \pi_{0}\Big(\underline{\mathrm{Hom}}_{\C}\big(X'\wedge (\underbrace{\mathbb{S}^1\wedge \dots \wedge \mathbb{S}^1}_{m}),Y'\big)\Big) \\
&\cong \; \pi_{0}\Big(\underline{\mathrm{Hom}}_{\C}\big(X'\wedge\mathbb{S}^m,Y'\big)\Big) \\
&\cong \; \pi_{0}\Big(\underline{\mathrm{Hom}}_{\simp_{*}}\big(\mathbb{S}^m,\underline{\mathrm{Hom}}^{\text{\tiny{$\text{pointé}$}}}_{\C}(X',Y')\big)\Big) \\
&\cong \; \pi_{m} \big(\mathbf{R}\mathrm{Hom}_{\C}(X,Y),\star_{X}^Y\big) 
\end{split}
\end{equation*}
\end{proof}

Montrons  :

\begin{theoreme}\label{cassi1}
Soit $0\leq n < \infty$ et $(\C,{\bf W}, {\bf cof},{\bf fib},\underline{\mathrm{Hom}}_{\C})$ une catégorie de modèles simpliciale pointée d'ordre $\infty$ telle que tous ses objets soient cofibrants; alors les énoncés suivants sont équi\-va\-lents:
\begin{enumerate}
\item $(\C,{\bf W}, {\bf cof},{\bf fib},\underline{\mathrm{Hom}}_{\C})$ est une catégorie de modèles simpliciale pointée d'ordre $n$ (\emph{i.e.} la structure de catégorie de modèles simpliciale d'ordre $\infty$ donnée, vérifie la propriété {\bf CMS2} pour les $n$-équivalences faibles).
\item Le foncteur $\Sigma^m\colon \xymatrix@C-10pt{\mathrm{Ho}(\C)\ar[r]&\mathrm{Ho}(\C)}$ est null pour tout $m>n$.
\item Le foncteur $\Omega^m\colon \xymatrix@C-10pt{\mathrm{Ho}(\C)\ar[r]&\mathrm{Ho}(\C)}$ est null pour tout $m>n$.
\item Si $X$ et $Y$ sont des objets de $\C$ quelconques, $\pi_{m} \big(\mathbf{R}\mathrm{Hom}_{\C}(X,Y),\star_{X}^Y\big)=0$ pour $m>n$.
\item Si $Y$ est un objet fibrant de $(\C,{\bf W}, {\bf cof},{\bf fib})$, l'ensemble simplicial $\underline{\mathrm{Hom}}_{\C}(X,Y)$ est un objet fibrant de la catégorie de modèles $(\simp,{\bf W}_{n}, {\bf mono},{\bf fib}_{n})$ (du Théorème \ref{ntypess}) pour tout objet $X$ de $\C$.
\end{enumerate} 

En fait les conditions (ii), (iii) et (iv) sont équivalentes pour toutes les catégories de modèles pointées.
\end{theoreme}
\begin{proof}
On n'a pas besoin de supposer que tous les objets de $\C$ soient cofibrants pour montrer (i) $\Rightarrow$ (ii). En effet remarquons que si $\C$ est une catégorie de modèles simpliciale pointée d'ordre $n$ quelconque, vu que le morphismes canonique $\xymatrix@C-8pt{\star\ar[r]& \mathbb{S}^{m+1}}$ est une $n$-équivalence faible d'ensembles simpliciaux pointés pour $n<m$, il se suit de la propriété {\bf CMS$2^{\text{\tiny{$\text{pointé}$}}}$} que le morphisme $\xymatrix{\star\cong (X\wedge\star)\ar[r] & X\wedge\mathbb{S}^m}$ est une équivalence faible de $\C$ pour tout $n<m$ et tout objet $X$ de $\C$ cofibrant. Donc, du Lemme \ref{classique} ci-dessus, on déduit que le foncteur $\Sigma^m$ est null pour tout $m>n$.

Pour montrer (iv)$\Leftrightarrow$(v) remarquons que si $Y$ est un objet fibrant d'une catégorie de modèles simpliciale (d'ordre $\infty$) dont tous les objets sont cofibrants $(\C,{\bf W}, {\bf cof},{\bf fib})$, alors pour tout objet $X$ de $\C$ l'ensemble simplicial $\underline{\mathrm{Hom}}_{\C}(X,Y)$ est un complexe de Kan dont le type d'homotopie est le complexe de fonction dérivé ${\bf R}\mathrm{Hom}_{\C}(X,Y)$ (voir le Lemme \ref{RHOMloc}). Donc $\underline{\mathrm{Hom}}_{\C}(X,Y)$ est un objet fibrant de la catégorie de modèles $(\simp,{\bf W}_{n}, {\bf mono},{\bf fib}_{n})$ du Théorème \ref{ntypess}, si et seulement si $\pi_{m} \big(\mathbf{R}\mathrm{Hom}_{\C}(X,Y),\star_{X}^Y\big)=0$ pour $m>n$.

D'un autre côté, on déduit de l'adjonction \eqref{adjlemme} que (ii) et (iii) sont des énoncés équivalents. En plus (iii)$\Leftrightarrow$(iv) parce que d'après le Lemme \ref{classique} il y a des isomorphismes naturels:
$$\big[X,\Omega^m Y\big]_{\C} \; \cong \; \pi_{m} \big(\mathbf{R}\mathrm{Hom}_{\C}(X,Y),\star_{X}^Y\big)\,,$$ 
quelques soient $X$ et $Y$. Il se suit en particulier que (ii)$\Leftrightarrow$(iii)$\Leftrightarrow$(iv) pour toues les catégories de modèles pointées (voir la Proposition \ref{adjsigma} et le Lemme 6.1.2 de \cite{hovey}). 

Supposons maintenant que $\C$ est une catégorie de modèles simpliciale pointé d'ordre $\infty$ dont les objets soient tous cofibrants et laquelle vérifie la propriété (iv).

Pour montrer que $\C$ vérifie la propriété (i) on doit montrer que si $\xymatrix@C-8pt{X\ar[r]^q&Y}$ est une cofibration de $\C$ et $\xymatrix@C-4pt{A\ar[r]^j&B}$ est un monomorphisme et une $n$-équivalence faible d'ensembles simpliciaux, la cofibration $\varphi$ dans le diagramme somme amalgamée suivant:
\begin{equation}\label{cocarproof}
\xymatrix@R-10pt@C-10pt{
X\otimes A \ar[dd]_{X\otimes \,j}\ar[rr]^{q\,\otimes A} && 
Y\otimes A \ar[dd]\ar@/^15pt/[rddd]^{Y\otimes \,j}&\\\\
X\otimes B\ar[rr]\ar@/_15pt/[rrrd]_{q\,\otimes B}&&
X\otimes B\underset{X\otimes A}{\bigsqcup}Y\otimes A\ar@{-->}[rd]|(.65)\varphi&\\
&&&Y\otimes B\,,}
\end{equation}
est une équivalence faible de $\C$.

Montrons premièrement qu'on peut supposer que $A$ et $B$ sont de complexes de Kan. Pour cela remarquons en premier lieu que dans le diagramme \eqref{cocarproof}, le morphisme $q\otimes A$ est une cofibration de $\C$ d'après {\bf CMS2} avec $j$ la cofibration canonique $\xymatrix@C-10pt{\emptyset\ar[r]&A}$ de $\simp$ (note que $X\otimes\emptyset$ et $Y\otimes\emptyset$ sont des objets initiaux de $\C$); en particulier, vu que tous les objets de $\C$ sont cofibrants, le carré cocartésien dans \eqref{cocarproof} est en fait un carré homotopiquement cocartésien d'après le Lemme \ref{fibrabra}. 

Si on considère maintenant un carré commutatif d'ensembles simpliciaux:
$$
\xymatrix{A\ar[r]^-{j}\ar[d]_-{q_{A}} & B \ar[d]^-{q_{B}} \\ A' \ar[r]_-{j'} &B'\,,}
$$
où les morphismes verticaux sont des $\infty$-équivalences faibles et des monomorphismes de but des complexes de Kan; on vérifie d'un côté que les morphismes: 
$$
\xymatrix@+5pt{X\otimes A\ar[r]^-{X\otimes q_{A}} & X\otimes A'}\,,\qquad
\xymatrix@+5pt{Y\otimes A\ar[r]^-{Y\otimes q_{A}} & Y\otimes A'}\,,
$$
$$
\xymatrix@+5pt{X\otimes B\ar[r]^-{X\otimes q_{B}} & X\otimes B'} \qquad \text{et}\qquad \xymatrix@+5pt{Y\otimes B \ar[r]^-{Y\otimes q_{B}} & Y \otimes B'\,,}
$$
sont des équivalences faibles de $\C$ (d'après {\bf CMS2}), et d'un autre côté que dans le diagramme:
\begin{equation}\label{cocarproof1}
\xymatrix@R-10pt@C-10pt{
X\otimes A' \ar[dd]_{X\otimes \,j'}\ar[rr]^{q\,\otimes A'} && 
Y\otimes A' \ar[dd]\ar@/^15pt/[rddd]^{Y\otimes \,j'}&\\\\
X\otimes B'\ar[rr]\ar@/_15pt/[rrrd]_{q\,\otimes B'}&&
X\otimes B'\underset{X\otimes A'}{\bigsqcup}Y\otimes A'\ar@{-->}[rd]|(.65){\varphi'}&\\
&&&Y\otimes B'\,,}
\end{equation}
le carré cocartésien est homotopiquement cocartésien (d'après le Lemme \ref{fibrabra} vu que $q\,\otimes A'$ est une cofibration).

On déduit du Lemme \ref{cubelele} que le morphisme $\varphi$ de \eqref{cocarproof} est une équivalence faible de $\C$, si et seulement si le morphisme $\varphi'$ de \eqref{cocarproof1} est une équivalence faible de $\C$. Autrement dit, pour montrer que dans le diagramme \eqref{cocarproof} le morphisme $\varphi$ de $\C$ soit une équivalence faible, on peut supposer que $j$ est un morphisme entre complexes de Kan. 

Deuxièmement, il se suit du fait que $A$ et $B$ soient de complexes de Kan, que le morphisme d'ensembles simpliciaux $\xymatrix@C+5pt{{\bf csq}_{n+1}(A)\ar[r]^-{{\bf csq}_{n+1}(j)}&{\bf csq}_{n+1}(B)}$ est un monomorphisme et une $\infty$-équivalence faible, parce que $j$ est un monomorphisme et une $n$-équivalence faible (voir le Corollaire \ref{lecoro}); donc, le morphisme $\psi$ dans le diagramme:
\begin{equation}\label{cocarproof2}
\xymatrix@R-10pt@C-10pt{
X\otimes {\bf csq}_{n+1}A \ar[dd]_{X\otimes \,{\bf csq}_{n+1} j}\ar[rr]^{q\,\otimes {\bf csq}_{n+1} A} && 
Y\otimes {\bf csq}_{n+1} A \ar[dd]\ar@/^15pt/[rddd]^{Y\otimes \,{\bf csq}_{n+1}j}&\\\\
X\otimes {\bf csq}_{n+1}B\ar[rr]\ar@/_15pt/[rrrd]_{q\,\otimes B}&&
X\otimes {\bf csq}_{n+1}B\underset{X\otimes {\bf csq}_{n+1} A}{\bigsqcup}Y\otimes {\bf csq}_{n+1} A\ar@{-->}[rd]|(.65)\psi&\\
&&&Y\otimes {\bf csq}_{n+1} B\,,}
\end{equation}
est bien une équivalence faible de $\C$.

Enfin, vu que le carré cocartésien de \eqref{cocarproof2} est un carré homotopiquement cocartésien, pour montrer que le morphisme $\varphi$ de \eqref{cocarproof} soit une équivalence faible de $\C$, d'après le Lemme \ref{cubelele} il suffit de démontrer que le morphisme $\xymatrix@C-3pt{Z\otimes K \ar[r]^-{Z\otimes \eta_{K}} & Z\otimes {\bf csq}_{n+1} K }$ induit du morphisme $\xymatrix@C-6pt{K \ar[r]^-{\eta_{K}} &  {\bf csq}_{n+1}K}$ où $\eta:\xymatrix@C-9pt{\mathrm{id}\ar@{=>}[r]&{\bf csq}_{n+1}}$ est une unité donnée, est une équivalence faible de $\C$ quelques soient $Z$ et $K$.

De fa\c con équivalente, il suffit de montrer que pour tout objet fibrant $W$ de $\C$, la fonction:
\begin{equation}\label{ilfaut}
\xymatrix@R=3pt@C+10pt{
[Z\otimes {\bf csq}_{n+1}K,W]_\C
\ar@{}[d]|-{\mathrel{\reflectbox{\rotatebox[origin=c]{90}{$\cong$}}}}
\ar[r] &  [Z\otimes K,W]_\C
\ar@{}[d]|-{\mathrel{\reflectbox{\rotatebox[origin=c]{90}{$\cong$}}}}\\
\pi_{0}\big( \underline{\mathrm{Hom}}_{\C}(Z\otimes {\bf csq}_{n+1}K,W)\big)\ar[r]_-{(Z\otimes\eta_{K})^\star} & 
\pi_{0}\big( \underline{\mathrm{Hom}}_{\C}(Z\otimes K,W)\big)\,,}
\end{equation}
induite du morphisme $\xymatrix@C-3pt{Z\otimes K \ar[r]^-{Z\otimes \eta_{K}} & Z\otimes {\bf csq}_{n+1} K }$, est bijective (Note que dans les isomorphismes de \eqref{ilfaut} on utilise à nouveau que tout objet de $\C$ est cofibrant).

Pour vérifier cette affirmation, remarquons que si $W$ est un objet fibrant de $\C$, alors $\underline{\mathrm{Hom}}_{\C}(Z,W)$ est un complexe de Kan pour tout objet $Z$ de $\C$ d'après la propriété duale à {\bf CMS2} (voir le Lemme 4.2.2 de \cite{hovey}); donc, vu qu'on a supposé que l'énoncé (iv) est vrai, le morphisme: 
\begin{equation}\label{loqbesoi}
\xymatrix{\underline{\mathrm{Hom}}_{\C}(Z,W) \ar[r]& {\bf csq}_{n+1}\underline{\mathrm{Hom}}_{\C}(Z,W)}
\end{equation}
est une $\infty$-équivalence faible entre complexes de Kan (voir le Corollaire \ref{lecoro}).

En particulier, pour tout ensemble simplicial $K$ on a un diagramme commutatif:
$$
\xymatrix@R=3pt{
\pi_{0}\big( \underline{\mathrm{Hom}}_{\C}(Z\otimes {\bf csq}_{n+1}K,W)\big)\ar[r]^-{(Z\otimes\eta_{K})^\star} 
\ar@{}[d]|-{\mathrel{\reflectbox{\rotatebox[origin=c]{90}{$\cong$}}}}& 
\pi_{0}\big( \underline{\mathrm{Hom}}_{\C}(Z\otimes K,W)\big)
\ar@{}[d]|-{\mathrel{\reflectbox{\rotatebox[origin=c]{90}{$\cong$}}}}\\
\pi_{0}\Big( \underline{\mathrm{Hom}}_{\simp}\big({\bf csq}_{n+1}K,\underline{\mathrm{Hom}}_{\C}(Z,W)\big)\Big) 
\ar@{}[d]|-{\mathrel{\reflectbox{\rotatebox[origin=c]{90}{$\cong$}}}}_{\eqref{loqbesoi}}& 
\pi_{0}\Big( \underline{\mathrm{Hom}}_{\simp}\big(K,\underline{\mathrm{Hom}}_{\C}(Z,W)\big)\Big)
\ar@{}[d]|-{\mathrel{\reflectbox{\rotatebox[origin=c]{90}{$\cong$}}}}^{\eqref{loqbesoi}}\\
\pi_{0}\Big( \underline{\mathrm{Hom}}_{\simp}\big({\bf csq}_{n+1}K,{\bf csq}_{n+1}\underline{\mathrm{Hom}}_{\C}(Z,W)\big)\Big)\ar[r]_-{\eta_{K}^{\;\star}}&\pi_{0}\Big( \underline{\mathrm{Hom}}_{\simp}\big(K,{\bf csq}_{n+1}\underline{\mathrm{Hom}}_{\C}(Z,W)\big)\Big)\,,}
$$
où la fonction horizontale inférieur est une bijection d'après le Lemme \ref{lemehomo}. Donc, \eqref{ilfaut} est une fonction bijective.
\end{proof}

\renewcommand{\thesubsection}{\S\thesection.\arabic{subsection}}
\subsection{}\;
\renewcommand{\thesubsection}{\thesection.\arabic{subsection}} 

Dans le présent paragraphe on va déduire du Théorème \ref{cassi1} et d'un résultat de D. Dugger (voir \cite{dugger}, \cite{rezk} et la Proposition \ref{duggerT} ci-dessous), l'énoncé suivant:

\begin{corollaire}\label{dugger2}
Soit $0\leq n\leq \infty$. Si $\C$ est une catégorie de modèles de Cisinski\footnote{Un \emph{catégorie de modèles de Cisinski} est une catégorie de modèles à engendrement cofibrant dont la catégorie sous-jacente est une catégorie de préfaisceaux d'ensembles et les cofibrations sont les monomorphismes.} alors les énoncés qui suivent sont équivalents:
\begin{enumerate}
\item Il existe une équivalence de Quillen:
$$
\C
\vcenter{\xymatrix@C+15pt{
\phantom{a}\ar@<-4pt>@/_10pt/[r]_-{G}\ar@{}[r]|{\perp}&\phantom{a}
\ar@<-4pt>@/_10pt/[l]_-{F}}}
\D
$$
où $\D$ est une catégorie de modèles simpliciale d'ordre $n$ et $F$ respect l'objet final $F(\star_\D)\cong\star_\C$.
\item $\C$ admet à équivalence de Quillen pointée près une structure de catégorie de modèles simpliciale d'ordre $n$, c'est-à-dire, il existe une chaîne d'équivalences de Quillen:
$$
\C
\vcenter{\xymatrix@C+15pt{
\phantom{a}\ar@<-4pt>@{-}@/_10pt/[r]_{\phantom{a}}&\phantom{a}
\ar@<-4pt>@{-}@/_10pt/[l]_{\phantom{a}}}}
\D_1
\vcenter{\xymatrix@C+15pt{
\phantom{a}\ar@<-4pt>@{-}@/_10pt/[r]&\phantom{a}
\ar@<-4pt>@{-}@/_10pt/[l]}}
\D_2
\vcenter{\xymatrix@C+15pt{
\phantom{a}\ar@{}[r]|-{\cdots}&\phantom{a}}}
\D_{i}
\vcenter{\xymatrix@C+15pt{
\phantom{a}\ar@<-4pt>@{-}@/_10pt/[r]&\phantom{a}
\ar@<-4pt>@{-}@/_10pt/[l]}}
\D
$$
où $\D$ est une catégorie de modèles simpliciale d'ordre $n$ et les foncteurs Quillen à gauche respectent l'objet final.
\item Le foncteur $\Sigma^m\colon \xymatrix@C-10pt{\mathrm{Ho}(\C_{\star})\ar[r]&\mathrm{Ho}(\C_{\star})}$ est null pour tout $m>n$, où $\C_\star$ note la catégorie des objets pointés de $\C$.
\item Le foncteur $\Omega^m\colon \xymatrix@C-10pt{\mathrm{Ho}(\C_{\star})\ar[r]&\mathrm{Ho}(\C_{\star})}$ est null pour tout $m>n$.
\item Si $X$ et $Y$ sont des objets de $\C_{\star}$ quelconques, $\pi_{m} \big(\mathbf{R}\mathrm{Hom}_{\C_{\star}}(X,Y),\star_{X}^Y\big)=0$ pour $m>n$.
\end{enumerate} 
\end{corollaire}

La version du résultat de D. Dugger dont ont aura besoin est la suivante:

\begin{proposition}\label{duggerT}
Si $\C$ est une catégorie de modèles de Cisinski, alors il existe une équivalence de Quillen:
$$
\C
\vcenter{\xymatrix@C+15pt{
\phantom{a}\ar@{}[r]|{\perp}\ar@<+4pt>@/^10pt/[r]^{F}&\phantom{a}
\ar@<+4pt>@/^10pt/[l]^-{G}}}
\D\,,
$$
vérifiant les propriétés:
\begin{enumerate}
\item $\D$ est une catégorie de modèles simpliciale d'ordre $\infty$.
\item Dans $\D_\star$ la catégorie de modèles des objets pointés de $\D$ (voir le Lemme \ref{lem1}) tous les objets sont cofibrants.
\item $F(\star_{\C})\cong\star_{\D}$, c'est-à-dire, l'image par $F$ d'un objet final de $\C$ est un objet final de $\D$.
\end{enumerate}
\end{proposition}
\begin{proof}
Rappelons qu'une catégorie de modèle de Cisinski est toujours propre à gauche (les cofibrations sont les monomorphismes) et combinatoire (à engendrement cofibrant et localement présentable); donc d'après \cite{dugger} (voir le Théorème 1.2 de \cite{dugger}) l'adjonction:
\begin{equation}\label{adjno}
\C
\vcenter{\xymatrix@C+15pt{
\phantom{a}\ar@{}[r]|{\perp}\ar@<+4pt>@/^10pt/[r]^{F \, = \, \text{préfaisceau constant}}&\phantom{a}
\ar@<+4pt>@/^10pt/[l]^-{G \, = \, \text{evaluation en $0$}}}}
\D \, = \, \C^{\Delta^{op}}\,,
\end{equation}
est une équivalence de Quillen si la catégorie de foncteurs $\D=\C^{\Delta^{op}}$ est munie de la structure de catégorie de modèles dont les cofibrations sont les cofibrations de Reedy (voir la Définition 15.3.3 de \cite{hirschhorn}) et les équivalences faibles sont les hocolim-équivalences (voir la Définition 5.1 de \cite{dugger}). En fait dans le Théorème 6.1 de \cite{dugger} il est montré que l'enrichissement canonique de $\C^{\Delta^{op}}$ dans $\simp$ (voir \S4.1 de \cite{dugger}), fait de $\D = \C^{\Delta^{op}}$ une catégorie de modèles simpliciale d'ordre $\infty$. 

Rappelons d'un autre côté qu'un objet $(X,x_0)$ de la catégorie des objets pointés $\D_\star = \big(\C^{\Delta^{op}}\big)_\star$ est cofibrant si et seulement si le morphisme $x_0\colon\xymatrix@C-7pt{\star\ar[r]&X}$ est une cofibration de Reedy de $\D = \C^{\Delta^{op}}$. Vu que les monomorphismes de la catégorie de préfaisceaux $\D=\C^{\Delta^{op}}$ sont exactement les cofibrations de Reedy, il se suit que tous les objets des catégories de modèles $\D$ et $\D_\star$ sont cofibrants.
 


Pour vérifier finalement la propriété (iii) remarquons que le foncteur $F$ de la adjonction \eqref{adjno} envoie un objet $X$ de $\C$ vers l'objet simplicial constant à valeurs $X$, donc $F$ envoie tout objet final de $\C$ dans un objet final de $\D$.
\end{proof}

Pour montrer le Corollaire \ref{dugger2} commençons par nous rappeler que si $\C$ est une catégorie munie d'un objet final $\star$, la catégorie $\C_{\star}$ des \emph{objets pointés de $\C$} est la catégorie des couples $X=(X,x)$ où $X$ est un objet de $\C$ et $x:\xymatrix@C-8pt{\star\ar[r]&X}$ est une flèche. Un morphisme de $(X,x)$ vers $(Y,y)$ est une flèche $f:\xymatrix@-6pt{X\ar[r]&Y}$ de $\C$ telle que $f x=y$.

Si on suppose de plus que $\C$ soit une catégorie admettant de coproduits finis, on montre facilement que le foncteur d'oubli $\pi :\xymatrix@-8pt{\C_{\star}\ar[r]&\C}$ défini par la règle $(X,x)\mapsto X$, admet un adjoint à  gauche $(\;\cdot\;)_{+}: \xymatrix@-5pt{\C\ar[r]&\C_{\star}}$, défini dans un objet $A$ de $\C$ comme le coproduit $A_{+}=A\sqcup\star$, pointé par le morphisme canonique $\xymatrix@C-8pt{\star\ar[r]&A\sqcup\star}$. 

L'énoncé suivant est la Proposition 1.1.8 de \cite{hovey}: 

\begin{lemme}\label{lem1}
Si $\C$ est une catégorie de modèles, alors $\C_{\star}$ admet une structure de catégorie de modèles, lorsque les équivalences faibles (resp. les cofibrations, les fibrations) sont les morphismes $f:\xymatrix@-14pt{(X,x)\ar[r]&(Y,y)}$ dont le morphisme sous-jacent $f:\xymatrix@-14pt{X\ar[r]&Y}$ est une équivalence faible (resp. une cofibration, une fibration) de $\C$. En particulier, 
\begin{equation}\label{olvidar}
\C_{\star}
\vcenter{\xymatrix@C+15pt{
\phantom{a}\ar@{}[r]|{\perp}\ar@<-4pt>@/_10pt/[r]_{\pi}&\phantom{a}
\ar@<-4pt>@/_10pt/[l]_-{(\,\cdot\,)_{+}}}}
\C\,,
\end{equation}
est une adjonction de Quillen.  
\end{lemme}

Remarquons maintenant que si $\C$ et $\D$ sont deux catégories munies des objets finals $\star_{\C}$ et $\star_{\D}$ respectivement, étant donnée une adjonction:
$$
\C
\vcenter{\xymatrix@C+15pt{
\phantom{a}\ar@{}[r]|{\perp}\ar@<+4pt>@/^10pt/[r]^{F}&\phantom{a}
\ar@<+4pt>@/^10pt/[l]^-{G}}}
\D\,,
$$
avec la propriété que le seul morphisme $\xymatrix@C-8pt{F(\star_{\C})\ar[r]^-{\Phi}_-{\cong}&\star_{\D}}$ soit un isomorphisme, on en déduit une adjonction:
$$
\C_{\star}
\vcenter{\xymatrix@C+15pt{
\phantom{a}\ar@{}[r]|{\perp}\ar@<+4pt>@/^10pt/[r]^{F_{\star}}&\phantom{a}
\ar@<+4pt>@/^10pt/[l]^-{G_{\star}}}}
\D_{\star}\,,
$$
où les foncteurs $F_\star$ et $G_\star$ sont définis par les formules: 
$$
F_{\star}(X,x_{0})=\big(FX,F(x_{0})\circ\Phi^{-1}\big) \qquad\text{et}\qquad  G_{\star}(A,a_{0})=\big(GA,G(a_{0})\circ\Psi^{-1}\big)\,,
$$
$\xymatrix@C-5pt{G(\star_{\D})\ar[r]^-{\Psi}_-{\cong}&\star_{\C}}$ étant le seul morphisme de but $\star_\C$ dans $\C$.

Pour montrer qu'on obtient effectivement une adjonction $F_\star \dashv G_\star$, considérons $\eta$ et $\epsilon$ une unité et une counité de l'adjonction $F \dashv G$ vérifiant les identités triangulaires. On constate aussi-tôt que les transformations naturelles $\eta_\star\colon\xymatrix@C-15pt{\mathrm{id}\ar@{=>}[r]&G_\star\circ F_\star}$ et $\epsilon_\star\colon\xymatrix@C-15pt{F_\star\circ G_\star\ar@{=>}[r]&\mathrm{id}}$ définies par les règles:
 \begin{align*}
\xymatrix@C+28pt{(X,x_0)\ar[r]^-{(\eta_\star)_{(X,x_0)}\,=\,\eta_X}&G_\star\circ F_\star(X,x_0)=\big(G\circ F(X), G\circ F (x_0)\circ G(\Phi)^{-1}\circ \Psi^{-1}\big)}\\
\text{et}\hspace{7cm}\\
\xymatrix@C+28pt{F_\star\circ G_\star(A,a_0)=\big(F\circ G(A), F\circ G (a_0)\circ F(\Psi)^{-1}\circ \Phi^{-1}\big)\ar[r]^-{(\epsilon_\star)_{(A,a_0)}\,=\,\epsilon_A}&(A,a_0)}\,,
 \end{align*}
lesquelles sont bien définies parce qu'on a des diagrammes commutatifs:
$$
\vcenter{\xymatrix@C-15pt@R-10pt{
X\ar[rr]^-{\eta_X} && GF(A)\\
\star_\C \ar[u]^-{x_0} \ar[rr]|-{\eta_{\star_\C}} &&GF(\star_\C) \ar[u]_-{GF(x_0)} \ar[dl]^-{G(\Phi)}\\
&G(\star_\D)\ar[ul]^-{\Psi}&}}
\qquad \text{et}\qquad 
\vcenter{\xymatrix@C-15pt@R-10pt{
FG(A)\ar[rr]^-{\epsilon_A} && A\\
FG(\star_D) \ar[u]^-{FG(a_0)} \ar[rr]|-{\epsilon_{\star_\D}} \ar[rd]_-{F(\Psi)} &&\star_\D \ar[u]_-{a_0} \\
&F(\star_\C)\ar[ur]^-{\Phi}&}}
$$ 
vérifient les identités triangulaires.

Montrons: 

\begin{lemme}\label{lem2}
Soit:
\begin{equation}\label{adjnini}
\C
\vcenter{\xymatrix@C+15pt{
\phantom{a}\ar@{}[r]|{\perp}\ar@<+4pt>@/^10pt/[r]^{F}&\phantom{a}
\ar@<+4pt>@/^10pt/[l]^-{G}}}
\D\,,
\end{equation}
une adjonction entre catégories de modèles vérifiant que le seul morphisme $\xymatrix@C-8pt{F(\star_{\C})\ar[r]^-{\Phi}_-{\cong}&\star_{\D}}$ soit un isomorphisme; si \eqref{adjnini} est une adjonction de Quillen (resp. une équivalence de Quillen) alors la adjonction induite:
\begin{equation}\label{adjninip}
\C_{\star}
\vcenter{\xymatrix@C+15pt{
\phantom{a}\ar@{}[r]|{\perp}\ar@<+4pt>@/^10pt/[r]^{F_{\star}}&\phantom{a}
\ar@<+4pt>@/^10pt/[l]^-{G_{\star}}}}
\D_{\star}\,,
\end{equation}
est une adjonction de Quillen (resp. une équivalence de Quillen) où $\C_\star$ et $\D_\star$ sont les catégories de modèles du Lemme \ref{lem1}.
\end{lemme}
\begin{proof}
Si $F$ est Quillen à gauche, on déduit que $F_\star$ l'est aussi de la définition de la structure de catégorie de modèles dans la catégorie des objets pointés d'une catégorie de modèles (voir le Lemme \ref{lem1}).

D'un autre, rappelons que un objet $(A,a_{0})$ est fibrant dans $\D_{\star}$ (resp. $(X,x_0)$ est cofibrant dans $\C_\star$) si et seulement si $A$ est un objet fibrant de $\D$ (resp.$\xymatrix@C-8pt{\star_{\C}\ar[r]^-{x_0}&X}$ est une cofibration de $\C$). Plus encore, on sais que:
$$
{\bf L}(F_{\star}) (X,x_0) = F_\star(LX,x'_0) \qquad \text{et}\qquad  {\bf R}(G_{\star}) (A,a_{0}) = G_\star \big( RA, p_{A}\circ a_{0}) 
$$
sont des foncteurs dérivés totals à gauche et à droite de $F_\star$ et $G_\star$ respectivement; où 
$$
\xymatrix@C-6pt{\star_\C\,\ar@{>->}[r]^-{x'_0} & LX\ar@{->>}[r]_-{\text{\rotatebox[origin=c]{360}{\Large $\widetilde{\phantom{w}}$}}}^-{q_{X}}& X}\,,\qquad\text{et}\qquad
\xymatrix@C-6pt{A\,\ar@{>->}[r]_-{\text{\rotatebox[origin=c]{360}{\Large $\widetilde{\phantom{w}}$}}}^-{p_{A}}& RA\ar@{->>}[r] & \star_{\D}}\,,
$$ 
sont des factorisations fonctorielles des morphismes $\xymatrix@C-8pt{\star_{\C}\ar[r]^-{x_0}&X}$ et $\xymatrix@C-8pt{A\ar[r]&\star_{\D}}$ respectivement.  

Dans ce cas une counité de l'adjonction ${\bf L}(F_{\star}) \dashv {\bf R}(G_{\star})$ dans l'objet $(A,a_{0})$:
$$
\xymatrix@C+5pt{{\bf L}(F_{\star})\,{\bf R}(G_{\star}) (A,a_{0})\ar[r]^-{\epsilon_{(A,a_0)}} & (A,a_{0})}
$$
est un morphisme de la catégorie homotopique $\mathrm{Ho}(\,\D_\star\,)$ dont le morphisme sous-jacent dans $\mathrm{Ho}(\,\D\,)$ est donné par la chaîne:
$$
\xymatrix@C+10pt{F \Big(L \big(G (RA)\big)\Big) \ar[r]^--{F(q_{G(RA)})} &F \big(G (RA)\big) \ar[r]^-{\epsilon_{RA}}  & RA&A
\ar[l]^-{\text{\rotatebox[origin=c]{360}{\Large $\widetilde{\phantom{w}}$}}}_-{p_A}}\,;
$$
où $\epsilon$ est une counité de l'adjonction $F \dashv G$. En particulier, si \eqref{adjnini} est une équivalence de Quillen alors \eqref{adjninip} l'est aussi.
\end{proof}

Rappelons maintenant que si $\C$ est une catégorie complète et cocomplète munie d'un enrichissement dans la catégorie cartésienne fermée $\simp$:
$$
\xymatrix@C+35pt{\C^{op} \times \C \ar[r]^-{\underline{\mathrm{Hom}}_{\C}(\,\cdot\,,\,\cdot\,)} & \simp},
$$
et munie des foncteurs adjoints:
\begin{equation}\label{lesdeux2}  
\C
\vcenter{\xymatrix@C+15pt{
\phantom{a}\ar@{}[r]|{\perp}\ar@<-4pt>@/_10pt/[r]_{\underline{\mathrm{Hom}}_{\C}(X,\,\cdot\,)}&\phantom{a}
\ar@<-4pt>@/_10pt/[l]_-{X\otimes\,\cdot\, }}}
\simp
\qquad\text{et}\qquad
\C^{op}
\vcenter{\xymatrix@C+15pt{
\phantom{a}\ar@{}[r]|{\perp}\ar@<-4pt>@/_10pt/[r]_{\underline{\mathrm{Hom}}_{\C}(\,\cdot\,,Y)}&\phantom{a}
\ar@<-4pt>@/_10pt/[l]_-{Y^{\,\cdot\,} }}}
\simp\,,
\end{equation}
pour $X$ et $Y$ des objets quelconques de $\C$; alors la catégorie $\C_{\star}$ des objets pointés de $\C$ admet un enrichissement canonique:
\begin{equation} \label{morstar}
\xymatrix@C+35pt{\C_{\star}^{op} \times \C_{\star} \ar[r]^-{\underline{\mathrm{Hom}}_{\C_{\star}}(\,\cdot\,,\,\cdot\,)} & \simp},
\end{equation}
muni des foncteurs adjoints:
\begin{equation}\label{adjstar}
\C_{\star}
\vcenter{\xymatrix@C+15pt{
\phantom{a}\ar@{}[r]|{\perp}\ar@<-4pt>@/_10pt/[r]_{\underline{\mathrm{Hom}}_{\C_{\star}}(X,\,\cdot\,)}&\phantom{a}
\ar@<-4pt>@/_10pt/[l]_-{X\wedge\,\cdot\, }}}
\simp
\qquad\text{et}\qquad
\C_{\star}^{op}
\vcenter{\xymatrix@C+15pt{
\phantom{a}\ar@{}[r]|{\perp}\ar@<-4pt>@/_10pt/[r]_{\underline{\mathrm{Hom}}_{\C_{\star}}(\,\cdot\,,Y)}&\phantom{a}
\ar@<-4pt>@/_10pt/[l]_-{Y^{\wedge\,\cdot\,} }}}
\simp
\end{equation}
pour $X=(X,x_0)$ et $Y=(Y,y_0)$ des objets pointés de $\C$.

En effet, si $X=(X,x_{0})$ est un objet pointé de $\C$, pour tout ensemble simplicial $A$  on définit l'objet pointé $X\wedge A$ de $\C$ par un carré cocartésien dans $\C$:
\begin{equation}
\xymatrix{
\star \ar[r]& X\wedge A\\ 
\star\otimes A \ar[r]_-{x_{0}\otimes A} \ar[u]  & X \otimes A\ar[u]\,.}
\end{equation}

Si on définit l'enrichissement de $\C_{\star}$ par la règle:
$$
\underline{\mathrm{Hom}}_{\C_{\star}}(X,Y)_{n} \, = \, \mathrm{Hom}_{\C_\star}(X\wedge\Delta^n,Y)\,,\qquad \text{où \,$n\geq 0$;}
$$
on vérifie sans peine que \eqref{adjstar} sont bien des adjonctions, lorsque $Y^{\wedge\, A}= \big(Y^{A},y_{0}^A\big)$ où $y_{0}^A$ est le morphisme $\xymatrix{\star\cong\star^A \ar[r]^-{y_{0}^A} & Y^A}$.

\begin{lemme}\label{deuxilemme}
Soit $0\leq n\leq \infty$. Si $\C$ est une catégorie de modèles munie d'un enrichissement $\underline{\mathrm{Hom}}_{\C}$ dans $\simp$ vérifiant l'existence des adjonctions:
$$
\C
\vcenter{\xymatrix@C+15pt{
\phantom{a}\ar@{}[r]|{\perp}\ar@<-4pt>@/_10pt/[r]_{\underline{\mathrm{Hom}}_{\C}(X,\,\cdot\,)}&\phantom{a}
\ar@<-4pt>@/_10pt/[l]_-{X\otimes\,\cdot\, }}}
\simp
\qquad\text{et}\qquad
\C^{op}
\vcenter{\xymatrix@C+15pt{
\phantom{a}\ar@{}[r]|{\perp}\ar@<-4pt>@/_10pt/[r]_{\underline{\mathrm{Hom}}_{\C}(\,\cdot\,,Y)}&\phantom{a}
\ar@<-4pt>@/_10pt/[l]_-{Y^{\,\cdot\,} }}}
\simp\,,
$$
pour $X$ et $Y$ des objets quelconques de $\C$; alors $(\C,\underline{\mathrm{Hom}}_{\C})$ est une catégorie de modèles simpliciale d'ordre $n$ si et seulement si, $(\C_{\star},\underline{\mathrm{Hom}}_{\C_{\star}})$ est une catégorie de modèles simpliciale pointée d'ordre $n$.
\end{lemme}
\begin{proof}
La preuve est très similaire à la preuve du Lemme \ref{equipointe}, mais c'est fois on utilise que $-\,\wedge A$ commute aux coproduits pour montrer que $(X\otimes A)_{+} \cong (X_+) \wedge A$. 
\end{proof}

\begin{proof}[Démonstration du Corollaire \ref{dugger2}]
On vérifie aussi-tôt que (i)$\Rightarrow$(ii). En plus, (ii)$\Rightarrow$(iii) d'après le Théorème \ref{cassi1} et les Lemmes \ref{quillenomegasus}, \ref{lem2} et \ref{deuxilemme}. D'un autre, on sait que (iii)$\Leftrightarrow$(iv)$\Leftrightarrow$(v) pour une catégorie de modèles pointée $\C$ quelconque; dans notre cas on peut utiliser le Lemme \ref{classique}.  

Enfin (iii)$\Rightarrow$(i) se suit du Théorème \ref{cassi1}, la Proposition \ref{duggerT} et des Lemmes \ref{quillenomegasus}, \ref{lem2} et \ref{deuxilemme}.
\end{proof}

\section{Pré-monoïdes de Segal}

\renewcommand{\thesubsection}{\S\thesection.\arabic{subsection}}
\subsection{}\;\label{ensreduitlacat}
\renewcommand{\thesubsection}{\thesection.\arabic{subsection}}

On désigne par $\simp_{0}$ la sous-catégorie pleine de $\simp$ dont les objets sont les \emph{ensembles simpliciaux réduits}, \emph{i.e.} les ensembles simpliciaux $X$ vérifiant que son ensemble des $0$-simplexes soit l'ensemble ponctuel $X_{0}=\star$. 

Remarquons que le foncteur d'inclusion de $\simp_{0}$ vers $\simp$ se factorise:
\begin{equation}\label{triang22}
\xymatrix@C-5pt@R-15pt{
\underset{\phantom{a}}{\simp_{0}}\,\ar@{^(->}@/_5pt/[rd] \ar@{^(->}[rr] & &\simp\,,\\
&\, \simp_{\star}\,\ar@/_5pt/[ru]_-{\pi =\, \textit{d'oubli}}&} 
\end{equation}
car tout ensemble simplicial réduit est canoniquement pointé. 

Plus encore, on vérifie qu'on a des adjonctions:
\begin{equation}\label{adjredu}
\simp_{0}
\xymatrix@C+15pt{
\phantom{a}\ar@<+10pt>@{}[r]|{\perp}\ar@{^(->}[r]&\phantom{a}
\ar@<-7pt>@/_12pt/[l]_-{\mathcal{F}}}
\simp
\qquad\text{et}\qquad
\simp_{0}
\vcenter{\xymatrix@C+15pt{
\phantom{a}\ar@<-10pt>@{}[r]|{\perp}\ar@<+10pt>@{}[r]|{\perp}\ar@{^(->}[r]&\phantom{a}
\ar@<-7pt>@/_12pt/[l]_-{\mathcal{G}}\ar@<+7pt>@/^12pt/[l]^-{\mathcal{H}}}}
\simp_{\star}\,,
\end{equation} 
définis pour tout ensemble simplicial $X$ et tout $0$-simplexe $x_{0}\in X_0$, par un carré cocartésien et un carré cartésien de $\simp$:
$$
\vcenter{\xymatrix{\mathcal{F} X =\mathcal{G} (X,x_{0}) & X\ar[l] \\\star\ar[u] & {\bf sq}_{0}(X)\ar[l]\ar[u]}} \qquad \text{et} \qquad 
\vcenter{\xymatrix{\mathcal{H} (X,x_{0}) \ar[d]\ar[r] & X \ar[d]\\ \star \ar[r]_-{x_{0}} & {\bf csq}_{0}(X)\,,}}
$$
respectivement\footnote{Rappelons que le $0$-squelette ${\bf sq}_{0}(X)$ de $X$ est l'ensemble simplicial constant à valeur $X_{0}$, l'ensemble des $0$-simplexes de $X$. Tandis que le $0$-cosquelette ${\bf csq}_{0}(X)$ est un ensemble simplicial dont l'ensemble de $n$-simplexes est les produit cartésien $X_0^n$}.

Autrement dit, $\mathcal{G}(X,x_{0})=\mathcal{F}(X)$ est l'ensemble simplicial quotient $X\big/{\bf sq}_{0}X$, c'est-à-dire on construit $\mathcal{G}(X,x_{0})_{n}=\mathcal{F}(X)_{n}$ pour $n\geq 0$ en identifiant dans l'ensemble $X_{n}$ tous les $n$-simplexes totalement dégénérés à un seul point. D'un autre côté, pour $n\geq 0$ l'ensemble $\mathcal{H}(X,x_{0})_{n}$ est le sous-ensemble des $n$-simplexes de $X$ ayant la propriété que tous ses $0$-sommets sont égaux à $x_{0}$.

\begin{lemme}\label{cocommm}
La catégorie $\simp_0$ admet des petites limites et colimites. Plus précisément, si $I$ est une petite catégorie et $\gamma\colon\xymatrix@C-10pt{I\ar[r]&\simp_0}$ est un foncteur, alors:
\begin{enumerate}
\item Une limite de $\gamma$ dans $\simp$ est une limite de $\gamma$ dans $\simp_0$.  
\item Un quotient $\underset{I}{\mathrm{colim}}\,\gamma\Big/{\bf sq}_{0} \big( \underset{I}{\mathrm{colim}}\,\gamma \big)$ dans $\simp$ est une colimite de $\gamma$ dans $\simp_0$, toujours  que $\underset{I}{\mathrm{colim}}\,\gamma$ soit une colimite de $\gamma$ dans $\simp$. En particulier si $I$ est une catégorie connexe, une colimite de $\gamma$ dans $\simp$ est une colimite de $\gamma$ dans $\simp_0$.
\end{enumerate}
\end{lemme}
\begin{proof}
Montrons seulement que le foncteur $\xymatrix@-5pt{\simp_0\,\ar@{^(->}[r]&\simp}$ commute aux colimites sur de catégories connexes. 

En effet si $\gamma\colon\xymatrix@C-10pt{I\ar[r]&\simp_0}$ est un foncteur dont la source est une petite catégorie connexe et on pose $\underset{i\in I}{\mathrm{colim}}\,\gamma(i) $ pour noter une colimite de $\gamma$ dans la catégorie $\simp$, on a des isomorphismes:
$$
\Big( \underset{i\in I}{\mathrm{colim}}\,\gamma(i)\Big)_0 \, \cong \,  \underset{i\in I}{\mathrm{colim}}\,\big(\gamma(i)_0\big) \, \cong \, \pi_0(I) \, \cong \, \star\,;
$$
où $\underset{i\in I}{\mathrm{colim}}\,\big(\gamma(i)_0\big)$ est une colimite dans la catégorie des ensembles du foncteur contant $i\longmapsto \star$. 
\end{proof}

D'un autre côté, remarquons que l'enrichissement de $\simp_\star$ dans $\simp$ définit dans \eqref{hompointe23}, induit l'enrichissement suivant de $\simp_0$ dans les ensembles simpliciaux:
\begin{equation}\label{homred}
\begin{split}
\xymatrix@C+40pt{\simp_{0} \times \simp_{0} \ar[r]^-{\underline{\mathrm{Hom}}_{\simp_{0}}(\,\cdot\,,\,\cdot\,)} & \simp}
\qquad\text{où}\qquad
\underline{\mathrm{Hom}}_{\simp_{0}}(X,Y)_{n} \, & = \, \mathrm{Hom}_{\simp_{0}}(X \wedge \Delta^n_{+},Y)\\
&= \, \mathrm{Hom}_{\simp_*}(X \wedge \Delta^n_{+},Y)\,.
\end{split}
\end{equation}
(note que si $X$ est un ensemble simplicial réduit, alors $X \wedge K_{+}=(X\times K) \big/ (\star\times K)$ l'est aussi).

On vérifie en particulier qu'on a des adjonctions:
\begin{equation}\label{adjreduit}
\simp_{0}
\vcenter{\xymatrix@C+15pt{
\phantom{a}\ar@{}[r]|{\perp}\ar@<-4pt>@/_10pt/[r]_{\underline{\mathrm{Hom}}_{\simp_{0}}(X,\,\cdot\,)}&\phantom{a}
\ar@<-4pt>@/_10pt/[l]_-{X\overset{\circ}{\wedge}(\,\cdot\,)}}}
\simp
\qquad\text{et}\qquad
\simp_{0}^{op}
\vcenter{\xymatrix@C+15pt{
\phantom{a}\ar@{}[r]|{\perp}\ar@<-4pt>@/_10pt/[r]_{\underline{\mathrm{Hom}}_{\simp_{0}}(\,\cdot\,,Y)}&\phantom{a}
\ar@<-4pt>@/_10pt/[l]_-{Y^{\overset{\circ}{\wedge}\,\cdot\,} }}}
\simp\,,
\end{equation}
\begin{align*}
\text{où}\quad \qquad  X\overset{\circ}{\wedge}(\,K\,)\, & = \, X\wedge(\,K\,)_{+} \qquad \text{et}\\
\big(Y^{\overset{\circ}{\wedge}K}\big)_{n} \, 
& = \, \mathcal{H}\big(\underline{\mathrm{Hom}}_{\simp_{*}}(K_{+},Y),\star\big)_{n}\\
& = \, \mathrm{Hom}_{\simp_{0}}\big( (K\times \Delta^n) \big/ (K\times {\bf sq}_{0}\Delta^n)  ,Y\big)\,.
\end{align*}

\begin{proposition}\label{modred}
Si $0\leq n\leq \infty$, la catégorie des ensembles simpliciaux réduits $\simp_{0}$ munie de l'enrichissement \eqref{homred} ci-dessous admet une structure de catégorie de modèles simpliciale pointée d'ordre $n$ (voir aussi le Corollaire \ref{dimminus}), lorsque:
\begin{align*}
\big\{\,\text{équivalences faibles}\,\big\} \quad &=\quad 
\big\{\, f:\xymatrix@-14pt{X\ar[r]&Y}\,\big|\; \text{$f$ est une $n$-équivalence faible dans $\simp$}\,\big\} \;\;\,= \;\,  {\bf W}_{n}^{red}\,,\\ 
\big\{\,\text{cofibrations}\,\big\} \quad &=\quad  \big\{\,\text{monomorphismes}\,\big\} \;\;\,= \;\,  {\bf mono}\,,\\
\text{et}\qquad\big\{\,\text{fibrations}\,\big\} \quad &=\quad  \big\{\,\text{morphismes avec la propriété de relèvement }\\
 &\phantom{=}\;\;\quad \text{\phantom{$\{\,$}à droite par rapport à ${\bf mono}\,\cap\,{\bf W}_{n}^{red}$}\,\big\}\\
 &=\quad \big\{\,\text{$n$-fibrations réduites}\,\big\}\;\,= \;\, {\bf fib}_{n}^{red}\,.
\end{align*}

Un ensemble simplicial réduit $Z$ est fibrant dans $(\simp_{0},{\bf W}^{red}_{n}, {\bf mono},{\bf fib}_{n}^{red})$ si et seulement si $Z$ est un complexe de Kan tel que $\pi_{m}(Z)=\pi_{m}(Z,\star)=0$ pour tout $m\geq n+1$. 
\end{proposition}
\begin{proof}
Commençons pour le cas $n=\infty$. Dans la Proposition 6.2 de \cite{jardineb} on montre que $(\simp_{0},{\bf W}^{red}_{\infty}, {\bf mono},{\bf fib}_{\infty}^{red})$ est une catégorie de modèles; plus encore, dans le Lemme 6.6 de \cite{jardineb} on vérifie que si $\xymatrix@C-8pt{X\ar[r]^-f&Y}$ est un morphisme d'ensembles simpliciaux réduits, les énoncés suivants sont équivalents:
\begin{enumerate}
\item $f$ est une fibration de la catégorie de modèles $(\simp_{0},{\bf W}^{red}_{\infty}, {\bf mono},{\bf fib}_{\infty}^{red})$ et $f$ vérifie la propriété de relèvement à droite par rapport au morphisme $\xymatrix@C-10pt{\star\ar[r]&\mathbb{S}^1}$.
\item $f$ est une fibration de Kan.
\end{enumerate}

Vu que tout morphisme $\xymatrix@C-10pt{X\ar[r]&\star}$ de but un ensemble simplicial réduit vérifie la propriété de relèvement à droite par rapport au morphisme $\xymatrix@C-10pt{\star\ar[r]&\mathbb{S}^1}$; il se suit qu'un ensemble simplicial réduit $X$ est fibrant dans la catégorie de modèles $(\simp_{0},{\bf W}^{red}_{\infty}, {\bf mono},{\bf fib}_{\infty}^{red})$, si et seulement si $X$ est un complexe de Kan. 

Remarquons maintenant que $(\simp_{0},{\bf W}^{red}_{\infty}, {\bf mono},{\bf fib}_{\infty}^{red})$ munie de l'enrichissement \eqref{homred} est une catégorie de modèle simpliciale pointée (d'ordre $\infty$). En effet, d'après le Lemme \ref{equipointe} on doit vérifier les propriétés {\bf CMS1} et {\bf CMS2}. Les adjonctions de la propriété {\bf CMS1} sont données dans \eqref{adjreduit}. D'un autre, pour vérifier la propriété {\bf CMS2} on doit montrer l'énoncé suivant:
\begin{quote}
Si $\xymatrix@C-8pt{A\ar[r]^j&B}$ est un monomorphisme d'ensembles simpliciaux et $\xymatrix@C-8pt{X\ar[r]^q&Y}$ est un monomorphisme d'ensembles simpliciaux réduits, alors le morphisme $\varphi$ dans le diagramme somme amalgamée de $\simp_{0}$ suivant:
\begin{equation}\label{carproof}
\xymatrix@R-10pt@C-10pt{
X\wedge A_{+} \ar[dd]_{X\wedge \,j_{+}}\ar[rr]^{q\,\wedge A_{+}} && 
Y\wedge A_{+} \ar[dd]\ar@/^15pt/[rddd]^{Y\wedge \,j_{+}}&\\\\
X\wedge B_{+}\ar[rr]\ar@/_15pt/[rrrd]_{q\,\wedge B_{+}}&&
X\wedge B_{+}\underset{X\wedge A_{+}}{\bigsqcup}Y\wedge A_{+}\ar@{-->}[rd]|(.65)\varphi&\\
&&&Y\wedge B_{+}\,,}
\end{equation}
est un monomorphisme d'ensembles simpliciaux réduits, lequel est une $\infty$-équivalence faible si $j$ est une $\infty$-équivalence faible d'ensembles simpliciaux, ou si $q$ est une $\infty$-équivalence faible entre ensembles simpliciaux réduits.
\end{quote}

Cette propriété c'est une conséquence du fait que le foncteur d'inclusion de $\simp_{0}$ vers $\simp_{*}$ commute aux colimites, et parce que  $(\simp_{\star},\pi^{-1}{\bf W}_{\infty}, {\bf mono}, \pi^{-1}{\bf fib}_{\infty},\underline{\mathrm{Hom}}_{\simp_*})$ est une catégorie de modèles simpliciale (d'ordre $\infty$).

D'un autre, remarquons que d'après un résultat de Smith (voir la Proposition 2.2 de \cite{barwick}) la catégorie de modèles $(\simp_{0},{\bf W}^{red}_{\infty}, {\bf mono},{\bf fib}_{\infty}^{red})$ est une catégorie de modèles combinatoire si et seulement si, $\simp_0$ est une catégorie localement présentable (voir \cite{adamek}) et il existe un ensemble $I$ de monomorphismes de $\simp_0$ tel que ${\bf cof}(I) = {\bf monos}$\footnote{${\bf cof}(I)$ note l'ensemble de tous les rétracts des compositions transfinies d'images directes d'éléments de $I$.}.

\begin{lemme}\label{lemm1}
La catégorie $\simp_0$ des ensembles simpliciaux réduits est localement présentable.
\end{lemme}
\begin{proof}
D'après la Proposition 1.46 de \cite{adamek} il suffit de montrer que $\simp_{0}$ est équivalente à une sous-catégorie pleine et réflective d'une catégorie de préfaisceaux $\mathbf{Ens}^{\mathcal{A}}$ laquelle soit fermée par des colimites directes (ou de façon équivalente par des colimites filtrantes).

Pour cela rappelons qu'on a une adjonction:
$$
\simp_{0}
\xymatrix@C+15pt{
\phantom{a}\ar@<+10pt>@{}[r]|{\perp}\ar@{^(->}[r]|-{\,\nu\,}&\phantom{a}
\ar@<-7pt>@/_12pt/[l]_-{\mathcal{F} = \,\cdot\, \big/ \mathbf{sq}_0(\,\cdot\,)}}
\simp\,,
$$
où le foncteur pleinement fidèle $\nu$ (l'inclusion canonique) commute aux colimites des diagrammes des foncteurs $\gamma:\xymatrix@C-8pt{I\ar[r]&\simp_{0}}$ où $I$ admet un objet initial (voir le Lemme \ref{cocommm}). 
\end{proof}

Montrons:

\begin{lemme}\label{lemm2}
Si $\xymatrix@C-7pt{X\ar[r]^\varphi&Y}$ es un monomorphisme d'ensembles simpliciaux réduits, alors $\varphi$ est la composition transfinie d'images directes dans la catégorie $\simp_0$ de morphismes dans l'ensemble:
$$
\bigg\{ \quad \xymatrix@+35pt{\partial\Delta^{n}/{\bf sq}_{0}\partial\Delta^n\; \ar@{^(->}[r]^-{\alpha^{n-1}/{\bf sq}_{0}\alpha^{n-1}}& \, \Delta^{n}/{\bf sq}_{0}\Delta^n} 
\quad \bigg| \quad n\geq 1 \quad \bigg\}\,.
$$
\end{lemme}
\begin{proof}
Si $\xymatrix@C-7pt{X\ar[r]^\varphi&Y}$ es un morphisme d'ensembles simpliciaux réduits, on va définir pour chaque $i\geq 0$ un ensemble simplicial réduit $X^{(i)}$ et un morphisme $\xymatrix@C-7pt{X^{(i)}\ar[r]^{g^i}& Y}$ tel que $X^0=X$ et $g^0=\varphi$; et si $i< j$ on va construire un morphisme $\xymatrix@C-7pt{X^{(i)}\ar[r]^{\varphi_i^j}&X^{(j)}}$ tel que:
$$
\xymatrix@R=6pt@C+18pt{
&X^{(i)}\ar[dd]|-{\varphi_i^j}\ar@/^6pt/[rd]^-{g^i}& \\
&&Y\;,\\
&X^{(j)}\ar@/_6pt/[ru]_-{g^j}&}
$$
soit un triangle commutatif et $\varphi_k^j \circ \varphi_i^k = \varphi_i^j$ toujours que $i< k< j$. 

Pour commencer on pose $X^{(0)}=X$ et on définit $\xymatrix@C+3pt{X^{(0)}\ar[r]^{g^0=\varphi}& Y}$ comme il est imposé. Si $i\geq 1$ on définit de façon inductive $X^{(i)}$ en prenant un carré cocartésien dans $\simp$ $\big($lequel est aussi un carré cocartésien dans $\simp_0$ $\big)$:
\begin{equation}\label{universele}
\xymatrix@R-10pt@C-5pt{
{\bf sq}_{i}(X^{(i-1)}) \ar[rr]\ar[dd]_-{{\bf sq}_{i}(g^{i-1})} && X^{(i-1)}\ar[dd]^-{\varphi_{i-1}^i} \\&&\\
{\bf sq}_{i}(Y) \ar[rr]& & X^{(i)}}
\end{equation}
le morphisme $g^{i}\colon\xymatrix@C-6pt{X^{(i)} \ar[r] &Y}$ étant donné par la propriété universelle:
\begin{equation}\label{universelle}
\xymatrix@R-10pt@C-5pt{
{\bf sq}_{i}(X^{(i-1)}) \ar[rr]\ar[dd]_-{{\bf sq}_{i}(g^{i-1})} && X^{(i-1)}\ar[dd]^-{\varphi_{i-1}^i} \ar@/^15pt/[dddr]^-{g^{i-1}}&\\&&&\\
{\bf sq}_{i}(Y)\ar@/_15pt/[rrrd] \ar[rr]& & X^{(i)}\ar@{-->}[rd]|-{g^{i}} &\\
&&&Y\,.}
\end{equation}

On définit finalement $\varphi_i^j=\varphi^j_{j-1}\circ\dots\circ\varphi_i^{i+1}$ pour $i<j$.

Remarquons que pour $i\geq 1$ le morphisme $g^{i}$ défini dans le diagramme \eqref{universelle} vérifie que les fonctions entre les ensembles des $n$-simplexes $g^{i}_n\colon\xymatrix@C-6pt{X^{(i)}_n \ar[r] &Y_n}$ sont bijectives pour $0\leq n\leq i$, vu que le morphisme $\xymatrix@C-6pt{{\bf sq}_{i}(Z)\ar[r]&Z}$ vérifie cette propriété pour tout ensemble simplicial $Z$. En fait, cela reste vrai pour $i\geq 0$ parce que $g^0=\varphi$ est un morphisme d'ensembles simpliciaux réduits. 

Si on suppose que $\xymatrix@C-7pt{X\ar[r]^\varphi&Y}$ soit un monomorphisme d'ensembles simpliciaux réduits, on vérifie par induction que les morphismes $g^{i}$ pour $i\geq 0$ sont aussi de monomorphismes. En particulier le morphisme $\xymatrix@C+10pt{\underset{i}{\mathrm{colim}}\,X^{(i)}\ar[r]^-{\mathrm{colim}\,g^i}& Y}$ dans le triangle commutatif:
$$
\xymatrix@C+10pt{\underset{i}{\mathrm{colim}}\,X^{(i)}\ar[r]^-{\mathrm{colim}\,g^i}& Y\\
X\ar[u]\ar@/_3pt/[ur]_-{\varphi}}
$$
où $\underset{i}{\mathrm{colim}}\,X^{(i)}$ est une colimite dans $\simp_0$ ou $\simp$ (elles sont isomorphes), est un isomorphisme d'ensembles simpliciaux réduits. Autrement dit, on vient de décrire le morphisme $\varphi$ comme la composition transfinie dans $\simp_0$ ou $\simp$ des morphismes $\big\{\varphi_{i-1}^{i}\big\}_{i\geq 1}$.

Plus encore, remarquons que des propriétés soulignées des morphismes $g^{i}$ on en déduit que pour $i\geq 1$ on a un carré cocartésien de la catégorie $\simp$:
$$
\xymatrix@R-10pt@C-5pt{
 \underset{\Omega_{i}}{\bigsqcup} \, \partial\Delta^{i}\ar[dd]_-{\bigsqcup \alpha^{i-1}} \ar[rr]&&{\bf sq}_{i}(X^{(i-1)})\ar[dd]^-{{\bf sq}_{i}(g^{i-1})}  \\&&\\
 \underset{\Omega_{i}}{\bigsqcup} \, \Delta^{i} \ar[rr]&&{\bf sq}_{i}(Y)}
$$
où $\Omega_{i}$ est l'ensemble des $i$-simplexes de $Y$ que ne sont pas dans l'image de la fonction $g^{i-1}_i$ (c'est-à-dire $Y_i\backslash X_i$). Donc on a un carré cocartésien dans la catégorie des ensembles simpliciaux réduits $\simp_0$:
\begin{equation}\label{universellle}
\xymatrix@R-10pt@C-5pt{
 \underset{\Omega_{i}}{\bigvee} \, \big(\partial\Delta^{i}/{\bf sq}_{0}\partial\Delta^{i}\big)\ar[dd]_-{\bigvee (\alpha^{i-1}/{\bf sq}_{0}\alpha^{i-1})} \ar[rr]&&{\bf sq}_{i}(X^{(i-1)})\ar[dd]^-{{\bf sq}_{i}(g^{i-1})}  \\&&\\
 \underset{\Omega_{i}}{\bigvee} \, \big(\Delta^{i}/{\bf sq}_{0}\Delta^n\big) \ar[rr]&&{\bf sq}_{i}(Y)\,,}
\end{equation}
où $\underset{k}{\bigvee} Z_k$ note une somme dans $\simp_0$ de la famille d'ensembles simpliciaux réduits $\big\{Z_k\big\}_k$, c'est-à-dire $\underset{k}{\bigvee} Z_k = \underset{k}{\bigsqcup} Z_k \big/ \big({\bf sq}_{0} \, \underset{k}{\bigsqcup} Z_k\big)$.

Il se suit des carré cocartésiens \eqref{universele} et \eqref{universellle} qu'on a un carré cocartésien de $\simp_0$:
$$
\xymatrix@R-10pt@C-5pt{
 \underset{\Omega_{i}}{\bigvee} \, \big(\partial\Delta^{i}/{\bf sq}_{0}\partial\Delta^{i}\big)\ar[dd]_-{\bigvee (\alpha^{i-1}/{\bf sq}_{0}\alpha^{i-1})} \ar[rr]&&X^{(i-1)}\ar[dd]^-{\varphi_{i-1}^{i}}  \\&&\\
 \underset{\Omega_{i}}{\bigvee} \, \big(\Delta^{i}/{\bf sq}_{0}\Delta^n\big) \ar[rr]&&X^{(i)}\,;}
$$
donc d'après les Lemmes 2.1.12 et 2.1.13 de \cite{hovey} $\varphi$ est la composition transfinie d'images directes dans $\simp_0$ d'éléments de l'ensemble:
$$
\bigg\{ \quad \xymatrix@+35pt{\partial\Delta^{n}/{\bf sq}_{0}\partial\Delta^n\; \ar@{^(->}[r]^-{\alpha^{n-1}/{\bf sq}_{0}\alpha^{n-1}}& \, \Delta^{n}/{\bf sq}_{0}\Delta^n} 
\quad \bigg| \quad n\geq 1 \quad \bigg\}\,.
$$
\end{proof}

On conclut que $(\simp_{0},{\bf W}^{red}_{\infty}, {\bf mono},{\bf fib}_{\infty}^{red})$ est bien une catégorie de modèles propre à gauche et combinatoire, donc d'après le Théorème 4.7 de \cite{barwick} on peut considérer sa localisation de Bousfield à gauche par rapport à l'ensemble des morphismes:
\begin{equation}\label{elsnred}
S_{n}^{red}\,=\,\Big\{ \quad \xymatrix@C+6pt{\star\, \ar@{^(->}[r]^-{}&\, \mathbb{S}^m \quad \Big| \quad m>n }\Big\}\,.
\end{equation}

Dans les Lemmes suivants on pose $[\,\cdot\,,\,\cdot\,]_{\infty}^{red}$ pour noter l'ensemble des morphismes dans la catégorie de fractions $\simp_0\big[({\bf W}^{red}_{\infty})^{-1}\big]$.

\begin{lemme}\label{lolocal} 
Soit $m\geq 1$ et $Z$ un ensemble simplicial réduit, alors $\pi_{m}(Z)=\pi_{m}(Z,\star)=0$ si et seulement si la fonction:
$$
\xymatrix@C+10pt{ [\mathbb{S}^{m},Z]_\infty^{\text{\tiny{$\text{red}$}}}\ar[r]^-{(\alpha^m)^*}&[\star,Z]_\infty^{\text{\tiny{$\text{red}$}}}}\,,
$$
induite par le morphisme $\xymatrix@C+6pt{\star\ar@{^(->}[r]^-{\alpha^m}&\mathbb{S}^{m}}$ dans la catégorie homotopique $\simp_0\big[({\bf W}^{red}_{\infty})^{-1}\big]$ est bijective.  
 
En particulier, si $n\geq 0$ un ensemble simplicial réduit $Z$ vérifie que $\pi_{m}(Z)=\pi_{m}(Z,\star)=0$ pour tout $m\geq n+1$ si et seulement si, $Z$ est un objet $S_{n}^{red}$-local de $(\simp_{0},{\bf W}^{red}_{\infty}, {\bf mono},{\bf fib}_{\infty}^{red})$.
\end{lemme}
\begin{proof}
Remarquons pour commencer que si $A$ et $B$ sont des ensembles simpliciaux réduits, d'après les définitions de $\underline{\mathrm{Hom}}_{\simp_{\star}}$ et $\underline{\mathrm{Hom}}_{\simp_{0}}$ dans \eqref{hompointe23} et \eqref{homred} respectivement, on a une bijection naturelle:
\begin{equation}\label{redpoin}
[A,B]_{\infty}^{\text{\tiny{$\text{red}$}}} \, \cong \, \pi_{0} \big(\underline{\mathrm{Hom}}_{\simp_{0}} (A,B')\big) \, = \, \pi_{0}\Big(\underline{\mathrm{Hom}}_{\simp_{\star}} \big( (A,\star),(B',\star)\big)\Big)
 \, \cong \, [(A,\star),(B,\star)]_{\infty}^{\text{\tiny{$\text{pointé}$}}}\,,
\end{equation}
où $B'$ est un remplacement fibrant de $B$ dans $(\simp_{0},{\bf W}^{red}_{\infty}, {\bf mono},{\bf fib}_{\infty}^{red})$; c'est-à-dire $B'$ est un ensemble simplicial réduit lequel est un complexe de Kan et en particulier $(B',\star)$ est un objet fibrant de $(\simp_{\star},\pi^{-1}{\bf W}_{\infty}, {\bf mono},\pi^{-1}{\bf fib}_{\infty})$ aussi.

Il se suit que si $Z$ est un ensemble simplicial réduit alors $[\mathbb{S}^{m},Z]_\infty^{\text{\tiny{$\text{red}$}}} \, \cong \, \pi_m(Z)$ (voir la bijection \eqref{isosi}); donc la fonction:
$$
\xymatrix@C+10pt{\pi_m(Z) \, \cong \, [\mathbb{S}^{m},Z]_\infty^{\text{\tiny{$\text{red}$}}}  \ar[r]^-{(\alpha^m)^*}&[\star,Z]_\infty^{\text{\tiny{$\text{red}$}}}\, \cong \, \star }\,,
$$
est bijective si et seulement si $\pi_m(Z)=0$. 
\end{proof}

\begin{lemme}\label{talvezsi}
Si $n\geq 0$, un morphisme $f:\xymatrix@-8pt{X\ar[r]&Y}$ d'ensembles simpliciaux réduits est une équivalence faible $S_{n}^{red}$-local, si et seulement si $f$ est une $n$-équivalence faible d'ensembles simpliciaux.
\end{lemme}
\begin{proof}
La preuve de ce Lemme est analogue à ceci du Lemme \ref{lemmelocal2} en considérant une version pointé du cosquelette ${\bf csq}^{\text{\tiny{$\text{pointé}$}}}_{n}(X,x_0) = \big({\bf csq}_{n}(X),x_0\big)$ (voir \ref{versionpointe}) et en remarquant que d'après le formules \eqref{hompointe23} et \eqref{homred} on a une égalité:
$$
\underline{\mathrm{Hom}}_{\simp_{0}} (A,B)\, = \,\underline{\mathrm{Hom}}_{\simp_{\star}} \big( (A,\star),(B,\star)\big)\,,
$$
pour des ensembles simpliciaux réduits quelconques $A$ et $B$.

Plus en détail: Soit $f\colon\xymatrix@C-8pt{X\ar[r]&Y}$un morphisme d'ensembles simpliciaux réduits et $Z$ un ensemble simplicial réduit. Si on considère de remplacements fibrants de $f$ et $Z$ dans  la catégorie de modèles $(\simp_0,{\bf W}^{red}_{\infty}, {\bf mono},{\bf fib}^{red}_{\infty})$, \emph{i.e.} un carré commutatif et une flèche:
$$
\vcenter{\xymatrix@R-6pt{X\ar[d]\ar[r]^f & Y\ar[d]\\ X' \ar[r]_{f'} & Y'}}\qquad\;\text{et}\qquad\;\vcenter{\xymatrix@R-6pt{Z\ar[d]\\Z'}}\,,
$$
où les morphismes verticaux sont des $\infty$-équivalences faibles entre ensembles simpliciaux réduits de but de complexes de Kan; alors la fonction:
\begin{equation}\label{inftyred}
\xymatrix@C+15pt{
[Y,Z]^{red}_{\infty}\ar[r]^-{f^*} & [X,Z]^{red}_{\infty}}
\end{equation}
entre les ensembles de morphismes de la catégorie $\simp_0\big[({\bf W}^{red}_{\infty})^{-1}\big]$, s'identifie avec la fonction: 
\begin{equation}\label{infty1ared}
\xymatrix@R=3pt@C+15pt{
\pi_{0}\Big(\underline{\mathrm{Hom}}_{\simp_0}\big(Y',Z'\big)\Big)\ar[r]^{(f')^*}\ar@{}[d]|-{\mathrel{\reflectbox{\rotatebox[origin=c]{90}{$=$}}}}& 
\pi_{0}\Big(\underline{\mathrm{Hom}}_{\simp_0}\big(X',Z'\big)\Big)\ar@{}[d]|-{\mathrel{\reflectbox{\rotatebox[origin=c]{90}{$=$}}}}\\
\pi_{0}\Big(\underline{\mathrm{Hom}}_{\simp_{\star}} \big((Y',\star),(Z',\star)\big)\Big)\ar[r]_{(f')^*}&\pi_{0}\Big(\underline{\mathrm{Hom}}_{\simp_{\star}} \big( (X',\star),(Z',\star) \big)\Big)\,.}
\end{equation}

Si on suppose que $Z$ soit un ensemble simplicial réduit $S^{red}_{n}$-local, il se suit du Lemme \ref{lolocal} et du Corollaire \ref{lecoro} que le morphisme:
$$
\xymatrix@C-8pt{Z'\ar[r] & {\bf csq}_{n+1}(Z')}
$$
est une $\infty$-équivalence faible entre ensembles simpliciaux réduits lesquels sont des objets fibrants de la catégorie de modèles $(\simp_{0},{\bf W}^{red}_{\infty}, {\bf mono},{\bf fib}^{red}_{\infty})$; en particulier:
$$
\xymatrix@C-8pt{(Z',\star)\ar[r] & \big({\bf csq}_{n+1}(Z'),\star\big)={\bf csq}^{\text{\tiny{$\text{pointé}$}}}_{n+1}(Z',\star)}\,,
$$
est une $\infty$-équivalence faible entre ensembles simpliciaux pointés lesquels sont des objets fibrants de la catégorie de modèles simpliciale $(\simp_\star,\pi^{-1}{\bf W}_{\infty}, {\bf mono},\pi^{-1}{\bf fib}_{\infty})$.

Donc, \eqref{infty1ared} s'identifie avec la fonction:
$$
\xymatrix@C+15pt{
\pi_{0}\Big(\underline{\mathrm{Hom}}_{\simp_\star}\big((Y',\star),{\bf csq}^{\text{\tiny{$\text{pointé}$}}}_{n+1}(Z',\star)\big)\Big)\ar[r]^-{(f')^*}& \pi_{0}\Big(\underline{\mathrm{Hom}}_{\simp_\star}\big((X',\star),{\bf csq}^{\text{\tiny{$\text{pointé}$}}}_{n+1}(Z',\star)\big)\Big)\,,}
$$
qui d'après le Lemme \ref{lemmepoi2} s'identifie avec la fonction: 
\begin{equation}\label{infty2ared}
\xymatrix@C-8pt@R=3pt{
\pi_{0}\Big(\underline{\mathrm{Hom}}_{\simp_\star}\big({\bf csq}^{\text{\tiny{$\text{pointé}$}}}_{n+1}(Y',\star),{\bf csq}^{\text{\tiny{$\text{pointé}$}}}_{n+1}(Z',\star)\big)\Big)\ar[r]\ar@{}[d]|-{\mathrel{\reflectbox{\rotatebox[origin=c]{90}{$=$}}}}& \pi_{0}\Big(\underline{\mathrm{Hom}}_{\simp_\star}\big({\bf csq}^{\text{\tiny{$\text{pointé}$}}}_{n+1}(X',\star),{\bf csq}^{\text{\tiny{$\text{pointé}$}}}_{n+1}(Z',\star)\big)\Big)\ar@{}[d]|-{\mathrel{\reflectbox{\rotatebox[origin=c]{90}{$=$}}}}\\
\pi_{0}\Big(\underline{\mathrm{Hom}}_{\simp_0}\big({\bf csq}_{n+1}Y',{\bf csq}_{n+1}X'\big)\Big)
\ar@{}[d]|-{\mathrel{\reflectbox{\rotatebox[origin=c]{90}{$\cong$}}}}
& \pi_{0}\Big(\underline{\mathrm{Hom}}_{\simp_0}\big({\bf csq}_{n+1}X',{\bf csq}_{n+1}X'\big)\Big)
\ar@{}[d]|-{\mathrel{\reflectbox{\rotatebox[origin=c]{90}{$\cong$}}}}\\
[{\bf csq}_{n+1}Y',{\bf csq}_{n+1}X'\big]^{red}_{\infty}\ar[r]_-{{\bf csq}_{n+1}(f')^*}& \big[{\bf csq}_{n+1}X',{\bf csq}_{n+1}X'\big]^{red}_{\infty}\,.}
\end{equation}

Si on suppose que $f$ soit une $n$-équivalence faible d'ensembles simpliciaux réduits, il se suit du Corollaire \ref{lecoro} que ${\bf csq}_{n+1}(f')$ est une $\infty$-équivalence faible, \emph{i.e.} \eqref{infty2ared} est une fonction bijective. Donc, si $f$ est une $n$-équivalence faible la fonction \eqref{inftyred} est bijective pour tout ensemble simplicial réduit $S^{red}_{n}$-local $Z$, c'est-à-dire, $f$ est une équivalence faible $S^{red}_{n}$-locale.

Soit maintenant $f\colon\xymatrix@C-8pt{X\ar[r]&Y}$ un morphisme d'ensembles simpliciaux réduits lequel soit une équivalence faible $S^{red}_{n}$-locale, et considérons un remplacement fibrant de $f$ dans la catégorie de modèles $(\simp_0,{\bf W}^{red}_{\infty}, {\bf mono},{\bf fib}^{red}_{\infty})$, \emph{i.e.} un carré commutatif:
$$
\vcenter{\xymatrix@R-6pt{X\ar[d]\ar[r]^f & Y\ar[d]\\ X' \ar[r]_{f'} & Y'\,,}}
$$
où les morphismes verticaux sont $\infty$-équivalences faibles entre ensembles simpliciaux réduits dont les buts sont des complexes de Kan.

Remarquons d'ailleurs que d'après le Corollaire \ref{lecoro} et le Lemme \ref{lemmelocala} les ensembles simpliciaux ${\bf csq}_{n+1}X'$ et ${\bf csq}_{n+1}Y'$ sont $S^{red}_{n}$-locaux; en particulier la fonction:
\begin{equation}\label{deuxeqared}
\xymatrix@R=5pt@C+15pt{
[Y,{\bf csq}_{n+1}X']^{red}_{\infty}\ar@{}[d]|-{\mathrel{\reflectbox{\rotatebox[origin=c]{90}{$\cong$}}}}
\ar[r]^-{f^*} & [X,{\bf csq}_{n+1}X']^{red}_{\infty}\ar@{}[d]|-{\mathrel{\reflectbox{\rotatebox[origin=c]{90}{$\cong$}}}}\\
\pi_{0}\Big(\underline{\mathrm{Hom}}_{\simp_0}\big(Y',{\bf csq}_{n+1}X'\big)\Big)\ar@{}[d]|-{\mathrel{\reflectbox{\rotatebox[origin=c]{90}{$=$}}}}
& \pi_{0}\Big(\underline{\mathrm{Hom}}_{\simp_0}\big(X',{\bf csq}_{n+1}X'\big)\Big)\ar@{}[d]|-{\mathrel{\reflectbox{\rotatebox[origin=c]{90}{$=$}}}}\\
\pi_{0}\Big(\underline{\mathrm{Hom}}_{\simp_\star}\big((Y',\star),{\bf csq}^{\text{\tiny{$\text{pointé}$}}}_{n+1}(X',\star)\big)\Big)
\ar[r]_-{(f')^*}& \pi_{0}\Big(\underline{\mathrm{Hom}}_{\simp_\star}\big((X',\star),{\bf csq}^{\text{\tiny{$\text{pointé}$}}}_{n+1}(X',\star)\big)\Big)
\,,}
\end{equation}
est bijective.

D'un autre côté, grâce au Lemme \ref{lemehomo} la fonction bijective \eqref{deuxeqared} s'identifie avec la fonction:
\begin{equation}\label{deuxeq2red}
\xymatrix@C+15pt@R=5pt{
\pi_{0}\Big(\underline{\mathrm{Hom}}_{\simp_\star}\big({\bf csq}^{\text{\tiny{$\text{pointé}$}}}_{n+1}(Y',\star),{\bf csq}^{\text{\tiny{$\text{pointé}$}}}_{n+1}(X',\star)\big)\Big)
\ar@{}[d]|-{\mathrel{\reflectbox{\rotatebox[origin=c]{90}{$=$}}}}
\ar[r]^{{\bf csq}^{\text{\tiny{$\text{pointé}$}}}_{n+1}(f')^*}& \pi_{0}\Big(\underline{\mathrm{Hom}}_{\simp_\star}\big({\bf csq}^{\text{\tiny{$\text{pointé}$}}}_{n+1}(X',\star),{\bf csq}^{\text{\tiny{$\text{pointé}$}}}_{n+1}(X',\star)\big)\Big)
\ar@{}[d]|-{\mathrel{\reflectbox{\rotatebox[origin=c]{90}{$=$}}}}\\
\pi_{0}\Big(\underline{\mathrm{Hom}}_{\simp_0}\big({\bf csq}_{n+1}Y',{\bf csq}_{n+1}X'\big)\Big)
\ar@{}[d]|-{\mathrel{\reflectbox{\rotatebox[origin=c]{90}{$\cong$}}}}
& \pi_{0}\Big(\underline{\mathrm{Hom}}_{\simp_0}\big({\bf csq}_{n+1}X',{\bf csq}_{n+1}X'\big)\Big)
\ar@{}[d]|-{\mathrel{\reflectbox{\rotatebox[origin=c]{90}{$\cong$}}}}\\
[{\bf csq}_{n+1}Y',{\bf csq}_{n+1}X'\big]^{red}_{\infty}\ar[r]_-{{\bf csq}_{n+1}(f')^*}& \big[{\bf csq}_{n+1}X',{\bf csq}_{n+1}X'\big]^{red}_{\infty}\,.}
\end{equation}

Soit $g:\xymatrix@C-8pt{{\bf csq}_{n+1}Y'\ar[r]&{\bf csq}_{n+1}X'}$ un morphisme dans la catégorie homotopique $\simp_0\big[({\bf W}^{red}_\infty)^{-1}\big]$, tel que le composé $g\circ {\bf csq}_{n+1}(f')$ soit le morphisme identité de l'objet ${\bf csq}_{n+1}X'$ dans cette catégorie de fractions.

On vérifie alors que la fonction composée: 
$$
\xymatrix@C+16pt{
 [{\bf csq}_{n+1}X',{\bf csq}_{n+1}Y']^{red}_{\infty} \ar[r]^-{g^*}&
 [{\bf csq}_{n+1}Y',{\bf csq}_{n+1}Y']^{red}_{\infty}\ar[r]^-{{\bf csq}_{n+1}(f')^*} & 
 [{\bf csq}_{n+1}X',{\bf csq}_{n+1}Y']^{red}_{\infty}\,,}
$$
est égale à la fonction identité; donc:
$$
\xymatrix@C+16pt{
 [{\bf csq}_{n+1}Y',{\bf csq}_{n+1}Y']^{red}_{\infty} \ar[r]^-{{\bf csq}_{n+1}(f')^*}&
 [{\bf csq}_{n+1}X',{\bf csq}_{n+1}Y']^{red}_{\infty}\ar[r]^-{g^*} & 
 [{\bf csq}_{n+1}Y',{\bf csq}_{n+1}Y']^{red}_{\infty}\,,}
$$
est aussi égale à la fonction identité, parce que la flèche:
$$
\xymatrix@C+16pt{
[{\bf csq}_{n+1}Y',{\bf csq}_{n+1}Y']^{red}_{\infty} \ar[r]^-{{\bf csq}_{n+1}(f')^*} & [{\bf csq}_{n+1}X',{\bf csq}_{n+1}Y']^{red}_{\infty}\,,}
$$
est bijective.

Il se suit que ${\bf csq}_{n+1}(f')\circ g$ est aussi le morphisme identité dans la catégorie $\simp_0\big[({\bf W}^{red}_\infty)^{-1}\big]$. Autrement dit, ${\bf csq}_{n+1}(f')$ est une $\infty$-équivalence faible. Donc, $f$ est une $n$-équivalence faible d'après  le Corollaire \ref{lecoro}.
\end{proof}

On déduit du Théorème 4.7 de \cite{barwick} et des Lemmes \ref{lolocal} et \ref{talvezsi} qu'on a bien une catégorie de modèles $(\simp_0,{\bf W}^{red}_n, {\bf mono},{\bf fib}^{red}_n)$ dont les objets fibrants sont les ensembles simpliciaux réduits lesquels sont les complexes de Kan $X$ tels que $\pi_m(X)=0$ pour tout $m\geq n+1$.

Finalement si $0\leq n < \infty$ on montre que $(\simp_{0},{\bf W}^{red}_{n}, {\bf mono},{\bf fib}_{n}^{red})$ munie de l'enrichissement \eqref{homred} est une catégorie de modèle simpliciale pointée d'ordre $n$ de façon analogue au cas $n=\infty$; cela étant possible parce que $(\simp_{\star},\pi^{-1}{\bf W}_{n}, {\bf mono}, \pi^{-1}{\bf fib}_{n},\underline{\mathrm{Hom}}_{\simp_*})$ est une catégorie de modèles simpliciale d'ordre $n$ (voir le Corollaire \ref{pointeco2} et les Lemmes \ref{equipointe} et \ref{deuxilemme}).
\end{proof}

On désigne par $\mathrm{Ho}_{n}\big(\,\simp_{0}\,\big)$ la \emph{catégorie des $n$-types d'homotopie réduits}, \emph{i.e.} la catégorie homotopique de la catégorie de modèles $(\simp_{0},{\bf W}^{red}_{\infty}, {\bf mono},{\bf fib}_{n}^{red})$, et on note $[\,\cdot\,,\,\cdot\,]_{n}^{red}$ l'ensemble des morphismes dans $\mathrm{Ho}_{n}(\simp_{0})$. Si $X$ est un ensemble simplicial réduit, la classe associée à $X$ dans l'ensemble des objets à isomorphisme près de la catégorie $\mathrm{Ho}_{n}\big(\,\simp_{0}\,\big)$, est appelée le \emph{$n$-type d'homotopie réduit de $X$}. 

Rappelons:

\begin{lemme}\label{hachelemme}
Si $0\leq n\leq \infty$, l'adjonction définie dans \eqref{adjredu} plus haut:
\begin{equation}\label{adjconne}
\simp_{0}
\xymatrix@C+15pt{
\phantom{a}\ar@<+10pt>@{}[r]|{\text{\rotatebox[origin=c]{180}{$\perp$}}}\ar@{^(->}[r]_-{\nu}&\phantom{a}
\ar@<-7pt>@/_12pt/[l]_-{\mathcal{H}}}
\simp_\star
\end{equation}
est une adjonction de Quillen par rapport aux structures de catégories de modèles de la Proposition \ref{modred} et du Corollaire \ref{pointeco2}.

Plus encore, l'adjonction induite entre les catégories homotopiques:
$$
\mathrm{Ho}_n\big(\simp_{0}\big)
\xymatrix@C+15pt{
\phantom{a}\ar@<+10pt>@{}[r]|{\text{\rotatebox[origin=c]{180}{$\perp$}}}\ar[r]_{{\bf L}\nu}&\phantom{a}
\ar@<-7pt>@/_12pt/[l]_-{{\bf R}\mathcal{H}}}
\mathrm{Ho}_n\big(\simp_\star\big)
$$
détermine une équivalence entre la catégorie des $n$-types d'homotopie réduits $\mathrm{Ho}_n\big(\simp_{0}\big)$ et la sous-catégorie pleine de $\mathrm{Ho}_n\big(\simp_\star\big)$ des $n$-types d'homotopie pointés connexes, \emph{i.e.} les $X$ tels que $\pi_0(X)=0$. 
\end{lemme}
\begin{proof}
Le foncteur d'inclusion $\simp_{0}\xymatrix@C-10pt{\phantom{a}\ar@{^(->}[r]&\phantom{a}}\simp_\star$ respecte évidemment les monomorphismes et les $n$-équivalences faibles entre des ensembles simpliciaux réduits. Donc \eqref{adjconne} est bien une adjonction de Quillen.

Si $A$ est un ensemble simplicial réduits lequel est un complexe de Kan, on vérifie facilement qu'une unité $\xymatrix{A\ar[r]&{\bf R}\mathcal{H}\circ {\bf L}\nu (A)}$ est un isomorphisme de $\mathrm{Ho}_n\big(\simp_{0}\big)$ pour n'importe quelle choix de foncteurs dérivé ${\bf R}\mathcal{H}$ et ${\bf L}\nu$. D'un autre, si $(X,x_0)$ est un ensemble simplicial pointé tel que $X$ est un complexe de Kan, on vérifie à l'aide de la Proposition \ref{kankan} que $\xymatrix{{\bf L}\nu \circ {\bf R}\mathcal{H}(X,x_0)\ar[r]&(X,x_0)}$ induit un isomorphisme de groupes d'homotopie $\pi_i$ pour $0<i<\infty$; Donc si $X$ est supposé connexe $\xymatrix{{\bf L}\nu \circ {\bf R}\mathcal{H}(X,x_0)\ar[r]&(X,x_0)}$ est un isomorphisme de la catégorie $\mathrm{Ho}_n\big(\simp_\star\big)$.
\end{proof}

On aura besoin du critère suivant qui permet d'identifier les fibrations entre des objets fibrants de la catégorie de modèles $(\simp_{0},{\bf W}^{red}_{n}, {\bf mono},{\bf fib}_{n}^{red})$ (voir le Lemme 2.1 de \cite{stancu}). 

\begin{proposition}\label{identfibra}
Si $0\leq n\leq \infty$, soient $X$ et $Y$ des objets fibrants de la catégorie de modèles $(\simp_{0},{\bf W}^{red}_{n}, {\bf mono},{\bf fib}_{n}^{red})$ de la Proposition \ref{modred}, c'est-à-dire $X$ et $Y$ sont des ensembles simpliciaux réduits lesquels sont de complexes de Kan tels que $\pi_{m}(X)=0=\pi_m(Y)$ pour tout $m\geq n+1$; alors un morphisme $f\colon\xymatrix@C-10pt{X\ar[r]&Y}$ est une fibration de $(\simp_{0},{\bf W}^{red}_{n}, {\bf mono},{\bf fib}_{n}^{red})$ si et seulement si $f$ vérifie la propriété de relèvement à droite par rapport à l'ensemble de morphismes:
$$
J_0\,=\,\bigg\{ \quad \xymatrix@+35pt{\Lambda^{m,k}/{\bf sq}_{0}\Lambda^{m,k}\; \ar@{^(->}[r]^-{\alpha^{m-1}/{\bf sq}_{0}\alpha^{m-1}}& \, \Delta^{m}/{\bf sq}_{0}\Delta^m} \quad \bigg| \quad  m\geq 2\,,\;\, 0\leq k\leq m \quad \bigg\}\,.
$$
(note que $m\geq 2$ dans la définition de $J_0$).
\end{proposition}
\begin{proof} 
Vu que $(\simp_{0},{\bf W}^{red}_{n}, {\bf mono},{\bf fib}_{n}^{red})$ est une localisation de Bousfield à gauche de la catégorie de modèles $(\simp_{0},{\bf W}^{red}_{\infty}, {\bf mono},{\bf fib}_{\infty}^{red})$, d'après la Proposition 3.4.16 de \cite{hirschhorn} si $X$ et $Y$ sont des objets fibrants de $(\simp_{0},{\bf W}^{red}_{n}, {\bf mono},{\bf fib}_{n}^{red})$ un morphisme $f\colon\xymatrix@C-10pt{X\ar[r]&Y}$ est une fibration de $(\simp_{0},{\bf W}^{red}_{n}, {\bf mono},{\bf fib}_{n}^{red})$ si et seulement si $f$ est une fibration de $(\simp_{0},{\bf W}^{red}_{\infty}, {\bf mono},{\bf fib}_{\infty}^{red})$.

Donc il suffit de montrer que si $X$ et $Y$ sont des ensembles simpliciaux réduits lesquels sont de complexes de Kan, alors un morphisme $f\colon\xymatrix@C-10pt{X\ar[r]&Y}$ est une fibration de la catégorie de modèles $(\simp_{0},{\bf W}^{red}_{\infty}, {\bf mono},{\bf fib}_{\infty}^{red})$, si et seulement si $f$ vérifie la propriété de relèvement à droite par rapport à l'ensemble de morphismes:
$$
J_0\,=\,\bigg\{ \quad \xymatrix@+35pt{\Lambda^{m,k}/{\bf sq}_{0}\Lambda^{m,k}\; \ar@{^(->}[r]^-{\alpha^{m-1}/{\bf sq}_{0}\alpha^{m-1}}& \, \Delta^{m}/{\bf sq}_{0}\Delta^m} \quad \bigg| \quad  m\geq 2\,,\;\, 0\leq k\leq m \quad \bigg\}\,.
$$

Vérifions d'abord:

\begin{lemme}\label{unnec}
Si $m\geq 2$ et $0\leq k\leq m$, le morphisme d'ensembles simpliciaux réduits:
$$
\xymatrix@+35pt{\Lambda^{m,k}/{\bf sq}_{0}\Lambda^{m,k}\; \ar@{^(->}[r]^-{\alpha^{m-1}/{\bf sq}_{0}\alpha^{m-1}}& \, \Delta^{m}/{\bf sq}_{0}\Delta^m} 
$$
est un monomorphisme et une $\infty$-équivalence faible. 
\end{lemme}
\begin{proof}
On vérifie facilement que si $\xymatrix@C-10pt{A\ar[r]^-{j}&B}$ est un monomorphisme d'ensembles simpliciaux, alors $\xymatrix@+5pt{A/{\bf sq}_0A\ar[r]^-{j/{\bf sq}_0j}&B/{\bf sq}_0B}$ est un monomorphisme d'ensembles simpliciaux réduits.

D'un autre côté en prenant de carrés cocartésiens:
$$
\vcenter{\xymatrix{
\star\ar[r] &\Delta^{m}/{\bf sq}_0\Delta^m \\
{\bf sq}_0\Delta^m\ar[u] \ar[r] & \Delta^m\ar[u]
}}\qquad \text{et}\qquad \vcenter{\xymatrix{
\star\ar[r] &\Lambda^{m,k}/{\bf sq}_0\Lambda^{m,k} \\
{\bf sq}_0\Lambda^{m,k}\ar[u] \ar[r] & \Lambda^{m,k}\ar[u]
}}
$$ 
lesquels d'après le Lemme \ref{fibrabra} sont en fait de carrés homotopiquement cocartésiens de la catégorie de modèles $(\simp,{\bf W}_\infty, {\bf mono},{\bf fib}_\infty)$; on construit un cube dans $\simp$ dont toutes les faces sont de carrés commutatifs:
$$
\xymatrix@-10pt{&{\bf sq}_0\Lambda^{m,k}  \ar[ld]\ar[rr] \ar'[d][dd] && \Lambda^{m,k} \ar[ld] \ar[dd]^{\alpha^{m,k}}\\ 
                           \star\ar@{=}[dd]  \ar@<-2pt>[rr]&&  \Lambda^{m,k} / {\bf sq}_0\Lambda^{m,k}\ar[dd] &\\
                           &{\bf sq}_0 \Delta^{m}  \ar[ld]\ar'[r][rr] && \Delta^{m}\ar[ld] \\ 
                          \star\ar[rr]&&  \Delta^{m}/{\bf sq}_0\Delta^m&}
$$                       
                          
Donc vu que $\xymatrix{\Lambda^{m,k}\ar[r]^{\alpha^{m,k}}&\Delta^{m}}$ est une $n$-équivalence faible d'ensembles simpliciaux si $m\geq 1$ et $0\leq k \leq m$ et $\xymatrix@+5pt{{\bf sq}_0\Lambda^{m,k}\ar[r]^{{\bf sq}_0\alpha^{m,k}}&{\bf sq}_0\Delta^{m}}$ est un isomorphisme d'ensembles simpliciaux si $m\geq 2$ et $0\leq k \leq m$; on déduit du Lemme \ref{cubelele} que le morphisme:
$$
\xymatrix@+35pt{\Lambda^{m,k}/{\bf sq}_{0}\Lambda^{m,k}\; \ar@{^(->}[r]^-{\alpha^{m-1}/{\bf sq}_{0}\alpha^{m-1}}& \, \Delta^{m}/{\bf sq}_{0}\Delta^m} \,,
$$
est une $\infty$-équivalence faible si $m\geq 2$ et $0\leq k \leq m$.
\end{proof}

Il se suit du Lemme \ref{unnec} qu'on vient de montrer que si $f\colon\xymatrix@C-10pt{X\ar[r]&Y}$ est une fibration de la catégorie de modèles $(\simp_{0},{\bf W}^{red}_{\infty}, {\bf mono},{\bf fib}_{\infty}^{red})$, alors $f$ vérifie la propriété de relèvement à droite par rapport à $J_0$. 

Montrons:

\begin{lemme}\label{deuxnec}
Si $j\colon\xymatrix@C-8pt{A\ar[r]& B}$ est un morphisme d'ensembles simpliciaux réduits lequel est un monomorphisme et une $\infty$-équivalence faible, alors $j$ vérifie la propriété de relèvement à gauche par rapport à l'ensemble de morphismes d'ensembles simpliciaux réduits $f\colon\xymatrix@C-5pt{X\ar[r]& Y}$ tels que $X$ et $Y$ sont de complexes de Kan et $f$ vérifie la propriété de relèvement à droite par rapport à l'ensemble de morphismes:
$$
J_0 \, = \, \bigg\{ \quad \xymatrix@+35pt{\Lambda^{m,k}/{\bf sq}_{0}\Lambda^{m,k}\; \ar@{^(->}[r]^-{\alpha^{m-1}/{\bf sq}_{0}\alpha^{m-1}}& \, \Delta^{m}/{\bf sq}_{0}\Delta^m} \quad \bigg| \quad  m\geq 2\,,\;\, 0\leq k\leq m \quad \bigg\}\,.
$$
\end{lemme}
\begin{proof}
Supposons qu'on s'est donné un carré commutatif d'ensembles simpliciaux réduits:
\begin{equation}\label{carrefinal}
\xymatrix{A\ar[d]_-{j}\ar[r]^-{\varphi} &  X \ar[d]^-{f}\\ 
B\ar[r]_{\psi} & Y }
\end{equation}
où $j$ est un monomorphisme et une $\infty$-équivalence faible, $X$ et $Y$ sont de complexes de Kan et $f$ vérifie la propriété de relèvement à droite par rapport à l'ensemble de morphismes:
$$
J_0\,=\,\bigg\{ \quad \xymatrix@+35pt{\Lambda^{m,k}/{\bf sq}_{0}\Lambda^{m,k}\; \ar@{^(->}[r]^-{\alpha^{m-1}/{\bf sq}_{0}\alpha^{m-1}}& \, \Delta^{m}/{\bf sq}_{0}\Delta^m} \quad \bigg| \quad  m\geq 2\,,\;\, 0\leq k\leq m \quad \bigg\}\,.
$$

On doit montrer que le carré \eqref{carrefinal} admet un relèvement, c'est-à-dire qu'il existe un morphisme d'ensembles simpliciaux réduits $F\colon\xymatrix@C-8pt{B\ar[r]&X}$ tel que $f\circ F=\psi$ et $F\circ j = \varphi$.

Dans ce but commençons pour prendre à l'aide de l'argument du petit objet une factorisation dans $\simp_0$ du morphisme $\psi$:
$$
\xymatrix{B  \ar@/^10pt/[rr]^-{\psi} \ar[r]_-{\psi_1}  & B' \ar[r]_-{\psi_2}  &Y\,,}
$$
où $\psi_2$ vérifie la propriété de relèvement à droite par rapport à l'ensemble de morphismes $J_0$ ci-dessus et $\psi_1$ est un morphisme dans ${\bf cell}\big(J_0\big)$.

En suit considérons un carré cartésien dans $\simp_0$: 
$$
\xymatrix{ B' \underset{Y}{\times} X \ar[d]_-{p_1} \ar[r]^-{p_2} & X \ar[d]^-{f}\\
B' \ar[r]_-{\psi_2} & Y \,,}
$$
et si $\xi\colon\xymatrix{A\ar[r] & B'\underset{Y}{\times} X}$ est le seul morphisme d'ensembles simpliciaux réduits tel que $p_2\circ\xi = \varphi$ et $p_1\circ\xi = \psi_1\circ j$, prenons une factorisation dans $\simp_0$:
$$
\xymatrix{A  \ar@/^10pt/[rr]^-{\xi} \ar[r]_-{\xi_1}  & A' \ar[r]_-{\xi_2}  & B' \underset{Y}{\times} X\,,}
$$ 
où $\psi_2$ vérifie la propriété de relèvement à droite par rapport à $J_0$ et $\psi_1$ est un morphisme dans ${\bf cell}\big(J_0\big)$.

Notons qu'on a un carré commutatif:
\begin{equation}\label{carrefinal2}
\xymatrix{
A\ar[d]_-{j} \ar[r]^-{\xi_1} & A'  \ar[d]^-{p_1\circ\xi_2} \\
B\ar[r]_-{\psi_1} & B'\,.}
\end{equation}
parce que $(p_1\circ \xi_2)\circ\xi_1=p_1\circ\xi=\psi_1\circ j$. 

En plus on vérifie que \eqref{carrefinal} admet un relèvement s'il est du même pour le carré \eqref{carrefinal2}. En effet, si $G\colon\xymatrix@C-8pt{B\ar[r]&A'}$ satisfait que $(p_1\circ\xi_2)\circ G=\psi_1$ et $G\circ j = \xi_1$, alors le morphisme $F\colon\xymatrix@C-8pt{B\ar[r]&X}$ défini par le composé $F=p_2\circ\xi_2\circ G$ vérifie que:
$$
f\circ F= f\circ p_2\circ\xi_2\circ G = \psi_2 \circ p_1\circ\xi_2\circ G = \psi_2 \circ \psi_1 = \psi
$$ 
$$\text{et}$$ 
$$
F\circ j = p_2\circ\xi_2\circ G\circ j = p_2\circ\xi_2\circ \xi_1= p_2\circ\xi = \varphi\,.
$$

Enfin, vu que $j\colon\xymatrix@C-8pt{A\ar[r]& B}$ est un monomorphisme d'ensembles simpliciaux réduits, pour montrer l'existence d'un relèvement du carré \eqref{carrefinal2} il suffit de montrer les points suivants:
\begin{enumerate}
\item[a)] $p_1\circ\xi_2$ est une $\infty$-équivalence faible.
\item[b)] $p_1\circ\xi_2$ est une fibration de Kan.
\end{enumerate}

Le point a) est déduit de la propriété de deux-sur-trois des $\infty$-équivalences faibles vu que $(p_1\circ \xi_2)\circ\xi_1 = \psi_1 \circ j$; et parce que $j$ est une $\infty$-équivalences faibles par hypothèse et les morphismes $\xi_1$ et $\psi_1$ sont aussi des $\infty$-équivalences faibles en tant que des éléments de ${\bf cell}(J_0)$ (voir le Lemme \ref{unnec}).

Pour montrer b) remarquons d'abord:
\begin{lemme}\label{fafalele}
Si $Z$ est un ensemble simplicial réduit, alors $\xymatrix@-10pt{Z\ar[r]&\star}$ vérifie la propriété de relèvement à droite par rapport à l'ensemble:
$$
\bigg\{ \quad \xymatrix@+35pt{\Lambda^{m,k}/{\bf sq}_{0}\Lambda^{m,k}\; \ar@{^(->}[r]^-{\alpha^{m-1}/{\bf sq}_{0}\alpha^{m-1}}& \, \Delta^{m}/{\bf sq}_{0}\Delta^m} \quad \bigg| \quad  m\geq 2\,,\;\, 0\leq k\leq m \quad \bigg\}\,.
$$
si et seulement $Z$ est un complexe de Kan (son image par le foncteur canonique $\xymatrix@C-3pt{\simp_0\ar@{^(->}[r]^-{\nu}&\simp}$).
\end{lemme}
\begin{proof}
Vu qu'on a une adjonction:
$$
\simp_{0}
\xymatrix@C+15pt{
\phantom{a}\ar@<+10pt>@{}[r]|{\perp}\ar@{^(->}[r]|-{\,\nu\,}&\phantom{a}
\ar@<-7pt>@/_12pt/[l]_-{\mathcal{F} = \,\cdot\, \big/ \mathbf{sq}_0(\,\cdot\,)}}
\simp\,,
$$
où le foncteur $\nu$ est normalement supprimé; il se suit qu'un ensemble simplicial réduit $X$ est un complexe de Kan si et seulement si le morphisme $\xymatrix@-10pt{X\ar[r]&\star}$ vérifie la propriété de relèvement à droite par rapport à l'ensemble:
$$
\bigg\{ \quad \xymatrix@+35pt{\Lambda^{m,k}/{\bf sq}_{0}\Lambda^{m,k}\; \ar@{^(->}[r]^-{\alpha^{m-1}/{\bf sq}_{0}\alpha^{m-1}}& \, \Delta^{m}/{\bf sq}_{0}\Delta^m} \quad \bigg| \quad  m\geq 1\,,\;\, 0\leq k\leq m \quad \bigg\}\,.
$$

Le Lemme se suit su fait que pour tout ensemble simplicial réduit $X$, le morphisme $\xymatrix@-10pt{X\ar[r]&\star}$ vérifie la propriété de relèvement à droite par rapport au morphisme $\xymatrix@-7pt{\star\ar@{^(->}[r]&\Delta^{1}/{\bf sq}_{0}\Delta^1 \, = \, \mathbb{S}^1}$.
\end{proof}

Montrons aussi:

\begin{lemme}\label{didifici}
Soit $Z$ et $W$ des ensembles simpliciaux réduits lesquels sont des complexes de Kan. Si $T\colon\xymatrix@C-8pt{Z\ar[r]&W}$ est un morphisme vérifiant la propriété de relèvement à droite par rapport à l'ensemble:
$$
J_0\,=\,
\bigg\{ \quad \xymatrix@+35pt{\Lambda^{m,k}/{\bf sq}_{0}\Lambda^{m,k}\; \ar@{^(->}[r]^-{\alpha^{m-1,k}/{\bf sq}_{0}\alpha^{m-1}}& \, \Delta^{m}/{\bf sq}_{0}\Delta^m} \quad \bigg| \quad  m\geq 2\,,\;\, 0\leq k\leq m \quad \bigg\}
$$ 
et tel que $T^\star\colon\xymatrix@C-8pt{\pi_1(Z)\ar[r]&\pi_1(W)}$ est une fonction surjective,  alors $T$ est une fibration de Kan (l'image de $f$ par le foncteur canonique $\xymatrix@C-3pt{\simp_0\ar@{^(->}[r]^-{\nu}&\simp}$). 
\end{lemme}
\begin{proof} Il nous reste à montrer que $T$ vérifie la propriété de relèvement à droite par rapport au morphisme:
$$
\xymatrix@+5pt{\star\; \ar@{^(->}[r]& \, \Delta^{1}/{\bf sq}_{0}\Delta^1}\,;
$$
c'est-à-dire que la fonction $T_1\colon\xymatrix@C-8pt{Z_1\ar[r]&W_1}$ est une fonction surjective.

Soit $w\in W_1$ un $1$-simplexe de l'ensemble simplicial réduit $W$. Vu que $T^\star\colon\xymatrix@C-8pt{\pi_1(Z)\ar[r]&\pi_1(W)}$ est une fonction surjective, il se suit de la Proposition \ref{kankan} qu'il existe un $1$-simplexe $z\in Z_1$ de $Z$ et un $2$-simplexe $\eta\in W_2$ de $W$ tels que $d_2(\eta)=T_1(z)$ et $d_1(\eta)=w$ et $d_0(\eta)=\star$.  

Considérons le morphisme $\eta'\colon\xymatrix{\Lambda^{2,1}/{\bf sq}_{0}\Lambda^{2,1}\ar[r] & Z}$ définit en imposant que: 
$$
\xymatrix{\Delta^1 \ar@/^15pt/[rrr]^-{z}\ar[r]_-{\delta_2} & \Lambda^{2,1} \ar[r]_-{\text{quotient}}& \Lambda^{2,1}/{\bf sq}_{0}\Lambda^{2,1}\ar[r]_-{\eta'} & Z}\qquad \text{et} \qquad 
\xymatrix{\Delta^1 \ar@/^15pt/[rrr]^-{\star}\ar[r]_-{\delta_0} & \Lambda^{2,1} \ar[r]_-{\text{quotient}}& \Lambda^{2,1}/{\bf sq}_{0}\Lambda^{2,1}\ar[r]_-{\eta'} & Z}\,;
$$
et remarquons qu'on a un carré commutatif:
\begin{equation}\label{relevecarr}
\xymatrix{
\Lambda^{2,1}/{\bf sq}_{0}\Lambda^{2,1} \ar[r]^-{\eta'} \ar[d]_-{\alpha^{2,1}/{\bf sq}_{0}\alpha^{2,1}} & Z \ar[d]^-{T}\\ 
\Delta^{2}/{\bf sq}_{0}\Delta^2 \ar[r]_-{\eta} & W}
\end{equation}

Si $\xymatrix{\Delta^{2}/{\bf sq}_{0}\Delta^2 \ar[r]^-{\xi} & Z}$ est un relèvement du carré \eqref{relevecarr}, on constate que:
$$T_1(d_1\xi) = d_1(T_2\xi) = d_1(\eta) = w\,.$$ 

Donc $T_1\colon\xymatrix@C-8pt{Z_1\ar[r]&W_1}$ est une fonction surjective.
\end{proof}

Remarquons enfin que d'après le Lemme \ref{fafalele} les ensembles simpliciaux $A'$ et $B'$ du diagramme \eqref{carrefinal2} sont des complexes de Kan. En effet, les morphismes $p_1$, $\xi_2$ et $\psi_2$ vérifient la propriété de relèvement à droite par rapport à $J_0$, et par hypothèse $X$ et $Y$ sont des complexes de Kan. Donc les morphismes $\xymatrix@C-12pt{A'\ar[r]&\star}$ et $\xymatrix@C-12pt{B'\ar[r]&\star}$ vérifient la propriété de relèvement à droite par rapport à $J_0$.

Donc, vu qu'on a déjà montré que le morphisme $p_1\circ\xi_2$ est une $\infty$-équivalence faible, on déduit du Lemme \ref{didifici} le point b) désiré, c'est-a-dire que $p_1\circ\xi_2$ est une fibration de Kan. 
\end{proof} 

Du Lemme \ref{deuxnec} qu'on vient de montrer on a que si $X$ et $Y$ sont des ensembles simpliciaux réduits lesquels sont de complexes de Kan et $f\colon\xymatrix@C-10pt{X\ar[r]&Y}$ est un morphismes vérifient la propriété de relèvement à droite par rapport à $J_0$, alors $f$ est une fibration de la catégorie de modèles $(\simp_{0},{\bf W}^{red}_{\infty}, {\bf mono},{\bf fib}_{\infty}^{red})$.
\end{proof}

\renewcommand{\thesubsection}{\S\thesection.\arabic{subsection}}
\subsection{}\;
\renewcommand{\thesubsection}{\thesection.\arabic{subsection}}

Dans le présent paragraphe on va montrer que si $1\leq n < \infty$ la catégorie de modèles simpliciale pointée $\big(\simp_{0},{\bf W}^{red}_{n}, {\bf mono},{\bf fib}_{n}^{red},\underline{\mathrm{Hom}}_{\simp_0}\big)$ de la Proposition \ref{modred} est une catégorie de modèles simpliciale pointée d'ordre $n-1$. 

D'après le Théorème \ref{cassi1} il suffit de construire un foncteur espace de lacets de la catégorie de modèles pointée $(\simp_{0},{\bf W}^{red}_{n}, {\bf mono},{\bf fib}_{n}^{red})$:
$$
\xymatrix{\simp_0\big[{\bf W}_n^{red}\big] \ar[r]^{\Omega} & \simp_0\big[{\bf W}_n^{red}\big]\,,}
$$
et noter que son $n$-ième itération $\Omega^{n}$ est le foncteur ponctuel.

Pour  cela considérons en premier lieu la suite de foncteurs:
\begin{equation}\label{suitedfo}
\xymatrix@C-8pt@R=1pt{
\;{\bf\Delta}\;\ar[r]^-{\mathrm{i}_{1,0}}&\;{\bf\Delta}\times  {\bf\Delta} \;\ar[r]^-\vartheta& \;{\bf\Delta}\;\\ 
\quad\text{\scriptsize{$[n]$}}\quad\ar@{|->}[r]&\quad\text{\scriptsize{$\big([n],[0]\big)$}}\quad\ar@{|->}[r]&\quad\text{\scriptsize{$[n+1]$}}\;\;}\,,
\end{equation}
où $\mathrm{i}_{1,0}$ est le foncteur d'inclusion de \eqref{iinclusss} et $\vartheta$ est la restriction à la catégorie des simplexes du produit $\otimes$ de la catégorie augmentée des simplexes ${\bf\Delta}_{+}$ défini dans \eqref{definidelta}.

On est intéressé au foncteur: 
$$
\xymatrix@C+15pt{
\simp \ar[r]^-{\mathbb{D}(\,-\,)} &\simp}\,
$$
défini par le composé des foncteurs induisent des foncteurs \eqref{suitedfo}:
$$
\xymatrix@C+15pt{
\simp \ar[r]^-{\vartheta^*} &\ssimp  \ar[r]^-{(\mathrm{i}_{1,0})^*}    &\simp }\,.
$$

Explicitement si $X$ est un ensemble simplicial, on trouve que $\mathbb{D}(X)_{n}=X_{n+1}$ et que les morphismes faces et dégénérescences de $\mathbb{D}(X)$ sont donnés par:
\begin{align*}
d_{i} : \mathbb{D}(X)_{n}=&\xymatrix@C+10pt{X_{n+1} \ar[r]^-{d_{i}}&X_{n}}=\mathbb{D}(X)_{n-1} \, \qquad \quad 0\leq i\leq n\\
&\qquad\qquad\quad\text{et}\\
s_{j}: \mathbb{D}(X)_{n-1}=&\xymatrix@C+10pt{X_{n}\ar[r]^-{s_{j}}&X_{n+1}}=\mathbb{D}(X)_{n}\, \qquad \quad 0\leq j\leq n-1\,.
\end{align*}

Considérons les transformations naturelles:
$$
\vcenter{\xymatrix@C+18pt{
\simp\ar@/^26pt/[rrr]^-{\text{\textit{Foncteur identité}}}\ar[rrr]|-{\mathbb{D}(-)}\ar@/_26pt/[rrr] _{{\bf sq}_0(\,-\,)}
&\rrtwocell<\omit>{<2> {\phantom{a}\alpha}}& \lltwocell<\omit>{<-2> {\beta\phantom{a}}}\lltwocell<\omit>{<2> {\varrho\phantom{a}}}& 
\simp }}\,,
$$
où les morphismes $\xymatrix@C+5pt{X  & \mathbb{D}(X) \ar@<-3pt>[r]_{\alpha_X}\ar[l]_-{\varrho_X}&{\bf sq}_0(X)\ar@<-5pt>[l]_{\beta_X}}$ sont définies par les fonctions:
$$
\xymatrix@C+15pt{X_n  & X_{n+1} \ar@<-5pt>[r]_{\underbrace{d_{0}\dots d_{0}}_{n+1}}\ar[l]_-{d_{n+1}}&X_0\ar@<-5pt>[l]_{\overbrace{s_{0}\dots s_0}^{n+1}}}\;,\quad\qquad\text{si}\quad \text{$n\geq 0$}
$$
pour tout ensemble simplicial $X$.

On vérifie qu'effectivement ces sont de transformations naturelles parce que d'après les identités \eqref{relsim} on a pour $0\leq i\leq n$ et $0\leq j \leq n-1$ que:
\begin{equation}
\begin{split}
d_n\circ d_i = d_i\circ d_{n+1}   & \qquad\quad \text{et} \qquad\quad  d_{n+1}\circ s_j = d_n\circ s_j\;,\\
\underbrace{d_0\circ\dots\circ d_0}_{n}\circ d_i = \underbrace{d_0\circ\dots\circ d_{0}}_{n+1} 
& \qquad\quad \text{et} \qquad\quad
\underbrace{d_0\circ \dots\circ d_{0}}_{n+1}\circ s_j = \underbrace{d_0\circ\dots\circ d_0}_n \;,\\
\underbrace{s_{0}\circ\dots\circ s_{0}}_{n} = d_i\circ \underbrace{s_{0}\circ\dots\circ s_{0}}_{n+1}
& \qquad\quad \text{et} \qquad\quad
s_j \circ \underbrace{s_{0}\circ\dots\circ s_{0}}_{n} = \underbrace{s_{0}\circ\dots\circ s_{0}}_{n+1}\,.
\end{split}
\end{equation}

\begin{lemme}\label{retractpd}
Si $X$ est un ensemble simplicial, les morphismes $\xymatrix@C+5pt{\mathbb{D}(X) \ar@<-3pt>[r]_{\alpha_X}&{\bf sq}_0(X)\ar@<-5pt>[l]_{\beta_X}}$ vérifient que $\alpha_X\circ\beta_X = \mathrm{id}_{{\bf sq}_0(X)}$ et il existe une homotopie $H:\beta_X\circ\alpha_X \simeq \mathrm{id}_{\mathbb{D}(X)}$; autrement dit ${\bf sq}_0(X)$ est un rétracte par déformation de $\mathbb{D}(X)$. En particulier, si $X$ est un ensemble simplicial réduit l'ensemble simplicial pointé $\mathbb{D}(X)$ de point base $\xymatrix{\star\ar[r]^-{\beta_X}&\mathbb{D}(X)}$ est contractile.
\end{lemme}
\begin{proof}
On vérifie que $\alpha_X\circ\beta_X = \mathrm{id}_{{\bf sq}_0(X)}$, c'est-à-dire que:
$$
\big( \alpha_X \big)_n \circ \big( \beta_X \big)_n \, = \, \underbrace{d_{0}\circ \dots \circ d_0}_{n+1}\circ \underbrace{s_0 \circ \dots \circ s_0}_{n+1} \, = \, \mathrm{id}_{X_0} \quad\qquad \text{si}\quad\;\, n\geq 0\,,
$$
parce que $d_0\circ s_{0} \, = \, \mathrm{id}_{X_{k}}$ si $k\geq 0$.

Rappelons par ailleurs:

\begin{lemme}\label{homotocombina}
Se donner une homotopie  $H\colon\xymatrix{Z\times\Delta^1\ar[r] & W}$ tel que $H:  f \simeq g$ équivaut à se donner une famille de fonctions:
$$
\xymatrix@C+5pt{Z_n\ar[r]^-{H(t)_n} & W_n} \qquad\quad \text{pour} \quad n\geq 0 \quad \text{et} \quad 0\leq t\leq n+1\,,
$$
vérifiant:
$$
d_i\circ H(t)_n \, = \, 
\begin{cases}
H(t)_{n-1}\circ d_i   & \quad \text{si \quad $0\leq t\leq i\leq n$} \\
H(t-1)_{n-1}\circ d_i  & \quad \text{si \quad $0\leq i < t \leq n+1$\,,} 
\end{cases}
$$
$$
s_j\circ H(t)_n \, = \, 
\begin{cases}
H(t)_{n+1}\circ s_j   & \quad \text{si \quad $0\leq t\leq j\leq n$} \\
H(t+1)_{n+1}\circ s_j  & \quad \text{si \quad $0\leq j < t\leq n+1$\,,} 
\end{cases}
$$
$$
H(0)_{n} \, = \, f_n  \qquad\text{et}\qquad  H(n+1)_{n} \, = \, g_n\,.
$$
\end{lemme}
\begin{proof}
Si $n\geq 0$ remarquons qu'on a une bijection:
\begin{equation}\label{n1homotobije}
\xymatrix@C+20pt@R=3pt{
\varphi   \ar@{}[r]|-{\longmapsto} &  t_{\varphi}  \\
\mathrm{Hom}_{\Delta}\big([n],[1]\big) \, = \, \Delta^1_n \quad \ar@{<->}[r] & \quad  \big\{\, t \,\big| \, 0\leq t \leq n+1\,\big\} \\
\varphi_t  \ar@{}[r]|(.55){\longleftarrow\!\!\shortmid} &  t  }
\end{equation}
où:
$$
t_\varphi=\begin{cases} \text{min}\{\, a \,|\, \varphi(a)=1 \,\}  \qquad &\text{si} \;\, \varphi \not\equiv 0\\ n+1 & \text{si} \;\, \varphi \equiv 0\end{cases}
\qquad\quad \text{et} \qquad\quad 
\varphi_t(a)=\begin{cases} 0 \qquad & \text{si} \;\,  0\leq a < t \\ 1 & \text{si} \;\,  t\leq a \leq n  \,. \end{cases}
$$

On déduit que les formules:
$$
H(t)_n(x) \, = \, H_n(x,\varphi_t) \qquad \text{et} \qquad H_n(x,\varphi) \, = \, H(t_\varphi)_n(x)
$$ 
définissent une bijection entre l'ensemble des familles de fonctions:
\begin{equation}\label{famille1sanstoi}
\Big\{ \, \xymatrix{Z_n\times \Delta^1_n \ar[r]^-{H_n}   & W_n} \, \Big\}_{n\geq 0}
\end{equation}
et l'ensemble des familles de fonctions:
\begin{equation}\label{famille2sanstoi}
\Big\{ \, \xymatrix{Z_n \ar[r]^-{H(t)_n}  & W_n} \, \Big\}_{n\geq 0,\,0\leq t\leq n+1}\,.
\end{equation}

En plus remarquons que si $0\leq j \leq n$ on a que:
$$
d_i\circ H(t)_n (x) = d_i\circ H_n(x,\varphi^t) = H_{n-1}(d_ix,\varphi_t\circ \delta_i) = H\big(t_{\varphi_t\circ\delta_i}\big)_{n-1}\circ d_i (x)
$$
où:
$$
s(i,t) \, = \, \begin{cases}  t \qquad & \text{si} \; 0\leq t \leq i\leq n \\ t-1 \qquad & \text{si} \; 0\leq i < t \leq n+1\,;  \end{cases}
$$
et si $0\leq j \leq n$ on a que:
$$
s_j\circ H(t)_n(x) = s_j\circ H_n(x,\varphi_t) = H_{n+1}(s_jx,\varphi_t\circ\sigma_j) = H\big(t_{\varphi_t\circ\sigma_j}\big)_{n+1}\circ s_j (x)
$$
où:
$$
s'(j,t) \, = \, \begin{cases}  t \qquad & \text{si} \;\, 0\leq t \leq j\leq n \\ t+1 \qquad & \text{si} \,\; 0\leq j < t \leq n+1 \,. \end{cases}
$$

D'un autre on a que:
$$
H_n\big(x,\delta_0\circ\underbrace{\sigma_0\circ\dots\circ\sigma_0}_{n}\big) \, = \, H_n(0)(x)
\qquad\text{et}\qquad
H_n\big(x,\delta_1\circ\underbrace{\sigma_0\circ\dots\circ\sigma_0}_{n}\big) \, = \, H_n(n+1)(x)
$$
pour tout $n$-simplexe $x$ de $Z$.
\end{proof}

On utilise le Lemme \ref{homotocombina} pour définir une homotopie $H:\beta_X\circ\alpha_X \simeq \mathrm{id}_{\mathbb{D}(X)}$:

Si $n\geq 0$ et $0\leq t\leq n+1$ on définie la fonction $H(t)_n \colon\xymatrix{X_{n+1}\ar[r] & X_{n+1}}$ par la formule:
$$
H(t)_n  \, = \, \underbrace{s_{n}\circ\dots\circ s_t}_{n+1-t} \circ \underbrace{d_t\circ\dots\circ d_{n}}_{n+1-t} \,.
$$

Remarquons alors que par définition $H(n+1)_n  = \mathrm{id}_{\mathbb{D}(X)_{n}}$ et qu'en plus:
$$
H(0)_n  \, = \, \underbrace{s_{n}\circ\dots\circ s_0}_{n+1} \circ \underbrace{d_0\circ\dots\circ d_{n}}_{n+1} \, = \, \underbrace{s_{0}\circ\dots\circ s_0}_{n+1} \circ \underbrace{d_0\circ\dots\circ d_{0}}_{n+1} \, = \, \big(\beta_X\big)_n\circ\big(\alpha_X\big)_n
$$
parce que   $s_k\circ s_l = s_{l+1}\circ s_k$ si $k\leq l$ et $d_k\circ d_l = d_{l-1}\circ d_k$ si $k<l$.

D'un autre côté si $0\leq t\leq i\leq n$ on a que:
\begin{align*}
d_i\circ H(t)_n 
\, =& \, d_i\circ \big(\underbrace{s_{n}\circ\dots\circ s_t}_{n+1-t} \circ \underbrace{d_t\circ\dots\circ d_{n}}_{n+1-t}\big) \\
\, =& \, \big(\underbrace{s_{n-1}\circ\dots s_i\circ(d_i\circ s_i)\circ s_{i-1}\dots\circ s_t}_{n+1-t+1}\big)\circ\big(\underbrace{d_t\circ\dots\circ d_{i-1}\circ d_i \circ d_{i+1}\circ\dots\circ d_{n}}_{n+1-t}\big)\\
\, =& \, \big(\underbrace{s_{n-1}\circ\dots\circ s_t}_{(n-1)+1-t}\big)\circ\big(\underbrace{d_t\circ\dots\circ d_{n-1}}_{(n-1)+1-t}\big)\circ d_i\\
\, =& \, H(t)_{n-1}\circ d_i\,
\end{align*}
et si $0\leq i < t \leq n+1$ on a que:
\begin{align*}
d_i\circ H(t)_n 
\, =& \, d_i\circ \big(\underbrace{s_{n}\circ\dots\circ s_t}_{n+1-t} \circ \underbrace{d_t\circ\dots\circ d_{n}}_{n+1-t}\big) \\
\, =& \,\big(\underbrace{s_{n-1}\circ\dots\circ s_{t-1}}_{n+1-t} \circ  d_i\circ  \underbrace{d_t\circ\dots\circ d_{n}}_{n+1-t}\big) \\
\, =& \,\big(\underbrace{s_{n-1}\circ\dots\circ s_{t-1}}_{(n-1)+1-(t-1)} \underbrace{d_{t-1}\circ\dots\circ d_{n-1}}_{(n-1)+1-(t-1)}\big) \circ d_i \\
\, =& \, H(t-1)_{n-1}\circ d_i\,;
\end{align*}
parce que $d_i\circ s_k = \begin{cases} s_{k-1}\circ d_i &\text{si\; $i<k$}\\\mathrm{id} & \text{si \; $i=k$}\end{cases}$ et $d_i\circ d_l = d_{l-1}\circ d_i$ si $i<l$.

Finalement notons que si $0\leq t\leq j \leq n$:
\begin{align*}
s_j\circ H(t)_n 
\, =& \, s_j\circ \big(\underbrace{s_{n}\circ\dots\circ s_t}_{n+1-t} \circ \underbrace{d_t\circ\dots\circ d_{n}}_{n+1-t}\big) \\
\, =& \, \big(\underbrace{s_{n+1}\circ\dots\circ s_t}_{(n+1)+1-t} \big)\circ \Big( (d_t\circ\dots\circ d_j) \circ (d_{j+1}\circ s_j)\circ (d_{j+1}\circ\dots\circ d_{n})\Big) \\
\, =& \,  \big(\underbrace{s_{n+1}\circ\dots\circ s_t}_{(n+1)+1-t} \circ \underbrace{d_t\circ\dots\circ d_{n+1}}_{(n+1)+1-t}\big) \circ s_j\\
\, =& \,  H(t)_{n+1}\circ s_j 
\end{align*}
et si $0\leq j < t \leq n+1$ on a que:
\begin{align*}
s_j\circ H(t)_n 
\, =& \, s_j\circ \big(\underbrace{s_{n}\circ\dots\circ s_t}_{n+1-t} \circ \underbrace{d_t\circ\dots\circ d_{n}}_{n+1-t}\big) \\
\, =& \, \big(\underbrace{s_{n+1}\circ\dots\circ s_{t+1}}_{n+1-t} \big)\circ s_j\circ \big(\underbrace{d_{t}\circ\dots\circ d_{n}}_{n+1-t}\big) \\
\, =& \, \big(\underbrace{s_{n+1}\circ\dots\circ s_{t+1}}_{(n+1)+1-(t+1)} \big)\big(\underbrace{d_{t+1}\circ\dots\circ d_{n+1}}_{(n+1)+1-(t+1)}\big)\circ s_j \\
\, =& \,  H(t+1)_{n+1}\circ s_j \,,
\end{align*}
parce que  $s_j\circ s_l = s_{l+1}\circ s_j$ si $j\leq l$, $d_{j+1}\circ s_j = \mathrm{id}$ et $s_{j} \circ d_{k} = d_{k+1}\circ s_{j}$ si $j<k$,
\end{proof}

Rappelons aussi:

\begin{lemme}\label{retractpd2}
Si $X$ est un ensemble simplicial lequel est un complexe de Kan, alors le morphisme d'ensembles simpliciaux $\xymatrix@C-8pt{\mathbb{D}(X) \ar[r]^-{\varrho_X}&X}$ est une fibration de Kan. 
\end{lemme}

Si $X$ est un ensemble simplicial pointé de point base $\xymatrix@C-12pt{\star\ar[r]^-{x_0}&X}$, on définit \emph{le espace de lacets simplicial de $X$ sur $x_0$} comme le sous-ensemble simplicial pointé $\mathbb{\Omega}_{x_0}(X)$ de $\mathbb{D}(X)$ où: 
\begin{equation}\label{lacetssimpisi}
\text{$\mathbb{\Omega}_{x_0}(X)_n= \Big\{\; \alpha\in \mathbb{D}(X)_n = X_{n+1}\;\;\Big|\;\; d_{n+1}(\alpha)=x_0 \;\Big\}\qquad$ si $n\geq 0$}\,.
\end{equation}

En particulier, on a un carré cartésien dans la catégorie $\simp$ ou la catégorie $\simp_\star$:
\begin{equation}\label{omegasimp}
\xymatrix{
\mathbb{\Omega}_{x_0}(X) \ar[r] \ar[d] & \mathbb{D}(X) \ar[d]^-{\rho_X}\\ 
\star \ar[r]_-{x_0} & X\,.}
\end{equation}

\begin{lemme}\label{minusone}
Soit $0 \leq n \leq \infty$. Si $X$ est un ensemble simplicial $n$-fibrant réduit c'est-à-dire $X$ est un objet fibrant de la catégorie de modèles $(\simp_{0},{\bf W}^{red}_{n}, {\bf mono},{\bf fib}_{n}^{red})$ de la Proposition \ref{modred}, ou de façon explicite $X$ est un complexe de Kan tel que $X_0=\star$ et $\pi_i(X)=0$ pour $i \geq n+1$, alors:
\begin{enumerate}
\item Le carré cartésien \eqref{omegasimp} est un carré homotopiquement cartésien dans la catégorie de modèles $(\simp_{\star},\pi^{-1}{\bf W}_{n}, {\bf mono},\pi^{-1}{\bf fib}_{n})$.
\item Le seul morphisme $\xymatrix@C-8pt{\mathbb{D}(X) \ar[r]^-{\alpha_X}&\star}$ est une $\infty$-équivalence faible.
\item Si $0 < n< \infty$, l'ensemble simplicial $\mathbb{\Omega}_\star(X)$ est un objet fibrant de la catégorie de modèles $(\simp_{\star},\pi^{-1}{\bf W}_{n-1}, {\bf mono},\pi^{-1}{\bf fib}_{n-1})$. De plus si $n=0$ le seul morphisme $\xymatrix@C-8pt{\mathbb{\Omega}_\star(X)\ar[r]^-{\alpha_X}&\star}$ est une $\infty$-équivalence faible.
\end{enumerate}
\end{lemme}
\begin{proof}
Étant donné $0\leq n\leq \infty$, supposons que $X$ est un complexe de Kan tel que $X_0=\star$ et $\pi_i(X)=0$ pour $i \geq n+1$.

Il se suit des Lemmes \ref{retractpd} et \ref{retractpd2} que le morphisme $\xymatrix@C-8pt{\mathbb{D}(X) \ar[r]^-{\varrho_X}&X}$ est une fibration de Kan et que $\xymatrix@C-8pt{\mathbb{D}(X) \ar[r]^-{\alpha_X}& {\bf sq}_0(X) = \star}$ une $\infty$-équivalence faible. En particulier les ensembles simpliciaux $X$ et $\mathbb{D}(X)$ sont des objets fibrants de la catégorie de modèles $(\simp_\star,\pi^{-1}{\bf W}_{n}, {\bf mono},\pi^{-1}{\bf fib}_{n})$ si on les considères comme des ensembles simpliciaux pointés par les $0$-simplexes $\star \in \{\star\} =X_0$ et $\star = s_0(\star)\in \mathbb{D}(X)_0=X_1$ respectivement. Donc la fibration de Kan $\xymatrix@C-8pt{\mathbb{D}(X) \ar[r]^-{\varrho_X}&X}$ est une fibration de la catégorie de modèles $(\simp_{\star},\pi^{-1}{\bf W}_{n}, {\bf mono},\pi^{-1}{\bf fib}_{n})$, et vu que \eqref{omegasimp} est un carré cartésien de $\simp_\star$, l'espace de lacet simplicial $\mathbb{\Omega}_\star(X)$ est un objet fibrant de $(\simp_{\star},\pi^{-1}{\bf W}_{n}, {\bf mono},\pi^{-1}{\bf fib}_{n})$.

D'après le Lemme \ref{fibrabra} le carré \eqref{omegasimp} est en fait un carré homotopiquement cartésien de la catégorie de modèles $(\simp_{\star},\pi^{-1}{\bf W}_{n}, {\bf mono},\pi^{-1}{\bf fib}_{n})$.

Remarquons d'un autre qu'un objet fibrant $Z$ de $(\simp_{\star},\pi^{-1}{\bf W}_{n}, {\bf mono},\pi^{-1}{\bf fib}_{n})$ est fibrant dans la catégorie de modèles $(\simp_{\star},\pi^{-1}{\bf W}_{n-1}, {\bf mono},\pi^{-1}{\bf fib}_{n-1})$ si et seulement si $\pi_{n}(Z)=0$. 

Il y a plusieurs arguments pour montrer que $\pi_{n}\big(\mathbb{\Omega}_\star(X)\big) \cong \pi_{n+1}\big(X\big)=0$. Par exemple, si $X$ est un complexe de Kan c'est une conséquence de la Proposition \ref{kankan} vu que:
$$
\Big\{ \alpha\in\mathbb{\Omega}_\star(X)_n \; \Big| \;  d_i(\alpha)=\star\quad 0\leq i\leq n\Big\} \; = \; 
\Big\{ \alpha\in X_{n+1} \; \Big|  \; d_i(\alpha)=\star\quad 0\leq i\leq n+1\Big\}
$$
$$\text{et}$$
$$
\Big\{ \sigma\in\mathbb{\Omega}_\star(X)_{n+1} \; \Big| \; d_j(\sigma)=\star \quad 2\leq j\leq n+1\Big\} \; = \; 
\Big\{ \sigma\in X_{n+2} \; \Big| \; d_j(\alpha)=\star\quad 2\leq j\leq n+2\Big\}\,.
$$
\end{proof}

Remarquons que d'après les Lemmes \ref{isoomega}, \ref{hachelemme} et \ref{minusone}, la restriction du foncteur:
$$
\xymatrix@C+10pt{\simp_\star\ar[r]^-{\mathbb{\Omega}}&\simp_\star}
$$ 
à la sous-catégorie pleine des ensembles simpliciaux $n$-fibrants réduits (les complexe de Kan $X$ tels que $X_0=\star$ et $\pi_i(X)=0$ pour $i \geq n+1$), admet un foncteur dérivé total par rapport aux $n$-équivalences faibles pointées $\pi^{-1}{\bf W}_{n}$. Plus encore un tel foncteur dérivé total est isomorphe à la restriction d'un foncteur espace de lacets de la catégorie de modèles $(\simp_{\star},\pi^{-1}{\bf W}_{n}, {\bf mono},\pi^{-1}{\bf fib}_{n})$ à la sous-catégorie des $n$-types d'homotopie pointés connexes.   

L'énoncé suivant est une version combinatoire du Lemme \ref{minusone} qu'on vient de montrer (voir aussi \cite{duskin}): 

\begin{lemme}\label{desgroupDD}
Soit $1\leq n< \infty$. Si $X$ est un $n$-groupoïde de Kan réduit (ce qu'on pourrait appeler un $n$-\emph{groupe de Kan}), alors $\mathbb{\Omega}_\star(X)$ est un $(n-1)$-groupoïde de Kan.
\end{lemme}
\begin{proof}
Il faut de montrer que si $X$ est un ensemble simplicial réduit lequel est un complexe de Kan qui satisfait la condition d'extension de Kan de façon stricte en dimension $m$ pour $m\geq n$ alors $\mathbb{\Omega}_\star(X)$ est un complexe de Kan qui satisfait la condition d'extension de Kan de façon stricte en dimension $m$ pour $m\geq n-1$.

Remarquons pour commencer que du plongement de Yoneda et de la définition de l'ensemble simplicial $\Lambda^{m+1,k}$, on peut identifier pour $0\leq k\leq m+1$ la fonction: 
$$
\xymatrix@C+20pt{\mathrm{Hom}_{\simp}\big(\Delta^{m+1},X\big)\ar[r]^-{\alpha^{m,k}_X} & \mathrm{Hom}_{\simp}\big(\Lambda^{m+1,k},X\big)}
$$
avec la fonction:
\begin{equation}\label{ide11}
\begin{split}
X_{m+1} &\vcenter{\xymatrix@R=3pt@C+10pt{\ar[r] & }}
\Bigg\{\text{\scriptsize{$\big(a_{0},\dots,a_{k-1},a_{k+1}\dots,a_{m+1}\big)\,\in\,\big(X_m\big)^{m+1}$}}\;\Bigg|\;
\vcenter{\xymatrix@R=1pt{\text{\scriptsize{$d_{i}a_{j}=d_{j-1}a_{i}$}}\\ \text{\scriptsize{si $0\leq i<j\leq m+1$ et 
$i,j\neq k$.}}}} \Bigg\}\,.\\
\eta \quad & \quad \longmapsto \quad \quad \qquad \; \big(d_{0}(\eta),\dots,d_{k-1}(\eta),d_{k+1}(\eta)\dots,d_{m+1}(\eta)\big) 
\end{split}
\end{equation}

Du même, si $0\leq k\leq m$ on identifie:
$$
\xymatrix@C+20pt{\mathrm{Hom}_{\simp}\big(\Delta^{m},\mathbb{\Omega}_\star(X)\big)\ar[r]^-{\alpha^{m-1,k}_{\mathbb{\Omega}_\star(X)}} & \mathrm{Hom}_{\simp}\big(\Lambda^{m,k},\mathbb{\Omega}_\star(X)\big)}
$$
avec la fonction suivante: 
\begin{equation}\label{ide22}
\begin{split}
\Big\{\alpha\in X_{m+1}\,\Big|\, d_{m+1}(\alpha)=\star\Big\}
&\vcenter{\xymatrix@R=3pt@C-10pt{\ar[r] & }}
\Bigg\{\text{\scriptsize{$\big(b_{0},\dots, b_{k-1},b_{k+1}\dots, b_{m}\big)\,\in\,\big(X_m\big)^{m}$}}\;\Bigg|\;
\vcenter{\xymatrix@R=1pt{\text{\scriptsize{$d_{m}(b_i)=\star$ pour tout $0\leq i \leq m$,\, $i\neq k$}}\\ \text{\scriptsize{et\; $d_{i}b_{j}=d_{j-1}b_{i}$ si $0\leq i<j\leq m$ et 
$i,j\neq k$.}}}} \Bigg\}\,.\\
\alpha \quad \qquad \quad & \quad \mapsto \quad \quad \qquad \; \big(d_{0}(\alpha),\dots,d_{k-1}(\alpha),d_{k+1}(\alpha)\dots,d_{m}(\alpha)\big) 
\end{split}
\end{equation}

Plus encore, on vérifie aussi-tôt que si $m\geq 0$ et $0\leq k\leq m$ on a un carré commutatif:
\begin{equation}\label{cacartgrp}
\begin{split}
& X_{m+1} \vcenter{\xymatrix@R=3pt@C+10pt{\ar[r]^-{\eqref{ide11}} & }}
\Bigg\{\text{\scriptsize{$\big(a_{0},\dots,a_{k-1},a_{k+1}\dots,a_{m+1}\big)\,\in\,\big(X_m\big)^{m+1}$}}\;\Bigg|\;
\vcenter{\xymatrix@R=1pt{\text{\scriptsize{$d_{i}a_{j}=d_{j-1}a_{i}$}}\\ \text{\scriptsize{si $0\leq i<j\leq m+1$ et 
$i,j\neq k$.}}}} \Bigg\}  \\
&\xymatrix{&\\\ar@<-10pt>[u]& \ar@<-145pt>[u]}\\
\Big\{\alpha\in X_{m+1}\,\Big|\, d&_{m+1}(\alpha)=\star\Big\}
\vcenter{\xymatrix@R=3pt@C-10pt{\ar[r]_-{\eqref{ide22}} & }}
\Bigg\{\text{\scriptsize{$\big(b_{0},\dots, b_{k-1},b_{k+1}\dots, b_{m}\big)\,\in\,\big(X_m\big)^{m}$}}\;\Bigg|\;
\vcenter{\xymatrix@R=1pt{\text{\scriptsize{$d_{m}(b_i)=\star$ pour tout $0\leq i \leq m$,\, $i\neq k$}}\\ \text{\scriptsize{et\; $d_{i}b_{j}=d_{j-1}b_{i}$ si $0\leq i<j\leq m$ et 
$i,j\neq k$.}}}} \Bigg\}\,,
\end{split}
\end{equation}
où la fonction qui monte à gauche est la contention et cela qui monte à droite est définie par la règle:
$$
\big(b_{0},\dots, b_{k-1},b_{k+1}\dots, b_{m}\big) \;\qquad \longmapsto \;\qquad
\big(b_{0},\dots, b_{k-1},b_{k+1}\dots, b_{m},\star\big)\,. 
$$

En fait on vérifie sans peine que le carré \eqref{cacartgrp} est un carré cartésien des ensembles. Donc si $X$ satisfait la condition d'extension de Kan en dimension $m$ (resp. il la satisfait de façon stricte en dimension $m$), alors $\mathbb{\Omega}_\star(X)$ satisfait la condition d'extension de Kan en dimension $m-1$ (resp. il la satisfait de façon stricte en dimension $m$).
\end{proof}

Si $X$ est un ensemble simplicial réduit considérons maintenant le morphisme d'ensembles simpliciaux réduits:
$$
\mathcal{H}\Big(\xymatrix{\mathbb{D}(X)\ar[r]^-{\varrho_X}&X}\Big)\;=\;\xymatrix{\mathbb{D}^{red}(X)\ar[r]^-{\varrho_X^{red}} &X}\,.
$$
image du morphisme $\xymatrix{\mathbb{D}(X)\ar[r]^-{\varrho_X}&X}$ par le foncteur adjoint à droite de l'adjonction:
$$
\simp_{0}
\xymatrix@C+15pt{
\phantom{a}\ar@<+10pt>@{}[r]|{\text{\rotatebox[origin=c]{180}{$\perp$}}}\ar@{^(->}[r]_-{\nu}&\phantom{a}
\ar@<-7pt>@/_12pt/[l]_-{\mathcal{H}}}
\simp_\star\,.
$$

Explicitement on a que $\mathbb{D}^{red}(X)_0=\{\star\}\subset X_1$ et:
$$
\mathbb{D}^{red}(X)_n = \Big\{ \;\;\alpha\in X_{n+1} \;\;\Big|\;\;  d_{i_n}\circ \dots \circ d_{i_1}(\alpha) = \star \in X_1 \quad \text{toujours \; que} \quad 0\leq i_j\leq  n+1-j \Big\}
$$
si $n\geq 1$.

On définit \emph{le espace de lacets simplicial réduit} d'un ensemble simplicial réduit $X$, comme le sous-ensemble simplicial réduit $\mathbb{\Omega}^{red}(X)$ de $\mathbb{D}^{red}(X)$ défini par: 
\begin{equation*}
\begin{split}
\mathbb{\Omega}^{red}(X)_n  \, = \; & \Big\{\alpha\in \mathbb{D}^{red}(X)_n \;\;\Big|\;\; d_{n+1}(\alpha)=\star\Big\}\\
 \, = \; & \Big\{\alpha\in X_{n+1} \;\;\Big|\;\; d_{n+1}(\alpha)=\star \in X_n \;\;\text{et}\;\; d_{i_n}\circ \dots \circ d_{i_1}(\alpha) = \star \in X_1\;\;\text{si}\;\; 0\leq i_j \leq n+1-j \Big\}
\end{split}
\end{equation*}
si $n\geq 1$.

En particulier on a un carré cartésien de $\simp_0$ (ou de $\simp$):
\begin{equation}\label{omegasimp2}
\xymatrix{
\mathbb{\Omega}^{red}(X) \ar[r] \ar[d] & \mathbb{D}^{red}(X) \ar[d]^-{\rho^{red}_X}\\ 
\star \ar[r] & X}
\end{equation}

\begin{lemme}\label{finomegared}
Soit $0 \leq n \leq \infty$. Si $X$ est un ensemble simplicial $n$-fibrant réduit c'est-à-dire $X$ est un objet fibrant de la catégorie de modèles $(\simp_{0},{\bf W}^{red}_{n}, {\bf mono},{\bf fib}_{n}^{red})$ de la Proposition \ref{modred}, ou de façon explicite $X$ est un complexe de Kan tel que $X_0=\star$ et $\pi_i(X)=0$ pour $i \geq n+1$, alors:
\begin{enumerate}
\item Le carré cartésien \eqref{omegasimp2} est un carré homotopiquement cartésien dans la catégorie de modèles  $(\simp_{0},{\bf W}^{red}_{n}, {\bf mono},{\bf fib}_{n}^{red})$.
\item Le seul morphisme $\xymatrix@C-8pt{\mathbb{D}^{red}(X) \ar[r]^-{\alpha_X}&\star}$ est une $\infty$-équivalence faible.
\item Si $1 \leq n  < \infty$ l'ensemble simplicial $\mathbb{\Omega}^{red}(X)$ est un objet fibrant de la catégorie de modèles $(\simp_{0},{\bf W}^{red}_{n-1}, {\bf mono},{\bf fib}_{n-1}^{red})$. En particulier $\xymatrix@C-8pt{\underbrace{\mathbb{\Omega}^{red}\circ\dots\circ\mathbb{\Omega}^{red}}_{n}(X)\ar[r] &\star}$ est une $\infty$-équivalence faible.
\end{enumerate}
\end{lemme}
\begin{proof}
Si $X$ est un objet fibrant de $(\simp_{0},{\bf W}^{red}_{n}, {\bf mono},{\bf fib}_{n}^{red})$ d'après le Lemme \ref{hachelemme} le morphisme:
$$
\mathcal{H}\Big(\xymatrix{\mathbb{D}(X)\ar[r]^-{\varrho_X}&X}\Big)\;=\;\xymatrix{\mathbb{D}^{red}(X)\ar[r]^-{\varrho_X^{red}} &X}\,,
$$ 
image du morphisme $\xymatrix@C-10pt{\mathbb{D}(X)\ar[r]^-{\varrho_X}&X}$ par le foncteur adjoint à droite de l'adjonction:
$$
\simp_{0}
\xymatrix@C+15pt{
\phantom{a}\ar@<+10pt>@{}[r]|{\text{\rotatebox[origin=c]{180}{$\perp$}}}\ar@{^(->}[r]_-{\nu}&\phantom{a}
\ar@<-7pt>@/_12pt/[l]_-{\mathcal{H}}}
\simp_\star\,
$$
est une fibration de la catégorie de modèles $(\simp_{0},{\bf W}^{red}_{n}, {\bf mono},{\bf fib}_{n}^{red})$.

En particulier les ensembles simpliciaux réduits $\mathbb{D}^{red}(X)$ et $\mathbb{\Omega}^{red}(X)$ sont des objets fibrants de la catégorie de modèles $(\simp_{0},{\bf W}^{red}_{n}, {\bf mono},{\bf fib}_{n}^{red})$.

Plus encore, vu que $\xymatrix@C-8pt{\mathbb{D}(X)\ar[r]^-{\alpha_X}&\star}$ est une $\infty$-équivalences faible entre complexes de Kan pointés:
$$
\mathcal{H}\Big(\xymatrix{\mathbb{D}(X)\ar[r]^-{\alpha_X}&\star}\Big)\;=\;\xymatrix{\mathbb{D}^{red}(X)\ar[r] &\star}
$$ 
est aussi une $\infty$-équivalence faible entre complexes de Kan réduits.

Si $X$ est un objet fibrant de la catégorie de modèles $(\simp_{0},{\bf W}^{red}_{n}, {\bf mono},{\bf fib}_{n}^{red})$ vu que le morphisme $\xymatrix@C-5pt{\mathbb{D}^{red}(X)\ar[r]^-{\varrho_X^{red}} &X}$ est une fibration de $(\simp_{0},{\bf W}^{red}_{n}, {\bf mono},{\bf fib}_{n}^{red})$, il se suit du Lemme \ref{fibrabra} que \eqref{omegasimp2} est un carré homotopiquement cartésien de $(\simp_{0},{\bf W}^{red}_{n}, {\bf mono},{\bf fib}_{n}^{red})$. 

D'un autre, rappelons qu'un objet fibrant $Z$ de $(\simp_{0},{\bf W}^{red}_{n}, {\bf mono},{\bf fib}_{n}^{red})$ est fibrant dans la catégorie de modèles $(\simp_{0},{\bf W}^{red}_{n-1}, {\bf mono},{\bf fib}_{n-1}^{red})$ si et seulement si $\pi_{n}(Z)=0$. On déduit que $\pi_{n}\big(\mathbb{\Omega}^{red}(X)\big) \cong \pi_{n+1}\big(X\big)=0$ de la Proposition \ref{kankan} parce que si $n\geq 1$:
\begin{equation}\label{losalfatos}
\begin{split}
&\Big\{\alpha\in\mathbb{\Omega}^{red}(X)_n \; \Big| \;  d_i(\alpha)=\star\quad 0\leq i\leq n\Big\} \\
\; =\; 
&\Big\{ \alpha\in X_{n+1} \; \Big|  \; d_i(\alpha)=\star\;\;\text{si}\;\; 0\leq i\leq n+1\;\;\text{et}\;\; d_{i_n}\circ \dots \circ d_{i_1}(\alpha) = \star \in X_1\;\;\text{si}\;\; 0\leq i_j \leq n+1-j \Big\} \\
= \; 
&\Big\{ \alpha\in X_{n+1} \; \Big|  \; d_i(\alpha)=\star\;\;\text{si}\;\; 0\leq i\leq n+1\;\Big\}
\end{split}
\end{equation}
et si $\alpha,\beta$ sont des éléments de l'ensemble \eqref{losalfatos} alors on a que:
\begin{equation*}
\begin{split}
&\Big\{w\in\mathbb{\Omega}^{red}(X)_{n+1} \; \Big| \;  d_0(w) = \alpha \,,\;\; d_1(w)=\beta\;\;\text{et}\;\;d_i(w)=\star\;\;\text{si}\;\;2\leq i\leq n+1\Big\} \\
\; =\; 
&\Big\{ w\in X_{n+2} \; \Big| \; d_0(w) = \alpha \,,\;\; d_1(w)=\beta\,,\;\;d_i(w)=\star\;\;\text{si}\;\;2\leq i\leq n+2\;\;\\
&\qquad \;\qquad\qquad \;\;\text{et}\;\; d_{i_{n+1}}\circ \dots \circ d_{i_1}(\alpha) = \star \in X_1\;\;\text{si}\;\; 0\leq i_j \leq n+2-j \Big\} \\
= \; 
&\Big\{ w\in X_{n+2} \; \Big| \; d_0(w) = \alpha \,,\;\; d_1(w)=\beta\,,\;\;d_i(w)=\star\;\;\text{si}\;\;2\leq i\leq n+2\;\Big\}\,.
\end{split}
\end{equation*}
\end{proof}

On déduit des Lemmes \ref{finomegared} et \ref{isoomega} l'énoncé suivant:

\begin{corollaire}\label{dimminus}
Si $0\leq n\leq  \infty$ le foncteur $\mathbb{\Omega}^{red}\colon\xymatrix@C-5pt{\simp_0\ar[r]&\simp_0}$ défini ci-dessus respecte les $n$-équivalences faibles entre les objets fibrants de la catégorie de modèles $(\simp_{0},{\bf W}^{red}_{n}, {\bf mono},{\bf fib}_{n}^{red})$ de la Proposition \ref{modred}. Plus encore un foncteur dérivé à droite de $\mathbb{\Omega}^{red}$ est un foncteur espace de lacet de $(\simp_{0},{\bf W}^{red}_{n}, {\bf mono},{\bf fib}_{n}^{red})$.

En particulier si $0< n< \infty$ la catégorie de modèles $(\simp_{0},{\bf W}^{red}_{n}, {\bf mono},{\bf fib}_{n}^{red})$ munie de l'enrichissement $\underline{\mathrm{Hom}}_{\simp_0}$ de \eqref{homred} est une catégorie de modèles simpliciale pointée d'ordre $n-1$ (voir le Théorème \ref{cassi1}).
\end{corollaire}

D'après le Corollaire \ref{dimminus} le foncteur dérivé des morphismes ${\bf R}\underline{\mathrm{Hom}}_{\simp_0}$ de la catégorie des $n$-types d'homotopie réduits $\mathrm{Ho}_{n}\big(\,\simp_{0}\,\big)$ c'est-à-dire de la catégorie homotopique $\simp_{0}\big[({\bf W}^{red}_{n})^{-1}\big]$ de la catégorie de modèles $(\simp_{0},{\bf W}^{red}_{n}, {\bf mono},{\bf fib}_{n}^{red})$, est un enrichissement tensoré et cotensoré sur la catégorie cartésienne fermée $\mathrm{Ho}_{n-1}\big(\,\simp\,\big)$ des $(n-1)$-types d'homotopie (la catégorie homotopique des $(n-1)$-groupoïdes). 

La catégorie des $n$-types d'homotopie réduits munie d'un tel enrichissement sur la catégorie homotopique des $(n-1)$-groupoïdes est dite parfois la \emph{catégorie homotopique des $n$-groupes}. 

On arrive à cette catégorie d'une façon équivalente comme suit: Posons ${\bf Fib}^{n}_0$ pour noter la sous-catégorie de $\simp_0$ dont les objets sont les ensembles simpliciaux $n$-fibrants réduits, c'est-à-dire les objets fibrants de la catégorie de modèles $(\simp_{0},{\bf W}^{red}_{n}, {\bf mono},{\bf fib}_{n}^{red})$ de la Proposition \ref{modred}, ou de façon explicite les complexe de Kan $X$ tels que $X_0=\star$ et $\pi_i(X)=0$ pour $i \geq n+1$.

On déduit du Corollaire \ref{dimminus}:

\begin{corollaire} 
Soit $0< n< \infty$. Si $X$ et $Y$ sont des ensembles simpliciaux $n$-fibrants réduits l'ensemble simplicial des morphismes $\underline{\mathrm{Hom}}_{\simp_0}(X,Y)$ de \eqref{homred} est un ensemble simplicial $(n-1)$-fibrant. Autrement dit le foncteur $\underline{\mathrm{Hom}}_{\simp_0}$ induit par restriction un enrichissement de la catégorie ${\bf Fib}_0^{n}$ sur la catégorie cartésienne fermée ${\bf Fib}^{n-1}$ (voir le Lemme \ref{nfibcartesiene}).
\end{corollaire}

Si $h{\bf Fib}_0^{n}$ note la catégorie des ensembles simpliciaux $n$-fibrants réduits dont l'ensemble des morphismes est l'ensemble des composantes connexes par arcs $\pi_0\big(\underline{\mathrm{Hom}}_{\simp_0}\big)$, il se suit que le foncteur d'inclusion de ${\bf Fib}^n_0$ vers $\simp_0$ induit une équivalence de catégories entre $h{\bf Fib}^n_0$ et la catégorie des $n$-types d'homotopie réduits $\mathrm{Ho}_{n}(\simp_0)$. 

Enfin remarquons que le foncteur $\underline{\mathrm{Hom}}_{\simp_0}$ est un enrichissement de $h{\bf Fib}_0^n$ sur la catégorie cartésienne fermée $h{\bf Fib}^{n-1}$, tandis que l'enrichissement de la catégorie des $n$-types d'homotopie réduits $\mathrm{Ho}_{n}(\simp_0)$ sur la catégorie cartésienne fermée des $(n-1)$-types d'homotopie $\mathrm{Ho}_{n-1}(\simp)$ est un foncteur dérivé des morphismes ${\bf R}\underline{\mathrm{Hom}}_{\simp_0}$.

\renewcommand{\thesubsection}{\S\thesection.\arabic{subsection}}
\subsection{}\;
\renewcommand{\thesubsection}{\thesection.\arabic{subsection}}

On va dire qu'une catégorie de modèles $\C$ est de \emph{dimension homotopique infinie} si pour tout entier $m > 0$ et tout foncteur espace de lacet ${\Omega}_{\C_\star}$ de la catégorie de modèles $\C_\star$ des objets pointés de $\C$ (voir \S\ref{lem1}), il existe un objet $X$ de $\C_\star$ tel que ${\Omega}_{\C_\star}^m(X)$ ne soit pas l'objet null de $\mathrm{Ho}\big(\C_\star\big)$. 

Une catégorie de modèles qui n'est pas de dimension homotopique infinie est dite de \emph{dimension homotopique finie}. Si $\C$ est une catégorie de modèles de dimension homotopique finie et $n\geq 0$ est un entier, on dit que $\C$ est de \emph{dimension homotopique $n$} s'il existe un foncteur espace de lacets ${\Omega}_{\C_\star}$ de la catégorie de modèles $\C_\star$ des objets pointés de $\C$ tel que:
\begin{enumerate}
\item ${\Omega}_{\C_\star}^{n}(X)$ n'est pas l'objet null de $\mathrm{Ho}\big(\C_\star\big)$ pour certain objets $X$ de $\C_\star$ (où $\Omega_{\C_\star}^0$ est le foncteur identité).
\item ${\Omega}_{\C_\star}^{n+1}(Y)$ est l'objet null de $\mathrm{Ho}\big(\C_\star\big)$ pour tout objet $Y$ de $\C_\star$ .
\end{enumerate}

Autrement dit la dimension homotopique d'une catégorie de modèle $\C$ est le entier le plus grand $n$ pour lequel la $n$-ième itération ${\Omega}_{\C_\star}\circ\cdots\circ{\Omega}_{\C_\star}$ d'un foncteur espace de lacets quelconques $\Omega_{\C_\star}$ de $\C_\star$ n'est pas le foncteur null. 

De façon équivalente la dimension homotopique d'une catégorie de modèle $\C$ mesure la taille de son espace des morphismes dérivé ${\bf R}\underline{\mathrm{Hom}}_{\C_\star}$. Plus précisément d'après le Théorème \ref{cassi1} une catégorie de modèles $\C$ est de dimension homotopique $n$ si et seulement si:
\begin{enumerate}
\item  Les groupes d'homotopie $\pi_{n+1}\big({\bf R}\underline{\mathrm{Hom}}_{\C_\star}(X,Y)\big)$ sont nuls pour tous les objets $X$ et $Y$ de $\C_\star$.
\item  Ils existent des objets $Z$ et $W$ de $\C_\star$ tels que $\pi_{n}\big({\bf R}\underline{\mathrm{Hom}}_{\C_\star}(Z,W)\big)$ n'est pas null.
\end{enumerate}

Plus encore d'après le Corollaire \ref{dugger2} la dimension homotopique d'une catégorie de modèles de Cisinski $\C$ est l'entier le plus petit $n\geq 0$ tel que $\C$ admette à équivalence de Quillen pointée près une structure de catégorie de modèles simplicial d'ordre $n$. 

Montrons:

\begin{corollaire}
$(\simp,{\bf W}_{\infty}, {\bf mono},{\bf fib}_\infty)$ et $(\simp_{0},{\bf W}^{red}_{\infty}, {\bf mono},{\bf fib}_{\infty}^{red})$ sont de catégories de modèles de dimension homotopique infinie.

D'un autre, si $n\geq 0$ la catégorie de modèles $(\simp,{\bf W}_{n}, {\bf mono},{\bf fib}_n)$ est de dimension homotopique $n$ et si $n\geq 1$ la catégorie de modèles $(\simp_{0},{\bf W}^{red}_{n}, {\bf mono},{\bf fib}_{n}^{red})$ est de dimension homotopique $n-1$.
\end{corollaire}
\begin{proof}
D'après le Lemme \ref{ntypessp} si $0\leq n \leq\infty$ la catégorie de modèles de Cisinski $(\simp,{\bf W}_{n}, {\bf mono},{\bf fib}_n)$ est une catégorie de modèles simplicial d'ordre $n$; donc $(\simp,{\bf W}_{n}, {\bf mono},{\bf fib}_n)$ est de dimension homotopique plus petite ou égale à $n$ (voir le Corollaire \ref{dugger2}). Du même si $0 < n \leq\infty$ d'après le Corollaire \ref{dimminus} la catégorie de modèles de Cisinski $(\simp_0,{\bf W}^{red}_{n}, {\bf mono},{\bf fib}^{red}_n)$ est de dimension homotopique plus petite ou égale à $n-1$.

Il se suit aussi-tôt que $(\simp,{\bf W}_{0}, {\bf mono},{\bf fib}_0)$ et $(\simp_0,{\bf W}^{red}_{1}, {\bf mono},{\bf fib}^{red}_1)$ sont des catégories de modèles de dimension homotopique $0$, parce que les foncteurs identité des catégories homotopiques $\simp_\star\big[(\pi^{-1}{\bf W}_0)^{-1}\big]$ et $\simp_0\Big[\big({\bf W}_1^{red}\big)^{-1}\Big]$ ne sont pas de foncteurs nuls.

Rappelons d'un autre que si $0\leq n \leq\infty$ d'après les Lemmes \ref{minusone}, \ref{finomegared} et \ref{isoomega} pour tout complexe de Kan réduit $X$ vérifiant que $\pi_i(X)=0$ si $i> n$, les ensembles simpliciaux $\mathbb{\Omega}_{\star}(X)$ et $\mathbb{\Omega}^{red}(X)$ sont des espaces de lacets fibrants de $X$ dans la catégorie de modèles des objets pointés de $(\simp,{\bf W}_{n}, {\bf mono},{\bf fib}_n)$ et dans la catégorie de modèles $(\simp_{0},{\bf W}^{red}_{n}, {\bf mono},{\bf fib}_{n}^{red})$ respectivement. Plus encore dans ce cas on a des isomorphismes (voir la preuve des Lemmes \ref{minusone} et \ref{finomegared}):
$$
\pi_i\big(\mathbb{\Omega}_\star^m(X)\big) \, \cong \, \pi_{m+i}(X) \qquad \text{si  $\quad m\geq 1$ et $i\geq 0$}
$$ 
$$\text{et}$$
$$
\pi_i\big((\mathbb{\Omega}^{red})^m(X)\big) \, \cong \, \pi_{m+i}(X) \qquad \text{si  $\quad m\geq 1$ et $i\geq 1$.}
$$ 

Si $m \geq 1$ considérons un complexe de Kan réduit $Z(m)$ tel que $\pi_m\big(Z(m)\big) \neq 0$ et $\pi_i\big(Z(m)\big)=0$ pour $i> m$ (prendre par exemple l'image par le foncteur adjoint à droite de l'adjonction \eqref{adjconne} d'un espace d'Eilenberg–Mac Lane connexe). 

On conclut que $(\simp,{\bf W}_{\infty}, {\bf mono},{\bf fib}_\infty)$ et $(\simp_{0},{\bf W}^{red}_{\infty}, {\bf mono},{\bf fib}_{\infty}^{red})$ sont de catégories de modèles de dimension homotopique infinie. En effet si $m\geq 1$ il existe un complexe de Kan pointé $X=Z(m)$ et un complexe de Kan réduit $Y=Z(m+1)$ tels que $\mathbb{\Omega}_\star(X)$ et $\mathbb{\Omega}^{red}_\star(Y)$ sont des espaces de lacets fibrants de $X$ et $Y$ dans $(\simp_\star,\pi^{-1}{\bf W}_{\infty}, {\bf mono},\pi^{-1}{\bf fib}_\infty)$ et $(\simp_{0},{\bf W}^{red}_{\infty}, {\bf mono},{\bf fib}_{\infty}^{red})$ respectivement, lesquels vérifient la propriété:
$$
\pi_0\big(\mathbb{\Omega}_\star^m(X)\big) \, \cong \, \pi_m\big(Z(m)\big) \neq 0 \qquad \text{et} \qquad 
\pi_1\big((\mathbb{\Omega}_\star^{red})^m(Y)\big) \, \cong \, \pi_{m+1}\big(Z(m+1)\big) \neq 0
$$
respectivement. Donc les espaces de lacets itérés $\mathbb{\Omega}_\star^m(X)$ et $(\mathbb{\Omega}^{red}_\star)^m(Y)$ ne sont pas l'objet null des catégories homotopiques $\simp_\star\big[(\pi^{-1}{\bf W}_\infty)^{-1}\big]$ et $\simp_0\Big[\big({\bf W}_\infty^{red}\big)^{-1}\Big]$ respectivement.

De la même manière, si $0 < n <\infty$ on considère le complexe de Kan pointé $X=Z(n)$ et le complexe de Kan réduit $Y=Z(n+1)$. On déduit que l'ensemble simplicial pointé $\mathbb{\Omega}_\star(X)$ et l'ensemble simplicial réduit $\mathbb{\Omega}^{red}_\star(Y)$ sont des espaces de lacets fibrants de $X$ et $Y$ dans les catégories de modèles $(\simp_\star,\pi^{-1}{\bf W}_{n}, {\bf mono},\pi^{-1}{\bf fib}_n)$ et $(\simp_{0},{\bf W}^{red}_{n+1}, {\bf mono},{\bf fib}_{n+1}^{red})$ respectivement, tels que:
$$
\pi_0\big(\mathbb{\Omega}_\star^n(X)\big) \, \cong \, \pi_n\big(Z(n)\big) \neq 0 \qquad \text{et} \qquad 
\pi_1\big((\mathbb{\Omega}_\star^{red})^n(Y)\big) \, \cong \, \pi_{n+1}\big(Z(n+1)\big) \neq 0
$$
respectivement. Autrement dit les espaces de lacets itérés $\mathbb{\Omega}_\star^n(X)$ et $(\mathbb{\Omega}^{red}_\star)^n(Y)$ ne sont pas l'objet null des catégories homotopiques $\simp_\star\big[(\pi^{-1}{\bf W}_n)^{-1}\big]$ et $\simp_0\Big[\big({\bf W}_{n+1}^{red}\big)^{-1}\Big]$ respectivement.

Donc $(\simp,{\bf W}_{n}, {\bf mono},{\bf fib}_n)$ et $(\simp_{0},{\bf W}^{red}_{n+1}, {\bf mono},{\bf fib}_{n+1}^{red})$ sont de catégories de modèles de dimension homotopique $n$. 
\end{proof}

\renewcommand{\thesubsection}{\S\thesection.\arabic{subsection}}
\subsection{}\;
\renewcommand{\thesubsection}{\thesection.\arabic{subsection}}

On désigne par $\ssimp$ la catégorie des ensembles bisimpliciaux; c'est-à-dire, des foncteurs:
$$\xymatrix@C+6pt@R=1pt{
\big({\bf \Delta}\times{\bf \Delta}\big)^{op}\ar[r]^-{X}&\ens\\
\text{\scriptsize{$\big([p],[q]\big)$}}\;\ar@{}[r]|-{\longmapsto}& \text{\scriptsize{$X_{p,q}$}}}$$ 
et des transformations naturelles entre eux. 

Les projections canoniques:
$$
\xymatrix@R=2pt{\Delta \times \Delta\ar[r]^-{\mathrm{p}_{1}} &\Delta \,\\  \text{\scriptsize{$\big([p],[q]\big)$}} \ar@{}[r]|-{\longmapsto} &\text{\scriptsize{$[p]$}}} 
\qquad\text{et}\qquad
\xymatrix@R=2pt{\Delta \times \Delta\ar[r]^-{\mathrm{p}_{2}} &\Delta \,,\\   \text{\scriptsize{$\big([p],[q]\big)$}} \ar@{}[r]|-{\longmapsto} &\text{\scriptsize{$[q]$}}} 
$$
induisent les foncteurs:
$$
\phantom{\text{et}\qquad}\xymatrix@C+10pt{\;\simp\;\ar[r]^-{\mathrm{p}_{1}^*} &\;\ssimp} \qquad\text{défini \;\; par}\qquad \mathrm{p}_{1}^*(K)_{p,q}\;=\; K_{p}\,,
$$
$$
\text{et}\qquad\xymatrix@C+10pt{\;\simp\;\ar[r]^-{\mathrm{p}_{2}^*} &\;\ssimp} \qquad\text{par}\qquad \mathrm{p}_{2}^*(K)_{p,q}\;=\; K_{q}\,;
$$
admettant des adjoints donnés par les formules:
$$
{\mathrm{p}_1\,}_{!} (X) = \underset{q}{\mathrm{colim}} \; X_{\bullet,q}\,,\qquad {\mathrm{p}_1\,}_{*} (X) = \underset{q}{\mathrm{lim}} \; X_{\bullet,q} \, \cong \, X_{\bullet, 0}\,,
$$
$$
{\mathrm{p}_2\,}_{!} (X) = \underset{p}{\mathrm{colim}} \, X_{p,\bullet}\qquad\text{et}\qquad\qquad{\mathrm{p}_2\,}_{*} (X) = \underset{p}{\mathrm{lim}} \; X_{p,\bullet} \, \cong \, X_{0,\bullet}\,.
$$

On définit le \emph{produit boite} des ensembles simpliciaux par la règle: 
$$
\vcenter{
\xymatrix@C+6pt@R=1pt{\simp\times\simp\ar[r]^-{\boxtimes}&\ssimp\,.
 \,\\  \text{\scriptsize{$\big(K,L\big)$}} \ar@{}[r]|-{\longmapsto} &\text{\scriptsize{$\mathrm{p}_{1}^*(K)\times \mathrm{p}_{2}^*L$}}}}
$$

Remarquons que si $k\geq 0$ les inclusions canoniques:
\begin{equation}\label{iinclusss}
\xymatrix@R=2pt{\Delta \ar[r]^-{\mathrm{i}_{1,k}} &\Delta \times \Delta\,\\ \text{\scriptsize{$[p]$}} \ar@{}[r]|-{\longmapsto} & \text{\scriptsize{$\big([p],[k]\big)$}} } 
\qquad\text{et}\qquad
\xymatrix@R=2pt{\Delta \ar[r]^-{\mathrm{i}_{2,k}} &\Delta \times \Delta\,,\\ \text{\scriptsize{$[q]$}} \ar@{}[r]|-{\longmapsto} & \text{\scriptsize{$\big([k],[q]\big)$}}} 
\end{equation}
induisent des adjonctions: 
$$
\ssimp\xymatrix@C+15pt{\phantom{a}\ar@/_12pt/[r]_{(\mathrm{i}_{1,k})^*}\ar@{}[r]|-{\perp}&\ar@/_12pt/[l]_{(\mathrm{i}_{1,k})_{!}}\phantom{a}}\simp
\qquad\text{et}\qquad
\ssimp\xymatrix@C+15pt{\phantom{a}\ar@/_12pt/[r]_{(\mathrm{i}_{2,k})^*}\ar@{}[r]|-{\perp}& \ar@/_12pt/[l]_{(\mathrm{i}_{2,k})_{!}}\phantom{a}}\simp\,,
$$
où on vérifie que ${(\mathrm{i}_{1,k})_{!}}\,A \, \cong \, A \boxtimes \Delta^k$ et ${(\mathrm{i}_{2,k})_{!}}\,A \, \cong \, \Delta^k \boxtimes A$ pour tout ensemble simplicial $A$.

Si $X$ est un ensemble bisimplicial:
$$
{(\mathrm{i}_{1,k})^*} X \,\,=\,\,  X_{\bullet,k} \;=\;\;
\xymatrix{
**[l]        X_{0,k} \ar@/_1pc/[r]|-{s^h_0}&
**[r]        X_{1,k} \,\ar@/_2pc/[l]|-{d^h_1}
                \ar@/_1pc/[l]|-{d^h_0} \ar@{}[r]|(.45){\dots\dots\dots\dots}&
**[l]        X_{n,k} \ar@/_1pc/[r]|-{s^{h}_{n}}
                \ar@/_2.8pc/[r]|-{s^{h}_0}^\vdots&
             X_{n+1,k} \ar@/_3pc/[l]|-{d^{h}_{n+1}}^{\vdots}
                \ar@/_1pc/[l]|-{d^{h}_0} \ar@{}[r]|-{\dots\dots\dots\dots}&}
$$
$$\text{et}$$
$$
{(\mathrm{i}_{2,k})^*} X \,\,=\,\, X_{k,\bullet} \;=\;\;
\xymatrix{
**[l]        X_{k,0} \ar@/_1pc/[r]|-{s^v_0}&
**[r]        X_{k,1} \,\ar@/_2pc/[l]|-{d^v_1}
                \ar@/_1pc/[l]|-{d^v_0} \ar@{}[r]|(.45){\dots\dots\dots\dots}&
**[l]        X_{k,n} \ar@/_1pc/[r]|-{s^{v}_{n}}
                \ar@/_2.8pc/[r]|-{s^{v}_0}^\vdots&
             X_{k,n+1} \ar@/_3pc/[l]|-{d^{v}_{n+1}}^{\vdots}
                \ar@/_1pc/[l]|-{d^{v}_0} \ar@{}[r]|-{\dots\dots\dots\dots}&},
$$
sont appelés respectivement l'\emph{ensemble simplicial horizontal} et \emph{vertical de degré $k$} de $X$.

Par exemple si $k=0$ on a que ${(\mathrm{i}_{1,0})_{!}} \cong \mathrm{p}_1^*$, ${(\mathrm{i}_{1,0})^*} \cong {\mathrm{p}_1\,}_*$, ${(\mathrm{i}_{2,0})_{!}} \cong \mathrm{p}_2^*$ et ${(\mathrm{i}_{2,0})^*} \cong {\mathrm{p}_2\,}_*$; en particulier:
\begin{equation}\label{hyv1}
\xymatrix@R=7pt{
 \mathrm{Hom}_{\simp}\big(\Delta^q,X_{p,\bullet}\big)\; \mbox{\Large $\cong$} \;
\mathrm{Hom}_{\ssimp}\big(\Delta^p\boxtimes \Delta^q,X\big)\; \mbox{\Large $\cong$} \;
\mathrm{Hom}_{\simp}\big(\Delta^p,X_{\bullet,q}\big)\,,\ar@{}[d]|-{\rotatebox[origin=c]{270}{{\Large $\cong$}}} \\
X_{p,q}}
\end{equation}
pour tout ensemble bisimplicial $X$.

Rappelons finalement que la catégorie $\ssimp$ est cartésienne fermée; en effet, pour tout ensemble bisimplicial $X$ on a une adjonction:
$$
\ssimp\xymatrix@C+15pt{\phantom{a}\ar@/^12pt/[r]^-{X \, \times \; \cdot\,}\ar@{}[r]|-{\perp}&\ar@/^12pt/[l]^-{\mathrm{hom}(X, \,\cdot\,)}\phantom{a}}\ssimp\,\qquad\;\text{où}\qquad\mathrm{hom}(X,Y)_{p,q}\,=\,\mathrm{Hom}_{\ssimp}\big(X\times(\Delta^p\boxtimes\Delta^q),Y\big)\,.
$$

On en déduit deux enrichissement de $\ssimp$ sur la catégorie cartésienne fermée des ensembles simpliciaux $\simp$: 
\begin{equation}\label{enrri12}
\begin{split}
\underline{\mathrm{Hom}}_{\ssimp}^{(1)}\big( X, Y\big)_{n}\,=&\;\mathrm{Hom}_{\ssimp}\big(X\times \mathrm{p}_1^*(\Delta^n),Y\big)\\
\text{et}&\\
\underline{\mathrm{Hom}}_{\ssimp}^{(2)}\big( X, Y\big)_{n}\,=&\;\mathrm{Hom}_{\ssimp}\big(X\times \mathrm{p}_2^*(\Delta^n),Y\big)
\end{split}
\end{equation}
munie des adjonctions:
$$
\simp\xymatrix@C+15pt{\phantom{a}\ar@/^12pt/[r]^-{X \, \times \, \mathrm{p}_{i}^* ( \cdot ) \,}\ar@{}[r]|-{\perp}&\ar@/^12pt/[l]^-{\underline{\mathrm{Hom}}_{\ssimp}^{(i)}\big(X,\,\cdot\,\big)}\phantom{a}}\ssimp
\qquad \quad \text{et} \quad \qquad
\simp\xymatrix@C+15pt{\phantom{a}\ar@/^12pt/[r]^-{\mathrm{hom}\big(\mathrm{p}_i^*(\cdot),Y\big)}\ar@{}[r]|-{\perp}&\ar@/^12pt/[l]^-{\underline{\mathrm{Hom}}_{\ssimp}^{(i)}\big(\,\cdot\,,Y\big)}\phantom{a}}\ssimp^{op}\,;
$$
pour $i=1$ ou $2$ respectivement.

\renewcommand{\thesubsection}{\S\thesection.\arabic{subsection}}
\subsubsection{}\;\label{premonsegal}
\renewcommand{\thesubsection}{\thesection.\arabic{subsection}}

Posons $\mathbf{MS}$ pour noter la sous-catégorie pleine de $\ssimp$ dont les objets sont les ensembles bisimpliciaux $X$ tels que $X_{p,0}=\star$ pour tout $p\geq 0$. De façon équivalente, $\mathbf{MS}$ est la catégorie des foncteurs:
$$\xymatrix@C+6pt@R=1pt{
{\bf \Delta}^{op}\ar[r]^-{X}&\simp_0\\
\text{\scriptsize{$[p]$}}\;\ar@{}[r]|-{\longmapsto}& \text{\scriptsize{$X_{p,\bullet}$}}}$$ 
et des transformations naturelles entre eux. On appelle les objets de la catégorie $\mathbf{MS}$ les \emph{pré-monoïdes de Segal} (les pré-catégories de Segal avec un seul objet).

Si $\big(\ssimp\big)_{\star}$ note la catégorie des ensembles bisimpliciaux pointés, le foncteur d'inclusion des pré-monoïdes de Segal $\mathbf{MS}$ vers $\ssimp$ se factorise:
\begin{equation}\label{factoMS}
\xymatrix@C-5pt@R-15pt{
\underset{\phantom{a}}{\mathbf{MS}}\,\ar@{^(->}@/_5pt/[rd] \ar@{^(->}[rr]^-{\textit{Foncteur d'inclusion}} & &\ssimp\,,\\
&\, \big(\ssimp\big)_{\star}\,\ar@/_5pt/[ru]_-{\textit{Foncteur d'oubli}}&} 
\end{equation}
autrement dit, tout pré-monoïde de Segal est canoniquement un ensemble bisimplicial pointé. 

Du même, on vérifie qu'on a des adjonctions:
\begin{equation}\label{adjMS}
\begin{split}
&\mathbf{MS}\,\cong\,\simp_0^{\Delta^{op}}
\xymatrix@C+15pt{
\phantom{a}\ar@<+10pt>@{}[r]|{\perp}\ar@{^(->}[r]&\phantom{a}
\ar@<-7pt>@/_12pt/[l]_-{\mathcal{F}^{\Delta^{op}}}}
\simp^{\Delta^{op}}\,\cong\,\ssimp\\&\phantom{a}\\
&\qquad\qquad\qquad\qquad\quad\text{et}\\&\phantom{a}\\
&\mathbf{MS}\,\cong\,\simp_0^{\Delta^{op}}
\vcenter{\xymatrix@C+15pt{
\phantom{a}\ar@<-10pt>@{}[r]|{\perp}\ar@<+10pt>@{}[r]|{\perp}\ar@{^(->}[r]&\phantom{a}
\ar@<-7pt>@/_12pt/[l]_-{\mathcal{G}^{\Delta^{op}}}\ar@<+7pt>@/^12pt/[l]^-{\mathcal{H}^{\Delta^{op}}}}}
\simp_\star^{\Delta^{op}}\,\cong\,\big(\ssimp\big)_{\star}\,;
\end{split}
\end{equation} 
définies à partir des adjonctions \eqref{adjredu}.

Les version pointées (voir \eqref{morstar} et \eqref{adjstar}) des enrichissements \eqref{enrri12} de la catégories des ensembles bisimpliciaux, induisent dans la catégorie $\mathbf{MS}$ les enrichissements simpliciaux suivants: 
\begin{equation}\label{MSHOM(1)}
\begin{split}
\underline{\mathrm{Hom}}_{\mathbf{MS}}^{(1)} \big( X, Y \big)_n
\,=&\,\mathrm{Hom}_{\mathbf{MS}}\Big( \big(X\times \mathrm{p}_1^*(\Delta^n)\big)\big/ \big(\star\times \mathrm{p}_1^*(\Delta^n) \big)  ,Y  \Big) \\
\,=&\,\mathrm{Hom}_{\ssimp_\star}\Big(\big(X\times \mathrm{p}_1^*(\Delta^n)\big)\big/ \big(\star\times \mathrm{p}_1^*(\Delta^n) \big) ,Y\Big) \\
\,=&\,\underline{\mathrm{Hom}}_{\ssimp_\star}^{(1)}\big( X, Y\big)_{n}\,,
\end{split}
\end{equation}
$$\text{et}$$
\begin{equation}
\begin{split}
\underline{\mathrm{Hom}}_{\mathbf{MS}}^{(2)} \big( X, Y \big)_n
\,=&\,\mathrm{Hom}_{\mathbf{MS}}\Big( \big(X\times \mathrm{p}_2^*(\Delta^n)\big)\big/ \big(\star\times \mathrm{p}_2^*(\Delta^n) \big)  ,Y  \Big) \\
\,=&\,\mathrm{Hom}_{\ssimp_\star}\Big(\big(X\times \mathrm{p}_2^*(\Delta^n)\big)\big/ \big(\star\times \mathrm{p}_2^*(\Delta^n) \big) ,Y\Big) \\
\,=&\,\underline{\mathrm{Hom}}_{\ssimp_\star}^{(2)}\big( X, Y\big)_{n}\,,
\end{split}
\end{equation}
munis des adjonctions:
\begin{equation}\label{adjoiii}
\simp\xymatrix@C+15pt{\phantom{a}\ar@/^12pt/[r]^-{X \, \overset{(i)}{\wedge} \, \cdot \,}\ar@{}[r]|-{\perp}&\ar@/^12pt/[l]^-{\underline{\mathrm{Hom}}^{(i)}_{\mathbf{MS}}\big(X,\,\cdot\,\big)}\phantom{a}}\mathbf{MS}
\qquad \quad \text{et} \quad \qquad
\simp\xymatrix@C+15pt{\phantom{a}\ar@/^12pt/[r]^-{Y^{\overset{(i)}{\wedge}\,\cdot\,}}\ar@{}[r]|-{\perp}&\ar@/^12pt/[l]^-{\underline{\mathrm{Hom}}^{(i)}_{\mathbf{MS}}\big(\,\cdot\,,Y\big)}\phantom{a}}\mathbf{MS}^{op}\,\;
\end{equation}
définis de la façon suivante:
$$
X \, \overset{(i)}{\wedge} \, K = \big(X\times \mathrm{p}_i^*(K) \big)\big/ \big( \star\times \, \mathrm{p}_i^*(K) \big)  
$$
$$\quad \text{et} \quad $$
$$
\big(Y^{\overset{(i)}{\wedge}K}\big)_{p,q} = \mathrm{Hom}_{\mathbf{MS}} 
\bigg(\Big(\mathrm{p}_i^* K \times (\Delta^p\boxtimes \Delta^q) \Big)\big/ \Big(\mathrm{p}_i^* K\times\big(\Delta^p\boxtimes {\bf sq}_0(\Delta^q)\big)\Big) , Y\bigg)
$$ 

Remarquons que si $X$ et $Y$ sont de pré-monoïdes de Segal quelconques, on déduit de la première adjonction de \eqref{adjMS} que:
\begin{equation}\label{homsimp}
\begin{split}
\underline{\mathrm{Hom}}^{(1)}_{\bf MS}\big(X,Y\big)_n \, & = \, \mathrm{Hom}_{\bf MS}\Big( \big(X\times \mathrm{p}_1^*(\Delta^n)\big)\big/ \big(\star\times \mathrm{p}_1^*(\Delta^n)\big)  ,Y  \Big) \\
& \cong \, \mathrm{Hom}_{\ssimp}\big(X\times \mathrm{p}_1^*(\Delta^n), Y \big) \\
& = \,  \underline{\mathrm{Hom}}^{(1)}_{\ssimp}\big(X,Y\big)_n\,.
\end{split}
\end{equation}

En revanche on a une fonction par nécessairement bijective:
$$
\xymatrix@R=3pt{
\underline{\mathrm{Hom}}^{(2)}_{\bf MS}\big(X,Y\big)_n  \ar@{}[d]|-{\rotatebox[origin=c]{270}{$=$}}\\ 
\mathrm{Hom}_{\bf MS}\big(X\times \mathrm{p}_2^*(\Delta^n) \big/ \big(\star\times \mathrm{p}_2^*(\Delta^n)\big)  , Y \big) \ar[dddd]   \\ \\\\\\
\mathrm{Hom}_{\bf MS}\Big( \big(X\times \mathrm{p}_2^*(\Delta^n)\big)\big/ \big(\star\times \mathrm{p}_2^*({\bf sq}_0\Delta^n)\big)  ,Y  \Big) 
\ar@{}[d]|-{\rotatebox[origin=c]{270}{$\cong$}}\\
\mathrm{Hom}_{\ssimp}\big(X\times \mathrm{p}_2^*(\Delta^n), Y \big)
\ar@{}[d]|-{\rotatebox[origin=c]{270}{$=$}}\\
\underline{\mathrm{Hom}}^{(2)}_{\ssimp}\big(X,Y\big)_n \,.}
$$

\begin{proposition}\label{reedyMS}
Soit $0\leq n\leq \infty$. La catégorie des pré-monoïdes de Segal ${\bf MS}$ munie de l'enrichissement $\underline{\mathrm{Hom}}^{(2)}_{{\bf MS}}$ admet une structure de catégorie de modèles simplicial d'ordre $n$, lorsque\footnote{Remarquons qu'un morphisme $f$ de $\mathbf{MS}$ est un monomorphisme si et seulement si  $f_{p,q}$ est une fonction injective pour tous $p,q\geq 0$}
\begin{align*}
\big\{\,\text{équivalences faibles}\,\big\} \quad &=\quad 
\big\{\, f:\xymatrix@-14pt{X\ar[r]&Y}\,\big|\; \text{$f_{p,\bullet}$ est une $n$-équivalence}\\
 &\phantom{=}\;\;\quad \text{\phantom{$\{\,$}faible d'ensembles simpliciaux  $\forall\; p\geq 0$}\,\big\}\;\,=\;\, {\bf W}_{n}^{(2)}\,,\\
\big\{\,\text{cofibrations}\,\big\} \quad &=\quad  \big\{\,\text{monomorphismes}\,\big\} \;\;\,= \;\,  {\bf mono}\,,\\
\big\{\,\text{fibrations}\,\big\} \quad &=\quad  \big\{\,\text{morphismes avec la propriété de relèvement }\\
 &\phantom{=}\;\;\quad \text{\phantom{$\{\,$}à droite par rapport à ${\bf mono}\,\cap\,{\bf W}_{n}^{(2)}$}\,\big\}\;\,= \;\, {\bf fib}_{n}^{(2)}\,.
\end{align*}


Plus encore, si $X$ est un pré-monoïde de Segal, alors $X$ est un objet fibrant de la catégorie de modèles $({\bf MS},{\bf W}^{(2)}_{n}, {\bf mono},{\bf fib}_{n}^{(2)})$ si et seulement si, le morphisme d'ensembles simpliciaux réduits:
\begin{equation}\label{segalfi}
\xymatrix{\underline{\mathrm{Hom}}^{(2)}_{\ssimp}\big(\mathrm{p}_1^* \Delta^p , X\big) \ar[r]^-{} &  \underline{\mathrm{Hom}}^{(2)}_{\ssimp}\big(\mathrm{p}_1^* \partial\Delta^p , X\big)}
\end{equation}
est une fibration de la catégorie de modèles $(\simp_{0},{\bf W}^{red}_{n}, {\bf mono},{\bf fib}_{n}^{red})$ pour tout $p\geq 0$. 
\end{proposition}
\begin{proof}
Vu que ${\bf MS}$ est isomorphe à la catégorie des foncteurs:
$$
\xymatrix@C+6pt@R=1pt{
{\bf \Delta}^{op}\ar[r]^-{X}&\simp_0\,,\\
\text{\scriptsize{$[p]$}}\;\ar@{}[r]|-{\longmapsto}& \text{\scriptsize{$X_{p,\bullet}$}}}
$$ 
et ${\bf \Delta}^{op}$ admet une structure de catégorie de Reedy (voir la Définition 6.1.1 de \cite{hirschhorn}) où: 
$$
{\bf \Delta}^{op}_{+}=\big\{ n\overset{f}{\to} m \in {\bf \Delta}^{op} \,\big|\, \text{$f$ est surjective dans ${\bf \Delta}$} \big\}
\;,\quad
{\bf \Delta}^{op}_{-}=\big\{ n\overset{f}{\to} m \in {\bf \Delta}^{op} \,\big|\, \text{$f$ est injective dans ${\bf \Delta}$} \big\}
$$
et la fonction degré $\mathrm{Ob}\big({\bf \Delta}^{op}\big)\to\mathbb{N}$ est la fonction identité; on déduit de la catégorie de modèles $(\simp_{0},{\bf W}^{red}_{n}, {\bf mono},{\bf fib}_{n}^{red})$ de la Proposition \ref{modred} une structure de catégorie de modèles dans ${\bf MS}$ dite \emph{de Reedy} (voir le Théorème 16.3.4 de \cite{hirschhorn}), dont les équivalences faibles sont précisément les éléments de ${\bf W}_n^{(2)}$.

Pour décrire les cofibration et les fibrations de Reedy, rappelons que si $X$ est un objet de ${\bf MS}$ on définie pour $k\geq 0$ les objets:
\begin{equation}\label{LM}
L_k (X) = \underset{\underset{\text{surj. non id.}}{k\to l}}{\mathrm{colim}} \, X_{l,\bullet}
\quad\qquad \text{et} \quad\qquad
M_k (X) = \underset{\underset{\text{inj. non id.}}{l\to k}}{\mathrm{lim}} \, X_{l,\bullet}\,;
\end{equation}
où la colimite et la limite sont pris dans $\simp_0$. 

En particulier on a des morphismes canoniques:
$$
\xymatrix{
L_k (X) \ar[r]& X_{k,\bullet} \ar[r]& M_k (X)}
$$

Un morphisme $f\colon\xymatrix@-5pt{X\ar[r]&Y}$ de ${\bf MS}$ est dit une cofibration (resp. une fibration) de Reedy si le morphisme $\varphi$ dans le diagramme somme amalgamée (resp. du produit fibré) de $\simp_0$:
\begin{equation} \label{LM2}
\vcenter{\xymatrix@-15pt{ 
L_k(X)\ar[ddd]_-{\mathrm{colim}\,f_{l,\bullet}} \ar[rrr]&&& X_{k,\bullet}\ar[ddd]\ar@/^15pt/[ddddr]^-{f_{k,\bullet}} &  \\ &&&&\\&&&&\\
L_k(Y)\ar[rrr]\ar@/_15pt/[rrrrd]&&&   L_k (Y)  \underset{L_k (X)}{\coprod} X_{k,\bullet} \ar@{-->}[rd]|-{\varphi}& \\  
&&&&Y_{k,\bullet} }}
\end{equation}
$$
\left(\text{resp.}\qquad\vcenter{\xymatrix@-15pt{ 
X_{k,\bullet} \ar@{-->}[rd]|-\varphi  \ar@/^15pt/[rrrrd] \ar@/_15pt/[ddddr]_-{f_{k,\bullet}} &&&&\\
& M_k (X)  \underset{M_k (Y)}{\times} Y_{k,\bullet} \ar[rrr]\ar[ddd]&&& M_k(X)\ar[ddd]^-{\mathrm{lim}\,f_{l,\bullet}}\\ 
&&&& \\&&&& \\
&Y_{k,\bullet} \ar[rrr]&&&M_k(Y)}}\right)
$$
est une cofibration (resp. une fibration) de $(\simp_{0},{\bf W}^{red}_{n}, {\bf mono},{\bf fib}_{n}^{red})$ pour tout $k\geq 0$.

Vu que d'après le Corollaire \ref{cocommm} les colimites (resp. limites) de \eqref{LM} et \eqref{LM2} sont de colimites (resp. limites) dans $\simp$, on déduit qu'un morphisme $f\colon\xymatrix@-5pt{X\ar[r]&Y}$ de ${\bf MS}$ est une cofibration de Reedy si et seulement si, son image dans la catégorie des ensembles bisimpliciaux $\ssimp\cong\simp^{\Delta^{op}}$ est une cofibration de Reedy. En particulier un morphisme $f\colon\xymatrix@-5pt{X\ar[r]&Y}$ de ${\bf MS}$ est une cofibration de Reedy si et seulement si les fonctions $f_{p,q}$ sont injectives pour $p,q\geq 0$. Donc $({\bf MS},{\bf W}^{(2)}_{n}, {\bf mono},{\bf fib}_{n}^{(2)})$ est bien une catégorie de modèles $\big($la catégorie de modèles de Reedy des objets simpliciaux de la catégorie de modèles $(\simp_{0},{\bf W}^{red}_{n}, {\bf mono},{\bf fib}_{n}^{red})$ $\big)$.

D'un autre côté, remarquons que d'après les isomorphismes:
\begin{equation}
\begin{split}
M_p (X)_q\,\cong\, \underset{\underset{\text{inj. non id.}}{l\to p}}{\mathrm{lim}} \, X_{l,q} \,&\cong\, 
\underset{\underset{\text{inj. non id.}}{l\to p}}{\mathrm{lim}} \, \bigg(\mathrm{Hom}_{\ssimp}\big(\Delta^l\boxtimes\Delta^q, X\big) \bigg)\\
\,&\cong\, 
 \mathrm{Hom}_{\ssimp}\Bigg(\Big(\underset{\underset{\text{inj. non id.}}{l\to p}}{\mathrm{colim}}\Delta^l\Big)\boxtimes\Delta^q, X\Bigg) \\
\,&\cong\, 
\mathrm{Hom}_{\ssimp}\big(\partial\Delta^p\boxtimes\Delta^q, X\big)\\
\,&=\,  \underline{\mathrm{Hom}}^{(2)}_{\ssimp}\big(\mathrm{p}_1^* \partial\Delta^p , X\big)_q\,,
\end{split}
\end{equation}
un objet $X$ de ${\bf MS}$ est un objet fibrant de Reedy si et seulement si, le morphisme d'ensembles simpliciaux réduits:
$$
\xymatrix{\underline{\mathrm{Hom}}^{(2)}_{\ssimp}\big(\mathrm{p}_1^* \Delta^p , X\big) \ar[r]^-{} &  \underline{\mathrm{Hom}}^{(2)}_{\ssimp}\big(\mathrm{p}_1^* \partial\Delta^p , X\big)}
$$
est une fibration de la catégorie de modèles $(\simp_{0},{\bf W}^{red}_{n}, {\bf mono},{\bf fib}_{n}^{red})$ pour tout $p\geq 0$. 

Finalement, vu qu'on a déjà les adjonctions \eqref{adjoiii}, pour montrer que $({\bf MS},{\bf W}^{(2)}_{n}, {\bf mono},{\bf fib}_{n}^{(2)})$ munie de l'enrichissement $\underline{\mathrm{Hom}}^{(2)}_{{\bf MS}}$ admette une structure de catégorie de modèles simplicial d'ordre $n$, il faut vérifier la propriété suivante:
\begin{quote}
Si $\xymatrix@C-8pt{A\ar[r]^j&B}$ est un monomorphisme d'ensembles simpliciaux et $\xymatrix@C-8pt{X\ar[r]^q&Y}$ est un monomorphisme de ${\bf MS}$, alors le morphisme $\varphi$ dans le diagramme somme amalgamée de ${\bf MS}$ suivant:
\begin{equation}
\xymatrix@R-10pt@C-10pt{
X\overset{(2)}{\wedge} A \ar[dd]_{X\overset{(2)}{\wedge} \,j}\ar[rr]^{q\,\overset{(2)}{\wedge} A} && 
Y \overset{(2)}{\wedge} A \ar[dd]\ar@/^15pt/[rddd]^{Y\overset{(2)}{\wedge} \,j}&\\\\
X\overset{(2)}{\wedge} B\ar[rr]\ar@/_15pt/[rrrd]_{q\,\overset{(2)}{\wedge} B}&&
X \overset{(2)}{\wedge} B\underset{X\overset{(2)}{\wedge} A}{\bigsqcup}Y\overset{(2)}{\wedge} A\ar@{-->}[rd]|(.65)\varphi&\\
&&&Y \overset{(2)}{\wedge} B\,,}
\end{equation}
soit un monomorphisme, lequel appartient à ${\bf W}_{n}^{(2)}$ si $j$ est une $\infty$-équivalence faible d'ensembles simpliciaux, ou si $q$ appartient à ${\bf W}_{n}^{(2)}$.
\end{quote}

Dans ce but remarquons que le foncteur d'inclusion canonique $\xymatrix@C-5pt{{\bf MS} \, \ar@{^(->}[r] & \ssimp_*}$ commute aux petites colimites et produits; donc d'après le Lemme \ref{deuxilemme} il suffit de montrer que la catégorie des ensembles bisimpliciaux $\ssimp$ munie de l'enrichissement $\underline{\mathrm{Hom}}_{\ssimp}^{(2)}$ défini dans \eqref{enrri12} plus haut, est une catégorie de modèles simpliciale d'ordre $n$ avec la structure de Reedy, c'est-à-dire lorsque:
\begin{align*}
\big\{\,\text{équivalences faibles}\,\big\} \quad &=\quad 
\big\{\, f:\xymatrix@-14pt{X\ar[r]&Y}\,\big|\; \text{$f_{p,\bullet}$ est une $n$-équivalence}\\
 &\phantom{=}\;\;\quad \text{\phantom{$\{\,$}faible d'ensembles simpliciaux  $\forall\; p\geq 0$}\,\big\}\,,\\
\big\{\,\text{cofibrations}\,\big\} \quad &=\quad  \big\{\,\text{monomorphismes}\,\big\}\,.
\end{align*}

Explicitement, il suffit de montrer la propriété suivante:
\begin{quote}
Si $\xymatrix@C-8pt{A\ar[r]^j&B}$ est un monomorphisme d'ensembles simpliciaux et $\xymatrix@C-8pt{X\ar[r]^q&Y}$ est un monomorphisme d'ensembles bisimpliciaux, alors le morphisme $\varphi$ dans le diagramme somme amalgamée de $\ssimp$ suivant:
\begin{equation}\label{wwwdiag}
\xymatrix@R-10pt@C-10pt{
X\times p_2^*A \ar[dd]_{X\times \,p_2^*j}\ar[rr]^{q\,\times p_2^*A} && 
Y\times p_2^*A \ar[dd]\ar@/^15pt/[rddd]^{Y\times \,p_2^*j}&\\\\
X\times p_2^*B\ar[rr]\ar@/_15pt/[rrrd]_{q\,\times p_2^*B}&&
X\times p_2^*B\underset{X\times p_2^*A}{\bigsqcup}Y\times p_2^*A\ar@{-->}[rd]|(.65)\varphi&\\
&&&Y\times p_2^*B\,,}
\end{equation}
est un monomorphisme d'ensembles bisimpliciaux, lequel appartient à: 
\begin{equation}\label{wwred}
\big\{\, f:\xymatrix@-14pt{X\ar[r]&Y}\,\big|\; \text{$f_{p,\bullet}$ est une $n$-équivalence faible d'ensembles simpliciaux  $\forall\; p\geq 0$}\,\big\}
\end{equation}
si $j$ est une $n$-équivalence faible d'ensembles simpliciaux, ou si $q$ appartient à \eqref{wwred}.
\end{quote}

Vu que $\ssimp$ est une catégorie de préfaisceaux, si $j$ et $q$ sont des monomorphismes tous les morphismes du diagramme \eqref{wwwdiag} sont aussi des monomorphismes. En plus remarquons que $p_2^*(j)$ appartient à \eqref{wwred} si $j$ est une $n$-équivalences faible des ensembles simpliciaux

On déduit la propriété désiré du fait que les monomorphismes qui appartient à \eqref{wwred} sont stables par cochangement de base et la famille \eqref{wwred} satisfait la propriété de deux-sur-trois.
\end{proof}

Dans l'énoncé suivant on donne des conditions suffisantes pour qu'un pré-monoïde de Segal soit un objet fibrant de la catégorie de modèles $({\bf MS},{\bf W}^{(2)}_{n}, {\bf mono},{\bf fib}_{n}^{(2)})$.

\begin{lemme}\label{fibrobj}
Si $0\leq n\leq \infty$ et $X$ est un objet de la catégorie ${\bf MS}$ vérifiant les trois conditions suivantes:
\begin{enumerate}
\item $X_{p,\bullet}$ est un $n$-groupoïde de Kan pour tout $p\geq 0$, c'est-à-dire si $p\geq 0$, $q\geq 2$ et $0\leq k\leq q$ la fonction:
\begin{equation}\label{ingrp}
\xymatrix{\mathrm{Hom}\big(\Delta^p\boxtimes\Delta^q,X\big)  \ar[r] &
\mathrm{Hom}\big(\Delta^p\boxtimes\Lambda^{q,k},X\big)}
\end{equation}
est surjective pour $2\leq q\leq n$ et bijective pour $q\geq n+1$.
\item Pour $p\geq 2$, $2\leq q\leq n$ et $0\leq k\leq q$ la fonction:
\begin{equation}\label{ingrpbord}
\xymatrix{\mathrm{Hom}\big(\partial\Delta^p\boxtimes\Delta^q,X\big)  \ar[r] &
\mathrm{Hom}\big(\partial\Delta^p\boxtimes\Lambda^{q,k},X\big)}
\end{equation}
est surjective.
\item Pour $p\geq 1$, $2\leq q\leq n$ et $0\leq k\leq q$ la fonction:
\begin{equation}
\xymatrix{\mathrm{Hom}\big(\Delta^p\boxtimes\Delta^q,X\big)  \ar[r] &
\mathrm{Hom}\big(\partial\Delta^p\boxtimes\Delta^{q},X\big) \underset{\mathrm{Hom}\big(\partial\Delta^p\boxtimes\Lambda^{q,k},X\big)}{\times}
\mathrm{Hom}\big(\Delta^p\boxtimes\Lambda^{q,k},X\big) }
\end{equation}
est surjective.
\end{enumerate}
alors $X$ est un objet fibrant de la catégorie de modèles $({\bf MS},{\bf W}^{(2)}_{n}, {\bf mono},{\bf fib}_{n}^{(2)})$ de la Proposition \ref{reedyMS}.
\end{lemme}
\begin{proof}
Soit $0\leq n\leq \infty$ et supposons que $X$ est un objet de la catégorie ${\bf MS}$ vérifiant les conditions (i), (ii) et (iii). Montrons qu'on a aussi:
\begin{enumerate}
\item[(iv)] Pour $2\leq q\leq n$ et $0\leq k\leq q$ la fonction:
$$
\xymatrix{\mathrm{Hom}\big(\partial\Delta^1\boxtimes\Delta^q,X\big)  \ar[r] &
\mathrm{Hom}\big(\partial\Delta^1\boxtimes\Lambda^{q,k},X\big)}
$$
est surjective.
\item[(v)] Pour $p\geq 1$ et $q\geq n+1$ et $0\leq k\leq q$ la fonction:
$$
\xymatrix{\mathrm{Hom}\big(\partial\Delta^p\boxtimes\Delta^q,X\big)  \ar[r] &
\mathrm{Hom}\big(\partial\Delta^p\boxtimes\Lambda^{q,k},X\big)}
$$
est bijective.
\item[(vi)] Pour $p\geq 1$, $q\geq n+1$ et $0\leq k\leq q$ la fonction:
$$
\xymatrix{\mathrm{Hom}\big(\Delta^p\boxtimes\Delta^q,X\big)  \ar[r] &
\mathrm{Hom}\big(\partial\Delta^p\boxtimes\Delta^{q},X\big) \underset{\mathrm{Hom}\big(\partial\Delta^p\boxtimes\Lambda^{q,k},X\big)}{\times}
\mathrm{Hom}\big(\Delta^p\boxtimes\Lambda^{q,k},X\big) }
$$
est bijective.
\end{enumerate}

En effet, vu que $\partial\Delta^1\cong \Delta^0\sqcup \Delta^0$ et la fonction \eqref{ingrp} est surjective si $p= 0$, $2\leq q\leq n$ et $0\leq k\leq q$, on déduit (iv). D'un autre vu que \eqref{ingrp} est bijective si $p\geq 0$, $q\geq n+1$ et $0\leq k\leq q$, il se suit (v). Finalement vu que d'après (i) et (v) les fonctions \eqref{ingrp} et \eqref{ingrpbord} sont bijectives si $p\geq 1$, $q\geq n+1$ et $0\leq k\leq q$, donc (vi).

Remarquons d'un autre côté que si $X$ vérifie les conditions (i)-(vi) ci-dessus, alors pour tout $p\geq 0$ les ensembles simpliciaux réduits 
$\underline{\mathrm{Hom}}^{(2)}_{\ssimp}\big(\mathrm{p}_1^* \Delta^p , X\big)$ et $\underline{\mathrm{Hom}}^{(2)}_{\ssimp}\big(\mathrm{p}_1^* \partial\Delta^p , X\big)$ sont des objets fibrants de la catégorie de modèles $(\simp_{0},{\bf W}^{red}_{n}, {\bf mono},{\bf fib}_{n}^{red})$ et le morphisme:
\begin{equation}\label{lemmadieu}
\xymatrix{\underline{\mathrm{Hom}}^{(2)}_{\ssimp}\big(\mathrm{p}_1^* \Delta^p , X\big) \ar[r]^-{} &  \underline{\mathrm{Hom}}^{(2)}_{\ssimp}\big(\mathrm{p}_1^* \partial\Delta^p , X\big)}
\end{equation} 
vérifie la propriété de relèvement à droite par rapport à l'ensemble de morphismes:
$$
\bigg\{ \quad \xymatrix@+45pt{\Lambda^{q,k}/{\bf sq}_{0}\Lambda^{q,k}\; \ar@{^(->}[r]^-{\alpha^{q-1,k}/{\bf sq}_{0}\alpha^{q-1,k}}& \, \Delta^{q}/{\bf sq}_{0}\Delta^q} \quad \bigg| \quad  q\geq 2\,,\;\, 0\leq k\leq q \quad \bigg\}\,.
$$

Donc d'après la Proposition \ref{identfibra}, le morphisme \eqref{lemmadieu} est une fibration de la catégorie de modèles $(\simp_{0},{\bf W}^{red}_{n}, {\bf mono},{\bf fib}_{n}^{red})$ pour tout $p\geq 0$; c'est-à-dire $X$ est un objet fibrant de la catégorie de modèles $({\bf MS},{\bf W}^{(2)}_{n}, {\bf mono},{\bf fib}_{n}^{(2)})$.
\end{proof}

\renewcommand{\thesubsection}{\S\thesection.\arabic{subsection}}
\subsubsection{}\;
\renewcommand{\thesubsection}{\thesection.\arabic{subsection}}

On va maintenant localiser la catégorie de modèles $({\bf MS},{\bf W}^{(2)}_{n}, {\bf mono},{\bf fib}_{n}^{(2)})$ de la Proposition \ref{reedyMS} par rapport aux équivalences faibles diagonales des ensembles bisimpliciaux; on obtiendra ainsi la catégorie de modèles simpliciale de Dugger (voir \cite{dugger}) des objets simpliciaux de la catégorie de modèles $(\simp_{0},{\bf W}^{red}_{n}, {\bf mono},{\bf fib}_{n}^{red})$ (voir la Proposition \ref{moduno} ci-dessous).

Rappelons pour commencer qu'à partir du foncteur diagonal de la catégorie des simplexes $\xymatrix@R=2pt{\Delta\ar[r]^-{\delta} &\Delta \times \Delta\,\\  \text{\scriptsize{$[n]$}} \ar@{}[r]|-{\longmapsto} &\text{\scriptsize{$\big([n],[n]\big)$}}}$, on construit une adjonction:
\begin{equation}\label{diago} 
\ssimp\xymatrix@C+15pt{\phantom{a}\ar@/^12pt/[r]^-{\delta^*}\ar@{}[r]|-{\perp}&\ar@/^12pt/[l]^-{\delta_*}\phantom{a}}
\simp\,,
\end{equation}
où $\delta^*(X) \, = \, X\circ\delta^{op}$ est dit \emph{l'ensemble simplicial diagonal} de $X$. 

Remarquons qu'il est possible de définir $\delta_*(A)_{p,q} = \mathrm{Hom}_{\simp}\big(\Delta^p\times \Delta^q , A \big)$. En effet le foncteur $\delta^*$ commute aux colimites et on a un isomorphisme naturel:
$$
\mathrm{Hom}_{\simp}\big(\delta^*(\Delta^p\boxtimes\Delta^q),A\big)\; = \;
\mathrm{Hom}_{\simp}\big(\Delta^p\times\Delta^q,A\big)\; = \;
\delta_*(A)_{p,q}\;\cong\;\mathrm{Hom}_{\ssimp}\big(\Delta^p\boxtimes\Delta^q,\delta_*A\big)
$$
pour tout ensemble simplicial $A$.

On déduit de \eqref{diago} une adjonction:
\begin{equation}\label{diagoadj}
{\bf MS} \xymatrix@C+15pt{\phantom{a}\ar@/^12pt/[r]^-{diag}\ar@{}[r]|-{\perp}&\ar@/^12pt/[l]^-{r}\phantom{a}}
\simp_0 \,;
\end{equation}
où $diag(X)_n=X_{n,n}$ et $r(A)_{p,q}=\mathrm{Hom}_{\simp_0}\Big(\big(\Delta^p\times\Delta^q\big)\big/\big(\Delta^p\times{\bf sq}_0(\Delta^q)\big),A\Big)$.

En particulier:
\begin{align*}
r(A)_{p,q} \, & = \, \mathrm{Hom}_{\simp_0}\bigg(\big(\Delta^p\times \Delta^q\big)\big/\big(\Delta^p\times {\bf sq}_0(\Delta^q)\big), \, A\bigg) & \\
& \cong \, \mathrm{Hom}_{\simp_0}\bigg(  \Big(\big(\Delta^q/{\bf sq}_0(\Delta^q)\big)\times \Delta^p\Big)\Big/\big(\star \times \Delta^p\big) , \, A\bigg) &\\
& = \, \underline{\mathrm{Hom}}_{\simp_0} \big(\Delta^q \big/ {\bf sq}_0(\Delta^q) , A  \big)_p & \text{(Voir la définition \eqref{homred})}\,;
\end{align*}
c'est-à-dire:
\begin{equation}\label{spacelaac}
r(A)_{p,\bullet} \, \cong \, A^{\overset{\circ}{\wedge}\Delta^p}
\qquad\text{et}\qquad
r(A)_{\bullet,q} \, \cong \, \underline{\mathrm{Hom}}_{\simp_0} \big(\Delta^q \big/ {\bf sq}_0(\Delta^q) , A  \big) \,.
\end{equation}


\begin{proposition}\label{moduno}
Si $0\leq n\leq \infty$ et $0\leq i\leq 1$, la catégorie des pré-monoïdes de Segal ${\bf MS}$ munie de l'enrichissement $\underline{\mathrm{Hom}}^{(i)}_{{\bf MS}}$ admet une structure de catégorie de modèles simplicial pointée d'ordre $n-1$, lorsque:
\begin{align*}
\big\{\,\text{équivalences faibles}\,\big\} \quad &=\quad 
\big\{\, f:\xymatrix@-14pt{X\ar[r]&Y}\,\big|\; \text{$\mathrm{diag}(f)$ est une $n$-équivalence}\\
 &\phantom{=}\;\;\quad \text{\phantom{$\{\,$}faible d'ensembles simpliciaux}\,\big\}\;\,=\;\, {\bf W}_{n}^{diag}\,,\\
\big\{\,\text{cofibrations}\,\big\} \quad &=\quad  \big\{\,\text{monomorphismes}\,\big\} \;\;\,= \;\,  {\bf mono}\,,\\
\big\{\,\text{fibrations}\,\big\} \quad &=\quad  \big\{\,\text{morphismes avec la propriété de relèvement }\\
 &\phantom{=}\;\;\quad \text{\phantom{$\{\,$}à droite par rapport à ${\bf mono}\,\cap\,{\bf W}_{n}^{diag}$}\,\big\}\;\,= \;\, {\bf fib}_{n}^{diag}\,.
\end{align*}

En plus, un pré-monoïde de Segal $X$ est fibrant dans $({\bf MS},{\bf W}^{diag}_{n}, {\bf mono},{\bf fib}_{n}^{diag})$ si et seulement si, $X$ est un objet fibrant de la catégorie de modèles $({\bf MS},{\bf W}^{(2)}_{n}, {\bf mono},{\bf fib}_{n}^{(2)})$ de la Proposition \ref{reedyMS} et tous les morphismes d'ensembles simpliciaux $\xymatrix@C-8pt{X_{k,\bullet} \ar[r]^{\varphi^*} & X_{l,\bullet}}$ induisent par les morphismes $\xymatrix@C-8pt{[l]\ar[r]^-{\varphi}&[k]}$ de $\Delta$ sont de $n$-équivalences faibles.
\end{proposition}
\begin{proof}
On sait que la catégorie de modèles $({\bf MS},{\bf W}^{(2)}_{n}, {\bf mono},{\bf fib}_{n}^{(2)})$ de la Proposition \ref{reedyMS} est une catégorie de modèles propre à gauche (les cofibrations sont les monomorphismes) et combinatoire (c'est la catégorie des diagrammes d'une catégorie combinatoire sur une catégorie de Reedy). Donc d'après le Théorème 4.7 de \cite{barwick} on peut considérer sa localisation de Bousfield à gauche par rapport à l'ensemble des morphismes:
\begin{equation}
S^{diag}_{n}\,=\,\Big\{ \quad 
\xymatrix@C+6pt{\big(\Delta^{l} \boxtimes \mathbb{S}^m\big)\big/\big(\Delta^{l} \boxtimes\star\big)\ar[r]^-{\varphi_*\boxtimes \mathrm{id}}&\big(\Delta^{k} \boxtimes \mathbb{S}^m\big)\big/\big(\Delta^{k} \boxtimes\star\big)} \quad \Big| \quad 1\leq m\leq n\,,\; [l]\overset{\varphi}{\to}[k] \,\in\,\Delta \quad \Big\}\,.
\end{equation}

Dans le deux Lemme suivants ont pose $[\,\cdot\,,\,\cdot\,]^{(2)}_{n}$ pour noter l'ensemble de morphismes dans la catégorie de fractions ${\bf MS}\big[\big({\bf W}^{(2)}_{n}\big)^{-1}\big]$.

\begin{lemme}\label{lolocal3}
Si $0\leq n\leq \infty$ un pré-monoïde de Segal $Z$ est un objet $S_{n}^{diag}$-local de la catégorie de modèles $({\bf MS},{\bf W}^{(2)}_{n}, {\bf mono},{\bf fib}_{n}^{(2)})$ si et seulement si, tous les morphismes d'ensembles simpliciaux $\xymatrix@C-8pt{Z_{k,\bullet} \ar[r]^{\varphi^*} & Z_{l,\bullet}}$ induisent par les morphismes $\xymatrix@C-8pt{[l]\ar[r]^-{\varphi}&[k]}$ de $\Delta$ sont de $n$-équivalences faibles.
\end{lemme}
\begin{proof}
Pour commencer remarquons que si $l\geq 0$, on a une adjonction:
$$
\xymatrix@C+15pt{
{\bf MS} \;\phantom{A}\ar@{}[r]|-{\perp\phantom{A}}
\ar@{<-}@/^14pt/[r]^-{(\Delta^l\boxtimes \cdot)/(\Delta^l \boxtimes\star)}
\ar@/_14pt/[r]_-{(\,\cdot\,)_{l,\bullet}}& \phantom{A}\simp_0}\,.;
$$
en particulier pour tout pré-monoïde de Segal $W$ on a une chaîne d'isomorphismes:
\begin{align*}
\underline{\mathrm{Hom}}_{\simp_0}\big(\mathbb{S}^m,W_{l,\bullet}\big)_n \, 
&= \, \mathrm{Hom}_{\simp_0}\Big((\mathbb{S}^m\times\Delta^n)\big/(\star\times\Delta^n),W_{l,\bullet}\Big)\\ 
&\cong \, \mathrm{Hom}_{{\bf MS}}\bigg(\Delta^l\boxtimes\Big[(\mathbb{S}^m\times\Delta^n)\big/(\star\times\Delta^n)\Big]\Big/(\Delta^l\boxtimes \star) ,W\bigg) \\ 
&\cong \, \mathrm{Hom}_{{\bf MS}}\bigg(\Big[(\Delta^l\boxtimes\mathbb{S}^m)\big/(\Delta^l\boxtimes \star)\Big]\times p_2^*\Delta^n\big/(\star\times p_2^*\Delta^n),W\bigg) \\ 
&\cong \, \underline{\mathrm{Hom}}_{{\bf MS}}^{(2)}\Big((\Delta^l\boxtimes\mathbb{S}^m)\big/(\Delta^l\boxtimes \star), W\Big)_n\,;
\end{align*}
pour tout $m\geq 1$ et $n,l\geq 0$.

Vu que le foncteur $(\,\cdot\,)_{l,\bullet}\colon\xymatrix@C-10pt{{\bf MS} \ar[r] & \simp_0}$ envoi les morphismes des familles ${\bf W}^{(2)}_{n}$ et ${\bf fib}_{n}^{(2)}$ vers les familles ${\bf W}^{red}_{n}$ et ${\bf fib}_{n}^{red}$ respectivement, on déduit que si $Z$ est un pré-monoïde de Segal arbitraire et $W$ est un remplacement fibrant de $Z$ dans la catégorie de modèles $({\bf MS},{\bf W}^{(2)}_{n}, {\bf mono},{\bf fib}_{n}^{(2)})$, alors $W_{l,\bullet}$ est un remplacement fibrant de $Z_{l,\bullet}$ dans $(\simp_{0},{\bf W}^{red}_{n}, {\bf mono},{\bf fib}_{n}^{red})$ pour tout $l\geq 0$.

Il se suit que si $\xymatrix@C-8pt{[l]\ar[r]^-{\varphi}&[k]}$ est un morphisme de $\Delta$ on a un diagramme commutatif:
$$
\xymatrix@C+10pt@R=3pt{ 
\big[(\Delta^k\boxtimes\mathbb{S}^m)\big/(\Delta^k\boxtimes \star),Z\big]_n^{(2)}
\ar[r]^-{(\varphi_*\boxtimes\mathrm{id})^*}\ar@{}[d]|-{\text{\rotatebox[origin=c]{270}{\Large $\cong$}}}&
\big[(\Delta^l\boxtimes\mathbb{S}^m)\big/(\Delta^l\boxtimes \star),Z\big]_n^{(2)}
\ar@{}[d]|-{\text{\rotatebox[origin=c]{270}{\Large $\cong$}}}\\
\pi_0\Big(\underline{\mathrm{Hom}}_{{\bf MS}}^{(2)}\big((\Delta^k\boxtimes\mathbb{S}^m)\big/(\Delta^k\boxtimes \star), W\big)\Big)
\ar@{}[d]|-{\text{\rotatebox[origin=c]{270}{\Large $\cong$}}}&
\pi_0\Big(\underline{\mathrm{Hom}}_{{\bf MS}}^{(2)}\big((\Delta^l\boxtimes\mathbb{S}^m)\big/(\Delta^l\boxtimes \star), W\big)\Big)
\ar@{}[d]|-{\text{\rotatebox[origin=c]{270}{\Large $\cong$}}}\\
\pi_0\Big(\underline{\mathrm{Hom}}_{\simp_0}\big(\mathbb{S}^m,W_{k,\bullet}\big)\Big)
\ar@{}[d]|-{\text{\rotatebox[origin=c]{270}{\Large $\cong$}}}&
\pi_0\Big(\underline{\mathrm{Hom}}_{\simp_0}\big(\mathbb{S}^m,W_{l,\bullet}\big)\Big)
\ar@{}[d]|-{\text{\rotatebox[origin=c]{270}{\Large $\cong$}}}\\
\big[\mathbb{S}^m,Z_{k,\bullet}\big]_n^{red}
\ar@{}[d]|-{\text{\rotatebox[origin=c]{270}{\Large $\cong$}}}&
\big[\mathbb{S}^m,Z_{l,\bullet}\big]_n^{red}
\ar@{}[d]|-{\text{\rotatebox[origin=c]{270}{\Large $\cong$}}}\\
\pi_m\big(Z_{k,\bullet}\big)
\ar[r]_-{\varphi^*}&
\pi_m\big(Z_{l,\bullet}\big)\,,
}
$$
ce qui montre l'énoncé désiré.
\end{proof}

Déterminons les équivalences faibles $S_{n}^{diag}$-locales: 

\begin{lemme} \label{lolocal4}
Si $0\leq n\leq \infty$, un morphisme de pré-monoïdes de Segal $f\colon\xymatrix@C-5pt{X\ar[r]&Y}$ est une équivalence faible $S_{n}^{diag}$-local de la catégorie de modèles $({\bf MS},{\bf W}^{(2)}_{n}, {\bf mono},{\bf fib}_{n}^{(2)})$, si et seulement si le morphisme d'ensembles simpliciaux réduits $\mathrm{diag}(f)\colon\xymatrix@C-5pt{\mathrm{diag}(X)\ar[r]&\mathrm{diag}(Y)}$ est une $n$-équivalence faible d'ensembles simpliciaux.
\end{lemme}
\begin{proof}
Commençons par noter que l'adjonction \eqref{diagoadj} est une adjonction de Quillen entre les catégories de modèles $({\bf MS},{\bf W}^{(2)}_{n}, {\bf mono},{\bf fib}_{n}^{(2)})$ et $(\simp_{0},{\bf W}^{red}_{n}, {\bf mono},{\bf fib}_{n}^{red})$. En effet, le foncteur $diag$ respect les monomorphismes évidemment. 

D'un autre pour montrer que $diag\big({\bf W}_{n}^{(2)}\big)\subset {\bf W}_n^{red}$ remarquons par ailleurs que si $f:\xymatrix@-14pt{X\ar[r]&Y}$ est un morphisme de ${\bf MS}$, on a que $f$ appartient à ${\bf W}^{(2)}_n$ si et seulement si son image par le foncteur d'inclusion $\xymatrix@C-8pt{{\bf MS}\ar[r]&\ssimp}$ appartient à l'ensemble des morphismes:
$$\big\{\, f:\xymatrix@-14pt{W\ar[r]&Z}\,\big|\; \text{$f_{p,\bullet}$ appartient à ${\bf W}_n$ pour tout $p\geq0$}\big\}\,,$$
et de façon analogue le morphisme d'ensembles simpliciaux réduits $diag(f)$ est une $n$-équivalence faible si et seulement si son image par le foncteur d'inclusion $\xymatrix@C-8pt{\simp_0\ar[r]&\simp}$ est une $n$-équivalence faible d'ensembles simpliciaux. 

Donc $diag\big({\bf W}_{n}^{(2)}\big)\subset {\bf W}_n^{red}$ parce que d'après un résultat de Denis-Charles Cisinski (voir le Corollaire 2.3.17 et le Théorème 1.4.3 de \cite{ci}),  pour tout $0\leq n\leq \infty$ l'image par le foncteur diagonal:
$$
\xymatrix@C+15pt{\ssimp\ar[r]^-{\delta^*}&\simp}
$$
de l'ensemble de morphismes $\big\{\, f:\xymatrix@-14pt{W\ar[r]&Z}\,\big|\; \text{$f_{p,\bullet}$ appartient à ${\bf W}_n$ pour tout $p\geq0$}\big\}$ est contenu dans l'ensemble des $n$-équivalences faibles. 

Remarquons d'un autre côté que pour tout ensemble simplicial réduit $A$ on a que:
$$
r(A)_{p,q} \,  = \, \mathrm{Hom}_{\simp_0}\bigg(\big(\Delta^p\times \Delta^q\big)\big/\big(\Delta^p\times {\bf sq}_0(\Delta^q)\big), \, A\bigg) \, = \, \Big( A^{\overset{\circ}{\wedge}\Delta^p}\Big)_q    \qquad \quad \text{(voir \eqref{adjreduit})\,;}
$$
alors si on pose $p^*\colon\xymatrix{\simp_0\ar[r]&\simp_0^{\Delta^{op}} \, = \, {\bf MS}}$ pour noter le foncteur préfaisceau constant $p^*(A)_{p,q} = A_{q} \cong \Big( A^{\overset{\circ}{\wedge}\Delta^0}\Big)_q$, on définit aussi-tôt une transformation naturelle $\xymatrix@C+8pt{\simp_0\rtwocell^{p^*}_r{\eta}&{\bf MS}}$ par le carré commutatif:
$$
\xymatrix@C+5pt@R=5pt{
p^*(A)_{p,q} \ar[r]^{(\eta_A)_{p,q}} \ar@{}[d]|-{\mathrel{\reflectbox{\rotatebox[origin=c]{90}{$\cong$}}}} & r(A)_{p,q} \ar@{}[d]|-{\mathrel{\reflectbox{\rotatebox[origin=c]{90}{$=$}}}}\\ 
\mathrm{Hom}_{\simp_0}\Big(\big(\Delta^0\times\Delta^q\big)\big/\big(\Delta^0\times {\bf sq}_0(\Delta^q)\big),A\Big)\ar@{}[d]|-{\mathrel{\reflectbox{\rotatebox[origin=c]{90}{$=$}}}}&
\mathrm{Hom}_{\simp_0}\Big(\big(\Delta^p\times\Delta^q\big)\big/\big(\Delta^p\times {\bf sq}_0(\Delta^q)\big),A\Big)\ar@{}[d]|-{\mathrel{\reflectbox{\rotatebox[origin=c]{90}{$=$}}}}
\\\big(A^{\overset{\circ}{\wedge}\Delta^0} \big)_q \ar[r]_{(A^{\overset{\circ}{\wedge}\varphi})_q}& \big(A^{\overset{\circ}{\wedge}\Delta^p}\big)_q\,,
}
$$
où $\varphi\colon\xymatrix@C-8pt{\Delta^p\ar[r]&\Delta^0}$ est le morphisme canonique, lequel est une $\infty$-équivalence faible. 

Il se suit que pour tout objet fibrant $A$ de $(\simp_{0},{\bf W}^{red}_{n}, {\bf mono},{\bf fib}_{n}^{red})$, le morphisme  d'ensembles simpliciaux réduits qu'on vient de définir $(\eta_A)_{p,\bullet}\colon\xymatrix@-3pt{p^*(A)_{p,\bullet} \ar[r] & r(A)_{p,\bullet}}$ est une $\infty$-équivalence faible pour tout $p\geq 0$, c'est-à-dire le morphisme de pré-monoïdes de Segal $\eta_A\colon\xymatrix@-5pt{p^*(A)\ar[r] & r(A)}$ est une équivalence faible de la catégorie de modèles $({\bf MS},{\bf W}^{(2)}_{n}, {\bf mono},{\bf fib}_{n}^{(2)})$. 

On conclut que si $A$ est un objet fibrant de $(\simp_{0},{\bf W}^{red}_{n}, {\bf mono},{\bf fib}_{n}^{red})$ et $X$ est un pré-monoïde de Segal quelconque, il y a des bijections naturelles: 
$$
\big[diag(X),A\big]^{red}_n \, \cong \, \big[{\bf L}diag(X),A\big]^{red}_n \, \cong \, \big[X,{\bf R}r(A)\big]^{(2)}_n \, \cong \, \big[X,p^*(A)\big]^{(2)}_n\,.
$$

Donc un morphisme de pré-monoïdes de Segal $f\colon\xymatrix@C-5pt{X\ar[r]&Y}$ satisfait que $diag(f)$ soit une $n$-équivalence faible d'ensembles simpliciaux réduits, si et seulement si pour tout objet fibrant $A$ de la catégorie de modèles $(\simp_{0},{\bf W}^{red}_{n}, {\bf mono},{\bf fib}_{n}^{red})$ la fonction:
$$
\xymatrix{\big[Y,p^*(A)\big]^{(2)}_n  \ar[r]^-{f^*} & \big[X,p^*(A)\big]^{(2)}_n}
$$
est bijective.

Vu que $p^*(A)$ est un objet $S^{diag}_n$-local pour tout ensemble simplicial réduit $A$, il se suit que si $f$ est une équivalence faible $S_{n}^{diag}$-local de la catégorie de modèles $({\bf MS},{\bf W}^{(2)}_{n}, {\bf mono},{\bf fib}_{n}^{(2)})$ alors $diag(f)$ est une $n$-équivalence faible d'ensembles simpliciaux réduits.

Pour montrer la réciproque remarquons que si $Z$ est un objet $S^{diag}_n$-local et fibrant de la catégorie de modèles $({\bf MS},{\bf W}^{(2)}_{n}, {\bf mono},{\bf fib}_{n}^{(2)})$, alors le morphisme $\xymatrix@-5pt{p^*(Z_{0,\bullet})\ar[r]&Z}$ défini par les fonctions $\xymatrix@C-5pt{Z_{0,q}\ar[r]^-{s_0\dots s_0}&Z_{p,q}}$ appartient à ${\bf W}^{(2)}_{n}$ où $Z_{0,\bullet}$ est un objet fibrant de $(\simp_{0},{\bf W}^{red}_{n}, {\bf mono},{\bf fib}_{n}^{red})$. 

Donc on a un carre commutatif:
$$
\xymatrix@R=5pt@C+25pt{
\big[Y,Z\big]^{(2)}_n \ar@{}[d]|-{\mathrel{\reflectbox{\rotatebox[origin=c]{90}{$\cong$}}}} \ar[r]^-{f^*} & \big[X,Z\big]^{(2)}_n\ar@{}[d]|-{\mathrel{\reflectbox{\rotatebox[origin=c]{90}{$\cong$}}}}\\
\big[Y,p^*(Z_{0,\bullet})\big]^{(2)}_n \ar@{}[d]|-{\mathrel{\reflectbox{\rotatebox[origin=c]{90}{$\cong$}}}}  &\big[X,p^*(Z_{0,\bullet})\big]^{(2)}_n\ar@{}[d]|-{\mathrel{\reflectbox{\rotatebox[origin=c]{90}{$\cong$}}}}\\
\big[diag(Y),Z_{0,\bullet}\big]^{red}_n  \ar[r]_-{diag(f)^*}&  \big[diag(X),Z_{0,\bullet}\big]^{red}_n}
$$
pour tout objet $S^{diag}_n$-local et fibrant $Z$ de la catégorie de modèles $({\bf MS},{\bf W}^{(2)}_{n}, {\bf mono},{\bf fib}_{n}^{(2)})$. Ce qui implique que si $diag(f)$ est une $n$-équivalence faible d'ensembles simpliciaux réduits, alors $f$ est une équivalence faible $S_{n}^{diag}$-local de la catégorie de modèles $({\bf MS},{\bf W}^{(2)}_{n}, {\bf mono},{\bf fib}_{n}^{(2)})$.
\end{proof}

Il se suit du Théorème 4.7 de \cite{barwick} et des Lemmes \ref{lolocal3} et \ref{lolocal4} ci-dessus qu'on a bien une catégorie de modèles $({\bf MS},{\bf W}^{diag}_{n}, {\bf mono},{\bf fib}_{n}^{diag})$, dont les objets fibrants $X$ sont les objets fibrants de la catégorie de modèles $({\bf MS},{\bf W}^{(2)}_{n}, {\bf mono},{\bf fib}_{n}^{(2)})$ de la Proposition \ref{reedyMS} vérifiant de plus que $\xymatrix@C-8pt{X_{k,\bullet} \ar[r]^{\varphi^*} & X_{l,\bullet}}$ soit une $n$-équivalence faible pour tout morphisme $\xymatrix@C-8pt{[l]\ar[r]^-{\varphi}&[k]}$ de $\Delta$ .

Si $0\leq i\leq 1$ pour montrer que la catégorie de modèles $({\bf MS},{\bf W}^{diag}_{n}, {\bf mono},{\bf fib}_{n}^{diag})$ avec les enrichissements $\underline{\mathrm{Hom}}^{(i)}_{{\bf MS}}$ est simpliciale d'ordre $n$, il faut vérifier la propriété suivante:
\begin{quote}
Si $\xymatrix@C-8pt{A\ar[r]^j&B}$ est un monomorphisme d'ensembles simpliciaux et $\xymatrix@C-8pt{X\ar[r]^q&Y}$ est un monomorphisme de ${\bf MS}$, alors le morphisme $\varphi$ dans le diagramme somme amalgamée de ${\bf MS}$ suivant:
\begin{equation}
\xymatrix@R-10pt@C-10pt{
X\overset{(i)}{\wedge} A \ar[dd]_{X\overset{(i)}{\wedge} \,j}\ar[rr]^{q\,\overset{(i)}{\wedge} A} && 
Y \overset{(i)}{\wedge} A \ar[dd]\ar@/^15pt/[rddd]^{Y\overset{(i)}{\wedge} \,j}&\\\\
X\overset{(i)}{\wedge} B\ar[rr]\ar@/_15pt/[rrrd]_{q\,\overset{(i)}{\wedge} B}&&
X \overset{(i)}{\wedge} B\underset{X\overset{(i)}{\wedge} A}{\bigsqcup}Y\overset{(i)}{\wedge} A\ar@{-->}[rd]|(.65)\varphi&\\
&&&Y \overset{(i)}{\wedge} B\,,}
\end{equation}
est un monomorphisme, lequel appartient à ${\bf W}_{n}^{diag}$ si $j$ est une $\infty$-équivalence faible d'ensembles simpliciaux, ou si $q$ appartient à ${\bf W}_{n}^{diag}$.
\end{quote}

Comme on l'a noté dans la Proposition \ref{reedyMS}, vu que le foncteur d'inclusion $\xymatrix@C-5pt{{\bf MS} \, \ar@{^(->}[r] & \ssimp_*}$ commute aux petites colimites et limites, d'après le Lemme \ref{deuxilemme} il suffit de montrer que la catégorie des ensembles bisimpliciaux $\ssimp$ munie de l'enrichissement $\underline{\mathrm{Hom}}_{\ssimp}^{(i)}$ défini dans \eqref{enrri12} plus haut, est une catégorie de modèles simpliciale d'ordre $n$ lorsque:
\begin{align*}
\big\{\,\text{équivalences faibles}\,\big\} \quad &=\quad 
\big\{\, f:\xymatrix@-14pt{X\ar[r]&Y}\,\big|\; \text{$\delta^*(f)$ est une $n$-équivalence}\\
 &\phantom{=}\;\;\quad \text{\phantom{$\{\,$}faible d'ensembles simpliciaux  $\forall\; p\geq 0$}\,\big\}\,,\\
\big\{\,\text{cofibrations}\,\big\} \quad &=\quad  \big\{\,\text{monomorphismes}\,\big\}\,;
\end{align*}
où $\xymatrix@C-5pt{\ssimp\ar[r]^-{\delta^*}&\simp}$ est le foncteur diagonal des ensembles bisimpliciaux.

Explicitement, il suffit de montrer la propriété suivante:
\begin{quote}
Si $\xymatrix@C-8pt{A\ar[r]^j&B}$ est un monomorphisme d'ensembles simpliciaux et $\xymatrix@C-8pt{X\ar[r]^q&Y}$ est un monomorphisme d'ensembles bisimpliciaux, alors le morphisme $\varphi$ dans le diagramme somme amalgamée de $\ssimp$ suivant:
\begin{equation}\label{wwwdiagMS}
\xymatrix@R-10pt@C-10pt{
X\times p_i^*A \ar[dd]_{X\times \,p_i^*j}\ar[rr]^{q\,\times p_i^*A} && 
Y\times p_i^*A \ar[dd]\ar@/^15pt/[rddd]^{Y\times \,p_i^*j}&\\\\
X\times p_i^*B\ar[rr]\ar@/_15pt/[rrrd]_{q\,\times p_i^*B}&&
X\times p_i^*B\underset{X\times p_i^*A}{\bigsqcup}Y\times p_i^*A\ar@{-->}[rd]|(.65)\varphi&\\
&&&Y\times p_i^*B\,,}
\end{equation}
est un monomorphisme d'ensembles bisimpliciaux, lequel appartient à: 
\begin{equation}\label{wwredMS}
\big\{\, f:\xymatrix@-14pt{X\ar[r]&Y}\,\big|\; \text{$\delta^*(f)$ est une $n$-équivalence faible d'ensembles simpliciaux  $\forall\; p\geq 0$}\,\big\}
\end{equation}
si $j$ est une $n$-équivalence faible d'ensembles simpliciaux, ou si $q$ appartient à \eqref{wwredMS}.
\end{quote}

Ce qui est simple à montrer vu que $\delta^*\big(p_i^*(j)\big) = j$ est une $n$-équivalence faible si $j$ l'est.

Finalement pour montrer que $({\bf MS},{\bf W}^{diag}_{n}, {\bf mono},{\bf fib}_{n}^{diag})$ avec les enrichissements $\underline{\mathrm{Hom}}^{(i)}_{{\bf MS}}$ est une catégorie de modèles simpliciale d'ordre $n-1$, d'après le Corollaire \ref{dimminus}, le Lemme \ref{quillenomegasus} et le Théorème \ref{cassi1}, il suffit de montrer que l'adjonction:
\begin{equation}\label{duggerQeq}
{\bf MS} \xymatrix@C+15pt{\phantom{a}\ar@/_12pt/[r]_-{(\,\cdot\,)_{0,\bullet}}\ar@{}[r]|-{\perp}&\ar@/_12pt/[l]_-{p^*\,=\,\text{préfaisceau constant}}\phantom{a}}
\simp_0 
\end{equation}
est une équivalence de Quillen entre $({\bf MS},{\bf W}^{diag}_{n}, {\bf mono},{\bf fib}_{n}^{diag})$  et $(\simp_{0},{\bf W}^{red}_{n}, {\bf mono},{\bf fib}_{n}^{red})$.

Pour commencer on vérifie sans peine que le foncteur $p^*$ est Quillen à gauche. D'un autre, si $A$ est un objet fibrant de la catégorie de modèles $(\simp_{0},{\bf W}^{red}_{n}, {\bf mono},{\bf fib}_{n}^{red})$, on vérifie aussi facilement que ${\bf R}(\;)_{0,\bullet}\circ {\bf L}p^* \,(A)$ est isomorphe à $A$. De la même façon si $X$ es un objet fibrant de la catégorie de modèles $({\bf MS},{\bf W}^{diag}_{n}, {\bf mono},{\bf fib}_{n}^{diag})$  on a que ${\bf L}p^*  \circ {\bf R}(\;)_{0,\bullet}\,(X)$ est isomorphe à $X$.
\end{proof}

On pose $\mathrm{Ho}_{n}({\bf MS})$ pour noter la catégorie homotopique ${\bf MS}\big[({\bf W}^{diag}_{n})^{-1}\big]$ de la catégorie de modèles $({\bf MS},{\bf W}^{diag}_{n}, {\bf mono},{\bf fib}_{n}^{diag})$ de la Proposition \ref{moduno}. Il se suit de la Proposition \ref{moduno} que les foncteurs dérivés des morphismes ${\bf R}\underline{\mathrm{Hom}}_{\bf MS}^{(i)}$ pour $0\leq i\leq 1$ sont des enrichissements forcement isomorphes de $\mathrm{Ho}_{n}({\bf MS})$ sur la catégorie cartésienne fermée des $(n-1)$-types d'homotopie $\mathrm{Ho}_{n-1}(\simp)$ (la catégories des $(n-1)$-groupoïdes). 

La catégorie $\mathrm{Ho}_{n}({\bf MS})$ ainsi enrichie sur la catégorie des $(n-1)$-groupoïdes est dit la \emph{catégorie homotopique des $n$-groupes de Segal}. Note que l'équivalence de Quillen \eqref{duggerQeq} induit une équivalence de la catégorie homotopique des $n$-groupes de Segal et la catégorie homotopique des $n$-groupes:
\begin{equation}\label{duggerQeqHO}
\mathrm{Ho}_{n}({\bf MS}) \xymatrix@C+15pt{\phantom{a}\ar@/_12pt/[r]_-{{\bf R}(\,\cdot\,)_{0,\bullet}}\ar@{}[r]|-{\simeq}&\ar@/_12pt/[l]_-{{\bf L}p^*}\phantom{a}}
\mathrm{Ho}_{n}(\simp_0)\,.
\end{equation}

Remarquons aussi que l'adjonction de Quillen \eqref{diagoadj} (note que le foncteur $diag$ est un foncteur de Quillen à gauche) est une équivalences de Quillen. En effet \ref{diagoadj} induit une équivalence de catégories:
\begin{equation}\label{diagoadjHO}
\mathrm{Ho}_{n}({\bf MS}) \xymatrix@C+15pt{\phantom{a}\ar@/_12pt/[r]_-{{\bf L}diag}\ar@{}[r]|-{\simeq}&\ar@/_12pt/[l]_-{{\bf R}r}\phantom{a}}
\mathrm{Ho}_{n}(\simp_0)
\end{equation}
parce que d'après la preuve du Lemme \ref{lolocal4} un foncteur dérivé ${\bf R}r$ est isomorphe à un foncteur dérivé ${\bf L}p^*$ (note que le foncteur $p^*$ envoi les $n$-équivalences faibles réduites vers des $n$-équivalences faibles diagonales).

L'énoncé suivant est déduit aussi-tôt du Lemme \ref{fibrobj} et la Proposition \ref{moduno}:

\begin{lemme}\label{fibrafifi}
Si $1\leq n\leq \infty$ et $X$ est un objet de la catégorie ${\bf MS}$ vérifiant les conditions suivantes:
\begin{enumerate}
\item Si $\xymatrix@C-10pt{[n]\ar[r]^{\varphi}&[m]}$ est un morphisme quelconque de la catégorie des simplexes $\Delta$, le morphisme induit $\xymatrix{X_{m,\bullet}\ar[r]^-{\varphi^*} & X_{n,\bullet}}$ est une $n$-équivalence faible d'ensembles simpliciaux.
\item $X_{p,\bullet}$ est un $n$-groupoïde de Kan pour tout $p\geq 0$, c'est-à-dire si $p\geq 0$, $q\geq 2$ et $0\leq k\leq q$ la fonction:
\begin{equation}
\xymatrix{\mathrm{Hom}\big(\Delta^p\boxtimes\Delta^q,X\big)  \ar[r] &
\mathrm{Hom}\big(\Delta^p\boxtimes\Lambda^{q,k},X\big)}
\end{equation}
est bijective pour $q\geq n+1$ et surjective pour $2\leq q\leq n$.
\item Si $2\leq p\leq n$, $2\leq q\leq n$ et $0\leq k\leq q$ la fonction:
\begin{equation}
\xymatrix{\mathrm{Hom}\big(\partial\Delta^p\boxtimes\Delta^q,X\big)  \ar[r] &
\mathrm{Hom}\big(\partial\Delta^p\boxtimes\Lambda^{q,k},X\big)}
\end{equation}
est surjective.
\item Si $1\leq p\leq n$, $2\leq q\leq n$ et $0\leq k\leq q$ la fonction:
\begin{equation}
\xymatrix{\mathrm{Hom}\big(\Delta^p\boxtimes\Delta^q,X\big)  \ar[r] &
\mathrm{Hom}\big(\partial\Delta^p\boxtimes\Delta^{q},X\big) \underset{\mathrm{Hom}\big(\partial\Delta^p\boxtimes\Lambda^{q,k},X\big)}{\bigtimes}
\mathrm{Hom}\big(\Delta^p\boxtimes\Lambda^{q,k},X\big) }
\end{equation}
est surjective.
\end{enumerate}
alors $X$ est un objet fibrant de la catégorie de modèles $({\bf MS},{\bf W}^{diag}_{n}, {\bf mono},{\bf fib}_{n}^{diag})$ de la Proposition \ref{moduno}.
\end{lemme}

\section{Ensembles simpliciaux tronqués}\label{tron}

Si $n\geq 0$, notons ${\bf\Delta}_{\leq n}$ la sous-catégorie pleine de la catégorie des simplexes ${\bf\Delta}$ dont les objets sont les catégories $[k]$, où $0 \leq k\leq n$. On appelle les objets de la catégorie des préfaisceaux $\simp_{\leq n}$, des \emph{ensembles simpliciaux $n$-tronqués}.

Considérons l'adjonction:
\begin{equation}\label{laa}
\ens^{\text{\scriptsize{$\Big({\bf\Delta}_{\leq n}\Big)^{op}$}}}\;=\,\vcenter{
\xymatrix@C+15pt{
\simp_{\leq n}\;\phantom{A}\ar@{}[r]|-{\perp\phantom{A}}
\ar@{<-}@/^18pt/[r]^-{\tau_{n}^{\,*}}
\ar@/_18pt/[r]_-{\tau_{n\;*}}& \phantom{A}\simp}}
\,=\;\ens^{\text{\scriptsize{${\bf\Delta}^{op}$}}}\,,
\end{equation}
induit par le foncteur d'inclusion canonique:
\begin{equation}\label{nun}
\xymatrix@C+10pt{{\bf\Delta}_{\leq n}\;\ar@{^(->}[r]^-{\tau_{n}}&{\bf\Delta}}.
\end{equation}

Si $X$ est un ensemble simplicial, on pose ${\bf csq}_{n}(X)$ pour noter l'ensemble simplicial $\tau_{n\;*}\tau_{n}^{*} \, (X)$, et on 
l'appelle le \emph{$n$-cosquelette de $X$}. Si le morphisme $\xymatrix@C-5pt{X\ar[r]^-{\eta_X}&{\bf csq}_{n}(X)}$ est un isomorphisme où $\eta$ est une unité quelconque de l'adjonction $\tau_{n\;*}  \dashv \tau_{n}^{*}$, on dit que $X$ est $n$-\emph{cosque\-lettique} ou qu'il est un $n$-\emph{cosquelette}. De fa\c con équivalente, $X$ est $n$-cosquelettique si pour tout ensemble simplicial $A$ la fonction canonique:
$$
\xymatrix@C+10pt{
\mathrm{Hom}_{\simp}\big(A,X\big) \ar[r]^-{\tau_{n}^{\;*}} &\mathrm{Hom}_{\simp_{\leq n}}\big(\tau_{n}^{\;*}A,\tau_{n}^{\;*}X\big)},
$$
est bijective. 

Plus encore, vu que le foncteur $\tau_{n\;*}$ de l'adjonction \eqref{laa} est pleinement fidèle, un ensemble simplicial $X$ est $n$-cosquelettique si et seulement si, $X$ est isomorphe à un objet dans l'image du foncteur $\tau_{n\;*}$.  

Montrons:

\begin{lemme}\label{cosque}
Si $X$ es un ensemble simplicial, les énoncés suivants sont équivalents:
\begin{enumerate}
\item $X$ est $n$-cosquelettique.
\item $X$ est $m$-cosquelettique pour tout $m\geq n$. 
\item Pour tout $m\geq n$, la fonction:
$$
\xymatrix@C+10pt{
\mathrm{Hom}_{\simp}\big(\Delta^{m+1},X\big) \ar[r]^{\alpha^{m}_{X}} &\mathrm{Hom}_{\simp}\big(\partial\Delta^{m+1},X\big)}
$$
induite du morphisme $\xymatrix@-10pt{\partial\Delta^{m+1} \ar[r]^-{\alpha^m} & \Delta^{m+1}}$ est bijective.
\end{enumerate}
\end{lemme}
\begin{proof}
Remarquons que si $\xymatrix@C-4pt{{\bf\Delta}_{\leq m}\;\ar@{^(->}[r]^-{j_{m}}&{\bf\Delta}_{\leq m+1}}$ note le foncteur d'inclusion canonique, on peut décomposer à isomorphisme près l'adjonction:
$$
\xymatrix@C+15pt{
\simp_{\leq m}\;\phantom{A}\ar@{}[r]|-{\perp\phantom{A}}
\ar@{<-}@/^18pt/[r]^-{\tau_{m}^{\phantom{A} *}}
\ar@/_18pt/[r]_-{\tau_{m\;*}}& \phantom{A}\simp}
$$
comme le composé suivant:
$$
\simp_{\leq m}
\xymatrix@C+18pt{
\phantom{a}\ar@{}[r]|-{\perp}
\ar@{<-}@/^18pt/[r]^-{j_{m}^{\phantom{a} *}}
\ar@/_18pt/[r]_-{j_{m\;*}}&
\phantom{a}}\simp_{\leq m+1}
\xymatrix@C+18pt{
\phantom{a}\ar@{}[r]|-{\perp}
\ar@{<-}@/^18pt/[r]^-{\tau_{m+1}^{\phantom{a} *}}
\ar@/_18pt/[r]_-{\tau_{m+1\;*}}
& \phantom{a}}\simp \,;
$$   
où les foncteurs $j_{m\;*}$ et $\tau_{m+1\;*}$ sont pleinement fidèles

On découle que si $X$ est un ensemble simplicial isomorphe à un objet dans l'image du foncteur $\tau_{m\;*}$, alors $X$ est aussi isomorphe à un objet dans l'image de $\tau_{m+1\;*}$. Donc, (i) implique (ii).

Pour montrer le reste de l'affirmation, observons par ailleurs que si $A$ est un ensemble simplicial $m$-tronqué, on peut décrire à isomorphisme près l'ensemble simplicial $(m+1)$-tronqué $j_{m\;*} (A)$ de la fa\c con suivante:
\begin{align*}
 j_{m\;*} (A)_{k} \, & = A_{k}\qquad\text{si}\;\;0\leq k\leq m \\
\text{et}\qquad j_{m\;*} (A)_{m+1} \, & =\mathrm{Hom}_{\simp_{\leq m}}\big(\tau_{m}^{\phantom{a} *}\partial\Delta^{m+1},A\big).
\end{align*}

Les morphismes faces et dégénérescences $\xymatrix{A_{k+1}\ar[r]^-{d^k_{i}} &A_{k}}$ et $\xymatrix{A_{k}\ar[r]^-{s^k_{i}} & A_{k+1}}$, sont pour $0\leq k< m$ ceux de $A$; tandis que le morphismes:
$$
\xymatrix@+10pt{
j_{m\;*} (A)_{m+1} \, = \, \mathrm{Hom}_{\simp_{\leq m}}\big(\tau_{m}^{\phantom{a} *}\partial\Delta^{m+1},A\big)
\ar@<+5pt>[r]^-{d^m_{i}} &\ar@<+5pt>[l]^-{s^m_{i}}
 \mathrm{Hom}_{\simp_{\leq m}}\big(\tau_{m}^{\phantom{a} *}\Delta^{m},A\big) \,\cong \,A_{m}},
$$ 
sont induisent des morphismes canoniques:
$$
\xymatrix@C-3pt{
\Delta^m \ar[rr]^-{\text{$i$-ième}}_(.38){\text{composante}}&& \underset{0\leq i\leq m+1}{\bigsqcup}\Delta^{m}\ar@{->>}[r]^{\sqcup \delta_{i}}& \partial\Delta^{m+1}}
\qquad\text{et}\qquad
\xymatrix@C-6pt{
\partial\Delta^{m+1} \ar[r]^-{\alpha^m}& \Delta^{m+1}\ar[r]^-{\sigma^m_{i}}& \Delta^{m}}.
$$

Autrement, vu qu'on a un isomorphisme:
\begin{equation}\label{trono}
j_{m\;*} (A)_{m+1}\; \cong \; \Bigg\{\text{\scriptsize{$\big(a_{0},\dots,a_{m+1}\big) \;\in\;\underset{0}{\overset{m+1}{\prod}} A_{m}$}}\;\Bigg|\;
\vcenter{\xymatrix@R=1pt{\text{\scriptsize{$d^{m-1}_{i}a_{j}=d^{m-1}_{j-1}a_{i}$}}\\ \text{\scriptsize{si $0\leq i<j\leq m+1$.}}}}\Bigg\},
\end{equation}
ils sont donnés par les règles:
$$
\xymatrix@C+15pt@R=1pt{
\quad (a_{0},\dots,a_{m+1})\quad  \ar@{|->}[r] & \quad a_{i}\quad\\
j_{m\;*} (A)_{m+1}   \ar@<+5pt>[r]|-{d^m_{i}} & \ar@<+5pt>[l]|-{s^m_{i}} A_{m}\\
\quad (b_{0},\dots,b_{m+1})\quad &\ar@{|->}[l] \quad a \quad }
$$
$$
\text{où}\qquad b_{k}\; = \,
\begin{cases}
 s_{i-1}^{m-1}d_{k}^{m-1}a  & 0\leq k\leq i-1\\
 a                                         &i\leq k\leq i+1 \\
 s_{i}^{m-1}d_{k-1}^{m-1}a     & i+2\leq k \leq m+1\,.
\end{cases}
$$

On peut maintenant décrire une unité de l'adjonction $j_{m\;*}  \dashv j_{m}^{\phantom{a} *}$ comme suit: Si $B$ est un ensemble simplicial $(m+1)$-tronqué, le morphisme $\xymatrix@C-5pt{B\ar[r]& j_{m\;*}  j_{m}^{\phantom{a} *}(B),}$ est défini pour $0\leq k\leq m$ comme la fonction identité, et pour $k=m+1$ comme la fonction:
$$
\xymatrix@C+10pt@R=12pt{
B_{m+1}\ar[r]\ar@{}[d]|-{\mathrel{\reflectbox{\rotatebox[origin=c]{90}{$\cong$}}}}& 
j_{m\;*}  j_{m}^{\phantom{a} *}(B)_{m+1}\ar@{}[d]|-{\mathrel{\reflectbox{\rotatebox[origin=c]{90}{$\cong$}}}}\\
\mathrm{Hom}_{\simp_{\leq m+1}}\big(\tau_{m+1}^{\phantom{a} *} \Delta^{m+1},B\big) \ar@/_10pt/[rd] & \mathrm{Hom}_{\simp_{\leq m}}\big(\tau_{m}^{\phantom{a} *} \partial\Delta^{m+1},j_{m}^{\phantom{a} *}B\big)\ar@{}[d]|-{\mathrel{\reflectbox{\rotatebox[origin=c]{90}{$\cong$}}}}\\
&\mathrm{Hom}_{\simp_{\leq m+1}}\big(\tau_{m+1}^{\phantom{a} *} \partial\Delta^{m+1},B\big)\,,}
$$
induite du morphisme $\xymatrix@C+10pt{ \partial\Delta^{m+1}\ar[r]^-{\alpha^m}&\Delta^{m+1}}$.

En particulier, vu que $j_{m\;*}$ est un foncteur pleinement fidèle, on vérifie que si $X$ est un ensemble simplicial quelconque, l'ensemble simplicial $(m+1)$-tronqué $\tau_{m+1}^{\phantom{a} *} (X)$ est isomorphe à un objet dans l'image du foncteur $j_{m\;*}$, si et seulement si la fonction:
$$
\xymatrix@C+10pt@R=12pt{
\mathrm{Hom}_{\simp_{\leq m+1}}\big(\tau_{m+1}^{\phantom{a} *} \Delta^{m+1},\tau_{m+1}^{\phantom{a} *}X\big) \ar[r]
\ar@{}[d]|-{\mathrel{\reflectbox{\rotatebox[origin=c]{90}{$\cong$}}}}&
\mathrm{Hom}_{\simp_{\leq m+1}}\big(\tau_{m+1}^{\phantom{a} *} \partial\Delta^{m+1},\tau_{m+1}^{\phantom{a} *}X\big)
\ar@{}[d]|-{\mathrel{\reflectbox{\rotatebox[origin=c]{90}{$\cong$}}}}\\
\mathrm{Hom}_{\simp}\big(\Delta^{m+1},X\big) \ar[r]_{\alpha^{m}_{X}} &\mathrm{Hom}_{\simp}\big(\partial\Delta^{m+1},X\big)},
$$
est bijective. Donc, les propriétés (ii) et (iii) sont équivalentes.
\end{proof}

Montrons aussi:

\begin{corollaire}\label{lecoro}
Soit $n\geq 0$. Si $X$ est un ensemble simplicial qui satisfait la condition d'extension de Kan en dimension $1\leq m\leq n+1$, c'est-à-dire si la fonction:
$$
\xymatrix{\mathrm{Hom}_{\simp}\big(\Delta^{m+1},X\big)\ar[r]^{\alpha^{m,k}_{X}}&\mathrm{Hom}_{\simp}\big(\Lambda^{m+1,k},X\big)}
$$
est surjective pour $1\leq m\leq n+1$ et $0\leq k\leq m+1$, alors ${\bf csq}_{n+1}(X)$ est un complexe de Kan tel que: 
$$
\text{$\pi_{m}\big({\bf csq}_{n+1}X,a\big) \, = \, 0$ \quad \text{pour \; tout}\quad $a\in {\bf csq}_{n+1}(X)_{0}=X_{0}$ \quad\text{et} \quad $m\geq n+1$;}
$$
autrement dit ${\bf csq}_{n+1}(X)$ est un ensemble simplicial $n$-fibrant \emph{i.e.} est un objet fibrant de la catégorie de modèles $(\simp,{\bf W}_n, {\bf mono},{\bf fib}_n)$ du Théorème \ref{ntypess}.

En particulier un ensemble simpliciaux $(n+1)$-cosquelettique $X$ est un complexe de Kan si et seulement si, il est un ensemble simplicial $n$-fibrant si et seulement si, il satisfait la condition d'extension de Kan en dimension $1\leq m\leq n+1$.

D'un autre si on suppose que $X$ soit un complexe de Kan et $\eta$ est une unité quelconque de l'adjonction $\tau_{n\;*}  \dashv \tau_{n}^{\phantom{a} *}$, alors le morphisme $\xymatrix@C-5pt{X\ar[r]^-{\eta_X}&{\bf csq}_{n+1}(X)}$ est une $n$-équivalence faible. 
\end{corollaire}
\begin{proof}
On déduit des Lemmes \ref{glennlemme} et \ref{cosque} que pour tout $X$ l'ensemble simplicial ${\bf csq}_{n+1}(X)$ satisfait la condition d'extension de Kan en dimension $m\geq n+2$.

Considérons maintenant le carré commutatif:
\begin{equation}\label{lecai}
\xymatrix@C+25pt{
\mathrm{Hom}_{\simp}\big(\Delta^{m+1},X\big)\ar[d]\ar[r]^{\alpha^{m,k}_{X}}&\mathrm{Hom}_{\simp}\big(\Lambda^{m+1,k},X\big)\ar[d]\\
\mathrm{Hom}_{\simp}\big(\Delta^{m+1},{\bf csq}_{n+1}X\big)\ar[r]_{\alpha^{m,k}_{{\bf csq}_{n+1}(X)}} &\mathrm{Hom}_{\simp}\big(\Lambda^{m+1,k},{\bf csq}_{n+1}X\big)\,,}
\end{equation}
induit des morphismes  $\xymatrix@C-5pt{X\ar[r]^-{\eta_X}& {\bf csq}_{n+1}X}$ et $\xymatrix@C-8pt{ \Lambda^{m+1,k}\ar[r]&\Delta^{m+1}}$.

Vu qu'on a un isomorphisme $(\eta_X)_p:X_p\cong {\bf csq}_{n+1}(X)_{p}$ pour $0\leq p\leq n+1$, on trouve que dans le carré \eqref{lecai} la flèche qui descende à droite est une bijection, pour tout $1\leq m\leq n+1$ et $0\leq k\leq m+1$  $\big(${de façon équivalente on peut remarquer que ${\bf sq}_{n+1}(\Lambda^{m+1,k}) \cong \Lambda^{m+1,k}$ si $1\leq m\leq n+1$}$\big)$; en particulier si on suppose que $\alpha_{X}^{m,k}$ soit surjective pour $1\leq m\leq n+1$ et $0\leq k\leq m+1$ alors la flèche $\alpha^{m,k}_{{\bf csq}_{n+1}(X)}$ est aussi surjective pour $1\leq m\leq n+1$ et $0\leq k\leq m+1$; c'est-à-dire, si $X$ satisfait la condition d'extension de Kan en dimension $1\leq m\leq n+1$ alors ${\bf csq}_{n+1}(X)$ aussi, donc ${\bf csq}_{n+1}(X)$ est un complexe de Kan.

Pour calculer maintenant les groupes d'homotopie de ${\bf csq}_{n+1}X$, remarquons par ailleurs que si $Y$ est un complexe de Kan quelconque dont la fonction:
$$
\xymatrix{\mathrm{Hom}_{\simp}\big(\Delta^{m+1},Y\big)\ar[r]^{\alpha^{m}_{Y}}&\mathrm{Hom}_{\simp}\big(\partial\Delta^{m+1},Y\big)}
$$ 
est injective pour $m\geq n+1$, il résulte de la Proposition \ref{kankan} que pour tout $a\in Y_{0}$ et $m\geq n+2$, le groupe $\pi_{m}(Y,a)$ est un quotient de l'ensemble ponctuel; donc, $\pi_{m}(Y,a)=0$ pour tout $a\in Y_{0}$ et $m\geq n+2$. En particulier, $\pi_{m}\big({\bf csq}_{n+1}X,a\big)=0$ pour tout $a\in X_{0}$ et $m\geq n+2$.

On déduit aussi de la Proposition \ref{kankan} que $\pi_{n+1}\big({\bf csq}_{n+1}X,a\big)=0$ pour tout $a\in X_{0}$. En effet, vu qu'on peut identifier:
$$
{\bf csq}_{n+1}(X)_{n+2}\;=\;
 \Bigg\{\text{\scriptsize{$\big(z_{0},\dots,z_{n+2}\big) \;\in\;\underset{0}{\overset{n+2}{\prod}} X_{n+1}$}}\;\Bigg|\;
\vcenter{\xymatrix@R=1pt{\text{\scriptsize{$d_{i}z_{j}=d_{j-1}z_{i}$}}\\ \text{\scriptsize{si $0\leq i<j\leq n+2$}}}}\Bigg\}\qquad\text{et}\qquad 
{\bf csq}_{n+1}(X)_{n+1}\;=\; X_{n+1}\,,
$$
de sort que les morphismes:
$$
\xymatrix@C+2pt{{\bf csq}_{n+1}(X)_{n+2}\ar[r]^-{d_{i}}&{\bf csq}_{n+1}(X)_{n+1}} \quad \text{pour}\quad 0\leq i \leq n+2,
$$
correspondent aux projections canoniques; étant donnés $x$ et $y$ de $(n+1)$-simplexes de ${\bf csq}_{n+1}(X)$ vérifiant que $d_{i}x=d_{i}y=a$ pour tout $0\leq i\leq n+1$, il est simple à voir que la donnée $(x,y,a,\dots,a)$ est une homotopie de $x$ vers $y$. 

Enfin, si $X$ est un complexe de Kan on peut faire appelle à la Proposition \ref{kankan} pour montrer que le morphisme $\xymatrix@C-5pt{X\ar[r]^-{\eta_X}&{\bf csq}_{n+1}(X)}$ est une $n$-équivalence faible parce que $(\eta_X)_p:X_p\cong {\bf csq}_{n+1}(X)_{p}$ est un isomorphisme pour $0\leq p\leq n+1$.
\end{proof}

Remarquons que si $X$ et $Y$ sont des ensembles simpliciaux $(n+1)$-cosquelettiques (resp. des complexes de Kan $(n+1)$-cosquelettiques), on déduit du Lemme \ref{cosque} (resp. les Lemmes \ref{cosque} et \ref{nfibcartesiene}) que le produit cartésien argument par argument $X\times Y$ est aussi un ensemble simplicial $(n+1)$-cosquelettiques (resp. un complexe de Kan $(n+1)$-cosquelettique). Montrons:

\begin{lemme}\label{lemehomo}
Soit $n\geq 0$. Si $X$ et $Y$ sont des ensembles simpliciaux alors $\underline{\mathrm{Hom}}_{\simp}\big(X,{\bf csq}_{n+1}Y\big)$ est un ensemble simplicial $(n+1)$-cosquelettique et le morphisme $\xymatrix@C-5pt{X\ar[r]^-{\eta_X}&{\bf csq}_{n+1}(X)}$ induit un isomorphisme d'ensembles simpliciaux:
\begin{equation}\label{lemcos1}
\xymatrix@C-8pt{
\underline{\mathrm{Hom}}_{\simp}\big({\bf csq}_{n+1}X,{\bf csq}_{n+1}Y\big) \ar[r]^-{\cong} &
\underline{\mathrm{Hom}}_{\simp}\big(X,{\bf csq}_{n+1}Y\big)}\,;
\end{equation}
où $\eta$ est une unité quelconque de l'adjonction $\tau_{n\;*}  \dashv \tau_{n}^{\phantom{a} *}$.

En particulier l'ensemble simplicial des morphismes $\underline{\mathrm{Hom}}(X,Y)$ est $(n+1)$-cosquelettiques (resp. un complexe de Kan $(n+1)$-cosquelettique) si $X$ et $Y$ sont des ensembles simpliciaux $(n+1)$-cosquelettiques (resp. des complexes de Kan $(n+1)$-cosquelettiques). Autrement dit, la sous-catégorie pleine de $\simp$ dont les objets sont les ensembles simpliciaux $(n+1)$-cosquelettique (resp. les complexes de Kan $(n+1)$-cosquelettiques) est une sous-catégorie cartésienne fermée.
\end{lemme}
\begin{proof}
Soient $X$ et  $Y$ des ensembles simpliciaux. Montrons que l'ensemble simplicial $\underline{\mathrm{Hom}}_{\simp}\big(X,{\bf csq}_{n+1}Y\big)$ est $(n+1)$-cosquelettique, \emph{i.e.} que le morphisme:
\begin{equation}\label{lemcos}
\xymatrix@C-8pt{
\mathrm{Hom}_{\simp}\Big(\Delta^{m+1},\underline{\mathrm{Hom}}_{\simp}\big(X,{\bf csq}_{n+1}Y\big)\Big)\ar[r]^{\alpha^{m}_{X}} &
\mathrm{Hom}_{\simp}\Big(\partial\Delta^{m+1},\underline{\mathrm{Hom}}_{\simp}\big(X,{\bf csq}_{n+1}Y\big)\Big)\,,}
\end{equation}
est une bijection si $m\geq n+1$.

Pour cela, on décompose \eqref{lemcos} comme la suite des bijectives suivantes:
\def\objectstyle{\scriptstyle}
\def\labelstyle{\scriptstyle}
\begin{align*}
\mathrm{Hom}_{\simp}\Big( \Delta^{m+1},\underline{\mathrm{Hom}}_{\simp}\big(X,{\bf csq}_{n+1}Y\big)\Big) \,& \cong \,
\mathrm{Hom}_{\simp}\big( \Delta^{m+1}\times X,{\bf csq}_{n+1}Y) &\text{(i)} \\
\,& \cong \, 
\mathrm{Hom}_{\simp_{\leq n+1}}\Big(\tau_{n+1}^*\big(\Delta^{m+1}\times X\big),\tau_{n+1}^*Y\Big) &\text{(ii)}\\
\,& \cong \, 
\mathrm{Hom}_{\simp_{\leq n+1}}\Big(\tau_{n+1}^*\big(\Delta^{m+1}\big)\times \tau_{n+1}^*\big(X\big),\tau_{n+1}^*Y\Big) &\text{(iii)}\\ 
\,& \cong \, 
\mathrm{Hom}_{\simp_{\leq n+1}}\Big(\tau_{n+1}^*\big(\partial\Delta^{m+1}\big)\times \tau_{n+1}^*\big(X\big),\tau_{n+1}^*Y\Big) &\text{(iv)}\\
\,& \cong \, 
\mathrm{Hom}_{\simp_{\leq n+1}}\Big(\tau_{n+1}^*\big(\partial\Delta^{m+1}\times X\big),\tau_{n+1}^*Y\Big) &\text{(v)}\\
\,& \cong \, 
\mathrm{Hom}_{\simp}\big(\partial\Delta^{m+1}\times X,{\bf csq}_{n+1}Y\big) &\text{(vi)}\\
\,& \cong \, 
\mathrm{Hom}_{\simp}\Big(\partial\Delta^{m+1},\underline{\mathrm{Hom}}_{\simp}\big(X,{\bf csq}_{n+1}Y\big)\Big)&\text{(vii)}
\end{align*}
où (i) et (vii) sont conséquence de l'adjonction \eqref{deltacarte}; (ii) et (vi) sont déduites de l'égalité ${\bf csq}_{n+1}=\tau_{n+1\;*}\tau_{n+1}^{\phantom{A}*}$ et parce que $\tau_{n+1\;*} \dashv \tau_{n+1}^{\phantom{A}*}$; (iii) et (v) sont des bijections, vu que $\tau_{n+1}^{\phantom{A}*}$ commute aux produits; et finalement on a (iv), car la $(n+1)$-troncation du morphisme $\xymatrix@C-9pt{\partial\Delta^{m+1}\ar@{^(->}[r]&\Delta^{m+1}}$ est un isomorphisme toujours que $m\geq n+1$. 

Montrons maintenant que le morphisme \eqref{lemcos1} est une bijection pour tout ensemble simplicial $Y$. Pour cela rappelons que le foncteur $\tau_{n+1}^{\phantom{A}*}$ commute aux produits, et que $\xymatrix{\tau_{n+1}^{\phantom{A}*}(X)\ar[r]^-{\tau_{n+1}^{\phantom{A}*}(\eta_X)}&\tau_{n+1}^{\phantom{A}*}\big({\bf csq}_{n+1}(X)\big)}$ est un isomorphisme d'ensembles simpliciaux tronqués. On vérifie alors que \eqref{lemcos1} se décompose comme la suite de bijections:
\begin{align*}
\underline{\mathrm{Hom}}_{\simp}\big(X,{\bf csq}_{n+1}Y\big)_{m} \,& = \,
\mathrm{Hom}_{\simp}\big( X\times\Delta^{m},{\bf csq}_{n+1}Y) \\
\,& \cong \, \mathrm{Hom}_{\simp_{\leq n+1}}\Big(\tau_{n+1}^{\phantom{A}*}\big(X\times\Delta^{m}\big),\tau_{n+1}^{\phantom{A}*}Y\Big) \\
\,& \cong \, \mathrm{Hom}_{\simp_{\leq n+1}}\big(\tau_{n+1}^{\phantom{A}*}(X)\times\tau_{n+1}^{\phantom{A}*}(\Delta^{m}),\tau_{n+1}^{\phantom{A}*}Y\Big) \\
\,& \cong \, \mathrm{Hom}_{\simp_{\leq n+1}}\Big(\tau_{n+1}^{\phantom{A}*} \big({\bf csq}_{n+1}(X)\big)\times\tau_{n+1}^{\phantom{A}*}(\Delta^{m}),\tau_{n+1}^{\phantom{A}*}Y\Big) \\
\,& \cong \, \mathrm{Hom}_{\simp_{\leq n+1}}\Big(\tau_{n+1}^{\phantom{A}*} \big({\bf csq}_{n+1}(X)\times\Delta^{m}\big),\tau_{n+1}^{\phantom{A}*}Y\Big) \\
\,& \cong \, \mathrm{Hom}_{\simp_{\leq n+1}}\Big({\bf csq}_{n+1}\big(X\big)\times\Delta^{m},{\bf csq}_{n+1}Y\Big) \\
\,& = \, \underline{\mathrm{Hom}}_{\simp}\big({\bf csq}_{n+1}X,{\bf csq}_{n+1}Y\big)_{m}\,.
\end{align*}
\end{proof}

Montrons une version réduit du Lemme \ref{lemehomo}:

\begin{lemme}\label{lemehomored}
Soient $X$ et $W$ des ensembles simpliciaux réduits. Si $W$ est $(n+1)$-cosquelettique alors $\underline{\mathrm{Hom}}_{\simp_0}\big(X,W\big)$ l'est aussi. 
\end{lemme}
\begin{proof}
Montrons que si $m\geq n+1$ la fonction:
\begin{equation}\label{lemcosred}
\xymatrix@C-8pt{
\mathrm{Hom}_{\simp}\Big(\Delta^{m+1},\underline{\mathrm{Hom}}_{\simp_0}\big(X,W\big)\Big)\ar[r]^{\alpha^{m}_{X}} &
\mathrm{Hom}_{\simp}\Big(\partial\Delta^{m+1},\underline{\mathrm{Hom}}_{\simp_0}\big(X,W\big)\Big)\,,}
\end{equation}
est bijective.

En effet, pour montrer cela on procède de façon analogue à la preuve du Lemme \ref{lemehomo} en décomposant la fonction \eqref{lemcosred} comme la suite des bijectives suivantes:
\begin{align*} 
\def\objectstyle{\scriptstyle}
\def\labelstyle{\scriptstyle}
\mathrm{Hom}_{\simp}\Big( \Delta^{m+1},\underline{\mathrm{Hom}}_{\simp_0}\big(X,W\big)\Big) \,& \cong \,
\mathrm{Hom}_{\simp}\Big( \Delta^{m+1},\underline{\mathrm{Hom}}_{\simp_0}\big(X,{\bf csq}_{n+1}W\big)\Big)&\text{(i)} \\
\,& \cong \,
\mathrm{Hom}_{\simp_0}\Big(  (X\times  \Delta^{m+1}) \big/ (\star\times  \Delta^{m+1}),{\bf csq}_{n+1}W  \Big) &\text{(ii)} \\
\,& \cong \,
\mathrm{Hom}_{\simp}\Big(  \big(X\times  \Delta^{m+1}\big) \big/ \big(\star\times  \Delta^{m+1}\big),{\bf csq}_{n+1}W  \Big) &\text{(iii)} \\
\,& \cong \,
\mathrm{Hom}_{\simp_{\leq n+1}}\bigg(  \tau_{n+1}^* \Big(\big(X\times  \Delta^{m+1}\big) \big/ \big(\star\times  \Delta^{m+1}\big)\Big), \tau_{n+1}^* W  \bigg) &\text{(iv)} \\
\,& \cong \, 
\mathrm{Hom}_{\simp_{\leq n+1}}\bigg(   \Big(\tau_{n+1}^* X\times \tau_{n+1}^* \Delta^{m+1}\Big) \big/ \Big(\star\times  \tau_{n+1}^* \Delta^{m+1}\Big), \tau_{n+1}^* W  \bigg) &\text{(v)} \\
\,& \cong \, 
\mathrm{Hom}_{\simp_{\leq n+1}}\bigg(   \Big(\tau_{n+1}^* X\times \tau_{n+1}^* \partial\Delta^{m+1}\Big) \big/ \Big(\star\times  \tau_{n+1}^* \partial\Delta^{m+1}\Big), \tau_{n+1}^* W  \bigg) &\text{(vi)} \\
\,& \cong \, 
\mathrm{Hom}_{\simp_{\leq n+1}}\bigg(  \tau_{n+1}^* \Big(\big(X\times  \partial\Delta^{m+1}\big) \big/ \big(\star\times  \partial\Delta^{m+1}\big)\Big), \tau_{n+1}^* W  \bigg) &\text{(vii)} \\
\,& \cong \, 
\mathrm{Hom}_{\simp}\Big(  \big(X\times  \partial\Delta^{m+1}\big) \big/ \big(\star\times  \partial\Delta^{m+1}\big),{\bf csq}_{n+1}W  \Big) &\text{(viii)}\\
\,& \cong \,
\mathrm{Hom}_{\simp_0}\Big(  (X\times  \partial\Delta^{m+1}) \big/ (\star\times  \partial\Delta^{m+1}),{\bf csq}_{n+1}W  \Big) &\text{(ix)} \\
 \,& \cong \,
\mathrm{Hom}_{\simp}\Big(\partial \Delta^{m+1},\underline{\mathrm{Hom}}_{\simp_0}\big(X,{\bf csq}_{n+1}W\big)\Big)&\text{(x)} \\
 \,& \cong \,
\mathrm{Hom}_{\simp}\Big( \partial \Delta^{m+1},\underline{\mathrm{Hom}}_{\simp_0}\big(X,W\big)\Big)&\text{(xi)} 
\end{align*}
où (i) et (xi) se suivent de l'isomorphisme $W\cong {\bf csq}_{n+1}W$; (ii) et (x) sont conséquence de l'adjonction:
$$
\mathrm{Hom}_{\simp_0}\Big(  (X\times  K) \big/ (\star\times K),Y  \Big)  \, \cong \,
\mathrm{Hom}_{\simp}\Big(K,\underline{\mathrm{Hom}}_{\simp_0}\big(X,Y\big)\Big)\,; 
$$
(iii) et (ix) sont bijections parce que $\xymatrix@C-10pt{\simp_0\ar[r]&\simp}$ est pleinement fidèle; (iv) et (viii) sont déduites de l'égalité ${\bf csq}_{n+1}=\tau_{n+1\;*}\tau_{n+1}^{\phantom{A}*}$ et parce que $\tau_{n+1\;*} \dashv \tau_{n+1}^{\phantom{A}*}$; (v) et (vii) sont des bijections vu que $\tau_{n+1}^{\phantom{A}*}$ commute aux petites limites et colimites; et finalement on a (vi) car la $(n+1)$-troncation du morphisme $\xymatrix@C-9pt{\partial\Delta^{m+1}\ar@{^(->}[r]&\Delta^{m+1}}$ est un isomorphisme toujours que $m\geq n+1$. 
\end{proof}

\renewcommand{\thesubsection}{\S\thesection.\arabic{subsection}}
\subsection{}\; \label{faiblecos}
\renewcommand{\thesubsection}{\thesection.\arabic{subsection}}

Si $n\geq 0$, un ensemble simplicial $X$ est appelé \emph{faiblement $n$-cosquelettique} si la fonction:
$$
\xymatrix@C+10pt{\mathrm{Hom}_{\simp}\big(\Delta^{m+1},X\big) \ar[r]^{\alpha^{m}_{X}} &\mathrm{Hom}_{\simp}\big(\partial\Delta^{m+1},X
\big)},
$$ 
induite de l'inclusion $\xymatrix@C-3pt{\partial\Delta^{m+1}\,\ar@{^(->}[r]^-{\alpha^{m}}&\Delta^{m+1}}$, est bijective pour $m>n$ et injective pour $m=n$. 

D'après le Lemme \ref{cosque}, on a que $X$ est faiblement $n$-cosquelettique si et seulement si, le morphisme canonique d'ensembles simpliciaux $\xymatrix@C-8pt{X\ar[r]&{\bf csq}_{m}(X)}$ est un monomorphisme pour $m=n$ et un isomorphisme pour $m\geq n+1$; \emph{i.e.} si $X$ est contenu dans son $n$-cosquelette et $X$ est $(n+1)$-cosquelettique\footnote{Voir la Définition 2.5 de \cite{duskin}, où les ensembles simpliciaux faiblement $n$-cosquelettiques sont dits complexes de Postnikov de dimension $n$}. 

Étant donné un ensemble simplicial $X$ et un nombre $n\geq 0$, on définit un ensemble simplicial faiblement $n$-cos\-que\-le\-ttique ${\bf csq}'_{n+1}(X)$ et une suite de morphismes:
$$
\xymatrix{X\ar[r] & {\bf csq}_{n+1}(X) \ar[r]& {\bf csq}'_{n+1}(X),}
$$
de la fa\c con suivante\footnote{L'ensemble simplicial qu'on note ${\bf csq}'_{n+1}(K)$ est désigné par $K^{(n)}$ dans la Définition 8.1 de \cite{may}.}: On considère la relation d'équivalence $\overset{n}{\sim}$ dans l'ensemble $X_{n+1}$ des $(n+1)$-simplexes de $X$:
\begin{equation}\label{relationn}
x\overset{n}{\sim} y \in X_{n+1}\quad\iff\quad d_{i}x=d_{i}y\quad\text{pour}\quad 0\leq i\leq n+1;
\end{equation}
et on note $\tau_{n+1}^{*}(X)'$ l'ensemble simplicial $(n+1)$-tronqué, qu'on déduit de la $(n+1)$-troncation $\tau_{n+1}^{*}(X)$ de $X$ en imposant la relation $\overset{n}{\sim}$ sur l'ensemble extrême $X_{n+1}$.

Si $0\leq m\leq n-1$, les morphismes faces et dégénérescences de $\tau_{n+1}^{*}(X)'$ sont ceux de $X$: 
$$
\vcenter{\xymatrix@C+12pt{ X_{m+1} \ar@{=}[r]\ar[d]_{d_{i}} & \tau_{n+1}^{*}(X)'_{m+1} \ar@{-->}[d]^{d_{i}}\\
X_{m} \ar@{=}[r] & \tau_{n+1}^{*}(X)'_{m}}}\qquad\text{et}\qquad
\vcenter{\xymatrix@C+12pt{ X_{m+1} \ar@{=}[r]& \tau_{n+1}^{*}(X)'_{m+1} \\
X_{m}\ar[u]^{s_{i}}  \ar@{=}[r] & \tau_{n+1}^{*}(X)'_{m} \ar@{-->}[u]_{s_{i}}}}\,;
$$  
tandis que pour $m=n$, ils sont définis par des carrés commutatifs:
$$
\vcenter{\xymatrix@C+12pt{ X_{n+1} \ar[r]^-{\text{quotient}}\ar[d]_{d_{i}} & \tau_{n+1}^{*}(X)'_{n+1} \ar@{-->}[d]^{d_{i}}\\
X_{n} \ar@{=}[r] & \tau_{n+1}^{*}(X)'_{n}}}\qquad\text{et}\qquad
\vcenter{\xymatrix@C+12pt{ X_{n+1} \ar[r]^-{\text{quotient}}& \tau_{n+1}^{*}(X)'_{n+1} \\
X_{n}\ar[u]^{s_{i}}  \ar@{=}[r] & \tau_{n+1}^{*}(X)'_{n} \ar@{-->}[u]_{s_{i}}}}\,.
$$

On définit ${\bf csq}'_{n+1}(X)$ comme l'ensemble simplicial $\tau_{n+1\,*}\big(\tau_{n+1}^{*}(X)' \big)$, où $\tau_{n+1\,*}$ est l'adjoint à droite du foncteur de troncation $\tau_{n+1}^*$ $\big($voir \eqref{laa}$\big)$. En appliquant $\tau_{n+1\,*}$ au morphisme quotient $\xymatrix@C-10pt{\tau_{n+1}^*(X)\ar[r]& \tau_{n+1}^{*}(X)'}$, on obtient un morphisme d'ensembles simpliciaux:
\begin{equation}\label{tildeco}
\xymatrix{{\bf csq}_{n+1}(X) \ar[r]& {\bf csq}'_{n+1}(X).}
\end{equation}


\begin{lemme}\label{dusl}
Si $X$ est un ensemble simplicial quelconque alors ${\bf csq}'_{n+1}(X)$ est un ensemble simplicial faiblement $n$-cosquelettique, et si $X$ est un ensemble simplicial faiblement $n$-cosquelettique alors le morphisme composé:
\begin{equation}\label{tildecoco}
\xymatrix{X\ar[r]&{\bf csq}_{n+1}(X) \ar[r]& {\bf csq}'_{n+1}(X),}
\end{equation}
est un isomorphismes d'ensembles simpliciaux. Autrement dit, un ensemble simplicial $X$ est faiblement $n$-cosquelettique si et seulement si le morphisme composé \eqref{tildecoco} est un isomorphismes d'ensembles simpliciaux.  

D'un autre, si $X$ satisfait la condition d'extension de Kan en dimension $1\leq m\leq n+1$, alors l'ensemble simplicial ${\bf csq}'_{n+1}(X)$ est un complexe de Kan et le morphisme \eqref{tildeco} ci-dessus est une $\infty$-équivalence faible d'ensembles simpliciaux. 
\end{lemme}
\begin{proof}
Remarquons que si $X$ est un ensemble simplicial quelconque, dans le carré commutatif induit du morphisme composé \eqref{tildecoco}:
\begin{equation}\label{lecarredu}
\xymatrix@C+25pt{
\mathrm{Hom}_{\simp}(\Delta^{n+1},X)\ar[d]_-{\gamma_{X}}\ar[r]^{\alpha^{n}_{X}}&\mathrm{Hom}_{\simp}(\partial\Delta^{n+1},X)\ar[d]\\
\mathrm{Hom}_{\simp}\big(\Delta^{n+1},{\bf csq}'_{n+1}(X)\big)\ar[r]_{\alpha^{n}_{{\bf csq}'_{n+1}X}} &\mathrm{Hom}_{\simp}\big(\partial\Delta^{n+1},{\bf csq}'_{n+1}(X)\big)\,,}
\end{equation}
la fonction $\gamma_{X}$ s'identifie au morphisme quotient $\xymatrix@C-10pt{X_{n+1}\ar[r] & X_{n+1}/\overset{n}{\sim}}$; et la fonction qui descende à droite est un isomorphisme (car les ensembles simpliciaux ${\bf csq}'_{n+1}(X)$ et $X$ ont la même $n$-troncation).  

D'un autre côté, il se suit de la définition de la fonction $\alpha_{X}^n$ que si $x$ et $y$ sont des éléments de l'ensemble $\mathrm{Hom}_{\simp}(\Delta^{n+1},X)$, alors $x\overset{n}{\sim}y$ si et seulement si $\alpha_{X}^n(x)=\alpha_{X}^n(y)$. Donc, $\gamma_{X}(x)=\gamma_{X}(y)$ si et seulement si $\alpha_{X}^n(x)=\alpha_{X}^n(y)$, et si $\alpha_{X}^n$ est injective alors $\gamma_{X}$ est un isomorphisme $\big($\emph{i.e.} {$\tau_{n+1}^{*}(X)' = \tau_{n+1}^{*}(X)$}$\big)$.

On en déduit en particulier que si $X$ est un ensemble simplicial quelconque, alors:
\begin{enumerate}
\item $\alpha^{n}_{{\bf csq}'_{n+1}X}$ est toujours injective.
\item Si la fonction $\alpha_{X}^n$ est injective, alors $\gamma_{X}$ est un isomorphisme, \emph{i.e.} $\tau_{n+1}^{*}(X)' = \tau_{n+1}^{*}(X)$.
\end{enumerate}

Vu que par définition ${\bf csq}'_{n+1}(X)$ est un ensemble simplicial $(n+1)$-cosquelettique, on en conclut:
\begin{enumerate} 
\item ${\bf csq}'_{n+1}(X)$ est toujours un ensemble simplicial faiblement $n$-cosquelettique.
\item Si $X$ est faiblement $n$-cosquelettique, alors \eqref{tildecoco} est un isomorphisme.
\end{enumerate}

Supposons maintenant que $X$ satisfait la condition d'extension de Kan en dimension $1\leq m\leq n+1$. D'après le Corollaire \ref{lecoro}, pour montrer que l'ensemble simplicial $(n+1)$-cosquelettique ${\bf csq}'_{n+1}(X)$ soit un complexe de Kan, il suffit de vérifier que ${\bf csq}'_{n+1}(X)$ satisfasse la condition d'extension de Kan en dimension $1\leq m\leq n+1$. 

Pour commencer, il se suit du fait que  ${\bf csq}'_{n+1}(X)$ et $X$ aient la même $n$-troncation, que ${\bf csq}'_{n+1}(X)$ satisfait la condition d'extension de Kan en dimension $1\leq m\leq n-1$. 

D'un autre côté, ${\bf csq}'_{n+1}(X)$ satisfait la condition d'extension de Kan en dimension $n$, vu que dans le carré commutatif:
$$
\xymatrix@C+25pt{
\mathrm{Hom}_{\simp}(\Delta^{n+1},X)\ar[d]\ar[r]^{\alpha^{n,k}_{X}}&\mathrm{Hom}_{\simp}(\Lambda^{n+1,k},X)\ar[d]\\
\mathrm{Hom}_{\simp}\big(\Delta^{n+1},{\bf csq}'_{n+1}(X)\big)\ar[r]_{\alpha^{n,k}_{{\bf csq}'_{n+1}X}} &\mathrm{Hom}_{\simp}\big(\Lambda^{n+1,k},{\bf csq}'_{n+1}(X)\big)\,,}
$$
la fonction $\alpha^{n,k}_{X}$ est surjective, et la fonction qui descende à droite est bijective. 

De la même fa\c con, pour montrer que ${\bf csq}'_{n+1}(X)$ satisfait la condition d'extension de Kan en dimension $n+1$, on considère le carré commutatif:
$$
\xymatrix@C+25pt{
\mathrm{Hom}_{\simp}(\Delta^{n+2},X)\ar[d]\ar[r]^{\alpha^{n+1,k}_{X}}&\mathrm{Hom}_{\simp}(\Lambda^{n+2,k},X)\ar[d]\\
\mathrm{Hom}_{\simp}\big(\Delta^{n+2},{\bf csq}'_{n+1}(X)\big)\ar[r]_{\alpha^{n+1,k}_{{\bf csq}'_{n+1}X}} &\mathrm{Hom}_{\simp}\big(\Lambda^{n+2,k},{\bf csq}'_{n+1}(X)\big)\,,}
$$
où $\alpha^{n+1,k}_{X}$ est une fonction surjective. 

Dans ce cas la fonction ${\alpha^{n+1,k}_{{\bf csq}'_{n+1}X}}$ est surjective, parce que:
$$
\xymatrix{\mathrm{Hom}_{\simp}(\Lambda^{n+2,k},X)\ar[r]& \mathrm{Hom}_{\simp}\big(\Lambda^{n+2,k},{\bf csq}'_{n+1}(X)\big)}
$$
est une fonction surjective, vu qu'elle s'insère dans un diagramme:
$$
\xymatrix@C-1pt{
\underset{i,j\neq k}{\underset{0\leq i<j\leq m+1}{\prod}}X_n\ar[d]|-{\textit{identité}}&&&
\underset{l\neq k}{\underset{0\leq l\leq m+1}{\prod}}X_{n+1}\ar[d]|-{\textit{quotient}}
\ar@<+7pt>[lll]^-{\underset{i<j}{\prod} \; d_{j-1}\circ\,{proj}_i}
\ar@<-7pt>[lll]_-{\underset{i<j}{\prod} \; d_{i} \circ\,{proj}_j}&
\ar[l]\mathrm{Hom}_{\simp}(\Lambda^{n+2,k},X)\ar[d]\\
\underset{i,j\neq k}{\underset{0\leq i<j\leq m+1}{\prod}}X_n&&&
\underset{l\neq k}{\underset{0\leq l\leq m+1}{\prod}}\big(X_{n+1}/\overset{n}{\sim}\big)
\ar@<+7pt>[lll]^-{\underset{i<j}{\prod} \; d_{j-1}\circ\,{proj}_i}
\ar@<-7pt>[lll]_-{\underset{i<j}{\prod} \; d_{i}\circ\,{proj}_j}&
\ar[l]\mathrm{Hom}_{\simp}\big(\Lambda^{n+2,k},{\bf csq}'_{n+1}(X)\big)}
$$ 
dont les lignes sont de noyaux de double flèches.

Enfin, vérifions que le morphisme \eqref{tildeco} est une $\infty$-équivalence si $X$ satisfait la condition d'extension de Kan en dimension $1\leq m\leq n+1$. En effet, vu que \eqref{tildeco} est un morphisme entre complexes de Kan $(n+1)$-cosquelettiques, d'après le Corollaire \ref{lecoro} il suffit de voir que \eqref{tildeco} est une $n$-équivalence faible; ce qui se suit de la Proposition \ref{kankan} car la fonction:
$$
\xymatrix{{\bf csq}_{n+1}(X)_{m} \ar[r]& {\bf csq}'_{n+1}(X)_{m}}
$$
est bijective si $0\leq m\leq n$ et surjective si $m=n+1$.
\end{proof}

Montrons:

\begin{lemme} \label{lemehomofaib}
Soit $n\geq 0$. Si $X$ et $Y$ sont des ensembles simpliciaux faiblement $n$-cosquelettiques alors le produit cartésien argument par argument $X\times Y$ et l'ensemble des morphismes $\underline{\mathrm{Hom}}_{\simp}\big(X,Y\big)$ sont des ensembles simpliciaux faiblement $n$-cosquelettiques.
\end{lemme}
\begin{proof}
Supposons que $X$ et $Y$ sont des ensembles simpliciaux faiblement $n$-cos\-que\-le\-ttiques. Il se suit des Lemmes \ref{cosque} et \ref{lemehomo} que les ensembles simpliciaux $X\times Y$ et $\underline{\mathrm{Hom}}_{\simp}\big(X,Y\big)$ sont $(n+1)$-cosquelettiques. 

L'ensemble simplicial $X\times Y$ est faiblement $n$-cosquelettique parce que la fonction:
\begin{equation}\label{injenosenueff1}
\xymatrix@C+10pt{ \mathrm{Hom}_{\simp}\big(\Delta^{n+1}, X\times Y\big) \ar[r] & \mathrm{Hom}_{\simp}\big(\partial\Delta^{n+1}, X\times Y\big)}
\end{equation}
induit de l'inclusion canonique $\xymatrix@C-8pt{\partial\Delta^{n+1}\ar@{^(->}[r]&\Delta^{n+1}}$ est à isomorphisme près le produit des fonctions injectives:
$$
\xymatrix@C-5pt{ \mathrm{Hom}_{\simp}\big(\Delta^{n+1}, X\big) \ar[r] & \mathrm{Hom}_{\simp}\big(\partial\Delta^{n+1}, X\big)}
\qquad\text{et}\qquad
\xymatrix@C-5pt{ \mathrm{Hom}_{\simp}\big(\Delta^{n+1}, Y\big) \ar[r] & \mathrm{Hom}_{\simp}\big(\partial\Delta^{n+1}, Y\big)}
$$
Donc \eqref{injenosenueff1} injective.  

Finalement remarquons que d'après le Lemme \ref{lemehomo} pour montrer que l'ensemble simplicial $\underline{\mathrm{Hom}}_{\simp}\big(X,Y\big)$ est faiblement $n$-cosquelettique, il suffit de vérifier que la fonction:
\begin{equation}\label{injenosenueff}
\xymatrix@C+10pt{
\mathrm{Hom}_{\simp_{\leq n+1}}\big(\tau_{n+1}^{\phantom{A}*}(X)\times\tau_{n+1}^{\phantom{A}*}(\Delta^{n+1}),\tau_{n+1}^{\phantom{A}*}Y\Big)
\ar[r] & 
\mathrm{Hom}_{\simp_{\leq n+1}}\big(\tau_{n+1}^{\phantom{A}*}(X)\times\tau_{n+1}^{\phantom{A}*}(\partial\Delta^{n+1}),\tau_{n+1}^{\phantom{A}*}Y\Big)}
\end{equation} 
induit de l'inclusion canonique $\xymatrix@C-8pt{\partial\Delta^{n+1}\ar@{^(->}[r]&\Delta^{n+1}}$ est injective.

De façon explicite il faut montrer que si $F,G\colon\xymatrix@C-4pt{\tau_{n+1}^{\phantom{A}*}(X)\times \tau_{n+1}^{\phantom{A}*}(\Delta^{n+1}) \ar[r] & \tau_{n+1}^{\phantom{A}*}(Y)}$ sont de morphismes d'ensembles simpliciaux tronqués tels que:
\begin{equation}\label{propriecuecuer}
\text{$F_{k}=G_{k}$ \; si \; $0\leq k\leq n$ \quad  et \quad $F_{n+1}(x,f) = G_{n+1}(x,f)$ \; si \; $x\in X_{n+1}$ \: et \; $f\in (\partial\Delta^{n+1})_{n+1}$;}
\end{equation}
alors $F_{n+1}(x,\mathrm{id}_{[n+1]}) = G_{n+1}(x,\mathrm{id}_{[n+1]})$ si $x\in X_{n+1}$ c'est-à-dire $F=G$.

Soient $F$ et $G$ deux morphismes vérifiant les propriétés \eqref{propriecuecuer}. Vu qu'on a supposé que $Y$ soit un ensemble simplicial faiblement $n$-cosquelettique, pour montrer que $F_{n+1}(x,\mathrm{id}_{[n]}) = G_{n+1}(x,\mathrm{id}_{[n]})$ pour tout $x\in X_{n+1}$ il suffit de vérifier que $d_i\circ F_{n+1}(x,\mathrm{id}_{[n]}) = d_i\circ G_{n+1}(x,\mathrm{id}_{[n]})$ pour tout $x\in X_{n+1}$ et $0\leq i\leq n+1$. Mais par hypothèse $F$ et $G$ vérifient:
$$
d_i\circ F_{n+1}(x,\mathrm{id}_{[n]}) \, = \, F_n(d_ix,\delta_i) \, = \, G_n(d_ix,\delta_i)  \, = \, d_i\circ G_{n+1}(x,\mathrm{id}_{[n]}).
$$ 
Donc $F=G$.
\end{proof}

\renewcommand{\thesubsection}{\S\thesection.\arabic{subsection}}
\subsection{}\;\label{groupoideDG}
\renewcommand{\thesubsection}{\thesection.\arabic{subsection}}

Soit $0\leq n\leq \infty$. Rappelons qu'un ensemble simplicial $X$ est appelé un \emph{$n$-groupoïde de Kan}, si $X$ est un complexe de Kan qui satisfait la condition d'extension de Kan de fa\c con stricte en dimension $m\geq n$. 

Explicitement, $X$ est un $n$-groupoïde de Kan si la fonction:
$$
\xymatrix@C+10pt{
\mathrm{Hom}_{\simp}\big(\Delta^{m+1},X\big) \ar[r]^{\alpha^{m,k}_{X}} &\mathrm{Hom}_{\simp}\big(\Lambda^{m+1,k},X\big)}
$$
induite de l'inclusion $\xymatrix@C-3pt{\Lambda^{m+1,k}\,\ar@{^(->}[r]^-{\alpha^{m,k}}&\Delta^{m+1}}$ est surjective pour $m\geq 1$ et $0\leq k\leq m+1$, et injective pour $m\geq n$ et $0\leq k\leq m+1$.

Montrons les équivalences suivantes:

\begin{lemme}\label{uneequ}
Si $X$ est un ensemble simplicial quelconque, les énoncés suivants sont équivalents:
\begin{enumerate}
\item $X$ est un $n$-groupoïde de Kan.
\item $X$ est un ensemble simplicial $(n+1)$-cosquelettique qui satisfait la condition d'extension de Kan en dimension $m$ pour $1\leq m\leq n+1$ et tel que la fonction:
$$
\xymatrix@C+10pt{
\mathrm{Hom}_{\simp}\big(\Delta^{n+1},X\big) \ar[r]^{\alpha^{n,k}_{X}} &\mathrm{Hom}_{\simp}\big(\Lambda^{n+1,k},X\big)}
$$
est injective pour tout $0\leq k\leq n+1$. 
\item $X$ est un ensemble simplicial faiblement $n$-cosquelettique qui satisfait la condition d'extension de Kan en dimension $m$ pour $1\leq m\leq n+1$ et tel que la fonction:
$$
\xymatrix@C+10pt{
\mathrm{Hom}_{\simp}\big(\Delta^{n+1},X\big) \ar[r]^{\alpha^{n,k}_{X}} &\mathrm{Hom}_{\simp}\big(\Lambda^{n+1,k},X\big)}
$$
est injective pour tout $0\leq k\leq n+1$. 
\end{enumerate}

En particulier, d'après le Corollaire \ref{lecoro} un $n$-groupoïde de Kan est un ensemble simplicial $n$-fibrant \emph{i.e.} est un objet fibrant de la catégorie de modèles $(\simp,{\bf W}_n, {\bf mono},{\bf fib}_n)$ du Théorème \ref{ntypess}. 
\end{lemme}
\begin{proof}
Rappelons qu'un ensemble simplicial faiblement $n$-cosquelettique est par définition $(n+1)$-cosquelettique. Donc (iii) $\Rightarrow$ (ii). D'un autre côté vu que $\alpha^{n,k}_{X}\,=\,\widetilde{\alpha}_X^{n,k}\circ \alpha_X^{n}$, si la fonction $\alpha^{n,k}_{X}$ est injective alors $\alpha_X^{n}$ l'est aussi; c'est-à-dire (ii) $\Rightarrow$ (iii).

Pour montrer que (i) $\Rightarrow$ (iii) considérons un $n$-groupoïde de Kan $X$. Remarquons pour commencer que par hypothèse $X$ satisfait la condition d'extension de Kan en dimension $m$ pour $1\leq m\leq n+1$ et que la fonction $\alpha^{n,k}_X$ est injective pour tout $0\leq k\leq n+1$. 
 
Plus encore, vu que par supposition la fonction $\alpha^{m,k}_{X}$ du diagramme:
\begin{equation}\label{trinLemmG}
\xymatrix@C+10pt{
\mathrm{Hom}_{\simp}\big(\Delta^{m+1},X\big)
\ar@<-3pt>@/_16pt/[rr]_{\alpha^{m,k}_{X}}\ar[r]^{\alpha^{m}_{X}}&
\mathrm{Hom}_{\simp}\big(\partial\Delta^{m+1},X\big)\ar[r]^{\widetilde{\alpha}^{m,k}_{X}}& 
\mathrm{Hom}_{\simp}\big(\Lambda^{m+1,k},X\big)}\,,
\end{equation}
est une fonction injective pour $m\geq n$ et $0\leq k\leq m+1$, on a que ${\alpha}^{m}_{X}$ est une fonction injective si $m\geq n$. Il se suit du Lemme \ref{glennlemme} que $\widetilde{\alpha}^{m,k}_{X}$ est une fonction injective pour $m\geq n+1$ et $0\leq k\leq m+1$; donc ${\alpha}^{m}_{X}$ est une fonction surjective si $m\geq n+1$ parce que ${\alpha}^{m,k}_{X}$ est surjective et $\widetilde{\alpha}^{m,k}_{X}$ injective toujours que $m\geq n+1$ et $0\leq k\leq m+1$. Autrement dit $X$ est un ensemble simplicial faiblement $n$-cosquelettique.

Finalement supposons que $X$ est un ensemble simplicial vérifiant les conditions de l'énoncé (iii). D'après le Lemme \ref{dusl} l'ensemble simplicial $X$ est un complexe de Kan, or pour montre que $X$ est un $n$-groupoïde de Kan il faut juste de montrer que la fonction $\alpha_X^{m,k}$ est injective si $m\geq n+1$ et $0\leq k\leq m+1$. Vu que par hypothèse la fonction $\alpha^m_X$ est une fonction injective si $m\geq n$, d'après le Lemme \ref{glennlemme} la fonction $\widetilde{\alpha}_X^{m,k}$ dans le triangle \eqref{trinLemmG} ci-dessus est aussi injective si $m\geq n+1$ et $0\leq k\leq m+1$. Donc $\alpha_X^{m,k}$ est une fonction injective si $m\geq n+1$ et $0\leq k\leq m+1$.
\end{proof}

On déduit du Lemme \ref{uneequ} que la sous-catégorie pleine de $\simp$ dont les objets sont les $0$-groupoïdes de Kan ${\bf Grpd}^{0}_{Kan}$, est équivalente à la catégorie des ensembles. 

En effet, rappelons qu'on a des adjonctions:
\begin{equation}\label{0cosque}
{\bf Ens}
\vcenter{\xymatrix@C+23pt{
\phantom{.}\ar[r]|-{\,[0]_!\,}&\phantom{.}
\ar@<+4pt>@/^12pt/[l]^-{[0]^*}_{\perp}
\ar@<-4pt>@/_12pt/[l]_-{\pi_0}^{\perp}
}}\simp
\end{equation}
où $\pi_0(X)$ est l'ensemble des composantes connexes par arcs de l'ensemble simplicial $X$, $[0]_!(A)$ est l'ensemble simplicial constant à valeurs l'ensemble $A$ et $[0]^*(X)=X_0$ est l'ensemble des $0$-simplexes de $X$.

Montrons:

\begin{corollaire}\label{0grpdclassi}
Si $X$ est un ensemble simplicial, ils sont équivalents:
\begin{enumerate}
\item $X$ est un $0$-groupoïde de Kan.
\item $X$ est un \emph{ensemble simplicial constant}, \emph{i.e.} les morphismes faces et dégénérescences de $X$ sont tous des fonctions bijectives.
\item $X$ est isomorphe à un ensemble simplicial dans l'image du foncteur $[0]_!$.
\item Si $\eta$ est une unité de l'adjonction $\pi_0 \dashv [0]_!$ le morphisme $\eta_X\colon \xymatrix@C-5pt{X\ar[r]&[0]_!\circ \pi_0 (X)}$ est un isomorphisme d'ensembles simpliciaux.
\item Si $\varepsilon$ est une counité de l'adjonction $[0]_! \dashv [0]^*$ le morphisme $\varepsilon_X\colon \xymatrix@C-5pt{[0]_!\circ [0]^* (X)\ar[r]&X}$ est un isomorphisme d'ensembles simpliciaux.
\end{enumerate}

En particulier, le foncteur ensemble simplicial constant $\xymatrix@C-2pt{{\bf Ens}\ar[r]^-{[0]_!}&\simp}$ induit une équivalence entre la catégorie des ensembles et la catégorie des $0$-groupoïdes de Kan.
\end{corollaire}
\begin{proof}
On déduit sans peine de la définition du foncteur $[0]_!$ que les énoncés (ii) et (iii) sont équivalents. D'un autre vu que $[0]_!$ est un foncteur pleinement fidèle, les énoncés (iii), (iv) et (v) sont équivalents.

D'un autre côté il se suit du Lemme \ref{uneequ} que si $X$ est un $0$-groupoïde de Kan alors $X$ est un ensemble simplicial $1$-cosquelettique et la fonction:
\begin{equation}	\label{ogrpdclass}
\xymatrix@C+10pt{
\mathrm{Hom}_{\simp}\big(\Delta^{1},X\big)  \ar[r]^{\alpha^{0,k}_X} & \mathrm{Hom}_{\simp}\big(\Lambda^{1,k},X\big)}
\end{equation}
est bijective si $0\leq k\leq 1$; c'est-à-dire les fonctions $d_0,d_1\colon\xymatrix@C-10pt{X_1\ar[r]&X_0}$ sont des bijections. On en déduit aussi-tôt que $X$ est un ensemble simplicial constant.

Réciproquement on vérifie sans peine qu'un ensemble simplicial constant est toujours un complexe de Kan, $1$-cosquelettique pour lequel la fonction \eqref{ogrpdclass} est bijective. Donc (i) $\Leftrightarrow$ (ii).
\end{proof}

Si $0\leq n\leq \infty$ montrons que la sous-catégorie pleine de $\simp$ dont les objets sont les $n$-groupoïdes de Kan noté ${\bf GrpdK}^{n}$ est une sous-catégorie cartésienne fermée:

\begin{lemme}\label{cartfermengrpKAN}
Soit $0\leq n\leq \infty$, si $X$ et $Y$ sont des $n$-groupoïdes de Kan l'ensemble simplicial produit $X\times Y$ et l'ensemble simplicial des morphismes $\underline{\mathrm{Hom}}_{\simp}(X,Y)$ sont de $n$-groupoïdes de Kan. Autrement dit la catégorie ${\bf GrpdK}^{n}$ des $n$-groupoïdes de Kan est une sous-catégorie cartésienne fermée de $\simp$.
\end{lemme}
\begin{proof}
Supposons que $X$ et $Y$ sont des $n$-groupoïdes de Kan. On vérifie aussi-tôt que l'ensemble simplicial produit $X\times Y$ est un $n$-groupoïde de Kan parce que la fonction:
$$
\xymatrix@C+10pt{
\mathrm{Hom}_{\simp}\big(\Delta^{m+1},X\times Y\big) \ar[r]^{\alpha^{m,k}_{X\times Y}} &\mathrm{Hom}_{\simp}\big(\Lambda^{m+1,k},X\times Y\big)}
$$
induite de l'inclusion $\xymatrix@C-3pt{\Lambda^{m+1,k}\,\ar@{^(->}[r]^-{\alpha^{m,k}}&\Delta^{m+1}}$ s'identifie à isomorphisme près à la fonction produit:
$$
\xymatrix@C+10pt{
\mathrm{Hom}_{\simp}\big(\Delta^{m+1},X\big) \times \mathrm{Hom}_{\simp}\big(\Delta^{m+1},Y\big) \ar[r]^{\alpha^{m,k}_{X}\times \alpha^{m,k}_{Y}} &\mathrm{Hom}_{\simp}\big(\Lambda^{m+1,k},X\big) \times \mathrm{Hom}_{\simp}\big(\Lambda^{m+1,k},Y\big)}
$$
laquelle est surjective pour $m\geq 1$ et $0\leq k\leq m+1$, et injective pour $m\geq n$ et $0\leq k\leq m+1$.

Remarquons d'un autre côté que d'après les Lemmes \ref{nfibcartesiene} et \ref{lemehomofaib} $\underline{\mathrm{Hom}}_{\simp}\big(X,Y\big)$ est un complexe de Kan faiblement $n$-cosquelettique. Il se suit des Lemmes \ref{uneequ} et \ref{lemehomo} que pour montrer que $\underline{\mathrm{Hom}}_{\simp}\big(X,Y\big)$ soit un $n$-groupoïde de Kan on doit juste vérifier que si $0\leq k\leq n+1$ la fonction:
\begin{equation}\label{injenosenueN}
\xymatrix@C+10pt{
\mathrm{Hom}_{\simp_{\leq n+1}}\big(\tau_{n+1}^{\phantom{A}*}(X)\times\tau_{n+1}^{\phantom{A}*}(\Delta^{n+1}),\tau_{n+1}^{\phantom{A}*}Y\Big)
\ar[r] & 
\mathrm{Hom}_{\simp_{\leq n+1}}\big(\tau_{n+1}^{\phantom{A}*}(X)\times\tau_{n+1}^{\phantom{A}*}(\Lambda^{n+1,k}),\tau_{n+1}^{\phantom{A}*}Y\Big)}
\end{equation} 
induit de l'inclusion canonique $\xymatrix@C-8pt{\Lambda^{n+1,k}\ar@{^(->}[r]&\Delta^{n+1}}$ est injective.

Soit $0\leq k\leq n+1$. Vu qu'on a les égalités des ensembles:
$$
\Delta^{n+1}_{n+1} \Big\backslash \Lambda^{n+1,k}_{n+1} \; =  \; 
\Big\{ \; \xymatrix@C-8pt{[n+1]\ar[r]^-{\mathrm{id}_{[n+1]}}&[n+1]}\,, \;\; \xymatrix@C-8pt{[n+1]\ar[r]^-{\sigma_i}&[n]\ar[r]^-{\delta_k}&[n+1]} \; \text{où \; $0\leq i\leq n$}\; \Big\} \,,
$$
$$
\Delta^{n+1}_{n} \Big\backslash \Lambda^{n+1,k}_{n} \; =  \; 
\Big\{ \; \xymatrix@C-8pt{[n]\ar[r]^-{\delta_k}&[n+1]} \; \; \Big\} \qquad \text{et} \qquad 
\Delta^{n+1}_{n-1} \Big\backslash \Lambda^{n+1,k}_{n-1} \; =  \;  \emptyset\,,
$$
il faut montrer que si $F,G\colon\xymatrix@C-7pt{\tau_{n+1}^{\phantom{A}*}(X)\times \tau_{n+1}^{\phantom{A}*}(\Delta^{n+1}) \ar[r] & \tau_{n+1}^{\phantom{A}*}(Y)}$ sont de morphismes d'ensembles simpliciaux tronqués tels que:
\begin{equation}\label{propriecuecuer2}
\text{$F_{k}=G_{k}$ \; si \; $0\leq k\leq n-1$, \quad  $F_{n}(a,f) = G_{n}(a,f)$ \; si \; $a\in X_{n}$ \; et \; $f\in \Lambda^{n+1,k}_{n}$}
\end{equation}
$$
\text{et \quad $F_{n+1}(x,f) = G_{n+1}(x,\varphi)$ \; si \; $x\in X_{n+1}$ \: et \; $\varphi\in \Lambda^{n+1,k}_{n+1}$;}
$$
alors:
\begin{equation}
\begin{split}
\text{$F_{n}(a,\delta_k) = G_{n}(a,\delta_k)$ \; si \; $a\in X_{n}$}  \hspace{4cm}\\
\text{$F_{n+1}(x,\mathrm{id}_{[n+1]}) = G_{n+1}(x,\mathrm{id}_{[n+1]})$ \; et \; $F_{n+1}(x,\delta_k\sigma_i) = G_{n+1}(x,\delta_k\sigma_i)$ \; si \; $x\in X_{n+1}$}\,.
\end{split}
\end{equation}

Pour montrer que $F_{n+1}(x,\mathrm{id}_{[n]}) = G_{n+1}(x,\mathrm{id}_{[n]})$ pour tout $x\in X_{n+1}$ il suffit de vérifier que $d_i\circ F_{n+1}(x,\mathrm{id}_{[n]}) = d_i\circ G_{n+1}(x,\mathrm{id}_{[n]})$ pour tout $x\in X_{n+1}$ et $0\leq i\leq n+1$ où $i\neq k$, vu que par hypothèse $Y$ est un $n$-groupoïde de Kan. Mais par hypothèse:
$$
d_i\circ F_{n+1}(x,\mathrm{id}_{[n]}) \, = \, F_n(d_ix,\delta_i) \, = \, G_n(d_ix,\delta_i)  \, = \, d_i\circ G_{n+1}(x,\mathrm{id}_{[n]})\,.
$$ 
pour tout $x\in X_{n+1}$ et $0\leq i\leq n+1$ où $i\neq k$. Donc $F_{n+1}(x,\mathrm{id}_{[n]}) = G_{n+1}(x,\mathrm{id}_{[n]})$ pour tout $x\in X_{n+1}$ .

Il se suit en particulier que pour tout $a\in X_n$:
$$
F_{n}(a,\delta_k) = d_k \circ F_{n+1}\big( s_k (a) ,\mathrm{id}_{[n]}\big) = d_k \circ G_{n+1}\big( s_k (a),\mathrm{id}_{[n]}) = G_{n}(a,\delta_k) \qquad \text{si \; $0\leq k\leq n$}
$$
$$\text{ou}$$
$$
F_{n}(a,\delta_k) = d_k \circ F_{n+1}\big( s_{n} (a) ,\mathrm{id}_{[n]}\big) = d_k \circ G_{n+1}\big( s_n (a),\mathrm{id}_{[n]}) = G_{n}(a,\delta_k) \qquad \text{si \; $k=n+1$\,.}
$$

Enfin on montre que $F_{n+1}(x,\delta_k\sigma_i) = G_{n+1}(x,\delta_k\sigma_i)$ si $x\in X_{n+1}$ parce que:
$$
d_j \circ F_{n+1}(x,\delta_k\sigma_i) = F_{n}\big(d_i(x),\delta_k\sigma_i\delta_j\big)  = G_{n}\big(d_i(x),\delta_k\sigma_i\delta_j\big) = d_j \circ G_{n+1}(x,\delta_k\sigma_i)
$$
pour tout $0\leq j\leq n+1$ et $Y$ est un $n$-groupoïde de Kan. Donc $F=G$.
\end{proof}

\renewcommand{\thesubsection}{\S\thesection.\arabic{subsection}}
\subsubsection{}\;  
\renewcommand{\thesubsection}{\thesection.\arabic{subsection}}

On va dire qu'un ensemble simplicial $X$ \emph{satisfait la condition de minimalité en dimension $m$} si pour tout $0\leq k\leq m+1$ la fonction:
$$
\xymatrix@C+10pt{
\mathrm{Hom}_{\simp}\big(\partial\Delta^{m+1},X\big) \ar[r]^{\widetilde{\alpha}^{m,k}_{X}} &\mathrm{Hom}_{\simp}\big(\Lambda^{m+1,k},X\big)},
$$ 
est injective lorsque on la restreinte à l'image de la fonction:
$$
\xymatrix@C+10pt{
\mathrm{Hom}_{\simp}\big(\Delta^{m+1},X\big) \ar[r]^{\alpha^{m,k}_{X}} &\mathrm{Hom}_{\simp}\big(\partial\Delta^{m+1,k},X\big)}\,.
$$

Si $n\geq 0$, un ensemble simplicial est dit \emph{$n$-minimal} si satisfait la condition de minimalité en dimension $m$ pour tout $m\geq n$ . Un ensemble simplicial $0$-minimal est appelé simplement un \emph{ensemble simplicial minimal} (voir II\S9 de \cite{may}).

\begin{corollaire}\label{groupidequi}
Si $n\geq 0$ et $X$ est un ensemble simplicial quelconque, les énoncés suivants sont équivalents:
\begin{enumerate}
\item $X$ est un $n$-groupoïde de Kan.
\item $X$ est un complexe de Kan $n$-minimal et faiblement $n$-cosquelettique.
\item $X$ est un ensemble simplicial faiblement $n$-cosquelettique, lequel satisfait la condition d'extension de Kan en dimension $m$ pour $1\leq m \leq n+1$ et satisfait la condition de minimalité en dimension $n$.
\end{enumerate}

Également on a les équivalences:
\begin{enumerate}
\item[(iv)] $X$ est un $n$-groupoïde de Kan minimal.
\item[(v)] $X$ est un complexe de Kan faiblement $n$-cosquelettique minimal.
\item[(vi)] $X$ est un ensemble simplicial faiblement $n$-cosquelettique, lequel satisfait la condition de minimalité en dimension $m$ pour $0\leq m\leq n$  et satisfait la condition d'extension de Kan en dimension $m$ pour $1\leq m \leq n+1$.
\end{enumerate}
\end{corollaire}
\begin{proof}
Si $X$ est un ensemble simplicial, du triangle commutatif:
$$
\xymatrix@C+10pt{
\mathrm{Hom}_{\simp}\big(\Delta^{n+1},X\big)
\ar@<-3pt>@/_16pt/[rr]_{\alpha^{n,k}_{X}}\ar[r]^{\alpha^{n}_{X}}&
\mathrm{Hom}_{\simp}\big(\partial\Delta^{n+1},X\big)\ar[r]^{\widetilde{\alpha}^{n,k}_{X}}& 
\mathrm{Hom}_{\simp}\big(\Lambda^{n+1,k},X\big)};
$$
on déduit que si la fonction $\alpha^{n}_{X}$ est  injective, par exemple si $X$ est un ensemble simplicial faiblement $n$-cosquelettique, alors $X$ satisfait la condition de minimalité en dimension $n$ si et seulement si la fonction $\alpha_{X}^{n,k}$ est injective. Il suit du Lemme \ref{uneequ} que les énoncés (i) et (iii) sont équivalents. 

D'après le Corollaire \ref{lecoro} pour montrer que (ii) $\Leftrightarrow$ (iii) il suffit de montrer l'énoncé suivant:
\begin{lemme}\label{dese}
Soit $n\geq 0$. Si $X$ est un ensemble simplicial $(n+1)$-cosquelettique (resp. faiblement $n$-cosquelettique), alors $X$ est satisfait la condition de minimalité en dimension $m$ pour $m\geq n+2$ (resp. $m\geq n+1$). 
\end{lemme}
\begin{proof}
Si $X$ est un ensemble simplicial $(n+1)$-cosquelettique, dans le triangle commutatif:
\begin{equation}
\xymatrix@C+10pt{
\mathrm{Hom}_{\simp}\big(\Delta^{m+1},X\big)
\ar@<-3pt>@/_16pt/[rr]_{\alpha^{m,k}_{X}}\ar[r]^{\alpha^{m}_{X}}&
\mathrm{Hom}_{\simp}\big(\partial\Delta^{m+1},X\big)\ar[r]^{\widetilde{\alpha}^{m,k}_{X}}& 
\mathrm{Hom}_{\simp}\big(\Lambda^{m+1,k},X\big)},
\end{equation}
la fonction $\alpha^m_{X}$ est bijective pour $m\geq n+1$. Il se suit de l'énoncé (ii) du Lemme \ref{glennlemme} que $\widetilde{\alpha}^{m,k}_{X}$ est une fonction bijective pour $m\geq n+2$ et $0\leq k \leq m+1$; en particulier, $X$ est $m$-minimal pour $m\geq n+2$.

D'un autre côté, si on suppose que $X$ soit faiblement $n$-cosquelettique, on a cette fois que la fonction $\alpha^m_{X}$ est injective pour $m\geq n$. Donc, par le même énoncé $X$ est un ensemble simplicial $m$-minimal pour $m\geq n+1$. 
\end{proof}

Enfin on déduit facilement les équivalences (iv) $\Leftrightarrow$ (v) $\Leftrightarrow$ (vi) à partir des équivalences (i) $\Leftrightarrow$ (ii) $\Leftrightarrow$ (iii) déjà montrées. 
\end{proof}

Rappelons de \cite{may}:

\begin{proposition}\label{propmay}
Soit $n\geq 0$. Si $X$ est un ensemble simplicial quelconque il existe un $n$-groupoïde de Kan minimal $W$ et un isomorphisme $\xymatrix@C-5pt{W\ar[r]&X}$ de la catégorie homotopique de la catégorie de modèles $(\simp,{\bf W}_n, {\bf mono},{\bf fib}_n)$ du Théorème \ref{ntypess}.

De façon analogue, si $X$ est un ensemble simplicial réduit il existe un $n$-groupoïde de Kan minimal $W$ lequel est un ensemble simplicial réduit (\emph{i.e.} un $n$-groupe de Kan minimal) et un isomorphisme $\xymatrix@C-5pt{W\ar[r]&X}$ de la catégorie homotopique de la catégorie de modèles $(\simp_0,{\bf W}^{red}_n, {\bf mono},{\bf fib}^{red}_n)$ de la Proposition \ref{modred}.
\end{proposition}
\begin{proof}
Si $X$ est un ensemble simplicial quelconque, on sait qu'il existe un complexe de Kan $Y$ et une $\infty$-équivalence faible d'ensembles simpliciaux $\xymatrix@C-5pt{X\ar[r]&Y}$.

D'un autre si $Y$ est un complexe de Kan, dans II\S9 de \cite{may} on peut trouver une preuve de l'existence d'un sous-ensemble simplicial $Z$ de $Y$ qui vérifie les propriétés suivantes:
\begin{enumerate}
\item $Z$ est un complexe de Kan minimal. 
\item Le morphisme d'inclusion $\xymatrix@C-10pt{Z\ar@{^(->}[r]&Y}$ est une $\infty$-équivalence faible.
\end{enumerate}

En particulier, d'après le Lemme \ref{dusl} on a un diagramme d'ensembles simpliciaux:
\begin{equation}\label{suisui} 
\xymatrix{
X\ar[d]&&\\Y&&\\ 
Z\ar[r]\ar@{^(->}[u]&{\bf csq}_{n+1}(Z)\ar[r]& {\bf csq}'_{n+1}(Z)}
\end{equation}
dont les morphismes sont des $n$-équivalences faibles et tel que ${\bf csq}'_{n+1}(Z)$ est un $n$-groupoïde de Kan.

 Montrons l'énoncé suivant:
\begin{lemme}\label{mininifaible}
Soit $n\geq 0$. Si $X$ est un ensemble simplicial lequel satisfait la condition de minimalité en dimension $m$ pour $0\leq m\leq n$, alors ${\bf csq}_{n+1}(X)'$ satisfait la condition de minimalité en dimension $m$ pour tout $m\geq 0$.
\end{lemme}
\begin{proof}
D'après le Lemme \ref{dese} il suffit de montrer que si $X$ satisfait la condition de minimalité en dimension $m$ pour $0\leq m\leq n$ alors ${\bf csq}_{n+1}(X)'$ aussi.

Premièrement ${\bf csq}_{n+1}'(X)$ est $m$-minimal pour $0\leq m\leq n-1$ car $X$ et ${\bf csq}_{n+1}'(X)$ ont la même $n$-troncations. Pour montrer que ${\bf csq}_{n+1}'(X)$ est $n$-minimal, considérons le diagramme commutatif:
$$
\xymatrix@R-3pt@C+7pt{
\mathrm{Hom}_{\simp}\big(\Delta^{n+1},X\big)  \ar[r]^-{\alpha^{n}_{X}}\ar[d]_-{f_{1}} & 
\mathrm{Hom}_{\simp}\big(\partial\Delta^{n+1},X\big) \ar[r]^{\widetilde{\alpha}^{n,k}_{X}}\ar[d]|-{f_{2}} & 
\mathrm{Hom}_{\simp}\big(\Lambda^{n+1,k},X\big)\ar[d]^-{f_{3}}\\
\mathrm{Hom}_{\simp}\big(\Delta^{n+1},{\bf csq}_{n+1}'X\big) \ar[r]_-{\alpha^{n}_{{\bf csq}_{n+1}'X}}& 
\mathrm{Hom}_{\simp}\big(\partial\Delta^{n+1},{\bf csq}_{n+1}'X\big) \ar[r]_-{\widetilde{\alpha}^{n,k}_{{\bf csq}_{n+1}'X}}&
\mathrm{Hom}_{\simp}\big(\Lambda^{n+1,k},{\bf csq}_{n+1}'X\big)\,,}
$$
où les fonctions $f_{i}$ sont induisent par $\xymatrix@C-7pt{X\ar[r]&{\bf csq}_{n+1}'X}$ le morphisme composé \eqref{tildecoco}. 

Remarquons que $f_{2}$ et $f_{3}$ sont des fonctions bijectives, car $X$ et ${\bf csq}_{n+1}'(X)$ ont la même $n$-troncation. D'un autre, vu qu'on peut identifier $f_{1}$ au morphisme quotient $\xymatrix@C-10pt{X_{n+1}\ar[r] & X_{n+1}/\overset{n}{\sim}}$, où $\overset{n}{\sim}$ est la relation d'équivalence \eqref{relationn}, $f_{1}$ est une fonction surjective. Donc, ${\bf csq}_{n+1}'(X)$ est un ensemble simplicial $n$-minimal si $X$ l'est. 
\end{proof}

On déduit du Lemme \ref{mininifaible} que le $n$-groupoïde de Kan ${\bf csq}'_{n+1}(Z)$ du diagramme \eqref{suisui} est un ensemble simplicial minimal. D'un autre vu que tous les morphismes de \eqref{suisui} sont des $n$-équivalences faibles, on en déduit un isomorphisme $\xymatrix@C-5pt{{\bf csq}'_{n+1}(Z)\ar[r]&X}$ de $\mathrm{Ho}_n(\simp)$ la catégorie homotopique de la catégorie de modèles $(\simp,{\bf W}_n, {\bf mono},{\bf fib}_n)$.

Finalement remarquons que si $X$ est un ensemble simplicial réduit on peut supposer que tous les ensembles simpliciaux du diagramme \eqref{suisui} sont aussi réduit. En effet pour commencer on choisi un remplacement fibrant $\xymatrix@C-5pt{X\ar[r]&Y}$ de $X$ dans la catégorie de modèles $(\simp_0,{\bf W}^{red}_n, {\bf mono},{\bf fib}^{red}_n)$. Puis on note que $Z$ est un sous-ensemble simplicial de $Y$, en particulier $Z$ est réduit. Enfin on note que les ensembles de $0$-simplexes de ${\bf csq}_{n+1}(Z)$ et ${\bf csq}_{n+1}(Z)'$ sont égaux à $Z_0=\star$.
\end{proof}

\renewcommand{\thesubsection}{\S\thesection.\arabic{subsection}}
\subsubsection{}\;  
\renewcommand{\thesubsection}{\thesection.\arabic{subsection}}

Si $n\geq 0$ rappelons que d'après les Lemmes \ref{cartfermengrpKAN}, \ref{lemehomo} et \ref{nfibcartesiene} les sous-catégories pleines ${\bf GrpdK}^n$, ${\bf Kcsq}^{n+1}$ et ${\bf Fib}^n$ de la catégorie des ensembles simpliciaux $\simp$ dont les objets sont les $n$-groupoïdes de Kan, les complexes de Kan $(n+1)$-cosquelettiques et les ensembles simpliciaux $n$-fibrants respectivement, sont de sous-catégories cartésiennes fermées de $\simp$.

Autrement dit, on a une chaîne des sous-catégories pleines:
\begin{equation}\label{chainenttyp}
\xymatrix{ 
\text{${\bf GrpdK}^n$} \; \ar@{^(->}[r]& \; \text{${\bf Kcsq}^{n+1}$}  \; \ar@{^(->}[r]& \; \text{${\bf Fib}^n$}   \; \ar@{^(->}[r]& \;  \simp}
\end{equation}
lesquelles sont stables par produits finis et par ensemble simplicial des morphismes $\underline{\mathrm{Hom}}_{\simp}$. 

Vu que dans $(\simp,{\bf W}_n, {\bf mono},{\bf fib}_n)$ tous les objets sont cofibrants et les objets fibrants sont par définition les ensembles simpliciaux $n$-fibrants, si on pose $h{\bf GrpdK}^n$, $h{\bf Kcsq}^{n+1}$ et $h{\bf Fib}^n$ pour noter les catégories qu'on obtient de ${\bf GrpdK}^n$, ${\bf Kcsq}^{n+1}$ et ${\bf Fib}^n$ respectivement en prenant l'ensemble des composantes connexes par arcs $\pi_0\big(\underline{\mathrm{Hom}}_{\simp}\big)$ de l'ensemble simplicial des morphismes $\underline{\mathrm{Hom}}_{\simp}$, les foncteurs d'inclusion de \eqref{chainenttyp} induisent une chaîne des foncteurs pleinement fidèles:
\begin{equation}\label{chainenttyph}
\xymatrix{ 
\text{$h{\bf GrpdK}^{n}$} \; \ar@{^(->}[r]& \; \text{$h{\bf Kcsq}^{n+1}$}  \; \ar@{^(->}[r]& \; \text{$h{\bf Fib}^n$}   \; \ar@{^(->}[r]& \;  \mathrm{Ho}_{n}(\simp)}
\end{equation}

La Proposition \ref{propmay} implique que ces foncteurs sont essentiellement surjectifs.

\begin{corollaire}\label{sonequivMOD}
La chaîne de foncteurs d'inclusion \eqref{chainenttyp} induite une chaîne des équivalences de catégories \eqref{chainenttyph}. En particulier la catégorie homotopique des $n$-groupoïdes de Kan $h{\bf GrpdK}^{n}$ est équivalente à la catégorie homotopique des $n$-types d'homotopie $\mathrm{Ho}_{n}(\simp)$.
\end{corollaire}

Remarquons:

\begin{corollaire}
Le foncteur canonique $\xymatrix@C-5pt{{\bf GrpdK}^0 \ar[r] & h{\bf GrpdK}^0}$ est un isomorphisme de catégories. En particulier le foncteur ensemble simplicial constant $\xymatrix@C-10pt{{\bf Ens}\ar[r]&\simp}$ induit une équivalence entre la catégorie cartésienne fermée des ensembles ${\bf Ens}$ et la catégorie cartésienne fermée des $0$-types d'homotopie $\mathrm{Ho}_0(\simp)$ \emph{i.e.} la catégorie homotopique des $0$-groupoïdes de Kan.
\end{corollaire}
\begin{proof}
Vu que le foncteur canonique $\xymatrix@C-5pt{{\bf GrpdK}^0 \ar[r] & h{\bf GrpdK}^0}$ est l'identité en objets, pour montrer que c'est un isomorphisme de catégories il suffit de vérifier que c'est un foncteur pleinement fidèle. 

Si $X$ et $Y$ sont de $0$-groupoïdes de Kan il se suit du Lemme \ref{cartfermengrpKAN} que l'ensemble simplicial $\underline{\mathrm{Hom}}_{\simp}(X,Y)$ est aussi un $0$-groupoïde de Kan. En particulier la fonction:
$$
\xymatrix{\underline{\mathrm{Hom}}_{\simp}(X,Y)_0 \ar[r] & \pi_0\big(\underline{\mathrm{Hom}}_{\simp}(X,Y)\big) }
$$
est bijective d'après le Corollaire \ref{0grpdclassi}. Autrement dit le foncteur $\xymatrix@C-5pt{{\bf GrpdK}^0 \ar[r] & h{\bf GrpdK}^0}$ est pleinement fidèle.

La deuxième assertion est une conséquence de ce qu'on vient de montrer, le Corollaire \ref{0grpdclassi} et le cas $n=0$ du Corollaire \ref{sonequivMOD}.
\end{proof}

\renewcommand{\thesubsection}{\S\thesection.\arabic{subsection}}
\subsubsection{}\;  
\renewcommand{\thesubsection}{\thesection.\arabic{subsection}}

Si $n\geq 1$ posons ${\bf GrpK}^n$, ${\bf Kcsq}_0^{n+1}$ et ${\bf Fib}_0^n$ pour noter les sous-catégories pleines de la catégorie des ensembles simpliciaux réduits$\simp_0$ dont les objets sont les $n$-groupoïdes de Kan, les complexes de Kan $(n+1)$-cosquelettiques et les ensembles simpliciaux $n$-fibrants respectivement. On appelle les objets de ${\bf GrpK}^n$ les $n$-\emph{groupes de Kan}.

Remarquons qu'on a une chaîne des foncteurs d'inclusion:
\begin{equation}\label{chainenttypred}
\xymatrix{ 
\text{${\bf GrpK}^n$} \; \ar@{^(->}[r]& \; \text{${\bf Kcsq}_0^{n+1}$}  \; \ar@{^(->}[r]& \; \text{${\bf Fib}_0^n$}   \; \ar@{^(->}[r]& \;  \simp}\,.
\end{equation}

Notons aussi $h{\bf GrpK}^n$, $h{\bf Kcsq}_0^{n+1}$ et $h{\bf Fib}_0^n$ les catégories qu'on obtient de ${\bf GrpK}^n$, ${\bf Kcsq}_0^{n+1}$ et ${\bf Fib}_0^n$ respectivement en prenant l'ensemble des composantes connexes par arcs $\pi_0\big(\underline{\mathrm{Hom}}_{\simp_0}\big)$ de l'ensemble simplicial des morphismes $\underline{\mathrm{Hom}}_{\simp_0}$ de $\simp_0$. Vu que dans la catégorie de modèles simpliciale $(\simp_0,{\bf W}^{red}_n, {\bf mono},{\bf fib}^{red}_n)$ tous les objets sont cofibrants et les objets fibrants sont par définition les ensembles simpliciaux $n$-fibrants réduits, les foncteurs d'inclusion de \eqref{chainenttypred} induisent une chaîne des foncteurs pleinement fidèles:
\begin{equation}\label{chainenttyphred}
\xymatrix{ 
\text{$h{\bf GrpK}^{n}$} \; \ar@{^(->}[r]& \; \text{$h{\bf Kcsq}_0^{n+1}$}  \; \ar@{^(->}[r]& \; \text{$h{\bf Fib}_0^n$}   \; \ar@{^(->}[r]& \;  \mathrm{Ho}_{n}(\simp_0)}
\end{equation}

\begin{corollaire}\label{sonequivMODred}
La chaîne de foncteurs d'inclusion \eqref{chainenttypred} induite une chaîne des équi\-va\-lences de catégories \eqref{chainenttyphred}. En particulier la catégorie homotopique des $n$-groupes de Kan $h{\bf GrpK}^{n}$ est équivalente à la catégorie homotopique des $n$-types d'homotopie réduits $\mathrm{Ho}_{n}(\simp_0)$.
\end{corollaire}
\begin{proof}
Voir la Proposition \ref{propmay}.
\end{proof}

Montrons:

\begin{corollaire}\label{equivgrpgrpK0}
Le foncteur canonique $\xymatrix@C-5pt{{\bf GrpK}^1 \ar[r] & h{\bf GrpK}^1\,\simeq\,\mathrm{Ho}_{1}(\simp_0)}$ est un isomorphisme de catégories (voir le Corollaire \ref{sonequivMODred} ci-dessus).
\end{corollaire}
\begin{proof}
Démontrons par ailleurs:

\begin{lemme} \label{HOMgrpdred}
Si $X$ et $W$ sont des ensembles simpliciaux réduits dont $W$ est un $1$-groupoïde de Kan, alors $\underline{\mathrm{Hom}}_{\simp_0}\big(X,W\big)$ est un $0$-groupoïde de Kan; \emph{i.e.} $\underline{\mathrm{Hom}}_{\simp_0}\big(X,W\big)$ est l'ensemble simplicial constant à valeurs $\mathrm{Hom}_{\simp_0}\big(X,W\big)$ (voir le Lemme \ref{0grpdclassi}). 

\end{lemme}
\begin{proof}
Rappelons que d'après le Lemme \ref{0grpdclassi} un ensemble simplicial $A$ est un $0$-groupoïde de Kan si et seulement si $A$ est un ensemble simplicial constant à valeurs $A_0$.

D'un autre remarquons que si $W$ est un $1$-groupoïde de Kan réduit alors $W$ est un objet fibrant de la catégorie de modèles $(\simp_0,{\bf W}^{red}_1, {\bf mono},{\bf fib}^{red}_1)$; en particulier l'ensemble simplicial $\underline{\mathrm{Hom}}_{\simp_0}\big(X,W\big)$ est un complexe de Kan pour tout ensemble simplicial $X$ réduit (voir le Lemme \ref{uneequ}).

Plus encore, vu que l'ensemble simplicial réduit $W$ est $2$-cosquelettique, $\underline{\mathrm{Hom}}_{\simp_0}\big(X,W\big)$ est aussi un ensemble simplicial $2$-cosquelettique d'après le Lemme \ref{lemehomored}.

Il nous reste à montrer que si $0\leq m\leq 1$ et $0\leq k\leq m+1$ la fonction (surjective):
$$
\xymatrix@R=3pt{
\mathrm{Hom}_{\simp}\Big( \Delta^{m+1},\underline{\mathrm{Hom}}_{\simp_0}\big(X,W\big)\Big)  \ar@{}[d]|-{\text{\rotatebox[origin=c]{90}{\Large $\cong$}}}\ar[r] &
\mathrm{Hom}_{\simp}\Big( \Lambda^{m+1,k},\underline{\mathrm{Hom}}_{\simp_0}\big(X,W\big)\Big) \ar@{}[d]|-{\text{\rotatebox[origin=c]{90}{\Large$\cong$}}}\\
\mathrm{Hom}_{\simp}\Big(  \big(X\times  \Delta^{m+1}\big) \big/ \big(\star\times  \Delta^{m+1}\big), W  \Big) & 
\mathrm{Hom}_{\simp}\Big(  \big(X\times  \Lambda^{m+1,k}\big) \big/ \big(\star\times  \Lambda^{m+1,k}\big), W  \Big)
}
$$
induit de l'inclusion $\xymatrix@C-5pt{\Lambda^{m+1,k}\ar[r] & \Delta^{m+1}}$ est injective.

\emph{Cas $m=1$:} Pour montrer qu'on a une bijection:
\begin{equation}\label{marininini}
\xymatrix{
\mathrm{Hom}_{\simp}\Big(  \big(X\times  \Delta^{2}\big) \big/ \big(\star\times  \Delta^{2}\big), W  \Big) \ar[r] &
\mathrm{Hom}_{\simp}\Big(  \big(X\times  \Lambda^{2,k}\big) \big/ \big(\star\times  \Lambda^{2,k}\big), W  \Big)}
\end{equation}
si $0\leq k\leq 2$, remarquons que par hypothèse:
$$
\xymatrix{
\mathrm{Hom}_{\simp}  \big(\star\times  \Delta^{2}, W  \big) \ar[r] &
\mathrm{Hom}_{\simp} \big(\star\times  \Lambda^{2,k}, W  \big)}
$$
est une fonction bijective et qu'on a déjà montré dans de Corollaire \ref{cartfermengrpKAN} que:
$$
\xymatrix{
\mathrm{Hom}_{\simp}  \big(X\times  \Delta^{2}, W  \big) \ar[r] &
\mathrm{Hom}_{\simp} \big(X\times  \Lambda^{2,k}, W  \big)}
$$
est aussi bijective si $W$ est un $1$-groupoïde. Donc \eqref{marininini} est une fonction bijective.

\emph{Cas $m=0$:}  Montrons que la fonction:
\begin{equation}\label{0grpreduit}
\xymatrix{
\mathrm{Hom}_{\simp}\Big(  \big(X\times  \Delta^{1}\big) \big/ \big(\star\times  \Delta^{1}\big), W  \Big) \ar[r] &
\mathrm{Hom}_{\simp}\Big(  \big(X\times  \Lambda^{1,k}\big) \big/ \big(\star\times  \Lambda^{1,k}\big), W  \Big)}
\end{equation}
est injective si $0\leq k\leq 1$. 

Explicitement, si $F,G\colon\xymatrix@C-8pt{X\times \Delta^1\ar[r]& W}$ sont deux morphismes d'ensembles simpliciaux tels que:
$$
F\;\Big|_{\star\times \Delta^{1}} \; = \; G\;\Big|_{\star\times \Delta^{1}} \; = \; \vcenter{\xymatrix@R=.5pt{\text{\scriptsize{Morphisme constant}}\\\text{\scriptsize{à valeurs $\star$}}}}
$$
$$
\text{et}\qquad \qquad F\;\Big|_{X\times \Lambda^{1,k}} \; = \; G\;\Big|_{X\times \Lambda^{1,k}}\quad\;\text{\scriptsize{(Remarquons que $\star\times  \Lambda^{1,k} \, \cong \, \star$)}}\,;
$$
on doit montrer que $F=G$.

Vu que $W$ est un ensemble simplicial faiblement $1$-cosquelettique et on a l'égalité:
$$
\Delta^1_1 \Big\backslash \Lambda^{1,k}_1 \; =  \; 
\Big\{ \; \xymatrix@C-8pt{[1]\ar[r]^-{\mathrm{id}_{[1]}}&[1]}\,, \;\; \xymatrix@C-8pt{[1]\ar[r]^-{\sigma_0}&[0]\ar[r]^-{\delta_k}&[1]} \; \Big\} \,,
$$
il suffit de démontrer que pour tout $1$-simplexe $a$ de $X$ on a que:
\begin{equation}\label{lesegagadesirr}
F_1\big(a,\mathrm{id}_{[1]}\big) \; = \; G_1\big(a,\mathrm{id}_{[1]}\big)\quad\text{et}\quad F_1\big(a,\delta_k\circ\sigma_0\big) \; = \; G_1\big(a,\delta_k\circ\sigma_0\big)\,.
\end{equation}

Si $a\in X_1$ considérons les $2$-simplexes $\big(s_1(a),\sigma_0\big)$ et $\big(s_0(a),\sigma_1\big)$ de $X$, et remarquons que:
\begin{align*}
d_0(s_1a,\sigma_0) \; = \; (s_0\star,\mathrm{id}_{[1]}) \qquad
d_1(s_1a,\sigma_0) \; & = \; (a,\mathrm{id}_{[1]}\big) \qquad
d_2(s_1a,\sigma_0) \; = \; (a,\delta_1\circ\sigma_0)
\\ & \text{et}\\
d_0(s_0a,\sigma_1)  \; = \; (a,\delta_0\circ\sigma_0) \qquad
d_1(s_0a,\sigma_1) \; & = \; (a,\mathrm{id}_{[1]}) \qquad
d_2(s_1a,\sigma_1) \; = \; (s_0\star,\mathrm{id}_{[1]})\,.
\end{align*}

On déduit que:
$$
F_2\big(s_1(a),\sigma_0\big) \; = \; G_2\big(s_1(a),\sigma_0\big) \qquad \text{et} \qquad F_2\big(s_0(a),\sigma_1\big) \; = \; G_2\big(s_0(a),\sigma_1\big)\,,
$$
parce que la fonction:
$$
\xymatrix{\mathrm{Hom}_{\simp}\big(\Delta^2,W\big)  \ar[r] & \mathrm{Hom}_{\simp}\big(\Lambda^{2},W\big)}
$$
induit de l'inclusion $\xymatrix@C-7pt{\Lambda^{2,m}\ar@{^(->}[r] & \Delta^{2}}$ est injective si $0\leq m\leq 2$ et on sait que:
$$
F_1\big(s_0\star,\mathrm{id}_{[1]}\big) \; = \; G_1\big(s_0\star,\mathrm{id}_{[1]}\big)\quad\text{et}\quad F_1\big(a,\delta_\ell\circ\sigma_0\big) \; = \; G_1\big(a,\delta_{\ell}\circ\sigma_0\big)
$$
si $\ell\neq k$.

En particulier: 
\begin{align*}
F_1(a,\delta_0\circ\sigma_0) \, = \,   d_1\circ F_2\big(s_0(a),\sigma_1\big) \, & = \, d_0\circ G_2\big(s_0(a),\sigma_1\big) \,  = \,  G_1(a,\delta_0\circ\sigma_0)  \\
F_1(a,\mathrm{id}_{[1]}) \, = \,   d_1\circ F_2\big(s_0(a),\sigma_1\big) \, & = \, d_1\circ G_2\big(s_0(a),\sigma_1\big) \,  = \,  G_1(a,\mathrm{id}_{[1]})   \\
F_1(a,\delta_1\circ\sigma_0) \, = \,   d_2\circ F_2\big(s_1(a),\sigma_0\big) \, & = \, d_2\circ G_2\big(s_1(a),\sigma_0\big) \,  = \,  G_1(a,\delta_1\circ\sigma_0)\,.
\end{align*}
\end{proof}

Si $X$ et $Y$ sont des $1$-groupes de Kan \emph{i.e.} des ensembles simpliciaux réduit lesquels sont des $1$-groupoïdes de Kan, il se suit du Lemme \ref{HOMgrpdred} que l'ensemble simplicial $\underline{\mathrm{Hom}}_{\simp_0}(X,Y)$ est un $0$-groupoïde de Kan. Donc:
$$
\xymatrix{\mathrm{Hom}_{\simp_0}(X,Y) \, \cong \, \underline{\mathrm{Hom}}_{\simp_0}(X,Y)_0 \ar[r] & \pi_0\big(\underline{\mathrm{Hom}}_{\simp_0}(X,Y)\big) }
$$
est une fonction bijective d'après le Corollaire \ref{0grpdclassi}. 

Autrement dit le foncteur $\xymatrix@C-5pt{{\bf GrpK}^1 \ar[r] & h{\bf GrpK}^1}$ est pleinement fidèle. Donc une équivalence de catégories.
\end{proof}

\renewcommand{\thesubsection}{\S\thesection.\arabic{subsection}}
\subsection{}\;\label{versionpointe}
\renewcommand{\thesubsection}{\thesection.\arabic{subsection}}

Dans le présent paragraphe, on va montrer des version pointée du Corollaire \ref{lecoro} et du Lemme \ref{lemehomo} qu'on a besoin pour montrer le Lemme \ref{talvezsi}. Pour cela considérons les adjonctions:
\begin{equation}\label{adjonctionpoin}
\ens_{\star}^{\text{\scriptsize{$\Big({\bf\Delta}_{\leq n+1}\Big)^{op}$}}}\;=\,\vcenter{
\xymatrix@C+15pt{
\simp_{\leq n+1,\,\star}\;\phantom{A}\ar@{}[r]|-{\perp\phantom{A}}
\ar@{<-}@/^18pt/[r]^-{\overline{\tau}_{n+1}^{\phantom{A} *}}
\ar@/_18pt/[r]_-{\overline{\tau}_{n+1\;*}}& \phantom{A}\simp_{\star}}}
\,\cong\;\ens_{\star}^{\text{\scriptsize{${\bf\Delta}^{op}$}}\,,}
\end{equation}
\begin{equation}\label{laa3}
\text{et}\qquad\ens^{\text{\scriptsize{$\Big({\bf\Delta}_{\leq n+1}\Big)^{op}$}}}\;=\,\vcenter{
\xymatrix@C+15pt{
\simp_{\leq n+1}\;\phantom{A}\ar@{}[r]|-{\perp\phantom{A}}
\ar@{<-}@/^18pt/[r]^-{\tau_{n+1}^{\phantom{A} *}}
\ar@/_18pt/[r]_-{\tau_{n+1\;*}}& \phantom{A}\simp}}
\,=\;\ens^{\text{\scriptsize{${\bf\Delta}^{op}$}}\,,}
\end{equation}
(voir le début de \S\ref{tron}) induisent par le foncteur d'inclusion canonique:
$$
\xymatrix@C+10pt{{\bf\Delta}_{\leq n+1}\;\ar@{^(->}[r]^-{\tau_{n+1}}&{\bf\Delta}}\,;
$$
et posons ${\bf csq}_{n+1}^{\text{\tiny{$\text{pointé}$}}} (X,x)=\overline{\tau}_{n+1\;*} \overline{\tau}_{n+1}^{\;*}(X,x)$, pour tout ensemble simplicial pointé $(X,x)$, de fa\c con analogue à la notation ${\bf csq}_{n+1} (X)=\tau_{n+1\;*} \tau_{n+1}^{\;*}(X)$, pour un ensemble simplicial $X$.

On peut trouver une description de ${\bf csq}_{n+1}^{\text{\tiny{$\text{pointé}$}}} (X,x)$ en fonction de l'ensemble simplicial ${\bf csq}_{n+1}(X)$, en considérant les adjonctions composées:
$$
\xymatrix@C+15pt{
\simp_{\leq n+1,\,\star}\;\phantom{A}\ar@{}[r]|-{\perp\phantom{A}}
\ar@{<-}@/^18pt/[r]^-{\overline{\tau}_{n+1}^{\phantom{A} *}}
\ar@/_18pt/[r]_-{\overline{\tau}_{n+1\;*}}& \phantom{A}
\simp_{\star}\phantom{a}
\ar@{}[r]|-{\perp}\ar@<-4pt>@/_18pt/[r]_-{\pi}&
\phantom{a}
\ar@<-4pt>@/_18pt/[l]_-{(\,\cdot\,)_{+}}
\simp}
$$
$$
\text{et}\qquad\xymatrix@C+15pt{
\simp_{\leq n+1,\,\star}
\;\phantom{A}
\ar@{}[r]|-{\perp}\ar@<-4pt>@/_18pt/[r]_-{\pi}&
\phantom{a}
\ar@<-4pt>@/_18pt/[l]_-{(\,\cdot\,)_{+}} 
\simp_{\leq n+1}
\;\phantom{A}
\ar@{}[r]|-{\perp\phantom{A}}
\ar@{<-}@/^18pt/[r]^-{\tau_{n+1}^{\phantom{A} *}}
\ar@/_18pt/[r]_-{\tau_{n+1\;*}}& 
\phantom{A}
\simp}
$$
qu'on construit à partir de \eqref{olvidar}, \eqref{adjonctionpoin} et \eqref{laa3}. 

En effet, vu que $\overline{\tau}_{n+1}^{\;*} (X_{+}) = \tau_{n+1}^*(X)_{+}$, pour tout ensemble simplicial $X$; il se suit que si $(A,a)$ est un ensemble simplicial pointé $(n+1)$-tronqué, l'ensemble simplicial sous-jacent à $\overline{\tau}_{n+1\;*}(A,a)$ est isomorphe à l'ensemble simplicial $\tau_{n+1\;*}(A)$. En particulier, il existe un isomorphisme naturel d'ensembles simpliciaux pointés:
\begin{equation}\label{legalite}
{\bf csq}_{n+1}^{\text{\tiny{$\text{pointé}$}}} (X,x)\,\cong\,\big({\bf csq}_{n+1}(X),\overline{x}\big)\,,
\end{equation}
où $\overline{x}$ est l'image de la flèche $\xymatrix@C-10pt{\star\ar[r]^-{x}&X}$ par le foncteur ${\bf csq}_{n+1}=\tau_{n+1\;*}\tau_{n+1}^{\;*}$; c'est-à-dire, $\overline{x}$ est le morphisme d'ensembles simpliciaux composé:
$$
\xymatrix{
\star \ar[r]^-{x} & 
X \ar[r]^-{\eta_{X}} & 
{\bf csq}_{n+1}(X)\,.}
$$ 

Finalement, remarquons que si $\overline{\eta}\colon \xymatrix@C-15pt{\mathrm{id}_{\simp_{\star}}\ar@{=>}[r] & {\bf csq}_{n+1}^{\text{\tiny{$\text{pointé}$}}}}$ est une unité de l'adjonction \eqref{adjonctionpoin}, la famille dont les éléments sont les morphismes images par le foncteur d'oubli $\xymatrix@C-5pt{\simp_\star \ar[r]^-{\pi} & \simp}$ des morphismes:
\begin{equation}
\xymatrix{(X,x)\ar[r]^-{\overline{\eta}_{X}} & {\bf csq}_{n+1}^{\text{\tiny{$\text{pointé}$}}}(X,x)\,,}\quad\text{pour $(X,x)$ un objet de $\simp_{*}$;}
\end{equation}
est une unité de l'adjonction \eqref{laa3}.

\begin{lemme}\label{lemmepoi1}
Soit $n\geq 0$ et $(X,x)$ un ensemble simplicial pointé. Si $(X,x)$ est un objet fibrant de la catégorie de modèles $(\simp_{\star},\pi^{-1}{\bf W}_{\infty}, {\bf mono},\pi^{-1}{\bf fib}_{\infty})$, alors ${\bf csq}^{\text{\tiny{$\text{pointé}$}}}_{n+1}(X,x)$ est un objet fibrant de la catégorie de modèles $(\simp_{\star},\pi^{-1}{\bf W}_{n}, {\bf mono},\pi^{-1}{\bf fib}_{n})$, et le morphisme canonique $\xymatrix@C-10pt{(X,x)\ar[r]&{\bf csq}^{\text{\tiny{$\text{pointé}$}}}_{n+1}(X,x)}$ appartient à $\pi^{-1}{\bf W}_{n}$.
\end{lemme}
\begin{proof} 
D'après le Corollaire \ref{pointeco2} un ensemble simplicial pointé $(X,x)$ est un objet fibrant de la catégorie de modèles $(\simp_{\star},\pi^{-1}{\bf W}_{\infty}, {\bf mono},\pi^{-1}{\bf fib}_{\infty})$, si et seulement si l'ensemble simplicial $X$ est un complexe de Kan; et ${\bf csq}^{\text{\tiny{$\text{pointé}$}}}_{n+1}(X)$ est un objet fibrant de la catégorie de modèles $(\simp_{\star},\pi^{-1}{\bf W}_{n}, {\bf mono},\pi^{-1}{\bf fib}_{n})$ si et seulement si ${\bf csq}_{n+1}(X)$ est un complexe de Kan tel que $\pi_{m}\big({\bf csq}_{n+1}(X)\big)=0$ si $m\geq n$. 

D'un autre dire que $\xymatrix@C-10pt{(X,x)\ar[r]&{\bf csq}^{\text{\tiny{$\text{pointé}$}}}_{n+1}(X,x)}$ appartienne à $\pi^{-1}{\bf W}_{n}$ équivaut à dire que $\xymatrix@C-10pt{X\ar[r]&{\bf csq}_{n+1}(X)}$ soit une $n$-équivalence faible. L'énoncé est alors une conséquence du Corollaire \ref{lecoro}.
\end{proof}

Montrons maintenant:

\begin{lemme}\label{lemmepoi2}
Soit $n\geq 0$. Si $X=(X,x)$ et $Y=(Y,y)$ sont des ensembles simpliciaux pointés, l'ensemble simplicial des morphismes $\underline{\mathrm{Hom}}_{\simp_{\star}}\big(X,{\bf csq}^{\text{\tiny{$\text{pointé}$}}}_{n+1}Y\big)$ est $(n+1)$-cosquelettique et le morphisme canonique $\xymatrix@C-10pt{X\ar[r]&{\bf csq}^{\text{\tiny{$\text{pointé}$}}}_{n+1}(X)}$ induit un isomorphisme d'ensembles simpliciaux:
\begin{equation}\label{lemcos3}
\xymatrix@C-8pt{
\underline{\mathrm{Hom}}_{\simp_{\star}}\big({\bf csq}^{\text{\tiny{$\text{pointé}$}}}_{n+1}X,{\bf csq}^{\text{\tiny{$\text{pointé}$}}}_{n+1}Y\big) \ar[r]^-{\cong} &
\underline{\mathrm{Hom}}_{\simp_{\star}}\big(X,{\bf csq}^{\text{\tiny{$\text{pointé}$}}}_{n+1}Y\big)}\,.
\end{equation}
\end{lemme}
\begin{proof}
La preuve est analogue à celle du Lemme \ref{lemehomo}. En effet vu que:
$$
\overline{\tau}_{n+1}^{\phantom{A}*}(A\wedge B)\cong \overline{\tau}_{n+1}^{\phantom{A}*}(B)\wedge \overline{\tau}_{n+1}^{\phantom{A}*}(B)
$$ 
parce que le foncteur $\overline{\tau}_{n+1}^{\phantom{A}*}$ commute aux limites et colimites, on considère les chaînes d'isomorphismes:
\begin{align*}
\mathrm{Hom}_{\simp}\Big( \Delta^{m+1},\underline{\mathrm{Hom}}_{\simp_\star}\big(X,{\bf csq}_{n+1}^{\text{\tiny{$\text{pointé}$}}}Y\big)\Big) \,& \cong\,
\mathrm{Hom}_{\simp_\star}\big( \Delta^{m+1}_+\wedge X,{\bf csq}_{n+1}^{\text{\tiny{$\text{pointé}$}}}Y) \\
\,& \cong \, 
\mathrm{Hom}_{\simp_{\leq n+1,\,\star}}\Big(\overline{\tau}_{n+1}^*\big(\Delta^{m+1}_+\wedge X\big),\overline{\tau}_{n+1}^*Y\Big) \\
\,& \cong \, 
\mathrm{Hom}_{\simp_{\leq n+1,\,\star}}\Big(\overline{\tau}_{n+1}^*\big(\Delta^{m+1}_+\big)\wedge \overline{\tau}_{n+1}^*\big(X\big),\overline{\tau}_{n+1}^*Y\Big) \\ 
\,& \cong \, 
\mathrm{Hom}_{\simp_{\leq n+1,\,\star}}\Big(\overline{\tau}_{n+1}^*\big(\partial\Delta^{m+1}_+\big)\wedge \overline{\tau}_{n+1}^*\big(X\big),\overline{\tau}_{n+1}^*Y\Big) \\
\,& \cong \, 
\mathrm{Hom}_{\simp_{\leq n+1,\,\star}}\Big(\overline{\tau}_{n+1}^*\big(\partial\Delta^{m+1}_+\wedge X\big),\overline{\tau}_{n+1}^*Y\Big) \\
\,& \cong \, 
\mathrm{Hom}_{\simp_\star}\big(\partial\Delta^{m+1}_+\wedge X,{\bf csq}_{n+1}^{\text{\tiny{$\text{pointé}$}}}Y\big) \\
\,& \cong \, 
\mathrm{Hom}_{\simp}\Big(\partial\Delta^{m+1},\underline{\mathrm{Hom}}_{\simp_\star}\big(X,{\bf csq}_{n+1}^{\text{\tiny{$\text{pointé}$}}}Y\big)\Big)
\end{align*}
où $m\geq n+1$ et:
\begin{align*}
\underline{\mathrm{Hom}}_{\simp_{\star}}\big(X,{\bf csq}^{\text{\tiny{$\text{pointé}$}}}_{n+1}Y\big)_{m} \,& = \,
\mathrm{Hom}_{\simp_{\star}}\big( X\wedge\Delta_{+}^{m+1},{\bf csq}^{\text{\tiny{$\text{pointé}$}}}_{n+1}Y) \\
\,& \cong \, \mathrm{Hom}_{\simp_{\leq n+1,\,\star}}\Big(\overline{\tau}_{n+1}^{\phantom{A}*}\big(X\wedge\Delta_{+}^{m}\big),\overline{\tau}_{n+1}^{\phantom{A}*}Y\Big) \\
\,& \cong \, \mathrm{Hom}_{\simp_{\leq n+1,\,\star}}\big(\overline{\tau}_{n+1}^{\phantom{A}*}(X)\wedge\overline{\tau}_{n+1}^{\phantom{A}*}(\Delta_{+}^{m}),\overline{\tau}_{n+1}^{\phantom{A}*}Y\Big) \\
\,& \cong \, \mathrm{Hom}_{\simp_{\leq n+1,\,\star}}\Big(\overline{\tau}_{n+1}^{\phantom{A}*} \big({\bf csq}^{\text{\tiny{$\text{pointé}$}}}_{n+1}(X)\big)\wedge\overline{\tau}_{n+1}^{\phantom{A}*}(\Delta_{+}^{m}),\overline{\tau}_{n+1}^{\phantom{A}*}Y\Big) \\
\,& \cong \, \mathrm{Hom}_{\simp_{\leq n+1,\,\star}}\Big(\overline{\tau}_{n+1}^{\phantom{A}*} \big({\bf csq}^{\text{\tiny{$\text{pointé}$}}}_{n+1}(X)\wedge\Delta_{+}^{m}\big),\overline{\tau}_{n+1}^{\phantom{A}*}Y\Big) \\
\,& \cong \, \mathrm{Hom}_{\simp_{\leq n+1,\,\star}}\Big({\bf csq}^{\text{\tiny{$\text{pointé}$}}}_{n+1}\big(X\big)\wedge\Delta_{+}^{m},{\bf csq}^{\text{\tiny{$\text{pointé}$}}}_{n+1}Y\Big) \\
\,& \cong \, \underline{\mathrm{Hom}}_{\simp_{\star}}\big({\bf csq}^{\text{\tiny{$\text{pointé}$}}}_{n+1}X,{\bf csq}^{\text{\tiny{$\text{pointé}$}}}_{n+1}Y\big)_{m}\,.
\end{align*}
\end{proof}

Finalement:

\begin{corollaire}\label{pointegroupdff}
Soit $n\geq 0$. Si $Y=(Y,y)$ est un ensemble simplicial pointés dont le ensemble simplicial sous-jacent est un $n$-groupoïde de Kan, alors l'ensemble simplicial des morphismes $\underline{\mathrm{Hom}}_{\simp_{\star}}\big(X,Y\big)$ est un $n$-groupoïde de Kan pour tout ensemble simplicial pointé $X=(X,x)$.
\end{corollaire}
\begin{proof}
L'ensemble simplicial $\underline{\mathrm{Hom}}_{\simp_{\star}}\big(X,Y\big)$ est $(n+1)$-cosquelettique d'après le Lemme \ref{lemmepoi2} et $\underline{\mathrm{Hom}}_{\simp_{\star}}\big(X,Y\big)$ est un complexe de Kan parce que l'ensemble simplicial sous-jacent à $Y$ est un complexe de Kan.

Il se suit du Lemme \ref{uneequ} qu'il faut montrer que la fonction:
$$
\xymatrix@C+10pt{
\mathrm{Hom}_{\simp_{\leq n+1,\,\star}}\big(\overline{\tau}_{n+1}^{\phantom{A}*}(X\wedge \Delta^{n+1}_+),\overline{\tau}_{n+1}^{\phantom{A}*}Y\Big)
\ar[r] & 
\mathrm{Hom}_{\simp_{\leq n+1,\,\star}}\big(\overline{\tau}_{n+1}^{\phantom{A}*}(X\wedge \Lambda^{n+1,k}_+),\overline{\tau}_{n+1}^{\phantom{A}*}Y\Big)}
$$
induit de l'inclusion canonique $\xymatrix@C-8pt{\Lambda^{n+1,k}\ar@{^(->}[r]&\Delta^{n+1}}$ est injective si $0\leq k\leq n+1$.

Soit $0\leq k\leq n+1$. Vu qu'on a les égalités des ensembles:
$$
\Delta^{n+1}_{n+1} \Big\backslash \Lambda^{n+1,k}_{n+1} \; =  \; 
\Big\{ \; \xymatrix@C-8pt{[n+1]\ar[r]^-{\mathrm{id}_{[n+1]}}&[n+1]}\,, \;\; \xymatrix@C-8pt{[n+1]\ar[r]^-{\sigma_i}&[n]\ar[r]^-{\delta_k}&[n+1]} \; \text{où \; $0\leq i\leq n$}\; \Big\} \,,
$$
$$
\Delta^{n+1}_{n} \Big\backslash \Lambda^{n+1,k}_{n} \; =  \; 
\Big\{ \; \xymatrix@C-8pt{[n]\ar[r]^-{\delta_k}&[n+1]} \; \; \Big\} \qquad \text{et} \qquad 
\Delta^{n+1}_{n-1} \Big\backslash \Lambda^{n+1,k}_{n-1} \; =  \;  \emptyset\,,
$$
si $F,G\colon\xymatrix@C-10pt{\tau_{n+1}^{\phantom{A}*}(X\times \Delta^{n+1})\ar[r] & \tau_{n+1}^{\phantom{A}*}(Y)}$ sont de morphismes d'ensembles simpliciaux tronqués tels que:
\begin{equation}\label{propriecuecuer3}
\text{$F_{k}=G_{k}$ \; si \; $0\leq k\leq n-1$, \quad  $F_{n}(a,f) = G_{n}(a,f)$ \; si \; $a\in X_{n}$ \; et \; $f\in \Lambda^{n+1,k}_{n}$\,,}
\end{equation}
$$
\text{$F_{n+1}(x,f) = G_{n+1}(x,\varphi)$ \; si \; $x\in X_{n+1}$ \: et \; $\varphi\in \Lambda^{n+1,k}_{n+1}$}
$$
$$
\text{et \qquad $F_i(x_0, \alpha) = y_0 = G_i(x_0, \alpha)$ \; si \; $\alpha \in \Delta_i^{n+1}$ \: et \; $0\leq i\leq n+1$\,,}
$$
on doit montrer que:
\begin{equation}
\begin{split}
\text{$F_{n}(a,\delta_k) = G_{n}(a,\delta_k)$ \; si \; $a\in X_{n}$}  \hspace{4cm}\\
\text{$F_{n+1}(x,\mathrm{id}_{[n+1]}) = G_{n+1}(x,\mathrm{id}_{[n+1]})$ \; et \; $F_{n+1}(x,\delta_k\sigma_i) = G_{n+1}(x,\delta_k\sigma_i)$ \; si \; $x\in X_{n+1}$}\,.
\end{split}
\end{equation}

La preuve est la même que cela du Lemme \ref{cartfermengrpKAN}.
\end{proof}

\renewcommand{\thesubsection}{\S\thesection.\arabic{subsection}}
\subsection{}
\renewcommand{\thesubsection}{\thesection.\arabic{subsection}}

Si $n\geq 0$ et ${\bf \Delta}_{\leq n}$ est la catégorie définie au début de \S\ref{tron}, considérons les catégories produit ${\bf\Delta}\times \big({\bf\Delta}_{\leq n}\big)$ et $ \big({\bf\Delta}_{\leq n}\big)\times{\bf\Delta}$. On appelle respectivement les objets des catégories de préfaisceaux:
$$
\ssimp_{\leq n} \,=\,  {\bf Ens}^{\text{\scriptsize{ $\Big({\bf\Delta}\times \big({\bf\Delta}_{\leq n}\big)\Big)^{op}$ }}}\qquad\text{et}\qquad
\widehat{\mathbf{{\bf\Delta}}_{\leq n}\times\mathbf{{\bf\Delta}}}
 \,=\,  {\bf Ens}^{\text{\scriptsize{ $\Big( \big({\bf\Delta}_{\leq n}\big)\times{\bf\Delta}\Big)^{op}$ }}}
$$ 
des ensembles bisimpliciaux \emph{verticalement $n$-tronqués} et \emph{horizontalement $n$-tronqués}.

Les foncteurs d'inclusion:
$$
\xymatrix@C+10pt{{\bf\Delta}\times\big({\bf\Delta}_{\leq n}\big)\;\ar@{^(->}[r]^-{\mathrm{id}\times\nu_{n}}&{\bf\Delta}\times{\bf\Delta}}
\qquad \text{et}\qquad
\xymatrix@C+10pt{\big({\bf\Delta}_{\leq n}\big)\times{\bf\Delta}\;\ar@{^(->}[r]^-{\nu_{n}\times\mathrm{id}}&{\bf\Delta}\times{\bf\Delta}}
$$ 
induisent d'adjonctions:
\begin{equation*}
\def\objectstyle{\scriptstyle}
\def\labelstyle{\scriptstyle}
\xymatrix@C+9pt{
\ssimp_{\leq n}\;\phantom{A}\ar@{}[r]|-{\perp\phantom{A}}
\ar@{<-}@/^15pt/[r]^-{(\mathrm{id}\times\nu_{n})^*}
\ar@/_15pt/[r]_-{(\mathrm{id}\times\nu_{n})_{*}}& \phantom{A}\ssimp }  \qquad\text{et}\qquad
\xymatrix@C+9pt{
\widehat{\mathbf{{\bf\Delta}}_{\leq n}\times\mathbf{{\bf\Delta}}}
\;\phantom{A}\ar@{}[r]|-{\perp\phantom{A}}
\ar@{<-}@/^15pt/[r]^-{(\nu_{n}\times\mathrm{id})^*}
\ar@/_15pt/[r]_-{(\nu_{n}\times\mathrm{id})_{*}}& \phantom{A}\ssimp }  \,,
\end{equation*}
respectivement.

Un ensemble bisimplicial $X$ est dit \emph{verticalement} $n$-\emph{cosque\-lettique}, si le morphisme canonique:
$$\xymatrix@C+4pt{X\ar[r]^-{\eta_{X}^n}&(\mathrm{id}\times\nu_{n})_{*}(\mathrm{id}\times\nu_{n})^* \, X }$$ 
est un isomorphisme; c'est-à-dire, si $X$ est isomorphe à un objet dans l'image du foncteur $(\mathrm{id}\times\nu_{n})_{*}$ ci-dessus. De fa\c con analogue, $X$ est dit \emph{horizontalement} $n$-\emph{cosque\-lettique} si $X$ est isomorphe à un objet dans l'image du foncteur $(\nu_{n}\times\mathrm{id})_{*}$.

\begin{lemme}\label{cosque2}
Pour tout ensemble bisimplicial $X$, les énoncés suivants sont équivalents:
\begin{enumerate}
\item $X$ est verticalement (resp. horizontalement) $n$-cosquelettique.  
\item $X$ est verticalement (resp. ho\-ri\-zon\-talement) $m$-cosquelettique pour $m\geq n$. 
\item La fonction:
$$
\def\objectstyle{\scriptstyle}
\def\labelstyle{\scriptstyle}
\xymatrix@C+20pt{
\mathrm{Hom}_{\ssimp}\Big(\Delta^p\boxtimes\Delta^{q},X\Big)
\ar[r]^-{\big(\mathrm{id}\boxtimes\alpha^{q-1}\big)^*} &
\mathrm{Hom}_{\ssimp}\Big(\Delta^p\boxtimes\partial\Delta^{q},X\Big)}
$$
$$
\bigg(\text{resp.}\quad\vcenter{
\def\objectstyle{\scriptstyle}
\def\labelstyle{\scriptstyle}
\xymatrix@C+20pt{
\mathrm{Hom}_{\ssimp}\Big(\Delta^p\boxtimes\Delta^{q},X\Big)
\ar[r]^-{\big(\alpha^{p-1}\boxtimes\mathrm{id}\big)^*} &
\mathrm{Hom}_{\ssimp}\Big(\partial\Delta^p\boxtimes\Delta^{q},X\Big)}}\bigg),
$$
induite du monomorphisme $\xymatrix@-10pt{\partial\Delta^{q} \ar[r]^-{\alpha^{q-1}} & \Delta^{q}}$ $\Big(\text{resp.}\;\, \xymatrix@-10pt{\partial\Delta^{p} \ar[r]^-{\alpha^{p-1}} & \Delta^{p}} \Big)$, est bijective pour tout $p\geq 0$ et $q\geq n+1$ (resp. $p\geq n+1$ et $q\geq 0$).
\item L'ensemble simplicial $X_{p,\bullet}$ (resp. $X_{\bullet,q}$) est $n$-cosquelettique pour tout $p\geq 0$ (resp. $q\geq 0$).
\item Les morphismes horizontaux (resp. verticaux) du carré suivant:
\begin{equation}\label{carcarte3}
\def\objectstyle{\scriptstyle}
\def\labelstyle{\scriptstyle}
\xymatrix@C+18pt{
\mathrm{Hom}_{\ssimp}\big(\Delta^{p}\boxtimes\,\Delta^{q},X\big) 
\ar[r]^-{(\mathrm{id}\boxtimes\alpha^{q-1})^*}\ar[d]_-{(\alpha^{p-1}\boxtimes\mathrm{id})^*}&
\mathrm{Hom}_{\ssimp}\big(\Delta^{p}\boxtimes\,\partial\Delta^{q},X\big)
\ar[d]^-{(\alpha^{p-1}\boxtimes\mathrm{id})^*}\,\\
\mathrm{Hom}_{\ssimp}\big(\partial\Delta^{p}\boxtimes\,\Delta^{q},X\big)
\ar[r]_-{(\mathrm{id}\boxtimes\alpha^{q-1})^*}&
\mathrm{Hom}_{\ssimp}\big(\partial\Delta^{p}\boxtimes\,\partial\Delta^{q},X\big)\,,}
\end{equation}
sont des isomorphismes pour tout $p\geq 0$ et $q\geq n+1$ (resp. $p\geq n+1$ et $q\geq 0$).
\item Pour tout $p\geq 0$ et $q\geq n+1$ (resp. $p\geq n+1$ et $q\geq 0$):
\begin{equation}\label{carcarte}
\def\objectstyle{\scriptstyle}
\def\labelstyle{\scriptstyle}
\xymatrix@C+18pt{
\mathrm{Hom}_{\ssimp}\big(\Delta^{p}\boxtimes\,\Delta^{q},X\big) 
\ar[r]^-{(\mathrm{id}\boxtimes\alpha^{q-1})^*}\ar[d]_-{(\alpha^{p-1}\boxtimes\mathrm{id})^*}&
\mathrm{Hom}_{\ssimp}\big(\Delta^{p}\boxtimes\,\partial\Delta^{q},X\big)
\ar[d]^-{(\alpha^{p-1}\boxtimes\mathrm{id})^*}\,\\
\mathrm{Hom}_{\ssimp}\big(\partial\Delta^{p}\boxtimes\,\Delta^{q},X\big)
\ar[r]_-{(\mathrm{id}\boxtimes\alpha^{q-1})^*}&
\mathrm{Hom}_{\ssimp}\big(\partial\Delta^{p}\boxtimes\,\partial\Delta^{q},X\big)\,,}
\end{equation}
est un carré cartésien.
\item Les morphismes horizontaux (resp. verticaux) du carré suivant: 
\begin{equation}\label{carcarte4}
\def\objectstyle{\scriptstyle}
\def\labelstyle{\scriptstyle}
\xymatrix@C+18pt{
\mathrm{Hom}_{\ssimp}\big(\Delta^{p}\boxtimes\,\Delta^{q},X\big) 
\ar[r]^-{(\mathrm{id}\boxtimes\alpha^{q-1})^*}\ar[d]_-{(\alpha^{p-1,k}\boxtimes\mathrm{id})^*}&
\mathrm{Hom}_{\ssimp}\big(\Delta^{p}\boxtimes\,\partial\Delta^{q},X\big)
\ar[d]^-{(\alpha^{p-1,k}\boxtimes\mathrm{id})^*}\,\\
\mathrm{Hom}_{\ssimp}\big(\Lambda^{p,k}\boxtimes\,\Delta^{q},X\big)
\ar[r]_-{(\mathrm{id}\boxtimes\alpha^{q-1})^*}&
\mathrm{Hom}_{\ssimp}\big(\Lambda^{p,k}\boxtimes\,\partial\Delta^{q},X\big)}
\end{equation}
$$
\left(\text{resp.}\qquad\quad
\def\objectstyle{\scriptstyle}
\def\labelstyle{\scriptstyle}
\vcenter{\xymatrix@C+18pt{
\mathrm{Hom}_{\ssimp}\big(\Delta^{p}\boxtimes\,\Delta^{q},X\big) 
\ar[r]^-{(\mathrm{id}\boxtimes\alpha^{q-1,k})^*}\ar[d]_-{(\alpha^{p-1}\boxtimes\mathrm{id})^*}&
\mathrm{Hom}_{\ssimp}\big(\Delta^{p}\boxtimes\,\Lambda^{q,k},X\big)
\ar[d]^-{(\alpha^{p-1}\boxtimes\mathrm{id})^*}\,\\
\mathrm{Hom}_{\ssimp}\big(\partial\Delta^{p}\boxtimes\,\Delta^{q},X\big)
\ar[r]_-{(\mathrm{id}\boxtimes\alpha^{q-1,k})^*}&
\mathrm{Hom}_{\ssimp}\big(\partial\Delta^{p}\boxtimes\,\Lambda^{q,k},X\big)}}\right)
$$
sont des isomorphismes si $p\geq 0$, $0\leq k\leq p$ et $q\geq n+1$ (resp. $p\geq n+1$, $0\leq k\leq p$ et $q\geq 0$).
\item Pour tout $p\geq 0$, $0\leq k\leq p$ et $q\geq n+1$ (resp. $p\geq n+1$, $0\leq k\leq p$ et $q\geq 0$):
\begin{equation}\label{carcarte2}
\def\objectstyle{\scriptstyle}
\def\labelstyle{\scriptstyle}
\xymatrix@C+18pt{
\mathrm{Hom}_{\ssimp}\big(\Delta^{p}\boxtimes\,\Delta^{q},X\big) 
\ar[r]^-{(\mathrm{id}\boxtimes\alpha^{q-1})^*}\ar[d]_-{(\alpha^{p-1,k}\boxtimes\mathrm{id})^*}&
\mathrm{Hom}_{\ssimp}\big(\Delta^{p}\boxtimes\,\partial\Delta^{q},X\big)
\ar[d]^-{(\alpha^{p-1,k}\boxtimes\mathrm{id})^*}\,\\
\mathrm{Hom}_{\ssimp}\big(\Lambda^{p,k}\boxtimes\,\Delta^{q},X\big)
\ar[r]_-{(\mathrm{id}\boxtimes\alpha^{q-1})^*}&
\mathrm{Hom}_{\ssimp}\big(\Lambda^{p,k}\boxtimes\,\partial\Delta^{q},X\big)}
\end{equation}
$$
\left(\text{resp.}\qquad\quad
\def\objectstyle{\scriptstyle}
\def\labelstyle{\scriptstyle}
\vcenter{\xymatrix@C+18pt{
\mathrm{Hom}_{\ssimp}\big(\Delta^{p}\boxtimes\,\Delta^{q},X\big) 
\ar[r]^-{(\mathrm{id}\boxtimes\alpha^{q-1,k})^*}\ar[d]_-{(\alpha^{p-1}\boxtimes\mathrm{id})^*}&
\mathrm{Hom}_{\ssimp}\big(\Delta^{p}\boxtimes\,\Lambda^{q,k},X\big)
\ar[d]^-{(\alpha^{p-1}\boxtimes\mathrm{id})^*}\,\\
\mathrm{Hom}_{\ssimp}\big(\partial\Delta^{p}\boxtimes\,\Delta^{q},X\big)
\ar[r]_-{(\mathrm{id}\boxtimes\alpha^{q-1,k})^*}&
\mathrm{Hom}_{\ssimp}\big(\partial\Delta^{p}\boxtimes\,\Lambda^{q,k},X\big)}}\right)
$$
est un carré cartésien.
\end{enumerate}
\end{lemme}
\begin{proof}
On montre l'équivalence des énoncés (i) et (ii)  comme dans la preuve du Lemme \ref{cosque}. D'un autre, d'après le Lemme \ref{cosque} les énoncés (iii) et (iv) sont équivalents parce qu'on a des isomorphismes:
$$
\mathrm{Hom}_{\ssimp}\big(\Delta^p\boxtimes W,X\big) \, \cong \, \mathrm{Hom}_{\simp}\big(W,X_{p,\bullet}\big)
\qquad\text{et}\qquad
\mathrm{Hom}_{\ssimp}\big(W\boxtimes \Delta^q,X\big) \, \cong \, \mathrm{Hom}_{\simp}\big(W,X_{\bullet,q}\big)
$$
pour tout ensemble simplicial $W$.

Pour montrer l'équivalence des énoncés (i) et (iii) dans le cas vertical considérons l'adjonction:
\begin{equation}\label{adjj}
\ssimp_{\leq m}
\xymatrix@C+18pt{
\phantom{a}\ar@{}[r]|-{\perp}
\ar@{<-}@/^18pt/[r]^-{(\mathrm{id}\times j_{m})^*}
\ar@/_18pt/[r]_-{(\mathrm{id}\times j_{m})_{*}}&
\phantom{a}}\ssimp_{\leq m+1}
\end{equation}
associée au foncteur d'inclusion $\xymatrix@C-4pt{{\bf\Delta}_{\leq m}\;\ar@{^(->}[r]^-{j_{m}}&{\bf\Delta}_{\leq m+1}}$.

On observe alors que si $A$ est un ensemble bisimplicial verticalement $m$-tronqué, l'ensemble bisimplicial verticalement $(m+1)$-tronqué $(\mathrm{id}\times j_{m})_{*} (A)$ est à isomorphisme près décrite comme suit: Pour tout $p\geq 0$, 
\begin{align*}
(\mathrm{id}\times j_{m})_{*} (A)_{p,q}\; & = \; A_{p,q} \qquad\text{si\; $0\leq q\leq m$,}\\
\text{et}\quad
(\mathrm{id}\times j_{m})_{*} (A)_{p,m+1}\; & = \; 
\mathrm{Hom}_{\ssimp_{\leq m}}
\Big( (\mathrm{id}\times\nu_{m})^*\big( \Delta^p \boxtimes\, \partial\Delta^{m+1} \big) ,A\Big)\\ 
 & \cong \;
\Bigg\{\text{\scriptsize{$\big(a_{0},\dots,a_{m+1}\big) \;\in\;\underset{0}{\overset{m+1}{\prod}} A_{p,m}$}}\;\Bigg|\;
\vcenter{\xymatrix@R=1pt{\text{\scriptsize{$d^{m-1}_{i}a_{j}=d^{m-1}_{j-1}a_{i}$}}\\ \text{\scriptsize{si $0\leq i<j\leq m+1$.}}}}\Bigg\}.
\end{align*}

Les morphismes faces et dégénérescences verticaux sont ceux de $A$ pour $p\geq 0$ et $0\leq q\leq m$; et:
 $$
\xymatrix{
(\mathrm{id}\times j_{m})_{*} (A)_{p,m+1} \ar@<-6pt>[d]_-{d_{i}^v}\\
A_{p,m}\ar@<-6pt>[u]_{s_{i}^v}}
$$
sont définis comme dans la preuve du Lemme \ref{cosque}. Enfin, les morphismes horizontaux:
$$
\xymatrix{
(\mathrm{id}\times j_{m})_{*} (A)_{p-1,m+1} \ar@<-6pt>[r]_-{s_{i}^h}&
(\mathrm{id}\times j_{m})_{*} (A)_{p,m+1} \ar@<-6pt>[l]_-{d_{i}^h}}
$$
sont définis argument par argument:
$$
\xymatrix@R=1pt{
\big(d_{i}^h a_{0},\dots,d_{i}^h a_{m+1}\big) &\big(a_{0},\dots,a_{m+1}\big) \ar@{|->}[l]_-{d_{i}^h}\\
\big( a_{0},\dots, a_{m+1}\big) \ar@{|->}[r]_-{s_{i}^h}&\big(s^h_{i}a_{0},\dots,s^h_{i}a_{m+1}\big)\,.}
$$

On peut construire une unité de l'adjonction \eqref{adjj}: Cela est définie dans un ensemble bisimplicial $(m+1)$-tronqué $B$ par un morphisme $\xymatrix@C-4pt{B\ar[r]&(\mathrm{id}\times j_{m})_{*}(\mathrm{id}\times j_{m})^* \, B}$, dont la partie non trivial s'identifie aux fonctions:
\begin{equation}\label{nont}
\def\objectstyle{\scriptstyle}
\def\labelstyle{\scriptstyle}
\xymatrix@C-5pt{
\mathrm{Hom}_{\ssimp_{\leq m+1}}\big((\mathrm{id}\times \nu_{m+1})^{*}
\big(\Delta^p\boxtimes\,\Delta^{m+1}\big),\,B\big)\ar[d]_-{(\mathrm{id}\boxtimes\alpha^m)^*} \,\cong\, B_{p,m+1}\\
\mathrm{Hom}_{\ssimp_{\leq m+1}}\big((\mathrm{id}\times \nu_{m+1})^{*} 
\big(\Delta^p\boxtimes\,\partial\Delta^{m+1}\big),\,B\big)\,\cong\, (\mathrm{id}\times j_{m})_{*}(\mathrm{id}\times j_{m})^* (B)_{p,m+1} ,}
\end{equation}
pour $p\geq 0$. Donc, $B$ est dans l'image du foncteur $(\mathrm{id}\times j_{m})_{*}$ si et seulement si, \eqref{nont} est une bijection pour tout $p\geq 0$.

On en conclut que $X$ est verticalement $n$-cosquelettique, si et seulement si la fonction:
$$
\def\objectstyle{\scriptstyle}
\def\labelstyle{\scriptstyle}
\xymatrix@C+25pt@R=6pt{
\mathrm{Hom}_{\ssimp}\Big(\Delta^p\boxtimes\,\Delta^{m+1},\,X\Big)
\ar[r]^-{(\mathrm{id}\boxtimes\,\alpha^m)^*} &
\mathrm{Hom}_{\ssimp}\Big(\Delta^p\boxtimes\,\partial\Delta^{m+1},\,X\Big),}
$$
est bijective pour tout $p\geq 0$ et $m\geq n$. Donc, les énoncés (i) et (iii) sont équivalents.

Les énoncés (iii) et (v) sont équivalents parce que la fonction:
$$
\def\objectstyle{\scriptstyle}
\def\labelstyle{\scriptstyle}
\xymatrix@C+20pt{
\mathrm{Hom}_{\ssimp}\Big(\Delta^p\boxtimes\Delta^{q},X\Big)
\ar[r]^-{\big(\mathrm{id}\boxtimes\alpha^{q-1}\big)^*} &
\mathrm{Hom}_{\ssimp}\Big(\Delta^p\boxtimes\partial\Delta^{q},X\Big)}
$$
est bijective pour tout $p\geq 0$, si et seulement si la fonction:
$$
\def\objectstyle{\scriptstyle}
\def\labelstyle{\scriptstyle}
\xymatrix@C+20pt{
\mathrm{Hom}_{\ssimp}\Big(W\boxtimes\Delta^{q},X\Big)
\ar[r]^-{\big(\mathrm{id}\boxtimes\alpha^{q-1}\big)^*} &
\mathrm{Hom}_{\ssimp}\Big(W\boxtimes\partial\Delta^{q},X\Big)}
$$
est bijective pour tout ensemble simplicial $W$.

Montrons que les énoncés (iii) et (vi) sont équivalents par un argument inductif sur $p\geq 0$. 

Si $p=0$, remarquons que dans le carré:
\begin{equation} \label{carcarte0}
\def\objectstyle{\scriptstyle}
\def\labelstyle{\scriptstyle}
\xymatrix@C+18pt{
\mathrm{Hom}_{\ssimp}\big(\Delta^{0}\boxtimes\,\Delta^{q},X\big) 
\ar[r]^-{(\mathrm{id}\boxtimes\alpha^{q-1})^*}\ar[d]_-{(\alpha^{-1}\boxtimes\mathrm{id})^*}&
\mathrm{Hom}_{\ssimp}\big(\Delta^{0}\boxtimes\,\partial\Delta^{q},X\big)
\ar[d]^-{(\alpha^{-1}\boxtimes\mathrm{id})^*}\,\\
\mathrm{Hom}_{\ssimp}\big(\partial\Delta^{0}\boxtimes\,\Delta^{q},X\big)
\ar[r]_-{(\mathrm{id}\boxtimes\alpha^{q-1})^*}&
\mathrm{Hom}_{\ssimp}\big(\partial\Delta^{0}\boxtimes\,\partial\Delta^{q},X\big)\,,}
\end{equation}
la fonction horizontale inférieur est l'identité de l'ensemble à un élément. Donc, puisqu'on a un carré commutatif:
$$
\def\objectstyle{\scriptstyle}
\def\labelstyle{\scriptstyle}
\xymatrix@C+25pt@R=6pt{
\mathrm{Hom}_{\simp}\big(\Delta^{q},\,X_{0,\bullet}\big)\ar[r]^-{\alpha^m_{X_{0,\bullet}}}
\ar@{}[d]|-{\mathrel{\reflectbox{\rotatebox[origin=c]{90}{$\cong$}}}}&
\mathrm{Hom}_{\simp}\big(\partial\Delta^{q},\,X_{0,\bullet}\big)\ar@{}[d]|-{\mathrel{\reflectbox{\rotatebox[origin=c]{90}{$\cong$}}}},\\
\mathrm{Hom}_{\ssimp}\Big(\Delta^0\boxtimes\,\Delta^q,\,X\Big)
\ar[r]_{(\mathrm{id}\boxtimes\,\alpha^{q-1})^*} &
\mathrm{Hom}_{\ssimp}\Big(\Delta^0\boxtimes\,\partial\Delta^{q},\,X\Big)\,,}
$$  
on en déduit que le carré \eqref{carcarte0} est cartésien pour $q\geq n+1$, si et seulement si $X_{0,\bullet}$ est un ensemble simplicial $n$-cosquelettique.

Supposons maintenant que les conditions (iii) et (vi) sont équivalentes pour $0\leq p\leq k$. Si $X$ est un ensemble bisimplicial tel que $X_{k,\bullet}$ est $n$-cosquelettique, on déduit que dans le carré:
\begin{equation} \label{carcartek}
\def\objectstyle{\scriptstyle}
\def\labelstyle{\scriptstyle}
\xymatrix@C+18pt{
\mathrm{Hom}_{\ssimp}\big(\Delta^{k+1}\boxtimes\,\Delta^{q},X\big) 
\ar[r]^-{(\mathrm{id}\boxtimes\alpha^{q-1})^*}\ar[d]_-{(\alpha^{k}\boxtimes\mathrm{id})^*}&
\mathrm{Hom}_{\ssimp}\big(\Delta^{k+1}\boxtimes\,\partial\Delta^{q},X\big)
\ar[d]^-{(\alpha^{k}\boxtimes\mathrm{id})^*}\,\\
\mathrm{Hom}_{\ssimp}\big(\partial\Delta^{k+1}\boxtimes\,\Delta^{q},X\big)
\ar[r]_-{(\mathrm{id}\boxtimes\alpha^{q-1})^*}&
\mathrm{Hom}_{\ssimp}\big(\partial\Delta^{k+1}\boxtimes\,\partial\Delta^{q},X\big)\,,}
\end{equation}
la fonction horizontale inférieur est un isomorphisme pour $q\geq n+1$, du fait qu'on peut l'identifier avec une colimite de fonctions du type:
$$
\xymatrix@C+18pt{
\mathrm{Hom}_{\simp}\big(\Delta^{q},X_{k,\bullet}\big)
\ar[r]^-{\alpha^{q-1}_{X}}&
\mathrm{Hom}_{\ssimp}\big(\partial\Delta^{q},X_{k,\bullet}\big)\,.}
$$

Donc, \eqref{carcartek} est cartésien pour $q\geq n+1$, si et seulement si $X_{k+1,\bullet}$ est $n$-cosquelettique. Autrement dit les conditions (iii) et (vi) sont équivalentes pour $0\leq p\leq k+1$.

On montre de façon analogue que les conditions (iii), (vii) et (viii) sont équivalentes.
\end{proof}

\renewcommand{\thesubsection}{\S\thesection.\arabic{subsection}}
\subsubsection{}\;
\renewcommand{\thesubsection}{\thesection.\arabic{subsection}} \label{menorn}

Considérons ${\bf\Delta}\underset{\leq n}{\times}{\bf\Delta}$, la sous-catégorie pleine de la catégorie produit ${\bf\Delta}\times{\bf\Delta}$ dont les objets sont les couples $\big([p],[q]\big)$ tels que $p+q\leq n$; et soit:
\begin{equation}\label{diagon}
\xymatrix@C+18pt{
\widehat{\mathbf{{\bf\Delta}}\underset{\leq n}{\times}\mathbf{{\bf\Delta}}}\phantom{A}\ar@{}[r]|-{\perp\phantom{A}}
\ar@{<-}@/^20pt/[r]^-{\mu_{n}^{\; *}}
\ar@/_20pt/[r]_-{\mu_{n\, *}}& \phantom{A}\underset{\phantom{\leq 2}}{\ssimp} }  \,,
\end{equation}
l'adjonction induite du foncteur d'inclusion canonique:
\begin{equation}\label{mun}
\xymatrix@C+10pt{{\bf\Delta}\underset{\leq n}{\times}{\bf\Delta}\;\ar@{^(->}[r]^-{\mu_{n}}&{\bf\Delta}\times{\bf\Delta}}.
\end{equation}

\begin{lemme}\label{derr}
Un ensemble bisimplicial $X$ est isomorphe à un objet dans l'image du foncteur $\mu_{n\, *}$, si et seulement si le carré suivant:
\begin{equation}\label{carcarteB}
\def\objectstyle{\scriptstyle}
\def\labelstyle{\scriptstyle}
\xymatrix@C+18pt{
\mathrm{Hom}_{\ssimp}\big(\Delta^{p}\boxtimes\,\Delta^{q},X\big) 
\ar[r]^-{(\mathrm{id}\boxtimes\alpha^{q-1})^*}\ar[d]_-{(\alpha^{p-1}\boxtimes\mathrm{id})^*}&
\mathrm{Hom}_{\ssimp}\big(\Delta^{p}\boxtimes\,\partial\Delta^{q},X\big)
\ar[d]^-{(\alpha^{p-1}\boxtimes\mathrm{id})^*}\,\\
\mathrm{Hom}_{\ssimp}\big(\partial\Delta^{p}\boxtimes\,\Delta^{q},X\big)
\ar[r]_-{(\mathrm{id}\boxtimes\alpha^{q-1})^*}&
\mathrm{Hom}_{\ssimp}\big(\partial\Delta^{p}\boxtimes\,\partial\Delta^{q},X\big)\,,}
\end{equation}
est cartésien pour tout $p,q\geq 0$ tels que $p+q\geq n+1$. 
\end{lemme}
\begin{proof}
Pour commencer rappelons que si $X$ est un ensemble bisimplicial arbitraire, le carré \eqref{carcarteB} induit simplement une fonction:
\begin{equation}\label{pfibre}
\xymatrix@R-8pt{
\mathrm{Hom}_{\ssimp}\big(\Delta^{p}\boxtimes\,\Delta^{q},X\big)  \ar[d]\\
\underset{\mathrm{Hom}_{\ssimp}\big(\partial\Delta^{p}\boxtimes\,\partial\Delta^{q},X\big)}{\mathrm{Hom}_{\ssimp}\big(\partial\Delta^{p}\boxtimes\,\Delta^{q},X\big)\;\,\bigtimes\;\,\mathrm{Hom}_{\ssimp}\big(\Delta^{p}\boxtimes\,\partial\Delta^{q},X\big)}\,,}
\end{equation}
vers le produit fibré des flèches:
$$
\def\objectstyle{\scriptstyle}
\def\labelstyle{\scriptstyle}
\xymatrix@C+18pt{
&\mathrm{Hom}_{\ssimp}\big(\Delta^{p}\boxtimes\,\partial\Delta^{q},X\big)
\ar[d]^-{(\alpha^{p-1}\boxtimes\mathrm{id})^*}\,\\
\mathrm{Hom}_{\ssimp}\big(\partial\Delta^{p}\boxtimes\,\Delta^{q},X\big)
\ar[r]_-{(\mathrm{id}\boxtimes\alpha^{q-1})^*}&
\mathrm{Hom}_{\ssimp}\big(\partial\Delta^{p}\boxtimes\,\partial\Delta^{q},X\big)\,.}
$$

Si par ailleurs on remarque que l'ensemble de départ de la fonction \eqref{pfibre} est isomorphe à $X_{p,q}$, ainsi que son ensemble d'arrivé est isomorphe à l'ensemble:
\begin{equation}\label{fifi}
\left\{\vcenter{\xymatrix@R=1pt{
\text{\scriptsize{$\big(a_{0},\dots,a_{p}; b_{0},\dots,b_{q}\big)$}} \\
\text{\scriptsize{$\in\;\,\underset{0}{\overset{p}{\prod}} X_{p-1,q}\;\times\;  \underset{0}{\overset{q}{\prod}} X_{p,q-1}$}}}}
\;\left| \;\vcenter{\xymatrix@R=1pt{
\text{\scriptsize{$d^{h}_{i}a_{i'}=d^{h}_{i'-1}a_{i}$}}, \;\, \text{\scriptsize{$d^{v}_{j}b_{j'}=d^{v}_{j'-1}b_{j}$}} \\
\text{\scriptsize{et}}\;\; \text{\scriptsize{$d^{v}_{j}a_{i}=d^{h}_{i}b_{j}$,}}\\
\text{\scriptsize{si\;\,$0\leq i<i' \leq p$\;\, et\;\, $0\leq j<j' \leq q$.}}}}\right.\right\}\footnote{Par définition $X_{p,-1}=X_{-1,q}=\star$ est l'ensemble à un élément.};
\end{equation}
on constat sans difficulté que la fonction \eqref{pfibre} s'identifie à l'association:
\begin{equation}\label{pfibre2}
\xymatrix@R=1pt@C+5pt{
X_{p,q}\ar[r]&\text{L'ensemble}\;\,\eqref{fifi}\\
a\ar@{}[r]|-{\mapsto} &\big(\,d^v_{0}a,\dots,d^v_{p}a\,; \,d^h_{0}a,\dots,d^h_{q}a\,\big)\,.}
\end{equation}

Donc, le carré \eqref{carcarteB} est cartésien si et seulement si la fonction \eqref{pfibre2} est bijective.

Pour montrer maintenant le Lemme, remarquons que si on décompose le foncteur d'inclusion canonique $\xymatrix@C+5pt{{\bf\Delta}\underset{\leq n}{\times}{\bf\Delta}\;\ar@{^(->}[r]^-{\mu_{n}}&{\bf\Delta}\times{\bf\Delta}}$ comme $\xymatrix@C+5pt{{\bf\Delta}\underset{\leq n}{\times}{\bf\Delta}\;\ar@{^(->}[r]^-{\rho_{n}}&{\bf\Delta}_{\leq n}\times{\bf\Delta}_{\leq n}\;\ar@{^(->}[r]^-{\nu_{n}\times\nu_{n}}&{\bf\Delta}\times{\bf\Delta}}$; on vérifie qu'un ensemble bisimplicial $X$ est isomorphe à un objet dans l'image de $\mu_{n\,*}\cong (\nu_{n}\times\nu_{n})_{*}\,\rho_{n\, *}$, si et seulement si $X$ vérifie aux propriétés:
\begin{enumerate}
\item[{\bf (A)}] $X$ est isomorphe à un objet dans l'image du foncteur $(\nu_{n}\times\nu_{n})_{*}$.
\item[{\bf (B)}] L'ensemble bisimplicial tronqué $(\nu_{n}\times\nu_{n})^*X$ est isomorphe à un objet dans l'image du foncteur $\rho_{n\, *}$. 
\end{enumerate}

En effet, cela provenant du fait que la counité de l'adjonction $(\nu_{n}\times\nu_{n})_{*} \dashv (\nu_{n}\times\nu_{n})^{*}$ est un isomorphisme. 

D'un autre côté, on découle du Lemme \ref{cosque2} qu'un ensemble bisimplicial $X$ satisfait la condition {\bf (A)}, si et seulement si le carré \eqref{carcarteB} est cartésien pour tout $p,q\geq 0$ tel que $\mathrm{max}(p,q)\geq n+1$. Il suffit alors de vérifier que $X$ satisfait la condition {\bf (B)}, si et seulement si la fonction \eqref{pfibre2} est bijective pour tous $0\leq p,q\leq n$ tels que $p+q \geq n$. 

Considérons dans ce but la décomposition:
$$
{\bf\Delta}\underset{\leq n}{\times}{\bf\Delta}\,=
\xymatrix@C-5pt{\A_{0}\ar@{^(->}[r]^-{F_{1}} & \A_{1}\ar@{^(->}[r]^-{F_{2}} &\cdots\ar@{^(->}[r]^-{F_{n}} &\A_{n}}
=\,{\bf\Delta}_{\leq n}\times{\bf\Delta}_{\leq n},
$$
de l'inclusion $\xymatrix{{\bf\Delta}\underset{\leq n}{\times}{\bf\Delta}\;\ar@{^(->}[r]^-{\rho_{n}}&{\bf\Delta}_{\leq n}\times{\bf\Delta}_{\leq n}}$, où $\A_{i}$ note la sous-catégorie pleine de ${\bf\Delta}_{\leq n}\times{\bf\Delta}_{\leq n}$, dont les objets sont les couples $\big([p],[q]\big)$ tels que $p+q\leq n+i$.

Pour chaque $1\leq i\leq n$, on peut décrire une unité de l'adjonction $F_{i\;*}  \dashv F_{i}^{\phantom{a} *}$ de la fa\c con suivante: 
Si $A$ est une préfaisceau d'ensembles sur $\A_{i}$, on vérifie qu'à isomorphisme près:
$$
F_{i\, *} F_{i}^{\;*}(A)_{p,q} \;=\,\;
\left\{\vcenter{\xymatrix@R=5pt{A_{p,q}  & \text{si\; \;$p+q < n+i$,} \\
\eqref{fifi} &\text{si\;\; $p+q=n+i$.}}}\right.
$$

Les morphismes faces et dégénérescences pour $p+q=n+i$:
$$
\xymatrix@C+10pt{A_{p-1,q}\ar@<-6pt>[r]_-{s_{i}^h} & 
F_{i\, *} F_{i}^{\;*}(A)_{p,q}\,, \ar@<-6pt>[d]_-{d_{i}^v}\ar@<-6pt>[l]_-{d_{i}^h}\\
&A_{p,q-1}\ar@<-6pt>[u]_-{s_{i}^v}}
$$
sont définis par les règles:
$$
\xymatrix@R=1pt@C-14pt{
a_{i} && \big(a_{0},\dots,a_{p}; b_{0},\dots,b_{q}\big) \ar@{|->}[dddd]^-{d_{i}^v}\ar@{|->}[ll]_-{d_{i}^h}&\text{et}&
a\ar@{|->}[rr]^-{s_{i}^h} && \big(a_{0},\dots,a_{p}; s_{i}^h d_{0}^v a,\dots,s_{i}^h d^v_{q} a\big) \\
&&&&&& \big(s_{i}^v d_{0}^h b,\dots,s_{i}^v d^h_{p} b; b_{0},\dots,b_{q}\big) \\
&&&&&&\\&&&&&&\\&&b_{i}&&&&\\
&&&&&&b\ar@{|->}[uuuu]_{s_{i}^v}\\&&&&&&}
$$
$$
\text{où}\quad 
a_{k}=\begin{cases}
s_{i-1}^hd_{k}^h a & \text{si\;$0\leq k\leq i-1$}\\
a & \text{si\;$i\leq k\leq i+1$}\\
s_{i}^hd_{k-1}^h a& \text{si\;$i+2\leq k\leq p$}
\end{cases}
\quad\text{et}\quad
b_{k}=\begin{cases}
s_{i-1}^vd_{k}^v b & \text{si\;$0\leq k\leq i-1$}\\
b & \text{si\;$i\leq k\leq i+1$}\\
s_{i}^vd_{k-1}^vb & \text{si\;$i+2\leq k\leq q$.}
\end{cases}
$$

Finalement, on définie le morphisme unité $\xymatrix@C-5pt{A\ar[r]&F_{i\, *} F_{i}^{\;*}(A)}$, comme la fonction identité si $p+q<n+i$ et comme la fonction \eqref{pfibre2} si $p+q=n+i$. En particulier, $A$ est isomorphe à un objet dans l'image du foncteur $F_{i}$, si et seulement si la fonction \eqref{pfibre2} est bijective, pour tous les couples $(p,q)$ tels que $p,q\geq 0$ et $p+q=n+i$.
 
Donc, $X$ satisfait la condition {\bf (B)}, si et seulement si la fonction \eqref{pfibre2} est bijective pour tout $0\leq p,q\leq n$ tel que $p+q \geq n$. 
\end{proof}


\chapter{Le $2$-groupe d'homotopie des espaces connexes}  
\setcounter{section}{11}

\section{Le nerf des petites catégories}\label{nerfptcat}

\renewcommand{\thesubsection}{\S\thesection.\arabic{subsection}}
\subsection{}\;
\renewcommand{\thesubsection}{\thesection.\arabic{subsection}}

Rappelons que le foncteur d'inclusion canonique $\xymatrix@C-10pt{\Delta\ar[r]&{\bf cat}}$ induit une adjonction:
\begin{equation}\label{nerfcategories}
\xymatrix@C+15pt{\mathbf{cat}\ar@/_12pt/[r]_{\mathrm{N}(\;\cdot\;)}\ar@{}[r]|-{\perp}&\simp \ar@/_12pt/[l]_{\mathrm{P}(\;\cdot\;)} \,,}
\end{equation}
où $\mathrm{N}(\,\cdot\,)$ est défini par la formule $\mathrm{N}(\A)_m=\mathrm{Hom}_{\bf cat} \big( [m], \A\big)$ et le foncteur $\mathrm{P}(\;\cdot\;)$ est une extension de Kan à gauche du foncteur $\xymatrix@C-10pt{\Delta\ar[r]&{\bf cat}}$ le long du plongement de Yoneda:
\begin{equation*}
\xymatrix@C+5pt@R=1pt{{\bf\Delta} \ar@{^(->}[r]&  \simp \\
[n] \ar@{}[r]|-{\longmapsto} & \Delta^n \,.}
\end{equation*}

On appelle $\mathrm{N}(\mathcal{A})$ l'\emph{nerf} (\emph{géométrique}) de la catégorie $\mathcal{A}$ et $\mathrm{P}(X)$ \emph{une catégorie des chemins} de l'ensemble simplicial $X$.

Si $\mathcal{A}$ est une petite catégorie l'ensemble simplicial tronqué $\tau_{2}^{\,*}\big(\mathrm{N}(\mathcal{A})\big)$ admet la description qui suit: Les ensembles $\mathrm{N}(\mathcal{A})_0$ et $\mathrm{N}(\mathcal{A})_1$ s'identifient aux ensembles des objets et des morphismes de $\mathcal{A}$ respectivement. L'image d'un objet $x$ de $\mathcal{A}$ par la fonction:
$$
\xymatrix@C-5pt{
\mathrm{N}(\mathcal{A})_0 \ar[r]^-{s_0} & \mathrm{N}(\mathcal{A})_1
}
$$
est le morphisme identité de $x$ dans $\mathcal{A}$. Les images d'un morphismes $f$ de $\mathcal{A}$ par les fonctions:
$$
\xymatrix@C-5pt{
\mathrm{N}(\mathcal{A})_1 \ar[r]^-{d_0} & \mathrm{N}(\mathcal{A})_0
}\quad\text{et}\quad
\xymatrix@C-5pt{
\mathrm{N}(\mathcal{A})_1 \ar[r]^-{d_1} & \mathrm{N}(\mathcal{A})_0
}
$$ 
sont le but et la source de $f$ respectivement. Finalement se donner un élément $\eta$ de l'ensemble $\mathrm{N}(\mathcal{A})_2$ équivaut à se donner trois morphismes de $\mathcal{A}$ dans un triangle commutatif:
$$
\vcenter{\xymatrix@R-10pt@C+3pt{
&A_1\ar[rd]^-{d_0\eta}&\\A_0\ar[ru]^-{d_2\eta} \ar[rr]_-{d_1\eta}& & A_2 \,.}}
$$

Montrons:

\begin{lemme}\label{nerfgrp2cos}
Le nerf $\mathrm{N}(\A)$ d'une petite catégorie $\A$ est un ensemble simplicial faiblement $1$-cos\-que\-le\-ttique (voir \ref{faiblecos}).
\end{lemme}
\begin{proof}
L'ensemble simplicial $\mathrm{N}(\A)$ est $2$-cosquelettique si la fonction: 
$$
\xymatrix{ \mathrm{Hom}_{\bf cat}\big([m+1],\A\big) \ar[r]^-{\underset{s}{\prod} \; \delta^*_{s}} & \underset{0\leq s\leq m+1}{\prod}\mathrm{Hom}_{\bf cat}\big([m],\A\big)}
$$ 
est le noyau dans ${\bf Ens}$ de la double flèche:
$$
\xymatrix@C-1pt{
\underset{0\leq s\leq q+1}{\prod}\mathrm{Hom}_{\bf cat}\big([m],\A\big)
\ar@<+7pt>[rrr]^-{\underset{i<j}{\prod} \;\delta^*_{j-1}\circ\mathrm{proj}_i}
\ar@<-7pt>[rrr]_-{\underset{i<j}{\prod} \;\delta^*_{i}\circ\mathrm{proj}_j}&&&
\underset{0\leq i<j\leq q+1}{\prod}\mathrm{Hom}_{\bf cat}\big([m-1],\A\big)}\,,
$$
pour toute petite catégorie $\A$ et tout entier $m\geq 2$.

Montrons l'énoncé équivalent:
\begin{quote}
Si $q\geq 2$ le morphisme $\xymatrix{\underset{0\leq s\leq q+1}{\bigsqcup}[q]\ar[r]^-{\underset{s}{\sqcup} \; \delta_{s}}& [q+1]}$ est le conoyau dans ${\bf cat}$ de la double flèche:
$$
\xymatrix@C-1pt{
\underset{0\leq i<j\leq q+1}{\bigsqcup}[q-1]
\ar@<+7pt>[rrr]^-{\underset{i<j}{\sqcup} \; \text{inc}_i\,\circ\,\delta_{j-1}}
\ar@<-7pt>[rrr]_-{\underset{i<j}{\sqcup} \;\text{inc}_{j}\,\circ\,\delta_{i}}&&&
\underset{0\leq s\leq q+1}{\bigsqcup}[q]}\,.
$$
\end{quote}

Commençons pour vérifier que pour $q\geq 1$ on a un carré cartésien dans ${\bf cat}$:
$$
\xymatrix{[q] \ar[r]^-{\delta_0} & [q+1] \\ [q-1] \ar[u]^-{\delta_q} \ar[r]_-{\delta_0}& [q]\ar[u]_-{\delta_{q+1}} }
$$

Pour cela remarquons que se donner deux foncteurs $\xymatrix@C-10pt{[q]\ar[r]^-{F}&\C}$ et $\xymatrix@C-10pt{[q]\ar[r]^-{G}&\C}$ tels que $F\circ\delta_q=G\circ\delta_0$ équivaut à se donner deux suites de morphismes dans $\C$:
$$
F \, = \, \Big(  \xymatrix@C-5pt{A_0 \ar[r]^-{f_1} & \cdots\cdots \ar[r]^-{f_{q}} &A_{q}}\Big)\qquad \text{et} \qquad 
G \, = \,\Big(  \xymatrix@C-5pt{B_0 \ar[r]^-{g_1} & \cdots\cdots \ar[r]^-{g_{q}} &B_{q}}\Big)
$$
telles que:
$$
\Big(  \xymatrix@C-5pt{A_0 \ar[r]^-{f_1} & \cdots\cdots \ar[r]^-{f_{q-1}} &A_{q-1}}\Big) \, = \, F\circ\delta_q \, = \, G\circ\delta_0 \, = \, \Big(  \xymatrix@C-5pt{B_1 \ar[r]^-{g_2} & \cdots\cdots \ar[r]^-{g_{q}} &B_{q}}\Big).
$$

Le seul foncteur $\xymatrix@C-10pt{[q+1]\ar[r]^-{H}&\C}$ tel que $H\circ\delta_0=F$ et $H\circ\delta_{q+1}=G$ est déterminé par la suite de morphismes:
$$
H \, = \,\left(  \vcenter{\xymatrix@C-5pt@R=2pt{
& A_0 \ar@{=}[d]\ar[r]^-{f_1} & \cdots\cdots \ar[r]^-{f_{q-1}} &A_{q-1}\ar@{=}[d]\ar[r]^-{f_q} & A_q\\
B_0 \ar[r]_-{g_1} &B_1 \ar[r]_-{g_2}& \cdots\cdots \ar[r]_-{g_{q}} &B_{q} & }}\right)\,.
$$

Montrons maintenant que si $\C$ est une petite catégorie et $\big\{ \xymatrix@C-10pt{[q]\ar[r]^{F_s}&\C }\big\}_{0\leq s\leq q+1}$ est une famille de foncteurs vérifiant $F_j\circ\delta_i=F_i\circ\delta_{j-1}$ si  $0\leq i < j \leq q+1$; il existe un seul foncteur $F:\xymatrix@C-8pt{[q+1]\ar[r]&\C }$ tel que $F\circ\delta_s=F_s$ pour $0\leq s \leq q+1$. En effet, vu que par hypothèse $F_{q+1}\circ\delta_0=F_0\circ\delta_q$, on sait qu'il existe un seul foncteur $\xymatrix@C-10pt{[q+1]\ar[r]^-{F}&\C}$ tel que $F\circ\delta_0=F_0$ et $F\circ\delta_{q+1}=F_{q+1}$. 

D'un autre côté on a que pour $0\leq j\leq q$: 
$$F\circ\delta_j\circ\delta_q = F\circ\delta_{q+1}\circ\delta_j = F_{q+1}\circ\delta_j = F_j\circ\delta_q$$ 
$$\text{et}$$
$$F\circ\delta_j\circ\delta_0 = F\circ\delta_{0}\circ\delta_{j-1} = F_{0}\circ\delta_{j-1} = F_j\circ\delta_0\,;$$
vu que $q\geq 2$ on déduit que $F\circ\delta_j=F_j$.

Donc $\mathrm{N}(\A)$ est un ensemble simplicial $2$-cosquelettique. Finalement montrons que la fonction:
\begin{equation} \label{NGfaibcoscat}
\xymatrix@C+10pt{
\mathrm{Hom}_{\simp}\Big(\Delta^{2},\mathcal{N}\big(\G\big)\Big)
\ar[r]&\mathrm{Hom}_{\simp}\Big(\partial\Delta^{2},\mathcal{N}\big(\G\big)\Big)}
\end{equation}
induite du morphisme d'inclusion $\xymatrix@C-10pt{\partial\Delta^2\ar@{^(->}[r]&\Delta^2}$ est injective, c'est-à-dire que $\mathrm{N}(\A)$ est un ensemble simplicial faiblement $1$-cosquelettique.

Pour cela notons que l'ensemble de départ de la fonction \eqref{NGfaibcoscat} s'identifie à l'ensemble des triangles commutatifs de $\G$:
\begin{equation}\label{trinerfcat}
\vcenter{\xymatrix@R-10pt@C+3pt{
&A_1\ar[rd]^-{f_0}&\\A_0\ar[ru]^-{f_2} \ar[rr]_-{f_1}& & A_2  }}
\end{equation}
et son ensemble d'arrivée s'identifie à l'ensemble des diagrammes de la même forme mais pas nécessairement commutatif. Vu que la fonction \eqref{NGfaibcoscat} oublie la commutativité du diagramme, elle est injective. 
\end{proof}

On déduit:

\begin{corollaire}\label{nerfplnmfidd} 
Le foncteur nerf $\mathrm{N}\colon\xymatrix@C-12pt{{\bf cat}\ar[r]&\simp}$ est pleinement fidèle.
\end{corollaire}
\begin{proof}
Le foncteur $\mathrm{N}$ est fidèle parce que si on se donne un foncteur $F\colon\xymatrix@C-8pt{\mathcal{A}\ar[r]&\mathcal{B}}$ entre petites catégories, les fonctions $\mathrm{N}(F)_0$ et $\mathrm{N}(F)_1$ sont les fonction qui définissent $F$ dans les objets et les morphismes respectivement. Donc si $\mathrm{N}(F)_0 = \mathrm{N}(G)_0$ et $\mathrm{N}(F)_1 = \mathrm{N}(G)_1$ on déduit que $F=G$.

D'un autre si $\mathcal{A}$ et $\mathcal{B}$ sont des petites catégories et $\varphi\colon\xymatrix@C-8pt{\mathrm{N}(\mathcal{A})\ar[r]&\mathrm{N}(\mathcal{B})}$ est un morphisme d'ensembles simpliciaux, on vérifie sans peine les fonctions $\varphi_0$ et $\varphi_1$ définissent un foncteur $F\colon\xymatrix@C-8pt{\mathcal{A}\ar[r]&\mathcal{B}}$ tel que $\mathrm{N}(F)_0=\varphi_0$ et $\mathrm{N}(F)_1=\varphi_1$. 

Il se suit que $\mathrm{N}(F)=\varphi$ vu que $\mathrm{N}(\mathcal{B})$ est un ensemble simplicial faiblement $1$-cosquelettique. Donc le foncteur $\mathrm{N}$ est plein.
\end{proof}

Si $\mathcal{A}$ et $\mathcal{B}$ sont deux petites catégories, notons $\mathcal{B}^{\mathcal{A}}$ la catégorie des foncteurs de $\mathcal{A}$ vers $\mathcal{B}$ et ses transformations naturelles. On vérifie sans peine qu'il y a des bijections naturelles:
$$
\mathrm{Hom}_{\bf cat}\big(\mathcal{A}\times \mathcal{C},\mathcal{B}\big) \, \cong \, \mathrm{Hom}_{\bf cat}\Big(\mathcal{C},\mathcal{B}^{\mathcal{A}}\Big)\,,
$$
pour des petites catégories $\mathcal{A}$, $\mathcal{B}$ et $\mathcal{C}$ quelconques. Autrement dit la catégorie ${\bf cat}$ est une catégorie cartésienne fermée.

Vu que le foncteur nerf $\mathrm{N}\colon\xymatrix@C-12pt{{\bf cat}\ar[r]&\simp}$ est pleinement fidèle et commute aux produits finis on déduit des bijections naturelles:
\begin{align*}
\mathrm{N}\Big(\mathcal{B}^{\mathcal{A}}\Big)_n \, & = \, 
\mathrm{Hom}_{\bf cat}\Big([n],\mathcal{B}^{\mathcal{A}}\Big) \, \cong \, \mathrm{Hom}_{\bf cat}\big(\mathcal{A}\times [n],\mathcal{B}\big)\\
&\cong \, \mathrm{Hom}_{\simp}\big(\mathrm{N}(\mathcal{A})\times \Delta^n,\mathrm{N}(\mathcal{B})\big) \, = \, \underline{\mathrm{Hom}}_{\simp}\big(\mathrm{N}(\mathcal{A}),\mathrm{N}(\mathcal{B})\big)_n\,,
\end{align*}
donc un isomorphisme naturel de foncteurs:
\begin{equation}\label{isonhomhom}
\vcenter{\xymatrix@C-10pt{
\mathrm{N}\Big((\cdot_{2})^{(\cdot_{1})}\Big)\ar@{=>}[r] & 
\underline{\mathrm{Hom}}_{\simp}\big(\mathrm{N}(\cdot_1),\mathrm{N}(\cdot_2)\big)}}\colon \vcenter{\xymatrix{{\bf cat}^{op}\times {\bf cat} \ar[r]& \simp}}\,.
\end{equation}

\renewcommand{\thesubsection}{\S\thesection.\arabic{subsection}}
\subsubsection{}\;\label{groupdnnerfe}
\renewcommand{\thesubsection}{\thesection.\arabic{subsection}}

Posons ${\bf Grpd}$ pour noter la sous-catégorie pleine de ${\bf cat}$ dont les objets sont les grou\-poïdes. Remarquons que ${\bf Grpd}$ est une sous-catégorie de ${\bf cat}$ stable par produits finis et objet de morphismes $\mathcal{H}^{\mathcal{G}}$, c'est-à-dire ${\bf Grpd}$ est une sous-catégorie cartésienne fermée de ${\bf cat}$. 

Montrons (voir les Corollaires \ref{nerfplnmfidd} et \ref{nerf1groupd}) que le foncteur nerf $\mathrm{N}\colon\xymatrix@C-12pt{{\bf cat}\ar[r]&\simp}$ induit par restriction une équivalence entre la catégorie cartésienne fermée ${\bf Grpd}$ et la catégorie cartésienne fermée ${\bf GrpdK}^1$ des $1$-groupoïdes de Kan (voir le Lemme \ref{cartfermengrpKAN}).

Rappelons par ailleurs: 

\begin{corollaire}\label{1grpKAN}
Si $\G$ est un groupoïde (\emph{i.e.} une petite catégorie dont toutes les flèches sont des isomorphismes) l'ensemble simplicial faiblement $1$-cos\-que\-le\-ttique $\mathrm{N}\big(\G\big)$ (voir le Lemme \ref{nerfgrp2cos}) satisfait la condition d'extension de Kan en dimension $1\leq m\leq 2$ et satisfait la condition de minimalité en dimension $1$, autrement dit $\mathrm{N}\big(\mathcal{G}\big)$ est un $1$-groupoïde de Kan (voir le Corollaire \ref{groupidequi}).

En particulier, si $\G$ est un groupoïde l'ensemble simplicial $\mathrm{N}\big(\mathcal{G}\big)$ est un objet fibrant de la catégorie de modèles $(\simp,{\bf W}_1, {\bf mono},{\bf fib}_1)$ du Théorème \ref{ntypess} (voir le Corollaire \ref{lecoro}); et si $G$ est un groupe (vu comme un groupoïde avec un seul objet) l'ensemble simplicial $\mathrm{N}\big(G\big)$ est un objet fibrant de la catégorie de modèles $(\simp_0,{\bf W}^{red}_1, {\bf mono},{\bf fib}^{red}_1)$ de la Proposition \ref{modred}.
\end{corollaire}
\begin{proof}
L'ensemble simplicial $\mathrm{N}\big(\G\big)$ est faiblement $1$-cosquelettique d'après la Lemme \ref{nerfgrp2cos}. Pour montrer que $\mathrm{N}(\G)$ satisfait en dimension $1$ la condition d'extension de Kan et la condition de minimalité, on simplement note que si on se donne dans le groupoïde $\G$ deux des trois morphismes d'un diagramme:
$$
\vcenter{\xymatrix@R-10pt@C+3pt{
&A_1\ar[rd]^-{f_0}&\\A_0\ar[ru]^-{f_2} \ar[rr]_-{f_1}& & A_2 \,, }}
$$
alors il existe un seul morphisme qui le complète dans un triangle commutatif. 

Vérifions que l'ensemble simplicial $\mathrm{N}(\G)$ satisfait la condition d'extension de Kan en dimension $2$, c'est-à-dire que la fonction:
\begin{equation} \label{NGfaibcoscat2}
\xymatrix@C+10pt{
\mathrm{Hom}_{\simp}\Big(\Delta^{3},\mathcal{N}\big(\G\big)\Big)
\ar[r]&\mathrm{Hom}_{\simp}\Big(\Lambda^{3,k},\mathcal{N}\big(\G\big)\Big)}
\end{equation}
induite du morphisme d'inclusion $\xymatrix@C-10pt{\Lambda^{3,k}\ar@{^(->}[r]&\Delta^3}$ est surjective pour tout $0\leq k\leq 3$. 

En effet, l'ensemble de départe de la fonction \eqref{NGfaibcoscat2} s'identifie à l'ensemble des diagrammes de $\G$ formés par quatre objets et six morphismes comme suit: 
\begin{equation}\label{thetracat}
\xymatrix@-12pt{
& X_1 \ar@/_12pt/[ldddd]\ar@/^15pt/[rrrddd]&&& \\ &&&&\\ 
& &&& \\
& X_0\ar[uuu] \ar[ld]\ar[rrr]&&&X_3\\
X_2\ar@/_12pt/[rrrru]&&&&&}
\end{equation}
tels que les quatre faces de \eqref{thetracat} soient de triangles commutatifs.

D'un autre l'ensemble d'arrivée de \eqref{NGfaibcoscat2} est égale à l'ensemble des diagrammes de la même forme, mais vérifiant seulement que trois de ses faces soient de triangles commutatifs. On note sans peine que cela suffit pour que la quatrième face soit aussi un triangle commutatif parce que les flèches sont des isomorphismes. Donc la fonction \eqref{NGfaibcoscat2} est bijective.
\end{proof}

Si $\mathcal{G}$ est un groupoïde rappelons que l'ensemble $\pi_0(\mathcal{G})$ des \emph{composantes connexes par arcs de $\mathcal{G}$} est par définition l'ensemble des objets à isomorphisme près du groupoïde $\mathcal{G}$. Si $a$ est un objet fixé de $\mathcal{G}$ le \emph{premier groupe d'homotopie de $\mathcal{G}$ basé en $a$} noté $\pi_1(\mathcal{G},a)$ est par définition $\mathrm{Hom}_{\mathcal{G}}(a,a)$, le groupe des automorphismes de $a$ dans $\mathcal{G}$. 

On définit sans peine des foncteurs:
\begin{equation}
\xymatrix@C+25pt{{\bf Grpd}\ar[r]^-{\pi_{0}(\,\cdot\,)} & {\bf Ens}}\quad \text{et}\quad
\xymatrix@C+25pt{{\bf Grpd}_{\star}\ar[r]^-{\pi_{1}(\,\cdot\,,\,\cdot\,)} & {\bf Grp}}\,,
\end{equation}
où ${\bf Grpd}_{\star}$ note la catégorie des groupoïdes pointés.

\begin{lemme}
Ils existent des isomorphismes canoniques de foncteurs:
\begin{equation}\label{isocanogroupd}
\xymatrix@C-8pt{\pi_0\ar@{=>}[r]^-{\alpha^0} & \pi_0\circ \mathrm{N}}\colon \xymatrix@C+5pt{{\bf Grpd} \ar[r] & {\bf Ens}}
\qquad \text{et} \qquad 
\xymatrix@C-8pt{\pi_1\ar@{=>}[r]^-{\alpha^1} & \pi_1\circ \mathrm{N}_{\star}}\colon \xymatrix@C+5pt{{\bf Grpd}_{\star} \ar[r] & {\bf Grp}}
\end{equation}
où $\pi_i$ d'un ensemble simplicial (pointé) $X$ est par définition le $i$-ème groupe d'homotopie de la réalisation géomé\-trique de $X$ (voir \ref{nequival}).
\end{lemme}
\begin{proof}
Si $\mathcal{G}$ est un groupoïde les ensembles $\pi_0(\mathcal{G})$ et $\pi_0\big(\mathrm{N}(\mathcal{G})\big)$ sont égaux à l'ensemble des objets de $\mathcal{G}$ soumis à la relation d'équivalence qu'identifie deux objets s'il existe un morphisme de $\mathcal{G}$ entre eux. On définit alors $\xymatrix@C-8pt{\pi_0\ar@{=>}[r]^-{\alpha^0} & \pi_0\circ \mathrm{N}}$ comme la transformation naturelle identité.

D'un autre d'après le Corollaire \ref{1grpKAN} et les Propositions \ref{kankan} et \ref{kankan2} on peut penser que le foncteur $\pi_1\circ \mathrm{N}_{\star}\colon \xymatrix@C+5pt{{\bf Grpd}_{\star} \ar[r] & {\bf Grp}}$ est le foncteur défini dans un groupoïde $\G$ muni d'un objet distingué $a$ par l'ensemble des morphismes $\mathrm{Hom}_{\mathcal{G}}(a,a)$ dont le produit est défini par la règle:
$$
f\,\star\, g \; = \; g\, \circ \, f \qquad \text{toujours \, que } \qquad f,g \, \in \, \mathrm{Hom}_{\mathcal{G}}(a,a)\,.
$$

On définit l'isomorphisme naturel $\alpha^1\colon\xymatrix@C-8pt{\pi_1\ar@{=>}[r] & \pi_1\circ \mathrm{N}_{\star}}$ dans un groupoïde $\G$ muni d'un objet distingué $a$ par la fonction:
$$
\vcenter{\xymatrix@C+13pt@R=1pt{
\pi_1\big(\G,a\big) \, = \, \mathrm{Hom}_{\mathcal{G}}(a,a) \ar[r]^-{\alpha^1_{(\G,a)}} & \mathrm{Hom}_{\mathcal{G}}(a,a)  \, = \, \pi_1\circ \mathrm{N}_{\star} \big(\G,a\big)\\
\phantom{espaciooo}f \ar@{}[r]|-{\longmapsto} & f^{-1}\phantom{espacioooespacio}
}}\,.
$$

Note que $\alpha^1_{(\G,a)}$ est bien un morphisme de groupes parce que:
$$
\alpha^1_{(\G,a)}\big(\,f\circ g\,\big) \, = \, \big(f\circ g\big)^{-1} \, = \, g^{-1}\circ f^{-1} \, = \, \alpha^1_{(\G,a)}\big(g\big) \circ \alpha^1_{(\G,a)}\big(f\big) \, = \, 
\alpha^1_{(\G,a)}\big(f\big) \star \alpha^1_{(\G,a)}\big(g\big) \,.
$$
\end{proof}

Rappelons que la \emph{catégorie homotopique des groupoïdes} $h{\bf Grpd}$ est la catégorie dont les objets sont les groupoïdes et les morphismes sont les classes à isomorphisme naturel près des foncteurs entre eux: 
$$
\mathrm{Hom}_{h{\bf Grpd}}\big( \G , \H\big) \, = \, \pi_0\big(\mathcal{H}^{\mathcal{G}}\big)\,.
$$ 

Note que $h{\bf Grpd}$ est une catégorie cartésienne fermée dont l'objet de morphismes est aussi le groupoïde des foncteurs $\G^{\H}$. De plus le foncteur canonique $\xymatrix@C-10pt{{\bf Grpd} \ar[r]&h{\bf Grpd}}$ respect les produits finis.  

Un foncteur  $F\colon \xymatrix@C-8pt{\mathcal{G} \ar[r] & \mathcal{H}}$ entre groupoïdes est dit une \emph{équivalence faible} s'il vérifie une des conditions équivalentes de l'énoncé qui suit:

\begin{corollaire}
Si $F\colon \xymatrix@C-8pt{\mathcal{G} \ar[r] & \mathcal{H}}$ est un foncteur entre groupoïdes, ils sont équivalents:
\begin{enumerate}
\item $\pi_0(F)$ et $\pi_1(F)$ sont de fonctions bijectives.
\item $F$ est une équivalence de catégories \emph{i.e.} $F$ est un foncteur pleinement fidèle et essentiellement surjectif. 
\item L'image de $F$ par le foncteur canonique $\xymatrix@C-10pt{{\bf Grpd} \ar[r]&h{\bf Grpd}}$ est un isomorphisme de $h{\bf Grpd}$.
\item $\mathrm{N}(F)$ est une $1$-équivalence homotopique faible d'ensembles simpliciaux.
\item $\mathrm{N}(F)$ est une $\infty$-équivalence homotopique faible d'ensembles simpliciaux.
\end{enumerate}
\end{corollaire}
\begin{proof}
On montre (i)$\Leftrightarrow$(ii)$\Leftrightarrow$(iii) sans difficulté.

D'un autre vu que les ensembles simpliciaux $\mathrm{N}(\mathcal{G})$ et $\mathrm{N}(\mathcal{H})$ sont des complexes de Kan dont les groupes d'homotopie $\pi_i\big(\mathrm{N}(\mathcal{G})\big)$ et $\pi_i\big(\mathrm{N}(\mathcal{H})\big)$ sont nuls pour $i\geq 2$ il se suit (iv)$\Leftrightarrow$(v).

Finalement (i)$\Leftrightarrow$(iv) d'après les isomorphismes \eqref{isocanogroupd}.  
\end{proof}

Montrons l'énoncé réciproque du Corollaire \ref{1grpKAN} (voir \cite{duskin}): 

\begin{corollaire}\label{nerf1groupd}
Si $W$ est un ensemble simplicial lequel est un $1$-groupoïde de Kan, alors il existe un groupoïde $\G$ et un isomorphisme d'ensembles simpliciaux $\mathrm{N}(\G)\,\cong\,W$. Autrement dit, un ensemble simplicial est un $1$-groupoïde de Kan si et seulement si il est isomorphe au nerf d'un groupoïde (voir le Corollaire \ref{1grpKAN}). 
\end{corollaire}
\begin{proof}
Soit $W$ un $1$-groupoïde de Kan. On définit le groupoïde $\G$ muni de l'isomorphisme désiré $\mathrm{N}(\G)\,\cong\,W$ comme suit:  L'ensemble des objets (resp. des morphismes) de $\G$ est égale à l'ensemble des $0$-simplexes (resp. des $1$-simplexes) de $W$. Si $f$ est un morphisme de $\G$, la source (resp. le but) de $f$ est l'objet $d_1(f)$ (resp. $d_0(f)$). En particulier l'ensemble $\mathrm{Hom}_{\simp}\big( \Lambda^{2,1}, W\big)$ s'identifie à l'ensemble des \emph{couples des morphismes composables} de $\G$ \emph{i.e.} des couples $(f,g)$ de $1$-simplexes de $W$ telles que $d_0(f)=d_1(g)$.  

On définit une loi de composition de $\G$ par le composé:
$$
\xymatrix@C+10pt{
\mathrm{Hom}_{\simp}\big(\Lambda^{2,1}, W\big) \ar[r]^-{(\alpha^{1,1}_W)^{-1}} & 
\mathrm{Hom}_{\simp}\big(\Delta^2,W\big) \ar[r]^-{\delta^*_1} &
\mathrm{Hom}_{\simp}\big(\Delta^1,W\big)\,.}
$$

Autrement dit, si $f,g$ sont des morphismes composables de $\G$ alors $g\circ f$ est égale au morphisme $d_1(\eta)$, où $\eta$ est le seul $2$-simplexe de $W$ tel que $d_2(\eta)=f$ et $d_0(\eta)=g$. 

Montrons d'abord que $\G$ ainsi défini est effectivement une catégorie: Si $f$ est un morphisme de $\G$ de source $d_1(f)$ et but $d_0(f)$, on a dans $\G$ que:
$$
f\circ \big(s_0(d_1f)\big) \, = \, f \qquad \text{et} \qquad \big(s_0(d_0f)\big) \circ f \, = \,  f\,;
$$
en effet $s_0(f)$ et $s_1(f)$ sont de $2$-simplexes de $W$ tels que:
$$
d_0\big(s_0(f)\big)  \, = \, f\, , \qquad d_1\big(s_0(f)\big)  \, = \, f  \qquad \text{et} \qquad d_2\big(s_0(f)\big)  \, = \, s_0\big(d_1(f)\big)\,,
$$ 
$$
d_0\big(s_1(f)\big)  \, = \, s_0\big(d_0(f)\big)\, , \qquad d_1\big(s_1(f)\big) \, = \, f \qquad \text{et} \qquad d_2\big(s_1(f)\big)  = f \,.
$$

Autrement dit, si $a$ est un objet de $\G$ le morphisme $s_0(a)$ est l'identité de $a$ dans $\G$.

Soient maintenant $f$, $g$ et $h$ trois morphismes de $\G$ composables:
$$
\xymatrix{a\ar[r]^-{f} & b \ar[r]^-{g} & c \ar[r]^-{h} & d\,;}
$$
c'est-à-dire, $f,g,h$ sont des $1$-simplexes de $W$ tels que $d_0(f) = d_1(g)$ et $d_0(g) = d_1(h)$.  

Pour montrer que $h\circ (g\circ f) = (h\circ g)\circ f$ considérons les seuls $2$-simplexes $\eta_0$, $\eta_2$ et $\eta_3$ de $W$ vérifiant les propriétés.
\begin{align*}
d_0(\eta_0) \, = & \, h\,     &d_2(\eta_0) \, &= \, g \\
d_0(\eta_2) \, =  & \, h\circ g \, &d_2(\eta_2) \, &= \, f\\
d_0(\eta_3) \, =  & \, g\,   &d_2(\eta_3) \, &= \, f\,;
\end{align*}
cela étant possible parce que $d_1(h) = d_0(g)$ et $d_1(h\circ g)=d_1(g)=d_0(f)$.

Vu que la fonction induite de l'inclusion $\xymatrix@C-8pt{\Lambda^{3,1}\ar[r] & \Delta^{3}}$:
$$
\xymatrix@C+10pt{
\mathrm{Hom}_{\simp}\big(\Delta^3,W\big) \ar[r] &
\mathrm{Hom}_{\simp}\big(\Lambda^{3,1},W\big)}
$$
est bijective et parce qu'on a les égalités: 
$$
d_0(\eta_2) \, = \, h\circ g \, = \, d_1(\eta_0)\,, \qquad d_0(\eta_3) \, = \,  g \, = \, d_2(\eta_0) \qquad \text{et} \qquad d_2(\eta_3) \, = \,  f \, = \, d_2(\eta_2)\,;
$$ 
il existe un seul $3$-simplexe $\xi$ de $W$ tel que $d_0(\xi)=\eta_0$, $d_2(\xi)=\eta_2$ et $d_3(\xi)=\eta_3$. En particulier:
\begin{align*}
h\circ (g\circ f) \, & = \, (d_0\eta_0)\circ (d_1\eta_3) \, = \, (d_0d_0\xi)\circ (d_1d_3\xi) \\ \, & = \, (d_0d_1\xi)\circ (d_2d_1\xi)  \, = \, d_1d_1\xi \, = \, d_1d_2\xi
\, = \, d_1\eta_2 \, = \, (h\circ g)\circ f\,.
\end{align*}

Donc $\G$ est bien une catégorie. 

Remarquons d'un autre côté que si $0\leq k\leq 2$, la fonction induite de l'inclusion $\xymatrix@C-8pt{\Lambda^{2,k}\ar[r] & \Delta^{2}}$:
$$
\xymatrix@C+10pt{
\mathrm{Hom}_{\simp}\big(\Delta^2,W\big) \ar[r] &
\mathrm{Hom}_{\simp}\big(\Lambda^{2,k},W\big)}
$$
est bijective par hypothèse. En particulier, si $f$ est un morphisme quelconque de $\G$ il existent de $2$-simplexes $\eta_1$ et $\eta_2$ de $W$ tels que:
\begin{align*}
d_1(\eta_1) \, = \, s_0(d_1f)\,, \qquad \; \qquad d_2(\eta_1) \, = \, f\,,\\
d_0(\eta_2) \, = \, f \qquad \text{et} \qquad d_1(\eta_2) \, = \,  s_0(d_0f) \,,
\end{align*}

Autrement dit, on a que $(d_0\eta_1) \circ f = s_0(d_1f)$ et $f \circ (d_2\eta_2) = s_0(d_0f)$; c'est-à-dire $f$ est un isomorphisme de $\G$. Donc $\G$ est un groupoïde.

Enfin, remarquons que l'isomorphisme d'ensembles simpliciaux tronqués $\varphi_{\bullet}\colon\xymatrix@C-8pt{\tau_{2}^{\,*}W \ar[r]^-{\cong} & \tau_{2}^{\,*}\circ  \mathrm{N}(\G)}$ défini par les règles: $\varphi_i=\mathrm{id}_{W_i}$ si $0\leq i\leq 1$ et $\varphi_2$ est le isomorphisme composé:
$$
\xymatrix@+5pt{W_2 \, \cong \, \mathrm{Hom}_{\simp}\big(\Delta^2,W\big) \ar[r]^-{(d_2,d_0)} & \mathrm{Hom}_{\simp}\big(\Lambda^{2,1},W\big) \, = \,  \mathrm{N}(\G)_2\,,} 
$$
induit un isomorphisme d'ensembles simpliciaux $\mathrm{N}(\G)\,\cong\,{\bf csq}_{2}\big(\mathrm{N}(\G)\big)\,\cong\,{\bf csq}_{2}(W)\,\cong\,W$.
\end{proof}

Il se suit des Corollaires \ref{nerfplnmfidd} et \ref{nerf1groupd} que le foncteur nerf $\mathrm{N}\colon\xymatrix@C-12pt{{\bf cat}\ar[r]&\simp}$ induit par restriction une équivalence entre la catégorie cartésienne fermée ${\bf Grpd}$ et la catégorie cartésienne fermée ${\bf GrpdK}^1$ des $1$-groupoïdes de Kan (voir le Lemme \ref{cartfermengrpKAN}). Plus encore d'après le Corollaire \ref{sonequivMOD} et les isomorphismes \eqref{isonhomhom} et \eqref{isocanogroupd} le foncteur nerf détermine aussi une équivalence entre la catégorie cartésienne fermées $h{\bf Grpd}$ (la catégorie homotopique des groupoïdes) et la catégorie cartésienne fermée des $1$-types d'homotopie $\mathrm{Ho}_1(\simp)$ (la catégorie homotopique des $1$-groupoïdes de Kan).

\renewcommand{\thesubsection}{\S\thesection.\arabic{subsection}}
\subsubsection{}\;\label{grpsontgrpKAN}
\renewcommand{\thesubsection}{\thesection.\arabic{subsection}}

Posons ${\bf Grp}$ pour noter la catégorie des groupes et les morphismes de groupes. Remarquons que d'après les Corollaires \ref{nerfplnmfidd} et \ref{nerf1groupd} le foncteur nerf $\mathrm{N}\colon\xymatrix@C-12pt{{\bf cat}\ar[r]&\simp}$ induit par restriction une équivalence entre la catégorie ${\bf Grp}$ des groupes et ${\bf GrpK}^1$ la sous-catégorie de $\simp_0$ dont les objets sont les $1$-groupoïdes de Kan réduits \emph{i.e.} les $1$-groupes de Kan. 

Plus encore dans ce cas l'isomorphisme $\alpha^1$ de \eqref{isocanogroupd} nous fourni d'un isomorphisme naturel: 
$$
\xymatrix{G\ar[r]^-{\alpha_G}_-{\cong} & \pi_1\big(\mathrm{N}(G)\big)} \; \text{\, pour tout groupe \, $G$}\,; 
$$
autrement dit la restriction du premier groupe d'homotopie des ensembles simpliciaux pointés aux $1$-groupes de Kan réduits est un inverse à isomorphisme près de la restriction du foncteur nerf  aux groupes:
\begin{equation}\label{eqpi1grps}
{\bf Grp}
\vcenter{\xymatrix@C+15pt{
\phantom{a}\ar@{}[r]|{\simeq}\ar@<-4pt>@/_10pt/[r]_{\mathrm{N}}&\phantom{a}
\ar@<-4pt>@/_10pt/[l]_-{\pi_1}}}
{\bf GrpK}^1\,.
\end{equation}

En particulier on a une équivalence entre la catégories des groupes ${\bf Grp}$ et la catégorie homotopique des $1$-groupes $\mathrm{Ho}_1(\simp_0)$:
\begin{equation}\label{eqpi1grps2}
{\bf Grp}
\vcenter{\xymatrix@C+15pt{
\phantom{a}\ar@{}[r]|{\simeq}\ar@<-4pt>@/_10pt/[r]_{\mathrm{N}}&\phantom{a}
\ar@<-4pt>@/_10pt/[l]_-{\pi_1}}}
\mathrm{Ho}_1(\simp_0)\,.
\end{equation}
(voir le Corollaire \ref{equivgrpgrpK0}).

Montrons:

\begin{corollaire}\label{adjonctiongrpss}
Le foncteur restriction du foncteur premier groupe d'homotopie des ensembles simpliciaux pointés aux ensembles simpliciaux réduits est un adjoint à gauche du foncteur nerf pour les groupes:
$$
{\bf Grp}
\vcenter{\xymatrix@C+15pt{
\phantom{a}\ar@{}[r]|{\perp}\ar@<-4pt>@/_10pt/[r]_{\mathrm{N}}&\phantom{a}
\ar@<-4pt>@/_10pt/[l]_-{\pi_1}}}
\simp_0\,,
$$
\end{corollaire}
\begin{proof}
Si $G$ est un groupe et $X$ est un ensemble simplicial réduit on a des isomorphismes naturels:
\begin{equation}\label{ferfer1}
\pi_0\Big(\underline{\mathrm{Hom}}_{\simp_0}\big(X,\mathrm{N}(G)\big)\Big) \; \cong \; \mathrm{Hom}_{\simp_0}\big(X,\mathrm{N}(G)\big)\,,
\end{equation}
\begin{equation}\label{ferfer2}
\pi_0\Big(\underline{\mathrm{Hom}}_{\simp_0}\big(X,\mathrm{N}(G)\big)\Big) \; \cong \; \big[X,\mathrm{N}(G)\big]_1^{red}
\end{equation}
\begin{equation} \label{ferfer3}
\text{et} \qquad \; \mathrm{Hom}_{\bf Grp}\big(\pi_1(X),G\big) \; \cong \; \big[X,\mathrm{N}(G)\big]_1^{red}
\end{equation}
où $\big[\,\cdot\,,\,\cdot\,\big]_1^{red}$ note l'ensemble des morphismes de la catégorie $\mathrm{Ho}_1(\simp_0)\,=\,\simp_0\big[\big({\bf W}^{red}_1\big)^{-1}\big]$ des $1$-types d'homotopie réduits.

En effet \eqref{ferfer1} se suit du Lemme \ref{HOMgrpdred} parce que $\mathrm{N}(G)$ est un $1$-groupe de Kan, \eqref{ferfer2} se suit du fait que $\mathrm{N}(G)$ soit un objet fibrant de la catégorie de modèles $(\simp_0,{\bf W}^{red}_1, {\bf mono},{\bf fib}^{red}_1)$  et \eqref{ferfer3} est une conséquence de l'équivalence de catégories \eqref{eqpi1grps2}.

Donc on a des isomorphismes naturels:
$$
\mathrm{Hom}_{\bf Grp}\big(\pi_1(X),G\big) \; \cong \; \mathrm{Hom}_{\simp_0}\big(X,\mathrm{N}(G)\big)\,.
$$
\end{proof}

Rappelons pour conclure que l'enrichissement $\mathrm{Hom}_{\bf Grp}$ de la catégorie des groupes ${\bf Grp}$ dans la catégorie des ensembles ${\bf Ens}$ est tensoré et cotensoré, c'est-à-dire si $G$ et $H$ sont de groupes on a de foncteurs adjoints:
$$
{\bf Grp}
\vcenter{\xymatrix@C+15pt{
\phantom{a}\ar@{}[r]|{\perp}\ar@<-4pt>@/_10pt/[r]_{\mathrm{Hom}_{{\bf Grp}}(G,\,\cdot\,)}&\phantom{a}
\ar@<-4pt>@/_10pt/[l]_-{G\otimes\,\cdot\, }}}
{\bf Ens}
\qquad\text{et}\qquad
{\bf Grp}^{op}
\vcenter{\xymatrix@C+15pt{
\phantom{a}\ar@{}[r]|{\perp}\ar@<-4pt>@/_10pt/[r]_{\mathrm{Hom}_{\bf Grp}(\,\cdot\,,H)}&\phantom{a}
\ar@<-4pt>@/_10pt/[l]_-{H^{\,\cdot\,} }}}
{\bf Ens}
$$

En effet si $A$ est un ensemble $H^A$ est le groupe des fonctions $\mathrm{Hom}_{\bf Ens}\big(A,H\big)$ dont le produit est défini argument par argument et $G\otimes A$ est le sous-groupe de $H^A$ des fonctions à support fini:
$$
G\otimes A \, = \, \big\{ \, \sigma\colon\xymatrix@C-10pt{A\ar[r]&G} \,\big|\, \text{$\sigma(a)=e_G$ pour tout $a\in A$ sauf pour un nombre fini} \, \big\}\,.
$$

\renewcommand{\thesubsection}{\S\thesection.\arabic{subsection}}
\subsection{}\;
\renewcommand{\thesubsection}{\thesection.\arabic{subsection}}

Soit $X$ un ensemble simplicial réduit et $G$ un groupe. Une \emph{fonction additive de $X$ à valeurs dans $G$} est une fonction de l'ensemble des $1$-simplexes de $X$ vers l'ensemble sous-jacent à $G$:
$$
\xymatrix@C+5pt{X_1\ar[r]^{D}  & G}
$$
vérifiant les propriétés:
\begin{enumerate}
\item $D(s_0\star) \, = \, e_G $, où $\star$ est le seul $0$-simplexe de $X$ et $e_G$ est l'élément neutre de $G$.
\item Si $\alpha\in X_2$ on a que $D\big(d_1\,\alpha\big) = D\big(d_2\,\alpha\big)\,\cdot\, D\big(d_0\,\alpha\big)$ dans $G$.
\end{enumerate}

On pose ${\bf add}_X(G)$ pour noter l'ensemble des fonctions additives de $X$ à valeurs dans $G$. Si ${\bf Grp}$ note la catégorie des groupes et les morphismes de groupes, on a un foncteur:
\begin{equation}\label{determmm0}
\xymatrix@R=5pt@C+15pt{
{\simp_0}^{op} \, \times \, \text{${\bf Grp}$} \ar[r]^-{{\bf det}}  &  {\bf Ens}\,,\\
\big(X, G\big) \;\, \ar@{}[r]|-{\longmapsto} & {\bf add}_X(G)}
\end{equation}
défini dans un morphisme d'ensembles simpliciaux réduits $f\colon\xymatrix@C-10pt{Y\ar[r]&X}$ et un morphisme de groupes  $\varphi\colon\xymatrix@C-10pt{G\ar[r]&H}$ par la fonction:
$$
\xymatrix@C+20pt@R=1pt{ 
{\bf add}_X(G)  \ar[r]^-{ {\bf add}_F(\varphi)} &  {\bf add}_Y(H)  \\
 D \quad \ar@{}[r]|-{\longmapsto} & \;  \varphi \circ D\circ F_1}
$$

Note que la fonction ${\bf add}_F(\varphi)$ est bien définie parce que: 
$$
\big( \varphi \circ D\circ F_1\big) (s_0 \star) = \varphi \Big(D\big(s_0 \, F_0(\star)\big)\Big) = \varphi \big(D (s_0\star)\big) = \varphi (e_G) = e_ H\,,
$$ 
et pour tout $\beta\in Y_2$ on a:
\begin{align*}
\big(\varphi\circ D \circ F_1\big) (d_1\,\beta) \,= \, & \varphi \Big(D\big(d_1\,(F_2 \beta)\big)\Big)  \\ 
 = \, & \varphi \Big(D\big(d_2\,(F_2\beta)\big)\,\cdot\, D\big(d_0\,(F_2\beta)\big)\Big)\\
 = \, & \varphi \Big(D\big(d_2\,(F_2\beta)\big)\Big)\,\cdot\, \varphi \Big( D\big(d_0\,(F_2\beta)\big)\Big)\\
 = \, & \big(\varphi\circ D \circ F_1 \big)(d_2\,\beta)\,\cdot\,\big(\varphi\circ D \circ F_1 \big)(d_0\,\beta)\,.
 \end{align*}

\begin{lemme}\label{lemmeinconudet}
Les foncteurs ${\bf add}_{\bullet_1}\big(\,\bullet_2\,\big)$ et $\mathrm{Hom}_{\simp_0} \big( \, \bullet_1\, , \, \mathrm{N}(\bullet_2) \, \big)$ de source la catégorie produit ${\simp_0}^{op}\times {\bf Grp}$ et but la catégorie des ensembles ${\bf Ens}$ sont naturellement isomorphes. 
\end{lemme}
\begin{proof}
Remarquons que si $G$ est un groupe et $X$ est un ensemble simplicial réduit, la fonction naturelle:
\begin{equation} \label{det0repress}
\xymatrix@C+5pt@R=1pt{
\mathrm{Hom}_{\simp_0} \big(X , \mathrm{N}(G)\big)\ar[r] & {\bf add}_{X}\big(G\big)\\
f_\bullet\quad \ar@{}[r]|-{\longmapsto} & \; f_1}
\end{equation}
est bien définie: En effet si $f\colon\xymatrix@C-8pt{X\ar[r]&\mathrm{N}(G)}$ est un morphisme d'ensembles simpliciaux, on a d'un côté que $f_1(s_0\star) = s_0(f_0\star) = s_0 (\star) = e_G$. D'un autre si $\alpha\in X_2$ on a que $f_2(\alpha)$ est un $2$-simplexe du nerf de $G$ tel que $d_i\circ f_2(\alpha) = f_1\circ d_i(\alpha)$, autrement dit on a que $f_1\big(d_1\,\alpha\big) = f_1\big(d_2\,\alpha\big)\,\cdot\, f_1\big(d_0\,\alpha\big)$ dans $G$.

D'après le Lemme \ref{nerfgrp2cos} l'ensemble simplicial $\mathrm{N}(G)$ est $2$-cosquelettique, donc pour montrer que la fonction \eqref{det0repress} est bijective il suffit de vérifier que la fonction: 
\begin{equation} \label{det0repress2}
\xymatrix@C+5pt@R=1pt{
\mathrm{Hom}_{\simp_{\leq 2}} \Big( \tau_{2}^{*} (X) , \tau_{2}^{*}\big(\mathrm{N}(G)\big)\Big)\ar[r] & {\bf add}_{X}\big(G\big)\\
f_\bullet\quad \ar@{}[r]|-{\longmapsto} & \; f_1}
\end{equation}
est bijective où $\tau_{2}^*\colon\xymatrix@C-10pt{\simp\ar[r]&\simp_{\leq 2}}$ est le foncteur de troncation.

\emph{La fonction \eqref{det0repress2} est injective:}

Soient $f,g\colon\xymatrix@C-8pt{\tau_{2}^*(X)\ar[r]&\tau_{2}^*\big(\mathrm{N}(G)\big)}$  deux morphismes d'ensembles simpliciaux tronqués tels que $f_1(a)=g_1(a)$ pour tout $a\in X_1$. 

On a carrément que $f_0=g_0$. D'un autre vu que $\mathrm{N}(G)$ est un ensemble simplicial faiblement $1$-cosquelettique et $d_i\big(f_2(\alpha)\big) = f_1\big(d_i(\alpha)\big) = g_1\big(d_i(\alpha)\big) = d_i\big(g_2(\alpha)\big)$ si $\alpha\in X_2$ il se suit que $f_2(\alpha)=g_2(\alpha)$ pour tout $\alpha\in X_2$.

\emph{La fonction \eqref{det0repress2} est surjective:}

Si $D\colon\xymatrix@C-8pt{X_1\ar[r] & G}$ est un $0$-déterminant de $X$ à valeurs dans $G$, on définit un morphisme d'ensembles simpliciaux tronqués $f^D\colon\xymatrix@C-8pt{\tau_{2}^*(X)\ar[r]&\tau_{2}^*\big(\mathrm{N}(G)\big)}$ comme suit: $f^D_0(\star)=\star$ si $X_0=\{\star\}$ et $f^D_1(a)=D(a)$ si $a\in X_1$. D'un autre, si $\alpha\in X_2$ vu que par hypothèse  $D(d_1\alpha)=D(d_2\alpha)\,\cdot\,D(d_0\alpha)$ on sait qu'il existe un seul $2$-simplexe $\eta$ du nerf de $G$ tel que $d_i(\eta)=D(d_i\alpha)$; on pose $f^D_2(\alpha)=\eta$ c'est-à-dire $f^D_2(\alpha)$ est le seul $2$-simplexe de $G$ tel que $d_i\big(f^D_2(\alpha)\big) = D(d_i\alpha) = f_1^D(d_i\alpha)$.

Vérifions que $f^D$ est un morphisme d'ensembles simpliciaux tronqués. Par ailleurs:
$$
f^D_0(d_i a) = f^D_0(\star) = \star = d_i (f^D_1a)   \quad \text{si \; $a\, \in \, X_1$\,,}\qquad \quad f^D_1(s_0\star) = D(s_0\star) = e_G =  s_0(\star) =  s_0(f^D_0\star)\,,
$$
$$
\text{et}\qquad f^D_1 (d_i \alpha) = D(d_i\alpha) = d_i(f_2^D\alpha)  \quad \text{si \; $\alpha \, \in \, X_2$\,.}
$$

En plus si $a\in X_1$ on a que $f^D_2(s_j a) = s_j(f_1^Da)$ parce que $d_i\big(f^D_2(s_j a)\big) = d_i \big(s_j(f_1^Da)\big)$, en effet: 
$$
d_i\big(f^D_2(s_j a)\big) = D\big(d_i(s_j a)\big)= 
\begin{cases} D(s_{j-1}d_i \, a)  = e_G & \text{si \, $i<j$} \\ D(a) & \text{si \, $j\leq i\leq j+1$} \\ D(s_{j}d_{i-1} \, a)  = e_G & \text{si \, $j+1<i$}\end{cases}
$$
$$
d_i \big(s_j(f_1^Da)\big) = 
\begin{cases} s_{j-1}d_i  (Da) = e_G & \text{si \, $i<j$} \\ D(a) & \text{si \, $j\leq i\leq j+1$} \\ s_{j}d_{i-1}  (Da)  = e_G & \text{si \, $j+1<i$\,.}\end{cases}
$$

Donc $f^D\colon\xymatrix@C-8pt{\tau_{2}^*(X)\ar[r]&\tau_{2}^*\big(\mathrm{N}(G)\big)}$ est bien un morphisme tel que $f^D_1=D$.
\end{proof}

Montrons finalement:
 
\begin{corollaire}\label{additivffm}
Si $G$ est un groupe et $X$ est un ensemble simplicial réduit les foncteurs:
$$
\vcenter{\xymatrix@R=5pt@C+15pt{
{\simp_0}^{op} \ar[r]^-{{\bf add}_{\bullet}(G)}  &  {\bf Ens}}}
\qquad\text{et}\qquad
\vcenter{\xymatrix@R=5pt@C+15pt{
\text{${\bf Grp}$} \ar[r]^-{{\bf add}_X(\,\bullet\,)}  &  {\bf Ens}}}
$$
sont représentables par $\mathrm{N}(G)$ et $\pi_1(X)$ respectivement. 

De plus si $f\colon\xymatrix@C-8pt{X\ar[r]&Y}$ est une $1$-équivalence faible des ensembles simpliciaux réduits, la fonction induite:
$$
\xymatrix@C+10pt{{\bf add}_{Y}\big(G\big) \ar[r]^-{{\bf add}_{f}\big(G\big)} & {\bf add}_{X}\big(G\big)}
$$
est bijective et le foncteur induit:
$$
\xymatrix@C+35pt{\mathrm{Ho}_1(\simp_0)^{op} \, = \, {\simp_0}\big[ \big({\bf W}^{red}_1\big)^{-1}\big]^{op}\ar[r]^-{h{\bf add}_{\bullet}(G)}&{\bf Ens}}
$$
est représentable aussi par l'ensemble simplicial réduit $\mathrm{N}(G)$.
\end{corollaire}
\begin{proof}

On déduit du Corollaire \ref{adjonctiongrpss} et du Lemme \ref{lemmeinconudet} que si $G$ est un groupe et $X$ un ensemble simplicial réduit les foncteurs:
$$
\vcenter{\xymatrix@R=5pt@C+15pt{
{\simp_0}^{op} \ar[r]^-{{\bf add}_{\bullet}(G)}  &  {\bf Ens}}}
\qquad\text{et}\qquad
\vcenter{\xymatrix@R=5pt@C+15pt{
\text{${\bf Grp}$} \ar[r]^-{{\bf add}_X(\,\bullet\,)}  &  {\bf Ens}}}
$$
sont représentables par $\mathrm{N}(G)$ et $\pi_1(X)$ respectivement. 

D'un autre il se suit des Corollaires \ref{1grpKAN} et \ref{adjonctiongrpss} et du Lemme \ref{HOMgrpdred} que si $f\colon\xymatrix@C-8pt{X\ar[r]&Y}$ est une $1$-équivalence faible des ensembles simpliciaux réduits la fonction induite:
$$
\xymatrix{{\bf add}_{Y}\big(G\big) \ar[r] & {\bf add}_{X}\big(G\big)}
$$
est bijective.

\end{proof}

\section{La catégorie des $2$-groupes}\label{2groupeseection}

\renewcommand{\thesubsection}{\S\thesection.\arabic{subsection}}
\subsection{}\;
\renewcommand{\thesubsection}{\thesection.\arabic{subsection}}

Rappelons qu'une \emph{catégorie monoïdale} (voir par exemple \cite{lane}) est la donnée d'une catégorie $\mathcal{M}$, un foncteur $\xymatrix@C-3pt{\mathcal{M}\times\mathcal{M}\ar[r]^-{\otimes}&\mathcal{M}}$, un objet distingué $\mathbb{1}$ et d'isomorphismes naturels: 
\begin{equation}\label{isos}
\text{$\xymatrix{ (X\otimes Y)\otimes Z \ar[r]^{a_{X,Y,Z}} & X\otimes (Y\otimes Z) }$, $\xymatrix@C-5pt{X \ar[r]^-{l_{X}} & \mathbb{1}\otimes X}$ et $\xymatrix@C-5pt{X\ar[r]^-{r_{X}} & X\otimes\mathbb{1}}$;}
\end{equation}
tels que les diagrammes suivants soient commutatifs:
\begin{equation}\label{pentagon} 
\def\objectstyle{\scriptstyle}
\def\labelstyle{\scriptstyle}
\vcenter{\xymatrix@C+18pt{
\big((W\otimes X)\otimes Y\big)\otimes Z\ar[r]^-{a_{W\otimes X,Y,Z}}\ar[d]_-{a_{X,Y,Z}\otimes Z}& (W\otimes X)\otimes (Y\otimes Z)  \ar[r]^-{a_{W,X,Y\otimes Z}}& W\otimes\big(X\otimes (Y\otimes Z)\big)\\
\big(W \otimes (X\otimes Y)\big)\otimes Z \ar[rr]_-{a_{W,X\otimes Y,Z}}&&W\otimes\big((X\otimes Y)\otimes Z\big)\ar[u]_-{W\otimes a_{X,Y,Z}}}}
\end{equation}
\begin{equation}\label{triangle} 
\text{et}\qquad
\def\objectstyle{\scriptstyle}
\def\labelstyle{\scriptstyle}
\vcenter{\xymatrix{
(X\otimes \mathbb{1}) \otimes Y  \ar[rr]^-{a_{X,\mathbb{1},Y}} &&    X\otimes (\mathbb{1}\otimes Y)\\
&X\otimes Y\ar[ru]_-{X\otimes l_Y}\ar[lu]^-{r_X\otimes Y}&}}
\end{equation}
pour $X$, $Y$, $Z$ et $W$ des objets quelconques de $\mathcal{M}$.

On vérifie:

\begin{lemme}\label{otrosdiagg}
Si $A$ et $B$ sont des objets d'une catégorie monoïdale quelconque, les triangles:
$$
\vcenter{\xymatrix@C-4pt@R+15pt{
(A\otimes B)\otimes \mathbb{1} \ar[rr]^-{a_{A,B,\mathbb{1}}}& & A\otimes (B\otimes \mathbb{1}) \\
&A\otimes B\ar[lu]^-{r_{A\otimes B}} \ar[ru]_-{A\otimes r_B}&}}\qquad et \qquad
\vcenter{\xymatrix@C-4pt@R+15pt{
& A\otimes B \ar[rd]^-{\ell_{A\otimes B}}\ar[dl]_-{\ell_A\otimes B}&\\
(\mathbb{1}\otimes A)\otimes B \ar[rr]_-{a_{\mathbb{1},A,B}}& & \mathbb{1}\otimes (A\otimes B)}}
$$
sont des diagrammes commutatifs. Plus encore on a que $\ell_\mathbb{1}=r_\mathbb{1}$.
\end{lemme}

Étant données deux catégories monoïdales $\mathcal{M}$ et $\mathcal{N}$, on définit un \emph{foncteur monoïdal lax et unitaire de $\mathcal{M}$ vers $\mathcal{N}$} (resp. \emph{foncteur monoïdal fort et unitaire}), comme la donnée d'un foncteur $\xymatrix@C-3pt{\mathcal{M}\ar[r]^-{F}&\mathcal{N}}$ tel que $F\big(\mathbb{1}_\mathcal{M}\big) = \mathbb{1}_\mathcal{N}$ et des morphismes naturels (resp. isomorphismes naturels): 
\begin{equation}\label{mu}
\xymatrix@C+8pt{ F(X)\otimes F(Y) \ar[r]^-{m_{X,Y}} & F(X\otimes Y) \,,}
\end{equation}
tels que les diagrammes suivants soient commutatifs:
\begin{equation}\label{hexa}
\def\objectstyle{\scriptstyle}
\def\labelstyle{\scriptstyle}
\vcenter{\xymatrix@C+18pt@R-3pt{
F\big((X\otimes Y)\otimes Z\big)   \ar[d]_-{Fa_{X,Y,Z}}&F(X\otimes Y)\otimes FZ  \ar[l]_-{m_{X\otimes Y, Z}}&(FX\otimes FY)\otimes FZ\ar[d]^-{a_{FX,FY,FZ}} \ar[l]_-{m_{X,Y}\otimes FZ}\\
F\big(X\otimes (Y\otimes Z)\big)  &FX\otimes F(Y\otimes Z) \ar[l]^-{m_{X,Y\otimes Z}} &FX\otimes (FY\otimes FZ)  \ar[l]^-{FX\otimes m_{Y,Z}}  }}
\end{equation}
$$\text{et}$$
\begin{equation}\label{hexa2}
\def\objectstyle{\scriptstyle}
\def\labelstyle{\scriptstyle}
\vcenter{\xymatrix@C-3pt@R-5pt{
\mathbb{1}\otimes FX \ar[dd]_-{m_{\mathbb{1},X}} &&FX\otimes \mathbb{1}  \ar[dd]^-{m_{X,\mathbb{1}}} \\
&FX\ar[lu]_{l_{FX}}\ar[ru]^{r_{FX}}\ar[rd]_{Fr_X}\ar[ld]^-{Fl_X}&\\
 F(\mathbb{1}\otimes X)   &&F(X\otimes\mathbb{1}) 
 }}
\end{equation}

Les catégories monoïdales et les foncteurs monoïdaux lax et unitaires (resp. foncteurs monoïdaux forts et unitaires) forment une catégorie pointée $\cat^{\otimes}_{lax,*}$ (resp. $\cat^{\otimes}_{fort,*}$) dont la composition est définie comme suit: Si $\xymatrix@C+8pt{\mathcal{M}\ar[r]^-{(F,m^F)}&\mathcal{M}'}$ et $\xymatrix@C+8pt{\mathcal{M}'\ar[r]^-{(G,m^G)}&\mathcal{M}''}$ sont de foncteur monoïdaux forts et unitaires (resp. foncteurs monoïdaux forts et unitaires) le composé $\xymatrix@C+25pt{\mathcal{M}\ar[r]^-{(G\circ F ,m^{G\circ F})}&\mathcal{M}''}$ est constitué du foncteur $G\circ F$ et de la transformation naturelle $m^{G\circ F}$ définie dans des objets $X$ et $Y$ de $\G$ par l'isomorphisme:
$$
\xymatrix@C+15pt{
G\circ F(X) \otimes G\circ F(Y) \ar[r]_-{m^G_{FX, FY}} \ar@/^18pt/[rr]^-{m^{G\circ F}_{X,Y}} & G\big(F(X)\otimes F(Y)\big) \ar[r]_-{G(m_{X,Y}^F)}  &  G\circ F (X\otimes Y)\,.
}
$$

Remarquons finalement qu'on a une $2$-catégorie $\underline{\cat}^{\otimes}_{lax,*}$ (resp. $\underline{\cat}^{\otimes}_{fort,*}$) dont la catégorie sous-jacente est $\cat^{\otimes}_{lax,*}$ (resp. $\cat^{\otimes}_{fort,*}$) et  les $2$-flèches sont les \emph{transformations} définies comme suit: Si $\xymatrix@C+3pt{\mathcal{M}\ar[r]^-{(F,m^F)}&\mathcal{N}}$ et $\xymatrix@C+3pt{\mathcal{M}\ar[r]^-{(G,m^G)}&\mathcal{N}}$ sont de foncteur monoïdaux forts et unitaires, une \emph{transformation de $(F,m^F)$ vers $(G,m^G)$} est la donnée d'une transformation naturelle de foncteurs $\xymatrix@C+3pt{\mathcal{M}\rtwocell^F_{G}{\eta}&\mathcal{N}}$ telle que $\eta_{\mathbb{1}_\mathcal{M}} = \mathrm{id}_{\mathbb{1}_\mathcal{N}}$ $\big(\text{rappelons que \,} F\mathbb{1}_\mathcal{M} = \mathbb{1}_\mathcal{N} = G\mathbb{1}_\mathcal{M}\big)$ et telle que le diagramme suivant soit commutatif:
\begin{equation}
\def\objectstyle{\scriptstyle}
\def\labelstyle{\scriptstyle}
\vcenter{\xymatrix@C+10pt@R-1pt{
FX\otimes FY \ar[d]_-{m_{X,Y}^F} \ar[r]^-{\eta_{X}\otimes\eta_{Y}} &GX\otimes GY   \ar[d]^-{m_{X,Y}^G}  \\
F(X\otimes Y)\ar[r]_-{\eta_{X\otimes Y}} & G(X\otimes Y) }}\,.
\end{equation}

Dans l'énoncé qui suit on montre que le $2$-foncteur d'oubli:
$$
\xymatrix@C-5pt{\underline{\cat}^{\otimes}_{lax,*} \ar[r] & \underline{\bf cat}}  
\qquad\quad
\big(\text{resp.}\;\; \xymatrix@C-5pt{\underline{\cat}^{\otimes}_{fort,*} \ar[r] & \underline{\bf cat}}\big) 
$$
reflète les équivalences internes et les isomorphismes.

\begin{lemme} \label{equiequi}
Si $\xymatrix@C-8pt{\mathcal{M}\ar[r]^-{F}&\mathcal{M}'}$ est un morphisme monoïdal lax et unitaire (resp. morphisme monoïdal fort et unitaire) entre catégories monoïdales dont le foncteur sous-jacent est pleinement fidèle et essentiellement surjectif, alors il existe un morphisme monoïdal lax et unitaire (resp. morphisme monoïdal fort et unitaire) $\xymatrix@C-8pt{\mathcal{M}'\ar[r]^-{G}&\mathcal{M}}$ et de transformations: 
$$
\xymatrix@C-15pt{F\,\circ\,G\;\ar@{=>}[r]^-{\alpha}& \; \mathrm{id}_{\mathcal{M}'}} \qquad \text{et} \qquad 
\xymatrix@C-15pt{G\,\circ\,F\;\ar@{=>}[r]^-{\beta}& \; \mathrm{id}_{\mathcal{M}}}
$$ 
telles que les morphismes $\alpha_A$ et $\beta_X$ sont des isomorphismes pour tous les objets $A$ et $X$.

Plus encore, si on suppose que $F$ est pleinement fidèle et bijectif sur les objets alors il existe un morphisme monoïdal lax et unitaire (resp. morphisme monoïdal fort et unitaire) $\xymatrix@C-8pt{\mathcal{G}'\ar[r]^-{G}&\mathcal{G}}$ tel que $G\circ F = \mathrm{id}_\mathcal{M}$ et $F\circ G = \mathrm{id}_{\mathcal{M}'}$.
\end{lemme}
\begin{proof}
Supposons que $\xymatrix@C-8pt{\mathcal{M}\ar[r]^-{F}&\mathcal{GM}'}$ est un morphisme monoïdal lax et unitaire (resp. morphisme monoïdal fort et unitaire) entre catégories monoïdales tel que le foncteur sous-jacent à $F$ soit pleinement fidèle et essentiellement surjectif. On va construit un morphisme monoïdal lax et unitaire (resp. morphisme monoïdal fort et unitaire) $\xymatrix@C-8pt{\mathcal{M}'\ar[r]^-{G}&\mathcal{M}}$ ainsi qu'une transformation $\xymatrix@C-15pt{F\,\circ\,G\;\ar@{=>}[r]^-{\alpha}& \; \mathrm{id}_{\mathcal{G}'}}$ de la façon suivante:

Étant donné un objet $A$ de $\mathcal{M}'$ choisissons un objet $G(A)$ de $\mathcal{M}$, tel que $F(GA)$ et $A$ soient isomorphes dans $\mathcal{M}'$ ($F$ est essentiellement surjectif). On peut alors choisir aussi un isomorphisme $\alpha_A\colon\xymatrix@-8pt{FG(A)\ar[r]&A}$ de $\mathcal{M}'$ et supposer de plus que  $G(\mathbb{1})=\mathbb{1}$ et $\alpha_\mathbb{1}=\mathrm{id}_{\mathbb{1}}$. Si $f\colon\xymatrix@-8pt{A\ar[r]&B}$ est un morphisme de $\mathcal{M}'$, on définit $G(f)\colon\xymatrix@-8pt{G(A)\ar[r]&G(B)}$ comme le seul morphisme de $\mathcal{M}$ tel que $\alpha_B \circ FG(f)= f\circ \alpha_A$. On vérifie sans peine qu'on obtient ainsi un foncteur $\xymatrix@C-8pt{\mathcal{M}'\ar[r]^-{G}&\mathcal{M}}$ entre les catégories sous-jacents aux catégories monoïdales et un isomorphisme naturelle de foncteurs $\xymatrix@C+3pt{\mathcal{M}'\rtwocell^{FG}_{\mathrm{id}}{\alpha}&\mathcal{M}'}$, tel que $G(\mathbb{1})=\mathbb{1}$ et $\alpha_\mathbb{1}=\mathrm{id}_{\mathbb{1}}$.

Complétons le foncteur $\xymatrix@C-8pt{\mathcal{M}'\ar[r]^-{G}&\mathcal{M}}$ à un morphisme monoïdal lax et unitaire (resp. morphisme monoïdal fort et unitaire) comme suit: Si $A$ et $B$ sont des objets de $\mathcal{M}'$ on pose:
$$
\xymatrix@+8pt{G(A)\otimes G(B)\ar[r]^-{m_{A,B}^{G}}&G(A\otimes B)}
$$ 
pour noter le seul morphisme (resp. isomorphisme) de $\mathcal{M}$ tel que le diagramme:
\begin{equation}\label{transequi}
\xymatrix@C+30pt{
 FG(A)\otimes FG(B) \ar[r]^-{m_{GA,GB}^F}  \ar[d]_-{\alpha_A\otimes \alpha_B}   &F\big( G(A)\otimes G(B)\big)\ar[r]^-{F(m_{A,B}^G)} & F\big(G(A\otimes B)\big)\ar[d]^-{\alpha_{A\otimes B}}  \\
A\otimes B \ar@{=}[rr]& & A\otimes B }
\end{equation} 
soit commutatif. En particulier $\xymatrix@C+3pt{\mathcal{M}'\rtwocell^{FG}_{\mathrm{id}}{\alpha}&\mathcal{M}'}$ serait une transformation entre morphisme de catégories monoïdales si $G$ était un morphisme monoïdal.

Montrons que $\xymatrix@C-8pt{\mathcal{M}'\ar[r]^-{G}&\mathcal{M}}$ est effectivement un morphisme monoïdal, c'est-à-dire montrons que les diagrammes qui suivent sont commutatifs:
\begin{equation}\label{laxo1proof}
\def\objectstyle{\scriptstyle}
\def\labelstyle{\scriptstyle}
\vcenter{\xymatrix@C+28pt@R-3pt{
G\big((A\otimes B)\otimes C\big)\ar[d]_-{Ga_{A,B,C}}&G(A\otimes B)\otimes GC\ar[l]_-{m^G_{A\otimes B, C}}&
(GA\otimes GB)\otimes GC\ar[d]^-{a_{GA,GB,GC}}\ar[l]_-{m^G_{A,B}\otimes GC}\\
G\big(A\otimes (B\otimes C)\big)  &GA\otimes G(B\otimes C)\ar[l]^-{m^G_{A,B\otimes C}}&
GA\otimes (GB\otimes GC)\ar[l]^-{GA\otimes m^G_{B,C}}}}
\end{equation}
$$\text{et}$$
\begin{equation}\label{laxo2proof}
\def\objectstyle{\scriptstyle}
\def\labelstyle{\scriptstyle}
\vcenter{\xymatrix@C+5pt@R-5pt{
\mathbb{1}\otimes GX \ar[dd]_-{m^G_{\mathbb{1},X}}&&GX\otimes \mathbb{1}\ar[dd]^-{m^G_{X,\mathbb{1}}} \\
&GX\ar[lu]_{l_{GX}}\ar[ru]^{r_{GX}}\ar[rd]_{Gr_X}\ar[ld]^-{Gl_X}&\\
 G(\mathbb{1}\otimes X)&&G(X\otimes\mathbb{1})
}}
\end{equation}

De façon équivalente montrons que les images par le foncteur fidèle $F$ des diagrammes \eqref{laxo1proof} et \eqref{laxo2proof} sont commutatifs. Pour commencer, on montrer que les triangles $F$\eqref{laxo2proof} sont commutatifs en considérant les diagrammes commutatifs:
$$
\xymatrix@C+25pt{
FG(X)  \ar[r]|-{l_{FGX}} \ar[d]_-{\alpha_X}\ar@/^18pt/[rr]^-{F(l_{GX})} & \mathbb{1}\otimes FG(X) \ar[r]|-{m^F_{\mathbb{1},GX}} \ar[d]|-{\mathbb{1}\otimes\alpha_X}& F(\mathbb{1}\otimes GX)\ar[r]^-{F(m^G_{\mathbb{1},X})}& FG(\mathbb{1}\otimes X) \ar[d]^-{\alpha_{\mathbb{1},X}}\\
X\ar[r]_-{l_X} &\mathbb{1}\otimes X   \ar[rr]_-{\mathrm{id}} &&\mathbb{1}\otimes X} 
$$
$$\text{et}$$
$$
\xymatrix@C+25pt{
FG(X)  \ar[r]|-{r_{FGX}} \ar[d]_-{\alpha_X}\ar@/^18pt/[rr]^-{F(r_{GX})} & FG(X)\otimes\mathbb{1}\ar[r]|-{m^F_{GX,\mathbb{1}}} \ar[d]|-{\alpha_X\otimes\mathbb{1}}& F(GX\otimes\mathbb{1})\ar[r]^-{F(m^G_{X,\mathbb{1}})}& FG(X\otimes\mathbb{1}) \ar[d]^-{\alpha_{X\otimes \mathbb{1}}}\\
X\ar[r]_-{r_X} &X\otimes\mathbb{1}\ar[rr]_-{\mathrm{id}} &&X\otimes\mathbb{1}} 
$$
et en se rappelant que: 
$$
\vcenter{\xymatrix@C+10pt{
FG(X)\ar[r]^-{FG(l_X)} \ar[d]_-{\alpha_X} & FG(\mathbb{1}\otimes X) \ar[d]^-{\alpha_{\mathbb{1}\otimes X}}\\
X\ar[r]_-{l_X} & \mathbb{1}\otimes X}} \qquad \text{et} \qquad \vcenter{\xymatrix@C+10pt{
FG(X)\ar[r]^-{FG(r_X)} \ar[d]_-{\alpha_X} & FG(X\otimes \mathbb{1}) \ar[d]^-{\alpha_{X\otimes\mathbb{1}}}\\
X\ar[r]_-{r_X} & X\otimes \mathbb{1}}}
$$ 
sont de diagrammes commutatifs parce que $\alpha$ est une transformation naturelle.

D'autre part, on vérifie que $F$\eqref{laxo1proof} est commutatif en l'insérant comme une face d'un cube, dont les autres faces sont les diagrammes commutatifs qui suit:
$$
\def\objectstyle{\scriptstyle}
\def\labelstyle{\scriptstyle}
\vcenter{\xymatrix@C+50pt@R+5pt{
FG\big((A\otimes B)\otimes C\big) \ar[d]_-{\alpha_{(A\otimes B)\otimes C}}& F\big(G(A\otimes B)\otimes GC\big) \ar[l]_-{F(m_{A\otimes B,C}^G)}& F\big((GA\otimes GB) \otimes GC\big) \ar[l]_-{F(m_{A,B}^G\otimes GC)}\\
(A\otimes B)\otimes C\ar@{=}[d] & FG(A\otimes B)\otimes FGC \ar[u]|-{m^F_{G(A\otimes B),GC}}\ar[l]|-{\alpha_{A\otimes B}\otimes \alpha_C}\ar[d]|-{\alpha_{A\otimes B}\otimes FGC}& F(GA\otimes GB) \otimes FGC\ar[l]|-{F(m^G_{A,B})\otimes FGC} \ar[u]_-{m^F_{GA\otimes GB,GC}}\\
(A\otimes B)\otimes C  & (A\otimes B)\otimes FGC\ar[l]^-{(A\otimes B)\otimes \alpha_C} & (FGA\otimes FGB)\otimes FGC\ar[u]_-{m^F_{GA,GB}\otimes FGC}\ar[l]^-{(\alpha_A\otimes \alpha_B)\otimes FGC}
}}\,,
$$
$$
\def\objectstyle{\scriptstyle}
\def\labelstyle{\scriptstyle}
\vcenter{\xymatrix@C+50pt@R+5pt{
FG\big(A\otimes (B\otimes C)\big) \ar[d]_-{\alpha_{A\otimes(B\otimes C)}}& F\big(G(A)\otimes G(B\otimes C)\big) \ar[l]_-{F(m_{A, B\otimes C}^G)}& F\big(GA\otimes (GB \otimes GC)\big) \ar[l]_-{F(GA\otimes m_{B,C}^G)}\\
A\otimes (B\otimes C)\ar@{=}[d] & FGA\otimes FG(B\otimes C) \ar[u]|-{m^F_{GA,G(B\otimes C)}}\ar[l]|-{\alpha_A\otimes \alpha_{B\otimes C}}\ar[d]|-{FGA\otimes \alpha_{B\otimes C}} & FG(A)\otimes F(GB\otimes GC)\ar[l]|-{FGA\otimes F(m^G_{B,C})} \ar[u]_-{m^F_{GA,GB\otimes GC}}\\
A\otimes (B\otimes C)  & FGA\otimes (B\otimes C)\ar[l]^-{\alpha_A\otimes (B\otimes C)} & FGA\otimes (FGB\otimes FGC)\ar[l]^-{FGA\otimes(\alpha_B\otimes\alpha_A)}\ar[u]_-{FGA\otimes m^F_{GB,GC}}\,,
}}\,,
$$
$$
\def\objectstyle{\scriptstyle}
\def\labelstyle{\scriptstyle}
\vcenter{\xymatrix@C+23pt@R-3pt{
(A\otimes B)\otimes C \ar[d]_-{a_{A,B,C}}\ar@{=}[r] & (A\otimes B)\otimes C \ar[d]|-{a_{A,B,C}} & FG\big((A\otimes B)\otimes C \big)\ar[l]_-{\alpha_{(A\otimes B)\otimes C}}\ar[d]^-{F(a_{A,B,C})}\\
A\otimes (B\otimes C) \ar@{=}[r] & A\otimes (B\otimes C)& FG\big(A\otimes (B\otimes C)\big)\ar[l]^-{\alpha_{A\otimes (B\otimes C)}}
}}\,,
$$
$$
\def\objectstyle{\scriptstyle}
\def\labelstyle{\scriptstyle}
\vcenter{\xymatrix@C+23pt@R-3pt{
(FGA\otimes FGB)\otimes FGC \ar[d]_-{a_{FGA,FGB,FGC}} \ar[r]^-{m^F_{GA,GB}\otimes FGC}&F(GA\otimes GB)\otimes FGC\ar[r]^-{m^F_{GA\otimes GB,GC}}& 
F\big((GA\otimes GB)\otimes GC \big)\ar[d]^-{F(a_{GA,GB,GC})}\\
FGA\otimes (FGB\otimes FGC)\ar[r]_-{FGA\otimes m^F_{GB,GC}}&FGA\otimes F(GB\otimes GC)\ar[r]_{m^F_{GA,GB\otimes GC}}& F\big(GA\otimes (GB\otimes GC)\big)
}}
$$
$$\text{et}$$
$$
\def\objectstyle{\scriptstyle}
\def\labelstyle{\scriptstyle}
\vcenter{\xymatrix@C+30pt@R+5pt{
(A\otimes B)\otimes C \ar[d]_-{a_{A,B,C}} & (A\otimes B)\otimes FGC\ar[l]_-{(A\otimes B)\otimes \alpha_C} & (FGA\otimes FGB)\otimes FGC\ar[l]_-{(\alpha_A\otimes \alpha_B)\otimes FGC}\ar[d]^-{a_{FGA,FGB,FGC}}\\
A\otimes (B\otimes C)  & FGA\otimes (B\otimes C)\ar[l]^-{\alpha_A\otimes (B\otimes C)} & FGA\otimes (FGB\otimes FGC)\ar[l]^-{FGA\otimes(\alpha_B\otimes\alpha_A)}
}}\,.
$$

Il manque juste de construire une transformation $\xymatrix@C-15pt{G\,\circ\,F\;\ar@{=>}[r]^-{\beta}& \; \mathrm{id}_{\mathcal{M}}}$. On vérifie sans peine que si on pose $\xymatrix@C-8pt{GF(X)\ar[r]^-{\beta_X}&X}$ pour noter le seul morphisme tel que $F(\beta_X)=\alpha_{FX}$,  alors $\beta_X$ est un isomorphisme de la catégorie monoïdale $\mathcal{M}$ et la famille $\beta=\{\beta_X\}_X$ est bien une transformation entre morphismes monoïdaux. 

Remarquons finalement que si $F$ était un foncteur bijectif sur les objets, dans la construction du foncteur $G$ ci-dessus on peut choisir pour tout objet $A$ de $\mathcal{M}'$ un (seul) objet $G(A)$ de $\mathcal{G}$ tel que $F(GA)=A$. En imposant aussi que $\alpha_A=\mathrm{id}_A$ on déduit que $G$ vérifie $F\circ G = \mathrm{id}$ et $G\circ F = \mathrm{id}$.
\end{proof}

\renewcommand{\thesubsection}{\S\thesection.\arabic{subsection}}
\subsection{}\;\label{2groupessection}
\renewcommand{\thesubsection}{\thesection.\arabic{subsection}}

Étant donné une catégorie monoïdale $\mathcal{M}$, si $X$ est un objet de $\mathcal{M}$ fixé on montre sans peine que les conditions suivantes sont équivalentes:
\begin{enumerate}
\item Les foncteurs:
$$
\xymatrix@C+10pt{\mathcal{M} \ar[r]^-{X\otimes -}& \mathcal{M}}
\qquad\text{et}\qquad 
\xymatrix@C+10pt{\mathcal{M}\ar[r]^-{-\otimes X} & \mathcal{M}}.
$$
sont des équivalences de catégories.
\item Ils existent des objets $X'$ et $X''$ de $\mathcal{M}$, et d'isomorphismes dans $\mathcal{M}$:
$$\xymatrix@C-5pt{X\otimes X' \ar[r]^-{\alpha_X}_-{\cong} & \mathbb{1}&\ar[l]_-{\beta_X}^-{\cong} X''\otimes X\,.}$$
\end{enumerate}

Un objet $X$ d'une catégorie monoïdale $\mathcal{M}$ qui satisfait à une de ceux deux conditions est appelé \emph{inversible}. Un \emph{2-groupe} (\emph{non commutatif}) est une catégorie monoïdale dont la catégorie sous-jacente est un groupoïde et tous ses objets sont inversibles. Un $2$-groupe est dit aussi une \emph{catégorie de Picard non commutatif}.

Notons $2$-${\bf Grp}$ la \emph{catégorie des} $2$-\emph{groupes} $\big(${resp. la $2$-\emph{catégorie des} $2$-\emph{groupes} $2$-$\underline{{\bf Grp}}$}$\big)$ c'est-à-dire la sous-catégorie pleine de $\cat^{\otimes}_{lax,*}$ $\big(${resp. la sous-2-catégorie pleine de $\underline{\bf cat}^{\otimes}_{lax,*}$}$\big)$ dont les objets sont les $2$-groupes. 

Montrons:

\begin{lemme}\label{inverses}
Si $\mathcal{M}$ est une catégorie monoïdale dont la catégorie sous-jacente est un groupoïde, alors $\mathcal{M}$ est un $2$-groupe si et seulement s'ils existent de foncteurs $\xymatrix@C-4pt{\mathcal{M} \ar[r]^-{\iota^d}& \mathcal{M}}$ et $\xymatrix@C-4pt{\mathcal{M} \ar[r]^-{\iota^g}& \mathcal{M}}$ et pour chaque objet $X$ de $\mathcal{M}$ d'isomorphismes:
$$\xymatrix@C+5pt{X\otimes \iota^d(X) \ar[r]^-{\alpha_X}_-{\cong} & \mathbb{1}&\ar[l]_-{\beta_X}^-{\cong} \iota^g(X)\otimes X}\,,$$
tels que si $f\colon\xymatrix@C-10pt{X\ar[r]&Y}$ est un morphisme de $\mathcal{M}$ on a des carrés commutatifs:
$$\xymatrix@C+5pt@R-3pt{
X\otimes \iota^d(X) \ar[r]^-{\alpha_X} \ar[d]_-{f\otimes \iota^d(f)} & \mathbb{1}\ar[d]|-{\mathrm{id}_{\mathbb{1}}}&\ar[l]_-{\beta_X}\iota^g(X)\otimes X\ar[d]^-{\iota^g(f)\otimes f}\\
Y\otimes \iota^d(Y) \ar[r]_-{\alpha_Y} & \mathbb{1}&\ar[l]^-{\beta_Y}\iota^g(Y)\otimes Y}$$
\end{lemme}
\begin{proof}
S'ils existent de foncteurs $\xymatrix@C-4pt{\mathcal{M} \ar[r]^-{\iota^d}& \mathcal{M}}$ et $\xymatrix@C-4pt{\mathcal{M} \ar[r]^-{\iota^g}& \mathcal{M}}$ et d'isomorphismes $\xymatrix@C+5pt{X\otimes \iota^d(X) \ar[r]^-{\alpha_X} & \mathbb{1}&\ar[l]_-{\beta_X} \iota^g(X)\otimes X}$ il est clair que tous les objets de $\mathcal{M}$ sont inversibles, c'est-à-dire $\mathcal{M}$ est un $2$-groupe.

Montrons que si $\mathcal{M}$ est un $2$-groupe alors il existe un foncteur $\xymatrix@C-4pt{\mathcal{M} \ar[r]^-{\iota^d}& \mathcal{M}}$ et pour chaque objet $X$ de $\mathcal{M}$ un isomorphisme $\xymatrix@C-5pt{X\otimes \iota^d(X) \ar[r]^-{\alpha_X}&\mathbb{1}}$ tel que pour tout morphisme $f\colon\xymatrix@C-10pt{X\ar[r]&Y}$ de $\mathcal{M}$:
$$
\xymatrix@R-10pt{
X\otimes \iota^d(X) \ar[dd]_-{f\otimes \iota^d(f)}\ar@/^7pt/[rd]^-{\alpha_X} & \\ &\mathbb{1} \\ Y\otimes \iota^d(Y) \ar@/_7pt/[ru]_-{\alpha_Y} & }
$$
soit un triangle commutatif. On vérifie de façon analogue l'existence du foncteur $\xymatrix@C-4pt{\mathcal{M} \ar[r]^-{\iota^g}& \mathcal{M}}$ et de l'isomorphisme $\xymatrix@C-5pt{\iota^g(X)\otimes X \ar[r]^-{\beta_X}&\mathbb{1}}$ avec la propriété désirée.

Choisissons pour chaque objet $X$ de $\mathcal{M}$ un objet $\iota^d(X)$ et un isomorphisme $\xymatrix@C-5pt{X\otimes \iota^d(X) \ar[r]^-{\alpha_X}& \mathbb{1}}$. Si $\xymatrix@C-5pt{X \ar[r]^-{f} & Y}$ est un morphisme de $\mathcal{M}$, vu que le foncteur $\xymatrix@C-5pt{\mathcal{M} \ar[r]^-{Y\otimes -}& \mathcal{M}}$ est une équivalence de catégories et $\mathcal{M}$ est un groupoïde, il existe un seul morphisme $\xymatrix@C-5pt{\iota^d(X) \ar[r]^-{\iota^d(f)} & \iota^d(Y)}$ de $\mathcal{M}$ tel que le carré suivant est commutatif:
$$
\xymatrix@C+20pt{
Y\otimes \iota^d(X) \ar@{-->}[d]_-{Y\otimes \iota^d(f)} & X\otimes \iota^d(X)  \ar[l]_-{f\otimes\iota^d(X)} \ar[d]^-{\alpha_X}\\
Y\otimes \iota^d(Y) \ar[r]_-{\alpha_Y} & \mathbb{1}}
$$

Vu que dans $\mathcal{M}$ on a un diagramme commutatif:
$$
\xymatrix@C+25pt@C+10pt{
X\otimes \iota^d(X) \ar[rd]|-{f\otimes \iota^d(f)} \ar[r]^-{f\otimes\iota^d(X)}\ar[d]_-{X\otimes \iota^d(f)} & Y\otimes \iota^d(X) \ar[d]^-{Y\otimes\iota^d(f)}\\
X\otimes \iota^d(Y) \ar[r]_-{f\otimes \iota^d(Y)} & Y\otimes \iota^d(Y)}
$$
associé aux morphismes $\xymatrix@C-5pt{X \ar[r]^-{f} & Y}$ et $\xymatrix@C-5pt{\iota^d(X) \ar[r]^-{\iota^d(f)} & \iota^d(Y)}$; le morphisme $\xymatrix@C-5pt{\iota^d(X) \ar[r]^-{\iota^d(f)} & \iota^d(Y)}$ est aussi le seul morphisme de $\mathcal{M}$ tel que le carré suivant est commutatif:
$$
\xymatrix@C+20pt{
X\otimes \iota^d(X) \ar[r]^-{\alpha_X}\ar@{-->}[d]_-{X\otimes \iota^d(f)} & \mathbb{1} \\
X\otimes \iota^d(Y) \ar[r]_-{f\otimes\iota^d(Y)} & Y\otimes \iota^d(Y)\ar[u]_-{\alpha_Y}}
$$
ou encore mieux le seul morphisme tel que le triangle suivant est commutatif:
$$
\xymatrix@R-10pt{
X\otimes \iota^d(X) \ar[dd]_-{f\otimes \iota^d(f)}\ar@/^7pt/[rd]^-{\alpha_X} & \\ &\mathbb{1} \\ Y\otimes \iota^d(Y) \ar@/_7pt/[ru]_-{\alpha_Y} & }
$$

On déduit de cette dernière caractérisation que $f\mapsto \iota^d(f)$ détermine un foncteur $\xymatrix@C-4pt{\mathcal{M} \ar[r]^-{\iota^d}& \mathcal{M}}$.
\end{proof}



Si $\mathcal{G}$ est un 2-groupe on pose $\pi_{0}(\mathcal{G})$ pour noter l'ensemble des composantes connexes par arcs du groupoïde $\mathcal{G}$ \emph{i.e.} l'ensemble des objets à isomorphisme près de $\mathcal{G}$. On vérifie que $\pi_{0}(\mathcal{G})$ est un ensemble muni d'une loi de composition induite par le foncteur $\otimes$ tel que $\pi_{0}(\mathcal{G})$ est un groupe (tout objet de $\mathcal{G}$ est inversible) dont l'élément neutre est la classe de l'objet $\mathbb{1}$. On appelle $\pi_{0}(\mathcal{G})$ \emph{le groupe de composantes connexes par arcs du $2$}-\emph{groupe $\mathcal{G}$}. Du même, notons $\pi_{1}(\mathcal{G})$ le groupe d'automorphismes $\mathrm{Hom}_{\mathcal{G}}\big(\mathbb{1},\mathbb{1}\big)$. Le groupe $\pi_{1}(\mathcal{G})$ est dit \emph{le groupe fondamental du $2$}-\emph{groupe $\mathcal{G}$}.

\begin{lemme}
Le groupe fondamental $\pi_{1}(\mathcal{G})$ d'un $2$-groupe $\mathcal{G}$ est un groupe commutatif.
\end{lemme}
\begin{proof}
Si $f,g\colon\xymatrix@C-6pt{\mathbb{1}\ar[r]&\mathbb{1}}$ sont des morphismes d'un $2$-groupe $\mathcal{G}$ on a un diagramme commutatif:
$$
\xymatrix@C+8pt{
\mathbb{1}\otimes\mathbb{1} \ar[d]_{\mathbb{1}\otimes g} \ar[rd]|-{f\otimes g}  \ar[r]^{f\otimes\mathbb{1}}  & \mathbb{1}\otimes\mathbb{1}\ar[d]^-{\mathbb{1}\otimes g} \\
\mathbb{1}\otimes\mathbb{1}\ar[r]_{f\otimes\mathbb{1}} & \mathbb{1}\otimes\mathbb{1}
}
$$

D'un autre on a aussi les diagrammes commutatifs:
$$
\vcenter{\xymatrix@C+6pt{
\mathbb{1}\otimes\mathbb{1}\ar[r]^-{f\otimes\mathbb{1}} & \mathbb{1}\otimes\mathbb{1} \ar[r]^-{\mathbb{1}\otimes g} & \mathbb{1}\otimes\mathbb{1} \\
\mathbb{1}\ar[u]^-{r_{\mathbb{1}}} \ar[r]_{f}&  \mathbb{1} \ar[u]|-{r_\mathbb{1} = \ell_\mathbb{1}} \ar[r]_-{g} & \mathbb{1} \ar[u]_-{\ell_{\mathbb{1}}}
}}
\qquad \text{et} \qquad 
\vcenter{\xymatrix{
\mathbb{1}\otimes\mathbb{1}\ar[r]^-{\mathbb{1}\otimes g} & \mathbb{1}\otimes\mathbb{1} \ar[r]^-{f\otimes \mathbb{1}} & \mathbb{1}\otimes\mathbb{1} \\
\mathbb{1}\ar[u]^-{\ell_{\mathbb{1}}} \ar[r]_{g}&  \mathbb{1} \ar[u]|-{\ell_\mathbb{1} = r_\mathbb{1}} \ar[r]_-{f} & \mathbb{1} \ar[u]_-{r_{\mathbb{1}}}
}}
$$

Donc $g\circ f \, = \, \ell_\mathbb{1}^{-1} \circ (\mathbb{1}\otimes g) \circ (f\otimes \mathbb{1})\circ r_{\mathbb{1}} 
\, = \, r_\mathbb{1}^{-1} \circ (f \otimes g) \circ \ell_{\mathbb{1}}
\, = \, r_\mathbb{1}^{-1} \circ (f\otimes \mathbb{1}) \circ (\mathbb{1}\otimes g)\circ \ell_{\mathbb{1}} \, = \, f\circ g$.
\end{proof}

On définit aussi-tôt des foncteurs:
\begin{equation}
\xymatrix@C+30pt{\text{$2$-${\bf Grp}$}\ar[r]^-{\pi_{0}(\,\cdot\,),\,\pi_{1}(\,\cdot\,)} & {\bf Grp}}\,.
\end{equation}

Un morphisme de 2-groupes $\xymatrix@C-8pt{\mathcal{G}\ar[r]^-{F}&\mathcal{G}'}$ est dit une \emph{équivalence faible} s'il vérifie une des propriétés équivalentes du Lemme suivant:

\begin{lemme}\label{eq2grps} 
Si $\xymatrix@C-8pt{\mathcal{G}\ar[r]^-{F}&\mathcal{G}'}$ est un morphisme monoïdal lax et unitaire entre 2-groupes, alors les énoncés suivants sont équivalents:
\begin{enumerate}
\item Les morphismes de groupes $\pi_{0}F$ et $\pi_{1}F$ sont d'isomorphismes.
\item Le foncteur sous-jacent au morphisme de $2$-groupes $\xymatrix@C-8pt{\mathcal{G}\ar[r]^-{F}&\mathcal{G}'}$ est pleinement fidèle et essentiellement surjectif, \emph{i.e.} l'image de $F$ par le $2$-foncteur d'oubli $\xymatrix@-8pt{\text{$2$-$\underline{\bf Grp}$}\ar[r]&\underline{\bf cat}}$ est une équivalence interne de la $2$-catégorie $\underline{\bf cat}$ des petites catégories.  
\item $F$ est une équivalence interne de la $2$-catégorie $2$-$\underline{\bf Grp}$ c'est-à-dire, il existe un morphisme monoïdal lax et unitaire $\xymatrix@C-8pt{\mathcal{G}'\ar[r]^-{G}&\mathcal{G}}$ et de transformations $\xymatrix@C-15pt{F\,\circ\,G\;\ar@{=>}[r]^-{\alpha}& \; \mathrm{id}_{\mathcal{G}'}}$ et $\xymatrix@C-15pt{G\,\circ\,F\;\ar@{=>}[r]^-{\beta}& \; \mathrm{id}_{\mathcal{G}}}$ (note que toutes les $2$-flèches de $2$-$\underline{\bf Grp}$ sont des isomorphismes).  
\end{enumerate}
\end{lemme}
\begin{proof}
Premièrement remarquons que si $\xymatrix@C-8pt{\mathcal{G}\ar[r]^-{F}&\mathcal{G}'}$ est un morphisme monoïdal lax et unitaire entre 2-groupes, alors la fonction $\pi_{0}F$ est surjective si et seulement si le foncteur sous-jacent à $F$ est essentiellement surjectif. Donc pour vérifier l'équivalence (i) $\Leftrightarrow$ (ii) il suffit de montrer que $F$ est un foncteur pleinement fidèle si et seulement si la fonction $\pi_{0}F$ est injective et $\pi_{1}F$ est bijective. 

D'un côté si le foncteur $F$ est pleinement fidèle on a en particulier que la fonction $\pi_{1}F$ est bijective. En plus si $X$ et $Y$ sont des objets de $\mathcal{G}$ tels que $F(X)$ et $F(Y)$ représentent la même classe dans l'ensemble $\pi_0(\mathcal{G}')$, alors il existe un isomorphisme $\varphi:\xymatrix@C-10pt{F(X)\ar[r]&F(Y)}$ vu que $\mathcal{G}'$ est un groupoïde. On a en plus que $\varphi=F(f)$ pour un morphisme $f:\xymatrix@C-10pt{X\ar[r]&Y}$ de $\mathcal{G}$ parce que $F$ est plein; donc la fonction $\pi_{0}F$ est bien injective.

Réciproquement, supposons que la fonction $\pi_{0}F$ est injective et $\pi_{1}F$ est bijective. Si $X$ et $Y$ sont des objets de $\mathcal{G}$ tels que l'ensemble $\mathrm{Hom}_\mathcal{G}\big(X,Y\big)$ est vide, alors l'ensemble de morphismes $\mathrm{Hom}_{\mathcal{G}'}\big(FX,FY\big)$ est aussi vide parce que la fonction $\pi_{0}F$ est injective et $\mathcal{G}$ est un groupoïde. D'un autre, s'il existe un morphisme $f\colon\xymatrix@C-8pt{X\ar[r]&Y}$ on considère le carré commutatif:
$$
\xymatrix@R=8pt@C+10pt{\mathrm{Hom}_\mathcal{G}(X,Y) \ar[r]^-{F_{X,Y}} & \mathrm{Hom}_{\mathcal{G}'}\big(F(X),F(Y)\big) \\
\mathrm{Hom}_\mathcal{G}(X,X) \ar[r]_-{F_{X,X}} \ar[u]^-{f\circ -}_-{\cong}& \mathrm{Hom}_{\mathcal{G}'}\big(F(X),F(X)\big) \ar[u]_-{F(f)\circ -}^{\cong}}
$$

On déduit que la function $F_{X,Y}$ est bijective si et seulement si $F_{X,X}$ l'est, vu que $f$ est un isomorphisme.

Montrons alors que pour un objet quelconque $X$, la fonction $F_{X,X}$ es bijective si et seulement $\pi_{1}F$ l'est. Pour cela remarquons qu'on a un diagramme composé de deux carré (plutôt pentagones) commutatifs:
\begin{equation}\label{picardhom}
\def\objectstyle{\scriptstyle}
\def\labelstyle{\scriptstyle}
\xymatrix@R-2pt@C+25pt{
\mathrm{Hom}_\mathcal{G}\big( \mathbb{1},\mathbb{1}\big) \ar[r]^-{X\otimes -}\ar[dd]_-{F_{\mathbb{1},\mathbb{1}}}  & 
\mathrm{Hom}_\mathcal{G}\big( X\otimes\mathbb{1}, X\otimes\mathbb{1}\big) \ar[d]|-{F_{X\otimes\mathbb{1},X\otimes\mathbb{1}}} &
\mathrm{Hom}_\mathcal{G}\big(X,X\big)  \ar[l]_-{r_X\circ \,-\, \circ r_X^{-1}}\ar[dd]^-{F_{X,X}}\\
& \mathrm{Hom}_\mathcal{G}\big( F(X\otimes\mathbb{1}),F(X\otimes\mathbb{1})\big)\ar[d]|-{m_{X,\mathbb{1}} \circ F(-) \circ m^{-1}_{X,\mathbb{1}}}   &\\
\mathrm{Hom}_{\mathcal{G}'}\big( \mathbb{1},\mathbb{1}\big)\ar[r]_-{FX\otimes -} & 
\mathrm{Hom}_{\mathcal{G}'}\big( FX\otimes\mathbb{1}, FX\otimes\mathbb{1}\big)&  
\mathrm{Hom}_{\mathcal{G}'}\big( FX,FX\big) \ar[l]^-{l_{FX}\circ - \circ l_{FX}^{-1}}\,.}
\end{equation}

Ceux pentagones sont effectivement commutatifs car si $\varphi\colon\xymatrix@C-8pt{\mathbb{1}\ar[r]&\mathbb{1}}$ et $\psi\colon\xymatrix@C-8pt{X\ar[r]&X}$ sont de morphismes quelconques de $\mathcal{G}$ on a un carré commutatif dans $\mathcal{G}'$:
$$
\def\objectstyle{\scriptstyle}
\def\labelstyle{\scriptstyle}
\xymatrix{
FX \otimes \mathbb{1}    \ar[r]^-{m_{X,\mathbb{1}}} \ar[d]_-{FX\otimes F\varphi} &   F (X \otimes \mathbb{1}) \ar[d]^-{F(X\otimes \varphi)}\\
FX \otimes \mathbb{1}    \ar[r]_-{m_{X,\mathbb{1}}} &   F (X \otimes \mathbb{1})}
$$
et un prisme de faces commutatives:
$$
\def\objectstyle{\scriptstyle}
\def\labelstyle{\scriptstyle}
\xymatrix@R-5pt{ 
F(X)\otimes\mathbb{1}\ar[dd]_-F\ar[rd]|-{m_{X,\mathbb{1}}} &&F(X)\ar[ll]_-{l_{FX}}\ar[ld]|-{F(l_X)} \ar[dd]^-{F\psi} \\
&\ar[dd]|(.31){F(\psi\otimes\mathbb{1})} F(X\otimes\mathbb{1})& \\
F(X)\otimes\mathbb{1}\ar[rd]_-{m_{X,\mathbb{1}}}  &&F(X)\ar'[l][ll]|-{l_{FX}}\ar[ld]^-{F(l_X)} \\
&F(X\otimes\mathbb{1})&}
$$

Donc on déduit du diagramme \eqref{picardhom} que $F_{X,X}$ est un fonction bijective si $F_{\mathbb{1},\mathbb{1}}$ l'est, car les foncteurs $X\otimes -$ et $FX\otimes -$ sont pleinement fidèles. 

L'équivalence (ii) $\Leftrightarrow$ (iii) est une conséquence du Lemme \ref{equiequi}. 
\end{proof}

Posons $2$-$h{\bf Grp}$ pour noter la \emph{catégorie homotopique des} $2$-\emph{groupes} c'est-à-dire la catégorie qu'on obtient de la $2$-catégorie $2$-$\underline{\bf Grp}$ en imposant la condition:
$$
\mathrm{Hom}_{\text{$2$-$h{\bf Grp}$}} \, = \, \pi_0\big(\underline{\mathrm{Hom}}_{\text{$2$-${\bf Grp}$}}\big)
$$

Vu que les catégories des morphismes $\underline{\mathrm{Hom}}_{\text{$2$-${\bf Grp}$}}$ de la $2$-catégorie $2$-$\underline{\bf Grp}$ sont de groupoïdes, on obtient la catégorie $2$-$h{\bf Grp}$ de la catégorie des $2$-groupes $2$-${\bf Grp}$ en imposant dans son ensemble de morphismes la relation d'équivalence qu'identifie deux morphismes de $2$-groupes $F,G\colon\xymatrix@C-10pt{\G\ar[r]&\H}$ s'il existe une transformation $\eta\colon\xymatrix@C-13pt{F\ar@{=>}[r]&G}$.

Remarquons que si $\pi\colon\xymatrix@C-10pt{\text{$2$-${\bf Grp}$} \ar[r] & \text{$2$-$h{\bf Grp}$}}$ note le foncteur quotient canonique, il se suit du Lemme \ref{eq2grps} que l'image d'un morphisme de $2$-groupes $F$ par le foncteur $\pi$ est un isomorphisme de la catégorie homotopique des $2$-groupes si et seulement si $F$ est une équivalence faible.


\renewcommand{\thesubsection}{\S\thesection.\arabic{subsection}}
\subsection{}\;\label{foncteurGn}
\renewcommand{\thesubsection}{\thesection.\arabic{subsection}}

Notons $\underline{\bf cat}$ la $2$-catégorie des petites catégories, foncteurs et transformations naturelles. Rappelons que le $2$-foncteur des morphismes:
\begin{equation}\label{eexpo}
\xymatrix@R=5pt@C+25pt{\underline{\bf cat}^{op}\times \underline{\bf cat} \ar[r]^-{\underline{\mathrm{Hom}}_{\bf cat}} & \underline{\bf cat} \,, }
\end{equation}
est défini par la règle:
$$
\bigg(\vcenter{\xymatrix{
\mathcal{A}\rtwocell\omit{\alpha} & \ltwocell^G_F{\omit}\mathcal{B}
}}\,,\;\,
\vcenter{\xymatrix{\C\rtwocell^\varphi_\psi{\beta}&\D}}\bigg) \quad \mapsto \quad 
\Bigg(
\vcenter{\xymatrix@C-18pt{
\underline{\mathrm{Hom}}_{\bf cat}(\A,\C)\ar@<+3pt>@/^15pt/[rrrrr]^{\varphi\circ-\circ F}\ar@<-3pt>@/_15pt/[rrrrr]_{\psi\circ-\circ G}&\rrtwocell<\omit>{\qquad\alpha\star-\star\beta}&&&&\underline{\mathrm{Hom}}_{\bf cat}(\B,\D)}}
\Bigg)\,.
$$

Si $\mathcal{A}$ et $\mathcal{C}$ sont de petites catégories, on pose plus-tôt $\C^{\A}$ pour noter $\underline{\mathrm{Hom}}_{\bf cat}(\A,\C)$. Si $\C=\G$ est un $2$-groupe on va définir dans la catégorie des foncteurs $\mathcal{G}^\mathcal{A}$ une structure de $2$-groupe. Commençons par la structure de catégorie monoïdale qu'on obtient en composant les isomorphismes canoniques des catégories: 
$$
\star^{\A} \cong \star \,,\qquad\quad
\big(\mathcal{G}\times \mathcal{G}\big)^{\mathcal{A}} \cong \mathcal{G}^{\mathcal{A}} \times \mathcal{G}^{\mathcal{A}} \qquad \text{et} \qquad 
\big(\mathcal{G}\times \mathcal{G} \times \mathcal{G}\big)^{\mathcal{A}}  \cong \mathcal{G}^{\mathcal{A}} \times \mathcal{G}^{\mathcal{A}} \times \mathcal{G}^{\mathcal{A}} 
$$
$\big($le foncteur $\xymatrix{{\bf cat} \ar[r]^-{-^{\mathcal{A}}}  & {\bf cat}}$ commute aux petites limites$\big)$ avec les images par le $2$-foncteur $\mathcal{G} \, \mapsto \, \mathcal{G}^{\mathcal{A}}$ des diagrammes:
$$
\xymatrix@C-5pt{\mathcal{G}\times \mathcal{G}  \ar[r]^-{\otimes}& \mathcal{G}}\,,
\quad
\xymatrix@C-5pt{\star\ar[r]^-{\mathbb{1}} & \mathcal{G}}\,,
\quad 
\vcenter{\xymatrix@R=4pt@C-6pt{ 
&&\mathcal{G}\times \mathcal{G}\ar@/^4pt/[rd]^\otimes&\\ 
\mathcal{G}\times \mathcal{G}\times \mathcal{G} \ar@/^5pt/[urr]^-{(p_1\otimes p_2,p_3)} \ar@/_5pt/[drr]_-{(p_1,p_2\otimes p_3)} &\rrtwocell<\omit>{a}   &  & \mathcal{G}\\
&&\mathcal{G}\times \mathcal{G} \ar@/_4pt/[ru]_\otimes& }}
\quad\text{et}\quad
\vcenter{\xymatrix@R=10pt{
&\mathcal{G}\times \mathcal{G}\ar@/^4pt/[rd]^\otimes&\\ 
\mathcal{G}  \ar[rr]|-{\mathrm{id}}
\ar@/^3pt/[ur]^-{(\mathbb{1},\mathrm{id})} \ar@/_4pt/[dr]_-{(\mathrm{id},\mathbb{1})} \rrtwocell<\omit>{<3>r}&&\mathcal{G}\lltwocell<\omit>{<3>l}\,.\\
&\mathcal{G}\times\mathcal{G} \ar@/_4pt/[ru]_\otimes&}}
$$

Explicitement la multiplication et l'unité de $\G^\A$ sont les foncteurs composés:
$$
\vcenter{\xymatrix@R=8pt{\G^\A\times\G^\A \ar@{}[rrd]|(.45){=}\ar@{}[d]|-{\text{\rotatebox[origin=c]{90}{$\cong$}}}  \ar@<+2pt>@/^8pt/[rrd]^-{\otimes}&&\\(\G\times\G)^\A \ar[rr]_-{\otimes^\A} &&\G^\A}}
\qquad\text{et}\qquad
\vcenter{\xymatrix@R=8pt{\star \ar@{}[rrd]|(.45){=}\ar@{}[d]|-{\text{\rotatebox[origin=c]{90}{$\cong$}}} \ar@<+2pt>@/^8pt/[rrd]^{\mathbb{1}}&&\\\star^\A\ar[rr]_-{\mathbb{1}^\A}&& \G^\A\,,}}
$$
et les contraintes d'associativité et d'unité sont les isomorphismes naturels:
$$
\xymatrix@C-15pt@R=8pt{
\ar@{}[rddddd]_-{=}&\G^\A\times\G^\A\times\G^\A\ar@{}[rrrrrd]|(.4){=}\ar@{}[d]|-{\text{\rotatebox[origin=c]{90}{$\cong$}}} \ar@<-25pt>@/_22pt/[ddddd]_{(p_1,p_2\otimes p_3)}   \ar@<+2pt>@/^17pt/[rrrrrd]^-{(p_1\otimes p_2,p_3)}&&&&&\\&
(\G\times\G\times\G)^{\A} \ar[ddd]|-{(p_1,p_2\otimes p_3)^\A}\ar[rrrr]|-{(p_1\otimes p_2,p_3)^{\A}} &&&& (\G\times\G)^{\A} \ar@{}[rddd]|(.45){=} \ar[ddd]|-{\otimes^{\A}} \ar@{}[r]|-\cong& 
\G^\A\times\G^\A\ar@/^12pt/[dddl]^-{\otimes}\\ 
&&\drtwocell<\omit>{\quad a^\A}&&&\\&&&&&&\\
&(\G\times\G)^{\A}\ar[rrrr]|-{\otimes^{\A}} \ar@{}[rrrd]|(.6){=} \ar@{}[d]|-{\text{\rotatebox[origin=c]{90}{$\cong$}}} &&&& \G^{\A} &\\
&\G^\A\times\G^\A \ar@/_15pt/[rrrru]_{\otimes} &&&& &}
$$
$$
\text{et}\qquad\qquad\qquad\vcenter{\xymatrix@C+10pt@R=5pt{
\ar@{}[rdd]_(.45){=}&\G^\A\times \G^\A \ar@<+5pt>@/^15pt/[rddd]^-{\otimes}\ar@{}[d]|-{\text{\rotatebox[origin=c]{90}{$\cong$}}} &\ar@{}[ldd]^(.45){=}\\
& (\G\times \G)^\A \ar@<+2pt>@/^2pt/[rdd]|-{\otimes^\A} &\\ &&\\
\G^\A \rrtwocell<\omit>{<3.5>\quad r^\A}
\ar@<+2pt>@/^2pt/[uur]|-{(\G,\mathbb{1})^\A}  \ar@<-2pt>@/_2pt/[ddr]|-{(\mathbb{1},\G)^\A} \ar@<+5pt>@/^15pt/[uuur]^{(\G^\A,\mathbb{1})} \ar@<-5pt>@/_15pt/[dddr]_{(\mathbb{1},\G^\A)}\ar[rr]|-{\text{identité}} &  & \G^{\A}\lltwocell<\omit>{<3.5>l^\A\quad}\,,\\&&\\
&(\G\times\G)^{\A}\ar@{}[d]|-{\text{\rotatebox[origin=c]{90}{$\cong$}}}\ar@<-2pt>@/_2pt/[ruu]|-{\otimes^\A}&\\
\ar@{}[ruu]^(.45){=}&\G^\A\times\G^\A\ar@<-5pt>@/_15pt/[ruuu]_-{\otimes}&\ar@{}[luu]_(.45){=}}}
$$
c'est-à-dire ce sont définis argument par argument. 

Il se suit qu'ils vérifient \eqref{pentagon} et \eqref{triangle} \emph{i.e.} étant donnés $\mathcal{X}$, $\mathcal{Y}$ et $\mathcal{Z}$ des objets de $\G^\A$ on a que:
$$
\def\objectstyle{\scriptstyle}
\def\labelstyle{\scriptstyle}
\vcenter{\xymatrix@C+18pt{
\big((\mathcal{W}\otimes \mathcal{X})\otimes \mathcal{Y}\big)\otimes \mathcal{Z}
\ar[r]^-{a_{\mathcal{W}\otimes \mathcal{X}, \mathcal{Y},\mathcal{Z}}}\ar[d]_-{a_{\mathcal{W},\mathcal{X},\mathcal{Y}}\otimes \mathcal{Z}}& 
(\mathcal{W}\otimes \mathcal{X})\otimes (\mathcal{Y}\otimes \mathcal{Z})  \ar[r]^-{a_{\mathcal{W},\mathcal{X},\mathcal{Y}\otimes \mathcal{Z}}}& \mathcal{W}\otimes\big(\mathcal{X}\otimes (\mathcal{Y}\otimes \mathcal{Z})\big)\\
\big(\mathcal{W} \otimes (\mathcal{X}\otimes \mathcal{Y})\big)\otimes \mathcal{Z} \ar[rr]_-{a_{\mathcal{W},\mathcal{X}\otimes \mathcal{Y},\mathcal{Z}}}&&\mathcal{W}\otimes\big((\mathcal{X}\otimes \mathcal{Y})\otimes \mathcal{Z}\big)\ar[u]_-{\mathcal{W}\otimes a_{\mathcal{X},\mathcal{Y},\mathcal{Z}}}}}
$$
$$
\text{et}\qquad
\def\objectstyle{\scriptstyle}
\def\labelstyle{\scriptstyle}
\vcenter{\xymatrix{
(\mathcal{X}\otimes \mathbb{1}) \otimes \mathcal{Y}  \ar[rr]^-{a_{\mathcal{X},\mathbb{1},\mathcal{Y}}} &&    \mathcal{X}\otimes (\mathbb{1}\otimes \mathcal{Y})\\
&\mathcal{X}\otimes \mathcal{Y}\ar[ru]_-{\mathcal{X}\otimes l_{\mathcal{Y}}}\ar[lu]^-{r_{\mathcal{X}\otimes \mathcal{Y}}}&}}
$$
sont de diagrammes commutatifs de $\G^\A$, parce que:
$$
\def\objectstyle{\scriptstyle}
\def\labelstyle{\scriptstyle}
\vcenter{\xymatrix@C+18pt{
\big((\mathcal{W}_a\otimes \mathcal{X}_a)\otimes \mathcal{Y}_a\big)\otimes \mathcal{Z}_a
\ar[r]^-{a_{\mathcal{W}_a\otimes \mathcal{X}_a, \mathcal{Y}_a,\mathcal{Z}_a}}\ar[d]_-{a_{\mathcal{W}_a,\mathcal{X}_a,\mathcal{Y}_a}\otimes \mathcal{Z}_a}& 
(\mathcal{W}_a\otimes \mathcal{X}_a)\otimes (\mathcal{Y}_a\otimes \mathcal{Z}_a)  \ar[r]^-{a_{\mathcal{W}_a,\mathcal{X}_a,\mathcal{Y}_a\otimes \mathcal{Z}_a}}& \mathcal{W}_a\otimes\big(\mathcal{X}_a\otimes (\mathcal{Y}_a\otimes \mathcal{Z}_a)\big)\\
\big(\mathcal{W}_a \otimes (\mathcal{X}_a\otimes \mathcal{Y}_a)\big)\otimes \mathcal{Z}_a \ar[rr]_-{a_{\mathcal{W}_a,\mathcal{X}_a\otimes \mathcal{Y}_a,\mathcal{Z}_a}}&&\mathcal{W}_a\otimes\big((\mathcal{X}_a\otimes \mathcal{Y}_a)\otimes \mathcal{Z}_a\big)\ar[u]_-{\mathcal{W}_a\otimes a_{\mathcal{X}_a,\mathcal{Y}_a,\mathcal{Z}_a}}}}
$$
$$
\text{et}\qquad
\def\objectstyle{\scriptstyle}
\def\labelstyle{\scriptstyle}
\vcenter{\xymatrix{
(\mathcal{X}_a\otimes \mathbb{1}) \otimes \mathcal{Y}_a  \ar[rr]^-{a_{\mathcal{X}_a,\mathbb{1},\mathcal{Y}_a}} &&    \mathcal{X}_a\otimes (\mathbb{1}\otimes \mathcal{Y}_a)\\
&\mathcal{X}_a\otimes \mathcal{Y}_a\ar[ru]_-{\mathcal{X}_a\otimes l_{\mathcal{Y}_a}}\ar[lu]^-{r_{\mathcal{X}_a\otimes \mathcal{Y}_a}}&}}
$$
sont de diagrammes commutatifs de $\G$ pour tout objet $a$ de $\A$.

D'un autre vu qu'on a supposé que $\mathcal{G}$ est un $2$-groupe, on trouve que la catégorie des foncteurs $\mathcal{G}^\mathcal{A}$ est un groupoïde et d'après le Lemme \ref{inverses} il est un $2$-groupe car on a des foncteurs:
$$
\xymatrix@C+10pt{\mathcal{G}^\A \ar[r]^-{(\iota^d)^\A}& \mathcal{G}^\A} 
\qquad \quad \text{et} \qquad \quad 
\xymatrix@C+10pt{\mathcal{G}^\A \ar[r]^-{(\iota^g)^\A}& \mathcal{G}^\A\,,}
$$
et des contraintes des inverses:
$$
\vcenter{\xymatrix@C+10pt@R=5pt{
\ar@{}[rdd]_(.45){=}&\G^\A\times \G^\A \ar@<+5pt>@/^15pt/[rddd]^-{\otimes}\ar@{}[d]|-{\text{\rotatebox[origin=c]{90}{$\cong$}}} &\ar@{}[ldd]^(.45){=}\\
& (\G\times \G)^\A \ar@<+2pt>@/^2pt/[rdd]|-{\otimes^\A} &\\ &&\\
\G^\A \rrtwocell<\omit>{<-3>\quad \alpha^\A}
\ar@<+2pt>@/^2pt/[uur]|-{(\G,\iota^d)^\A}  \ar@<-2pt>@/_2pt/[ddr]|-{(\iota^g,\G)^\A} \ar@<+5pt>@/^15pt/[uuur]^-{(\G^\A,\iota^d)} \ar@<-5pt>@/_15pt/[dddr]_{(\iota^g,\G^\A)}\ar[rr]|-{\text{identité}} &  & \G^{\A}\lltwocell<\omit>{<-3>\beta^\A\quad}\,.\\&&\\
&(\G\times\G)^{\A}\ar@{}[d]|-{\text{\rotatebox[origin=c]{90}{$\cong$}}}\ar@<-2pt>@/_2pt/[ruu]|-{\otimes^\A}&\\
\ar@{}[ruu]^(.45){=}&\G^\A\times\G^\A\ar@<-5pt>@/_15pt/[ruuu]_-{\otimes}&\ar@{}[luu]_(.45){=}}}
$$

On va étendre la fonction $(\A,\G)\longmapsto \G^\A$ qu'on vient de définir à un $2$-foncteur qui complète le carré commutatif:
\begin{equation}\label{homodeG}
\xymatrix@R=1pt@C+10pt{
(\A,\G)\ar@{}[r]|-{\longmapsto}&\G^\A\\
\underline{\bf cat}^{op}\times \text{$2$-$\underline{\bf Grp}$} \ar@{-->}[r] \ar[ddddd]_-{\underline{\bf cat}^{op}\times \pi}& \text{$2$-$\underline{\bf Grp}$} \ar[ddddd]^-{\pi}\\&\\&\\&\\&\\
\underline{\bf cat}^{op}\times \underline{\bf cat} \ar[r]_-{\underline{\mathrm{Hom}}_{\bf cat}} & \underline{\bf cat}\,,}
\end{equation}
où $\pi\colon\xymatrix@C-5pt{\text{$2$-$\underline{\bf Grp}$} \ar[r] & \underline{\bf cat}}$ est le $2$-foncteur d'oubli.

Si $F\colon\xymatrix@C-4pt{\B\ar[r]&\A}$ est un foncteur et $(\varphi,m^\varphi)\colon\xymatrix@C-4pt{\G\ar[r]&\mathcal{H}}$ est un morphisme de $2$-groupes, on définit le morphisme de $2$-groupes $(\varphi^F,m^{\varphi^F})\colon\xymatrix@C-4pt{\G^\A\ar[r]&\mathcal{H}^\B}$ par le foncteur $\varphi^F=\varphi \circ - \circ F$ et l'isomorphisme naturelle:
$$
\vcenter{\xymatrix@C=20pt@R=20pt{
\G^\A\times\G^\A \ar[d]_-{\varphi^F\times\varphi^F}\ar[r]^-{\otimes} \drtwocell<\omit>{\; m^{\varphi^F}} & \G^\A \ar[d]^{\varphi^F}\\
\H^\B\times\H^\B\ar[r]_{\otimes}& \H^\B
}}\qquad = \qquad
\vcenter{\xymatrix@C=5pt@R+5pt{
\G^\A\times\G^\A \ar@{}[r]|-{\cong} \ar[dd]_-{\varphi^F\times\varphi^F}\ar@/^20pt/[rrrr]^{\otimes}  &  (\G\times\G)^\A  \ar[rrr]|-{\otimes^\A} \ar[d]|-{(\varphi\times\varphi)^\A} & \drtwocell<\omit>{(m^\varphi)^\A}& & \G^\A  \ar[d]^-{\varphi^\A} \ar@<+3pt>@/^18pt/[dd]^-{\varphi^F}\\
 & (\mathcal{H}\times\mathcal{H})^\A \ar[rrr]|-{\otimes^\A}\ar[d]|-{(\mathcal{H}\times\mathcal{H})^\F}&&& \mathcal{H}^\A \ar[d]^{\mathcal{H}^F}\\ 
 \mathcal{H}^\B\times\mathcal{H}^\B \ar@{}[r]|-{\cong}\ar@/_20pt/[rrrr]_-{\otimes} & (\mathcal{H}\times\mathcal{H})^\B\ar[rrr]|-{\otimes^\B}&&& \mathcal{H}^\B\,.}}
$$ 

La couple $(\varphi^F,m^{\varphi^F})$ est effectivement un morphisme lax et unitaire de $2$-groupes parce que: 
$$
\varphi^F(\,\text{foncteur contant à valeurs $\mathbb{1}_\G$}\,)\;=\;\text{foncteur contant à valeurs $\mathbb{1}_\H$}
$$ 
et parce qu'on a des diagrammes commutatifs de $\G$:
\begin{equation*}
\def\objectstyle{\scriptstyle}
\def\labelstyle{\scriptstyle}
\vcenter{\xymatrix@C+18pt@R-3pt{
\varphi\big((\mathcal{X}_{Fb}\otimes \mathcal{Y}_{Fb})\otimes \mathcal{Z}_{Fb}\big)\ar[d]_-{\varphi a}&\varphi (\mathcal{X}_{Fb}\otimes \mathcal{Y}_{Fb})\otimes \varphi \mathcal{Z}_{Fb}\ar[l]_-{m^\varphi}&(\varphi \mathcal{X}_{Fb}\otimes \varphi \mathcal{Y}_{Fb})\otimes \varphi \mathcal{Z}_{Fb}\ar[d]^-{a}\ar[l]_-{m^\varphi\otimes \varphi \mathcal{Z}_{Fb}}\\
\varphi\big(\mathcal{X}_{Fb}\otimes (\mathcal{Y}_{Fb}\otimes \mathcal{Z}_{Fb})\big)  &\varphi\mathcal{X}_{Fb}\otimes \varphi (\mathcal{Y}_{Fb}\otimes \mathcal{Z}_{Fb})\ar[l]^-{m^\varphi}&\varphi \mathcal{X}_{Fb}\otimes (\varphi \mathcal{Y}_{Fb}\otimes \varphi \mathcal{Z}_{Fb})\ar[l]^-{\varphi \mathcal{X}_{Fb}\otimes m^\varphi}}}
\end{equation*}
$$\text{et}$$
\begin{equation*}
\def\objectstyle{\scriptstyle}
\def\labelstyle{\scriptstyle}
\vcenter{\xymatrix@C-3pt@R-5pt{
\mathbb{1}\otimes \varphi \mathcal{X}_{Fb} \ar[dd]_-{m^\varphi_{\mathbb{1},\mathcal{X}_{Fb}}}&&\varphi \mathcal{X}_{Fb}\otimes \mathbb{1}\ar[dd]^-{m^\varphi_{\mathcal{X}_{Fb},\mathbb{1}}} \\
&\varphi \mathcal{X}_{Fb}\ar[lu]_{l_{\varphi \mathcal{X}_{Fb}}}\ar[ru]^{r_{\varphi \mathcal{X}_{Fb}}}\ar[rd]_{\varphi r_{\mathcal{X}_{Fb}}}\ar[ld]^-{\varphi l_{\mathcal{X}_{Fb}}}&\\
 \varphi(\mathbb{1}\otimes \mathcal{X}_{Fb})&&\varphi(\mathcal{X}_{Fb}\otimes\mathbb{1})
 }}
\end{equation*}
si $\mathcal{X},\mathcal{Y},\mathcal{Z}$ sont des objets de $\G^{\A}$ et $b$ un objet de $\B$.

Remarquons aussi que si on a de foncteurs: 
$$
\xymatrix@C+5pt{\C\ar[r]^-{G} & \B \ar[r]^-{F} & \A }
$$ 
et de morphismes de $2$-groupes:
$$
\xymatrix@C+10pt{\G\ar[r]^-{(\varphi,m^\varphi)}&\mathcal{H} \ar[r]^-{(\psi,m^\psi)} & \mathcal{K}\,,}
$$
alors $\big(\psi^G,m^{\psi^G}\big)\circ\big(\varphi^F,m^{\varphi^F}\big)=\big((\psi\circ\varphi)^{F\circ G},m^{(\psi\circ \varphi)^{F\circ G} }\big)$. 

En effet, on voit facilement que $(\psi\circ\varphi)^{F\circ G}=\psi^G\circ\varphi^F$ et qu'on a un triangle commutatif:
$$
\xymatrix{
(\psi\circ\varphi)^{F\circ G}(\mathcal{X})\otimes (\psi\circ\varphi)^{F\circ G}(\mathcal{Y}) 
\ar[rr]^-{m^{(\psi\circ\varphi)^{F\circ G}}_{\mathcal{X},\mathcal{Y}}} \ar[rd]_-{m^{\psi^G}_{\varphi^F\mathcal{X},\varphi^F\mathcal{Y}}} &  & 
(\psi\circ\varphi)^{F\circ G}(\mathcal{X}\otimes \mathcal{Y})\\
& \psi^G\big(\varphi^F(\mathcal{X})\otimes\varphi^F(\mathcal{Y})\big)
\ar[ru]_-{\psi^G(m^{\varphi^F}_{\mathcal{X},\mathcal{Y}})}& }
$$
pour $\mathcal{X}$ et $\mathcal{Y}$ des objets de $\G^\A$; parce que si $c$ est un objet de $\C$ on a un triangle commutatif:
$$
\xymatrix{
(\psi\circ\varphi)(\mathcal{X}_{FGc})\otimes (\psi\circ\varphi)(\mathcal{Y}_{FGc}) 
\ar[rr]^-{m^{(\psi\circ\varphi)}_{\mathcal{X}_{FGc},\mathcal{Y}_{FGc}}} \ar[rd]_-{m^{\psi}_{\varphi \mathcal{X}_{FGc},\varphi\mathcal{Y}_{FGc}}} &  & 
(\psi\circ\varphi)(\mathcal{X}_{FGc}\otimes \mathcal{Y}_{FGc})\\
& \psi\big(\varphi(\mathcal{X}_{FGc})\otimes\varphi(\mathcal{Y}_{FGc})\big)
\ar[ru]_-{\psi(m^{\varphi}_{\mathcal{X}_{FGc},\mathcal{Y}_{FGc}})}& }
$$

Finalement considérons une transformation naturelle entre foncteurs de petites catégories et une transformation entre morphismes de $2$-groupes: 
$$
\xymatrix@C+12pt{\mathcal{B}\rtwocell<5>^F_G{\, \alpha} &\mathcal{A}}
\qquad\qquad \text{et}\qquad\qquad 
\xymatrix@C+12pt{\mathcal{G}\rtwocell<5>^{(\varphi,m^\varphi)}_{(\psi,m^\psi)}{\, \eta} &\mathcal{H}\,,}
$$ 
respectivement; on définit la transformation $\xymatrix@C+16pt{\G^\mathcal{A}\rtwocell<5>^{(\varphi^F,m^{\varphi^F})}_{(\psi^G,m^{\psi^G})}{\;\;\;\eta^\alpha} &\G^\mathcal{B}}$ par le composé horizontal:
$$
\xymatrix@C+12pt{\G^\mathcal{A}\rtwocell<5>^{\varphi^\A}_{\psi^\A}{\;\,\; \eta^\A} &\H^\mathcal{A} \rtwocell<5>^{\H^F}_{\H^G}{\;\;\, \H^\alpha} &\H^\mathcal{B}}
\qquad = \qquad 
\xymatrix@C+12pt{\G^\mathcal{A}\rtwocell<5>^{\G^F}_{\G^G}{\;\;\, \G^\alpha} &\G^\mathcal{B} \rtwocell<5>^{\varphi^\B}_{\psi^\B}{\;\;\, \eta^\B} &\H^\mathcal{B}\,.}
$$

On vérifie que si $\mathcal{X}$ et $\mathcal{Y}$ sont des objets de $\G^\A$ alors:
\begin{equation*}
\xymatrix@C+25pt@R+25pt{
\varphi^F\mathcal{X}\otimes \varphi^F\mathcal{Y}\ar[r]^-{\eta^\alpha_{\mathcal{X}}\otimes \eta^\alpha_{\mathcal{Y}}} \ar[d]_{m^{\varphi^F}_{\mathcal{X},\mathcal{Y}}}&
\psi^F\mathcal{X}\otimes \psi^F\mathcal{Y}\ar[d]^{m^{\psi^G}_{\mathcal{X},\mathcal{Y}}}\\
\varphi^F(\mathcal{X}\otimes \mathcal{Y}) \ar[r]_-{\eta^{\alpha}_{\mathcal{X}\otimes\mathcal{Y}}}  & 
\psi^G(\mathcal{X}\otimes \mathcal{Y})  }
\end{equation*}
est un diagramme commutatif de $\H^\B$ de la façon suivante: Pour un objet quelconque $a$ de $\A$ on décompose le diagramme:
\begin{equation*}
\xymatrix@C+30pt@R+25pt{
\varphi\mathcal{X}_{Fa}\otimes \varphi\mathcal{Y}_{Fa}\ar[r]^-{(\eta^\alpha_{\mathcal{X}})_a\otimes (\eta^\alpha_{\mathcal{Y}})_a} \ar[d]_{m^{\varphi}_{\mathcal{X}_{Fa},\mathcal{Y}_{Fa}}}&\psi\mathcal{X}_{Fa}\otimes \psi\mathcal{Y}_{Fa}\ar[d]^{m^{\psi}_{\mathcal{X}_{Ga},\mathcal{Y}_{Ga}}}\\
\varphi(\mathcal{X}_{Fa}\otimes \mathcal{Y}_{Fa}) \ar[r]_-{(\eta^{\alpha}_{\mathcal{X}\otimes\mathcal{Y}})_a}  & 
\psi(\mathcal{X}_{Ga}\otimes \mathcal{Y}_{Ga})  }
\end{equation*}
par les sous-diagrammes:
$$
\xymatrix@C+20pt{
\varphi (\mathcal{X}_{Fa}) \otimes \varphi (\mathcal{Y}_{Fa})\ar[rd]|-{\varphi (\mathcal{X}_{\alpha}) \otimes \varphi (\mathcal{Y}_{\alpha})} 
\ar[dddd]_-{m^\varphi_{\mathcal{X}_{Fa},\mathcal{Y}_{Fa}}}
\ar[rr]^-{(\eta^\alpha_\mathcal{X})_a\otimes (\eta^\alpha_\mathcal{Y})_a}\ar@{}[ddddr]|-{\text{(III)}}&\ar@{}[ddddr]|-{\text{(IV)}}\ar@{}[d]|-{\text{(I)}}&
\psi (\mathcal{X}_{Ga})\otimes \psi (\mathcal{X}_{Ga})\ar[dddd]^-{m^\psi_{\mathcal{X}_{Ga},\mathcal{Y}_{Ga}}}\\
&  \varphi (\mathcal{X}_{Ga}) \otimes \varphi (\mathcal{Y}_{Ga})\ar[ru]|-{\eta_{\mathcal{X}_{Ga}}\otimes \eta_{\mathcal{Y}_{Ga}}} \ar[dd]|-{m^\varphi_{\mathcal{X}_{Ga},\mathcal{Y}_{Ga}}}&\\&&\\& \varphi (\mathcal{X}_{Ga}\otimes \mathcal{Y}_{Ga}) \ar[rd]|-{\eta_{\mathcal{X}_{Ga} \otimes \mathcal{Y}_{Ga}}} \ar@{}[d]|-{\text{(II)}}&\\
\varphi (\mathcal{X}_{Fa}\otimes \mathcal{Y}_{Fa})\ar[rr]_-{(\eta^{\alpha}_{\mathcal{X}\otimes\mathcal{Y}})_a} \ar[ru]|-{\varphi(\mathcal{X}_\alpha\otimes\mathcal{Y}_\alpha)}
&&\psi (\mathcal{X}_{Ga}\otimes \mathcal{Y}_{Ga}) }
$$

Puis on note que (I) et (II) sont commutatifs par définition de $\eta^\alpha$, (III) est commutatif parce que $m^\varphi$ est une transformation naturelle de foncteurs et (IV) commute vu que $\eta$ est une transformation de morphismes de $2$-groupes.

\renewcommand{\thesubsection}{\S\thesection.\arabic{subsection}}
\subsubsection{}\;
\renewcommand{\thesubsection}{\thesection.\arabic{subsection}}

Posons $\underline{\Delta}$ pour noter la sous-$2$-catégorie pleine de $\underline{\bf cat}$ dont les objets sont les catégories $[n]=\{0<\dots<n\}$ pour $n\geq 0$. Autrement dit, $\underline{\Delta}$ est la $2$-catégorie dont la catégorie sous-jacente est la catégories des simplexes $\Delta$ et où il existe une seul $2$-flèche $\xymatrix@C+5pt{[n]\rtwocell^{\varphi}_{\psi}&[m]}$ si et seulement si $\varphi(i)\leq \psi(i)$ pour tout $0\leq i\leq n$.

Si $\G$ est un $2$-groupe, considérons le $2$-foncteur induit du $2$-foncteur \eqref{homodeG}:
\begin{equation}\label{2grpsimpi}
\xymatrix@R=5pt{\underline{\Delta}^{op}\ar[r]^-{\mathcal{G}^{\bullet}} & \text{$2$-$\underline{\bf Grp}$} \\ [n]\ar@{}[r]|-{\longmapsto} &\mathcal{G}^{[n]}}\,.
\end{equation}

Montrons que l'objet simplicial $\G^{[\bullet]}\colon\xymatrix@C-5pt{\Delta^{op}\ar[r]&\text{$2$-${\bf Grp}$}}$ sous-jacent au $2$-foncteur \eqref{2grpsimpi} est une "résolution simpliciale" du $2$-groupe $\G$. Plus précisément:

\begin{lemme}\label{porporfin}
Si $\G$ est un $2$-groupe et $\xymatrix@C-10pt{[n]\ar[r]^{\varphi}&[m]}$ est un morphisme quelconque de la catégorie des simplexes $\Delta$, le morphisme induit $\xymatrix@C-3pt{\mathcal{G}^{[m]}\ar[r]^-{\varphi^\star}&\mathcal{G}^{[n]}}$ est une équivalence faible de $2$-groupes (voir le Lemme \ref{eq2grps}).
\end{lemme}
\begin{proof}
On sait que tout morphisme $\xymatrix@C-10pt{[n]\ar[r]^{\varphi}&[m]}$ de la catégorie $\Delta$ est égal à un composé de la forme $\varphi = \delta_{j_{m-n+k}}\circ\dots \circ \delta_{j_1}\circ\sigma_{i_k}\circ\dots\circ\sigma_{i_1}$ où  $\delta_{j}$ et $\sigma_{i}$ notent les morphismes faces et dégénérescences. En particulier, il suffit de vérifier que pour tout $0\leq i\leq k$ et $0\leq j \leq k+1$ les morphismes $[k]\xymatrix@C-3pt{\ar@/^3pt/@<+3pt>[r]^-{\delta_j}&\ar@/^3pt/@<+3pt>[l]^-{\sigma_{i}}}[k+1]$ induisent des équivalences faibles des $2$-groupes $\mathcal{G}^{[k+1]}\xymatrix@C-3pt{\ar@/^3pt/@<+3pt>[r]^-{\delta^\star_j}&\ar@/^3pt/@<+3pt>[l]^-{\sigma^\star_{i}}}\mathcal{G}^{[k]}$.

Remarquons que si $0\leq i\leq k$ alors on a une égalité $\sigma_i\circ\delta_{i} = \mathrm{id}_{[k]}$. D'un autre, il y a une $2$-flèche $\mathrm{id}_{[k+1]}\Rightarrow \delta_{i}\circ\sigma_i$ de la $2$-catégorie $\underline{\Delta}$ pour $0\leq i\leq k$, vu que $a\leq \delta_{i}\circ\sigma_i(a)$ si $0\leq a\leq n+1$; en effet on note que:
$$
\delta_{i}\circ\sigma_i(a)=\begin{cases}a\quad &a\neq i \\ i+1 \quad &a=i \,.\end{cases}
$$ 

On en conclut que $\delta^\star_i\circ\sigma^{\star}_{i} = \mathrm{id}_{\mathcal{G}^{[k]}}$ et qu'il existe une $2$-flèche $\mathrm{id}_{\mathcal{G}^{[k+1]}}\Rightarrow \sigma^\star_{i}\circ\delta^\star_i$ de la $2$-catégorie $2$-$\underline{\bf Grp}$. Donc de la propriété (iii) du Lemme \ref{eq2grps} il se suit que $\mathcal{G}^{[k+1]}\xymatrix@C-3pt{\ar@/^3pt/@<+3pt>[r]^-{\delta^\star_j}&\ar@/^3pt/@<+3pt>[l]^-{\sigma^\star_{i}}}\mathcal{G}^{[k]}$ sont bien des équivalences faibles des $2$-groupes si $0\leq i,j\leq k$.

Finalement pour montrer que $\xymatrix{\mathcal{G}^{[k+1]} \ar[r]^-{\delta_{k+1}^\star} & \mathcal{G}^{[k]}}$ est une équivalence faible de $2$-groupes, on note de la même façon que $\sigma_k\circ\delta_{k+1} = \mathrm{id}_{[k]}$ et qu'il y a une $2$-flèche $\delta_{k+1}\circ\sigma_k\Rightarrow \mathrm{id}_{[k+1]} $ de la $2$-catégorie $\underline{\Delta}$ vu que: 
$$
\delta_{k+1}\circ\sigma_k(a)=\begin{cases}a\quad &0 \leq a\leq k \\ k \quad &a=k+1\,. \end{cases}
$$
\end{proof}

\section{Le nerf des $2$-groupes}


Si $\M$ est une catégorie monoïdale et $q\geq 0$ est un entier, un \emph{$q$-simplexe de $\M$} est par définition un couple $(X,\alpha)$ où $X$ est un ensemble d'objets de $\M$:
$$
X \, = \, \Big\{  X_{ij}  \,\Big|\,  0\leq i<j\leq q \Big\}
$$
et $\alpha$ est un ensemble de morphismes de $\M$:
$$
\alpha \, = \, \Big\{  \alpha_{ijk}\colon\xymatrix@-10pt{X_{ij}\otimes X_{jk} \ar[r]& X_{ik}}  \; \Big| \;  0\leq i<j<k\leq q \Big\};
$$
tels que si $0\leq i<j<k<l\leq q$ on a un diagramme commutatif:
\begin{equation}\label{ccohecohe}
\xymatrix@R=3pt@C+15pt{
X_{il}& \ar[l]_-{\alpha_{ijl}} X_{ij}\otimes X_{jl}  \\&\\&\\
 &  X_{ij}\otimes \big( X_{jk}\otimes X_{kl} \big) \ar[uuu]_-{X_{ij}\otimes\alpha_{jkl}} \ar@{}[d]_-{a}|-{\text{\rotatebox[origin=c]{90}{$\cong$}}} \\ 
X_{ik}\otimes X_{kl} \ar[uuuu]^-{\alpha_{ikl}}  & \big(X_{ij}\otimes X_{jk}\big)\otimes X_{kl}  \ar[l]^-{\alpha_{ijk} \otimes X_{kl}} \,.}
\end{equation}

Posons $\mathcal{\M}_q$ pour note l'ensemble des $q$-simplexes de $\M$. On définit un foncteur:
\begin{equation}\label{nervio22}
\xymatrix@R=3pt{{\bf cat}_{lax,\star}^{\otimes}\times\Delta^{op} \ar[rr] && {\bf Ens}\\ \big(\M,[q]\big) \ar@{}[rr]|-{\longmapsto} && \M_q}
\end{equation}
comme suit: Si $\varphi\colon\xymatrix@-5pt{[q]\ar[r]&[q']}$ est un morphisme de $\Delta$, la fonction induite:
$$
\xymatrix@R=3pt{
\M_{q'} \ar[rr]^-{\varphi^*} &&\M_{q}\\
(X,\alpha) \ar@{}[rr]|-{\longmapsto}&& (Y,\beta)}
$$
est donné par les règles:
\begin{equation}\label{lasyy}
Y_{ij} \, = \, \begin{cases}
\mathbb{1}_{\G} & \text{si \; $\varphi i=\varphi j$} \\
X_{\varphi i\varphi j}& \text{si \; $\varphi i<\varphi j$} 
\end{cases}
\end{equation}
toujours que $0\leq i< j \leq q$ et:
\begin{equation}\label{lasbebe}
\beta_{ijk} \, = \, \begin{cases}
\ell^{-1}_{X_{\varphi i\varphi k}}\colon\xymatrix@-10pt{\mathbb{1} \otimes X_{\varphi i\varphi k} \ar[r]& X_{\varphi i\varphi k} } & \text{si \; $\varphi i=\varphi j\leq\varphi k$}\\
r^{-1}_{X_{\varphi i\varphi k}}\colon\xymatrix@-10pt{X_{\varphi i\varphi k}\otimes \mathbb{1}\ar[r]& X_{\varphi i\varphi k}}  & \text{si \; $\varphi i\leq\varphi j=\varphi k$}\\
\alpha_{\varphi i\varphi j\varphi k}\colon\xymatrix@-10pt{X_{\varphi i\varphi j}\otimes X_{\varphi j\varphi k}\ar[r]& X_{\varphi i\varphi k}} & \text{si \; $\varphi i<\varphi j<\varphi k$}
\end{cases}
\end{equation}
pour $0\leq i< j < k \leq q'$.

Il faut vérifier que la couple $(Y,\beta)$ ainsi définie est un $q$-simplexe de $\M$, \emph{i.e.} il faut montrer qu'avec les définition \eqref{lasyy} et \eqref{lasbebe} on a un diagramme commutatif:
\begin{equation}\label{ccohecohe2}
\xymatrix@R=3pt@C+15pt{
Y_{il} & Y_{ij}\otimes Y_{jl} \ar[l]_-{\beta_{ijl}} \\&\\&\\
 &  Y_{ij}\otimes \big( Y_{jk}\otimes Y_{kl} \big) \ar[uuu]_-{Y_{ij}\otimes\beta_{jkl}}  \ar@{}[d]_-{a}|-{\text{\rotatebox[origin=c]{90}{$\cong$}}} \\ 
Y_{ik}\otimes Y_{kl} \ar[uuuu]^-{\beta_{ikl}} & \big(Y_{ij}\otimes Y_{jk}\big)\otimes Y_{kl}\ar[l]^-{\beta_{ijk} \otimes Y_{kl}} }
\end{equation}
toujours que $0\leq i<j<k<l\leq q$.

On considère plusieurs cas: Dans le cas où $0\leq \varphi i<\varphi j<\varphi k<\varphi l\leq q$ le diagramme \eqref{ccohecohe2} est un diagramme de la forme \eqref{ccohecohe}; donc il est commutatif. Si d'un autre $0\leq \varphi i= \varphi j<\varphi k<\varphi l\leq q$ alors le diagramme \eqref{ccohecohe2} devient:
$$
\xymatrix@R=3pt@C+15pt{
X_{\varphi i\varphi l}& \mathbb{1}\otimes X_{\varphi j\varphi l} \ar[l]_-{\ell^{-1}_X} \\&\\&\\
 &  \mathbb{1}\otimes \big( X_{\varphi i\varphi k}\otimes X_{\varphi k\varphi l} \big) \ar[uuu]_-{\mathbb{1}\otimes\alpha_{\varphi i\varphi k\varphi l}} 
 \ar@/_6pt/@{-->}[ld]|-{\ell^{-1}_{X\otimes X}} \ar@{}[d]_-{a}|-{\text{\rotatebox[origin=c]{90}{$\cong$}}} \\ 
X_{\varphi i\varphi k}\otimes X_{\varphi k\varphi l} \ar[uuuu]^-{\alpha_{\varphi i\varphi k\varphi l}}  & \big(\mathbb{1}\otimes X_{\varphi j\varphi k}\big)\otimes X_{\varphi k\varphi l} \ar[l]^-{\ell^{-1}_X \otimes X_{\varphi k\varphi l}}}
$$
qui est commutatif d'après la naturalité de $\ell$ et le Lemme \ref{otrosdiagg}.

Si $0\leq \varphi i < \varphi j = \varphi k<\varphi l\leq q$ il se suit de \eqref{triangle} qu'on a un diagramme commutatif:
$$
\xymatrix@R=3pt@C+18pt{
X_{\varphi i\varphi l}  & X_{\varphi i\varphi j}\otimes X_{\varphi j\varphi l} \ar[l]_-{\alpha_{\varphi i\varphi j\varphi l}}  \\&\\&\\
 &  X_{\varphi i\varphi j}\otimes \big( \mathbb{1} \otimes X_{\varphi j\varphi l} \big) \ar[uuu]_-{X_{\varphi i\varphi j}\otimes\ell^{-1}_X} \ar@{}[d]_-{a}|-{\text{\rotatebox[origin=c]{90}{$\cong$}}} \\ 
X_{\varphi i\varphi j}\otimes X_{\varphi j\varphi l} \ar[uuuu]^-{\alpha_{\varphi i\varphi j\varphi l}} \ar@{=}[ruuuu]& \big(X_{\varphi i\varphi j}\otimes \mathbb{1}\big)\otimes X_{\varphi j\varphi l}\ar[l]^-{r_X^{-1} \otimes X_{\varphi j\varphi l}}\,;}
$$
et si $0\leq \varphi i < \varphi j < \varphi k = \varphi l\leq q$ il se suit de la naturalité de $r$ et du Lemme \ref{otrosdiagg}, qu'on a un diagramme commutatif: 
$$
\xymatrix@R=3pt@C+30pt{
X_{\varphi i\varphi k}  & X_{\varphi i\varphi j}\otimes X_{\varphi j\varphi k}\ar[l]_-{\alpha_{\varphi i\varphi j\varphi k}}\\&\\&\\
 &  X_{\varphi i\varphi j}\otimes \big( X_{\varphi j\varphi k} \otimes \mathbb{1}\big)  \ar[uuu]_-{X_{\varphi i\varphi j}\otimes r_X^{-1}}  \ar@{}[d]_-{a}|-{\text{\rotatebox[origin=c]{90}{$\cong$}}} \\ 
X_{\varphi i\varphi k}\otimes\mathbb{1} \ar[uuuu]^-{r^{-1}_X}& \big(X_{\varphi i\varphi j}\otimes X_{\varphi j\varphi k}\big)\otimes \mathbb{1}
\ar@/^28pt/@<+30pt>@{-->}[uuuu]^-{r^{-1}_{X_{\varphi i\varphi j}\otimes X_{\varphi j\varphi k}}}\ar[l]^-{\alpha_{\varphi i\varphi j\varphi k} \otimes \mathbb{1}} \,.}
$$

D'un autre, si $0\leq \varphi i = \varphi j = \varphi k < \varphi l\leq q$ le diagramme \eqref{ccohecohe2} est de la forme:
$$
\xymatrix@R=3pt@C+15pt{
X_{\varphi i\varphi l}& \mathbb{1}\otimes X_{\varphi i\varphi l} \ar[l]_-{\ell^{-1}_X}\\&\\&\\
 &  \mathbb{1}\otimes \big( \mathbb{1}\otimes X_{\varphi i\varphi l} \big)\ar[uuu]_-{\mathbb{1}\otimes\ell^{-1}_X}  \ar@{}[d]_-{a}|-{\text{\rotatebox[origin=c]{90}{$\cong$}}} \\ 
\mathbb{1}\otimes X_{\varphi i\varphi l}\ar@{=}[uuuur] \ar[uuuu]^-{\ell^{-1}_X}  & \big(\mathbb{1}\otimes \mathbb{1}\big)\otimes X_{\varphi i\varphi l}\ar[l]^-{\ell^{-1}_{\mathbb{1}} \otimes X_{\varphi i\varphi l}}\,,}
$$
si $0\leq \varphi i < \varphi j = \varphi k = \varphi l\leq q$ il est de la forme:
$$
\xymatrix@R=3pt@C+10pt{
X_{\varphi i\varphi j}& X_{\varphi i\varphi j}\otimes \mathbb{1}  \ar[l]_-{r^{-1}_X} \\&\\&\\
 &  X_{\varphi i \varphi j}\otimes \big( \mathbb{1}\otimes \mathbb{1} \big) \ar[uuu]_-{X_{\varphi i \varphi j}\otimes \ell^{-1}_\mathbb{1}}
 \ar@{}[d]_-{a}|-{\text{\rotatebox[origin=c]{90}{$\cong$}}} \\ 
X_{\varphi i \varphi j}\otimes \mathbb{1}  \ar[uuuu]^-{r^{-1}_X}  \ar@{=}[uuuur]& \big(X_{\varphi i \varphi j}\otimes \mathbb{1}\big)\otimes \mathbb{1}\ar[l]^-{r^{-1}_{X} \otimes \mathbb{1}} }
$$
et si $0\leq \varphi i = \varphi j = \varphi k = \varphi l\leq q$ de la forme:
$$
\xymatrix@R=3pt@C+10pt{
\mathbb{1} & \mathbb{1}\otimes \mathbb{1} \ar[l]_-{\ell^{-1}_{\mathbb{1}}} \\&\\&\\
 &  \mathbb{1}\otimes \big( \mathbb{1}\otimes \mathbb{1} \big) \ar[uuu]_-{\mathbb{1}\otimes \ell^{-1}_\mathbb{1}}
 \ar@{}[d]_-{a}|-{\text{\rotatebox[origin=c]{90}{$\cong$}}} \\ 
\mathbb{1}\otimes \mathbb{1}\ar@{=}[uuuur]  \ar[uuuu]^-{\ell^{-1}_{\mathbb{1}}} & \big(\mathbb{1}\otimes \mathbb{1}\big)\otimes \mathbb{1}
\ar[l]^-{\ell^{-1}_{\mathbb{1}} \otimes \mathbb{1}}\,;}
$$
donc ils sont commutatifs d'après le Lemme \ref{otrosdiagg}.

Supposons maintenant que $\xymatrix@C+8pt{\M\ar[r]^-{(F,m^F)}&\M'}$ est un morphisme lax et unitaire entre catégories monoïdaux et $q\geq 0$. On définit la fonction:
$$
\xymatrix@R=3pt{
\M_{q} \ar[rr]^-{(F,m^F)_q} &&\M'_{q}\\
(X,\alpha) \ar@{}[rr]|-{\longmapsto}&& (Y,\beta)}
$$
par les règles: $Y_{ij} = F(X_{ij})$ si $0\leq i<j\leq q$ et $\beta_{ijk}$ est le composé:
$$
\beta_{ijk}\; =\; 
\bigg(\xymatrix@C-18pt{
Y_{ij}\otimes Y_{jk}\ar@{=}[r] & 
F(X_{ij})\otimes F(X_{jk}) \ar[rrrr]^-{m^F_{X_{ij},X_{jk}}} &&&& 
F(X_{ij}\otimes X_{jk}) \ar[rrrr]^-{F(\alpha_{ijk})} &&&& 
F(X_{ik})\ar@{=}[r] & 
Y_{ik}}\bigg)
$$
si $0\leq i<j<k\leq q$.

On montre sans peine que la couple ainsi définie $(Y,\beta)$ est bien un $q$-simplexe de $\mathcal{M}'$.


Le foncteur adjoint du foncteur \eqref{nervio22} qu'on vient de définir est dit le \emph{foncteur nerf} (\emph{géométrique}) des catégories monoïdales:
\begin{equation}\label{nervio22a}
\xymatrix@R=3pt{{\bf cat}_{lax,\star}^\otimes \ar[rr]^-{\mathcal{N}(\,\cdot\,)} && \simp \,;\\ \M \ar@{}[rr]|-{\longmapsto} && \M_\bullet}
\end{equation}
si $\M$ est une catégorie monoïdale on appelle l'ensemble simplicial $\mathcal{N}\big(\M\big)$ \emph{le nerf de $\M$}.

\begin{lemme}\label{3cosqull}
Le nerf $\mathcal{N}\big(\M\big)$ d'une catégorie monoïdale $\M$ est un ensemble simplicial faiblement $2$-cosquelettique (voir \ref{faiblecos}).
\end{lemme}
\begin{proof}
On va déduire cette affirmation de l'énoncé suivant:
\begin{lemme}\label{combinat}
Si $q\geq 2$ et $0\leq n\leq q$ la fonction:
\begin{align*}
\Big\{ \; \big(x_0,\dots,x_n\big) \,\in\, \mathbb{N}^{n+1}\;\, &\Big| \;\, 0\leq x_0<\dots<x_n\leq q+1\;\, \Big\}\\
\xymatrix{\\\ar[u]^-{\varphi}}&\\
\Big\{ \; \big(i_0,\dots,i_n ;s\big) \,\in\, \mathbb{N}^{n+2}\;\, &\Big| \;\, 0\leq i_0<\dots<i_n\leq q  \quad\text{et}\quad 0\leq s\leq q+1 \;\, \Big\} 
\end{align*}
définie par $\varphi(i_0,\dots,i_n;s) = \big(\delta_{s}i_0,\dots,\delta_{s}i_n\big)$ est le conoyau de la double flèche:
\begin{align*}
\Big\{ \; \big(i_0,\dots,i_n ;s\big) \,\in\, \mathbb{N}^{n+2}\;\, \Big|& \;\, 0\leq i_0<\dots<i_n\leq q  \quad\text{et}\quad 0\leq s\leq q+1 \;\, \Big\} \\
&\xymatrix{\\\ar@<+5pt>[u]^-{\varphi_1}\ar@<-5pt>[u]_-{\varphi_2} }\\
\Big\{ \; \big(a_0,\dots,a_n ;s,s'\big) \,\in\, \mathbb{N}^{n+3}\;\, \Big|& \;\, 0\leq a_0<\dots<a_n\leq q-1  \quad\text{et}\quad 0\leq s<s'\leq q+1 \;\, \Big\}
\end{align*}
où $\varphi_1(a_0,\dots,a_n;s,s') = \big(\delta_{s}a_0,\dots,\delta_{s}a_n;s'\big)$ et $\varphi_2(a_0,\dots,a_n;s,s') = \big(\delta_{s'-1}a_0,\dots,\delta_{s'-1}a_n;s\big)$.
\end{lemme}
\begin{proof}
Si $0\leq x_0<\dots<x_n\leq q+1$ son des entiers naturels, vu que par hypothèse $0\leq n \leq q$ il existe un entier $0\leq s\leq q+1$ tel que: Soit $x_k < s < x_{k+1}$ pour certain $0\leq k\leq n$, soit $x_n<s\leq q+1$ ou soit $0\leq s<x_0$. Dans ceux cas on a:
$$
\varphi(x_0,\dots,x_k,x_{k+1}-1,\dots,x_n-1 ;s)=(x_0,\dots,x_n)\qquad \text{où} \qquad 0\leq x_0<\dots<x_k<x_{k+1}-1<\dots<x_n-1\leq q\,,
$$
$$
\varphi(x_0,\dots,x_n ;s)=(x_0,\dots,x_n)\qquad \text{où} \qquad 0\leq x_0<\dots<x_n\leq q \qquad \text{ou}
$$
$$
\varphi(x_0,\dots,x_n ;s)=(x_0-1,\dots,x_n-1)\qquad \text{où} \qquad 0\leq x_0-1<\dots<x_n-1\leq q\,, \qquad \text{respectivement.}
$$

Autrement dit $\varphi$ est une fonction surjective. 

D'un autre côté on sait que $\delta_{s'}\circ\delta_s = \delta_s\circ\delta_{s'-1}$ si $0\leq s<s'\leq q+1$, alors:
\begin{align*}
\varphi\circ\varphi_1 (a_0,\dots,a_n;s,s') \, &= \, \varphi\big( \delta_sa_0,\dots,\delta_sa_n;s'\big) \\
\, &= \, (\delta_{s'}\circ\delta_sa_0,\dots,\delta_{s'}\circ\delta_sa_n)\\
\, &= \, (\delta_s\circ\delta_{s'-1}a_0,\dots,\delta_s\circ\delta_{s'-1}a_n)\\
\, &= \, \varphi\big( \delta_{s'-1}a_0,\dots,\delta_{s'-1}a_n;s\big) \\
\, &= \, \varphi\circ\varphi_2 (a_0,\dots,a_n;s,s')\,,
\end{align*}
toujours que $0\leq a_0<\dots<a_n\leq q-1$. 

On va montrer que si $\varphi(i_0,\dots,i_n ;s)\,=\,\varphi(i'_0,\dots,i'_n ;s')$ où $(i_0,\dots,i_n ;s)\,\neq\,(i'_0,\dots,i'_n ;s')$ et $s\leq s'$, alors $s\neq s'$ et ils existent $0\leq a_0<\dots<a_n\leq q-1$ tels que:
$$
\varphi_1(a_0,\dots,a_n;s,s') \, = \, (i_0,\dots,i_n ;s)
\qquad\text{et}\qquad
\varphi_2(a_0,\dots,a_n;s,s') \, = \, (i'_0,\dots,i'_n ;s')\,.
$$

En effet supposons que $\varphi(i_0,\dots,i_n ;s)\,=\,\varphi(i'_0,\dots,i'_n ;s')$ où $s\leq s'$. Il se suit de la définition de la fonction $\varphi$ que: 
$$
0 \, \leq \, \delta_s(i_0)\,=\,\delta_{s'}(i'_0) \,<\, \dots\, \,<\, \delta_s(i_n)\,=\,\delta_{s'}(i'_n) \, \leq \, q+1\,.
$$

Si $s=s'$ on a que $\delta_s(i_k)=\delta_{s}(i'_k)$ pour tout $0\leq k\leq n$; donc $(i_0,\dots,i_n ;s)\,=\,(i'_0,\dots,i'_n ;s')$ vu que $\delta_s$ est une fonction injective. 

Si $s<s'$ on va considérer deux cas, le cas où aucun des entiers $\delta_s(i_0)\,=\,\delta_{s'}(i'_0),\cdots,\delta_s(i_n)\,=\,\delta_{s'}(i'_n)$ est compris entre $s$ et $s' $ et le cas où il existe $0\leq l\leq n$ tel que $s < \delta_s(i_{l})\,=\,\delta_{s'}(i'_{l})  < s'$.

Dans le premier cas on peut dire que $\delta_s(i_k)\,=\,\delta_{s'}(i'_k)<s<s'<\delta_s(i_{k+1})\,=\,\delta_{s'}(i'_{k+1})$ alors:
$$
\vcenter{\xymatrix@R=2pt{\delta_s(i_0)\,  =\,i_0 \\\vdots\\\delta_s(i_k)\,  =\,i_k \\
\delta_s(i_{k+1})\, =\,i_{k+1}+1 \\\vdots\\\delta_s(i_{n})\,  =\,i_{n}+1}}
\qquad\vcenter{\xymatrix{\text{et}}}\qquad
\vcenter{\xymatrix@R=2pt{\delta_{s'}(i'_0)\,  =\,i'_0 \\\vdots\\\delta_s(i'_k)\,  =\,i'_k \\
\delta_{s'}(i'_{k+1})\, =\,i'_{k+1}+1 \\\vdots\\\delta_{s'}(i'_{n})\,  =\,i'_{n}+1\,;}}
$$
donc $i_{0}=i'_{0}$, $i_{1}=i'_{1}$, $\dots$, $i_{n}=i'_{n}$ et $i_k<s<s'\leq i_{k+1}$. 

On déduit que les entiers $i_0,\dots,i_k,i_{k+1}-1,\dots,i_{n}-1$ vérifient:
$$
0\leq i_0<\dots<i_k<i_{k+1}-1<\dots<i_n -1\leq q-1
$$ 
et en plus:
\begin{align*}
\varphi_1(i_0,\dots,i_k,i_{k+1}-1,\dots,i_n -1;s,s') = &\, \big(\delta_si_0,\dots,\delta_si_k,\delta_s(i_{k+1}-1),\dots,\delta_s(i_n -1);s'\big) \\
 = &\, (i_0,\dots,i_k,i_{k+1},\dots,i_n;s') \\
 = &\, (i'_0,\dots,i'_n;s') 
\end{align*}
$$\text{et}$$
\begin{align*}
\varphi_2(i_0,\dots,i_k,i_{k+1}-1,\dots,i_n -1;s,s') = &\, \big(\delta_{s'-1}i_0,\dots,\delta_{s'-1}i_k,\delta_{s'-1}(i_{k+1}-1),\dots,\delta_{s'-1}(i_n -1);s\big) \\
 = &\, (i_0,,\dots,i_n;s)
\end{align*}
parce que $i_k<s\leq i_{k+1}-1$ et $i_k<s'-1\leq i_{k+1}-1$.


Enfin si on a:
$$
\delta_s(i_k)\,=\,\delta_{s'}(i'_k) < s < \delta_s(i_{k+1})\,=\,\delta_{s'}(i'_{k+1}) < \dots < \delta_s(i_l)\,=\,\delta_{s'}(i'_l) < s' < \delta_s(i_{l+1})\,=\,\delta_{s'}(i'_{l+1})\,,
$$
alors on déduit:
$$
\vcenter{\xymatrix@R=2pt{\delta_s(i_0)\,  =\,i_0 \\\vdots\\\delta_s(i_k)\,  =\,i_k \\
\delta_s(i_{k+1})\, =\,i_{k+1}+1\\\vdots\\\delta_s(i_l)\,  =\,i_l+1 \\
\delta_s(i_{k+1})\, =\,i_{l+1}+1 \\\vdots\\\delta_s(i_{n})\,  =\,i_{n}+1}}
\qquad\vcenter{\xymatrix{\text{et}}}\qquad
\vcenter{\xymatrix@R=2pt{\delta_{s'}(i'_0)\,  =\,i'_0 \\\vdots\\\delta_{s'}(i'_k)\,  =\,i'_k \\
\delta_{s'}(i'_{k+1})\, =\,i'_{k+1}\\\vdots\\\delta_{s'}(i'_l)\,  =\,i'_l \\
\delta_{s'}(i'_{l+1})\, =\,i'_{l+1}+1 \\\vdots\\\delta_{s'}(i'_{n})\,  =\,i'_{n}+1\,;}}
$$
donc:
$$
i_{0}=i'_{0}< \dots < i_{k}=i'_{k} < s <  i_{k+1}+1=i'_{k+1} < \dots < i_{l}+1=i'_{l} < s' < i_{l+1}+1=i'_{l+1}+1 < \dots < i_{n}+1=i'_{n}+1\,.
$$ 

En particulier:
$$
0\leq i'_0 < \dots < i'_k<s\leq i'_{k+1}-1 < \dots < i'_{n}-1 \leq q-1 \qquad \text{et} \qquad
0 \leq i_0 < \dots < i_l<s'-1\leq i_{l+1}-1 < \dots <  i_{n}-1 \leq q-1\,;
$$ 
donc:
\begin{align*}
\,\varphi_1 & \, (i'_0,\dots,i'_k,i'_{k+1}-1,\dots,i'_n -1;s,s') \\
= &\, \big(\delta_s(i'_0),\dots,\delta_s(i'_k),\delta_s(i'_{k+1}-1),\dots,\delta_s(i'_n -1);s'\big) \\
 = &\, \big(i'_0,\dots,i'_k,i'_{k+1},\dots,i'_n ;s'\big) 
\end{align*}
$$\text{et}$$
\begin{align*}
\varphi_2 & \, (i'_0,\dots,i'_k,i'_{k+1}-1,\dots,i'_n -1;s,s') \\
= &\, \varphi_2(i_0,\dots,i_{l},i_{l+1}-1,\dots,i_n -1;s,s') \\
= &\, \big(\delta_{s'-1}(i_0),\dots,\delta_{s'-1}(i_{l}),\delta_{s'-1}(i_{l+1}-1),\dots,\delta_{s'-1}(i_n -1);s\big) \\
= &\, \big(i_0,\dots,i_{l},i_{l+1},\dots,i_n ;s\big)\,.
\end{align*}
\end{proof}

Si $\M$ est une catégorie monoïdale montrons que $\mathcal{N}\big(\M\big)$ est un ensemble simplicial $3$-cosquelettique, c'est-à-dire montrons que si $q\geq 3$ la fonction:
$$
\xymatrix{\M_{q+1}\ar[r]^-{\underset{s}{\prod} \; \delta_{s}}&\underset{0\leq s\leq q+1}{\prod}\M_q}
$$ 
est le noyau dans ${\bf Ens}$ de la double flèche:
$$
\xymatrix@C-1pt{
\underset{0\leq s\leq q+1}{\prod}\M_q
\ar@<+7pt>[rrr]^-{\underset{s<s'}{\prod} \; \delta^*_s\,\circ\, \text{proj}_{s'}}
\ar@<-7pt>[rrr]_-{\underset{s<s'}{\prod} \; \delta^*_{s'-1}\,\circ\, \text{proj}_{s}}&&&
\underset{0\leq s<s'\leq q+1}{\prod}\M_{q-1}}\,.
$$

De façon explicite étant donné un ensemble d'objets de $\M$:
$$
\Big\{  X^s_{ij}  \,\Big|\,  0\leq i<j\leq q\,, \quad 0\leq s\leq q+1\Big\}
$$
et un ensemble de morphismes de $\M$:
$$
\Big\{  \alpha^s_{ijk}\colon\xymatrix@-10pt{ X^s_{ij}\otimes X^s_{jk} \ar[r]& X^s_{ik} }  \; \Big| \;  0\leq i<j<k\leq q\,, \quad 0\leq s\leq q+1\Big\};
$$
tels que:
\begin{enumerate}
\item Si $0\leq i<j<k<l\leq q$ et $0\leq s \leq q+1$ on a un diagramme commutatif:
$$
\xymatrix@R=3pt@C+15pt{
X^s_{il}  & X^s_{ij}\otimes X^s_{jl}  \ar[l]_-{\alpha^s_{ijl}} \\&\\&\\
 &  X^s_{ij}\otimes \big( X^s_{jk}\otimes X^s_{kl} \big) \ar[uuu]_-{X^s_{ij}\otimes\alpha^s_{jkl}}
\ar@{}[d]_-{a}|-{\text{\rotatebox[origin=c]{90}{$\cong$}}} \\ 
X^s_{ik}\otimes X^s_{kl} \ar[uuuu]^-{\alpha^s_{ikl}} & \big(X^s_{ij}\otimes X^s_{jk}\big)\otimes X^s_{kl} \ar[l]^-{\alpha^s_{ijk} \otimes X^s_{kl}}\,,}
$$
\item Si $0\leq s < s' \leq q+1$ on a les égalités:
\begin{align*}
X^{s'}_{\delta_s \, a\delta_sb} \, &= \, X^{s}_{\delta_{s'-1} a \, \delta_{s'-1}b} \qquad &&\text{si}\quad \, 0\leq a< b\leq q-1\\
\alpha^{s'}_{\delta_s a \, \delta_s b \, \delta_s c}  \, &= \, \alpha^{s}_{\delta_{s'-1} a \, \delta_{s'-1}b \, \delta_{s'-1}c}  \qquad &&\text{si}\quad \, 0\leq a< b<c\leq q-1\,;
\end{align*}
\end{enumerate} 
montons qu'il existe un seul ensemble d'objets de $\M$:
$$
\Big\{  Y_{xy}  \,\Big|\,  0\leq x<y\leq q+1\,\Big\}
$$
et un seul ensemble de morphismes de $\M$:
$$
\Big\{  \beta_{xyz}\colon\xymatrix@-10pt{ Y_{xy}\otimes Y_{yz}\ar[r]& Y_{xz} }  \; \Big| \;  0\leq x<y<z\leq q+1\,\Big\};
$$
tels que:
\begin{enumerate}
\item $X_{ij}^s=Y_{\delta_s i \, \delta_s j}$ si $0\leq i<j\leq q$ et $0\leq s \leq q+1$.
\item $\alpha_{ijk}^s=\beta_{\delta_s i \, \delta_s j \, \delta_s k}$ si $0\leq i<j<k\leq q$ et $0\leq s \leq q+1$
\item On a un diagramme commutatif:
$$
\xymatrix@R=3pt@C+15pt{
Y_{il}  & Y_{ij}\otimes Y_{jl} \ar[l]_-{\beta_{ijl}} \\&\\&\\
 &  Y_{ij}\otimes \big( Y_{jk}\otimes Y_{kl} \big) \ar[uuu]_-{Y_{ij}\otimes \beta_{jkl}} 
\ar@{}[d]_-{a}|-{\text{\rotatebox[origin=c]{90}{$\cong$}}} \\ 
 Y_{ik}\otimes Y_{kl} \ar[uuuu]^-{\beta_{ikl}} & \big(Y_{ij}\otimes Y_{jk}\big)\otimes Y_{kl} \ar[l]^-{\beta_{ijk} \otimes Y_{kl}} }
$$
si $0\leq i<j<k<l\leq q+1$.
\end{enumerate}

Vu que par hypothèse $q\geq 3>2$ et on a que:
$$
X^{s'}_{\delta_s a \, \delta_s b} = X^{s}_{\delta_{s'-1} a \, \delta_{s'-1} b} \qquad \text{si} \quad 0\leq a< b\leq q-1 \;\, et \;\, 0\leq s < s' \leq q+1\,, 
$$
il se suit du Lemme \ref{combinat} (pour $n=1$) qu'il existe un seul ensemble d'objets de $\M$:
 $$
\Big\{  Y_{xy}  \,\Big|\,  0\leq x<y\leq q+1\,\Big\}
$$
tels que $X_{ij}^s=Y_{\delta_s i \, \delta_s j}$ si $0\leq i<j\leq q$ et $0\leq s \leq q+1$. 

D'un autre, on déduit du Lemme \ref{combinat} (pour $n=2$) que si $0\leq x<y<z\leq q+1$ alors ils existent $0\leq s \leq q+1$ et $0\leq i<j<k\leq q$ tels que $\delta_s i = x$, $\delta_s j = y$ et $\delta_s k = z$. En particulier on peut définir $\beta_{xyz} = \beta_{\delta_s i \, \delta_s j \, \delta_s k} = \alpha_{ijk}^s$.  

Remarquons que si on a des autres $0\leq s_0 \leq q+1$ et $0\leq i_0<j_0<k_0\leq q$ tels que $\delta_{s_0} i_0 = x$, $\delta_{s_0} j_0 = y$ et $\delta_{s_0} k_0 = z$, on sait de la preuve du Lemme \ref{combinat} qu'il a deux cas: Soit $s=s_0$ et alors $i=i_0$, $j=j_0$ et $k=k_0$; soit disons $s_0<s$ et alors il existent $0\leq a<b<c\leq q-1$ tels que:
$$
\varphi_1(a,b,c;s_0,s) \, = \, (i,j,k;s) \qquad \text{et} \qquad \varphi_2(a,b,c;s_0,s) \, = \, (i_0,j_0,k_0;s_0)
$$
c'est-à-dire tels que:
\begin{align*}
\delta_{s_0} (a) \, =& \, i & \delta_{s-1} (a) \, =& \, i_0 \\
\delta_{s_0} (b) \, =& \, j & \delta_{s-1} (b) \, =& \, j_0 \\
\delta_{s_0} (c) \, =& \, k & \delta_{s-1} (c) \, =& \, k_0 \,.
\end{align*}

Donc on a que $\alpha_{ijk}^s = \alpha_{\delta_{s_0}a \, \delta_{s_0} b \, \delta_{s_0} c }^s = \alpha_{\delta_{s-1}a \, \delta_{s-1} b \, \delta_{s-1} c }^{s_0}=  \alpha_{i_0 \, j_0 \, k_0}^{s_0}$. Autrement dit, $\beta_{xyz} $ est bien défini pour $0\leq x<y<z\leq q+1$ et on a que $\beta_{\delta_s i \, \delta_s j \, \delta_s k}=\alpha_{ijk}^s$ toujours que $0\leq i<j<k\leq q$ et $0\leq s \leq q+1$.

Il nous reste à montrer que si  $0\leq x<y<z<w\leq q+1$ on a carré commutatif:
\begin{equation}\label{lecacarrr}
\xymatrix@R=3pt@C+15pt{
Y_{xw}& Y_{xy}\otimes Y_{yw}  \ar[l]_-{\beta_{xyw}} \\&\\&\\
 &  Y_{xy}\otimes \big( Y_{yz}\otimes Y_{zw} \big) \ar[uuu]_-{Y_{xy}\otimes \beta_{yzw}} 
 \ar@{}[d]_-{a}|-{\text{\rotatebox[origin=c]{90}{$\cong$}}} \\ 
 Y_{xz}\otimes Y_{zw} \ar[uuuu]^-{\beta_{xzw}} & \big(Y_{xy}\otimes Y_{yz}\big)\otimes Y_{zw} \ar[l]^-{\beta_{xyz} \otimes Y_{zw}} }
\end{equation}

En effet si $0\leq x<y<z<w\leq q+1$, vu que $q\geq 3$ ils existent $0\leq s \leq q+1$ et $0\leq i<j<k<l\leq q$ tels que $\delta_s i = x$, $\delta_s j = y$, $\delta_s k = z$ et $\delta_s k = w$; alors le diagramme \eqref{lecacarrr} est le diagramme commutatif:
$$
\xymatrix@R=3pt@C+15pt{
X^s_{il}& X^s_{ij}\otimes X^s_{jl}  \ar[l]_-{\alpha^s_{ijl}}\\&\\&\\
 &  X^s_{ij}\otimes \big( X^s_{jk}\otimes X^s_{kl} \big) \ar[uuu]_-{X^s_{ij}\otimes \alpha^s_{jkl}}
 \ar@{}[d]_-{a}|-{\text{\rotatebox[origin=c]{90}{$\cong$}}} \\ 
 X^s_{ik}\otimes X^s_{kl} \ar[uuuu]^-{\alpha^s_{ikl}} & \big(X^s_{ij}\otimes X^s_{jk}\big)\otimes X^s_{kl} \ar[l]^-{\alpha^s_{ijk} \otimes X^s_{kl}}  }
$$

Donc $\mathcal{N}\big(\mathcal{M}\big)$ est un ensemble simplicial $3$-cosquelettique. 

Finalement remarquons que $\mathcal{N}\big(\mathcal{M}\big)$ est un ensemble simplicial faiblement $2$-cosquelettique, c'est-à-dire notons que la fonction:
\begin{equation} \label{NGfaibcos}
\xymatrix@C+10pt{
\mathrm{Hom}_{\simp}\Big(\Delta^{3},\mathcal{N}\big(\mathcal{M}\big)\Big)
\ar[r]&\mathrm{Hom}_{\simp}\Big(\partial\Delta^{3},\mathcal{N}\big(\mathcal{M}\big)\Big)}
\end{equation}
induite du morphisme d'inclusion $\xymatrix@C-10pt{\partial\Delta^3\ar@{^(->}[r]&\Delta^3}$ est injective. 

En effet on a d'un côté que l'ensemble de départ de \eqref{NGfaibcos} s'identifie à l'ensemble des diagrammes commutatifs de $\M$ de la forme:
\begin{equation}
\xymatrix@R=3pt@C+15pt{
X_{03}& \ar[l]_-{\alpha_{2}} X_{01}\otimes X_{13}  \\&\\&\\
 &  X_{01}\otimes \big( X_{12}\otimes X_{23} \big) \ar[uuu]_-{X_{12}\otimes\alpha_{0}} \ar@{}[d]_-{a}|-{\text{\rotatebox[origin=c]{90}{$\cong$}}} \\ 
X_{02}\otimes X_{23} \ar[uuuu]^-{\alpha_{1}}  & \big(X_{01}\otimes X_{12}\big)\otimes X_{23}  \ar[l]^-{\alpha_{3} \otimes X_{23}} }
\end{equation}
et son ensemble d'arrivée s'identifie à l'ensemble des diagrammes de la même forme mais pas nécessairement commutatif. La fonction \eqref{NGfaibcos} oublie la commutativité; donc elle est injective. 
\end{proof}

Si $\mathcal{M}$ est une catégorie monoïdale, l'ensemble simplicial tronqué $\mathcal{N}\big(\mathcal{M}\big)_q=\M_q$ pour $0\leq q\leq 3$ est décrit de la façon suivante:
\begin{equation}\label{nerf}
\begin{split}
\def\objectstyle{\scriptstyle}
\def\labelstyle{\scriptstyle}
\left\{\left.\vcenter{
\xymatrix@R=3pt@C+15pt{
X_{03}& \ar[l]_-{\alpha_{2}} X_{01}\otimes X_{13}  \\&\\&\\
 &  X_{01}\otimes \big( X_{12}\otimes X_{23} \big) \ar[uuu]_-{X_{12}\otimes\alpha_{0}} \ar@{}[d]_-{a}|-{\text{\rotatebox[origin=c]{90}{$\cong$}}} \\ 
X_{02}\otimes X_{23} \ar[uuuu]^-{\alpha_{1}}  & \big(X_{01}\otimes X_{12}\big)\otimes X_{23}  \ar[l]^-{\alpha_{3} \otimes X_{23}}}}
\right|\vcenter{\xymatrix@R=3pt{\text{Les $\xi_{i}$ sont de}\\\text{morphismes de $\mathcal{M}$}\\
\text{faissant commutatif}\\\text{le diagramme}}}\right\}\\
\xymatrix{\ar@<-32pt>@/_4pt/[d]_-{d_{0}}\ar@<-10pt>@/_1pt/[d]|-{d_{1}}\ar@<10pt>@/^1pt/[d]|-{d_{2}}\ar@<32pt>@/^4pt/[d]^-{d_{3}}\\
\ar@<21pt>@/^2.5pt/[u]|-{s_{0}}\ar[u]|-{s_{1}}\ar@<-21pt>@/_2.5pt/[u]|-{s_{2}}\,}\qquad\qquad\qquad\\
\def\objectstyle{\scriptstyle}
\def\labelstyle{\scriptstyle}
\Big\{  \xymatrix{ X_{2}\otimes X_{0} \ar[r]^-{\xi} & X_{1}}   
\Big|   \vcenter{\xymatrix{\text{$\xi$ est un morphisme de $\mathcal{M}$}}} \Big.\Big\} \qquad\\
\xymatrix{\ar@<-20pt>@/_3.5pt/[d]_-{d_{0}}\ar[d]|-{d_{1}}\ar@<20pt>@/^3.5pt/[d]^-{d_{2}}\\
\ar@<10pt>@/^2pt/[u]|-{s_{0}}\ar@<-10pt>@/_2pt/[u]|-{s_{1}}\,}\qquad\qquad\qquad\quad\,\\
\def\objectstyle{\scriptstyle}
\def\labelstyle{\scriptstyle}
\big\{ \xymatrix{X} \big| \xymatrix{\text{$X$ est un objet de $\mathcal{M}$}}\big\}\qquad\qquad\qquad\\
\xymatrix{\ar@<-12pt>@/_2pt/[d]_-{d_{0}}\ar@<12pt>@/^2pt/[d]^-{d_{1}}\\\ar[u]|-{s_{0}}\,}\qquad\qquad\qquad\qquad\\
\star\qquad\qquad\qquad\qquad\qquad\;\,\,
\end{split}
\end{equation}

Explicitement il existe un seul $0$-simplexes de $\M$, le $0$-simplexe vide. L'ensemble des $1$-simplexes s'identifie à l'ensemble des objets de $\M$, parce qu'il existe une seul couple d'entiers $(i,j)$ tels que $0\leq i<j\leq 1$ et il n'y a pas des triplets $(i,j,k)$ des entiers tels que $0\leq i<j<k\leq 1$. D'un autre un $2$-simplexe est constitué de trois objets  $X_{0}$, $X_{1}$ et $X_{2}$, un pour chacune des couples des entiers $(1,2)$, $(0,2)$ et $(0,1)$ qui vérifient $0\leq i<j\leq 2$ respectivement, et un morphisme: 
\begin{equation}\label{2simplele}
\xi \colon \xymatrix{X_2\otimes X_0 \ar[r] & X_1}
\end{equation} 
pour le seul triplet $(0,1,2)$ des entiers tels que $0\leq i<j<k\leq 2$.

Enfin un $3$-simplexe de $\M$ est déterminé par six objets $\big\{ X_{ij} \,\big|\, 0\leq i<j\leq 3\big\}$ et quatre morphismes de $\M$ :
$$
\Big\{\;\xi_l=\alpha_{ijk}\colon \xymatrix@C-10pt{ X_{ij}\otimes X_{jk} \ar[r] & X_{ik} } \;\,\Big|\quad 0\leq i<j<k\leq 3 \;\;\Big\}\qquad(\text{\tiny{où $0 \leq l \leq 3$ est le seul entier avec $l\neq i,j,k$}})\,;
$$
tels que le diagramme suivant soit commutatif:
\begin{equation}\label{carr}
\eta \quad\;=\;\quad\vcenter{\xymatrix@R=3pt@C+15pt{
X_{03}& \ar[l]_-{\alpha_{2}} X_{01}\otimes X_{13}  \\&\\&\\
 &  X_{01}\otimes \big( X_{12}\otimes X_{23} \big) \ar[uuu]_-{X_{12}\otimes\alpha_{0}} \ar@{}[d]_-{a}|-{\text{\rotatebox[origin=c]{90}{$\cong$}}} \\ 
X_{02}\otimes X_{23} \ar[uuuu]^-{\alpha_{1}}  & \big(X_{01}\otimes X_{12}\big)\otimes X_{23}  \ar[l]^-{\alpha_{3} \otimes X_{23}}\,\text{\normalsize{.}}}}
\end{equation} 

On dit par abus que le morphisme \eqref{2simplele} est un $2$-simplexe de $\M$ et que le diagramme commutatif \eqref{carr} est un $3$-simplexes de $\M$.

Remarquons aussi que les morphismes faces de l'ensemble simplicial tronqué \eqref{nerf} sont donnés par les règles suivantes: Si $\eta$ note le $3$-simplexe \eqref{carr}, alors $d_{i}(\eta)=\xi_{i}$; et si \eqref{2simplele} est un $2$-simplexe, $d_{i}(\xi)=X_{i}$.

D'un autre, les morphismes dégénérescences sont donnés par: $s_0(\star)=\mathbb{1}$,
$$
s_0 (X) \,=\;\Big(\xymatrix@C-5pt{\mathbb{1}\otimes X \ar[r]^-{l_{X}^{-1}} & X}\Big)\,,\quad\quad\quad
s_1 (X) \,=\;\Big(\xymatrix@C-5pt{X\otimes \mathbb{1} \ar[r]^-{r_{X}^{-1}} & X}\Big),
$$
$$
\def\objectstyle{\scriptstyle}
\def\labelstyle{\scriptstyle}
s_{0}(\xi)\,=\,\vcenter{\xymatrix@C+5pt@R-16pt{
X_{1} &\mathbb{1}\otimes X_1 \ar[l]_-{\ell^{-1}} \\&\\
& \mathbb{1}\otimes (X_{2}\otimes X_{0})  \ar[uu]_-{\mathbb{1}\otimes\xi}  \ar@{}[d]_-{a}|-{\text{\rotatebox[origin=c]{90}{$\cong$}}}\\
X_{2}\otimes X_{0}\ar[uuu]^-{\xi}& (\mathbb{1}\otimes X_2)\otimes X_0 \ar[l]^-{\ell^{-1}\otimes X_0}  \,\text{\normalsize{,}}}}\;\;
\qquad 
s_{1}(\xi)\,=\,\vcenter{\xymatrix@C+6pt@R-16pt{
X_{1}&X_{2}\otimes X_{0}  \ar[l]_-{\xi} \\&\\
& X_2(\mathbb{1}\otimes X_0) \ar[uu]_-{X_2\otimes \ell^{-1}} \ar@{}[d]_-{a}|-{\text{\rotatebox[origin=c]{90}{$\cong$}}}\\
X_{2}\otimes X_{0}   \ar[uuu]^-{\xi} &(X_2\otimes \mathbb{1}) \otimes X_0   \ar[l]^-{r^{-1}\otimes X_0}}}\;
$$
$$
\def\objectstyle{\scriptstyle}
\def\labelstyle{\scriptstyle}
\text{et}\qquad s_{2}(\xi)\,=\,\vcenter{\xymatrix@C+6pt@R-16pt{
X_{1}  &X_{2}\otimes X_{0}  \ar[l]_-{\xi} \\&\\
&X_2 \otimes (X_{0}\otimes \mathbb{1}) \ar[uu]_-{X_2\otimes r^{-1}}\ar@{}[d]_-{a}|-{\text{\rotatebox[origin=c]{90}{$\cong$}}}\\
X_{1}\otimes \mathbb{1}  \ar[uuu]^-{\ell^{-1}} &(X_{2} \otimes X_{0})\otimes \mathbb{1}  \ar[l]^-{\xi\otimes \mathbb{1}}\;\text{\normalsize{.}}}}
$$

Montrons premièrement:

\begin{corollaire}\label{plfimonlaxun}
Le foncteur nerf de la catégorie des catégories monoïdales et les foncteurs lax et unitaires $\mathcal{N}\colon\xymatrix{{\bf cat}_{lax,\star}^\otimes \ar[r] & \simp}$ est pleinement fidèle.
\end{corollaire}
\begin{proof}
Pour montrer que $\mathcal{N}$ est un foncteur fidèle considérons $F,G\colon\xymatrix@C-7pt{\mathcal{M}\ar[r]&\mathcal{M}'}$ deux morphismes lax et unitaires entre catégories monoïdales tels que $\mathcal{N}(F)=\mathcal{N}(G)$. Si $X$ est un objet de $\mathcal{M}$ il se suit que $F(X) = \mathcal{N}(F)_1(X) = \mathcal{N}(G)_1(X) = G(X)$. D'un autre si $X$ et $Y$ sont deux objets de $\mathcal{M}$ on que:
$$
m^F_{X,Y}\, = \, \mathcal{N}(F)_2(\mathrm{id}_{X\otimes Y}) \, = \, \mathcal{N}(G)_2(\mathrm{id}_{X\otimes Y}) \, = \, m^G_{X,Y}\,.
$$

Enfin si $f\colon\xymatrix@C-8pt{X\ar[r]&Y}$ est un morphisme de $\mathcal{M}$ remarquons par ailleurs que:
$$
F(f\otimes \mathbb{1}) \circ m^F_{Y,\mathbb{1}}\, = \, \mathcal{N}(F)_2(f\otimes \mathbb{1}) \, = \, \mathcal{N}(G)_2(f\otimes \mathbb{1}) \, = \, G(f\otimes \mathbb{1}) \circ m^G_{Y,\mathbb{1}}\,.
$$

Vu qu'on a de diagrammes commutatifs:
$$
\xymatrix@C+30pt@R-2pt{
F(X\otimes \mathbb{1}) \ar[r]^-{F(f\otimes \mathbb{1})}&F(Y \otimes \mathbb{1})   &
G(X\otimes \mathbb{1}) \ar[r]^-{G(f\otimes \mathbb{1})}&G(Y \otimes \mathbb{1})\\ 
F(X)\otimes F(\mathbb{1})\ar[u]^-{m^F_{X,\mathbb{1}}} \ar[r]|-{F(f)\otimes F(\mathbb{1})}& F(Y)\otimes F(\mathbb{1})\ar[u]_-{m^F_{Y,\mathbb{1}}} &
G(X)\otimes G(\mathbb{1})\ar[u]^-{m^G_{X,\mathbb{1}}} \ar[r]|-{G(f)\otimes G(\mathbb{1})}& G(Y)\otimes G(\mathbb{1})\ar[u]_-{m^G_{Y,\mathbb{1}}}\\
F(X)\ar[u]^-{r_{FX}}\ar[r]_-{F(f)}& F(Y)\ar[u]_{r_{FY}} &
G(X)\ar[u]^-{r_{GX}}\ar[r]_-{G(f)}& G(Y)\ar[u]_{r_{GY}}}
$$ 
on déduit que $F(f)=G(f)$, parce que $r_{FY} = r_{GY}$ et  $m^F_{Y,\mathbb{1}}= m^G_{Y,\mathbb{1}}$ sont des isomorphismes de $\mathcal{M}'$: 
$$
\vcenter{\xymatrix@C+6pt@R=4pt{
F(Y\otimes \mathbb{1})  & \\
&   F(Y)  \ar@<-3pt>[lu]_-{\ell_{FY}}^-{\cong}  \ar@<+3pt>[ld]^-{F(\ell_Y)}_-{\cong}  \,. \\
F(Y)\otimes\mathbb{1}\ar[uu]^-{m^F_{Y,\mathbb{1}}}    & 
}}
$$
Donc $\mathcal{N}$ est un foncteur fidèle.

D'un autre si $\mathcal{M}$ et $\mathcal{M}'$ sont de catégories monoïdales, soit $\varphi\colon\xymatrix@C-5pt{\mathcal{N}(\mathcal{M}) \ar[r] &\mathcal{N}(\mathcal{M}')}$ un morphisme d'ensembles simpliciaux. On définit un morphisme lax et unitaire de catégories monoïdales $F\colon\xymatrix@C-7pt{\mathcal{M}\ar[r]&\mathcal{M}'}$ tel que $\mathcal{N}(F)=\varphi$ de la façon suivante: Si $X$ est un objet de $\mathcal{M}$ on pose $F(X)=\varphi_1(X)$. Il se suit que: 
$$
F(\mathbb{1})=\varphi_1\circ s_0 (\star) = s_0\circ \varphi_0 (\star) = s_0 (\star) = \mathbb{1}\,. 
$$

Si $f\colon\xymatrix@C-8pt{X\ar[r]&Y}$ est un morphisme de $\mathcal{M}$ on définit $F(f)\colon\xymatrix@C-8pt{F(X)\ar[r]&F(Y)}$ comme le morphisme composé de $\mathcal{M}'$:
$$
\xymatrix@C+13pt{ 
F(X) \ar[r]^-{r_{FX}} &F(X)\otimes \mathbb{1} \ar[r]^-{\varphi_2(f \circ r_X^{-1})} & F(Y)\,.}
$$

On a en particulier que:
\begin{align*}
F(\mathrm{id}_X) \, &= \, \varphi_2(\mathrm{id}_X\circ r_X^{-1}) \circ r_{FX} \, = \, \varphi_2(r_X^{-1}) \circ r_{FX} \, = \,  \big(\varphi_2\circ s_1(X)\big) \circ r_{FX}\\ 
\, &= \, \big(s_1 \circ \varphi_1(X)\big) \circ r_{FX}  \, = \, r_{FX}^{-1} \circ r_{FX}   \, = \, \mathrm{id}_{FX}\,.
\end{align*}

D'un autre si $f\colon\xymatrix@C-8pt{X\ar[r]&Y}$ et $g\colon\xymatrix@C-8pt{Y\ar[r]&Z}$ sont deux morphismes de $\mathcal{M}$, remarquons par ailleurs que le diagramme: 
\begin{equation}\label{3simpFER}
\xymatrix@R-16pt@C+18pt{
Z &X\otimes \mathbb{1}  \ar[l]_-{(g\circ f)\circ r_X^{-1}} \\&\\
&X \otimes (\mathbb{1}\otimes \mathbb{1}) \ar[uu]_-{X\otimes r^{-1}_{\mathbb{1}} }  \ar@{}[d]_-{a}|-{\text{\rotatebox[origin=c]{90}{$\cong$}}}\\
Y\otimes \mathbb{1} \ar[uuu]^-{g\circ r_Y^{-1}} &(X\otimes \mathbb{1})\otimes \mathbb{1}   \ar[l]^-{(f\circ r_X^{-1})\otimes \mathbb{1}} \,\text{\normalsize{,}}   }
\end{equation}
est commutatif. En effet on a que dans le diagramme:
$$
\vcenter{\xymatrix@R-13pt@C+30pt{
Y&X\otimes \mathbb{1}  \ar[l]_-{(g\circ f)\circ r_X^{-1}} & \\
&&\\
&X \otimes (\mathbb{1}\otimes \mathbb{1}) \ar@{}[rddd]_{\text{(II)}} \ar[uu]|-{X\otimes r^{-1}_{\mathbb{1}}}  \ar@{}[d]_-{a}|-{\text{\rotatebox[origin=c]{90}{$\cong$}}} &\\
Y\otimes \mathbb{1} \ar[uuu]^-{g\circ r_Y^{-1}}  &\big(X\otimes \mathbb{1}\big)\otimes \mathbb{1}  \ar[l]|-{(f\circ r_X^{-1})\otimes \mathbb{1}} &\\
&&\\
Y \ar@{}[uur]|-{\text{(I)}} \ar[uu]^-{r_{Y}}&X\otimes \mathbb{1} \ar[l]^-{f\circ r_X^{-1}} \ar[uu]|-{r_{X\otimes\mathbb{1}}}& X\ar[l]^-{r_{X}}  \ar@/_14pt/@<-6pt>[luuuuu]_{r_{X}}}}
$$
le carré (I) est commutatif parce que $r_{?}$ est une transformation naturelle et le diagramme (II) est commutatif d'après le Corollaire \ref{otrosdiagg}. Donc vu que:
$$
\big((g\circ f)\circ r_X^{-1}\big) \circ r_X \, = \, g\circ f \, = \,  \big(g\circ r_Y^{-1}\big) \circ r_Y \circ \big(f\circ r_X^{-1}\big) \circ r_X  \,,
$$
on en conclut que:
$$
\big((g\circ f)\circ r_X^{-1}\big) \circ \big(X\otimes r_{\mathbb{1}}^{-1}\big) \circ a_{X,\mathbb{1},\mathbb{1}} \circ r_{X\otimes \mathbb{1}}\circ r_X 
\, = \,  \big(g\circ r_Y^{-1}\big) \circ \big((f\circ r_X^{-1} )\otimes \mathbb{1}\big) \circ r_{X\otimes \mathbb{1}} \circ r_X  \,.
$$

Il se suit que \eqref{3simpFER} est un diagramme commutatif parce que $r_{X\otimes \mathbb{1}} \circ r_X$ est un isomorphisme de $\mathcal{M}$.

Considérons maintenant l'image du $3$-simplexe \eqref{3simpFER} par le morphisme $\varphi$: 
$$
\xymatrix@R-16pt@C+23pt{
F(Z) &F(X)\otimes \mathbb{1}  \ar[l]_-{\varphi_2\big((g\circ f)\circ r_X^{-1}\big)} \\&\\
&F(X) \otimes (\mathbb{1}\otimes \mathbb{1}) \ar[uu]_-{F(X)\otimes r^{-1}_{\mathbb{1}}}  \ar@{}[d]_-{a}|-{\text{\rotatebox[origin=c]{90}{$\cong$}}}\\
F(Y)\otimes \mathbb{1} \ar[uuu]^-{\varphi_2\big(g\circ r_Y^{-1}\big)} &\big(F(X)\otimes \mathbb{1}\big)\otimes \mathbb{1}   \ar[l]^-{\varphi_2(f\circ r_X^{-1})\otimes \mathbb{1}} \,\text{\normalsize{.}}   }
$$

Vu qu'on a des diagramme commutatifs:
$$
\vcenter{\xymatrix@C+18pt{
F(Y)\otimes \mathbb{1} &\big(F(X)\otimes \mathbb{1}\big)\otimes \mathbb{1}   \ar[l]_-{\varphi_2(f\circ r_X^{-1})\otimes \mathbb{1}} \\
F(Y)  \ar[u]^-{r_{FY}}  & F(X)\otimes \mathbb{1} \ar[l]^-{\varphi_2\big(f\circ r_Y^{-1}\big)}\ar[u]_{r_{FX\otimes\mathbb{1}}}
}}\qquad \;\text{et}\;\qquad
\vcenter{\xymatrix@R-16pt@C+20pt{
F(X)\otimes \mathbb{1}  & \\
&\\
F(X) \otimes (\mathbb{1}\otimes \mathbb{1}) \ar[uu]^-{F(X)\otimes r^{-1}_{\mathbb{1}}}  \ar@{}[d]_-{a}|-{\text{\rotatebox[origin=c]{90}{$\cong$}}} &\\
\big(F(X)\otimes \mathbb{1}\big)\otimes \mathbb{1}&\\
&\\
F(X)\otimes \mathbb{1} \ar[uu]^-{r_{FX\otimes\mathbb{1}}}& F(X)\ar[l]^-{r_{FX}}  \ar@/_14pt/@<-6pt>[luuuuu]_{r_{FX}}}}
$$
il se suit que:
$$
F(g\circ f) \, = \, \varphi_2\big((g\circ f)\circ r_X^{-1}\big) \circ \big(r_{FX}\big) \, = \,   \big(\varphi_2( g\circ r_X^{-1}) \circ (r_{FY})\big) \circ \big(\varphi_2( f\circ r_X^{-1}) \circ (r_{FX})\big) \, = \, F(g)\circ F(f)\,.
$$

Donc $F\colon\xymatrix@C-7pt{\mathcal{M}\ar[r]&\mathcal{M}'}$ est bien un foncteur tel que $F(\mathbb{1})=\mathbb{1}$.

D'un autre côté si $X$ et $Y$ sont des objets de $\mathcal{M}$ on définit $m^F_{X,Y}\colon\xymatrix@C-5pt{F(X)\otimes F(Y)\ar[r]& F(X\otimes Y)}$ comme le morphisme de $\mathcal{M}'$:
$$
\xymatrix@C+25pt{
F(X)\otimes F(Y) \ar[r]^-{\varphi_2(\mathrm{id}_{X\otimes Y})} & F(X\otimes Y)\,.}
$$

Montrons que $m^F$ est bien une transformation naturelle: De façon explicite si $f\colon\xymatrix@C-8pt{X\ar[r]&X'}$ et $g\colon\xymatrix@C-8pt{Y\ar[r]&Y'}$ sont deux morphismes de $\mathcal{M}$ vérifions que:
$$
F(f\otimes g) \circ m^F_{X,Y} \, = \, m^F_{X',Y'}\circ \big(F(f)\otimes F(g)\big)\,,
$$
c'est-à-dire montrons que:
\begin{equation}\label{unofielferfer}
\begin{split}
\text{\scriptsize{$
\Big( \varphi_2\big((f\otimes g)\circ r_{X\otimes Y}^{-1}\big)\Big) \circ \big(r_{F(X\otimes Y)}\big) \circ \big(\varphi_2(\mathrm{id}_{X\otimes Y})\big) \, = \, $}}\qquad \quad \quad\\
\text{\scriptsize{$\big(\varphi_2(\mathrm{id}_{X'\otimes Y'})\big) \circ \Big( \big(\varphi_2(f\circ r_{X}^{-1})\circ r_{FX}\big) \otimes F(Y')\Big) \circ 
\Big(F(X) \otimes \big(\varphi_2(g\circ r_{Y}^{-1})\circ r_{FY}\big)\Big)$}}\,.
\end{split}
\end{equation}

Pour commencer remarquons que le carré commutatif:
$$
\xymatrix@C+35pt@R+5pt{
F(X\otimes Y) \otimes \mathbb{1} & \big( F(X) \otimes F(Y) \big)\otimes \mathbb{1}\ar[l]_-{\varphi_2(\mathrm{id}_{X\otimes Y}) \otimes \mathbb{1}}\\
F(X\otimes Y) \ar[u]^-{r_{F(X\otimes Y)}} & F(X)\otimes F(Y) \ar[l]^-{\varphi_2(\mathrm{id}_{X\otimes Y}) }
\ar[u]_-{r_{F(X)\otimes F(Y)}}
}
$$
implique que:
\begin{equation}\label{dosfielferfer}
\text{\scriptsize{$
\Big( \varphi_2\big((f\otimes g)\circ r_{X\otimes Y}^{-1}\big)\Big) \circ \big(r_{F(X\otimes Y)}\big) \circ \big(\varphi_2(\mathrm{id}_{X\otimes Y})\big) \, = \, 
\Big( \varphi_2\big((f\otimes g)\circ r_{X\otimes Y}^{-1}\big)\Big) \circ \big(\varphi_2(\mathrm{id}_{X\otimes Y})\otimes \mathbb{1}\big) \circ \big(r_{FX\otimes FY}\big)$}}\,.
\end{equation}

D'un autre, vu qu'on a les diagrammes commutatifs de $\mathcal{M}$:
$$
\vcenter{\xymatrix@R=3pt@C+22pt{
X'\otimes Y' & X\otimes Y'  \ar[l]_-{f\otimes Y'} \\&\\&\\
 &  X\otimes \big( Y\otimes \mathbb{1} \big) \ar[uuu]_-{X\otimes (g\circ r_Y^{-1})} 
  \ar@{}[d]_-{a}|-{\text{\rotatebox[origin=c]{90}{$\cong$}}} \\ 
(X\otimes Y)\otimes \mathbb{1} \ar[uuuu]^-{(f\otimes g)\circ r_{X\otimes Y}^{-1}}  & \big(X\otimes Y\big)\otimes \mathbb{1} \ar[l]^-{\mathrm{id}_{X\otimes Y} \otimes \mathbb{1}}
}}\qquad \; \text{et}\qquad \;
\vcenter{\xymatrix@R=3pt@C+25pt{
X'\otimes Y' & X\otimes Y'  \ar[l]_-{f\otimes Y'} \\&\\&\\
 &  X\otimes \big( \mathbb{1}\otimes Y' \big) \ar[uuu]_-{X\otimes (\ell_{Y'}^{-1})} 
  \ar@{}[d]_-{a}|-{\text{\rotatebox[origin=c]{90}{$\cong$}}} \\ 
X'\otimes Y' \ar[uuuu]^-{\mathrm{id}_{X'\otimes Y'}}  & \big(X\otimes \mathbb{1}\big)\otimes Y' \ar[l]^-{(f\circ r_{X}^{-1})\otimes Y'}
\ar@/^28pt/@<+16pt>@{-->}[uuuu]^-{r^{-1}_{X}\otimes Y'} }}
$$
il se suit qu'on a les diagrammes commutatifs de $\mathcal{M}'$:
$$
\xymatrix@R=3pt@C+28pt{
F(X'\otimes Y') & F(X)\otimes F(Y')  \ar[l]_-{\varphi_2(f\otimes Y')} \\&\\&\\
 &  F(X)\otimes \big( F(Y)\otimes \mathbb{1} \big) \ar[uuu]_-{FX\otimes \varphi_2(g\circ r_Y^{-1})} 
  \ar@{}[d]_-{a}|-{\text{\rotatebox[origin=c]{90}{$\cong$}}} \\ 
F(X\otimes Y)\otimes \mathbb{1} \ar[uuuu]^-{\varphi_2((f\otimes g)\circ r_{X\otimes Y}^{-1})}  & \big(F(X)\otimes F(Y)\big)\otimes \mathbb{1} \ar[l]^-{\varphi_2(\mathrm{id}_{X\otimes Y}) \otimes \mathbb{1}}
}
$$
$$
\text{et}\qquad \; 
\vcenter{\xymatrix@R=3pt@C+32pt{
F(X'\otimes Y') & F(X)\otimes F(Y')  \ar[l]_-{\varphi_2(f\otimes Y')} \\&\\&\\
 &  F(X)\otimes \big( \mathbb{1}\otimes F(Y') \big) \ar[uuu]_-{FX\otimes (\ell_{FY'}^{-1})} 
  \ar@{}[d]_-{a}|-{\text{\rotatebox[origin=c]{90}{$\cong$}}} \\ 
F(X')\otimes F(Y') \ar[uuuu]^-{\varphi_2(\mathrm{id}_{X'\otimes Y'})}  & \big(F(X)\otimes \mathbb{1}\big)\otimes F(Y') \ar[l]^-{\varphi_2(f\circ r_{X}^{-1})\otimes F(Y')}
\ar@/^28pt/@<+25pt>@{-->}[uuuu]^-{r^{-1}_{FX}\otimes FY'} \,;}}
$$
en particulier: 
\begin{equation}\label{trefielferfer}
\begin{split}
\text{\scriptsize{$
\Big( \varphi_2\big((f\otimes g)\circ r_{X\otimes Y}^{-1}\big)\Big) \circ \big(\varphi_2(\mathrm{id}_{X\otimes Y})\otimes \mathbb{1}\big) \circ \big(r_{FX\otimes FY}\big) \, = \, $}}
\qquad\qquad\qquad\qquad\;\\
\text{\scriptsize{$ \big( \varphi_2(f\otimes Y') \big) \circ \big( F(X)\otimes \varphi_2 (g\circ r_Y^{-1})\big)\circ \big(a_{FX,FY,\mathbb{1}}\big)\circ \big(r_{FX\otimes FY}\big)  \, = \, $}}
\qquad\qquad\qquad\;\\
\text{\scriptsize{$\big( \varphi_2 (\mathrm{id}_{X'\otimes Y'})\big) \circ \big(\varphi_2 (f\circ r_X^{-1})\otimes F(Y')\big) \circ \big(r_{FX}\otimes F(Y')\big)\circ \big( F(X)\otimes \varphi_2 (g\circ r_Y^{-1})\big) \circ \big(a_{FX,FY,\mathbb{1}}\big)\circ \big(r_{FX\otimes FY}\big)$}}\,.
\end{split}
\end{equation}

Il se suit de \eqref{dosfielferfer}, \eqref{trefielferfer} et du triangle commutatif:
$$
\xymatrix@C-5pt@R-5pt{
\big(F(X)\otimes F(Y)\big)\otimes \mathbb{1} \ar[rr]^-{a_{FX,FY,\mathbb{1}}} && F(X)\otimes \big(F(Y)\otimes \mathbb{1}\big)\\
&F(X)\otimes F(Y)\ar[ru]_-{F(X)\otimes r_{FY}} \ar[lu]^-{r_{FX\otimes FY}}& }
$$
qu'on a bien l'égalité \eqref{unofielferfer}.

Donc $m^F$ est une transformation naturelle. Remarquons d'un autre que:
$$
\def\objectstyle{\scriptstyle}
\def\labelstyle{\scriptstyle}
\vcenter{\xymatrix@C-3pt@R-5pt{
&&FX\otimes \mathbb{1}  \ar[dd]^-{m^F_{X,\mathbb{1}}} \\
&FX\ar[ru]^{r_{FX}}\ar[rd]_{Fr_X}&\\
&&F(X\otimes\mathbb{1}) 
 }}
$$
est un triangle commutatif parce que:
$$
F(r_X) \, = \, \varphi_2(\mathrm{id}_{X\otimes \mathbb{1}}) \circ r_{FX} \, = \, m^F_{X,\mathbb{1}} \circ r_{FX}\,.
$$

Pour montrer qu'on a aussi un triangle commutatif:
$$
\def\objectstyle{\scriptstyle}
\def\labelstyle{\scriptstyle}
\vcenter{\xymatrix@C-3pt@R-5pt{
\mathbb{1}\otimes FX \ar[dd]_-{m^F_{\mathbb{1},X}} && \\
&FX\ar[lu]_{l_{FX}}\ar[ld]^-{Fl_X}&\\
 F(\mathbb{1}\otimes X)   &&
 }}
$$
montrons que si $f\colon\xymatrix@C-8pt{X\ar[r]&Y}$ est un morphisme de $\mathcal{M}$ alors le morphisme $F(f)=\varphi_2(f \circ r_X^{-1})\circ r_{FX}$ est égale au composé:
$$
\xymatrix@C+13pt{ 
F(X) \ar[r]^-{\ell_{FX}} &\mathbb{1}\otimes F(X) \ar[r]^-{\varphi_2(f \circ \ell_X^{-1})} & F(Y)}
$$

Pour cela remarquons que d'après les diagrammes commutatifs:
$$
\vcenter{\xymatrix@C+8pt@R+3pt{
X\otimes\mathbb{1} \ar[r]^-{\ell_X\otimes\mathbb{1}}  & (\mathbb{1}\otimes X) \otimes \mathbb{1} \\
X\ar[r]_-{\ell_X} \ar[u]^-{r_X}&\mathbb{1}\otimes X\ar[u]_-{r_{\mathbb{1}\otimes X}}
}}\qquad \text{et} \qquad
\vcenter{\xymatrix@R-5pt@C-8pt{
(\mathbb{1}\otimes X) \otimes \mathbb{1} \ar[rr]^-{a_{\mathbb{1},X,\mathbb{1}}} && \mathbb{1}\otimes (X\otimes \mathbb{1})\\
 &\mathbb{1}\otimes X\ar[ul]^-{r_{\mathbb{1}\otimes X}}\ar[ru]_-{\mathbb{1}\otimes r_X}&   
}}
$$
on a les égalités des morphismes de $\mathcal{M}'$ suivantes:
\begin{align*}
\big(r_X\big)\circ\big(\ell_X^{-1}\big)\circ\big(\mathbb{1}\otimes r_X^{-1}\big)\circ\big(a_{\mathbb{1},X,\mathbb{1}}\big) \, &= \,
\big(r_X\big)\circ\big(\ell_X^{-1}\big)\circ\big(\mathbb{1}\otimes r_X\big)^{-1}\circ\big(a_{\mathbb{1},X,\mathbb{1}}\big) \\
\, &= \, \big(r_X\big)\circ\big(\ell_X^{-1}\big)\circ\big(r_{\mathbb{1}\otimes X}\big)^{-1} \, = \,  \big(\ell_X\otimes \mathbb{1}\big)^{-1} \, = \, \ell_{X}^{-1}\otimes \mathbb{1}\,;
\end{align*}
autrement dit on a un diagramme commutatif:
\begin{equation}\label{3simplexecoqu}
\xymatrix@R=3pt@C+15pt{
X& \mathbb{1}\otimes X  \ar[l]_-{\ell^{-1}_X} \\&\\&\\
 &  \mathbb{1}\otimes \big( X\otimes \mathbb{1} \big) \ar[uuu]_-{\mathbb{1}\otimes r_X^{-1}} 
  \ar@{}[d]_-{a}|-{\text{\rotatebox[origin=c]{90}{$\cong$}}} \\ 
X\otimes \mathbb{1} \ar[uuuu]^-{r_X^{-1}}  & \big(\mathbb{1}\otimes X\big)\otimes \mathbb{1} \ar[l]^-{\ell^{-1}_X \otimes \mathbb{1}}
}
\end{equation}

Considérons le $3$-simplexe de $\mathcal{M}'$ obtenu de façon analogue:
$$
\xymatrix@R=3pt@C+15pt{
FX& \mathbb{1}\otimes FX  \ar[l]_-{\ell^{-1}_{FX}} \\&\\&\\
 &  \mathbb{1}\otimes \big( FX\otimes \mathbb{1} \big) \ar[uuu]_-{\mathbb{1}\otimes r_{FX}^{-1}} 
  \ar@{}[d]_-{a}|-{\text{\rotatebox[origin=c]{90}{$\cong$}}} \\ 
FX\otimes \mathbb{1} \ar[uuuu]^-{r_{FX}^{-1}}  & \big(\mathbb{1}\otimes FX\big)\otimes \mathbb{1} \ar[l]^-{\ell^{-1}_{FX} \otimes \mathbb{1}}
}
$$
ainsi que le $3$-simplexe de $\mathcal{M}'$ image de \eqref{3simplexecoqu} par la fonction $\varphi_3$:
$$
\xymatrix@R=3pt@C+25pt{
FY& \mathbb{1}\otimes FX  \ar[l]_-{\varphi_2(f\circ \ell^{-1}_{X})} \\&\\&\\
 &  \mathbb{1}\otimes \big( FX\otimes \mathbb{1} \big) \ar[uuu]_-{\mathbb{1}\otimes  r_{FX}^{-1}} 
  \ar@{}[d]_-{a}|-{\text{\rotatebox[origin=c]{90}{$\cong$}}} \\ 
FX\otimes \mathbb{1} \ar[uuuu]^-{\varphi_2(f\circ r_{X}^{-1})}  & \big(\mathbb{1}\otimes FX\big)\otimes \mathbb{1} \ar[l]^-{\ell^{-1}_{FX}\otimes \mathbb{1}}
}
$$ 

Il se suit que:
\begin{align*}
F(f) \, &= \, \varphi_2(f \circ r_X^{-1})\circ r_{FX} \\
&= \, \varphi_2(f \circ \ell_X^{-1})\circ\big(\mathbb{1}\otimes r^{-1}_{FX}\big)\circ\big(a_{\mathbb{1},FX,\mathbb{1}}\big)\circ \big( \ell_{FX}^{-1}\otimes \mathbb{1}\big)^{-1}\circ r_{FX}\\
&= \, \varphi_2(f \circ \ell_X^{-1})\circ \ell_{FX}\,.
\end{align*}

Finalement pour montrer qu'on a un diagramme commutatif:
\begin{equation}\label{ferrrrrr}
\def\objectstyle{\scriptstyle}
\def\labelstyle{\scriptstyle}
\vcenter{\xymatrix@C+23pt@R-3pt{
F\big((X\otimes Y)\otimes Z\big)   \ar[d]_-{Fa_{X,Y,Z}}&F(X\otimes Y)\otimes FZ  \ar[l]_-{m^F_{X\otimes Y, Z}}&(FX\otimes FY)\otimes FZ\ar[d]^-{a_{FX,FY,FZ}} \ar[l]_-{m^F_{X,Y}\otimes FZ}\\
F\big(X\otimes (Y\otimes Z)\big)  &FX\otimes F(Y\otimes Z) \ar[l]^-{m^F_{X,Y\otimes Z}} &FX\otimes (FY\otimes FZ)  \ar[l]^-{FX\otimes m^F_{Y,Z}}  }}
\end{equation}
considérons les $3$-simplexes de $\mathcal{M}$:
$$
\vcenter{\xymatrix@R=3pt@C+30pt{
X\otimes (Y\otimes Z)& X \otimes (Y\otimes Z) \ar[l]_-{\mathrm{id}_{X\otimes (Y\otimes Z)}} \\&\\&\\
 &  X\otimes (Y\otimes Z) \ar[uuu]_-{X\otimes \mathrm{id}_{Y\otimes Z}} 
  \ar@{}[d]_-{a}|-{\text{\rotatebox[origin=c]{90}{$\cong$}}} \\ 
(X\otimes Y)\otimes Z  \ar[uuuu]^-{a_{X,Y,Z}}  & 
(X\otimes Y)\otimes Z  \ar[l]^-{\mathrm{id}_{X\otimes Y}\otimes Z}
}}
$$
$$
\text{et}\qquad\;\vcenter{\xymatrix@R=3pt@C+30pt{
X\otimes (Y\otimes Z)& (X\otimes Y)\otimes Z \ar[l]_-{a_{X,Y,Z}} \\&\\&\\
 &  (X\otimes Y)\otimes \big( Z\otimes \mathbb{1} \big) \ar[uuu]_-{(X\otimes Y)\otimes r_Z^{-1}} 
  \ar@{}[d]_-{a}|-{\text{\rotatebox[origin=c]{90}{$\cong$}}} \\ 
\big((X\otimes Y)\otimes Z\big)\otimes \mathbb{1} \ar[uuuu]^-{a_{X,Y,Z}\circ r_{(X\otimes Y)\otimes Z}^{-1}}  & 
\big((X\otimes Y)\otimes Z\big)\otimes \mathbb{1} \ar[l]^-{\mathrm{id}_{(X\otimes Y)\otimes Z}\otimes \mathbb{1}}\,,
\ar@/^28pt/@<+30pt>@{-->}[uuuu]^-{r^{-1}_{(X\otimes Y)\otimes Z}}}}
$$
ainsi que ses images par le morphismes $\varphi$:
\begin{equation}\label{1erimagine}
\xymatrix@R=3pt@C+30pt{
F\big(X\otimes (Y\otimes Z)\big)& F(X) \otimes F(Y\otimes Z) \ar[l]_-{\varphi_2\big(\mathrm{id}_{X\otimes (Y\otimes Z)}\big)} \\&\\&\\
 &  F(X)\otimes \big(F(Y)\otimes F(Z)\big) \ar[uuu]_-{FX\otimes \varphi_2(\mathrm{id}_{Y\otimes Z})} 
  \ar@{}[d]_-{a}|-{\text{\rotatebox[origin=c]{90}{$\cong$}}} \\ 
F(X\otimes Y)\otimes F(Z)  \ar[uuuu]^-{\varphi_2(a_{X,Y,Z})}  & 
\big(F(X)\otimes F(Y)\big)\otimes F(Z)  \ar[l]^-{\varphi_2(\mathrm{id}_{X\otimes Y})\otimes FZ}}
\end{equation}
\begin{equation} \label{2erimagine}
\text{et}\qquad\vcenter{\xymatrix@R=3pt@C+30pt{
F\big(X\otimes (Y\otimes Z)\big)& F(X\otimes Y)\otimes FZ \ar[l]_-{\varphi_2(a_{X,Y,Z})} \\&\\&\\
 &  F(X\otimes Y)\otimes \big( FZ\otimes \mathbb{1} \big) \ar[uuu]_-{F(X\otimes Y)\otimes r_{FZ}^{-1}} 
  \ar@{}[d]_-{a}|-{\text{\rotatebox[origin=c]{90}{$\cong$}}} \\ 
F\big((X\otimes Y)\otimes Z\big)\otimes \mathbb{1}\ar[uuuu]^-{\varphi_2\big(a_{X,Y,Z}\circ r_{(X\otimes Y)\otimes Z}^{-1}\big)}  & 
\big(F(X\otimes Y)\otimes FZ\big)\otimes \mathbb{1}  \ar[l]^-{\varphi_2(\mathrm{id}_{(X\otimes Y)\otimes Z})\otimes \mathbb{1}}    
\ar@/^28pt/@<+30pt>@{-->}[uuuu]^-{r^{-1}_{F(X\otimes Y)\otimes FZ}}} }
\end{equation}

Il se suit des diagrammes \eqref{1erimagine} et \eqref{2erimagine} que:
\begin{align*}
\big(\varphi_2(\mathrm{id}_{X\otimes(Y\otimes Z)})\big)\circ\big( FX\otimes \varphi_2(\mathrm{id}_{Y\otimes Z})\big)\circ\big(a_{FX,FY,FZ}\big) \, = \, 
\big(\varphi_2(a_{X,Y,Z})\big)\circ\big(\varphi_2(\mathrm{id}_{X\otimes Y})\otimes FZ\big) \\
 \, = \, 
 \big(\varphi_2(a_{X,Y,Z}\circ r^{-1}_{(X\otimes Y)\otimes Z})\big)\circ    \big(\varphi_2(\mathrm{id}_{(X\otimes Y)\otimes Z} \otimes \mathbb{1})\big)\circ \big(r_{F(X\otimes Y)\otimes FZ}\big)\circ   \big(\varphi_2(\mathrm{id}_{X\otimes Y})\otimes FZ\big) 
\end{align*}

Plus encore vu que $r_{?}$ est une transformation naturelle on a un carré commutatif:
$$
\xymatrix@C+35pt@R+5pt{
F\big((X\otimes Y) \otimes Z\big) \otimes \mathbb{1} & \big( F(X\otimes Y)\otimes FZ\big)\otimes \mathbb{1}\ar[l]_-{\varphi_2(\mathrm{id}_{(X\otimes Y)\otimes Z}) \otimes \mathbb{1}}\\
F\big((X\otimes Y) \otimes Z\big)\ar[u]^-{r_{F((X\otimes Y)\otimes Z)}} & \big( F(X\otimes Y)\otimes FZ\big) \ar[l]^-{\varphi_2(\mathrm{id}_{(X\otimes Y)\otimes Z}) }
\ar[u]_-{r_{F(X\otimes Y)\otimes FZ}}
}
$$
c'est-à-dire on a que:
\begin{align*}
\big(\varphi_2(a_{X,Y,Z}\circ r^{-1}_{(X\otimes Y)\otimes Z})\big)\circ    \big(\varphi_2(\mathrm{id}_{(X\otimes Y)\otimes Z} \otimes \mathbb{1})\big)\circ \big(r_{F(X\otimes Y)\otimes FZ}\big)\circ   \big(\varphi_2(\mathrm{id}_{X\otimes Y})\otimes FZ\big)  \, = \,  \\
\big(\varphi_2(a_{X,Y,Z}\circ r^{-1}_{(X\otimes Y)\otimes Z})\big)\circ  \big(r_{F\big((X\otimes Y)\otimes Z\big)}\big)\circ \big(\varphi_2(\mathrm{id}_{(X\otimes Y)\otimes Z})\big)\circ  \big(\varphi_2(\mathrm{id}_{X\otimes Y})\otimes FZ\big)
\end{align*}
autrement dit \eqref{ferrrrrr} est un diagramme commutatif. 

Donc $(F,m^F)\colon\xymatrix@C-7pt{\mathcal{M}\ar[r]&\mathcal{N}'}$ est un morphisme lax et unitaire entre $2$-groupes. 

Enfin vu que $\mathcal{N}(\mathcal{M}')$ est un ensemble simplicial faiblement $2$-cosquelettique,  pour montrer que $\mathcal{N}(F) = \varphi$ il suffit de vérifier que les morphismes des ensembles simpliciaux tronqués $\tau_{2}^{*}\big(\mathcal{N}(F)\big)$ et $\tau_{2}^{*}(\varphi)$ sont égaux.

On a évidement que $\tau_{1}^{*}\big(\mathcal{N}(F)\big) = \tau_{1}^{*}(\varphi)$. Montrons que si:
$$
\xymatrix@C+5pt{ A\otimes B \ar[r]^-{\xi} & C\,;}
$$
est un $2$-simplexe de $\mathcal{N}(\mathcal{M})$ alors $\mathcal{N}(F)_2\big(\xi\big) = \varphi_2(\xi)$.
 
Vu qu'on a un carré commutatif:
$$
\xymatrix@C+35pt@R+5pt{
F(A\otimes B) \otimes \mathbb{1} & \big( F(A) \otimes F(B) \big)\otimes \mathbb{1}\ar[l]_-{\varphi_2(\mathrm{id}_{A\otimes B}) \otimes \mathbb{1}}\\
F(A\otimes B) \ar[u]^-{r_{F(A\otimes B)}} & F(A)\otimes F(B) \ar[l]^-{\varphi_2(\mathrm{id}_{A\otimes B}) }
\ar[u]_-{r_{F(A)\otimes F(B)}}
}
$$
on a que:
\begin{equation*}
\begin{split}
\mathcal{N}(F)_2\big(\xi\big) \, & = \, F(\xi)\circ m^F_{A\otimes B} \, = \, \varphi_2\big(\xi\circ r^{-1}_{A\otimes B}\big)\circ r_{F(A\otimes B)}\circ \varphi_2(\mathrm{id}_{A\otimes B})\\
\, & = \, \varphi_2\big(\xi\circ r^{-1}_{A\otimes B}\big)\circ \big(\varphi_2(\mathrm{id}_{A\otimes B})\otimes \mathbb{1}\big)\circ r_{FA\otimes FB}
\end{split}
\end{equation*}

Si d'un autre on considère l'image par $\varphi$ du $3$-simplexe de $\mathcal{N}(\mathcal{M})$:
$$
\xymatrix@R=3pt@C+25pt{
C&  A\otimes B\ar[l]_-{\xi} \\&\\&\\
 &  A\otimes (B\otimes \mathbb{1}) \ar[uuu]_-{A\otimes r^{-1}_{B}} 
  \ar@{}[d]_-{a}|-{\text{\rotatebox[origin=c]{90}{$\cong$}}} \\ 
(A\otimes B)\otimes \mathbb{1} \ar[uuuu]^-{\xi \circ r^{-1}_{A\otimes B}}  & \big(A\otimes B\big)\otimes \mathbb{1} \ar[l]^-{(\mathrm{id}_{A\otimes B}) \otimes \mathbb{1}}
\ar@/^28pt/@<+23pt>@{-->}[uuuu]^-{r^{-1}_{A\otimes B}}}
$$
on obtient un diagramme commutatif de $\mathcal{N}(\mathcal{M}')$:
$$
\xymatrix@R=3pt@C+25pt{
F(C)&  F(A)\otimes F(B)\ar[l]_-{\varphi_2(\xi)} \\&\\&\\
 &  F(A)\otimes (F(B)\otimes \mathbb{1}) \ar[uuu]_-{FA\otimes r^{-1}_{FB}} 
  \ar@{}[d]_-{a}|-{\text{\rotatebox[origin=c]{90}{$\cong$}}} \\ 
F(A\otimes B)\otimes \mathbb{1} \ar[uuuu]^-{\varphi_2(\xi \circ r^{-1}_{A\otimes B})}  & \big(F(A)\otimes F(B)\big)\otimes \mathbb{1} \ar[l]^-{\varphi_2(\mathrm{id}_{A\otimes B}) \otimes \mathbb{1}}\ar@/^25pt/@<+30pt>@{-->}[uuuu]^-{r^{-1}_{FA\otimes FB}}
}
$$

En particulier on déduit que:
$$
\varphi_2\big(\xi\circ r^{-1}_{A\otimes B}\big)\circ \big(\varphi_2(\mathrm{id}_{A\otimes B})\otimes \mathbb{1}\big)\circ r_{FA\otimes FB} \, = \, \varphi_2(\xi)\,;
$$
autrement dit $\mathcal{N}(F)_2\big(\xi\big) = \varphi_2(\xi)$. 
\end{proof}

Montrons l'énoncé suivante:
\begin{proposition} \label{ilest} 
Si $\G$ est un $2$-groupe l'ensemble simplicial faiblement $2$-cosquelettique $\mathcal{N}\big(\mathcal{G}\big)$ (voir le Lemme \ref{3cosqull}) satisfait la condition d'extension de Kan en dimension $1\leq m\leq 3$ et satisfait la condition de minimalité en dimension $2$, autrement dit $\mathcal{N}\big(\mathcal{G}\big)$ est un $2$-groupoïde de Kan (voir le Corollaire \ref{groupidequi}).

En particulier, l'ensemble simplicial réduit $\mathcal{N}\big(\mathcal{G}\big)$ est un objet fibrant de la catégorie de modèles $(\simp_0,{\bf W}^{red}_2, {\bf mono},{\bf fib}^{red}_2)$ de la Proposition \ref{modred}.
\end{proposition}
\begin{proof}
Rappelons que l'ensemble simplicial $\mathcal{N}\big(\mathcal{G}\big)$ est faiblement $2$-cosquelettique d'après la Lemme \ref{3cosqull}. Montrons que $\mathcal{N}\big(\mathcal{G}\big)$ satisfait la condition d'extension de Kan en dimension $1\leq m\leq 3$ et satisfait la condition de minimalité en dimension $2$.

\emph{Condition d'extension de Kan en dimension $1$}: 

Montrons que si $X_{0}$, $X_{1}$ et $X_{2}$ sont des objets de $\G$ ils existent de morphismes de $\mathcal{G}$ de la forme:
$$
\xymatrix@C-8pt{X_2\otimes A_0 \ar[r]^-{\xi_0} & X_1}
\qquad \;\,
\xymatrix@C-8pt{X_2\otimes X_0 \ar[r]^-{\xi_1} & A_1}
\qquad \text{et} \qquad 
\xymatrix@C-8pt{A_2\otimes X _0\ar[r]^-{\xi_2} & X_1}\,,
$$
pour certains objets $A_0$, $A_1$ et $A_2$ de $\mathcal{G}$.

En effet si on prend des isomorphismes $\xymatrix@C-8pt{X_{2}\otimes X_{2}'\ar[r]^-{\varphi}&\mathbb{1}}$ et $\xymatrix@C-8pt{X_{0}'\otimes X_{0}\ar[r]^-{\psi}&\mathbb{1}}$, on définit: 
$$
\xi_{0}:  \bigg(\vcenter{\xymatrix@C-2pt{
X_{2}\otimes\big(X_{2}'\otimes X_{1}\big) \ar[rr]^-{a_{X_2,X_2',X_1}^{-1}} && \big(X_{2}\otimes X_{2}'\big)\otimes X_1 \ar[r]^-{\varphi\otimes X_1}&  \mathbb{1}\otimes X_{1}  \ar[r]^-{\ell_{X_1}^{-1}} & X_{1}
}} \bigg)\,,
$$
$$
\xi_{1}:  \bigg(\vcenter{\xymatrix@C-4pt{X_{2}\otimes X_{0} \ar[r]^-{\mathrm{id}} & X_{2}\otimes X_{0} }}\bigg)  \qquad\qquad \text{et}
$$
$$
\xi_{2}:   \bigg(\vcenter{\xymatrix@C-6pt{
\big(X_{1}\otimes X_{0}'\big)\otimes X_0 \ar[rr]^-{a_{X_1,X_0',X_0}} &&   X_{1}\otimes  \big(X_{0}'\otimes X_{0}\big) \ar[r]^-{X_1\otimes\psi}  & X_{1}\otimes \mathbb{1}  \ar[r]^-{r_{X_1}^{-1}} & X_{1}
 }} \bigg)  \,.
$$

\emph{Condition d'extension de Kan en dimension $2$ et de minimalité en dimension $2$}:

Pour montrer ces deux conditions, on doit vérifier que si on se donne trois des quatre morphismes $\xi_{0}$, $\xi_{1}$, $\xi_{2}$ et $\xi_{3}$ dans un diagramme de la forme:
\begin{equation}\label{cocol}
\vcenter{\xymatrix@R-16pt@C+8pt{
X_{03}  &X_{01}\otimes X_{13}   \ar[l]_-{\xi_{2}}   \\&\\
&X_{01} \otimes (X_{12}\otimes X_{23})  \ar[uu]_-{X_{01}\otimes \xi_{0} }  \ar@{}[d]_-{a}|-{\text{\rotatebox[origin=c]{90}{$\cong$}}}\\
X_{02}\otimes X_{23}  \ar[uuu]^-{\xi_{1}}   &(X_{01}\otimes X_{12})\otimes X_{23}  \ar[l]^-{\xi_{3}\otimes X_{23}} \,,}}
\end{equation}
on peut déterminer de façon unique un quatrième en faisant commutatif le diagramme \eqref{cocol}. Ceci est une conséquence du fait que $\mathcal{G}$ est un groupoïde et que les foncteurs de la forme $(A\otimes\,-\,)$ et $(\,-\,\otimes B)$ sont des équivalences de catégories.

\emph{Condition d'extension de Kan en dimension $3$}: 

Vu que $\mathcal{N}\big(\mathcal{G}\big)$ est un ensemble simplicial $3$-cosquelettique, l'ensemble des $4$-simplexes de $\G$ s'identifie canoniquement à l'ensemble des $5$-uplets $(\eta^{0},\eta^1,\eta^2,\eta^3,\eta^{4})$ de $3$-simplexes de $\G$:
\begin{equation}\label{33}
\eta^i \quad\;=\;\quad
\vcenter{\xymatrix@R-16pt@C+8pt{
X^i_{03}  &X^i_{01}\otimes X^i_{13}   \ar[l]_-{\xi^i_{2}}   \\&\\
&X^i_{01} \otimes (X^i_{12}\otimes X^i_{23})  \ar[uu]_-{X^i_{01}\otimes \xi^i_{0} }  \ar@{}[d]_-{a}|-{\text{\rotatebox[origin=c]{90}{$\cong$}}}\\
X^i_{02}\otimes X^i_{23}  \ar[uuu]^-{\xi^i_{1}}   &(X^i_{01}\otimes X^i_{12})\otimes X^i_{23}  \ar[l]^-{\xi^i_{3}\otimes X^i_{23}} \,,}}
\end{equation}
tels que:
\begin{equation}\label{condd}
\text{$d_{i}\eta^{j}=d_{j-1}\eta^{i}$\qquad  \text{toujours} \; {que} \qquad $0\leq i<j\leq 4$\,.}
\end{equation}

On vérifie alors qu'un $5$-uplets $(\eta^{0},\eta^1,\eta^2,\eta^3,\eta^{4})$ vérifiant la propriété \eqref{condd}, peut se ranger dans un cube:
\begin{equation}\label{cube}
\vcenter{\xymatrix@R-5pt@C+1pt{
&\centerdot\ar@{}[ddddl]|-{\eta_{2}}\ar@{}[rrrddd]|-{\eta_1}&&&\centerdot   \ar[lll]\\
\centerdot \ar[ur] &&&&\\&&&&\centerdot \ar[uu] \\
\centerdot \ar[uu] &\centerdot\ar@{}[rrd]|-{\mathrm{id}\otimes\eta_4} \ar[uuu]&&&\centerdot  \ar[lll]\ar[u]\\
\centerdot \ar[ru]\ar[u]&&\centerdot\ar[ll]&\centerdot\ar[l]\ar[ur]&}}
\qquad\vcenter{\xymatrix@R-18.5pt@C+1pt{
&\centerdot\ar@{}[rrdd]|-{\eta_3}&&&\centerdot\ar[lll]\ar@{}[ldddddd]|-{IV}\\ 
&&&&\\
\centerdot\ar@{}[rrdddd]|-{\eta_0\otimes\mathrm{id}} \ar[ruu]&&\centerdot\ar@{}[rddd]|-{II} \ar[ll] &\centerdot    \ar[l]\ar[ruu] &\\
&&&&\\
&&&&\centerdot  \ar@{}[ldddd]|-{V} \ar[uuuu] \\
&&\centerdot  \ar[uuu]&&\\
\centerdot\ar@{}[rrdd]|--{I}\ar[uuuu]&&\centerdot\ar@{}[rdd]|-{III} \ar[ll]\ar[u] &\centerdot  \ar[l]\ar[uuuu] \ar[uur]&\centerdot \ar[uu]\\
&&&&\\
\centerdot  \ar[uu]&&\centerdot \ar[ll] \ar[uu]&\centerdot \ar[l]\ar[uu]\ar[ruu]&}}
\end{equation}
où les carrés $I$, $II$, $III$, $IV$ et $V$, sont de diagrammes de la forme:
$$ 
\def\objectstyle{\scriptstyle}
\def\labelstyle{\scriptstyle}
\vcenter{\xymatrix@C-8pt@R+3pt{
(A\otimes B)\otimes C\ar@{}[rrd]|-{I}\ar[d]_-{a}\ar[rr]^-{(\mathrm{id}\otimes\xi)\otimes\mathrm{id}}&&\big(A\otimes (X\otimes Y)\big)\otimes C\ar[d]^-a\\A\otimes(B\otimes C)\ar[rr]_-{\mathrm{id}\otimes(\xi\otimes\mathrm{id})}&&A\otimes\big((X\otimes Y)\otimes C\big)\,\text{\normalsize{,}} }}
\qquad
\vcenter{\xymatrix@C+8pt@R+3pt{
\big( A\otimes (B \otimes C)\big)\otimes Y\ar[d]_-{a}\ar@{}[rd]|-{III}\ar[r]^-{a(a^{-1}\otimes \mathrm{id})}&(A\otimes B)\otimes (C\otimes D)\ar[d]^-{a}\\A\otimes \big((B\otimes C)\otimes D\big)\ar[r]_{\mathrm{id}\otimes a}&A\otimes \big(B\otimes (C\otimes D)\big)\,\text{\normalsize{,}}}}
$$
$$
\def\objectstyle{\scriptstyle}
\def\labelstyle{\scriptstyle}
\vcenter{\xymatrix@C-8pt@R+3pt{
A\otimes B\ar@{}[rd]|-{IV}\ar[d]_{\xi'\otimes\mathrm{id}}\ar[r]^-{\mathrm{id}\otimes \xi}& A\otimes (X\otimes Y)\ar[d]^-{\xi'\otimes\mathrm{id}}\\
(Z\otimes W)\otimes B\ar[r]_-{\mathrm{id}\otimes\xi}&(Z\otimes W)\otimes(X\otimes Y)\,\text{\normalsize{,}}}}
\qquad\vcenter{\xymatrix@C+9pt@R+3pt{
(A\otimes B)\otimes C \ar@{}[rd]|-{V}\ar[r]^-{(\mathrm{id}\otimes\mathrm{id})\otimes\xi}\ar[d]_-{a} & (A\otimes B)\otimes (X\otimes Y)\ar[d]^-{a} \\
A\otimes (B\otimes C) \ar[r]_-{\mathrm{id}\otimes(\mathrm{id}\otimes \xi)} &A\otimes \big(B\otimes(X\otimes Y)\big)}}
$$
$$
\def\objectstyle{\scriptstyle}
\def\labelstyle{\scriptstyle}
\text{et}\qquad\quad\vcenter{\xymatrix@C+3pt@R-8pt{
(A\otimes B)\otimes C\ar@{}[rddd]|-{II}\ar[dd]_{(\xi\otimes\mathrm{id})\otimes\mathrm{id}}\ar[r]^-{a}&
A\otimes (B\otimes C)\ar[ddd]^-{\xi\otimes(\mathrm{id}\otimes\mathrm{id})}\\&\\
\big((X\otimes Y)\otimes B\big)\otimes C\ar[d]_-{a\otimes\mathrm{id}}&\\
\big( X\otimes (Y\otimes B)\big)\otimes C\ar[r]_-{a(a^{-1}\otimes\mathrm{id})}&(X\otimes Y)\otimes(B\otimes C)\,\text{\normalsize{.}}}}
$$

On vérifie alors que tous les diagrammes dans les faces de \eqref{cube}, sont de diagrammes commutatifs. En effet, par hypothèse $\eta_{0}\otimes\mathrm{id}$, $\eta_{1}$, $\eta_{2}$, $\eta_{3}$ et $\mathrm{id}\otimes \eta_{4}$ sont commutatifs. Ensuite, puisque l'isomorphisme d'associativité $a$ est naturel, $I$, $II$ et $V$ sont commutatifs. D'un autre, $III$ est un diagramme de la forme \eqref{pentagon}, qu'on a supposé commutatif dans la définition d'une catégorie monoïdale. Enfin, $IV$ est commutatif car $\otimes$ est un foncteur en deux variables.

Pour montrer maintenant que $\mathcal{N}\big(\mathcal{G}\big)$ satisfait la condition d'extension de Kan en dimension $3$, remarquons que si on se donne seulement quatre des $3$-simplexes de $(\eta_{0},\dots,\eta_{4})$ sauf disons $\eta_{k}$, et on suppose que:
$$
\text{$d^2_{i}\eta_{j}=d^2_{j-1}\eta_{i}$\qquad  pour\quad $0\leq i<j\leq 4$ \quad et \quad $i,j\neq k$;}
$$ 
on peut toujours construire le cube \eqref{cube}, dont tous les diagrammes dans les faces sont commutatifs, sauf \emph{a priori} celui qui corresponde au $\eta_{k}$ manquant. 

Donc, le problème d'extension a une solution unique, car les foncteurs $(\mathrm{id}\otimes\,-\,)$ et $(\,-\,\otimes\mathrm{id})$ sont pleins et fidèles.
\end{proof}

Remarquons:

\begin{corollaire}\label{isoGAMMA}
Il existe un isomorphisme naturel de foncteurs:
\begin{equation}
\xymatrix@C+12pt{
\text{$2$-${\bf Grp}$} \drtwocell<\omit>{\phantom{aa}\Gamma}  \ar[r]^{\mathcal{N}} \ar[d]_{s} & \simp_0 \ar[d]^{\mathbb{\Omega}_\star} \\
{\bf Grpd} \ar[r]_{\mathrm{N}} & \simp }
\end{equation}
où $s\colon\xymatrix@C-10pt{\text{$2$-${\bf Grp}$} \ar[r] & {\bf Grpd}}$ est le foncteur groupoïde sous-jacent et $\mathbb{\Omega}_\star\colon\xymatrix@C-10pt{\simp_0\ar[r] & \simp}$ est le foncteur espace de lacet simplicial sur $\star$ de \eqref{lacetssimpisi}.
\end{corollaire}
\begin{proof}
Notons pour commencer que:
$$
\mathrm{N}\big(s(\mathcal{G})\big)_0 \, = \, 
\Big\{  \;  A \; \Big| \; \text{$A$ est un objet de $\mathcal{G}$} \; \Big\} \, = \, 
\mathbb{\Omega}_\star \big(\mathcal{N}(\mathcal{G})\big)_0\;,
$$
$$
\mathrm{N}\big(s(\mathcal{G})\big)_1 \, = \, \Big\{  \;  f\colon\xymatrix@C-10pt{A\ar[r]&B} \; \Big| \; \text{$f$ est un morphisme de $\mathcal{G}$} \; \Big\}
$$
$$
\text{et}
$$
$$
\mathbb{\Omega}_\star \big(\mathcal{N}(\mathcal{G})\big)_1 \, = \, \Big\{ \; \xi\colon\xymatrix@C-10pt{\mathbb{1}\otimes X\ar[r]&Y} \; \Big| \; \text{$\xi$ est un morphisme de $\mathcal{G}$} \; \Big\}\,,
$$
où les morphismes faces et dégénérescences de $\mathrm{N}\big(s(\mathcal{G})\big)$ sont donnés par: 
$$
\text{$d_0(f)=B$, \;  $d_1(f)=A$ \; et \; $s_0(A)=\mathrm{id}_A$}
$$
 et ceux de $\mathbb{\Omega}_\star \big(\mathcal{N}(\mathcal{G})\big)$ par 
$$
\text{$d_0(\xi)=X$, \; $d_1(\xi)=Y$ \; et \; $s_0(A)=\ell^{-1}_A$}\,.
$$

On définit la fonction:
$$
\xymatrix@C+18pt{ \mathbb{\Omega}_\star \big(\mathcal{N}(\mathcal{G})\big)_0 \ar[r]^{(\Gamma_\mathcal{G})_0} & \mathrm{N}\big(s(\mathcal{G})\big)_0}
$$
comme l'identité et la fonction:
$$
\xymatrix@C+18pt{ \mathbb{\Omega}_\star \big(\mathcal{N}(\mathcal{G})\big)_1 \ar[r]^{(\Gamma_\mathcal{G})_1} & \mathrm{N}\big(s(\mathcal{G})\big)_1}
$$
par la règle:
$$
\bigg(\xymatrix@C-10pt{\mathbb{1}\otimes X\ar[r]^-{\xi}&Y}\bigg)\qquad \longmapsto \qquad  \bigg(\xymatrix@C-4pt{Y\ar[r]^{\xi^{-1}}&\mathbb{1} \otimes X \ar[r]^{\ell_X^{-1}} &X}\bigg)\,.
$$

On obtient ainsi un morphisme d'ensembles simpliciaux tronqués:
\begin{equation}\label{troncationferferencore}
\xymatrix@C+10pt{\tau_{1}^{*}\Big(\mathbb{\Omega}_\star \big(\mathcal{N}(\mathcal{G})\big)\Big)\ar[r] &\tau_{1}^{*} \Big(\mathrm{N}\big(s(\mathcal{G})\big)\Big)}
\end{equation}
où $\tau_{1}^{*}\colon\xymatrix@C-10pt{\simp\ar[r]&\simp_{\leq 1}}$ est le foncteur de troncation. 

Vu que d'après la Proposition \ref{ilest}, le Lemmes \ref{desgroupDD} et le Corollaire \ref{1grpKAN} les ensembles simpliciaux $\mathbb{\Omega}_\star \big(\mathcal{N}(\mathcal{G})\big)$ et $\mathrm{N}\big(s(\mathcal{G})$ sont de $1$-groupoïdes de Kan, $\mathbb{\Omega}_\star \big(\mathcal{N}(\mathcal{G})\big)$ et $\mathrm{N}\big(s(\mathcal{G})$ sont des ensembles simpliciaux faiblement $1$-cosquelettiques. Donc pour étendre le morphisme \eqref{troncationferferencore} de façon unique en un morphisme d'ensembles simpliciaux: 
$$
\xymatrix@C+8pt{ \mathbb{\Omega}_\star \big(\mathcal{N}(\mathcal{G})\big)\ar[r]^-{\Gamma_\mathcal{G}} & \mathrm{N}\big(s(\mathcal{G})\big)}\,,
$$ 
il suffit de vérifier que si on se donne un $2$-simplexe de $\mathbb{\Omega}_\star \big(\mathcal{N}(\mathcal{G})\big)$:
\begin{equation}\label{ultimasssss}
\vcenter{\xymatrix@R-16pt@C+8pt{
Z  &\mathbb{1}\otimes Y \ar[l]_-{\xi_{2}}   \\&\\
&\mathbb{1} \otimes (\mathbb{1}\otimes X)  \ar[uu]_-{\mathbb{1}\otimes \xi_{0} }  \ar@{}[d]_-{a}|-{\text{\rotatebox[origin=c]{90}{$\cong$}}}\\
\mathbb{1}\otimes X  \ar[uuu]^-{\xi_{1}}   &(\mathbb{1}\otimes \mathbb{1})\otimes X  \ar[l]^-{\ell^{-1}_{\mathbb{1}}\otimes X} \,,}}
\end{equation}
alors on a un diagramme commutatif de $\mathcal{G}$:
$$
\xymatrix@C-1pt@R-10pt{
  && Y \ar[rd]^{\xi_0^{-1}}&& \\
  &\mathbb{1}\otimes Y\ar[ur]^-{\ell_Y^{-1}}&   &\mathbb{1}\otimes X\ar[rd]^-{\ell_X^{-1}}& \\
  Z\ar[ru]^{\xi_2^{-1}}\ar[rr]_-{\xi_1^{-1}}&& \mathbb{1}\otimes X\ar[rr]_-{\ell_X^{-1}}&& X\,.
 }
$$

Pour cela remarquons que d'après les diagrammes commutatifs:
$$
\vcenter{\xymatrix@C+8pt{
\mathbb{1}\otimes X \ar[r]^-{\xi_0} \ar[d]_{\ell_{\mathbb{1}\otimes X}}& Y \ar[d]^-{\ell_Y} \\
 \mathbb{1}\otimes (\mathbb{1}\otimes X) \ar[r]_-{\mathbb{1}\otimes \xi_0} & \mathbb{1} \otimes X \,,} }
\qquad 
\vcenter{\xymatrix@C-4pt@R+12pt{
& \mathbb{1}\otimes X \ar[rd]^-{\ell_{\mathbb{1}\otimes X}}\ar[dl]_-{\ell_\mathbb{1}\otimes X}&\\
(\mathbb{1}\otimes \mathbb{1})\otimes X \ar[rr]_-{a_{\mathbb{1},\mathbb{1},X}}& & \mathbb{1}\otimes (\mathbb{1}\otimes X)}}
$$
et \eqref{ultimasssss}, on a que:
\begin{align*}
\ell_X^{-1}\circ \xi_0^{-1}\circ \ell_Y^{-1}\circ \xi_2^{-1} \, &= \, (\xi_2\circ\ell_Y\circ\xi_0\circ\ell_X)^{-1} \\
\, &= \, \big(\xi_2\circ (\mathbb{1}\otimes \xi_0) \circ (a_{\mathbb{1},\mathbb{1},X}) \circ (\ell_\mathbb{1}\otimes X) \circ\ell_X\big)^{-1} \\
\, &= \, (\xi_1\circ \ell_X)^{-1} \, = \, \ell_X^{-1}\circ \xi_1^{-1}\,.
\end{align*}

\emph{$\Gamma$ est une transformation naturelle:}

Si $F\colon\xymatrix@C-10pt{\mathcal{G}\ar[r]& \mathcal{H}}$ est un morphisme de $2$-groupes il suffit de montrer qu'on a un carré commutatif:
$$
\xymatrix@C+15pt@R-3pt{ 
\mathbb{\Omega}_\star \big(\mathcal{N}(\mathcal{G})\big)_n\ar[r]^-{(\Gamma_\mathcal{G})_n} \ar[d]_-{\mathbb{\Omega}_\star\,\mathcal{N}(F)_n}& \mathrm{N}\big(s(\mathcal{G})\big)_n \ar[d]^-{\mathrm{N} \, s(F)_n}\\
\mathbb{\Omega}_\star \big(\mathcal{N}(\mathcal{H})\big)_n\ar[r]_-{(\Gamma_\mathcal{H})_n} & \mathrm{N}\big(s(\mathcal{H})\big)_n\,,}
$$ 
pour $0\leq n\leq 1$.

En effet si $n=0$ le carré commute évidement. D'un autre si $n=1$ et $\xi\colon\xymatrix@C-10pt{\mathbb{1}\otimes X\ar[r]& Y}$ est un $1$-simplexe de $\mathbb{\Omega}_{\star}(\mathcal{N} \mathcal{H})$ on a que les morphismes:
\begin{align*}
\big(\Gamma_{\mathcal{H}}\big)_1\circ \big( \mathbb{\Omega}_{\star}(\mathcal{N} \, F)\big)_1\Big(\xymatrix@C-10pt{\mathbb{1}\otimes X\ar[r]^-{\xi} & Y}\Big) \, & = \, 
\big(\Gamma_{\mathcal{H}}\big)_1 \Big(\xymatrix@C-6pt{\mathbb{1}\otimes F(X)\ar[r]^-{m_{\mathbb{1},X}^F} & F(\mathbb{1}\otimes X) \ar[r]^-{F(\xi)} & FY}\Big) \\
\, & = \, \Big(\xymatrix@C-3pt{
FY \ar[r]^-{F(\xi)^{-1}} &  F(\mathbb{1}\otimes X)  \ar[r]^-{(m_{\mathbb{1},X}^F)^{-1}}  &   \mathbb{1}\otimes F(X)  \ar[r]^-{\ell_{FX}^{-1}}& FX }\Big) 
\end{align*}
$$\text{et}$$
\begin{align*}
\big(\mathrm{N}(s\, F)\big)_1 \circ  \big(\Gamma_{\mathcal{H}}\big)_1 \Big(\xymatrix@C-10pt{\mathbb{1}\otimes X\ar[r]^-{\xi} & Y}\Big) \, & = \, 
\big(\mathrm{N}(s\, F)\big)_1 \Big(\xymatrix@C-6pt{Y\ar[r]^{\xi^{-1}} & \mathbb{1} \otimes X \ar[r]^{\ell_X^{-1}} & X}\Big) \\
\, & = \, \Big(\xymatrix@C-3pt{
FY \ar[r]^-{F(\xi)^{-1}} &  F(\mathbb{1}\otimes X)  \ar[r]^-{F(\ell_X)^{-1}}  &   FX }\Big)   
\end{align*}
sont égaux parce qu'on a un triangle commutatif:
$$
\xymatrix@R-5pt{
F(\mathbb{1}\otimes X) \ar[rr]^{m_{\mathbb{1},X}^F} && \mathbb{1}\otimes FX \\ 
&FX\ar[ul]^-{F(\ell_X)}\ar[ru]_-{\ell_{FX}}&}
$$
\end{proof}

L'énoncé réciproque de la Proposition \ref{ilest} a été montré par Duskin dans \cite{duskin}:

\begin{lemme} \label{ilestres}
Si $Z$ est un ensemble simplicial réduit lequel est un $2$-groupoïde de Kan, alors il existe un $2$-groupe $\G$ et un isomorphisme d'ensembles simpliciaux $Z\,\cong\,\mathcal{N}\big(\mathcal{G}\big)$. En particulier un ensemble simplicial réduit $Z$ est un $2$-groupoïde de Kan si et seulement si $Z$ est isomorphe au nerf d'un $2$-groupe (voir la Proposition \ref{ilest}).
\end{lemme}

Montrons finalement:

\begin{corollaire}\label{grho}
Si $0\leq i\leq 1$ on a des isomorphismes naturels de foncteurs:
$$
\xymatrix@C-8pt{
\big(\pi_i  \ar@{=>}[r]^-{\alpha^i} & \pi_{i+1}\circ \mathcal{N}\big) \colon \text{$2$-${\bf Grp}$} \ar[rr]& & {\bf Grp}\,.
}
$$

Donc un morphisme de $2$-groupes $F$ est une équivalence faible de $2$-groupes si et seulement si, $\mathcal{N}(F)$ est une $2$-équivalence faible d'ensembles simpliciaux réduits si et seulement si, $\mathcal{N}(F)$ est une $\infty$-équivalence faible d'ensembles simpliciaux.

En particulier, le foncteur nerf pour les $2$-groupes $\mathcal{N}\colon\xymatrix@C-10pt{\text{$2$-${\bf Grp}$} \ar[r] & \simp_0}$ induit un foncteur essentiellement surjectif:
$$
\xymatrix@C+10pt{\text{$2$-$h{\bf Grp}$} \ar[r]^{h\mathcal{N}} & \mathrm{Ho}_2\big(\simp_0\big)}
$$
de la catégorie homotopique des $2$-groupes vers la catégorie des $2$-types d'homotopie réduits (la catégorie homotopique des $2$-groupes).
\end{corollaire}
\begin{proof} D'après la Proposition \ref{ilest} si $\mathcal{G}$ est un $2$-groupe l'ensemble simplicial $\mathcal{N}\mathcal{G}$ est un complexe de Kan. Donc si $0\leq i\leq 1$ on peut définir les groupes $\pi_{i+1}\big(\mathcal{N}\mathcal{G},\star\big)$ de façon combinatoire par les règles des Propositions \ref{kankan} et \ref{kankan2}. 

{\bf $\alpha^0$}: L'ensemble sous-jacent au groupe $\pi_{1}\big(\mathcal{N}\mathcal{G},\star\big)$ est égal à l'ensemble d'objets de $\mathcal{G}$ soumis à la relation d'équivalence:
\begin{equation}\label{rela}
X\sim Y \qquad\text{s'il existe un morphisme de $\mathcal{G}$ de la forme} \;\;\xymatrix@C-5pt{X\otimes \mathbb{1}\ar[r]&Y}\;.
\end{equation}

Remarquons que si on considère la relation à isomorphisme près dans l'ensemble des objets du groupoïde $\mathcal{G}$:
\begin{equation}
X\sim' Y \qquad\text{s'il existe un morphisme de $\mathcal{G}$ de la forme} \;\;\xymatrix@C-5pt{X\ar[r]&Y}\;, 
\end{equation}
il se suit que $\sim \, = \, \sim'$ parce que on a l'isomorphisme naturel $\ell_{?}$. Donc $\pi_0(\mathcal{G}) = \pi_{1}\big(\mathcal{N}\mathcal{G}\big)$ en tant que ensembles.

D'un autre si on pose $[A]$ pour note la classe d'un objet de $\mathcal{G}$ à isomorphisme près on a que $[A]\cdot [B] = [C]$ dans le groupe $\pi_{1}\big(\mathcal{N}\mathcal{G}\big)$ si et seulement s'il existe un morphisme de $\mathcal{G}$ de la forme:
$$
\xymatrix{A\otimes B \ar[r]^-{\xi} & C\,.} 
$$

Vu qu'on a des morphismes:
$$
\xymatrix{A\otimes \mathbb{1} \ar[r]^-{\ell_A^{-1}} & A}   \qquad \text{et}\qquad \xymatrix{\mathbb{1}\otimes A \ar[r]^-{r_A^{-1}} & A\,,} 
$$
pour $A$ un objet quelconque de $\mathcal{G}$, il se suit que l'élément neutre du groupe $\pi_{1}\big(\mathcal{N}\mathcal{G}\big)$ est la classe de l'objet $s_{0}(\star)=\mathbb{1}$. De la même façon si $X$ et $Y$ sont des objets de $\mathcal{G}$ on a que $[X]\cdot [Y] = [X\otimes Y]$ parce qu'on a un morphisme:
$$
\xymatrix@C-4pt{(X)\otimes (Y)\ar[r]^-{\mathrm{id}}&(X\otimes Y)}\,.
$$

Donc $\pi_0(\mathcal{G}) = \pi_{1}\big(\mathcal{N}\mathcal{G}\big)$ en tant que groupes, c'est-à-dire:
$$
\xymatrix@C-8pt{
\big(\pi_0  \ar@{=>}[r]^-{\alpha^0} & \pi_{1}\circ \mathcal{N}\big) \colon \text{$2$-${\bf Grp}$} \ar[rr]& & {\bf Grp}\,
}
$$
est la transformation identité.

{\bf $\alpha^1$}: Remarquons que l'ensemble sous-jacent au groupe $\pi_{2}\big(\mathcal{N}\mathcal{G}\big)$ est égal à l'ensemble des morphismes de $\mathcal{G}$ de la forme $\xi:\xymatrix@C-8pt{\mathbb{1}\otimes \mathbb{1}\ar[r]&\mathbb{1}}$, soumis à la relation d'équivalence: 
$$
\text{$\xi\sim \xi'$\qquad s'il existe}\qquad
\def\objectstyle{\scriptstyle}
\def\labelstyle{\scriptstyle}
\vcenter{\xymatrix@R=3pt@C+15pt{
\mathbb{1}& \mathbb{1}\otimes \mathbb{1}  \ar[l]_-{\xi'} \\&\\&\\
 &  \mathbb{1}\otimes \big( \mathbb{1}\otimes \mathbb{1} \big) \ar[uuu]_-{\mathbb{1}\otimes \ell_\mathbb{1}^{-1}} 
  \ar@{}[d]_-{a}|-{\text{\rotatebox[origin=c]{90}{$\cong$}}} \\ 
\mathbb{1}\otimes \mathbb{1} \ar[uuuu]^-{\ell_\mathbb{1}^{-1}}  & \big(\mathbb{1}\otimes \mathbb{1}\big)\otimes \mathbb{1} \ar[l]^-{\xi \otimes \mathbb{1}}}}
\qquad\text{commutatif.}
$$

Vu que $(\mathbb{1}\otimes \ell_{\mathbb{1}}^{-1})\circ a=\ell_{\mathbb{1}\otimes\mathbb{1}}^{-1}$, il se suit de la naturalité de l'isomorphisme $\ell_{?}$, que $\xi\sim\xi'$ si et seulement si $\xi\otimes\mathbb{1}=\xi'\otimes\mathbb{1}$. Donc, $\xi\sim\xi'$ si et seulement si $\xi=\xi'$ parce que le foncteur $-\otimes\mathbb{1}$ est pleinement fidèle. Autrement dit l'ensemble sous-jacent à $\pi_{2}\big(\mathcal{N}\mathcal{G}\big)$ est simplement l'ensemble des morphismes de $\mathcal{G}$ de la forme $\xi:\xymatrix@C-5pt{\mathbb{1}\otimes\mathbb{1}\ar[r]&\mathbb{1}}$.

Remarquons aussi que si $\xi,\xi',\xi'':\xymatrix@C-5pt{\mathbb{1}\otimes\mathbb{1}\ar[r]&\mathbb{1}}$ sont des morphismes de $\mathcal{G}$ alors $\xi\,\cdot\,\xi'=\xi''$ dans $\pi_{2}\big(\mathcal{N}\mathcal{G}\big)$ s'il existe un diagramme commutatif de $\mathcal{G}$:
$$
\vcenter{\xymatrix@R=3pt@C+15pt{
\mathbb{1}& \mathbb{1}\otimes \mathbb{1}  \ar[l]_-{\xi''} \\&\\&\\
 &  \mathbb{1}\otimes \big( \mathbb{1}\otimes \mathbb{1} \big) \ar[uuu]_-{\mathbb{1}\otimes \ell_\mathbb{1}^{-1}} 
  \ar@{}[d]_-{a}|-{\text{\rotatebox[origin=c]{90}{$\cong$}}} \\ 
\mathbb{1}\otimes \mathbb{1} \ar[uuuu]^-{\xi'}  & \big(\mathbb{1}\otimes \mathbb{1}\big)\otimes \mathbb{1} \ar[l]^-{\xi \otimes \mathbb{1}}}}
$$

Montrons que la transformation:
$$
\xymatrix@C-8pt{
\big(\pi_1  \ar@{=>}[r]^-{\alpha^1} & \pi_{2}\circ \mathcal{N}\big) \colon \text{$2$-${\bf Grp}$} \ar[rr]& & {\bf Grp}\,
}
$$
définie par la famille des fonctions:
$$
\vcenter{\xymatrix@R=1pt{
\pi_{1}\big(\mathcal{G}\big)\ar[r]^-{\alpha^1_{\mathcal{G}}}  &\pi_{2}\big(\mathcal{N}\mathcal{G}\big)\\
\qquad \varphi\qquad \ar@{|->}[r]&\qquad \varphi\circ\ell^{-1}_{\mathbb{1}}= \varphi\circ r_{\mathbb{1}}^{-1} \qquad }}\qquad\quad\text{si \; $\mathcal{G}$ \; est \; un \; $2$-groupe\,,}
$$
est bien un isomorphisme naturel de groupes.

En effet pour commencer on note que $\alpha_\mathcal{G}^{1}$ est une bijection dont l'inverse est définie par la règle: 
$$
\big(\alpha^1_{\mathcal{G}}\big)^{-1}(\xi)=\xi\circ\ell_{\mathbb{1}}=\xi\circ r_{\mathbb{1}}\,.
$$

D'un autre si $\varphi$ et $\psi$ sont des éléments de $\pi_1(\mathcal{G})$ remarquons qu'on a un diagramme commutatif:
$$
\vcenter{\xymatrix@R=3pt@C+15pt{
\mathbb{1}& \mathbb{1}\otimes \mathbb{1}  \ar[l]_-{(\psi\circ\varphi)\circ\ell_\mathbb{1}^{-1}} \\&\\&\\
 &  \mathbb{1}\otimes \big( \mathbb{1}\otimes \mathbb{1} \big) \ar[uuu]_-{\mathbb{1}\otimes \ell_\mathbb{1}^{-1}} 
  \ar@{}[d]_-{a}|-{\text{\rotatebox[origin=c]{90}{$\cong$}}} \\ 
\mathbb{1}\otimes \mathbb{1} \ar[uuuu]^-{\psi\circ\ell_{\mathbb{1}}^{-1}}  & \big(\mathbb{1}\otimes \mathbb{1}\big)\otimes \mathbb{1} \ar[l]^-{(\varphi\circ\ell^{-1}_{\mathbb{1}}) \otimes \mathbb{1}}}}
$$
parce qu'on a les diagrammes commutatifs de $\mathcal{G}$:
$$
\vcenter{\xymatrix{
\mathbb{1} \ar[r]^-{\varphi} & \mathbb{1} \\
\mathbb{1}\otimes \mathbb{1}\ar[r]_-{\varphi\otimes \mathbb{1}}\ar[u]^-{\ell_{\mathbb{1}}^{-1}}&\mathbb{1}\otimes \mathbb{1}\ar[u]_-{\ell_{\mathbb{1}}^{-1}}
}}\;,\qquad\;
\vcenter{\xymatrix@C-12pt{
(\mathbb{1}\otimes \mathbb{1})\otimes \mathbb{1} \ar[rd]_-{\ell^{-1}_{\mathbb{1}} \otimes \mathbb{1} } \ar[rr]^-{(\varphi\circ\ell^{-1}_{\mathbb{1}})\otimes\mathbb{1}} && 
\mathbb{1}\otimes\mathbb{1} \\&\mathbb{1}\otimes\mathbb{1}\ar[ru]_{\varphi\otimes\mathbb{1}}&}}
\quad\,\text{et} \quad \,
\vcenter{\xymatrix@C-18pt{
(\mathbb{1}\otimes \mathbb{1})\otimes \mathbb{1} \ar[rd]_-{\ell^{-1}_{\mathbb{1}} \otimes \mathbb{1} } \ar[rr]^-{a} && 
\mathbb{1}\otimes(\mathbb{1}\otimes\mathbb{1}) \ar[ld]^-{\mathbb{1}\otimes\ell_{\mathbb{1}}^{-1} } \\&\mathbb{1}\otimes\mathbb{1}&}}\,.
$$

Autrement dit on a que:
$$
\alpha^{1}_{\mathcal{G}}(\varphi)\,\cdot\, \alpha^{1}_{\mathcal{G}}(\psi) \, = \, \alpha^{1}_{\mathcal{G}}(\psi\circ \varphi) 
$$

Donc $\alpha^1_{\mathcal{G}}$ est un isomorphisme de groupes. 

Vérifions que $\alpha^1$ est une transformation naturelle: Si $F: \xymatrix@C-9pt{\mathcal{G}\ar[r]&\mathcal{H}}$ est un morphisme de $2$-groupes on a par définition que:
$$
\vcenter{\xymatrix@R=1pt@C-5pt{
\pi_{2}\big(\mathcal{N}\mathcal{G}\big)\ar[r]^-{F_{*}}  &\pi_{2}\big(\mathcal{N}\mathcal{H},\big)\\
\quad\qquad \xi\quad\qquad \ar@{|->}[r]&\qquad F(\xi)\circ m^F_{\mathbb{1},\mathbb{1}} \qquad }}
\qquad\text{et}\qquad
\vcenter{\xymatrix@R=1pt@C-5pt{
\pi_{1}\big(\mathcal{G}\big)\ar[r]^-{F_{*}}  &\pi_{1}\big(\mathcal{H},\big)\,.\\
\quad\qquad \varphi\quad\qquad \ar@{|->}[r]&\qquad F(\varphi)\qquad }}
$$

Il se suit qu'on a un carré commutatif: 
$$
\xymatrix{
\pi_1(\mathcal{G}) \ar[r]^-{\alpha^1_{\mathcal{G}}} \ar[d]_-{F_*} & \pi_2\big(\mathcal{N}(\mathcal{G})\big)\ar[d]^-{F_*}\\
\pi_1(\mathcal{H}) \ar[r]_-{\alpha^1_{\mathcal{G}}} & \pi_2\big(\mathcal{N}(\mathcal{H})\big)
}
$$
parce qu'on a un triangle commutatif:
$$
\xymatrix@R-18pt@C-6pt{&\mathbb{1}\ar[rd]^-{F(l_{\mathbb{1}})}\ar[ld]_-{l_{\mathbb{1}}}&\\
F(\mathbb{1}\otimes\mathbb{1})\ar[rr]_-{m_{\mathbb{1},\mathbb{1}}}&&\mathbb{1}\otimes\mathbb{1}\,.}
$$

On déduit de ce qu'on a montré qu'un morphisme de $2$-groupes $F$ est une équivalence faible de $2$-groupes si et seulement si $\mathcal{N}(F)$ est une $2$-équivalence faible d'ensembles simpliciaux réduits. Pour montrer que $\mathcal{N}(F)$ est une $2$-équivalence faible d'ensembles simpliciaux réduits si et seulement si $\mathcal{N}(F)$ est une $\infty$-équivalence faible d'ensembles simpliciaux remarquons que les groupes $\pi_{i}\big(\mathcal{N}\mathcal{G}\big)$ sont nuls pour $i=0$ et $i\geq 3$.

Donc le foncteur nerf pour les $2$-groupes $\mathcal{N}\colon\xymatrix@C-10pt{\text{$2$-${\bf Grp}$} \ar[r] & \simp_0}$ induit un foncteur:
$$
\xymatrix@C+10pt{\text{$2$-$h{\bf Grp}$} \ar[r]^{h\mathcal{N}} & \mathrm{Ho}_2\big(\simp_0\big)}\,.
$$
lequel est essentiellement surjectif d'après la Proposition \ref{propmay} et le Lemme \ref{ilestres}.
\end{proof}

\renewcommand{\thesubsection}{\S\thesection.\arabic{subsection}}
\subsection{}\;
\renewcommand{\thesubsection}{\thesection.\arabic{subsection}}

Si $\mathcal{G}$ est un $2$-groupe d'après  la Proposition \ref{ilest} l'ensemble simplicial $\mathcal{N}(\G)$ est un objet fibrant de la catégorie de modèles $(\simp_{0},{\bf W}^{red}_{2}, {\bf mono},{\bf fib}_{2}^{red})$, il se suit du Corollaire \ref{dimminus} que si $X$ est un ensemble simplicial quelconque l'ensemble simplicial $\underline{\mathrm{Hom}}_{\simp_0}\big(X ,\mathcal{N}(\G)\big)$ est un objet fibrant de la catégorie de modèles $(\simp,{\bf W}_{1}, {\bf mono},{\bf fib}_{1})$ \emph{i.e.} un complexe de Kan dont les groupes d'homotopie $\pi_i$ sont nuls pour $i\geq 2$. 

On est tenter a penser que l'ensemble simplicial $\underline{\mathrm{Hom}}_{\simp_0}\big(X,\mathcal{N}(\G)\big)$ est un $1$-groupoïde de Kan \emph{i.e.} le nerf d'un groupoïde (voir le Lemme \ref{HOMgrpdred}), par contre c'est pas le cas en général:

\begin{lemme}\label{elelemmafer}
Si $\G$ est un $2$-groupe vérifiant la propriété:
\begin{quote}
Pour tout objet $X$ de $\mathcal{G}$ il existe au moins un morphisme non identité $f$ de source $X$.
\end{quote}
l'ensemble simplicial $\underline{\mathrm{Hom}}_{\simp_0}\big(\Delta^1\big/ {\bf sq}_0\Delta^1 ,\mathcal{N}(\G)\big)$ n'est pas un $1$-groupoïde de Kan.
\end{lemme}
\begin{proof}
Remarquons que l'ensemble $\underline{\mathrm{Hom}}_{\simp_0}\big(\Delta^1\big/ {\bf sq}_0\Delta^1 ,\mathcal{N}(\G)\big)_n$ s'identifie à l'ensemble:
\begin{equation}\label{delta1deltan}
\mathrm{Hom}_{\simp} \Big(\big(\Delta^1\times  \Delta^n\big) \big/ \big({\bf sq}_0\Delta^1 \times  \Delta^n \big), \mathcal{N}(\G)\Big)\,.
\end{equation}

Pour comprendre la fonction:
\begin{equation} 
\xymatrix@C+10pt{
\mathrm{Hom}_{\simp}\bigg(\Delta^{2},\underline{\mathrm{Hom}}_{\simp_0}\big(\Delta^1\big/ {\bf sq}_0\Delta^1 ,\mathcal{N}(\G)\big)\bigg) \ar[r]^{\alpha^{1,k}} &\mathrm{Hom}_{\simp}\bigg(\Lambda^{2,k},\underline{\mathrm{Hom}}_{\simp_0}\big(\Delta^1\big/ {\bf sq}_0\Delta^1 ,\mathcal{N}(\G)\big)\bigg)}
\end{equation}
on va décrire l'ensemble \eqref{delta1deltan} pour $0\leq n \leq 2$.

{\bf $n=0$}: L'ensemble simplicial $\big(\Delta^1\times  \Delta^0\big) \big/ \big({\bf sq}_0\Delta^1 \times  \Delta^0 \big)$ est isomorphe à l'ensemble simplicial $\Delta^1\big/ {\bf sq}_0\Delta^1$. Un morphisme de $\Delta^1\big/ {\bf sq}_0\Delta^1$ vers $\mathcal{N}(\G)$ est déterminé par un $1$-simplexe de $\mathcal{N}(\G)$ c'est-à-dire par un objet de $\G$.

{\bf $n=1$}: On veut décrire l'ensemble simplicial $\big(\Delta^1\times  \Delta^1\big) \big/ \big({\bf sq}_0\Delta^1 \times  \Delta^1 \big)$; pour cela remarquons que:
$$
\Big(\Delta^1_0\times  \Delta^1_0\Big) \Big\backslash \Big(\big({\bf sq}_0\Delta^1\big)_0 \times  \Delta^1_0\Big) \; = \; \emptyset \;, 
$$
\begin{align*}
\Big(\Delta^1_1\times  \Delta^1_1\Big) \Big\backslash \Big(\big({\bf sq}_0\Delta^1\big)_1 \times  \Delta^1_1\Big) \; = \; 
\Bigg\{ \, \big(\mathrm{id}_{[1]},\mathrm{id}_{[1]}\big), \; \big(\mathrm{id}_{[1]},\delta_0\sigma_0\big), \; \big(\mathrm{id}_{[1]},\delta_1\sigma_0\big) \;  \Bigg\}
\end{align*}
$$\text{et}$$
\begin{align*}
\Big(\Delta^1_2\times  \Delta^1_2\Big) \Big\backslash \Big(\big({\bf sq}_0\Delta^1\big)_2 \times  \Delta^1_2\Big) \; = \;  \Bigg\{ \, 
&\big(\sigma_0,\sigma_0\big), \; \big(\sigma_0,\sigma_1\big), \; \big(\sigma_0,\delta_0\sigma_0\sigma_0\big), \; \big(\sigma_0,\delta_1\sigma_0\sigma_0\big) \;  \\
&\big(\sigma_1,\sigma_0\big), \; \big(\sigma_1,\sigma_1\big), \; \big(\sigma_1,\delta_0\sigma_0\sigma_0\big), \; \big(\sigma_1,\delta_1\sigma_0\sigma_0\big) \; \Bigg\}\,;
\end{align*}
en plus dans l'ensemble simplicial $\Delta^1\times\Delta^1$ on a que:
\begin{align*}
s_0\big(\mathrm{id}_{[1]},\mathrm{id}_{[1]}\big) \; &= \;  \big(\sigma_0,\sigma_0\big) & 
s_1\big(\mathrm{id}_{[1]},\mathrm{id}_{[1]}\big) \; &= \;  \big(\sigma_1,\sigma_1\big) \\
s_0\big(\mathrm{id}_{[1]},\delta_0\sigma_0\big) \; &= \;  \big(\sigma_0,\delta_0\sigma_0\sigma_0\big) & 
s_1\big(\mathrm{id}_{[1]},\delta_0\sigma_0\big) \; &= \;  \big(\sigma_1,\delta_0\sigma_0\sigma_0\big) \\
s_0\big(\mathrm{id}_{[1]},\delta_1\sigma_0\big) \; &= \;  \big(\sigma_0,\delta_1\sigma_0\sigma_0\big) &
s_1\big(\mathrm{id}_{[1]},\delta_1\sigma_0\big) \; &= \;  \big(\sigma_1,\delta_1\sigma_0\sigma_0\big)\,. \\
\end{align*}

Plus encore, 
\begin{lemme}\label{squeletiqueee}
$\big(\Delta^1\times  \Delta^1\big) \big/ \big({\bf sq}_0\Delta^1 \times  \Delta^1 \big)$ est un ensemble simplicial $2$-squelettique, c'est-à-dire tous ses $n$-simplexes sont dégénérés pour $n\geq 3$. 
\end{lemme}
\begin{proof}
En effet, il suffit de montrer cette affirmation pour l'ensemble simplicial $\Delta^1\times  \Delta^1$: Remarquons que si $m\geq 3$ les éléments de $\Delta^1_m\times  \Delta^1_m$ sont les couples:
$$
\big(\varphi\sigma_{i_1}\dots\sigma_{i_{m-1}},\psi\sigma_{j_1}\dots\sigma_{j_{m-1}}\big)
$$ 
où $\varphi$ et $\psi$ sont des morphismes $\xymatrix@C-8pt{[1]\ar[r]&[1]}$ de $\Delta$ et $0 \leq i_1< \dots < i_{m-1}\leq 2$, $0 \leq j_1< \dots < j_{m-1} \leq 2$ sont de suites des entiers.   Il se suit en particulier que $m-2\leq i_{m-1},j_{m-1}\leq m-1$.

Dans le cas où $i_{m-1}=j_{m-1}=k$ on a que:
$$
s_k\big(\varphi\sigma_{i_1}\dots\sigma_{i_{m-2}},\psi\sigma_{j_1}\dots\sigma_{j_{m-2}}\big) \, = \, \big(\varphi\sigma_{i_1}\dots\sigma_{i_{m-1}},\psi\sigma_{j_1}\dots\sigma_{j_{m-1}}\big)
$$ 

Supposons que $m-2 = i_{m-1} < j_{m-1} = m-1$. On vérifie d'abord que:
$$
i_1=0, \, i_2=1, \, \dots \, i_{m-2}=m-3, \, i_{m-1}=m-2 \qquad\quad  \text{et}  \qquad\quad  m-3 \leq j_{m-2}\leq m-2\,.
$$

Vu que $\sigma_k\sigma_l = \sigma_{l-1}\sigma_{k}$ toujours que $k< l$, dans le cas où $j_{m-2}=m-2$ on a que:

\begin{align*}
\big(\varphi\sigma_{i_1}\dots\sigma_{i_{m-2}}\sigma_{i_{m-1}},\psi\sigma_{j_1}\dots\sigma_{j_{m-2}}\sigma_{j_{m-1}}\big) \; & = \; 
\big(\varphi\sigma_{0}\dots\sigma_{m-3}\sigma_{m-2},\psi\sigma_{j_1}\dots\sigma_{m-2}\sigma_{m-1}\big)  \\ 
\; & = \;  \big(\varphi\sigma_{0}\dots\sigma_{m-3}\sigma_{m-2},\psi\sigma_{j_1}\dots\sigma_{m-2}\sigma_{m-2}\big)\\
\; & = \;  s_{m-2}\big(\varphi\sigma_{0}\dots\sigma_{m-3},\psi\sigma_{j_1}\dots\sigma_{m-2}\big)\,,
\end{align*}
et dans le cas où $j_{m-2}=m-3$ on a que:
\begin{align*}
\big(\varphi\sigma_{i_1}\dots\sigma_{i_{m-2}}\sigma_{i_{m-1}},\psi\sigma_{j_1}\dots\sigma_{j_{m-2}}\sigma_{j_{m-1}}\big) \; & = \; 
\big(\varphi\sigma_{0}\dots\sigma_{m-3}\sigma_{m-2},\psi\sigma_{j_1}\dots\sigma_{m-3}\sigma_{m-1}\big)  \\ 
\; & = \;  \big(\varphi\sigma_{0}\dots\sigma_{m-3}\sigma_{m-3},\psi\sigma_{j_1}\dots\sigma_{m-2}\sigma_{m-3}\big)\\
\; & = \;  s_{m-3}\big(\varphi\sigma_{0}\dots\sigma_{m-3},\psi\sigma_{j_1}\dots\sigma_{m-2}\big)\,.
\end{align*}

Donc $\Delta^1\times\Delta^1$ n'a pas de $n$-simplexes non-dégénérés si $n\geq 3$.
\end{proof}

On déduit que l'ensemble simplicial quotient $\big(\Delta^1\times  \Delta^1\big) \big/ \big({\bf sq}_0\Delta^1 \times  \Delta^1 \big)$ a un seul $0$-simplexe, trois $1$-simplexes non-dégénérés associés aux $1$-simplexes de $\Delta^1\times  \Delta^1$:
$$
\big(\mathrm{id}_{[1]},\mathrm{id}_{[1]}\big), \; \big(\mathrm{id}_{[1]},\delta_0\sigma_0\big), \; \big(\mathrm{id}_{[1]},\delta_1\sigma_0\big)
$$
et deux $2$-simplexes non-dégénérés:
$$
\big(\sigma_0,\sigma_1\big), \; \big(\sigma_1,\sigma_0\big)\,.
$$

Vu que dans l'ensemble simplicial $\Delta^1\times\Delta^1$ on a que: 
\begin{align*}
d_0\big(\sigma_0,\sigma_1\big) \; &= \; \big(\mathrm{id}_{[1]},\delta_0\sigma_0\big) \qquad 
d_1\big(\sigma_0,\sigma_1\big) \; = \; \big(\mathrm{id}_{[1]},\mathrm{id}_{[1]}\big) \qquad 
d_2\big(\sigma_0,\sigma_1\big) \; = \; \big(\delta_1\sigma_0,\mathrm{id}_{[1]}\big) \\
d_0\big(\sigma_1,\sigma_0\big) \; &= \; \big(\delta_0\sigma_0,\mathrm{id}_{[1]}\big) \qquad 
d_1\big(\sigma_1,\sigma_0\big) \; = \; \big(\mathrm{id}_{[1]},\mathrm{id}_{[1]}\big) \qquad 
d_2\big(\sigma_1,\sigma_0\big) \; = \; \big(\mathrm{id}_{[1]},\delta_1\sigma_0\big)\;.
\end{align*}
il se suit que se donner un élément de l'ensemble:
$$
\mathrm{Hom}_{\simp} \Big(\big(\Delta^1\times  \Delta^1\big) \big/ \big({\bf sq}_0\Delta^1 \times  \Delta^1 \big), \mathcal{N}(\G)\Big)\,,
$$
équivaut à se donner trois objets et deux morphismes de $\G$:
\begin{equation}
\xymatrix@C+10pt{
\mathbb{1} \otimes X_{(\mathrm{id}_{[1]},\delta_0\sigma_0)} \ar[r]^-{\xi_{(\sigma_0,\sigma_1)}}  & X_{(\mathrm{id}_{[1]},\mathrm{id}_{[1]})} 
& X_{(\mathrm{id}_{[1]},\delta_1\sigma_0)}\otimes \mathbb{1} \ar[l]_-{\xi_{(\sigma_1,\sigma_0)}} 
}
\end{equation}

$n=2$: Pour décrire l'ensemble simplicial $\big(\Delta^1\times  \Delta^2\big) \big/ \big({\bf sq}_0\Delta^1 \times  \Delta^2 \big)$ remarquons que:
$$
\Big(\Delta^1_0\times  \Delta^2_0\Big) \Big\backslash \Big(\big({\bf sq}_0\Delta^1\big)_0 \times  \Delta^2_0\Big) \; = \; \emptyset \;, 
$$
\begin{align*}
\Big(\Delta^1_1\times  \Delta^2_1\Big) \Big\backslash \Big(\big({\bf sq}_0\Delta^1\big)_1 \times  \Delta^2_1\Big) \; = \;  \Bigg\{ \, 
&\big(\mathrm{id}_{[1]}, \varphi\big) \; \Bigg| \; \varphi \in \big\{ \delta_0, \, \delta_1, \, \delta_2, \, \delta_1\delta_0\sigma_0, \, \delta_2\delta_0\sigma_0, \, \delta_2\delta_1\sigma_0\big\} \; \Bigg\}
\end{align*}

\begin{align*}
\Big(\Delta^1_2\times  \Delta^2_2\Big) \Big\backslash \Big(\big({\bf sq}_0\Delta^1\big)_2 \times  \Delta^2_2\Big) \; = \;  \Bigg\{ \, 
\big(\psi, \varphi\big) \; \Bigg| \; &\psi \in \Big\{ \sigma_0, \, \sigma_1 \Big\} \quad\text{et}\quad  
\varphi \in \Big\{ \mathrm{id}_{[2]}, \delta_0\sigma_0, \, \delta_1\sigma_0, \, \delta_2\sigma_0, \, \\
                            &\delta_0\sigma_1, \, \delta_1\sigma_1, \, \delta_2\sigma_1, \, \delta_1\delta_0\sigma_0\sigma_0, \, \delta_2\delta_0\sigma_0\sigma_0, \, \delta_2\delta_1\sigma_0\sigma_0\Big\} \; \Bigg\}
\end{align*}
$$\text{et}$$
\begin{align*}
\Big(\Delta^1_3\times  \Delta^2_3\Big) \Big\backslash \Big(\big({\bf sq}_0\Delta^1\big)_3 \times  \Delta^2_3\Big) \; = \;  \Bigg\{ \, 
\big(\psi, \varphi\big) \; \Bigg| \; &\psi \in \Big\{ \sigma_0\sigma_0, \, \sigma_1\sigma_0, \, \sigma_1\sigma_1 \Big\} \quad\text{et}\quad  
\varphi \in \Big\{ \sigma_0, \, \sigma_1, \, \sigma_2, \, \delta_0\sigma_0\sigma_0, \, \\
                            & \delta_0\sigma_1\sigma_0, \,  \delta_0\sigma_1\sigma_1, \, \delta_1\sigma_0\sigma_0, \, \delta_1\sigma_1\sigma_0, \, \delta_1\sigma_1\sigma_1, \, 
                            \delta_2\sigma_0\sigma_0, \,  \\
                            &\delta_2\sigma_1\sigma_0, \, \delta_2\sigma_1\sigma_1, \, \delta_1\delta_0\sigma_0\sigma_0\sigma_0, \,  \delta_2\delta_0\sigma_0\sigma_0\sigma_0, \, 
                             \delta_2\delta_1\sigma_0\sigma_0\sigma_0\big\} \; \Bigg\}\,.
\end{align*}

D'un autre on montre que $\big(\Delta^1\times  \Delta^2\big) \big/ \big({\bf sq}_0\Delta^1 \times  \Delta^2 \big)$ est un ensemble simplicial $3$-squelettique par la même méthode que dans la preuve du Lemme \ref{squeletiqueee} et on vérifie que:
{\footnotesize{
\begin{align*}
s_0\big(\mathrm{id}_{[1]},\delta_0\big) \; &= \;  \big(\sigma_0,\delta_0\sigma_0\big) & 
s_0\big(\mathrm{id}_{[1]},\delta_1\delta_0\sigma_0\big) \; &= \;  \big(\sigma_0,\delta_1\delta_0\sigma_0\sigma_0\big)\\
s_0\big(\mathrm{id}_{[1]},\delta_1\big) \; &= \;  \big(\sigma_0,\delta_1\sigma_0\big) &
s_0\big(\mathrm{id}_{[1]},\delta_2\delta_0\sigma_0\big) \; &= \;  \big(\sigma_0,\delta_2\delta_0\sigma_0\sigma_0\big)\\
s_0\big(\mathrm{id}_{[1]},\delta_2\big) \; &= \;  \big(\sigma_0,\delta_2\sigma_0\big) &
s_0\big(\mathrm{id}_{[1]},\delta_2\delta_1\sigma_0\big) \; &= \;  \big(\sigma_0,\delta_2\delta_1\sigma_0\sigma_0\big)\\
&&&\\
s_1\big(\mathrm{id}_{[1]},\delta_0\big) \; &= \;  \big(\sigma_1,\delta_0\sigma_1\big) &
s_1\big(\mathrm{id}_{[1]},\delta_1\delta_0\sigma_0\big) \; &= \;  \big(\sigma_1,\delta_1\delta_0\sigma_0\sigma_0\big)\\
s_1\big(\mathrm{id}_{[1]},\delta_1\big) \; &= \;  \big(\sigma_1,\delta_1\sigma_1\big) &
s_1\big(\mathrm{id}_{[1]},\delta_2\delta_0\sigma_0\big) \; &= \;  \big(\sigma_1,\delta_2\delta_0\sigma_0\sigma_0\big)\\
s_1\big(\mathrm{id}_{[1]},\delta_2\big) \; &= \;  \big(\sigma_1,\delta_2\sigma_1\big) &
s_1\big(\mathrm{id}_{[1]},\delta_2\delta_1\sigma_0\big) \; &= \;  \big(\sigma_1,\delta_2\delta_1\sigma_0\sigma_0\big)
\end{align*}}}

{\footnotesize{
\begin{align*}
s_0\big(\sigma_0,\mathrm{id}_{[2]}\big) \; &= \;  \big(\sigma_0\sigma_0,\sigma_0\big) & 
s_0\big(\sigma_1,\mathrm{id}_{[2]}\big) \; &= \;  \big(\sigma_1\sigma_0,\sigma_0\big) \\
s_0\big(\sigma_0,\delta_0\sigma_0\big) \; &= \;  \big(\sigma_0\sigma_0,\delta_0\sigma_0\sigma_0\big) & 
s_0\big(\sigma_1,\delta_0\sigma_0\big) \; &= \;  \big( \sigma_1\sigma_0,\delta_0\sigma_0\sigma_0\big) \\
s_0\big(\sigma_0,\delta_1\sigma_0\big) \; &= \;  \big(\sigma_0\sigma_0,\delta_1\sigma_0\sigma_0\big) & 
s_0\big(\sigma_1,\delta_1\sigma_0\big) \; &= \;  \big(\sigma_1\sigma_0,\delta_1\sigma_0\sigma_0\big) \\
s_0\big(\sigma_0, \delta_2\sigma_0\big) \; &= \;  \big(\sigma_0\sigma_0, \delta_2\sigma_0\sigma_0\big) & 
s_0\big(\sigma_1, \delta_2\sigma_0\big) \; &= \;  \big(\sigma_1\sigma_0, \delta_2\sigma_0\sigma_0\big) \\
s_0\big(\sigma_0,\delta_0\sigma_1\big) \; &= \;  \big(\sigma_0\sigma_0,\delta_0\sigma_1\sigma_0\big) & 
s_0\big(\sigma_1,\delta_0\sigma_1\big) \; &= \;  \big(\sigma_1\sigma_0,\delta_0\sigma_1\sigma_0\big) \\
s_0\big(\sigma_0,\delta_1\sigma_1\big) \; &= \;  \big(\sigma_0\sigma_0,\delta_1\sigma_1\sigma_0\big) & 
s_0\big(\sigma_1,\delta_1\sigma_1\big) \; &= \;  \big(\sigma_1\sigma_0,\delta_1\sigma_1\sigma_0\big) \\
s_0\big(\sigma_0,\delta_2\sigma_1\big) \; &= \;  \big(\sigma_0\sigma_0,\delta_2\sigma_1\sigma_0\big) & 
s_0\big(\sigma_1,\delta_2\sigma_1\big) \; &= \;  \big(\sigma_1\sigma_0,\delta_2\sigma_1\sigma_0\big) \\
s_0\big(\sigma_0,\delta_1\delta_0\sigma_0\sigma_0\big) \; &= \;  \big( \sigma_0\sigma_0,\delta_1\delta_0\sigma_0\sigma_0\sigma_0 \big) & 
s_0\big(\sigma_1,\delta_1\delta_0\sigma_0\sigma_0\big) \; &= \;  \big(\sigma_1\sigma_0,\delta_1\delta_0\sigma_0\sigma_0\sigma_0\big) \\
s_0\big(\sigma_0,\delta_2\delta_0\sigma_0\sigma_0\big) \; &= \;  \big(\sigma_0\sigma_0,\delta_2\delta_0\sigma_0\sigma_0\sigma_0\big) & 
s_0\big(\sigma_1,\delta_2\delta_0\sigma_0\sigma_0\big) \; &= \;  \big(\sigma_1\sigma_0,\delta_2\delta_0\sigma_0\sigma_0\sigma_0\big) \\
s_0\big(\sigma_0,\delta_2\delta_1\sigma_0\sigma_0\big) \; &= \;  \big(\sigma_0\sigma_0,\delta_2\delta_1\sigma_0\sigma_0\sigma_0\big) & 
s_0\big(\sigma_1,\delta_2\delta_1\sigma_0\sigma_0\big) \; &= \;  \big(\sigma_1\sigma_0,\delta_2\delta_1\sigma_0\sigma_0\sigma_0\big) 
\end{align*}}}

{\footnotesize{
\begin{align*}
s_1\big(\sigma_0,\mathrm{id}_{[2]}\big) \; &= \;  \big(\sigma_0\sigma_0,\sigma_1\big) & 
s_1\big(\sigma_1,\mathrm{id}_{[2]}\big) \; &= \;  \big(\sigma_1\sigma_1,\sigma_1\big) \\
s_1\big(\sigma_0,\delta_0\sigma_0\big) \; &= \;  \big(\sigma_0\sigma_0,\delta_0\sigma_0\sigma_0\big) & 
s_1\big(\sigma_1,\delta_0\sigma_0\big) \; &= \;  \big( \sigma_1\sigma_1,\delta_0\sigma_0\sigma_0\big) \\
s_1\big(\sigma_0,\delta_1\sigma_0\big) \; &= \;  \big(\sigma_0\sigma_0,\delta_1\sigma_0\sigma_0\big) & 
s_1\big(\sigma_1,\delta_1\sigma_0\big) \; &= \;  \big(\sigma_1\sigma_1,\delta_1\sigma_0\sigma_0\big) \\
s_1\big(\sigma_0, \delta_2\sigma_0\big) \; &= \;  \big(\sigma_0\sigma_0, \delta_2\sigma_0\sigma_0\big) & 
s_1\big(\sigma_1, \delta_2\sigma_0\big) \; &= \;  \big(\sigma_1\sigma_1, \delta_2\sigma_0\sigma_0\big) \\
s_1\big(\sigma_0,\delta_0\sigma_1\big) \; &= \;  \big(\sigma_0\sigma_0,\delta_0\sigma_1\sigma_1\big) & 
s_1\big(\sigma_1,\delta_0\sigma_1\big) \; &= \;  \big(\sigma_1\sigma_1,\delta_0\sigma_1\sigma_1\big) \\
s_1\big(\sigma_0,\delta_1\sigma_1\big) \; &= \;  \big(\sigma_0\sigma_0,\delta_1\sigma_1\sigma_1\big) & 
s_1\big(\sigma_1,\delta_1\sigma_1\big) \; &= \;  \big(\sigma_1\sigma_1,\delta_1\sigma_1\sigma_1\big) \\
s_1\big(\sigma_0,\delta_2\sigma_1\big) \; &= \;  \big(\sigma_0\sigma_0,\delta_2\sigma_1\sigma_1\big) & 
s_1\big(\sigma_1,\delta_2\sigma_1\big) \; &= \;  \big(\sigma_1\sigma_1,\delta_2\sigma_1\sigma_1\big) \\
s_1\big(\sigma_0,\delta_1\delta_0\sigma_0\sigma_0\big) \; &= \;  \big( \sigma_0\sigma_0,\delta_1\delta_0\sigma_0\sigma_0\sigma_0 \big) & 
s_1\big(\sigma_1,\delta_1\delta_0\sigma_0\sigma_0\big) \; &= \;  \big(\sigma_1\sigma_1,\delta_1\delta_0\sigma_0\sigma_0\sigma_0\big) \\
s_1\big(\sigma_0,\delta_2\delta_0\sigma_0\sigma_0\big) \; &= \;  \big(\sigma_0\sigma_0,\delta_2\delta_0\sigma_0\sigma_0\sigma_0\big) & 
s_1\big(\sigma_1,\delta_2\delta_0\sigma_0\sigma_0\big) \; &= \;  \big(\sigma_1\sigma_1,\delta_2\delta_0\sigma_0\sigma_0\sigma_0\big) \\
s_1\big(\sigma_0,\delta_2\delta_1\sigma_0\sigma_0\big) \; &= \;  \big(\sigma_0\sigma_0,\delta_2\delta_1\sigma_0\sigma_0\sigma_0\big) & 
s_1\big(\sigma_1,\delta_2\delta_1\sigma_0\sigma_0\big) \; &= \;  \big(\sigma_1\sigma_1,\delta_2\delta_1\sigma_0\sigma_0\sigma_0\big) 
\end{align*}
}}

{\footnotesize{
\begin{align*}
s_2\big(\sigma_0,\mathrm{id}_{[2]}\big) \; &= \;  \big(\sigma_1\sigma_0,\sigma_2\big) & 
s_2\big(\sigma_1,\mathrm{id}_{[2]}\big) \; &= \;  \big(\sigma_1\sigma_1,\sigma_2\big) \\
s_2\big(\sigma_0,\delta_0\sigma_0\big) \; &= \;  \big(\sigma_1\sigma_0,\delta_0\sigma_1\sigma_0\big) & 
s_2\big(\sigma_1,\delta_0\sigma_0\big) \; &= \;  \big( \sigma_1\sigma_1,\delta_0\sigma_1\sigma_0\big) \\
s_2\big(\sigma_0,\delta_1\sigma_0\big) \; &= \;  \big(\sigma_1\sigma_0,\delta_1\sigma_1\sigma_0\big) & 
s_2\big(\sigma_1,\delta_1\sigma_0\big) \; &= \;  \big(\sigma_1\sigma_1,\delta_1\sigma_1\sigma_0\big) \\
s_2\big(\sigma_0, \delta_2\sigma_0\big) \; &= \;  \big(\sigma_1\sigma_0, \delta_2\sigma_1\sigma_0\big) & 
s_2\big(\sigma_1, \delta_2\sigma_0\big) \; &= \;  \big(\sigma_1\sigma_1, \delta_2\sigma_1\sigma_0\big) \\
s_2\big(\sigma_0,\delta_0\sigma_1\big) \; &= \;  \big(\sigma_1\sigma_0,\delta_0\sigma_1\sigma_1\big) & 
s_2\big(\sigma_1,\delta_0\sigma_1\big) \; &= \;  \big(\sigma_1\sigma_1,\delta_0\sigma_1\sigma_1\big) \\
s_2\big(\sigma_0,\delta_1\sigma_1\big) \; &= \;  \big(\sigma_1\sigma_0,\delta_1\sigma_1\sigma_1\big) & 
s_2\big(\sigma_1,\delta_1\sigma_1\big) \; &= \;  \big(\sigma_1\sigma_1,\delta_1\sigma_1\sigma_1\big) \\
s_2\big(\sigma_0,\delta_2\sigma_1\big) \; &= \;  \big(\sigma_1\sigma_0,\delta_2\sigma_1\sigma_1\big) & 
s_2\big(\sigma_1,\delta_2\sigma_1\big) \; &= \;  \big(\sigma_1\sigma_1,\delta_2\sigma_1\sigma_1\big) \\
s_2\big(\sigma_0,\delta_1\delta_0\sigma_0\sigma_0\big) \; &= \;  \big( \sigma_1\sigma_0,\delta_1\delta_0\sigma_0\sigma_0\sigma_0 \big) & 
s_2\big(\sigma_1,\delta_1\delta_0\sigma_0\sigma_0\big) \; &= \;  \big(\sigma_1\sigma_1,\delta_1\delta_0\sigma_0\sigma_0\sigma_0\big) \\
s_2\big(\sigma_0,\delta_2\delta_0\sigma_0\sigma_0\big) \; &= \;  \big(\sigma_1\sigma_0,\delta_2\delta_0\sigma_0\sigma_0\sigma_0\big) & 
s_2\big(\sigma_1,\delta_2\delta_0\sigma_0\sigma_0\big) \; &= \;  \big(\sigma_1\sigma_1,\delta_2\delta_0\sigma_0\sigma_0\sigma_0\big) \\
s_2\big(\sigma_0,\delta_2\delta_1\sigma_0\sigma_0\big) \; &= \;  \big(\sigma_1\sigma_0,\delta_2\delta_1\sigma_0\sigma_0\sigma_0\big) & 
s_2\big(\sigma_1,\delta_2\delta_1\sigma_0\sigma_0\big) \; &= \;  \big(\sigma_1\sigma_1,\delta_2\delta_1\sigma_0\sigma_0\sigma_0\big) \,.
\end{align*}
}}

Autrement dit l'ensemble simplicial quotient $\big(\Delta^1\times  \Delta^2\big) \big/ \big({\bf sq}_0\Delta^1 \times  \Delta^2\big)$ a un seul $0$-simplexe, six $1$-simplexes non-dégénérés associés aux $1$-simplexes de $\Delta^1\times  \Delta^2$ suivants: 
$$
\big(\mathrm{id}_{[1]},\delta_0\big),\;\;
\big(\mathrm{id}_{[1]},\delta_1\big) ,\; \;
\big(\mathrm{id}_{[1]},\delta_2\big),\;\;
\big(\mathrm{id}_{[1]},\delta_1\delta_0\sigma_0\big) ,\;\;
\big(\mathrm{id}_{[1]},\delta_2\delta_0\sigma_0\big)\quad \text{et}\quad
\big(\mathrm{id}_{[1]},\delta_2\delta_1\sigma_0\big);\,
$$
huit $2$-simplexes non-dégénérés:
$$
\big(\sigma_0,\mathrm{id}_{[1]}\big),\;\;
\big(\sigma_1,\mathrm{id}_{[1]}\big) ,\;\;
\big(\sigma_0,\delta_0\sigma_1\big),\;\;
\big(\sigma_0,\delta_1\sigma_1\big),\;\;
\big(\sigma_0,\delta_2\sigma_1\big),\;\;
\big(\sigma_1,\delta_0\sigma_0\big),\;\;
\big(\sigma_1,\delta_1\sigma_0\big)\quad \text{et}\quad
\big(\sigma_1,\delta_2\sigma_0\big);\;
$$
et trois $3$-simplexes non-dégénérés:
$$
\big(\sigma_0\sigma_0,\sigma_2\big),\;\;
\big(\sigma_1\sigma_0,\sigma_1\big)\quad \text{et}\quad
\big(\sigma_1\sigma_1,\sigma_0\big)\,.
$$

Vu qu'on a les égalités:
{\footnotesize{
\begin{align*}
d_0\big(\sigma_0,\mathrm{id}_{[1]}\big) &= \big(\mathrm{id}_{[1]},\delta_0\big),\;\;
d_1\big(\sigma_0,\mathrm{id}_{[1]}\big) = \big(\mathrm{id}_{[1]},\delta_1\big)\quad\text{et}\quad
d_2\big(\sigma_0,\mathrm{id}_{[1]}\big) = \big(\delta_1\sigma_0,\delta_2\big),\\
d_0\big(\sigma_1,\mathrm{id}_{[1]}\big) &= \big(\delta_0\sigma_0,\delta_0\big),\;\;
d_1\big(\sigma_1,\mathrm{id}_{[1]}\big) = \big(\mathrm{id}_{[1]},\delta_1\big)\quad\text{et}\quad
d_2\big(\sigma_1,\mathrm{id}_{[1]}\big) = \big(\mathrm{id}_{[1]},\delta_2\big), \\
d_0\big(\sigma_0,\delta_0\sigma_1\big) &= \big(\mathrm{id}_{[1]},\delta_1\delta_0\sigma_0\big),\;\;
d_1\big(\sigma_0,\delta_0\sigma_1\big) = \big(\mathrm{id}_{[1]},\delta_0\big)\quad\text{et}\quad
d_2\big(\sigma_0,\delta_0\sigma_1\big) = \big(\delta_1\sigma_0,\delta_0\big),\\
d_0\big(\sigma_0,\delta_1\sigma_1\big) &= \big(\mathrm{id}_{[1]},\delta_1\delta_0\sigma_0\big),\;\;
d_1\big(\sigma_0,\delta_1\sigma_1\big) = \big(\mathrm{id}_{[1]},\delta_1\big)\quad\text{et}\quad
d_2\big(\sigma_0,\delta_1\sigma_1\big) = \big(\delta_1\sigma_0,\delta_1\big),\\
d_0\big(\sigma_0,\delta_2\sigma_1\big) &= \big(\mathrm{id}_{[1]},\delta_2\delta_0\sigma_0\big),\;\;
d_1\big(\sigma_0,\delta_2\sigma_1\big) = \big(\mathrm{id}_{[1]},\delta_2\big)\quad\text{et}\quad
d_2\big(\sigma_0,\delta_2\sigma_1\big) = \big(\delta_1\delta_0,\delta_2\big),\\
d_0\big(\sigma_1,\delta_0\sigma_0\big) &= \big(\delta_0\sigma_0,\delta_0\big),\;\;
d_1\big(\sigma_1,\delta_0\sigma_0\big) = \big(\mathrm{id}_{[1]},\delta_0\big)\quad\text{et}\quad
d_2\big(\sigma_1,\delta_0\sigma_0\big) = \big(\mathrm{id}_{[1]},\delta_2\delta_0\sigma_0\big),\\
d_0\big(\sigma_1,\delta_1\sigma_0\big) &= \big(\delta_0\sigma_0,\delta_1\big),\;\;
d_1\big(\sigma_1,\delta_1\sigma_0\big) = \big(\mathrm{id}_{[1]},\delta_1\big)\quad\text{et}\quad
d_2\big(\sigma_1,\delta_1\sigma_0\big) = \big(\mathrm{id}_{[1]},\delta_2\delta_1\sigma_0\big),\\
d_0\big(\sigma_1,\delta_2\sigma_0\big) &= \big(\delta_0\sigma_0,\delta_2\big),\;\;
d_1\big(\sigma_1,\delta_2\sigma_0\big) = \big(\mathrm{id}_{[1]},\delta_2\big)\quad\text{et}\quad
d_2\big(\sigma_1,\delta_2\sigma_0\big) = \big(\mathrm{id}_{[1]},\delta_2\delta_1\sigma_0\big),
\end{align*}
}}

{\footnotesize{
\begin{align*}
d_0\big(\sigma_0\sigma_0,\sigma_2\big) = \big(\sigma_0,\delta_0\sigma_1\big) 
d_1\big(\sigma_0\sigma_0,\sigma_2\big) = \big(\sigma_0,\delta_1\sigma_1\big)
d_2\big(\sigma_0\sigma_0,\sigma_2\big) = \big(\sigma_0,\mathrm{id}_{[1]}\big)
d_3\big(\sigma_0\sigma_0,\sigma_2\big) = \big(\delta_1\sigma_0\sigma_0,\mathrm{id}_{[1]}\big)
\end{align*}
}}

{\footnotesize{
\begin{align*}
d_0\big(\sigma_1\sigma_0,\sigma_1\big) = \big(\sigma_1,\delta_0\sigma_0\big) 
d_1\big(\sigma_1\sigma_0,\sigma_1\big) = \big(\sigma_1,\mathrm{id}_{[1]}\big)
d_2\big(\sigma_1\sigma_0,\sigma_1\big) = \big(\sigma_0,\mathrm{id}_{[1]}\big)
d_3\big(\sigma_1\sigma_0,\sigma_1\big) = \big(\sigma_0,\delta_2\sigma_1\big)
\end{align*}
}}
$$\text{et}$$
{\footnotesize{
\begin{align*}
d_0\big(\sigma_1\sigma_1,\sigma_0\big) = \big(\delta_0\sigma_0\sigma_0,\mathrm{id}\big) 
d_1\big(\sigma_1\sigma_1,\sigma_0\big) = \big(\sigma_1,\mathrm{id}_{[1]}\big)
d_2\big(\sigma_1\sigma_1,\sigma_0\big) = \big(\sigma_1,\delta_1\sigma_0\big)
d_3\big(\sigma_1\sigma_1,\sigma_0\big) = \big(\sigma_1,\delta_2\sigma_0\big)\,{\normalsize{;}}
\end{align*}
}}
il se suit que se donner un élément de l'ensemble:
$$
\mathrm{Hom}_{\simp} \Big(\big(\Delta^1\times  \Delta^2\big) \big/ \big({\bf sq}_0\Delta^1 \times  \Delta^2 \big), \mathcal{N}(\G)\Big)\,,
$$
équivaut à se donner un diagramme en $\mathcal{G}$ de la forme:
\begin{equation}\label{nononoferfer}
{\footnotesize{\xymatrix@R=5pt@C+42pt{
\big(X_{(\mathrm{id}_{[1]},\delta_2\delta_1\sigma_0)}\otimes\mathbb{1}\big)\otimes\mathbb{1} \ar[r]^-{\xi_{(\sigma_1,\delta_2\sigma_0)}\otimes \mathbb{1}} 
 \ar@{}[rdddd]|-{\eta_{(\sigma_1\sigma_1,\sigma_0)}} &
X_{(\mathrm{id}_{[1]},\delta_2)}\otimes \mathbb{1}\ar[dddd]|-{\xi_{(\sigma_1,\mathrm{id}_{[1]})}}  \ar@{}[rdddd]|-{\eta_{(\sigma_1\sigma_0,\sigma_1)}}  & 
\big(\mathbb{1}\otimes X_{(\mathrm{id}_{[1]},\delta_2\delta_0\sigma_0)}\big)\otimes \mathbb{1}\ar[l]_-{\xi_{(\sigma_0,\delta_2\sigma_1)} \otimes \mathbb{1}}\\
X_{(\mathrm{id}_{[1]},\delta_2\delta_1\sigma_0)}\otimes\big(\mathbb{1}\otimes\mathbb{1}\big)   
\ar[ddd]_-{X_{(\mathrm{id}_{[1]},\delta_2\delta_1\sigma_0)\otimes r_{\mathbb{1}}^{-1}}}\ar@{}[u]^-{a}|-{\text{\rotatebox[origin=c]{90}{$\cong$}}}& & 
\mathbb{1}\otimes \big(X_{(\mathrm{id}_{[1]},\delta_2\delta_0\sigma_0)}\otimes \mathbb{1} \big) \ar[ddd]^-{\mathbb{1}\otimes \xi_{(\sigma_1,\delta_0\sigma_0)}} 
  \ar@{}[u]_-{a}|-{\text{\rotatebox[origin=c]{90}{$\cong$}}} \\ &&\\&&\\
 X_{(\mathrm{id}_{[1]},\delta_2\delta_1\sigma_0)}\otimes\mathbb{1}\ar[r]_-{\xi_{(\sigma_1,\delta_1\delta_0)}}
&X_{(\mathrm{id}_{[1]},\delta_1)}  \ar@{}[rdddd]|-{\eta_{(\sigma_0\sigma_0,\sigma_2)}} & 
\mathbb{1}\otimes X_{(\mathrm{id}_{[1]},\delta_0)}  \ar[l]|-{\xi_{(\sigma_0,\mathrm{id}_{[1]})}} \\&&\\&&\\
& &  \mathbb{1}\otimes \big(\mathbb{1} \otimes X_{(\mathrm{id}_{[1]},\delta_1\delta_0\sigma_0)}\big) \ar[uuu]_-{\mathbb{1}\otimes \xi_{(\sigma_0,\delta_0\sigma_1)}} 
  \ar@{}[d]^-{a}|-{\text{\rotatebox[origin=c]{90}{$\cong$}}} \\ 
&\mathbb{1} \otimes X_{(\mathrm{id}_{[1]},\delta_1\delta_0\sigma_0)}\ar[uuuu]^-{\xi_{(\sigma_0,\delta_1\sigma_1)}}  & \big(\mathbb{1}\otimes \mathbb{1}\big)\otimes X_{(\mathrm{id}_{[1]},\delta_1\delta_0\sigma_0)}\ar[l]^-{r_{\mathbb{1}}^{-1}\otimes X_{(\mathrm{id}_{[1]},\delta_1\delta_0\sigma_0)}}
}}}
\end{equation}

Finalement remarquons que dans l'ensemble simplicial $\underline{\mathrm{Hom}}_{\simp_0}\big(\Delta^1\big/ {\bf sq}_0\Delta^1 ,\mathcal{N}(\G)\big)$ les images par les fonctions $d_0$, $d_1$ et $d_2$ du $2$-simplexe \eqref{nononoferfer} sont les $1$-simplexes:
$$
\xymatrix@C+10pt{
\mathbb{1} \otimes X_{(\mathrm{id}_{[1]},\delta_1\delta_0\sigma_0)}  \ar[r]^-{\xi_{(\sigma_0,\delta_0\sigma_1)}} & X_{(\mathrm{id}_{[1]},\delta_0)}  
 & X_{(\mathrm{id}_{[1]},\delta_2\delta_0\sigma_0)}\otimes \mathbb{1}  \ar[l]_-{\xi_{(\sigma_1,\delta_0\sigma_0)}}\,,
}
$$
$$
\xymatrix@C+10pt{
\mathbb{1} \otimes X_{(\mathrm{id}_{[1]},\delta_1\delta_0\sigma_0)}  \ar[r]^-{\xi_{(\sigma_0,\delta_1\sigma_1)}} & X_{(\mathrm{id}_{[1]},\delta_1)}  
 & X_{(\mathrm{id}_{[1]},\delta_2\delta_1\sigma_0)}\otimes \mathbb{1}  \ar[l]_-{\xi_{(\sigma_1,\delta_1\sigma_0)}}\qquad \text{et}
}
$$
$$
\xymatrix@C+10pt{
\mathbb{1} \otimes X_{(\mathrm{id}_{[1]},\delta_2\delta_0\sigma_0)}  \ar[r]^-{\xi_{(\sigma_0,\delta_2\sigma_1)}} & X_{(\mathrm{id}_{[1]},\delta_0)}  
 & X_{(\mathrm{id}_{[1]},\delta_2\delta_1\sigma_0)}\otimes \mathbb{1}  \ar[l]_-{\xi_{(\sigma_1,\delta_2\sigma_0)}}\,,
}
$$
respectivement.

On déduit de cette description que si $0\leq k\leq 2$ et $\mathcal{G}$ est un $2$-groupe vérifiant les hypothèse de l'énoncé du Lemme \ref{elelemmafer} alors la fonction:
\begin{equation} \label{lafontitifeeer}
\xymatrix@C+10pt{
\mathrm{Hom}_{\simp}\bigg(\Delta^{2},\underline{\mathrm{Hom}}_{\simp_0}\big(\Delta^1\big/ {\bf sq}_0\Delta^1 ,\mathcal{N}(\G)\big)\bigg) \ar[r]^{\alpha^{1,k}} &\mathrm{Hom}_{\simp}\bigg(\Lambda^{2,k},\underline{\mathrm{Hom}}_{\simp_0}\big(\Delta^1\big/ {\bf sq}_0\Delta^1 ,\mathcal{N}(\G)\big)\bigg)}
\end{equation}
induite de l'inclusion $\xymatrix@C-3pt{\Lambda^{2,k}\,\ar@{^(->}[r]^-{\alpha^{m,k}}&\Delta^{2}}$ n'est pas injective; ce qu'implique que l'ensemble simplicial $\underline{\mathrm{Hom}}_{\simp_0}\big(\Delta^1\big/ {\bf sq}_0\Delta^1 ,\mathcal{N}(\G)\big)$ n'est pas un $1$-groupoïde de Kan.
\end{proof}

\renewcommand{\thesubsection}{\S\thesection.\arabic{subsection}}
\subsection{}\;\label{detrereduit}
\renewcommand{\thesubsection}{\thesection.\arabic{subsection}}

Si $X$ est un ensemble simplicial réduit et $\G$ est un $2$-groupe, \emph{un déterminant de $X$ à valeurs dans $\G$} est par définition un couple de fonctions $D=(D,T)$:
$$
\xymatrix@C+10pt{X_1\ar[r]^-{D} & \big\{\text{Objets de $\G$}\big\}} \qquad\text{et}\qquad \xymatrix{X_2\ar[r]^-{T} & \big\{\text{Morphismes de $\G$}\big\}}
$$
vérifiant les propriétés:
\begin{enumerate}
\item  (Compatibilité) Si $\xi\in X_2$ alors $T(\xi)$ est un morphisme de $\G$ de la forme:
$$
\xymatrix@+10pt{
D(d_2\xi)\otimes D(d_0 \xi) \ar[r]^-{T(\xi)} &  D(d_1\xi)   \,.}
$$
\item (Unitaire) On a que $D(s_0 \, \star) = \mathbb{1}$ et $T(s_0\circ s_0 \, \star) = \ell_{\mathbb{1}}^{-1} = r_{\mathbb{1}}^{-1}$.
\item (Associativité)\, Si $\eta\in X_{3}$ on a un diagramme commutatif:
$$
\vcenter{\xymatrix@R-10pt@C+35pt{
D(A_{03})&D(A_{01})\otimes D(A_{13})\ar[l]_-{T(d_2\,\eta)}\\&\\
&D(A_{01}) \otimes (D(A_{12})\otimes D(A_{23})) \ar@{}[d]_-{a}|-{\text{\rotatebox[origin=c]{90}{$\cong$}}}\ar[uu]_-{D(A_{01})\otimes T(d_0\, \eta) }\\
D(A_{02})\otimes D(A_{23}) \ar[uuu]^-{T(d_1\,\eta)}&(D(A_{01})\otimes D(A_{12}))\otimes D(A_{23}) \ar[l]^-{T(d_3\,\eta)\otimes D(A_{23})}\,,}}
$$
\begin{align*}
\text{où} \qquad 
A_{03} \, = \, & d_{1} d_{1} \eta \, = \, d_{1} d_{2} \eta \qquad
A_{01} \, = \,  d_{2} d_{2} \eta \, = \, d_{2} d_{3} \eta \qquad
A_{13} \, = \,  d_{1} d_{0} \eta \, = \, d_{0} d_{2} \eta \\
A_{02} \, = \, & d_{2} d_{1} \eta \, = \, d_{1} d_{3}  \eta \qquad
A_{23} \, = \,  d_{0} d_{0} \eta \, = \, d_{0} d_{1} \eta \qquad
A_{12} \, = \,  d_{2} d_{0} \eta \, = \, d_{0} d_{3} \eta \,.
\end{align*}
\end{enumerate}

On pose ${\bf det}_X(\G)$ pour noter l'ensemble des déterminants d'un ensemble simplicial réduit $X$ à valeurs dans un $2$-groupe $\G$. Observons qu'on a un foncteur: 
\begin{equation}\label{tututumared}
\xymatrix@R=5pt@C+15pt{
{\simp_0}^{op} \, \times \, \text{$2$-${\bf Grp}$} \ar[r]^-{{\bf det}}  &  {\bf Ens}\,,\\
\big(X,\G\big) \;\, \ar@{}[r]|-{\longmapsto} & {\bf det}_X(\G)}
\end{equation}
défini dans un morphisme d'ensembles simpliciaux réduits $f\colon\xymatrix@C-10pt{Y\ar[r]&X}$ et un morphisme de $2$-groupes $F=(F,m^{F})\colon\xymatrix@C-10pt{\G\ar[r]&\H}$ par la fonction:
$$
\xymatrix@C+34pt@R=1pt{
{\bf det}_X(\G) \ar[r]^-{{\bf det}_f(F)} & {\bf det}_Y(\H)\\
(D,T) \ar@{}[r]|-{\longmapsto} & (\overline{D},\overline{T})}
$$
où:
\begin{equation}\label{overlineDT}
\xymatrix@C+10pt@R=3pt{
Y_1\ar[r]^-{\overline{D}} & \big\{\text{Objets de $\H$}\big\}\\
B\ar@{}[r]|-{\mapsto} & F\big( D ( f_1\, B)\big)} 
\qquad\text{et}\qquad 
\xymatrix@C+10pt@R=3pt{Y_2\ar[r]^-{\overline{T}} & \big\{\text{Morphismes de $\H$}\big\}\\
\tau\ar@{}[r]|-{\mapsto} &  F \big( T (f_2\tau)\big)\circ m^F}
\end{equation}

Montrons que la couple $(\overline{D},\overline{T})$ est effectivement un déterminant de $Y$ à valeurs dans $\mathcal{H}$: En effet pour commencer remarquons que si $\tau$ est un $2$-simplexe de $Y$ alors $T(f_2\tau)$ est un morphisme de $\mathcal{G}$ de la forme:
$$
\xymatrix{
D(d_2\,f_2 \, \tau)\otimes D(d_0 \, f_2\, \tau) \ar[r]^-{T(f_2\,\tau)} &  D(d_1 \, f_2\,\tau)\,; }
$$
alors $\overline{T}(\tau)$ est un morphisme de $\mathcal{H}$ de la forme:
$$
\xymatrix{
\overline{D}(d_2 \, \tau)\otimes \overline{D}(d_0 \, \tau) \ar[r]^-{\overline{T}(\tau)} &  \overline{D}(d_1 \,\tau) \,,}
$$
parce que si $0\leq i\leq 2$ on a que:
$$
F\big(D(d_i \circ f_2\,\tau)\big) \, = \, F\big(D\circ f_1 (d_i\,\tau)\big) \, = \, \overline{D}(d_i\,\tau)\,.
$$

D'un autre côté on a que:
$$
\overline{D}(s_0\,\star) \, = \, F\big(D(f_1\circ s_0 \, \star)\big) \, = \, F\big(D(s_0\circ f_1 \, \star)\big) \, = \, F\big(D(s_0 \, \star)\big) \, = \, F(\mathbb{1}) \, = \, \mathbb{1}
$$
$$\text{et}$$
$$
\overline{T}(s_0s_0\, \star) \, = \, F\big(T(f_2s_0s_0\,\star)\big)\circ m^F \, = \,  F\big(T(s_0s_0\,\star)\big)\circ m^F \, = \,  F\big(\ell^{-1}_\mathbb{1}\big)\circ m^F_{\mathbb{1},\mathbb{1}} \, = \, \ell_{\mathbb{1}}^{-1}\,. 
$$

On vérifie sans peine que \eqref{tututumared} est un foncteur, par exemple si on se donne des morphismes des ensembles simpliciaux réduits:
$$
\xymatrix@C-10pt{Z\ar[r]^g&Y\ar[r]^{f} &X}
$$
et des morphismes de $2$-groupes:
$$
\xymatrix@C-10pt{\G\ar[r]^F&\H\ar[r]^-{G}&\mathcal{K}}
$$ 
pour montrer que ${\bf det}_{g}(G)\circ {\bf det}_{f}(F) = {\bf det}_{f\circ g}(G\circ F)$ il suffit de noter que pour tout déterminant $(D,T)$ de $X$ à valeurs dans $\mathcal{G}$ et tout $2$-simplexe $\tau$ de $Z$ on a un triangle commutatif de $\mathcal{K}$:
\begin{equation}
\def\objectstyle{\scriptstyle}
\def\labelstyle{\scriptstyle}
\xymatrix@C-20pt{
G\circ F\big(Df_1g_1(d_2\tau)\big) \otimes G\circ F\big(Df_1g_1(d_0\tau)\big) \ar[rd]_-{m^G} \ar[rr]^-{m^{G\circ F}} &   &  G\circ F \big(Df_1g_1(d_2\tau)\otimes Df_1g_1(d_0\tau)\big)\\
& G\big(F\big(Df_1g_1(d_2\tau)\big)\otimes F\big(Df_1g_1(d_0\tau)\big)\big) \ar[ru]_-{G(m^F)}&
}
\end{equation}


Montrons:

\begin{proposition}\label{detpeqNsF}
Les foncteurs ${\bf det}_{\bullet_1}\big(\,\bullet_2\,\big)$ et $\mathrm{Hom}_{\simp_0} \big( \, \bullet_1\, , \, \mathcal{N}(\bullet_2) \, \big)$ de source la catégorie produit ${\simp_0}^{op}\times\text{$2$-${\bf Grp}$}$ et but la catégorie des ensembles ${\bf Ens}$ sont naturellement isomorphes.
\end{proposition}
\begin{proof}
Si $\mathcal{G}$ est un $2$-groupe et $X$ est un ensemble simplicial réduit, remarquons que la fonction naturelle:
\begin{equation} \label{det1repress}
\xymatrix@C+5pt@R=1pt{
\mathrm{Hom}_{\simp_0} \big(X , \mathcal{N}(\mathcal{G})\big)\ar[r] & {\bf det}_{X}\big(\mathcal{G}\big)\\
f_\bullet\quad \ar@{}[r]|-{\longmapsto} & \; (f_1,f_2)}
\end{equation}
est bien définie: En effet si $f\colon\xymatrix@C-8pt{X\ar[r]&\mathcal{N}(\mathcal{G})}$ est un morphisme d'ensembles simpliciaux, on a par ailleurs que:
$$
f_1(s_0\star) = s_0(f_0\star) = s_0 (\star) = \mathbb{1} \qquad \; \text{et} \; \qquad f_2(s_0\circ s_0 \, \star) = s_0\circ s_0 (f_2\star) = s_0\circ s_0 (\star) = \ell_\mathrm{1}^{-1}\,.
$$ 

D'un autre si $\xi\in X_2$ alors $f_2(\xi)$ est un $2$-simplexe du nerf de $\mathcal{G}$ tel que $d_i\circ f_2(\xi) = f_1\circ d_i(\xi)$, autrement dit $f_2(\xi)$ est un morphisme de $\mathcal{G}$ de la forme:
$$
\xymatrix@C+5pt{ f_1(d_2\,\xi) \otimes f_1(d_0\,\xi) \ar[r]^-{f_2(\xi)} & f_1(d_1\,\xi)\,.}
$$

En fin si $\eta$ est un $3$-simplexe de $X$ alors $f_3(\eta)$ est un $3$-simplexe du nerf de $\mathcal{G}$ tel que $d_i\circ f_3(\eta) = f_2 \circ d_i(\eta)$, c'est-à-dire on a diagramme commutatif:
$$
\vcenter{\xymatrix@R-10pt@C+35pt{
f_1(A_{03})&f_1(A_{01})\otimes f_1(A_{13})\ar[l]_-{f_2(d_2\,\eta)}\\&\\
&f_1(A_{01}) \otimes (f_1(A_{12})\otimes f_1(A_{23})) \ar@{}[d]_-{a}|-{\text{\rotatebox[origin=c]{90}{$\cong$}}}\ar[uu]_-{f_1(A_{01})\otimes f_2(d_0\, \eta) }\\
f_1(A_{02})\otimes f_1(A_{23}) \ar[uuu]^-{f_2(d_1\,\eta)}&(f_1(A_{01})\otimes f_1(A_{12}))\otimes f_1(A_{23}) \ar[l]^-{f_2(d_3\,\eta)\otimes f_1(A_{23})}\,,}}
$$
\begin{align*}
\text{où} \qquad 
A_{03} \, = \, & d_{1} d_{1} \eta \, = \, d_{1} d_{2} \eta \qquad
A_{01} \, = \,  d_{2} d_{2} \eta \, = \, d_{2} d_{3} \eta \qquad
A_{13} \, = \,  d_{1} d_{0} \eta \, = \, d_{0} d_{2} \eta \\
A_{02} \, = \, & d_{2} d_{1} \eta \, = \, d_{1} d_{3}  \eta \qquad
A_{23} \, = \,  d_{0} d_{0} \eta \, = \, d_{0} d_{1} \eta \qquad
A_{12} \, = \,  d_{2} d_{0} \eta \, = \, d_{0} d_{3} \eta \,.
\end{align*}

Donc si $f\colon\xymatrix@C-8pt{X\ar[r]&\mathcal{N}(\mathcal{G})}$ est un morphisme d'ensembles simpliciaux, la couple $(f_1,f_2)$ est bien un déterminant de $X$ à valeurs dans $\G$.

Remarquons d'un autre que d'après le Lemme \ref{3cosqull} l'ensemble simplicial $\mathcal{N}(\mathcal{G})$ est $3$-cosque\-le\-ttique, donc pour montrer que la fonction \eqref{det1repress} est bijective il suffit de vérifier que la fonction: 
\begin{equation} \label{det1repress2}
\xymatrix@C+5pt@R=1pt{
\mathrm{Hom}_{\simp_{\leq 3}} \Big( \tau_{3}^{*} (X) , \tau_{3}^{*}\big(\mathcal{N}(\mathcal{G})\big)\Big)\ar[r] & {\bf det}_{X}\big(\mathcal{G}\big)\\
f_\bullet\quad \ar@{}[r]|-{\longmapsto} & \; (f_1,f_2)}
\end{equation}
est bijective où $\tau_{3}^*\colon\xymatrix@C-10pt{\simp\ar[r]&\simp_{\leq 3}}$ est le foncteur de troncation.

\emph{La fonction \eqref{det1repress2} est injective:}  

Supposons que $f,g\colon\xymatrix@C-10pt{\tau_{3}^*X\ar[r]&\tau_{3}^*\big(\mathcal{N}(\mathcal{G})\big)}$ sont deux morphismes des ensembles simpliciaux tronqués tels que $f_1=g_1$ et $f_2=g_2$. On a carrément que $f_0=g_0$. D'un autre vu que l'ensemble simplicial $\mathcal{N}(\mathcal{G})$ est faiblement $2$-cosquelettique, pour montrer que $f_3=g_3$ il suffit de noter que pour tout $\alpha\in X_3$ et $0\leq i\leq 3$ on a que:
$$
d_i(f_3\alpha) \, = \, f_2 \circ d_i(\alpha) \, = \, g_2 \circ d_i(\alpha) \, = \, d_i(g_3\alpha)\,.
$$

\emph{La fonction \eqref{det1repress2} est surjective:}

Si $(D,T)$ est un déterminant de $X$ à valeurs dans $\mathcal{G}$ on va définir un morphisme d'ensembles simpliciaux $f\colon\xymatrix@C-10pt{\tau_{3}^*X\ar[r]&\tau_{3}^*\big(\mathcal{N}(\mathcal{G})\big)}$ tel que $f_1=D$ et $f_2=T$. 

Pour commencer on définit $g\colon\xymatrix@C-10pt{\tau_{2}^*X\ar[r]&\tau_{2}^*\big(\mathcal{N}(\mathcal{G})\big)}$ par les règles $g_2=T$, $g_1=D$ et $g_0=\star$. $\tau_{2}^*f$ est effectivement un morphisme d'ensembles simpliciaux tronqués d'après les propriétés (i) et (ii) d'un déterminant et du Lemme qui suit:

\begin{lemme}
Soit $X$ un ensemble simplicial réduit et $\G$ un $2$-groupe. Si $(D,T)$ est un déterminant de $X$ à valeurs dans $\G$ pour tout élément $A$ de $X_{1}$ on a que $T\big(s_i(A)\big)=s_i\big(D(A)\big)$ si $0\leq i\leq 1$, c'est-à-dire:
$$
\xymatrix@C+48pt{\mathbb{1}\otimes D(A) \ar[r]^-{T(s_0(A)) \, = \, \ell^{-1}_{D(A)}} & D(A)}
\qquad\text{et}\qquad 
\xymatrix@C+48pt{D(A)\otimes\mathbb{1} \ar[r]^-{T(s_1(A)) \, = \, r^{-1}_{D(A)}} & D(A)}\,.
$$
\end{lemme}
\begin{proof}
Si $A\in X_{1}$ on va montrer que $T(s_0(A)) = \ell^{-1}_{D(A)}$. On vérifie l'égalité $T(s_1(A)) = r^{-1}_{D(A)}$ de façon analogue. 

Par ailleurs considérons $s_0 s_0(A)\in X_{3}$. D'après les propriétés (ii) et (iii) qui vérifie le déterminant $(D,T)$ on a un diagramme commutatif:
$$
\vcenter{\xymatrix@R-10pt@C+25pt{
D(A)&\mathbb{1}\otimes D(A)  \ar[l]_-{T(s_0\,A)}\\&\\
&\mathbb{1} \otimes (\mathbb{1}\otimes D(A)) \ar@{}[d]_-{a}|-{\text{\rotatebox[origin=c]{90}{$\cong$}}} \ar[uu]_-{\mathbb{1}\otimes T(s_0\, A) }\\
\mathbb{1}\otimes D(A) \ar[uuu]^-{T(s_0\,A)}&(\mathbb{1}\otimes \mathbb{1})\otimes D(A)  \ar[l]^-{\ell_\mathbb{1}^{-1}\otimes D(A)} \,;}}
$$
autrement dit on a un triangle commutatif:
$$
\vcenter{\xymatrix{
\big(\mathbb{1}\otimes \mathbb{1}\big) \otimes D(A) \ar[rr]^-{a_{\mathbb{1},\mathbb{1},D(A)}} \ar[rd]_-{\ell_\mathbb{1}^{-1}\otimes D(A)} &&    
\mathbb{1}\otimes \big(\mathbb{1}\otimes D(A)\big) \ar[ld]^-{\mathbb{1}\otimes T(s_0 \, A)}\\
&\mathbb{1}\otimes D(A)&}}
$$

Vu que $l_\mathbb{1}=r_\mathbb{1}$ et on a un triangle commutatif:
$$
\vcenter{\xymatrix{
\big(\mathbb{1}\otimes \mathbb{1}\big) \otimes D(A) \ar[rr]^-{a_{\mathbb{1},\mathbb{1},D(A)}} \ar[rd]_-{r_\mathbb{1}^{-1}\otimes D(A)} &&    
\mathbb{1}\otimes \big(\mathbb{1}\otimes D(A)\big) \ar[ld]^-{\mathbb{1}\otimes \ell^{-1}_{D(A)}}\\
&\mathbb{1}\otimes D(A)&}}
$$
il se suit que $\mathbb{1}\otimes T(s_0\, A) \, = \, \mathbb{1}\otimes \ell^{-1}_{D(A)}$ .

Donc $T\big(s_0 (A)\big) \, = \, \ell^{-1}_{D(A)}$ parce que $\mathbb{1} \otimes \, \cdot$ est une équivalence de catégories.
\end{proof}

Enfin vu que l'ensemble simplicial $\mathcal{N}\big(\mathcal{G}\big)$ est faiblement $2$-cosquelettique, il se suit de la propriété (iii) d'un déterminant à valeurs dans un $2$-groupe qu'il existe un seul morphisme d'ensembles simpliciaux tronqués $f\colon\xymatrix@C-10pt{\tau_{3}^*X\ar[r]&\tau_{3}^*\big(\mathcal{N}(\mathcal{G})\big)}$ tel que $\tau_{2}^*(f)=g$, en particulier $f_1=D$ et $f_2=T$. 
\end{proof}

Si $D=(D,T)$ et $D'=(D',T')$ sont deux déterminants d'un ensemble simplicial réduit $X$ à valeurs dans le $2$-groupe $\G$, un \emph{morphisme de déterminants} $H\colon\xymatrix@C-10pt{D\ar[r]&D'}$ est une fonction:
$$
\xymatrix@C+10pt{X_1\ar[r]^-{H} & \big\{\text{Morphismes de $\G$}\big\}}
$$
vérifiant les propriétés:
\begin{enumerate}
\item $H(A)$ est un morphisme en $\G$ de la forme:
$$
\xymatrix@C+10pt{D(A) \ar[r]^-{H(A)} & D'(A)} 
$$
pour tout $1$-simplexe $A$ de $X$.
\item On a que $H(s_0 \star) = \mathrm{id}_{\mathbb{1}}$.
\item On a un diagramme commutatif:
$$
\xymatrix@+10pt{
D(d_2\xi)\otimes D(d_0 \xi)\ar[d]_-{H(d_2\xi)\otimes H(d_0\xi)}\ar[r]^-{T(\xi)} & D(d_1\xi)  \ar[d]^-{H(d_1\xi)}  \\
D'(d_2\xi)\otimes D'(d_0 \xi)  \ar[r]_-{T'(\xi)} &  D'(d_1\xi)\,.}
$$
pour tout $2$-simplexe $\xi$ de $X$. 
\end{enumerate}

On vérifie sans peine qu'il y a un groupoïde canonique $\underline{{\bf det}}_X(\G)$ dont les objets sont les déterminants de $X$ à valeurs dans $\mathcal{G}$ et les morphismes son les morphismes des déterminants. 

Remarquons que le foncteur \eqref{tututumared} s'étend en un foncteur:
\begin{equation}\label{tututumaredhomgrp}
\vcenter{\xymatrix@R=5pt@C+15pt{
{\simp_0}^{op} \, \times \, \text{$2$-${\bf Grp}$} \ar[r]^-{\underline{\bf det}}  &  {\bf Grpd}\\
\big(X,\G\big) \;\, \ar@{}[r]|-{\longmapsto} & \underline{{\bf det}}_X(\G)} }\,,
\end{equation}
défini dans un morphisme d'ensembles simpliciaux réduits $f\colon\xymatrix@C-10pt{Y\ar[r]&X}$ et un morphisme de $2$-groupes $F=(F,m^{F})\colon\xymatrix@C-10pt{\G\ar[r]&\H}$ par le foncteur:
$$
\xymatrix@C+34pt@R=1pt{
\underline{\bf det}_X(\G) \ar[r]^-{\underline{\bf det}_f(F)} & \underline{\bf det}_Y(\H)\\
(D,T)  \ar[dddd]_-{H} & (\overline{D},\overline{T})  \ar[dddd]^-{\overline{H}} \\  & \\  \ar@{}[r]|-{\longmapsto} & \\    & \\  
(D',T') & (\overline{D}',\overline{T}')}
$$
où $\overline{H}$ est défini dans un $1$-simplexe $B$ de $Y$ comme le morphisme de $\H$:
$$
\vcenter{\xymatrix@C+25pt@R=2pt{ 
\overline{D}(B) \ar@{=}[d] \ar[r]^-{\overline{H}(B)}     & \overline{D}'(B) \ar@{=}[d]\\
F\big(D(f_1\,B)\big)  \ar[r]_-{F\big(H(f_1\,B)\big)}  & F\big(D'(f_1\,B)\big)
}}\,.
$$

Dans l'énoncé qui suit on note $\pi_0\colon\xymatrix@C-10pt{{\bf Grpd}\ar[r]&{\bf Ens}}$ le foncteur ensemble des objets à isomorphismes près. 

\begin{proposition} \label{detpeqNsF2}
Les foncteurs $\pi_0\big(\underline{\bf det}_{\bullet_1}(\,\bullet_2\,)\big)$ et $\pi_0\Big(\underline{\mathrm{Hom}}_{\simp_0} \big( \, \bullet_1\, , \, \mathcal{N}(\bullet_2) \, \big)\Big)$ de source la catégorie produit ${\simp_0}^{op}\times\text{$2$-${\bf Grp}$}$ et but la catégorie des ensembles ${\bf Ens}$ sont naturellement isomorphes.
\end{proposition}

\begin{proof}
D'après la Proposition \ref{detpeqNsF} il suffit de montrer que si $(D,T)$ et $(D',T')$ sont deux déterminants d'un ensemble simplicial réduit $X$ à valeurs dans un $2$-groupe $\G$, alors il existe un morphisme de déterminants $H\colon\xymatrix@C-10pt{(D,T)\ar[r]&(D',T')}$ si et seulement s'il existe un morphisme d'ensembles simpliciaux $F\colon\xymatrix@C-10pt{X\times\Delta^1\ar[r]&\mathcal{N}(\G)}$ tel que le morphisme composé:
$$
\xymatrix@C+10pt{\star\times\Delta^1 \ar[r]^-{s_0\times\Delta^1} & X\times\Delta^1\ar[r]^-{F}&\mathcal{N}(\G)}
$$
est le morphisme constant à valeurs $\mathbb{1}$ et le diagramme:
$$
\xymatrix@R=2pt@C+10pt{
X \ar[rd]|-{\nu_1} \ar@/^12pt/[rrd]^-{F_{(D,T)}}& &\\
& X\times \Delta^1 \ar[r]|-{F} & \mathrm{N}(\G) \\
X \ar[ru]|-{\nu_0} \ar@/_12pt/[rru]_-{F_{(D',T')}} & &
}
$$
est un diagramme commutatif, où $\nu_i$ est le morphisme composé:
$$
\xymatrix@C+12pt{X_{\bullet,1}\,\cong \, X_{\bullet,1}\times \Delta^0 \ar[r]^-{X_{\bullet,1}\times \delta_i} & X_{\bullet,1}\times \Delta^1}\;.
$$

Vu que $\mathcal{N}(\G)$ est un ensemble simplicial faiblement $2$-cosquelettique on vérifie sans difficulté que se donner un morphisme $F\colon\xymatrix@C-10pt{X\times\Delta^1\ar[r]&\mathcal{N}(\G)}$ vérifiant les propriétés désirées équivaut à se donner deux fonctions:
\begin{equation}\label{lesfoncctt}
\vcenter{\xymatrix@C+8pt{X_1\ar[r]^-{\sigma} & \Big\{ \text{Objets  de  $\G$}\Big\} }}
\qquad\text{et}\qquad
\vcenter{ \xymatrix@C+15pt{X_2\ar[r]^-{(S,S')} & \Big\{ \text{Morphismes  de  $\G$}\Big\}^2}}
 \end{equation}
telles que:
\begin{enumerate}
\item[(a)] Pour tout $\xi\in X_2$:
$$
\xymatrix@C+10pt{
D(d_2\xi) \otimes \sigma(d_0\xi) \ar[r]^-{S_\xi} & \sigma(d_1\xi) & \ar[l]_-{S'_\xi} \sigma(d_2\xi)\otimes D'(d_0\xi)
}
$$
\item[(b)]  $\sigma(s_0\star)=\mathbb{1}$  et $S_{s_0s_0\star} = r_{\mathbb{1}}^{-1} = l_{\mathbb{1}}^{-1} = S'_{s_0s_0\star}$
\item[(c)] Pour tout $\eta\in X_3$ on a le diagramme commutatif de $\G$:
$$
{\footnotesize{\xymatrix@R=5pt@C+42pt{
\big(\sigma(A_{01})\otimes D'(A_{12})\big)\otimes D'(A_{23}) \ar@{}[rdddd]|-{{\bf (I)}} \ar[r]^-{S'_{d_3\eta}\otimes D'(A_{23})}  &
\sigma(A_{02}) \otimes D'(A_{23}) \ar[dddd]|-{S'_{d_1\eta}}  \ar@{}[rdddd]|-{{\bf (II)}}  & 
\big(D(A_{01})\otimes \sigma(A_{12})\big)\otimes D'(A_{23})\ar[l]_-{S_{d_3\eta} \otimes D'(A_{23})}\\
\sigma(A_{01})\otimes\big(D'(A_{12})\otimes D'(A_{23})\big)  \ar[ddd]_-{\sigma(A_{01})\otimes T'_{d_0\eta} }
\ar@{}[u]^-{a}|-{\text{\rotatebox[origin=c]{90}{$\cong$}}}& & 
D(A_{01})\otimes \big(\sigma(A_{12})\otimes D'(A_{23}) \big) \ar[ddd]^-{D(A_{01})\otimes S'_{d_0\eta}} 
  \ar@{}[u]_-{a}|-{\text{\rotatebox[origin=c]{90}{$\cong$}}} \\ &&\\&&\\
 \sigma(A_{01})\otimes D'(A_{13})\ar[r]_-{S'_{d_2\eta}}
&\sigma(A_{03}) \ar@{}[rdddd]|-{{\bf (III)}} & D(A_{01})\otimes \sigma(A_{13}) \ar[l]|-{S_{d_2\eta}} \\&&\\&&\\
& &  D(A_{01})\otimes\big(D(A_{12})\otimes \sigma(A_{23})\big)  \ar[uuu]_-{D(A_{01})\otimes S_{d_0\eta}}
  \ar@{}[d]^-{a}|-{\text{\rotatebox[origin=c]{90}{$\cong$}}} \\ 
&D(A_{02}) \otimes \sigma(A_{23})\ar[uuuu]^-{S_{d_1\eta}}  & \big(D(A_{02})\otimes D(A_{12})\big)\otimes \sigma(A_{23})\ar[l]^-{T_{d_3\eta}\otimes \sigma(A_{23})}
}}}
$$
\begin{align*}
\text{où} \qquad 
A_{03} \, = \, & d_{1} d_{1} \eta \, = \, d_{1} d_{2} \eta \qquad
A_{01} \, = \,  d_{2} d_{2} \eta \, = \, d_{2} d_{3} \eta \qquad
A_{13} \, = \,  d_{1} d_{0} \eta \, = \, d_{0} d_{2} \eta \\
A_{02} \, = \, & d_{2} d_{1} \eta \, = \, d_{1} d_{3}  \eta \qquad
A_{23} \, = \,  d_{0} d_{0} \eta \, = \, d_{0} d_{1} \eta \qquad
A_{12} \, = \,  d_{2} d_{0} \eta \, = \, d_{0} d_{3} \eta \,.
\end{align*}
\end{enumerate}

Étant donné de fonctions \eqref{lesfoncctt} vérifiant les propriétés (a)-(c) on définit un morphisme de déterminants $H\colon\xymatrix@C-10pt{(D,T)\ar[r]&(D',T')}$ par la fonction:
$$
\xymatrix@C+10pt@R=3pt{
X_1\ar[r]^-{H} & \big\{\text{Morphismes de $\G$}\;,\big\}\\
A\ar@{}[r]|-{\longmapsto} &    l^{-1}_{DA}\circ (S'_{s_0A})^{-1}\circ S_{s_1A}\circ r_{DA}
}
$$
c'est-à-dire si $A$ est un $1$-simplexe de $X$ le morphisme $H(A)$ est le composé:
$$
\xymatrix@C+10pt{
D(A) \ar[r]^-{r_{D(A)}} & D(A) \otimes \sigma(s_0\star) \ar[r]^-{S_{s_1A}} & \sigma(A) \ar[r]^-{(S'_{s_0A})^{-1}}&  \sigma(s_0\star)\otimes D'(A) \ar[r]^-{l^{-1}_{DA}} & D'(A) \;.
}
$$

Il se suit en particulier que $H(s_0\star) = \mathrm{id}_\mathbb{1}$. En plus si $\xi$ est un $2$-simplexe de $X$ on montre que:
$$
\xymatrix@+10pt{
D(d_2\xi)\otimes D(d_0 \xi)\ar[d]_-{H(d_2\xi)\otimes H(d_0\xi)}\ar[r]^-{T(\xi)} & D(d_1\xi)  \ar[d]^-{H(d_1\xi)}  \\
D'(d_2\xi)\otimes D'(d_0 \xi)  \ar[r]_-{T'(\xi)} &  D'(d_1\xi)\,.}
$$
est un carré commutatif en considérant le diagramme commutatif ${\bf (I)}$ pour $\eta=s_0\xi$, le diagramme commutatif ${\bf (II)}$ pour $\eta=s_1\xi$ et le diagramme commutatif ${\bf (III)}$ pour $\eta=s_2\xi$.

Enfin si d'un autre $H\colon\xymatrix@C-10pt{(D,T)\ar[r]&(D',T')}$ est un morphisme de déterminants, les fonctions \eqref{lesfoncctt} vérifiant les propriétés (a)-(c) sont définies par les règles:
$$
\vcenter{\xymatrix@C+18pt@R=3pt{
X_1\ar[r]^-{\sigma \, = \, D} & \Big\{ \text{Objets  de  $\G$}\Big\}\;, \\
A\ar@{}[r]|-{\longmapsto} & D(A)
}}
\qquad \;
\vcenter{\xymatrix@C+18pt@R=3pt{
X_2\ar[r]^-{S \, = \, T} & \Big\{ \text{Morphismes  de  $\G$}\Big\}\;, \\
\xi\ar@{}[r]|-{\longmapsto} & T_\xi
}}
$$
$$
\text{et}\qquad\;
\vcenter{\xymatrix@C+13pt@R=3pt{
X_2\ar[r]^-{S'} & \Big\{ \text{Morphismes  de  $\G$}\Big\}\;. \\
\xi\ar@{}[r]|-{\longmapsto} &  T_\xi \circ \big(D(d_2\xi)\otimes H(d_0\xi)^{-1}\big)
}}
$$
\end{proof}

On déduit des Propositions \ref{modred}, \ref{ilest} et \ref{detpeqNsF2} et du Corollaire \ref{grho}:  

\begin{corollaire}\label{repsteouni2}
Si $X$ est un ensemble simplicial réduit et $(\varphi,m^{\varphi})\colon \xymatrix@C-10pt{\G\ar[r]& \H}$ est une $2$-équivalence faible de $2$-groupes la fonction
$$
\xymatrix@C+45pt{
 \pi_0\big(\underline{{\bf det}}_X(\G)\big) \ar[r]^-{\pi_0\big(\underline{{\bf det}}_X(\varphi,m^{\varphi})\big)} &  \pi_0\big(\underline{{\bf det}}_X(\H)\big)
}
$$
est bijective.
\end{corollaire}

\renewcommand{\thesubsection}{\S\thesection.\arabic{subsection}}
\subsection{}\;\label{lenerfSegalss}
\renewcommand{\thesubsection}{\thesection.\arabic{subsection}}

On définit le foncteur \emph{nerf de Segal} pour les $2$-groupes:
\begin{equation}
\xymatrix{ \text{$2$-${\bf Grp}$}\ar[r]^-{\mathcal{N}_{\mathcal{S}}} & {\bf MS}\,,}
\end{equation}
par le composé:
\begin{equation}\label{nP22}
\xymatrix@C+20pt@R=5pt{
\text{$2$-${\bf Grp}$} \ar[r]& \text{$2$-${\bf Grp}$}^{\Delta^{op}} \ar[r] & \simp_0^{\Delta^{op}}  \, \cong \, {\bf MS}\,,\\
 \G \ar@{}[r]|-{\longmapsto}& \G^{[\bullet]} \ar@{}[r]|-{\longmapsto}&  \mathcal{N}\big(\G^{[\bullet_1]}\big)_{\bullet_2}}
\end{equation}
où $\G\longmapsto\G^{[\bullet]}$ est déduit du $2$-foncteur \eqref{homodeG} et $\mathcal{N}$ est le foncteur nerf des $2$-groupes \eqref{nervio22a}.

\begin{proposition}\label{fibrachii}
Si $\G$ est un $2$-groupe le pré-monoïde de Segal $\mathcal{N}_{\mathcal{S}}(\G)$ est un objet fibrant de la catégorie de modèles $({\bf MS},{\bf W}^{diag}_{2}, {\bf mono},{\bf fib}_{2}^{diag})$ de la Proposition \ref{moduno}.
\end{proposition}
\begin{proof}
Si $\G$ est un $2$-groupe, on va montrer que l'ensemble bisimplicial $\mathcal{N}_{\mathcal{S}}(\G)$ vérifie les propriétés (i)-(iv) du Lemme \ref{fibrafifi}. 

Pour commencer remarquons que d'après le Lemme \ref{porporfin} et le Corollaire \ref{grho} l'ensemble simplicial $\mathcal{N}_{\mathcal{S}}(\G)$ vérifie la propriété (i) du Lemme \ref{fibrafifi}. Du même la propriété (ii) est déduit de la Proposition \ref{ilest}.

Dans le reste de la preuve on va appeler \emph{un prisme de $\G$} la donnée de neuf objets de $\G$ et douze morphismes arrangés dans la forme:
\begin{equation}\label{prismes}
\xymatrix@C-10pt{  
&A_2\otimes A_0 \ar'[d][dd]|(.35){\tau_2\otimes\tau_0} \ar[rrr]^f \ar[dl]_{\varphi_2\otimes \varphi_0} &&&A_1 \ar[dl]|{\varphi_1} \ar[dd]^-{\tau_1}\\
B_2\otimes B_0 \ar@<-2pt>[rrr]|(.65)g \ar[rd]_{\psi_2\otimes \psi_0} &&& B_1 \ar[rd]|-{\psi_1}&\\
&C_2\otimes C_0\ar[rrr]_-h&&&C_1\,,}
\end{equation}
sans demander que les faces du prisme soient des diagrammes commutatifs.

Pour montrer la propriété (iii) du Lemme \ref{fibrafifi} il faut montrer que si $0\leq k\leq 2$ la fonction:
\begin{equation}\label{iiipp}
\xymatrix{
\mathrm{Hom}_{\ssimp}\Big(\partial\Delta^2\boxtimes\Delta^2, \mathcal{N}_{\mathcal{S}}(\G)\Big)\ar[rr]&&
\mathrm{Hom}_{\ssimp}\Big(\partial\Delta^2\boxtimes\Lambda^{2,k}, \mathcal{N}_{\mathcal{S}}(\G)\Big)}
\end{equation}
induite par le morphisme d'inclusion canonique $\xymatrix@C-5pt{\Lambda^{2,k}\ar@{^(->}[r]&\Delta^2}$ est surjective.

En premier lieu remarquons que l'ensemble $\mathrm{Hom}_{\ssimp}\big(\partial\Delta^2\boxtimes\Delta^2, \mathcal{N}_{\mathcal{S}}(\G)\big)$ s'identifie à l'ensemble des prismes \eqref{prismes} dont les faces en forme de carrés sont de diagrammes commutatifs, mais pas nécessairement les faces en forme de triangles.

D'un autre l'ensemble $\mathrm{Hom}_{\ssimp}\big(\partial\Delta^2\boxtimes\Lambda^{2,k}, \mathcal{N}_{\mathcal{S}}(\G)\big)$ s'identifie à l'ensemble des couples de triangles pas nécessairement commutatifs dans $\G$:
$$
\left(\vcenter{
\xymatrix@R-5pt@C-10pt{
&X\ar[dd]\ar[dl]\\
Y\ar@{}[r]|(.6){?}\ar[rd]&\\
&Z}}\quad,
\qquad
\vcenter{\xymatrix@R-5pt@C-10pt{
&X'\ar[dd]\ar[dl]\\
Y'\ar@{}[r]|(.6){?}\ar[rd]&\\
&Z'}
}\right)\,.
$$

L'image par la fonction \eqref{iiipp} d'un prisme \eqref{prismes} dans son ensemble de départe est égale à la couple des triangles pas nécessairement commutatifs:
$$
\left(\vcenter{\xymatrix@R-5pt@C-10pt{
&A_i\ar[dd]\ar[dl]_-{\varphi_i}\ar[dd]^-{\tau_i}\\
B_i\ar[rd]_-{\psi_i}\ar@{}[r]|(.6){?}&\\
&C_i}}\quad,
\qquad
\vcenter{\xymatrix@R-5pt@C-10pt{
&A_j\ar[dd]\ar[dl]_-{\varphi_i}\ar[dd]^-{\tau_j}\\
B_j\ar[rd]_-{\psi_i}\ar@{}[r]|(.6){?}&\\
&C_j}}\right)
$$
où $0\leq i < j \leq 2$ avec $i,j\neq k$.

Pour montrer que \eqref{iiipp} est une fonction surjective rappelons que d'après le Lemme \ref{inverses} ils existent de foncteurs $\xymatrix@C-4pt{\mathcal{G} \ar[r]^-{\iota^d}& \mathcal{G}}$ et $\xymatrix@C-4pt{\mathcal{G} \ar[r]^-{\iota^g}& \mathcal{G}}$ et d'isomorphismes naturels:
$$\xymatrix@C+5pt{X\otimes \iota^d(X) \ar[r]^-{\alpha_X}_-{\cong} & \mathbb{1}&\ar[l]_-{\beta_X}^-{\cong} \iota^g(X)\otimes X\,.}$$

Supposons alors que $k=2$ et considérons la couple de triangles pas nécessairement commutatifs de $\G$:
$$
\left(\vcenter{\xymatrix@R-5pt@C-10pt{
&A_0\ar[dd]\ar[dl]_-{\varphi_0}\ar[dd]^-{\tau_0}\\
B_0\ar[rd]_-{\psi_0}\ar@{}[r]|(.6){?}&\\
&C_0}}\quad,
\qquad\;
\vcenter{\xymatrix@R-5pt@C-10pt{
&A_1\ar[dd]\ar[dl]_-{\varphi_1}\ar[dd]^-{\tau_1}\\
B_1\ar[rd]_-{\psi_1}\ar@{}[r]|(.6){?}&\\
&C_1}}\right)\,.
$$

On définit le triangle manquant du prisme \eqref{prismes} par:
$$
\left(\vcenter{\xymatrix@R-5pt@C-10pt{
&A_2\ar[dd]\ar[dl]_-{\varphi_2}\ar[dd]^-{\tau_2}\\
B_2\ar[rd]_-{\psi_2}\ar@{}[r]|(.6){?}&\\
&C_2}}\right)
\qquad = \qquad
\left(\vcenter{\xymatrix@R-5pt@C-10pt{
&A_1\otimes \iota^g(A_0)\ar[dd]^-{\tau_1\otimes \iota^g(\tau_0)}\ar[dl]_-{\varphi_1\otimes \iota^g(\varphi_0)}\\
B_1\otimes \iota^g(B_0)\ar[rd]_-{\psi_1\otimes \iota^g(\psi_0)}\ar@{}[r]|(.6){?}&\\
&C_1\otimes \iota^g(C_0)}}\right)
$$
et les morphismes $f$, $g$ et $h$ comme étant les morphismes composés de $\G$ suivants:
$$
f\colon\bigg(\xymatrix@C+5pt{
\big(A_1\otimes \iota^g(A_0) \big)\otimes A_0 \ar[r]^{a} & A_1\otimes \big( \iota^g(A_0) \otimes A_0 \big) \ar[r]^-{A_1\otimes \beta_{A_0}}& A_1\otimes \mathbb{1} \ar[r]^-{r_{A_1}^{-1}}& A_1
}\bigg)\,,
$$
$$
g\colon\bigg(\xymatrix@C+5pt{
\big(B_1\otimes \iota^g(B_0) \big)\otimes B_0 \ar[r]^{a} & B_1\otimes \big( \iota^g(B_0) \otimes B_0 \big) \ar[r]^-{B_1\otimes \beta_{B_0}}& B_1\otimes \mathbb{1} \ar[r]^-{r_{B_1}^{-1}}& B_1
}\bigg)
$$
$$\text{et}$$
$$
h\colon\bigg(\xymatrix@C+5pt{
\big(C_1\otimes \iota^g(C_0) \big)\otimes C_0 \ar[r]^{a} & C_1\otimes \big( \iota^g(C_0) \otimes C_0 \big) \ar[r]^-{C_1\otimes \beta_{C_0}}& C_1\otimes \mathbb{1} \ar[r]^-{r_{C_1}^{-1}}& C_1
}\bigg)\,.
$$

On obtient alors des carrés commutatifs:
$$
\vcenter{\xymatrix@R-5pt@C+5pt{
\big(A_1\otimes \iota^g(A_0) \big)\otimes A_0 \ar[d]_-{\big(\varphi_1\otimes \iota^g(\varphi_0)\big)\otimes \varphi_0}   \ar[r]^-f &  A_1\ar[d]^{\varphi_1} \\
\big(B_1\otimes \iota^g(B_0) \big)\otimes B_0 \ar[r]_-g & B_1 \,,
}}\qquad\qquad
\vcenter{\xymatrix@R-5pt@C+5pt{
 \big(B_1\otimes \iota^g(B_0) \big)\otimes B_0 \ar[d]_-{\big(\psi_1\otimes \iota^g(\psi_0)\big)\otimes \psi_0}  \ar[r]^-g & B_1\ar[d]^-{\psi_1}\\
\big(C_1\otimes \iota^g(C_0) \big)\otimes C_0 \ar[r]_-h & C_1 
}}
$$
$$
\text{et}\qquad\;
\vcenter{\xymatrix@R-5pt@C+5pt{
\big(A_1\otimes \iota^g(A_0) \big)\otimes A_0 \ar[d]_-{\big(\tau_1\otimes \iota^g(\tau_0)\big)\otimes \tau_0} \ar[r]^-f &  A_1\ar[d]^-{\tau_1}  \\
\big(C_1\otimes \iota^g(C_0) \big)\otimes C_0 \ar[r]_-h & C_1\,;
}}
$$
parce que si $\varphi\colon\xymatrix@C-5pt{X\ar[r]&Y}$ et $\varphi'\colon\xymatrix@C-5pt{X'\ar[r]&Y'}$ sont de morphismes quelconques de $\G$, on a des carrés commutatifs:
$$
\xymatrix{
X\ar[r]^-{\varphi} \ar[d]_-{r_X} & Y \ar[d]^-{r_Y} \\
X\otimes \mathbb{1} \ar[r]_-{\varphi\otimes \mathbb{1}} & Y\otimes \mathbb{1}}
\qquad 
\xymatrix@C+20pt{
(X\otimes \iota^gX')\otimes X' \ar[d]_-{a}\ar[r]^-{(\varphi\otimes\iota^g\varphi')\otimes\varphi'} & (Y\otimes \iota^gY')\otimes Y'\ar[d]^-{a}\\
X\otimes (\iota^gX'\otimes X')\ar[r]_-{\varphi\otimes(\iota^g\varphi'\otimes\varphi')} & Y\otimes (\iota^gY'\otimes Y')}
$$
$$
\text{et}\qquad\;
\xymatrix@C+10pt{
X\otimes \big( \iota^gX'\otimes X' \big) \ar[d]_-{\varphi \otimes (\iota^g\varphi'\otimes\varphi')}\ar[r]^-{X\otimes\beta_{X'}} & X\otimes\mathbb{1} \ar[d]^-{\varphi\otimes\mathbb{1}}\\
Y\otimes \big( \iota^gY'\otimes Y' \big) \ar[r]_-{Y\otimes\beta_{Y'}} & Y\otimes\mathbb{1} \,.}
$$

Si $k=0$ on montre de façon analogue que la fonction \eqref{iiipp} est surjective. Supposons finalement que $k=1$ et considérons la couple de triangles pas nécessairement commutatifs de $\G$:
$$
\left(\vcenter{\xymatrix@R-5pt@C-10pt{
&A_0\ar[dd]\ar[dl]_-{\varphi_0}\ar[dd]^-{\tau_0}\\
B_0\ar[rd]_-{\psi_0}\ar@{}[r]|(.6){?}&\\
&C_0}}\quad,
\qquad\;
\vcenter{\xymatrix@R-5pt@C-10pt{
&A_2\ar[dd]\ar[dl]_-{\varphi_2}\ar[dd]^-{\tau_2}\\
B_2\ar[rd]_-{\psi_2}\ar@{}[r]|(.6){?}&\\
&C_2}}\right)\,.
$$

On définit le triangle manquant par:
$$
\left(\vcenter{\xymatrix@R-5pt@C-10pt{
&A_1\ar[dd]\ar[dl]_-{\varphi_1}\ar[dd]^-{\tau_1}\\
B_1\ar[rd]_-{\psi_1}\ar@{}[r]|(.6){?}&\\
&C_1}}\right)
\qquad = \qquad
\left(\vcenter{\xymatrix@R-5pt@C-10pt{
&A_0\otimes A_2 \ar[dd]^-{\tau_0\otimes\tau_2}\ar[dl]_-{\varphi_0\otimes\varphi_2}\\
B_1\otimes B_2\ar[rd]_-{\psi_0\otimes\psi_2}\ar@{}[r]|(.6){?}&\\
&C_1\otimes C_2}}\right)
$$
et les morphismes $f$, $g$ et $h$ comme étant des identités. Il se suit immédiatement que dans le prisme résultant \eqref{prismes} les faces en forme de carrés sont de diagrammes commutatifs.

Enfin, montrons la propriété (iv) du Lemme \ref{fibrafifi}. Commençons pour vérifier que si $0\leq k\leq 2$ la fonction:
\begin{equation}\label{fonctionfinalll1}
\def\objectstyle{\scriptstyle}
\def\labelstyle{\scriptstyle}
\xymatrix{
\mathrm{Hom}_{\ssimp}\big(\Delta^2\boxtimes\Delta^2, \mathcal{N}_{\mathcal{S}}(\G)\big)\ar[d]\\
\mathrm{Hom}_{\ssimp}\big(\partial\Delta^2\boxtimes\Delta^{2}, \mathcal{N}_{\mathcal{S}}(\G)\big)
\underset{\mathrm{Hom}_{\ssimp}(\partial\Delta^2\boxtimes\Lambda^{2,k}, \mathcal{N}_{\mathcal{S}}(\G))}{\times}
\mathrm{Hom}_{\ssimp}\big(\Delta^2\boxtimes\Lambda^{2,k}, \mathcal{N}_{\mathcal{S}}(\G)\big)}
\end{equation}
induite du carré commutatif:
$$
\xymatrix@-5pt{
\Delta^2\times\Delta^2 & \ar@{_(->}[l]\; \Delta^2\times\Lambda^{2,k}\\
\ar@{^(->}[u] \partial\Delta^2\times\Delta^2&\ar@{_(->}[l]\ar@{^(->}[u] \; \partial\Delta^2\times\Lambda^{2,k}}
$$
est surjective.

Notons en premier lieu que l'ensemble de départe de la fonction \eqref{fonctionfinalll1} est isomorphe à l'ensemble des prismes \eqref{prismes} vérifiant les deux propriétés: Les faces en forme des carrés sont des diagrammes commutatifs et les triangles:
$$
\vcenter{\xymatrix@R-5pt@C-10pt{
&A_0\ar[dd]\ar[dl]_-{\varphi_0}\ar[dd]^-{\tau_0}\\
B_0\ar[rd]_-{\psi_0}\ar@{}[r]|(.6){}&\\
&C_0}}\qquad\;
\vcenter{\xymatrix@R-5pt@C-10pt{
&A_1\ar[dd]\ar[dl]_-{\varphi_1}\ar[dd]^-{\tau_1}\\
B_1\ar[rd]_-{\psi_1}\ar@{}[r]|(.6){}&\\
&C_1}}\qquad\text{et}
\qquad
\vcenter{\xymatrix@R-5pt@C-10pt{
&A_2\ar[dd]\ar[dl]_-{\varphi_2}\ar[dd]^-{\tau_2}\\
B_2\ar[rd]_-{\psi_2}\ar@{}[r]|(.6){}&\\
&C_2}}\,,
$$
sont commutatifs. En particulier les deux faces en forme de triangles de \eqref{prismes} sont de diagrammes commutatifs.

D'un autre l'ensemble d'arrivée de \eqref{fonctionfinalll1} s'identifie à l'ensemble des prismes \eqref{prismes} tels que ses faces en forme des carrés sont des diagrammes commutatifs et seulement les deux triangles:
$$
\vcenter{\xymatrix@R-5pt@C-10pt{
&A_i\ar[dd]\ar[dl]_-{\varphi_i}\ar[dd]^-{\tau_i}\\
B_i\ar[rd]_-{\psi_i}&\\
&C_i}}\qquad\text{et}
\qquad
\vcenter{\xymatrix@R-5pt@C-10pt{
&A_j\ar[dd]\ar[dl]_-{\varphi_i}\ar[dd]^-{\tau_j}\\
B_j\ar[rd]_-{\psi_i}&\\
&C_j}}\,,
$$
où $0\leq i<j\leq 2$ et $i,j\neq k$, sont commutatifs. 

La fonction \eqref{fonctionfinalll1} est alors la fonction qu'oublie que le triangle suivant est commutatif:
$$
\vcenter{\xymatrix@R-5pt@C-10pt{
&A_k\ar[dd]\ar[dl]_-{\varphi_k}\ar[dd]^-{\tau_k}\\
B_k\ar[rd]_-{\psi_k}\ar@{}[r]|(.6){?}&\\
&C_k}}\,.
$$

On montre qu'en fait la fonction \eqref{fonctionfinalll1} est bijective. En effet si $k=1$ c'est une conséquence du fait que $\G$ est un groupoïde. On déduit les cas $k=0,2$ de l'énoncé suivant:
\begin{lemme}\label{fifidelite}
Si $f\colon\xymatrix@C-5pt{X\ar[r]&Y}$ et $\xi\colon\xymatrix@C-5pt{X\otimes X'\ar[r]&Y\otimes Y'}$ sont des morphismes d'un $2$-groupe $\G$, il existe un seul morphisme $f'\colon\xymatrix@C-5pt{X'\ar[r]&Y'}$ tel que $f\otimes f' = \xi$.
\end{lemme}
\begin{proof}
Vu que tous les objets de $\G$ sont inversibles le foncteur $X\otimes \,-\, \colon \xymatrix@C-5pt{\G\ar[r] & \G}$ est une équivalence de catégories, donc il existe un seul morphisme $f'\colon\xymatrix@C-10pt{X'\ar[r]&Y'}$ tel que $X\otimes f' \, = \, (f\otimes Y')^{-1}\circ \xi $:
$$
\xymatrix{
X\otimes X' \ar[r]^-{X\otimes f'} \ar[rd]_{\xi}& X\otimes Y' \ar[d]^-{f\otimes Y'}\\
& Y\otimes Y'}
$$

Autrement dit il existe un seul morphisme $f'\colon\xymatrix@C-10pt{X'\ar[r]&Y'}$ tel que: 
$$f\otimes f' \, = \, (f \otimes Y') \circ (X \otimes f') \, = \, \xi\,.$$
\end{proof}

Montrons pour conclure que la fonction:
\begin{equation}\label{fonctionfinalll}
\def\objectstyle{\scriptstyle}
\def\labelstyle{\scriptstyle}
\xymatrix{
\mathrm{Hom}_{\ssimp}\big(\Delta^1\boxtimes\Delta^2, \mathcal{N}_{\mathcal{S}}(\G)\big)\ar[d]\\
\mathrm{Hom}_{\ssimp}\big(\partial\Delta^1\boxtimes\Delta^{2}, \mathcal{N}_{\mathcal{S}}(\G)\big)
\underset{\mathrm{Hom}_{\ssimp}(\partial\Delta^1\boxtimes\Lambda^{2,k}, \mathcal{N}_{\mathcal{S}}(\G))}{\times}
\mathrm{Hom}_{\ssimp}\big(\Delta^1\boxtimes\Lambda^{2,k}, \mathcal{N}_{\mathcal{S}}(\G)\big)}
\end{equation}
induite du carré commutatif:
$$
\xymatrix@-5pt{
\Delta^1\times\Delta^2 & \ar@{_(->}[l]\; \Delta^1\times\Lambda^{2,k}\\
\ar@{^(->}[u] \partial\Delta^1\times\Delta^2&\ar@{_(->}[l]\ar@{^(->}[u] \; \partial\Delta^1\times\Lambda^{2,k}}
$$
est bijective, en particulier \eqref{fonctionfinalll} est surjective.

Pour cela remarquons que l'ensemble $\mathrm{Hom}_{\ssimp}\big(\Delta^1\boxtimes\Delta^2, \mathcal{N}_{\mathcal{S}}(\G)\big)$ est isomorphe à l'ensemble des carrés commutatifs de la forme:
\begin{equation}\label{fibrantcarr}
\xymatrix@C+8pt@R-2pt{
X_2\otimes X_0 \ar[d]_-{\varphi_2\otimes\varphi_0}   \ar[r]^-{f} &  X_1 \ar[d]^-{\varphi_1} \\
Y_2\otimes Y_0 \ar[r]_-g & Y_1\,;}
\end{equation}
tandis que l'ensemble d'arrivée de \eqref{fonctionfinalll} s'identifie à l'ensemble dont les objets sont déterminés par six objets et quatre morphismes comme ci-dessous:
$$
\left(
\vcenter{\xymatrix@C+5pt@R-2pt{
X_2\otimes X_0 \ar[r]^-{f} & X_1\\
Y_2\otimes Y_0 \ar[r]_-g  & Y_1}}\,;\quad
\vcenter{\xymatrix@C-5pt@R-2pt{
X_i \ar[d]_-{\varphi_i} \\
Y_i}}\,,
\vcenter{\xymatrix@C-5pt@R-2pt{
 X_j\ar[d]^-{\varphi_j} \\
Y_j}}\right)
$$
où $0\leq i<j\leq 2$ et $i,j\neq k$.

La fonction \eqref{fonctionfinalll} est alors la fonctions qu'oublie le morphisme $\varphi_k$; en particulier on montre qu'elle est bijective pour $0\leq k\leq 2$. En effet, si $k=1$ alors $\varphi_1= g \circ (\varphi_2 \otimes \varphi_0) \circ f^{-1}$; et si $k=0$ ou $k=2$ alors $\varphi_k$ est le seul morphisme de $\G$ tel que $\varphi_2\otimes \varphi_1 = g^{-1} \circ \varphi_1 \circ f$ (voir le Lemme \ref{fifidelite}).
\end{proof}

\renewcommand{\thesubsection}{\S\thesection.\arabic{subsection}}
\subsubsection{}\;     
\renewcommand{\thesubsection}{\thesection.\arabic{subsection}}

Voyons une construction équivalente du foncteur nerf de Segal ${\mathcal{N}_{\mathcal{S}}}$: Soit $\G$ un $2$-groupe et $q\geq 0$ un entier. Si $(X,\alpha)$ et $(Y,\beta)$ sont de $q$-simplexes de $\G$, on définit un \emph{morphisme de $q$-simplexes} $f\colon\xymatrix@C-5pt{(X,\alpha)\ar[r]&(Y,\beta)}$ comme étant la donnée d'une famille de morphismes de $\G$:
$$
\Big\{  f_{ij}\colon\xymatrix{X_{ij}\ar[r]& Y_{ij}}  \,\Big|\,  0\leq i<j\leq q \Big\}
$$
tels que si $ 0\leq i<j<k\leq q$, on a un carré commutatif:
$$
\xymatrix@R-3pt@C+6pt{
  X_{ij}\otimes X_{jk}\ar[r]^-{\alpha_{ijk}}\ar[d]_-{f_{ij}\otimes f_{jk}}&  X_{ik} \ar[d]^-{f_{ik}}\\
  Y_{ij}\otimes Y_{jk}\ar[r]_-{\beta_{ijk}}&  Y_{ik} \,.} 
$$

L'ensemble des $q$-simplexes et des morphismes de $q$-simplexes de $\G$ forment un groupoïde $\underline{\G_q}$ dont la composition est définit argument par argument.  

On va définir un foncteur:
\begin{equation}\label{nervio22cat}
\xymatrix@R=3pt{\text{$2$-${\bf Grp}$}\times\Delta^{op} \ar[rr] && {\bf Grpd}\\ \big(\G,[q]\big) \ar@{}[rr]|-{\longmapsto} && \underline{\G_q}}
\end{equation}
dont le composé avec le foncteur "ensemble des objets" $\xymatrix@C-5pt{{\bf Grpd}\ar[r]&{\bf Ens}}$ est la restriction du foncteur \eqref{nervio22} aux $2$-groupes.

Si $\varphi\colon\xymatrix@-6pt{[q]\ar[r]&[q']}$ est un morphisme de $\Delta$ et $F=(F,m^F)\colon\xymatrix@C-6pt{\G\ar[r]&\H}$ est un morphisme lax et unitaire de $2$-groupes, le foncteur:
\begin{equation}\label{qsimpfun}
\xymatrix@R=2pt{
\G_{q'} \ar[rr]^-{F^\varphi} &&\H_{q}\\
(X,\alpha) \ar[dddd]_-{f}&& (X',\alpha') \ar[dddd]^-{f'} \\&&\\\ar@{}[rr]|-{\longmapsto}&&\\&&\\(Y,\beta) && (Y',\beta')}
\end{equation}
a été déjà défini en objets par les formules:

\begin{equation}\label{lasyy2}
X'_{ij} \, = \, \begin{cases}
\mathbb{1}_{\H} & \text{si \; $\varphi i=\varphi j$} \\
FX_{\varphi i\varphi j}& \text{si \; $\varphi i<\varphi j$} 
\end{cases}
\end{equation}
toujours que $0\leq i< j \leq q$ et:
\begin{equation}\label{lasbebe2}
\alpha'_{ijk} \, = \, \begin{cases}
\ell_{FX_{\varphi i\varphi k}}^{-1}\colon\xymatrix@-10pt{  \mathbb{1} \otimes FX_{\varphi i\varphi k} \ar[r]&   FX_{\varphi i\varphi k}} & \text{si \; $\varphi i=\varphi j\leq\varphi k$}\\
r_{FX_{\varphi i\varphi k}}^{-1}\colon\xymatrix@-10pt{ FX_{\varphi i\varphi k}\otimes \mathbb{1} \ar[r]& FX_{\varphi i\varphi k}}  & \text{si \; $\varphi i\leq\varphi j=\varphi k$}\\
 F(\alpha_{\varphi i\varphi j\varphi k})\circ m^F\colon\xymatrix@-10pt{ FX_{\varphi i\varphi j}\otimes FX_{\varphi j\varphi k} \ar[r]& FX_{\varphi i\varphi k}} & \text{si \; $\varphi i<\varphi j<\varphi k$}
\end{cases}
\end{equation}
pour $0\leq i< j < k \leq q$. 

On définit le foncteur \eqref{qsimpfun} en morphismes par la règle:
$$
f'_{ij} \, = \, 
\begin{cases}
\mathrm{id}_{\mathbb{1}} & \varphi(i)=\varphi(j)\\
F(f_{\varphi i\,\varphi j}) & \varphi(i)<\varphi(j)
\end{cases}
$$
si $0\leq i<j\leq q$.

Pour montrer que $f'$ ainsi défini est effectivement un morphisme de $q$-simplexes, on doit vérifier que:
$$
\xymatrix@R-3pt@C+6pt{
  X'_{ij}\otimes X'_{jk}\ar[r]^-{\alpha'_{ijk}}\ar[d]_-{f'_{ij}\otimes f'_{jk}}&  X'_{ik} \ar[d]^-{f'_{ik}}\\
  Y'_{ij}\otimes Y'_{jk}\ar[r]_-{\beta'_{ijk}}&  Y'_{ik} \,.} 
$$
est un carré commutatif de $\H$ toujours que $0\leq i<j<k\leq q$. Pour cela remarquons que si $0\leq i'<j'<k'\leq q'$ on a un diagramme commutatif:
$$
\xymatrix@C+18pt{
F(X_{i'j'})\otimes F(X_{j'k'})  \ar[r]^-{m^F} \ar[d]_-{F(f_{i'j'})\otimes F(f_{j'k'})}  & 
F(X_{i'j'}\otimes X_{j'k'}) \ar[d]|-{F(f_{i'j'}\otimes f_{j'k'})}  \ar[r]^-{F(\alpha_{i'j'k'})}  &
F(X_{i'k'}) \ar[d]^-{F(f_{i'k'})}   \\
F(Y_{i'j'})\otimes F(Y_{j'k'})  \ar[r]_-{m^F} & F(Y_{i'j'}\otimes Y_{j'k'})  \ar[r]_-{F(\beta_{i'j'k'})} & 
F(Y_{i'k'})  }
$$
et si $0\leq i'<k'\leq q'$ les carrés:
$$
\vcenter{\xymatrix@C+12pt{
F(X_{i'k'}) \ar[r]^-{\ell_{F(X_{i'k'})}}\ar[d]_-{F(f_{i'k'})} & \mathbb{1}\otimes F(X_{i'k'}) \ar[d]^-{\mathbb{1}\otimes F(f_{i'k'})}\\
F(Y_{i'k'}) \ar[r]_-{\ell_{F(Y_{i'k'})}} & \mathbb{1}\otimes F(Y_{i'k'})}}
\qquad \text{et} \qquad 
\vcenter{\xymatrix@C+12pt{
F(X_{i'k'}) \ar[r]^-{r_{F(X_{i'k'})}}\ar[d]_-{F(f_{i'k'})} &  F(X_{i'k'}) \otimes \mathbb{1} \ar[d]^-{F(f_{i'k'})\otimes \mathbb{1}} \\
F(Y_{i'k'}) \ar[r]_-{r_{F(Y_{i'k'})}} &  F(Y_{i'k'})\otimes \mathbb{1} }}
$$
son commutatifs.

Remarquons que le foncteur adjoint:
$$
\xymatrix@R=3pt{\text{$2$-${\bf Grp}$}\ar[rr] && \ssimp \\ 
\G \ar@{}[rr]|-{\longmapsto} && \mathrm{N}\big(\underline{\G_{\bullet_2}}\big)_{\bullet_1}}
$$
du foncteur composé:
$$
\xymatrix@R=3pt{\text{$2$-${\bf Grp}$}\times\Delta^{op} \ar[rr]^-{\eqref{nervio22cat}}&& {\bf Grpd} \ar[rr]^-{\mathrm{N}} && \simp \\ 
\big(\G,[q]\big) \ar@{}[rr]|-{\longmapsto} && \underline{\G_q}\ar@{}[rr]|-{\longmapsto} && \mathrm{N}\big(\underline{\G_q}\big)}
$$
est isomorphe au foncteur nerf de Segal pour le $2$-groupes.

En effet si $p,q\geq 0$ et $\G$ est un $2$-groupe, on vérifie sans peine que les ensembles:
$$
\mathrm{N}\big(\underline{\G_{q}}\big)_{p}\,=\,\mathrm{Hom}_{\bf cat}\big([p],\underline{\G_q}\big) \, = \, \, \text{$p$-simplexes du groupoïde}\;  \underline{\G_q}  
$$ 
$$\text{et}$$
$$
\mathcal{N}_{\mathcal{S}}\big(\G\big)_{p,q} \, = \, \mathcal{N}\big(\G^{[p]}\big)_{q} \, = \, \, \text{$q$-simplexes du $2$-groupe}\;  \G^{[p]}  
$$ 
s'identifient à l'ensemble dont les objets sont les couples formées par une famille de suites de morphismes de $\G$:
\begin{equation}\label{unlemmepart}
\Big\{ \xymatrix{X_{ij}^0 \ar[r]^{f_{ij}^1}&\dots\ar[r]^{f_{ij}^p}&X_{ij}^p}  \;\, \Big| \,\;  0\leq i<j\leq q \;\Big\}
\end{equation}
et une famille de morphismes de la forme:
\begin{equation}\label{doslemmepart}
\Big\{ \xymatrix{ X^s_{ij}\otimes X^s_{jk} \ar[r]^-{\alpha^s_{ijk}}& X_{ik}^s}  \;\, \Big| \,\;  0\leq i<j<k\leq q \; \; 0\leq s\leq p \;\Big\}
\end{equation}
telles que si $0\leq i<j<k\leq q$ et $0\leq s\leq p-1$ on a un carré commutatif:
$$
\xymatrix@C+15pt{
X_{ik}^s\ar[r]^-{f_{ik}^{s+1}}  & X_{ik}^{s+1} \\
X_{ij}^s\otimes X_{jk}^s \ar[r]_-{f^s_{ij}\otimes f^s_{jk}} \ar[u]^-{\alpha^s_{ijk}}  & X_{ij}^{s+1}\otimes X_{jk}^{s+1} \ar[u]_-{\alpha_{ijk}^{s+1}} }
$$
et si $0\leq i<j<k<l\leq q$ et $0\leq s\leq p$ on a un diagramme commutatif:
$$
\xymatrix@R=3pt@C+15pt{
X^s_{il} & X^s_{ij}\otimes X^s_{jl} \ar[l]_-{\alpha^s_{ijl}} \\&\\&\\
 &  X^s_{ij}\otimes \big( X^s_{jk}\otimes X^s_{kl} \big) \ar@{}[d]_-{a}|-{\text{\rotatebox[origin=c]{90}{$\cong$}}} \ar[uuu]_-{X^s_{ij}\otimes\alpha^s_{jkl}}  \\ 
X^s_{ik}\otimes X^s_{kl}  \ar[uuuu]^-{\alpha^s_{ikl}}  & \big(X^s_{ij}\otimes X^s_{jk}\big)\otimes X^s_{kl} \ar[l]^-{\alpha^s_{ijk} \otimes X^s_{kl}}\,.}
$$

Notons aussi que si $\varphi\colon\xymatrix@C-10pt{[p']\ar[r] & [p]}$ et $\psi\colon\xymatrix@C-10pt{[q']\ar[r] & [q]}$ sont deux morphismes de la catégorie simplicial $\Delta$ et $F\colon\xymatrix@C-8pt{\G\ar[r]&\H}$ est un morphisme lax et unitaire de $2$-groupes, l'image par la fonction induite:
$$
\xymatrix{\mathrm{N}\big(\underline{\G_{q}}\big)_{p}\,\cong\,\mathcal{N}_{\mathcal{S}}\big(\G\big)_{p,q}\ar[r] & 
\mathcal{N}_{\mathcal{S}}\big(\H\big)_{p',q'}\,\cong\,\mathrm{N}\big(\underline{\H_{q'}}\big)_{p'}}
$$
d'une couple \eqref{unlemmepart} et \eqref{doslemmepart} vérifiant les propriétés ci-dessus est la couple des familles:
$$
\Big\{\vcenter{\xymatrix{Y_{ab}^0 \ar[r]^{g_{ab}^0}&\dots\ar[r]^{g_{ab}^{p'}}&Y_{ab}^{p'}}} \;\, \Big| \,\;  0\leq a<b\leq q' \;\Big\}
$$
$$\text{et}$$
$$
\Big\{ \xymatrix{Y^t_{ab}\otimes Y^t_{bc}\ar[r]^-{\beta^t_{abc}}& Y_{ac}^t}  \;\, \Big| \,\;  0\leq a<b<c\leq q' \; \; 0\leq t\leq p' \;\Big\}
$$
définies par les formules:
$$
Y_{ab}^t \, = \, 
\begin{cases}
\mathbb{1} \qquad & \psi (a)=\psi (b) \\
F X_{\psi a\psi b}^{\varphi t} & 0\leq \psi (a) < \psi (b) \leq q
\end{cases}  \qquad\;\text{et}
$$
$$
g_{ab}^t \, = \, 
\begin{cases}
\mathrm{id}_{\mathbb{1}} \qquad & \psi (a)=\psi (b) \\
\mathrm{id}_{FX^{\varphi t}_{\psi a \psi b}} \qquad & 0\leq \psi (a) < \psi (b) \leq q \quad \text{et}\quad \varphi(t-1) = \varphi(t) \\
F\big(f_{\psi a\psi b}^{\varphi t}\circ \dots \circ f_{\psi a\psi b}^{\varphi (t-1)+1}\big) \qquad & 0\leq \psi (a) < \psi (b) \leq q \quad \text{et}\quad 0\leq \varphi(t-1) < \varphi(t)\leq p
\end{cases}
$$
si $0\leq a <b\leq q'$ et $0\leq t \leq p'$; et:
$$
\beta^t_{abc} \; = \;
\begin{cases}
\ell_{FX^{\varphi t}_{\psi a\psi c}}^{-1}\colon\xymatrix@-10pt{ \mathbb{1} \otimes FX^{\varphi t}_{\psi a\psi c}  \ar[r]& FX^{\varphi t}_{\psi a\psi c} } & \text{si \; $\psi a=\psi b\leq\psi c$}\\
r^{-1}_{FX^{\varphi t}_{\psi a\psi c}}\colon\xymatrix@-10pt{ FX^{\varphi t}_{\psi a\psi c}\otimes \mathbb{1} \ar[r]& FX^{\varphi t}_{\psi a \psi c}}  & \text{si \; $\psi a\leq\psi b=\psi c$}\\
F\alpha^{\varphi t}_{\psi a\psi b\psi c}\circ m^F\colon\xymatrix@-10pt{FX^{\varphi t}_{\psi a\psi b}\otimes FX^{\varphi t}_{\psi b\psi c}\ar[r]& FX^{\varphi t}_{\psi a\psi c}} & \text{si \; $\psi a<\psi b<\psi c$}
\end{cases}
$$
si $0\leq a <b<c\leq q'$ et $0\leq t \leq p'$.

On a montré:
 
\begin{lemme}\label{sonisomm}
Si $p,q\geq 0$ et $\G$ est un $2$-groupe, il existe une bijection naturelle entre l'ensemble des $q$-simplexes du $2$-groupe $\G^{[p]}$ et l'ensemble des $p$-simplexes du nerf du groupoïde $\underline{G_q}$ des $q$-simplexes du $2$-groupe $\G$:
$$
\mathcal{N}_{\mathcal{S}}\big(\G\big)_{p,q} \, = \, \mathcal{N}\big(\G^{[p]}\big)_{q}  \, \cong \,  \mathrm{N}\big(\underline{\G_{q}}\big)_{p}\,.
$$
\end{lemme}

\renewcommand{\thesubsection}{\S\thesection.\arabic{subsection}}
\subsubsection{}\;
\renewcommand{\thesubsection}{\thesection.\arabic{subsection}}

Posons ${\bf\Delta}\underset{\leq 3}{\times}{\bf\Delta}$ pour noter la sous-catégorie pleine de la catégorie produit ${\bf\Delta}\times{\bf\Delta}$ dont les objets sont les couples $\big([p],[q]\big)$ tels que $p+q\leq 3$ (voir \ref{menorn}) et considérons l'adjonction:
\begin{equation}\label{diagon3}
\xymatrix@C+18pt{
\widehat{\mathbf{{\bf\Delta}}\underset{\leq 3}{\times}\mathbf{{\bf\Delta}}}\phantom{A}\ar@{}[r]|-{\perp\phantom{A}}
\ar@{<-}@/^20pt/[r]^-{\mu_{3}^{\; *}}
\ar@/_20pt/[r]_-{\mu_{3\, *}}& \phantom{A}\underset{\phantom{\leq 2}}{\ssimp} }  \,,
\end{equation}
induite du foncteur d'inclusion canonique:
\begin{equation}\label{mun3}
\xymatrix@C+10pt{{\bf\Delta}\underset{\leq 3}{\times}{\bf\Delta}\;\ar@{^(->}[r]^-{\mu_{3}}&{\bf\Delta}\times{\bf\Delta}}.
\end{equation}

Remarquons que si $A$ est un objet de la catégorie $\widehat{\mathbf{{\bf\Delta}}\underset{\leq 3}{\times}\mathbf{{\bf\Delta}}}$ tel que $A_{p,0}=\star$ si $0\leq p \leq 3$, alors l'ensemble simplicial $\mu_{3\, *}(A)$ vérifie que $\mu_{3\, *}(A)_{p,0}=\star$ pour tout $p\geq 0$. Donc, si on note ${\bf MS}_{\leq 3}$ la sous-catégorie pleine de $\widehat{\mathbf{{\bf\Delta}}\underset{\leq 3}{\times}\mathbf{{\bf\Delta}}}$ dont les objets sont les ensembles bisimpliciaux tronqués $A$ tels que $A_{p,0}=\star$ pour $0\leq p \leq 3$. L'adjonction \eqref{diagon3} induit par restriction une adjonction:
\begin{equation}\label{diagon3re}
\xymatrix@C+18pt{
{\bf MS}_{\leq 3}\phantom{A}\ar@{}[r]|-{\perp\phantom{A}}
\ar@{<-}@/^20pt/[r]^-{\mu_{3}^{\; *}}
\ar@/_20pt/[r]_-{\mu_{3\, *}}& \phantom{A}{\bf MS} }  \,.
\end{equation}

\begin{corollaire}\label{fernandaadnanref}
Le nerf de Segal $\mathcal{N}_{\mathcal{S}}\big(\G\big)$ d'un $2$-groupe $\G$ est dans l'image essentielle du foncteur $\mu_{3\, *}$ de l'adjonction \eqref{diagon3re}. En particulier si $X$ est un objet de ${\bf MS}$ la fonction de restriction: 
$$
\xymatrix{
\mathrm{Hom}_{\bf MS} \big( X, \mathcal{N}_{\mathcal{S}}(\G)\big)\ar[r]  & 
\mathrm{Hom}_{{\bf MS}_{\leq 3}} \big( \mu_{3}^{\; *}(X), \mu_{3}^{\; *}\circ \mathcal{N}_{\mathcal{S}}(\G)\big)\,.}
$$
est bijective.
\end{corollaire}
\begin{proof}
Soit $\G$ un $2$-groupe. Si $p,q\geq 0$ d'après les Lemmes \ref{sonisomm} et \ref{3cosqull} l'ensemble simplicial $\mathcal{N}_{\mathcal{S}}\big(\G\big)_{\bullet,q}$ est $3$-cosquelettique et l'ensemble simpliciaux $\mathcal{N}_{\mathcal{S}}\big(\G\big)_{p,\bullet}$ est $2$-cosquelettique. Il se suit des Lemmes \ref{cosque2} et \ref{derr} qu'il suffit de montrer que le carré:
\begin{equation}\label{carcarteBproof}
\def\objectstyle{\scriptstyle}
\def\labelstyle{\scriptstyle}
\xymatrix@C+18pt{
\mathrm{Hom}_{\ssimp}\big(\Delta^{p}\boxtimes\,\Delta^{q},\mathcal{N}_{\mathcal{S}}(\G)\big) 
\ar[r]\ar[d]&
\mathrm{Hom}_{\ssimp}\big(\partial\Delta^{p}\boxtimes\,\Delta^{q},\mathcal{N}_{\mathcal{S}}(\G)\big)
\ar[d]\,\\
\mathrm{Hom}_{\ssimp}\big(\Delta^{p}\boxtimes\,\partial\Delta^{q},\mathcal{N}_{\mathcal{S}}(\G)\big)
\ar[r]&
\mathrm{Hom}_{\ssimp}\big(\partial\Delta^{p}\boxtimes\,\partial\Delta^{q},\mathcal{N}_{\mathcal{S}}(\G)\big)}
\end{equation}
induit du carré d'ensembles bisimpliciaux:
$$
\xymatrix@-5pt{
\Delta^p\boxtimes\Delta^q & \ar@{_(->}[l]\; \partial\Delta^p\boxtimes\Delta^{q}\\
\ar@{^(->}[u] \Delta^p\boxtimes\partial\Delta^q&\ar@{_(->}[l]\ar@{^(->}[u] \; \partial\Delta^p\boxtimes\partial\Delta^{q}\,,}
$$
est cartésien pour tout $0\leq p\leq 2$ et $0\leq q\leq 3$ tels que $p+q\geq 4$, c'est-á-dire pour les couples $(p,q)\, \in \, \Big\{(2,2), \; (2,3), \; (1,3)\Big\}$.

Pour montrer les cas $(2,2)$ et $(2,3)$ considérons le diagramme:
\begin{equation}\label{doscuadrr}
\def\objectstyle{\scriptstyle}
\def\labelstyle{\scriptstyle}
\xymatrix@C-4pt{
\mathrm{Hom}_{\ssimp}\big(\Delta^{2}\boxtimes\,\Delta^{q},\mathcal{N}_{\mathcal{S}}(\G)\big) 
\ar[r]\ar[d]&
\mathrm{Hom}_{\ssimp}\big(\partial\Delta^{2}\boxtimes\,\Delta^{q},\mathcal{N}_{\mathcal{S}}(\G)\big)
\ar[r]\ar[d]&
\mathrm{Hom}_{\ssimp}\big(\Lambda^{2,1}\boxtimes\,\Delta^{q},\mathcal{N}_{\mathcal{S}}(\G)\big)
\ar[d]
\\
\mathrm{Hom}_{\ssimp}\big(\Delta^{2}\boxtimes\,\partial\Delta^{q},\mathcal{N}_{\mathcal{S}}(\G)\big)
\ar[r]&
\mathrm{Hom}_{\ssimp}\big(\partial\Delta^{2}\boxtimes\,\partial\Delta^{q},\mathcal{N}_{\mathcal{S}}(\G)\big)
\ar[r]&
\mathrm{Hom}_{\ssimp}\big(\Lambda^{2,1}\boxtimes\,\partial\Delta^{q},\mathcal{N}_{\mathcal{S}}(\G)\big)\,,}
\end{equation}
induit du diagramme d'ensembles bisimpliciaux:
$$
\xymatrix@-5pt{
\Delta^2\boxtimes\Delta^q & \ar@{_(->}[l]\; \partial\Delta^2\boxtimes\Delta^{q} & \ar@{_(->}[l]\; \Lambda^{2,1}\boxtimes\Delta^{q}\\
\ar@{^(->}[u] \Delta^2\boxtimes\partial\Delta^q&\ar@{_(->}[l]\ar@{^(->}[u] \; \partial\Delta^2\boxtimes\partial\Delta^{q}
&\ar@{_(->}[l]\ar@{^(->}[u] \; \Lambda^{2,1}\boxtimes\partial\Delta^{q}}
$$

Vu que $\mathcal{N}_{\mathcal{S}}(\G)_{\bullet,q}$ est un $1$-groupoïde de Kan pour $q\geq 0$, on sait que les morphismes composés horizontaux de \eqref{doscuadrr} sont des isomorphismes. En particulier, si $p=2$ et $2\leq q\leq 3$ le carré \eqref{carcarteBproof} est cartésien toujours que la fonction:
$$
\def\objectstyle{\scriptstyle}
\def\labelstyle{\scriptstyle}
\xymatrix@C-5pt{
\mathrm{Hom}_{\ssimp}\big(\partial\Delta^{2}\boxtimes\,\Delta^{q},\mathcal{N}_{\mathcal{S}}(\G)\big) \ar[r] & 
\mathrm{Hom}_{\ssimp}\big(\partial\Delta^{2}\boxtimes\,\partial\Delta^{q},\mathcal{N}_{\mathcal{S}}(\G)\big) \, \times \,
\mathrm{Hom}_{\ssimp}\big(\Lambda^{2,1}\boxtimes\,\Delta^{q},\mathcal{N}_{\mathcal{S}}(\G)\big) 
}
$$
soit injective.

Si $q=3$ remarquons qu'en fait la fonction: 
$$
\def\objectstyle{\scriptstyle}
\def\labelstyle{\scriptstyle}
\xymatrix@C-5pt{
\mathrm{Hom}_{\ssimp}\big(\partial\Delta^{2}\boxtimes\,\Delta^{3},\mathcal{N}_{\mathcal{S}}(\G)\big) \ar[r] & 
\mathrm{Hom}_{\ssimp}\big(\partial\Delta^{2}\boxtimes\,\partial\Delta^{3},\mathcal{N}_{\mathcal{S}}(\G)\big)
}
$$
est injective. En effet, on sait que $\mathcal{N}_{\mathcal{S}}(\G)_{p,\bullet}$ est un $2$-groupoïde de Kan si $p\geq 0$, en particulier la fonction composée:
$$
\def\objectstyle{\scriptstyle}
\def\labelstyle{\scriptstyle}
\xymatrix@C-5pt{
\mathrm{Hom}_{\ssimp}\big(\partial\Delta^{2}\boxtimes\,\Delta^{3},\mathcal{N}_{\mathcal{S}}(\G)\big) \ar[r] & 
\mathrm{Hom}_{\ssimp}\big(\partial\Delta^{2}\boxtimes\,\partial\Delta^{3},\mathcal{N}_{\mathcal{S}}(\G)\big)\ar[r] & 
\mathrm{Hom}_{\ssimp}\big(\partial\Delta^{2}\boxtimes\,\Lambda^{3,1},\mathcal{N}_{\mathcal{S}}(\G)\big)
}
$$
est bijective.

Si $q=2$ montrons que la fonction:
\begin{equation}\label{fonctprodd}
\def\objectstyle{\scriptstyle}
\def\labelstyle{\scriptstyle}
\xymatrix@C-5pt{
\mathrm{Hom}_{\ssimp}\big(\partial\Delta^{2}\boxtimes\,\Delta^{2},\mathcal{N}_{\mathcal{S}}(\G)\big) \ar[r] & 
\mathrm{Hom}_{\ssimp}\big(\partial\Delta^{2}\boxtimes\,\partial\Delta^{2},\mathcal{N}_{\mathcal{S}}(\G)\big) \, \times \,
\mathrm{Hom}_{\ssimp}\big(\Lambda^{2,1}\boxtimes\,\Delta^{2},\mathcal{N}_{\mathcal{S}}(\G)\big) 
}
\end{equation}
est injective. 

Pour commencer remarquons que l'ensemble $\mathrm{Hom}_{\ssimp}\big(\partial\Delta^{2}\boxtimes\,\Delta^{2},\mathcal{N}_{\mathcal{S}}(\G)\big)$ s'identifie à l'ensemble des diagrammes:
\begin{equation}\label{prismes2}
\xymatrix@C-10pt{  
&A_2\otimes A_0 \ar'[d][dd]|(.35){\tau_2\otimes\tau_0} \ar[rrr]^f \ar[dl]_{\varphi_2\otimes \varphi_0} &&&A_1 \ar[dl]|{\varphi_1} \ar[dd]^-{\tau_1}\\
B_2\otimes B_0 \ar@<-2pt>[rrr]|(.65)g \ar[rd]_{\psi_2\otimes \psi_0} &&& B_1 \ar[rd]|-{\psi_1}&\\
&C_2\otimes C_0\ar[rrr]_-h&&&C_1}
\end{equation}
dont les faces en forme de carrés sont de diagrammes commutatifs, mais pas nécessairement les faces en forme de triangles.

D'un autre l'ensemble $\mathrm{Hom}_{\ssimp}\big(\partial\Delta^{2}\boxtimes\,\partial\Delta^{2},\mathcal{N}_{\mathcal{S}}(\G)\big)$ s'identifie à l'ensemble des triplets de triangles pas  nécessairement commutatifs:
$$
\vcenter{\xymatrix@R-5pt@C-10pt{
&A_0\ar[dl]_-{\varphi_0}\ar[dd]^-{\tau_0}\\
B_0\ar[rd]_-{\psi_0}\ar@{}[r]|(.6){?}&\\
&C_0}}\qquad\;
\vcenter{\xymatrix@R-5pt@C-10pt{
&A_1\ar[dl]_-{\varphi_1}\ar[dd]^-{\tau_1}\\
B_1\ar[rd]_-{\psi_1}\ar@{}[r]|(.6){?}&\\
&C_1}}\qquad\text{et}
\qquad
\vcenter{\xymatrix@R-5pt@C-10pt{
&A_2\ar[dl]_-{\varphi_2}\ar[dd]^-{\tau_2}\\
B_2\ar[rd]_-{\psi_2}\ar@{}[r]|(.6){?}&\\
&C_2}}\,,
$$
et l'ensemble $\mathrm{Hom}_{\ssimp}\big(\Lambda^{2,1}\boxtimes\,\Delta^{2},\mathcal{N}_{\mathcal{S}}(\G)\big)$ s'identifie à l'ensemble des diagrammes de la forme:
$$
\xymatrix@C-10pt{  
&A_2\otimes A_0 \ar[rrr]^f \ar[dl]_{\varphi_2\otimes \varphi_0} &&&A_1 \ar[dl]|{\varphi_1} \\
B_2\otimes B_0 \ar[rrr]|-g \ar[rd]_{\psi_2\otimes \psi_0} &&& B_1 \ar[rd]|-{\psi_1}&\\
&C_2\otimes C_0\ar[rrr]_-h&&&C_1\,.} 
$$
où les carrés sont commutatifs.

On déduit sans peine (voir le Lemme \ref{fifidelite}) que la fonction \eqref{fonctprodd} est bijective.

Finalement pour montrer que le carré:
\begin{equation}\label{feramorreal}
\def\objectstyle{\scriptstyle}
\def\labelstyle{\scriptstyle}
\xymatrix@C+18pt{
\mathrm{Hom}_{\ssimp}\big(\Delta^{1}\boxtimes\,\Delta^{3},\mathcal{N}_{\mathcal{S}}(\G)\big) 
\ar[r]\ar[d]&
\mathrm{Hom}_{\ssimp}\big(\partial\Delta^{1}\boxtimes\,\Delta^{3},\mathcal{N}_{\mathcal{S}}(\G)\big)
\ar[d]\,\\
\mathrm{Hom}_{\ssimp}\big(\Delta^{1}\boxtimes\,\partial\Delta^{3},\mathcal{N}_{\mathcal{S}}(\G)\big)
\ar[r]&
\mathrm{Hom}_{\ssimp}\big(\partial\Delta^{1}\boxtimes\,\partial\Delta^{3},\mathcal{N}_{\mathcal{S}}(\G)\big)}
\end{equation}
est cartésien; remarquons que dans le diagramme:
$$
\def\objectstyle{\scriptstyle}
\def\labelstyle{\scriptstyle}
\xymatrix@C+18pt{
\mathrm{Hom}_{\ssimp}\big(\Delta^{1}\boxtimes\,\Delta^{3},\mathcal{N}_{\mathcal{S}}(\G)\big) 
\ar[r]\ar[d]&
\mathrm{Hom}_{\ssimp}\big(\partial\Delta^{1}\boxtimes\,\Delta^{3},\mathcal{N}_{\mathcal{S}}(\G)\big)
\ar[d]\,\\
\mathrm{Hom}_{\ssimp}\big(\Delta^{1}\boxtimes\,\partial\Delta^{3},\mathcal{N}_{\mathcal{S}}(\G)\big)
\ar[r]\ar[d]&
\mathrm{Hom}_{\ssimp}\big(\partial\Delta^{1}\boxtimes\,\partial\Delta^{3},\mathcal{N}_{\mathcal{S}}(\G)\big)\ar[d]\\
\mathrm{Hom}_{\ssimp}\big(\Delta^{1}\boxtimes\,\Lambda^{3,k},\mathcal{N}_{\mathcal{S}}(\G)\big)
\ar[r]&
\mathrm{Hom}_{\ssimp}\big(\partial\Delta^{1}\boxtimes\,\Lambda^{3,k},\mathcal{N}_{\mathcal{S}}(\G)\big)
}
$$
induit du diagramme d'ensembles bisimpliciaux:
$$
\xymatrix@-5pt{
\Delta^1\boxtimes\Delta^3 & \ar@{_(->}[l]\; \partial\Delta^1\boxtimes\Delta^{3} \\
\ar@{^(->}[u] \Delta^1\boxtimes\partial\Delta^3&\ar@{_(->}[l]\ar@{^(->}[u] \; \partial\Delta^1\boxtimes\partial\Delta^{3}\\
\Delta^{1}\boxtimes\Lambda^{3,k}  \ar@{^(->}[u] & \ar@{_(->}[l]\; \partial\Delta^{1}\boxtimes\Lambda^{3,k} \ar@{^(->}[u]\,,
}
$$
les morphismes composés verticaux sont des bijections parce que $\mathcal{N}_{\mathcal{S}}(\G)_{p,\bullet}$ est un $2$-groupoïde de Kan pour $p\geq 0$.

Donc le carré \eqref{feramorreal} est cartésien si la fonction:
\begin{equation} \label{feramorreal2}
\def\objectstyle{\scriptstyle}
\def\labelstyle{\scriptstyle}
\xymatrix@C-5pt{
\mathrm{Hom}_{\ssimp}\big(\Delta^{1}\boxtimes\,\partial\Delta^{3},\mathcal{N}_{\mathcal{S}}(\G)\big) \ar[r] & 
\mathrm{Hom}_{\ssimp}\big(\Delta^{1}\boxtimes\,\Lambda^{3,k},\mathcal{N}_{\mathcal{S}}(\G)\big) \, \times \,
\mathrm{Hom}_{\ssimp}\big(\partial\Delta^{1}\boxtimes\,\partial\Delta^{3},\mathcal{N}_{\mathcal{S}}(\G)\big) 
}
\end{equation} 
est injective. 

Pour montrer cela remarquons que se donner un élément de la source de la fonction \eqref{feramorreal2} équivaut à se donner deux diagrammes pas nécessairement commutatifs de $\mathcal{G}$ de la forme:
\begin{equation}\label{sourceferamorreal2}
\vcenter{\xymatrix@R-16pt@C+8pt{
X_{03}  &X_{01}\otimes X_{13}   \ar[l]_-{\xi_{2}}   \\&\\
&X_{01} \otimes (X_{12}\otimes X_{23})  \ar[uu]_-{X_{01}\otimes \xi_{0} }  \ar@{}[d]_-{a}|-{\text{\rotatebox[origin=c]{90}{$\cong$}}}\\
X_{02}\otimes X_{23}  \ar[uuu]^-{\xi_{1}}   &(X_{01}\otimes X_{12})\otimes X_{23}  \ar[l]^-{\xi_{3}\otimes X_{23}}}}
\qquad\text{et}\qquad
\vcenter{\xymatrix@R-16pt@C+8pt{
X'_{03}  &X'_{01}\otimes X'_{13}   \ar[l]_-{\xi'_{2}}   \\&\\
&X'_{01} \otimes (X'_{12}\otimes X'_{23})  \ar[uu]_-{X'_{01}\otimes \xi'_{0} }  \ar@{}[d]_-{a}|-{\text{\rotatebox[origin=c]{90}{$\cong$}}}\\
X'_{02}\otimes X'_{23}  \ar[uuu]^-{\xi'_{1}}   &(X'_{01}\otimes X'_{12})\otimes X'_{23}  \ar[l]^-{\xi'_{3}\otimes X'_{23}}}}\,,
\end{equation}
et quatre diagrammes commutatifs de $\mathcal{G}$:
\begin{equation} \label{sourceferamorreal22}
\vcenter{\xymatrix{ 
X_{12}\otimes X_{23} \ar[r]^-{\xi_0}  \ar[d]_-{f_{12}\otimes f_{23}}& X_{13} \ar[d]^-{f_{13}} \\
X'_{12}\otimes X'_{23} \ar[r]_-{\xi'_0} & X'_{13}}}\;,
\qquad \qquad
\vcenter{\xymatrix{ 
X_{02}\otimes X_{23} \ar[r]^-{\xi_1}  \ar[d]_-{f_{02}\otimes f_{23}}& X_{03} \ar[d]^-{f_{03}} \\
X'_{02}\otimes X'_{23} \ar[r]_-{\xi'_1} & X'_{03}}}\;,
\end{equation}
\begin{equation}\label{sourceferamorreal23}
\vcenter{\xymatrix{ 
X_{01}\otimes X_{13} \ar[r]^-{\xi_2}  \ar[d]_-{f_{01}\otimes f_{13}}& X_{03} \ar[d]^-{f_{03}} \\
X'_{01}\otimes X'_{13} \ar[r]_-{\xi'_2} & X'_{03}}}
\qquad \text{et} \qquad
\vcenter{\xymatrix{ 
X_{01}\otimes X_{12} \ar[r]^-{\xi_3}  \ar[d]_-{f_{01}\otimes f_{12}}& X_{02} \ar[d]^-{f_{02}} \\
X'_{01}\otimes X'_{12} \ar[r]_-{\xi'_3} & X'_{02}}}\;;
\end{equation} 
ce qui se range dans un cube de $\mathcal{G}$ dont quatre des six faces sont commutatifs.

Plus encore la première coordonné de l'image de \eqref{sourceferamorreal2}, \eqref{sourceferamorreal22} et \eqref{sourceferamorreal23} par la fonction \eqref{feramorreal2} oubli les morphismes $\xi_k$ et $\xi'_k$, et la deuxième coordonné oubli la commutativité des quatre diagrammes \eqref{sourceferamorreal22} et \eqref{sourceferamorreal23}. Donc \eqref{feramorreal2} est bien une fonction injective.
\end{proof}

Montrons:

\begin{corollaire}\label{pfNerSegal}
Le foncteur nerf de Segal pour les $2$-groupes $\mathcal{N}_{\mathcal{S}}\colon\xymatrix{ \text{$2$-${\bf Grp}$}\ar[r] & {\bf MS}}$ est pleinement fidèle.
\end{corollaire}
\begin{proof}
Soient $F, G\colon \xymatrix@C-8pt{\mathcal{G}\ar[r] & \mathcal{H}}$ deux morphisme de $2$-groupes tels que $\mathcal{N}_{\mathcal{S}}(F) \, = \, \mathcal{N}_{\mathcal{S}}(G)$. Vu que le foncteur nerf $\mathcal{N}\colon\xymatrix{ \text{$2$-${\bf Grp}$}\ar[r] & \simp}$ est fidèle (voir le Corollaire \ref{plfimonlaxun}) et on a que:
$$
\mathcal{N}(F) \, = \, \mathcal{N}_{\mathcal{S}}(F)_{0,\bullet} \, = \, \mathcal{N}_{\mathcal{S}}(G)_{0,\bullet} \, = \, \mathcal{N}(G)\,;
$$
il se suit que $F=G$

Si d'un autre $f\colon \xymatrix@C-8pt{\mathcal{N}_{\mathcal{S}}(\mathcal{G}) \ar[r] & \mathcal{N}_{\mathcal{S}}(\mathcal{H})}$ est un morphisme d'ensembles bisimpliciaux, on sait qu'il existe un morphisme de $2$-groupes $F\colon \xymatrix@C-8pt{\mathcal{G}\ar[r] & \mathcal{H}}$ tel que $f_{0,\bullet} \, = \, \mathcal{N}(F) \, = \, \mathcal{N}_{\mathcal{S}}(F)_{0,\bullet}$, parce que le foncteur $\mathcal{N}\colon\xymatrix{ \text{$2$-${\bf Grp}$}\ar[r] & \simp}$ est plein.

D'un autre vu qu'on a les égalités des foncteurs: 
$$
\xymatrix@C+27pt{\text{$2$-${\bf Grp}$}\ar[r]^-{(\mathcal{N}_{\mathcal{S}})_{0,\bullet} \, = \, \mathcal{N}}&\simp_0}
\qquad\text{et}\qquad 
\xymatrix@C+27pt{\text{$2$-${\bf Grp}$}\ar[r]^-{(\mathcal{N}_{\mathcal{S}})_{\bullet,1} \, = \, \mathrm{N}\circ s}&\simp}
$$ 
où $s\colon\xymatrix@C-12pt{\text{$2$-${\bf Grp}$} \ar[r] & {\bf Grpd}}$ note le foncteur groupoïde sous-jacent, il se suit du Corollaire \ref{isoGAMMA} qu'on a un carré commutatif:
$$
\xymatrix@C+10pt{
\mathbb{\Omega}_\star\big((\mathcal{N}_{\mathcal{S}}\mathcal{G})_{0,\bullet} \big) \ar[d]_-{\mathbb{\Omega}_\star\big((\mathcal{N}_{\mathcal{S}}F)_{0,\bullet}\big)}\ar[r]^-{\Gamma_\mathcal{G}} & 
\big(\mathcal{N}_{\mathcal{S}}\mathcal{G}\big)_{\bullet,1}\ar[d]^-{(\mathcal{N}_{\mathcal{S}} F)_{\bullet,1}}\\
\mathbb{\Omega}_\star\big((\mathcal{N}_{\mathcal{S}}\mathcal{H})_{0,\bullet} \big) \ar[r]_-{\Gamma_\mathcal{H}} & \big(\mathcal{N}_{\mathcal{S}}\mathcal{H}\big)_{\bullet,1}
}
$$

Donc pour montrer qu'on a $f_{\bullet,1} \, = \,  \mathcal{N}_{\mathcal{S}}(F)_{\bullet,1}$ il suffit de montrer qu'on a un carré commutatif pour $0\leq n\leq 1$:
\begin{equation}\label{naturaliteexten}
\xymatrix@C+15pt{
\mathbb{\Omega}_\star\big((\mathcal{N}_{\mathcal{S}}\mathcal{G})_{0,\bullet} \big)_n \ar[d]_-{\mathbb{\Omega}_\star(f_{0,\bullet})_n}\ar[r]^-{(\Gamma_\mathcal{G})_n} & 
\big(\mathcal{N}_{\mathcal{S}}\mathcal{G}\big)_{n,1}\ar[d]^-{f_{n,1}}\\
\mathbb{\Omega}_\star\big((\mathcal{N}_{\mathcal{S}}\mathcal{H})_{0,\bullet} \big)_n \ar[r]_-{(\Gamma_\mathcal{H})_n} & \big(\mathcal{N}_{\mathcal{S}}\mathcal{H}\big)_{n,1}}
\end{equation}
car $f_{0,\bullet} \, = \,  \mathcal{N}_{\mathcal{S}}(F)_{0,\bullet}$ et $\big(\mathcal{N}_{\mathcal{S}}\mathcal{H}\big)_{\bullet,1}$ est un ensemble simplicial faiblement $1$-cosquelettique.

Si $n=0$ le carré \eqref{naturaliteexten} est commutatif parce que $(\Gamma_\mathcal{G})_0$ et $(\Gamma_\mathcal{H})_0$ sont les fonctions identité entre les ensembles des objets de $\mathcal{G}$ et $\mathcal{H}$ respectivement, tandis que $f_{0,1}\, = \, \mathbb{\Omega}_\star(f_{0,\bullet})_0$ est la fonction induit par $F$ de l'ensemble des objets $\mathcal{G}$ vers l'ensemble des objets de $\mathcal{H}$.

Si $n=1$ montrons que le carré \eqref{naturaliteexten} est commutatif: Si $\xi\colon\xymatrix@C-10pt{\mathbb{1}\otimes X\ar[r]& Y}$ est un $1$-simplexe de l'ensemble simplicial $\mathbb{\Omega}_\star\big((\mathcal{N}_{\mathcal{S}}\mathcal{G})_{0,\bullet} \big)$, considérons le $(1,2)$-simplexe suivant de l'ensemble simplicial $\mathcal{N}_{\mathcal{S}}(\mathcal{G})$:
$$
\xymatrix@C+13pt{
\mathbb{1}\otimes X \ar[r]^-{\ell_X^{-1}} & X \\
\mathbb{1}\otimes X \ar[u]^-{\mathrm{id}\,=\,\mathbb{1}\otimes X}\ar[r]_-{\xi} & Y\ar[u]_-{\ell_X^{-1}\circ \xi^{-1}} 
}
$$
et son image par la fonction $f_{1,2}$:
$$
\xymatrix@C+40pt{
\mathbb{1}\otimes f_{0,1}(X) \ar[r]^-{f_{0,2}(\ell_X^{-1}) \, = \, \ell^{-1}_{f_{0,1}(X)}} & f_{0,1}(X) \\
\mathbb{1}\otimes f_{0,1}(X) \ar[u]^-{\mathrm{id}\,=\,\mathbb{1}\otimes f_{0,1}(X)}\ar[r]_-{f_{0,2}(\xi)} & f_{0,1}(Y)\ar[u]_-{f_{1,1}(\ell_X^{-1}\circ \xi^{-1})} 
}
$$

On déduit que dans le $2$-groupe $\mathcal{H}$ on a que:
$$
f_{1,1}\big(\ell_X^{-1}\circ \xi^{-1}\big) \; = \; \ell^{-1}_{f_{0,1}(X)} \circ f_{0,2}(\xi)^{-1}\,;
$$
autrement dit \eqref{naturaliteexten} est un carré commutatif pour $n=1$:
\begin{align*}
f_{1,1}\circ \big(\Gamma_{\mathcal{G}}\big)_1 (\xi) \, &= \, f_{1,1}\big(\ell_X^{-1}\circ \xi^{-1})   \\
\, &= \,    \ell^{-1}_{f_{0,1}(X)} \circ f_{0,2}(\xi)^{-1} \\
\, &= \,    \big(\Gamma_{\mathcal{H}}\big)_1\big(f_{0,2}(\xi)\big)    \\
\, &= \,   \big(\Gamma_{\mathcal{H}}\big)_1\circ \mathbb{\Omega}_\star\big(f_{0,\bullet}\big)_1 (\xi)\,.
\end{align*}

Donc $f_{\bullet,1} = \mathcal{N}_{\mathcal{S}}(F)_{\bullet,1}$.

Remarquons finalement que $f_{1,2} \, = \,  \mathcal{N}_{\mathcal{S}}(F)_{1,2}$ parce que \eqref{fonctionfinalll} est une fonction bijective pour tout $2$-groupe $\G$ et on a déjà montré que $f_{p,q} = \mathcal{N}_{\mathcal{S}}(F)_{p,q}$ si $(p,q)\,\in\, \big\{ \, (0,2), \, (0,1), \, (1,1)\, \big\}$.
\end{proof}

\renewcommand{\thesubsection}{\S\thesection.\arabic{subsection}}
\subsubsection{}\;\label{deuxfoncteurnerfSEGAL}     
\renewcommand{\thesubsection}{\thesection.\arabic{subsection}}

Rappelons que la catégorie ${\bf Grpd}^{\Delta^{op}}$ des foncteurs de $\Delta^{op}$ vers la catégorie des groupoïdes ${\bf Grpd}$ admet un structure canonique de $2$-catégorie $\underline{{\bf Grpd}}^{\Delta^{op}}$ dont les $2$-flèches son définies comme suit: Si $X$ et $Y$ sont des objets de ${\bf Grpd}^{\Delta^{op}}$ considérons deux flèche de $X$ vers $Y$ dans ${\bf Grpd}^{\Delta^{op}}$:
$$
\xymatrix@C+5pt{
\Delta^{op} \; \ar@<+3pt>@/^11pt/[rr]^-{X}  \ar@<-3pt>@/_11pt/[rr]_-{Y} &{\Downarrow F,\; G}& \, {\bf Grpd} 
}
$$

Une $2$-flèche $\Gamma$ de $F$ vers $G$ dans $\underline{{\bf Grpd}}^{\Delta^{op}}$ est une famille de transformations naturelles:
$$
\alpha \, = \, \left\{
\vcenter{\xymatrix@C-8pt{ 
X_q \ar@<+2pt>@/^10pt/[rr]^-{F_q}  \ar@<-2pt>@/_10pt/[rr]_-{G_q} &{\Downarrow \Gamma_q}& Y_q
}}\;\right\}_{n \in \mathbb{N}}
$$
telle que:
$$
\left(\vcenter{\xymatrix@C-11pt{ 
X_{q+1} \ar@<+2pt>@/^10pt/[rr]^-{F_{q+1}}  \ar@<-2pt>@/_10pt/[rr]_-{G_{q+1}} &{\Downarrow \Gamma_{q+1}}& Y_{q+1}\ar[rr]^-{d_i} && Y_{q} 
}}\right)\;=\;
\left(\vcenter{\xymatrix@C-11pt{ 
X_{q+1} \ar[rr]^-{d_i}& &X_q \ar@<+2pt>@/^10pt/[rr]^-{F_q}  \ar@<-2pt>@/_10pt/[rr]_-{G_q} &{\Downarrow  \Gamma_q}& Y_q
}}\right)
$$
si $0\leq i\leq q+1$ et:
$$
\left(\vcenter{\xymatrix@C-11pt{ 
X_{q} \ar@<+2pt>@/^10pt/[rr]^-{F_{q}}  \ar@<-2pt>@/_10pt/[rr]_-{G_{q}} &{\Downarrow \Gamma_q}& Y_{q}\ar[rr]^-{s_j} && Y_{q+1} 
}}\right)\;=\;
\left(\vcenter{\xymatrix@C-11pt{ 
X_{q} \ar[rr]^-{s_j}&& X_{q+1} \ar@<+2pt>@/^10pt/[rr]^-{F_{q+1}}  \ar@<-2pt>@/_10pt/[rr]_-{G_{q+1}} &{\Downarrow \Gamma_{q+1}}& Y_{q+1}
}}\right)
$$
si $0\leq j\leq q$.

On va définir un $2$-foncteur:
\begin{equation}\label{nervio22catenrri}
\xymatrix@R=3pt{\text{$2$-$\underline{{\bf Grp}}$} \ar[rr]^-{\underline{\mathcal{N}}_{\mathcal{S}}} && \underline{{\bf Grpd}}^{\Delta^{op}}\\ \G \ar@{}[rr]|-{\longmapsto} && \underline{\G_\bullet}}
\end{equation}
dont le foncteur sous-jacent est l'adjoint du foncteur \eqref{nervio22cat}. 

Si $\xymatrix@C-2pt{\mathcal{G}\rtwocell^{F}_G{\eta}& \mathcal{H}}$ est une transformation entre morphismes de $2$-groupes on définit la $2$-flèches:
$$
\xymatrix@C-12pt{
\mathcal{N}_{\mathcal{S}}(\mathcal{G}) \ar@<+2pt>@/^10pt/[rr]^-{\mathcal{N}_{\mathcal{S}}(F)}  \ar@<-2pt>@/_10pt/[rr]_-{\mathcal{N}_{\mathcal{S}}(G)} &{\phantom{.}\quad\Downarrow \;\; \text{\scriptsize $\underline{\mathcal{N}}_{\mathcal{S}}(\eta)$}}& \mathcal{N}_{\mathcal{S}}(\mathcal{H})} \; = \;
\left\{ \; \vcenter{\xymatrix@C+20pt{\underline{\mathcal{G}}_q\rtwocell^{F_q}_{G_q}{\,\;\eta_q}& \underline{\mathcal{H}}_q}} \; \right\}_{q\geq 0}
$$
de la façon suivante: Si $q\geq 0$ et $(X,\alpha)$ est un $q$-simplexe de $\mathcal{G}$ le morphisme de $q$-simplexes de $\mathcal{H}$ : 
$$
\xymatrix@C+20pt{F_q(X,\alpha) \ar[r]^-{(\eta_q)_{(X,\alpha)}}& G_q(X,\alpha)}
$$ 
est égal à la famille de morphismes:
\begin{equation}\label{definietatodos}
(\eta_q)_{(X,\alpha)}  \, = \, \Bigg\{\;\xymatrix@C+25pt{F(X_{ij})\ar[r]^{\eta_{X_{ij}} } & G(X_{ij})}\;\Bigg|\; 0\leq i < j \leq q \;\Bigg\}\,;
\end{equation}
laquelle est bien un morphisme de $q$-simplexes parce qu'on a le diagramme commutatif:
$$
\xymatrix@C+10pt{
F(X_{ij})\otimes F(X_{jk})  \ar[d]_-{\eta_{X_{ij}} \otimes \eta_{X_{jk}}} \ar[r]^-{m^F} & F(X_{ij}\otimes X_{jk})  \ar[r]^-{F(\alpha_{ijk})}  
\ar[d]|-{\eta_{X_{ij}\otimes X_{jk}}} & F(X_{ik}) \ar[d]^-{\eta_{X_{ik}}}\\
G(X_{ij})\otimes G(X_{jk}) \ar[r]_-{m^G} & G(X_{ij}\otimes X_{jk})  \ar[r]_-{G(\alpha_{ijk})}  & G(X_{ik}) }
$$
toujours que $0\leq i<j<k\leq q$.

On vérifie sans peine que \eqref{nervio22catenrri} ainsi défini est bien un $2$-foncteur.

\begin{corollaire}\label{2pleinefideleNS}
Le $2$-foncteur nerf de Segal $\underline{\mathcal{N}}_{\mathcal{S}}\colon\xymatrix@C-10pt{\text{$2$-$\underline{\bf Grp}$}\ar[r]&\underline{\bf Grpd}^{\Delta^{op}}}$ est pleinement fidèle dans le sens que si $\mathcal{G}$ et $\mathcal{H}$ sont des $2$-groupes le foncteur:
\begin{equation}\label{isogrpdSegal}
\xymatrix@C+20pt
{\underline{\mathrm{Hom}}_{\text{$2$-${\bf Grp}$}}\big(\mathcal{G},\mathcal{H}\big)\ar[r]&
\underline{\mathrm{Hom}}_{{\bf Grpd}^{\Delta^{op}}}\big(\mathcal{N}_{\mathcal{S}}(\mathcal{G}),\mathcal{N}_{\mathcal{S}}(\mathcal{H})\big)
}
\end{equation}
est un isomorphisme de groupoïdes.
\end{corollaire}
\begin{proof}
D'après le Corollaire \ref{pfNerSegal} le foncteur \eqref{isogrpdSegal} est une bijection entre les ensembles des objets. Montrons que c'est un foncteur pleinement fidèle. 

Si $\xymatrix@C+2pt{\mathcal{G}\rtwocell^{F}_G{\eta}& \mathcal{H}}$ et $\xymatrix@C+2pt{\mathcal{G}\rtwocell^{F}_G{\eta'}& \mathcal{H}}$ sont deux transformations entre morphismes de $2$-groupes telles que $\underline{\mathcal{N}}_{\mathcal{S}}(\eta) \, = \, \underline{\mathcal{N}}_{\mathcal{S}}(\eta')$, il se suit que pour tout $1$-simplexe $X$ de $\mathcal{G}$ \emph{i.e.} pour tout objet $X$ de $\mathcal{G}$ on a que: 
$$
\eta_X \, = \, \big(\underline{\mathcal{N}}_{\mathcal{S}}(\eta)_1\big)_X \, = \, \big(\underline{\mathcal{N}}_{\mathcal{S}}(\eta')_1\big)_X  \, = \, \eta'_X\,.
$$ 

Donc $\eta=\eta'$. 

Considérons maintenant deux morphismes de $2$-groupes $F, G\colon \xymatrix@C-8pt{\mathcal{G}\ar[r] & \mathcal{H}}$ et supposons qu'on s'est donné une $2$-flèche de la $2$-catégorie $\underline{\bf Grpd}^{\Delta^{op}}$:
$$
\xymatrix{\mathcal{N}_{\mathcal{S}}(\mathcal{G}) \ar@<+2pt>@/^10pt/[rr]^-{\mathcal{N}_{\mathcal{S}}(F)}  \ar@<-2pt>@/_10pt/[rr]_-{\mathcal{N}_{\mathcal{S}}(G)} &{\phantom{.}\quad\Downarrow \;\; \text{\scriptsize $\Gamma$}} & \mathcal{N}_{\mathcal{S}}(\mathcal{H})}  \; = \;
\left\{ \; \vcenter{\xymatrix@C+20pt{\underline{\mathcal{G}}_q\rtwocell^{F_q}_{G_q}{\,\;\Gamma_q}& \underline{\mathcal{H}}_q}} \; \right\}_{q\geq 0}\,.
$$

On va définir une $2$-flèche $\xymatrix@C+2pt{\mathcal{G}\rtwocell^{F}_G{\eta}& \mathcal{H}}$ de la $2$-catégorie $2$-$\underline{\bf Grp}$ telle que $\underline{\mathcal{N}}_{\mathcal{S}}(\eta)\, = \, \Gamma$.

Si $X$ est un objet de $\mathcal{G}$ \emph{i.e.} si $X$ est un $1$-simplexe de $\mathcal{G}$, on définit le morphisme $\eta_X$ de $\mathcal{H}$ comme suit:
$$
\xymatrix@C+28pt{F(X)\ar[r]^-{\eta_X\,=\,(\Gamma_1)_X} & G(X)}
$$

\emph{Montrons que $\eta_\mathbb{1}=\mathrm{id}_\mathbb{1}$:}

Vu qu'on a l'égalité de transformations naturelles:
$$
\left(\vcenter{\xymatrix@C-11pt{ 
\underline{\mathcal{G}}_{0} \ar@<+2pt>@/^10pt/[rr]^-{F_{0}}  \ar@<-2pt>@/_10pt/[rr]_-{G_{0}} &{\Downarrow \text{\scriptsize $\Gamma_0$}}& \underline{\mathcal{H}}_{0} \ar[rr]^-{s_0} && \underline{\mathcal{H}}_{1} 
}}\right)\;=\;
\left(\vcenter{\xymatrix@C-11pt{ 
\underline{\mathcal{G}}_{0} \ar[rr]^-{s_0}&& \underline{\mathcal{G}}_{1}  \ar@<+2pt>@/^10pt/[rr]^-{F_{1}}  \ar@<-2pt>@/_10pt/[rr]_-{G_{1}} &{\Downarrow \text{\scriptsize $\Gamma_1$}}& \underline{\mathcal{H}}_{1} 
}}\right)
$$
où $\Gamma_0$ est la seul transformation naturelle entre le foncteur identité $F_0=G_0$ de la catégorie ponctuelle $\underline{\mathcal{G}}_0=\underline{\mathcal{H}}_0$, on déduit que $(\Gamma_1)_{\mathbb{1}}= \mathrm{id}_{\mathbb{1}}$. 

\emph{Montrons que $\eta$ est une transformation naturelle:}

Soit $f\colon\xymatrix@C-5pt{X \ar[r]& Y}$ un morphisme de $\mathcal{G}$. Si on considère les $2$-simplexes $r_X^{-1}\colon \xymatrix@C-5pt{X\otimes \mathbb{1} \ar[r] & X }$ et $f\circ r_X^{-1}\colon\xymatrix@C-5pt{X\otimes \mathbb{1} \ar[r]& Y}$ de $\mathcal{G}$ on déduit qu'on a les diagrammes commutatifs de $\mathcal{H}$:
$$
\vcenter{\xymatrix@C+15pt{
F(X)\otimes\mathbb{1}\ar@<+2pt>@/^20pt/[rr]^-{r_{FX}^{-1}}\ar[d]_-{(\Gamma_1)_X\otimes (\Gamma_1)_\mathbb{1}} \ar[r]|-{m^F} & F(X\otimes \mathbb{1}) \ar[r]|-{F(r_X^{-1})} & F(X) \ar[d]^-{(\Gamma_1)_X}\\
G(X)\otimes\mathbb{1} \ar[r]|-{m^G} \ar@<-2pt>@/_20pt/[rr]_-{r_{GX}^{-1}} & G(X\otimes\mathbb{1}) \ar[r]|-{G(r_X^{-1})} & G(X)}}\qquad \qquad\text{et}
$$
$$
\xymatrix@C+15pt{
F(X)\otimes\mathbb{1}\ar@<+2pt>@/^20pt/[rr]^-{r_{FX}^{-1}}\ar[d]_-{(\Gamma_1)_X\otimes (\Gamma_1)_\mathbb{1}} \ar[r]|-{m^F} & F(X\otimes \mathbb{1}) \ar[r]|-{F(r_X^{-1})} & F(X)  \ar[r]^-{Ff} & F(Y)\ar[d]^-{(\Gamma_1)_Y}\\
G(X)\otimes\mathbb{1} \ar[r]|-{m^G} \ar@<-2pt>@/_20pt/[rr]_-{r_{GX}^{-1}} & G(X\otimes\mathbb{1}) \ar[r]|-{G(r_X^{-1})} & G(X) \ar[r]_-{Gf} & G(Y)\,.}
$$

Donc on a un carré commutatif: 
$$
\xymatrix@C+10pt{
F(X) \ar[r]^-{Ff} \ar[d]_-{(\Gamma_1)_X \,=\, \eta_X} & F(Y) \ar[d]^-{(\Gamma_1)_Y \,=\, \eta_Y}\\
G(X) \ar[r]_-{Gf} & G(Y)
}
$$

\emph{Montrons que $\eta$ est une transformation entre morphismes de $2$-groupes:}

Si $X$ et $Y$ sont des objets de $\mathcal{G}$ considérons le $2$-simplexe $\mathrm{id}_{X\otimes Y}\xymatrix@C-8pt{X\otimes Y\ar[r] & X\otimes Y}$ de $\mathcal{G}$ . On déduit qu'on a un diagramme commutatif de $\mathcal{H}$: 
$$
\xymatrix@C+22pt{
F(X)\otimes F(Y)\ar[d]_-{(\Gamma_1)_X\otimes (\Gamma_1)_Y} \ar[r]^-{m_{X,Y}^F} & F(X\otimes Y) \ar@{=}[r]^-{F(\mathrm{id}_{X\otimes Y})} & F(X\otimes Y) \ar[d]^-{(\Gamma_1)_{X\otimes Y}}\\
G(X)\otimes G(Y) \ar[r]_-{m^G_{X,Y}} & G(X\otimes Y) \ar@{=}[r]_-{G(\mathrm{id}_{X\otimes Y})} & G(X\otimes Y)\,;}
$$
c'est-à-dire on a un carré commutatif:
$$
\xymatrix@C+10pt{
F(X)\otimes F(Y) \ar[r]^-{m^F_{X,Y}} \ar[d]_-{(\Gamma_1)_X \otimes  (\Gamma_1)_Y \,=\, \eta_X\otimes \eta_Y} & F(X\otimes Y)\ar[d]^-{\eta_{X\otimes Y} \, = \, (\Gamma_1)_{X\otimes Y}} \\
G(X)\otimes G(Y) \ar[r]_-{m^G_{X,Y}} & G(X\otimes Y)\,.
}
$$

\emph{Montrons que $\underline{\mathcal{N}}_{\mathcal{S}}(\eta)\, =\, \Gamma$:}

On a que $\underline{\mathcal{N}}_{\mathcal{S}}(\eta)_0\, =\, \Gamma_0$ est la seule transformation naturelle du foncteur identité $F_0=G_0$ de la catégorie ponctuelle $\underline{\mathcal{G}}_0 = \underline{\mathcal{H}}_0$. En plus par définition $\underline{\mathcal{N}}_{\mathcal{S}}(\eta)_1\, = \, \eta \, =\, \Gamma_1$. 

D'un autre remarquons que si $q\geq 2$ on a l'égalité des transformations naturelles:
$$
\left(\vcenter{\xymatrix@C-11pt{ 
\underline{\mathcal{G}}_{q} \ar@<+2pt>@/^10pt/[rr]^-{F_{q}}  \ar@<-2pt>@/_10pt/[rr]_-{G_{q}} &{\Downarrow \text{\scriptsize $\Gamma_q$}}& \underline{\mathcal{H}}_{q} \ar[rr]^-{d_{i_{q-1}}\dots d_{i_1}} && \underline{\mathcal{H}}_{1} 
}}\right)\;=\;
\left(\vcenter{\xymatrix@C-11pt{ 
\underline{\mathcal{G}}_{q} \ar[rr]^-{d_{i_{q-1}}\dots d_{i_1}}&& \underline{\mathcal{G}}_{1}  \ar@<+2pt>@/^10pt/[rr]^-{F_{1}}  \ar@<-2pt>@/_10pt/[rr]_-{G_{1}} &{\Downarrow \text{\scriptsize $\Gamma_1$}}& \underline{\mathcal{H}}_{1} 
}}\right)\,;
$$
alors si $(X,\alpha)$ est un $q$-simplexe de $\G$ il se suit que le morphisme de $q$-simplexes de $\H$:
$$
\xymatrix@C+25pt{F_q(X,\alpha) \ar[r]^-{(\Gamma_q)_{(X,\alpha)}}& G_q(X,\alpha)}
$$
est défini par la famille de morphismes de $\H$:
\begin{equation}\label{definietatodos2}
(\Gamma_q)_{(X,\alpha)}  \, = \, \Bigg\{\;\xymatrix@C+25pt{F(X_{ij})\ar[r]^{(\Gamma_1)_{X_{ij}} } & G(X_{ij})}\;\Bigg|\; 0\leq i < j \leq q \;\Bigg\}\,.
\end{equation}

Donc $\underline{\mathcal{N}}_{\mathcal{S}}(\eta)_q\, =\, \Gamma_q$ si $q\geq 2$ d'après la définition \eqref{definietatodos}.
\end{proof}

Considérons maintenant le foncteur pleinement fidèle: 
\begin{equation}\label{grpdcatsimp}
\vcenter{\xymatrix@C+18pt@R=2pt{
{\bf Grpd}^{\Delta^{op}} \ar[r]^{\mathrm{N}^{\Delta^{op}}} & \ssimp \\
\G_{\bullet} \ar@{}[r]|-{\longmapsto} & \mathrm{Hom}_{\bf cat}\big(\Delta^{\bullet_2},\G_{\bullet_1}\big)
}}\,,
\end{equation}
déduit du foncteur nerf des petites catégories $\mathrm{N}\colon\xymatrix@C-10pt{{\bf cat}\ar[r]& \simp}$. 

On vérifie sans peine que l'isomorphisme de foncteurs:
$$
\vcenter{\xymatrix@C-10pt{
\mathrm{Hom}_{{\bf Grpd}^{\Delta^{op}}}\big(\bullet_1,\bullet_2\big)\ar@{=>}[r] & 
\mathrm{Hom}_{\ssimp}\Big(\mathrm{N}^{\Delta^{op}}(\bullet_1),\mathrm{N}^{\Delta^{op}}(\bullet_2)\Big)}}\colon 
\vcenter{\xymatrix{\big({{\bf Grpd}^{\Delta^{op}}}\big)^{op}\times {\bf Grpd}^{\Delta^{op}}\ar[r]& {\bf Ens}}}
$$
défini par le foncteur \eqref{grpdcatsimp} s'étend  en un isomorphisme naturel entre foncteurs de la catégorie $\big({{\bf Grpd}^{\Delta^{op}}}\big)^{op}\times {\bf Grpd}^{\Delta^{op}}$ vers la catégorie des ensembles simpliciaux $\simp$: 
\begin{equation} \label{elisoquebesofer}
\vcenter{\xymatrix@C-10pt{
\mathrm{N}\Big(\underline{\mathrm{Hom}}_{{\bf Grpd}^{\Delta^{op}}}\big(\bullet_1,\bullet_2\big)\Big)      \ar@{=>}[r] & 
\underline{\mathrm{Hom}}_{\ssimp}^{(1)}\Big(\mathrm{N}^{\Delta^{op}}(\bullet_1),\mathrm{N}^{\Delta^{op}}(\bullet_2)\Big)}}\,,
\end{equation}
vu que d'après \eqref{isonhomhom} on a un isomorphisme:
$$
\vcenter{\xymatrix@C-10pt{
\mathrm{N}\Big(\underline{\mathrm{Hom}}_{\bf Grpd}\big(\bullet_1,\bullet_2\big)\Big)\ar@{=>}[r] & 
\underline{\mathrm{Hom}}_{\simp}\big(\mathrm{N}(\bullet_1),\mathrm{N}(\bullet_2)\big)}}\colon \vcenter{\xymatrix{{\bf Grpd}^{op}\times {\bf Grpd} \ar[r]& \simp}}\,.
$$

On déduit du Corollaire \ref{2pleinefideleNS} et de l'isomorphisme \eqref{elisoquebesofer}:

\begin{corollaire}\label{pleinffMS1}
L'isomorphisme naturel de foncteurs:
$$
\vcenter{\xymatrix@C-10pt{
\mathrm{Hom}_{\text{$2$-${\bf Grp}$}}\big(\bullet_1,\bullet_2\big)\ar@{=>}[r] & 
\mathrm{Hom}_{\bf MS}\Big(\mathcal{N}_{\mathcal{S}}(\bullet_1),\mathcal{N}_{\mathcal{S}}(\bullet_2)\Big)}}\colon 
\vcenter{\xymatrix{\text{$2$-${\bf Grp}$}^{op}\times\text{$2$-${\bf Grp}$}\ar[r]& {\bf Ens}}}
$$
défini par le foncteur nerf de Segal $\mathcal{N}_{\mathcal{S}}\colon\xymatrix@C-12pt{\text{$2$-${\bf Grp}$} \ar[r] & {\bf MS}}$ (voir le Corollaire \ref{pfNerSegal}) s'étend  en un isomorphisme naturel entre foncteurs de la catégorie $\text{$2$-${\bf Grp}$}^{op}\times\text{$2$-${\bf Grp}$}$ vers la catégorie des ensembles simpliciaux $\simp$: 
\begin{equation} \label{eqpleinffMS1}
\vcenter{\xymatrix@C-10pt{
\mathrm{N}\Big(\underline{\mathrm{Hom}}_{\text{$2$-${\bf Grp}$}}\big(\bullet_1,\bullet_2\big)\Big)      \ar@{=>}[r] & 
\underline{\mathrm{Hom}}_{\bf MS}^{(1)}\Big(\mathcal{N}_{\mathcal{S}}(\bullet_1),\mathcal{N}_{\mathcal{S}}(\bullet_2)\Big)}}\,.
\end{equation}
\end{corollaire}

Rappelons que la catégorie homotopique des $2$-groupes $2$-$h{\bf Grp}$ (voir \ref{2groupessection}) est la catégorie dont les objets sont les $2$-groupes et les morphismes sont les classes à isomorphismes près des morphismes de $2$-groupes:
$$
\mathrm{Hom}_{\text{$2$-$h{\bf Grp}$}} \, = \, \pi_0\big(\underline{\mathrm{Hom}}_{\text{$2$-${\bf Grp}$}}\big).
$$

On déduit du Corollaire \ref{pleinffMS1} que si $\G$ et $\H$ sont des $2$-groupes, le foncteur nerf de Segal induit un carré commutatif:

$$
\xymatrix@C+20pt@R-5pt{
\mathrm{Hom}_{\text{$2$-${\bf Grp}$}}\big(\G,\H\big) \ar[r]^-{\mathcal{N}_{\mathcal{S}}}_-{\cong} \ar[d]&   \mathrm{Hom}_{\bf MS}\big(\mathcal{N}_{\mathcal{S}}\G,\mathcal{N}_{\mathcal{S}}\H\big)\ar[d]\\
\pi_0\Big(\underline{\mathrm{Hom}}_{\text{$2$-${\bf Grp}$}}\big(\G,\H\big)\Big) \ar[r]_-{\pi_0\eqref{eqpleinffMS1}} & \pi_0\Big(\underline{\mathrm{Hom}}^{(1)}_{\bf MS}\big(\mathcal{N}_{\mathcal{S}}\G,\mathcal{N}_{\mathcal{S}}\H\big)\Big)}
$$
dont les fonctions horizontales sont bijectives. 

En particulier, vu que d'après la Proposition \ref{fibrachii} l'ensemble bisimplicial $\mathcal{N}_{\mathcal{S}}(\H)$ est un objet fibrant de la catégorie de modèles simplicial pointée $({\bf MS},{\bf W}^{diag}_{2}, {\bf mono},{\bf fib}_{2}^{diag})$ de la Proposition \ref{moduno}, le foncteur nerf de Segal $\mathcal{N}_{\mathcal{S}}\colon\xymatrix@C-12pt{\text{$2$-${\bf Grp}$} \ar[r] & {\bf MS}}$ induit un foncteur pleinement fidèle:
$$
\xymatrix@C+15pt{
\text{$2$-$h{\bf Grp}$} \ar[r]^-{h\mathcal{N}_{\mathcal{S}}} & \mathrm{Ho}_2\big({\bf MS}\big)\,,
}
$$ 

Rappelons d'un autre côté que la adjonction \eqref{duggerQeq} de la page \pageref{duggerQeq}:
$$
{\bf MS} \xymatrix@C+15pt{\phantom{a}\ar@/_12pt/[r]_-{(\,\cdot\,)_{0,\bullet}}\ar@{}[r]|-{\perp}&\ar@/_12pt/[l]_-{p^*\,=\,\text{préfaisceau constant}}\phantom{a}}
\simp_0 
$$
est une équivalence de Quillen entre $({\bf MS},{\bf W}^{diag}_{2}, {\bf mono},{\bf fib}_{2}^{diag})$ et $(\simp_{0},{\bf W}^{red}_{2}, {\bf mono},{\bf fib}_{2}^{red})$. Donc le diagramme commutatif de foncteurs:
$$
\xymatrix@C+18pt@R=3pt{
& \simp_0\\
\text{$2$-${\bf Grp}$}\ar@/^10pt/[ru]^-{\mathcal{N}} \ar@/_10pt/[rd]_{\mathcal{N}_{\mathcal{S}}}&\\
& {\bf MS} \ar[uu]_-{(\,\cdot\,)_{0,\bullet}}
}
$$
induit un diagramme commutatif à isomorphisme près:
$$
\xymatrix@C+18pt@R=3pt{
& \mathrm{Ho}_2\big(\simp_0\big)\\
\text{$2$-$h{\bf Grp}$}\ar@/^10pt/[ru]^-{h\mathcal{N}} \ar@/_10pt/[rd]_{h\mathcal{N}_{\mathcal{S}}}\ar@{}[r]|(.55){\cong}&\\
& \mathrm{Ho}_2\big({\bf MS}\big) \ar[uu]_-{{\bf R}(\,\cdot\,)_{0,\bullet}}
}
$$
où $h\mathcal{N}$ est un foncteur essentiellement surjectif (voir le Corollaire \ref{grho}), $h\mathcal{N}_{\mathcal{S}}$ est pleinement fidèle et ${\bf R}(\,\cdot\,)_{0,\bullet}$ est une équivalence de catégories.

Donc:

\begin{corollaire}\label{equiv2grpMS2}
Le foncteur nerf $\mathcal{N}\colon\xymatrix@C-12pt{\text{$2$-${\bf Grp}$} \ar[r] & \simp_0}$ induit une équivalence de catégories:
\begin{equation}\label{equiv2grpMS2foncteur00}
\xymatrix@C+15pt{
\text{$2$-$h{\bf Grp}$} \ar[r]^-{h\mathcal{N}} & \mathrm{Ho}_2\big(\simp_0\big)
}
\end{equation}
de la catégorie homotopique des $2$-groupes vers a catégorie des $2$-types d'homotopie réduits, et le foncteur nerf de Segal $\mathcal{N}_{\mathcal{S}}\colon\xymatrix@C-12pt{\text{$2$-${\bf Grp}$} \ar[r] & {\bf MS}}$ induit une équivalence de catégories:
\begin{equation}\label{equiv2grpMS2foncteur}
\xymatrix@C+15pt{
\text{$2$-$h{\bf Grp}$} \ar[r]^-{h\mathcal{N}_{\mathcal{S}}} & \mathrm{Ho}_2\big({\bf MS}\big)
}
\end{equation}
de la catégorie homotopique des $2$-groupes vers la catégorie homotopique des $2$-groupes de Segal.
\end{corollaire}

Si $X$ est un ensemble simplicial réduit, un $2$-groupe $\G$ tel que $h\mathcal{N}(\G)$ soit isomorphe au $2$-type d'homotopie réduit de $X$ dans la catégorie $\mathrm{Ho}_2\big(\simp_0\big)$ est dit un $2$-\emph{groupe d'homotopie} de $X$. De façon analogue un $2$-\emph{groupe d'homotopie} d'un pré-monoïde de Segal $X$ est un $2$-groupe $\G$ tel que $h\mathcal{N}_{\mathcal{S}}(\G)$ soit isomorphe à $X$ dans la catégorie des $2$-groupes de Segal $\mathrm{Ho}_2\big({\bf MS}\big)$. D'après le Corollaire \ref{equiv2grpMS2} tout ensemble simplicial réduit (resp. tout pré-monoïde de Segal) admet un $2$-groupe d'homotopie et ce $2$-groupe est unique à équivalence près dans la $2$-catégorie $\text{$2$-$\underline{\bf Grp}$}$. (Voir par exemple \cite{moerdijk}).

On définit le (ou un) $2$-groupe d'homotopie d'un ensemble simplicial pointé $X$ comme le (ou un) $2$-groupe d'homotopie de l'ensemble simplicial réduit ${\bf R}\mathcal{H}(X)$ (voir le Lemme \ref{hachelemme}). 

\renewcommand{\thesubsection}{\S\thesection.\arabic{subsection}}
\subsection{}\;
\renewcommand{\thesubsection}{\thesection.\arabic{subsection}}

Si $X$ est un pré-monoïde de Segal (\emph{i.e.} un objet de ${\bf MS}$) et $\G$ est un $2$-groupe, un \emph{déterminant de $X$ à valeurs $\G$} est par définition une couple $D=(D,T)$ composée par un morphisme d'ensemble simpliciaux (vers le nerf du groupoïde sous-jacent au $2$-groupe $\G$) 
$$
\xymatrix{X_{\bullet, 1} \ar[r]^-{D} & \mathrm{N}(\G)}
$$
et une fonction d'ensembles:
$$
\xymatrix{X_{0,2} \ar[r]^-{T} & \{\text{Morphismes de $\G$}\}}\,,
$$
tels que:
\begin{enumerate}
\item (Compatibilité) \, Si $\xi\in X_{0,2}$, $T(\xi)$ est un morphisme de $\G$ de la forme: 
$$
\xymatrix@C+8pt{D_0(d^v_2\xi) \otimes D_0(d^v_0\xi)  \ar[r]^-{T(\xi)} & D_0(d^v_1\xi) \,;}
$$
\item (Unitaire) \, On a que $D_0\big(s_0^v(\star)\big)\,=\,\mathbb{1}$ et $T\big(s_0^vs_0^v(\star)\big)\,=\,l_{\mathbb{1}}^{-1}\,=\,r_{\mathbb{1}}^{-1}$.
\item (Naturalité)\, Si $\zeta \in X_{1,2}$ on a un carré commutatif:
$$
\xymatrix@C+10pt{
D_0(d^v_2 d^h_1 \, \zeta) \otimes D_0(d^v_0 d^h_1 \, \zeta)\ar[d]_-{D_1(d^v_2 \, \zeta)\otimes D_1(d^v_0 \, \zeta)} \ar[r]^-{T(d^h_1 \, \zeta)} & D_0(d^v_1 d^h_1 \, \zeta)   \ar[d]^-{D_1(d^v_1 \, \zeta)}\\
D_0(d^v_2 d^h_0 \, \zeta) \otimes D_0(d^v_0 d^h_0 \, \zeta) \ar[r]_-{T(d^h_0 \, \zeta)} & D_0(d^v_1 d^h_0 \, \zeta) 
}
$$

(Rappelons que $d^h_i d^v_j \, = \, d^v_j d^h_i$).

\item (Associativité)\, Si $\eta\in X_{0,3}$ on a un diagramme commutatif:
\begin{equation}\label{pentafinipli}
\vcenter{\xymatrix@R-10pt@C+35pt{
D_0(A_{03}) &D_0(A_{01})\otimes D_0(A_{13}) \ar[l]_-{T(d_2^v\,\eta)}\\&\\
&D_0(A_{01}) \otimes \big(D_0(A_{12})\otimes D_0(A_{23})\big) \ar[uu]_-{D_0(A_{01})\otimes T(d_0^v\, \eta) }\ar@{}[d]_-{a}|-{\text{\rotatebox[origin=c]{90}{$\cong$}}}\\
D_0(A_{02})\otimes D_0(A_{23}) \ar[uuu]^-{T(d_1^v\,\eta)} &\big(D_0(A_{01})\otimes D_0(A_{12}))\otimes D_0(A_{23}\big)\ar[l]^-{T(d_3^v\,\eta)\otimes D_0(A_{23})}\,,}}
\end{equation}
\begin{align*}
\text{où} \qquad 
A_{03} \, = \, & d^v_{1} d^v_{1} \eta \, = \, d^v_{1} d^v_{2} \eta \qquad
A_{01} \, = \,  d^v_{2} d^v_{2} \eta \, = \, d^v_{2} d^v_{3} \eta \qquad
A_{13} \, = \,  d^v_{1} d^v_{0} \eta \, = \, d^v_{0} d^v_{2} \eta \\
A_{02} \, = \, & d^v_{2} d^v_{1} \eta \, = \, d^v_{1} d^v_{3}  \eta \qquad
A_{23} \, = \,  d^v_{0} d^v_{0} \eta \, = \, d^v_{0} d^v_{1} \eta \qquad
A_{12} \, = \,  d^v_{2} d^v_{0} \eta \, = \, d^v_{0} d^v_{3} \eta \,.
\end{align*}
\end{enumerate}

Notons ${\bf det}_X(\G)$ l'ensemble des déterminants d'un pré-monoïdes de Segal $X$ à valeurs dans un $2$-groupe $\G$. On définit un foncteur:  
\begin{equation}\label{tututuma}
\xymatrix@R=5pt@C+15pt{
{\bf MS}^{op} \, \times \, \text{$2$-${\bf Grp}$} \ar[r]  &  {\bf Ens}\,,\\
\big(X,\G\big) \;\, \ar@{}[r]|-{\longmapsto} & {\bf det}_X(\G)}
\end{equation}
comme suit: Si $f\colon\xymatrix@C-10pt{Y\ar[r]&X}$ est un morphisme de ${\bf MS}$ et $(\varphi,m^{\varphi})\colon\xymatrix@C-10pt{\G\ar[r]&\H}$ est un morphisme de $2$-groupes, on définit la fonction:
$$
\xymatrix@C+26pt@R=5pt{
{\bf det}_X(\G) \ar[r]^-{{\bf det}_f(\varphi,m^\varphi)} & {\bf det}_Y(\H)
}
$$
par la formule ${\bf det}_f\big(\varphi,m^\varphi\big)(D,T)\,=\,(\overline{D},\overline{T})$, où $\overline{D}$ et $\overline{T}$ sont les composés:
$$
\xymatrix@C-5pt{
Y_{\bullet,1} \ar[r]^-{f_{\bullet,1}} & X_{\bullet,1} \ar[r]^-{D} & \mathcal{N}_{\mathcal{S}}\big(\G\big)_{\bullet,1} 
\ar[rr]^-{\mathcal{N}_\mathcal{S}(\varphi,m^\varphi)_{\bullet,1}} && \mathcal{N}_{\mathcal{S}}\big(\H\big)_{\bullet,1}}
$$
$$\text{et}$$
$$
\xymatrix@C-5pt{
Y_{0,2} \ar[r]^-{f_{0,2}} & X_{0,2} \ar[r]^-{T} & \mathcal{N}_{\mathcal{S}}\big(\G\big)_{0,2} 
\ar[rr]^-{\mathcal{N}_\mathcal{S}(\varphi,m^\varphi)_{0,2}} && \mathcal{N}_{\mathcal{S}}\big(\H\big)_{0,2}\,.}
$$

Si la couple $(\overline{D},\overline{T})$ ainsi définie est un déterminant de $Y$ à valeurs dans $\H$, on vérifie sans difficulté que ${\bf det}_f(\varphi,m^\varphi)$ est effectivement un foncteur.

Montrons que $(\overline{D},\overline{T})$ est un déterminant. Pour montrer la propriété (i) remarquons que on a un diagramme:
$$
\xymatrix@C+5pt{
Y_{0,1} \ar@{}[rd]|-{\text{(I)}}\ar[r]^-{f_{0,1}} & X_{0,1} \ar[r]^-{D_0} \ar@{}[rd]|-{\text{(II)}}& \mathcal{N}_{\mathcal{S}}\big(\G\big)_{0,1} 
\ar[rr]^-{\mathcal{N}_\mathcal{S}(\varphi,m^\varphi)_{0,1}} \ar@{}[rrd]|-{\text{(III)}}&& \mathcal{N}_{\mathcal{S}}\big(\H\big)_{0,1}\\
Y_{0,2} \ar[u]^-{d_i^v}\ar[r]_-{f_{0,2}} & X_{0,2} \ar[r]_-{T} \ar[u]|-{d_i^v}& \mathcal{N}_{\mathcal{S}}\big(\G\big)_{0,2} 
\ar[rr]_-{\mathcal{N}_\mathcal{S}(\varphi,m^\varphi)_{0,2}} \ar[u]|-{d_i^v} && \mathcal{N}_{\mathcal{S}}\big(\H\big)_{0,2} \ar[u]_-{d_i^v} }
$$
où (I) et (III) sont de carrés commutatifs parce que $f$ et $(\varphi,m^\varphi)$ sont morphismes et (II) est commutatif vu que $(D,T)$ est un déterminant. Alors si $\xi'\in Y_{0,2}$ on a que $d_i^v\,T'(\xi')\,=\,\overline{D}_0\, d_i^v(\xi')$.

On vérifie la propriété (ii) c'est-à-dire que $\overline{D}_0\big(s_0^v(\star)\big)=\mathbb{1}$ et $T\big(s_0^vs_0^v(\star)\big)\,=\,l_{\mathbb{1}}^{-1}\,=\,r_{\mathbb{1}}^{-1}$ en remarquant qu'on a de diagrammes commutatifs:
$$
\xymatrix@C+5pt{
Y_{0,1} \ar[r]^-{f_{0,1}} & X_{0,1} \ar[r]^-{D_0} & \mathcal{N}_{\mathcal{S}}\big(\G\big)_{0,1} 
\ar[rr]^-{\mathcal{N}_\mathcal{S}(\varphi,m^\varphi)_{0,1}} && \mathcal{N}_{\mathcal{S}}\big(\H\big)_{0,1}\\
\star \ar[u]^-{s_0^v}\ar@{=}[r] & \star \ar@{=}[r] \ar[u]|-{s_0^v}& \star
\ar@{=}[rr] \ar[u]|-{s_0^v} && \star \ar[u]_-{s_0^v}}
$$
$$\text{et}$$
$$
\xymatrix@C+5pt{
Y_{0,2} \ar[r]^-{f_{0,2}} & X_{0,2} \ar[r]^-{T} & \mathcal{N}_{\mathcal{S}}\big(\G\big)_{0,2} 
\ar[rr]^-{\mathcal{N}_\mathcal{S}(\varphi,m^\varphi)_{0,2}} && \mathcal{N}_{\mathcal{S}}\big(\H\big)_{0,2}\\
\star \ar[u]^-{s_0^vs_0^v}\ar@{=}[r] & \star \ar@{=}[r] \ar[u]|-{s_0^vs_0^v}& \star
\ar@{=}[rr] \ar[u]|-{s_0^vs_0^v} && \star \ar[u]_-{s_0^vs_0^v}\,.}
$$

D'un autre côté, si $\zeta' \in Y_{1,2}$ on a un carré commutatif de $\G$:
$$
\xymatrix@C+25pt{
  D_0(d^v_2 d^h_1 \, f_{1,2}\zeta') \otimes D_0(d^v_0 d^h_1 \, f_{1,2}\zeta')  \ar[d]_-{D_1(d^v_2 \, f_{1,2}\zeta')\otimes D_1(d^v_0 \, f_{1,2}\zeta')}  \ar[r]^-{T(d^h_1 \, f_{1,2}\zeta')} & 
  D_0(d^v_1 d^h_1 \, f_{1,2}\zeta') \ar[d]^-{D_1(d^v_1 \, f_{1,2}\zeta')}  \\
D_0(d^v_2 d^h_0 \, f_{1,2}\zeta') \otimes D_0(d^v_0 d^h_0 \, f_{1,2}\zeta')  \ar[r]_-{T(d^h_0 \, f_{1,2}\zeta')} & D_0(d^v_1 d^h_0 \, f_{1,2}\zeta') \,;
}
$$
donc en appliquant le morphisme $(\varphi,m^{\varphi})$ on obtient le carré commutatif de $\H$:
$$
\xymatrix@C+15pt{
\overline{D}_0(d^v_2 d^h_1 \, \zeta') \otimes \overline{D}_0(d^v_0 d^h_1 \, \zeta')\ar[d]_-{\overline{D}_1(d^v_2 \, \zeta')\otimes \overline{D}_1(d^v_0 \, \zeta')} \ar[r]^-{\overline{T}(d^h_1 \, \zeta')} & \overline{D}_0(d^v_1 d^h_1 \, \zeta')   \ar[d]^-{\overline{D}_1(d^v_1 \, \zeta')}\\
\overline{D}_0(d^v_2 d^h_0 \, \zeta') \otimes \overline{D}_0(d^v_0 d^h_0 \, \zeta') \ar[r]_-{\overline{T}(d^h_0 \, \zeta')} & \overline{D}_0(d^v_1 d^h_0 \, \zeta') 
}
$$

On montre de façon analogue la propriété d'associativité.

\begin{proposition}\label{representpri}
Les foncteurs ${\bf det}_{\bullet_1}\big(\,\bullet_2\,\big)$ et $\mathrm{Hom}_{\bf MS} \big( \, \bullet_1\, , \, \mathcal{N}_{\mathcal{S}}(\bullet_2) \, \big)$ de la catégorie produit ${\bf MS}^{op}\times\text{$2$-${\bf Grp}$}$ vers la catégorie des ensembles ${\bf Ens}$ sont naturellement isomorphes.
\end{proposition}
\begin{proof}
Si $X$ est un pré-monoïde de Segal et $\G$ est un $2$-groupe on vérifie sans peine que la fonction naturelle:
\begin{equation}\label{transnaturrr}
\vcenter{\xymatrix@R=2pt@C+5pt{
\mathrm{Hom}_{\bf MS} \big( X, \mathcal{N}_{\mathcal{S}}(\G)\big)\ar[r] & {\bf det}_{X}\big(\G\big)\\
F\ar@{}[r]|-{\longmapsto} & (F_{\bullet,1},F_{0,2})
}}
\end{equation}
est bien définie; autrement dit on vérifie aussi-tôt que si $F\colon\xymatrix@C-5pt{X\ar[r]& \mathcal{N}_{\mathcal{S}}(\G)}$ est un morphisme de la catégorie ${\bf MS}$ la couple:
$$
\xymatrix{X_{\bullet, 1} \ar[r]^-{F_{\bullet,1}} & \mathcal{N}_{\mathcal{S}}(\G)_{\bullet,1}}\,=\,\mathrm{N}(\G)\qquad\text{et} \qquad \xymatrix{X_{0,2} \ar[r]^-{F_{0,2}} & \mathcal{N}_{\mathcal{S}}(\G)_{0,2} \, \subset \, \big\{\text{Morphismes \, de \, $\G$}\big\}}
$$
est bien un déterminant de $X$ à valeurs dans $\G$.

\vspace{.5cm}

\emph{La fonction \eqref{transnaturrr} est injective:}

Soient $F,G\colon\xymatrix@C-5pt{X\ar[r]& \mathcal{N}_{\mathcal{S}}(\G)}$ deux morphismes de pré-monoïdes de Segal tels que $F_{\bullet,1}= G_{\bullet,1}$ et $F_{0,2}= G_{0,2}$. Remarquons pour commencer que $F_{\bullet,0}=G_{\bullet,0}$ parce que $\mathcal{N}_{\mathcal{S}}(\G)_{\bullet,0}$ est l'ensemble simplicial constant $\star$.  En plus $F_{0,3}=G_{0,3}$ parce que $\mathcal{N}_{\mathcal{S}}(\G)_{0,\bullet} \, = \, \mathcal{N}(\G)$ est un ensemble simplicial faiblement $2$-cosquelettique (voir le Lemme \ref{3cosqull}) et $F_{0,2}= G_{0,2}$. Enfin on a que $F_{1,2}=G_{1,2}$ parce que la fonction \eqref{fonctionfinalll} est bijective et on a que $F_{0,2}= G_{0,2}$ et $F_{1,1}= G_{1,1}$.

Donc $F_{p,q}= G_{p,q}$ si $p,q\geq 0$ et $p+q \leq 3$. Il se suit que $F=G$ d'après le Corollaire \ref{fernandaadnanref}.

\vspace{.5cm}

\emph{La fonction \eqref{transnaturrr} est surjective:}

Montrons par ailleurs:

\begin{lemme}\label{pttlmsoufer}
Soit $X$ un objet de ${\bf MS}$ et $\G$ un $2$-groupe. Si $(D,T)$ est un déterminant de $X$ à valeurs dans $\G$ pour tout élément $A$ de $X_{0,1}$ on a une égalité de morphismes $T\big(s_i^v(A)\big)=s_i^v\big(D_0(A)\big)$ pour $0\leq i\leq 1$; c'est-à-dire:
$$
\xymatrix@C+40pt{\mathbb{1} \otimes D_0(A) \ar[r]^-{T(s_0^v(A)) = l^{-1}_{D_0(A)}} & D_0(A)}
\qquad\text{et}\qquad 
\xymatrix@C+40pt{D_0(A)\otimes\mathbb{1} \ar[r]^-{T(s_1^v(A)) = r^{-1}_{D_0(A)}} & D_0(A)}\,.
$$
\end{lemme}
\begin{proof}
Si $A\in X_{0,1}$ montrons que $T(s_0^v(A)) = l^{-1}_{D_0(A)}$; on vérifie $T(s_1^v(A)) = r^{-1}_{D_0(A)}$ de façon analogue. 

Pour commencer on considère l'élément $s_0^vs_0^v(A)$ de $X_{0,3}$, alors d'après les propriétés (ii) et (iv) ci-dessus, on a le diagramme commutatif:
$$
\vcenter{\xymatrix@R-10pt@C+25pt{
D_0(A) &\mathbb{1}\otimes D_0(A) \ar[l]_-{T(s_0^v\,A)} \\&\\
&\mathbb{1} \otimes (\mathbb{1}\otimes D_0(A))  \ar[uu]_-{\mathbb{1}\otimes T(s_0^v\, A) } \ar@{}[d]_-{a}|-{\text{\rotatebox[origin=c]{90}{$\cong$}}}\\
\mathbb{1}\otimes D_0(A)  \ar[uuu]^-{T(s_0^v\,A)}&(\mathbb{1}\otimes \mathbb{1})\otimes D_0(A) \ar[l]^-{l^{-1}_\mathbb{1}\otimes D_0(A)} \,;}}
$$
autrement dit on a un triangle commutatif:
$$
\vcenter{\xymatrix{
\big(\mathbb{1}\otimes \mathbb{1}\big) \otimes D_0(A) \ar[rr]^-{a} \ar[rd]_-{l^{-1}_\mathbb{1}\otimes D_0(A)} && 
   \mathbb{1}\otimes \big(\mathbb{1}\otimes D_0(A)\big) \ar[ld]^-{\mathbb{1}\otimes T(s_0^v \, A)}\\
&\mathbb{1}\otimes D_0(A)&}}
$$

D'un autre vu que $l_\mathbb{1}=r_\mathbb{1}$ (voir le Lemme \ref{otrosdiagg}), d'après \eqref{triangle} on a le triangle commutatif:
$$
\vcenter{\xymatrix{
\big(\mathbb{1}\otimes \mathbb{1}\big) \otimes D_0(A) \ar[rr]^-{a} &&    \mathbb{1}\otimes \big(\mathbb{1}\otimes D_0(A)\big)\\
&\mathbb{1}\otimes D_0(A)\ar[ru]_-{\mathbb{1}\otimes l_{D_0(A)}}\ar[lu]^-{l_\mathbb{1}\otimes D_0(A)}&}}
$$

Donc $T\big(s_0^v (A)\big) \, = \, l^{-1}_{D_0(A)}$.
\end{proof}

D'après le Corollaire \ref{fernandaadnanref}  $\big(\text{voir aussi l'adjonction \eqref{diagon3re}}\big)$ il suffit de se donner un morphisme $F\colon\xymatrix@C-5pt{\mu_{3}^{\; *}(X)\ar[r]& \mu_{3}^{\; *}\circ \mathcal{N}_{\mathcal{S}}(\G)}$ de pré-monoïdes de Segal tronqués tel que: 
$$
\text{$F_{p,1}\,=\,D_p$ \; si \; $0\leq p\leq 2$ \qquad et \qquad $F_{0,2}\,=\,T\,.$}
$$

Pour commencer on définit:
\begin{equation}\label{poquitasfer}
\begin{split}
F_{p,0} \, = \, \mathrm{id}_\star \qquad  \text{si \;  $0\leq p\leq 3$}\,,  \\
F_{p,1} \, = \, D_p   \qquad  \text{si \;  $0\leq p\leq 2$}\phantom{\,,}  \\ 
\text{et} \qquad  F_{0,2} \, = \, T  \,.   \qquad  \phantom{\text{si \;  $0\leq p\leq 2$}\;,}
\end{split}
\end{equation}

Il se suit du Lemme \ref{pttlmsoufer} et des propriétés (i)-(ii) qui vérifie le déterminant $(D,T)$ que les fonctions $F_{p,q}$ définies dans \eqref{poquitasfer} commutent avec les opérateurs faces et dégénérescence.

Enfin d'après la définition de l'ensemble bisimplicial $\mathcal{N}_{\mathcal{S}}(\G)$ et la propriété (iii) $\big(\text{resp. (iv)}\big)$ qui vérifie $(D,T)$ il existe une seule fonction $F_{1,2}$ $\big(\text{resp. $F_{0,3}$}\big)$ laquelle commutent avec les opérateurs faces et dégénérescence.
\end{proof}

Si $(D,T)$ et $(D',T')$ sont deux déterminants d'un pré-monoïde de Segal $X$ à valeurs dans un $2$-groupe $\G$, un \emph{morphisme de déterminants} $H\colon\xymatrix@C-10pt{(D,T)\ar[r]&(D',T')}$ est un morphisme d'ensembles simpliciaux (vers le nerf du groupoïde sous-jacent à $\G$):
$$
\xymatrix@C+13pt{
X_{\bullet,1} \times \Delta^1 \ar[r]^-{H} & \mathrm{N}(\G) 
}
$$
vérifiant les propriétés: 
\begin{enumerate}
\item $H$ est une homotopie de $D$ vers $D'$, c'est-à-dire on a un diagramme commutatif de morphismes d'ensembles simpliciaux:
$$
\xymatrix@R=2pt@C+10pt{
X_{\bullet,1} \ar[rd]|-{\nu_1} \ar@/^12pt/[rrd]^-D& &\\
& X_{\bullet,1}\times \Delta^1 \ar[r]|-{H} & \mathrm{N}(\G) \\
X_{\bullet,1} \ar[ru]|-{\nu_0} \ar@/_12pt/[rru]_-{D'} & &
}
$$
où $\nu_i$ est le morphisme composé:
$$
\xymatrix@C+12pt{X_{\bullet,1}\,\cong \, X_{\bullet,1}\times \Delta^0 \ar[r]^-{X_{\bullet,1}\times \delta_i} & X_{\bullet,1}\times \Delta^1}\,.
$$
\item $H$ est une homotopie pointée \emph{i.e.} le morphisme d'ensembles simpliciaux composé: 
$$
\xymatrix@C+15pt{
\Delta^1 \, \cong \, \star \times \Delta^1 \, = \, X_{\bullet,0} \times \Delta^1 \ar[r]^-{s_0^v \times \Delta^1} & X_{\bullet,1} \times \Delta^1 \ar[r]^-{H} & \mathrm{N}(\G)
}
$$
est le morphisme constant à valeurs $\mathbb{1}$. De façon équivalente $H_1\big(\star,\mathrm{id}_{[1]}\big) = \mathrm{id}_{\mathbb{1}}$.
\item Si $\chi\in X_{1,2}$ on a un diagramme commutatif de $\G$:
$$
\xymatrix@+10pt{
D_0(d^v_2d_1^h\chi)\otimes D_0(d^v_0 d_1^h\chi)\ar[d]_-{H_1( d^v_2\,\chi,\mathrm{id}_{[1]})\otimes H_1( d^v_0\,\chi,\mathrm{id}_{[1]})}\ar[r]^-{T(d_1^h\chi)} & D_0(d^v_1d_1^h\chi)  \ar[d]^-{H_1( d^v_1\,\chi,\mathrm{id}_{[1]})}  \\
D'_0(d^v_2d_0^h\chi)\otimes D'_0(d^v_0 d_0^h\chi)  \ar[r]_-{T'(d_0^h\chi)} &  D'_0(d^v_1d_0^h\chi)\,.}
$$ 
\end{enumerate}

Rappelons que d'après la définition \eqref{MSHOM(1)} et l'adjonction \eqref{adjMS} on a un isomorphisme naturel:
\begin{equation*}
\begin{split}
\underline{\mathrm{Hom}}_{\mathbf{MS}}^{(1)} \big( X, \mathcal{N}_{\mathcal{S}}(\G)\big)_{1}
\,=&\,\mathrm{Hom}_{\mathbf{MS}}\Big( \big(X\times \mathrm{p}_1^*(\Delta^1)\big)\big/ \big(\star\times \mathrm{p}_1^*(\Delta^1) \big)  ,\mathcal{N}_{\mathcal{S}}(\G)  \Big)\\
\,\cong&\,\mathrm{Hom}_{\ssimp}\big(X\times \mathrm{p}_1^*(\Delta^1),\mathcal{N}_{\mathcal{S}}(\G)\big)\,.
\end{split}
\end{equation*}

On définit ainsi une fonction:
\begin{equation}\label{derderHom1}
\vcenter{\xymatrix@R=2pt@C+5pt{
\underline{\mathrm{Hom}}^{(1)}_{\bf MS} \big( X, \mathcal{N}_{\mathcal{S}}(\G)\big)_{1}\ar[r] & \left\{\overset{\text{Morphismes de}}{\underset{\text{à valeurs dans $\G$}}{\text{\tiny{$\text{déterminants de \; $X$}$}}}}\right\}\\  
\mathcal{F}\ar@{}[r]|-{\longmapsto} & H_\mathcal{F}
}}
\end{equation}
où $H_\mathcal{F}\colon \xymatrix@C-10pt{(D_\mathcal{F},T_\mathcal{F})\ar[r] & (D'_\mathcal{F},T'_\mathcal{F})}$ est le morphisme du déterminant: 
\begin{equation*}
\def\objectstyle{\scriptstyle}
\def\labelstyle{\scriptstyle}
\xymatrix{(D_\mathcal{F},T_\mathcal{F})}  =
\Big(
\def\objectstyle{\scriptstyle}
\def\labelstyle{\scriptstyle}
\xymatrix@C-10pt{
X_{\bullet,1}\,\cong \, X_{\bullet,1}\times \Delta^0 \ar[rr]^-{X_{\bullet,1}\times \delta_1} && X_{\bullet,1}\times \Delta^1\ar[r]^-{\mathcal{F}_{\bullet,1}}&\mathrm{N}(\G)
},
\def\objectstyle{\scriptstyle}
\def\labelstyle{\scriptstyle}
\xymatrix@C-10pt{
X_{0,2} \,\cong \, X_{0,2}\times \Delta^0_0 \ar[rr]^-{X_{0,2}\times \delta_1} &&X_{0,2}\times \Delta^1_0\ar[r]^-{\mathcal{F}_{0,2}} & \mathcal{N}_{\mathcal{S}}(\G)_{0,2}
}
\Big)
\end{equation*}
vers le déterminant:
\begin{equation*}
\def\objectstyle{\scriptstyle}
\def\labelstyle{\scriptstyle}
\xymatrix{(D'_\mathcal{F},T'_\mathcal{F})}  =
\Big(
\def\objectstyle{\scriptstyle}
\def\labelstyle{\scriptstyle}
\xymatrix@C-10pt{
X_{\bullet,1}\,\cong \, X_{\bullet,1}\times \Delta^0 \ar[rr]^-{X_{\bullet,1}\times \delta_0} && X_{\bullet,1}\times \Delta^1\ar[r]^-{\mathcal{F}_{\bullet,1}}&\mathrm{N}(\G)
},
\def\objectstyle{\scriptstyle}
\def\labelstyle{\scriptstyle}
\xymatrix@C-10pt{
X_{0,2} \,\cong \, X_{0,2}\times \Delta^0_0 \ar[rr]^-{X_{0,2}\times \delta_0} &&X_{0,2}\times \Delta^1_0\ar[r]^-{\mathcal{F}_{0,2}} & \mathcal{N}_{\mathcal{S}}(\G)_{0,2}
}
\Big)
\end{equation*}
défini par le morphisme d'ensembles simpliciaux:
$$
\xymatrix@C+13pt{
X_{\bullet,1} \times \Delta^1 \ar[r]^-{\mathcal{F}_{\bullet,1}} & \mathrm{N}(\G) 
}\,.
$$

Notons que $H_\mathcal{F}$ vérifie la propriété (iii) parce que d'après la définition de $\mathcal{N}_{\mathcal{S}}(\G)$ si $\chi\in X_{1,2}$ on a que $\mathcal{F}_{1,2}(\chi,\mathrm{id}_{[1]})$ est le diagramme commutatif de $\G$ qui suit:
$$
\xymatrix@+18pt{
\mathcal{F}_{0,1}(d^v_2d_1^h\chi,\delta_1)\otimes \mathcal{F}_{0,1}(d^v_0 d_1^h\chi,\delta_1)\ar[d]_-{\mathcal{F}_{1,1}(d^v_2\,\chi,\mathrm{id}_{[1]})\otimes \mathcal{F}_{1,1}(d^v_0\,\chi,\mathrm{id}_{[1]})}\ar[r]^-{\mathcal{F}_{0,2}(d_1^h\chi,\delta_1)} & \mathcal{F}_{0,1}(d^v_1\chi,\delta_1)  \ar[d]^-{\mathcal{F}_{1,1}(d^v_1\,\chi,\mathrm{id}_{[1]})}  \\
\mathcal{F}_{0,1}(d^v_2d_0^h\chi,\delta_0)\otimes \mathcal{F}_{0,1}(d^v_0 d_0^h\chi,\delta_0)  \ar[r]_-{\mathcal{F}_{0,2}(d_0^h\chi,\delta_0)} &  \mathcal{F}_{0,1}(d^v_1d_0^h\chi,\delta_0)\,.}
$$ 

Constatons:

\begin{lemme}\label{igrphom1MS00}
La fonction \eqref{derderHom1} est bijective.
\end{lemme}
\begin{proof}
Soient $\mathcal{F},\mathcal{F'}\colon\xymatrix@C-8pt{X\times \mathrm{p}_1^*(\Delta^1)\ar[r]&\mathcal{N}_{\mathcal{S}}(\G) }$ deux morphismes d'ensembles bisimpliciaux tels que $H_\mathcal{F} = H_{\mathcal{F}'}$.

Il se suit de la définition de \eqref{derderHom1} que $\mathcal{F}_{\bullet,1}=\mathcal{F}'_{\bullet,1}$ et $\mathcal{F}_{0,2}=\mathcal{F}'_{0,2}$. D'un autre $\mathcal{F}_{\bullet,0}=\mathcal{F}'_{\bullet,0}$ parce que $\mathcal{N}_{\mathcal{S}}(\G) _{\bullet,0}$ est l'ensemble simplicial constant à valeurs $\star$. Vu que $\mathcal{N}_{\mathcal{S}}(\G) _{0,\bullet} = \mathcal{N}(\G)$ est un ensemble simplicial faiblement $2$-cosquelettique et $\mathcal{F}_{0,2}=\mathcal{F}'_{0,2}$ on a que $\mathcal{F}_{0,3}=\mathcal{F}'_{0,3}$. En plus $\mathcal{F}_{1,2}=\mathcal{F}'_{1,2}$ parce que $\mathcal{F}_{0,2}=\mathcal{F}'_{0,2}$, $\mathcal{F}_{1,1}=\mathcal{F}'_{1,1}$ et la fonction \eqref{fonctionfinalll} est injective.

Donc $\mathcal{F}=\mathcal{F'}$ d'après le Corollaire \ref{fernandaadnanref}. 

Supposons maintenant que $H\colon\xymatrix@C-10pt{(D,T)\ar[r]&(D',T')}$ est un morphisme de déterminants de $X$ à valeurs dans $\G$. On montre qu'il existe un morphisme d'ensembles bisimpliciaux: 
$$
\xymatrix@C+15pt{X\times \mathrm{p}_1^*(\Delta^1)\ar[r]^-{\mathcal{F}}&\mathcal{N}_{\mathcal{S}}(\G) }
$$ 
tel que $H_{\mathcal{F}}=H$ de la façon suivante: Par ailleurs il se suit de la Proposition \ref{representpri} qu'il existe un morphisme d'ensembles simpliciaux:
$$
\xymatrix@C+12pt{X_{0,\bullet}\times \Delta^1_0\ar[r]^-{\mathcal{F}_{0,\bullet}}&\mathcal{N}_{\mathcal{S}}(\G)_{0,\bullet} = \mathcal{N}(\G)}
$$ 
tel que:
$$
\mathcal{F}_{0,1} (A,\delta_0) \, = \,  D'_0(A)\;,  \qquad
\mathcal{F}_{0,1} (A,\delta_1) \, = \,  D_0(A) \qquad \text{si \; \, $A\in X_{0,1}$ \,.}
$$
$$\text{et}$$
$$
\mathcal{F}_{0,2} (\xi,\delta_0) \, = \,  T'(\xi)\;,  \qquad 
\mathcal{F}_{0,2} (\xi,\delta_1) \, = \,  T(\xi) \qquad \text{si \; \, $\xi\in X_{0,2}$ \,.}
$$

D'un autre le morphisme d'ensembles simpliciaux $\mathcal{F}_{\bullet,0}$ est par définition le morphisme constant sur $\mathcal{N}_{\mathcal{S}}(\G)_{\bullet,0}$ l'ensemble simplicial constant à valeurs $\star$ et $\mathcal{F}_{\bullet,1}$ est égale au morphisme d'ensembles simpliciaux:
$$
\xymatrix@C-8pt{
X_{\bullet,1}\times \Delta^1\ar[r]^-{H}&\mathcal{N}_{\mathcal{S}}(\G)_{\bullet,1} \, = \, \mathrm{N}(\G)\,.
}
$$

Notons que d'après la propriété (ii) qui vérifie le morphisme de déterminants $H$ on a que $\mathcal{F}_{\bullet,1} \circ s_0^v = s_0^v \circ \mathcal{F}_{\bullet,0}$. En plus d'après la propriété (iii) on a pour tout $\chi\in X_{1,2}$ un diagramme commutatif:
$$
\vcenter{\xymatrix@+18pt{
\mathcal{F}_{0,1}(d^v_2d_1^h\chi,\delta_1)\otimes \mathcal{F}_{0,1}(d^v_0 d_1^h\chi,\delta_1)\ar[d]_-{\mathcal{F}_{1,1}(d^v_2\,\chi,\mathrm{id}_{[1]})\otimes \mathcal{F}_{1,1}(d^v_0\,\chi,\mathrm{id}_{[1]})}\ar[r]^-{\mathcal{F}_{0,2}(d_1^h\chi,\delta_1)} & \mathcal{F}_{0,1}(d^v_1\chi,\delta_1)  \ar[d]^-{\mathcal{F}_{1,1}(d^v_1\,\chi,\mathrm{id}_{[1]})}  \\
\mathcal{F}_{0,1}(d^v_2d_0^h\chi,\delta_0)\otimes \mathcal{F}_{0,1}(d^v_0 d_0^h\chi,\delta_0)  \ar[r]_-{\mathcal{F}_{0,2}(d_0^h\chi,\delta_0)} &  \mathcal{F}_{0,1}(d^v_1d_0^h\chi,\delta_0)}}\;,
$$ 
or une fonction $\mathcal{F}_{1,2}$ commutant aux morphismes faces et dégénérescences qui le concerne.

Il se suit du Corollaire \ref{fernandaadnanref} que la fonction \eqref{derderHom1} est surjective.
 \end{proof}

Montrons:

\begin{lemme}\label{igrphom1MS}
Si $X$ est un pré-monoïde de Segal et $\G$ est un $2$-groupe, l'ensemble simplicial $\underline{\mathrm{Hom}}^{(1)}_{\bf MS} \big( X , \mathcal{N}_{\mathcal{S}}(\G) \big)$ est un $1$-groupoïde de Kan (voir le Lemme \ref{elelemmafer}).
\end{lemme}
\begin{proof}
D'après le Corollaire \ref{fernandaadnanref} $\underline{\mathrm{Hom}}^{(1)}_{\bf MS} \big( X , \mathcal{N}_{\mathcal{S}}(\G) \big)$ est un ensemble simplicial $2$-cosquelettique. En plus $\underline{\mathrm{Hom}}^{(1)}_{\bf MS} \big( X , \mathcal{N}_{\mathcal{S}}(\G) \big)$ est un complexe de Kan parce que $\mathcal{N}_{\mathcal{S}}(\G)$ est un objet fibrant de la catégorie de modèles simplicial pointée $({\bf MS},{\bf W}^{diag}_{2}, {\bf mono},{\bf fib}_{2}^{diag})$ de la Proposition \ref{moduno}.

Il se suit du Lemme \ref{uneequ} qu'il faut juste montrer que la fonction:
\begin{equation}\label{ferteeexxxtr}
\xymatrix@C+15pt{
\mathrm{Hom}_{\simp}\Big(\Delta^{2},\underline{\mathrm{Hom}}^{(1)}_{\bf MS} \big( X , \mathcal{N}_{\mathcal{S}}(\G) \big)\Big) \ar[r]^{\alpha^{1,k}_{X}} &\mathrm{Hom}_{\simp}\Big(\Lambda^{2,k},\underline{\mathrm{Hom}}^{(1)}_{\bf MS} \big( X , \mathcal{N}_{\mathcal{S}}(\G) \big)\Big)}
\end{equation}
est injective pour tout $0\leq k\leq 2$. 

Remarquons que d'après le Corollaire \ref{fernandaadnanref} et le Lemme \ref{igrphom1MS00} se donner un élément de la source de la fonction \eqref{ferteeexxxtr} équivaut à se donner trois morphismes de déterminants de $X$ à valeurs dans $\G$ comme dans le diagramme:
$$
\vcenter{\xymatrix@R=3pt@C+5pt{
&(D',T')\ar@/^3pt/[rd]^-{H^0}&\\
(D,T)\ar@/^3pt/[ru]^-{H^2}\ar[rr]_{H^1}&&(D'',T'') \,,
}}
$$
plus un morphisme d'ensembles simpliciaux:
$$
\vcenter{\xymatrix@C+10pt{
X_{\bullet,1} \times \Delta^2 \ar[r]^-{\Phi} & \mathcal{N}_{\mathcal{S}}\big(\G\big)_{\bullet,1} \, = \, \mathrm{N}(\G)
}}\,,
$$
tel que si $0\leq i\leq 2$ le morphisme $H^i\colon\xymatrix@C-10pt{X_{\bullet,1} \times \Delta^1\ar[r]& \mathrm{N}(\G)}$ est égale au composé:
$$
\vcenter{\xymatrix@C+10pt{
X_{\bullet,1} \times \Delta^1  \ar[r]^-{X_{\bullet,1} \times \delta_i} & X_{\bullet,1} \times \Delta^2 \ar[r]^-{\Phi} & \mathcal{N}_{\mathcal{S}}\big(\G\big)_{\bullet,1} \, = \, \mathrm{N}(\G)
}}\,.
$$

Or la fonction \eqref{ferteeexxxtr} est injective parce que $\underline{\mathrm{Hom}}_{\simp} \big( X_{\bullet,1} , \mathrm{N}(\G) \big)$ est un $1$-groupoïde de Kan.
\end{proof}

Si $X$ est un pré-monoïde de Segal et $\G$ est un $2$-groupe, il se suit de la Proposition \ref{representpri} et des Lemmes \ref{igrphom1MS00} et \ref{igrphom1MS} que les déterminants de $X$ à valeurs dans $\G$ et les morphismes entre eux forment un groupoïde qu'on note $\underline{{\bf det}}_X(\G)$.

De façon explicite si:
$$
\xymatrix@C+10pt{(D,T)\ar[r]^-{H}&(D',T')\ar[r]^-{H'}&(D'',T'')}
$$
sont de morphismes de déterminants de $X$ à valeurs dans $\G$ on définit le morphisme de déterminants composé $H'\circ H\colon\xymatrix@C-10pt{(D,T)\ar[r]&(D'',T'')}$ comme suit: Vu que d'après le Corollaire \ref{pointegroupdff} l'ensemble simplicial $\underline{\mathrm{Hom}}_{\simp_\star}\big(X_{\bullet,1},\mathrm{N}(\G)\big)$ est un $1$-groupoïde de Kan, il existe un seul morphisme:
$$
\xymatrix@C+13pt{
X_{\bullet,1} \times \Delta^2 \ar[r]^-{\Phi} & \mathrm{N}(\G) 
}
$$
tel que le morphisme d'ensembles simpliciaux composé: 
$$
\xymatrix@C+12pt{X_{\bullet,1}\times \Delta^1 \ar[r]^-{X_{\bullet,1}\times \delta_0} & X_{\bullet,1}\times \Delta^2 \ar[r]^-{\Phi} & \mathrm{N}(\G) }
$$
est égale au morphisme $H'$, le morphisme d'ensembles simpliciaux composé: 
$$
\xymatrix@C+12pt{X_{\bullet,1}\times \Delta^1 \ar[r]^-{X_{\bullet,1}\times \delta_2} & X_{\bullet,1}\times \Delta^2 \ar[r]^-{\Phi} & \mathrm{N}(\G) }
$$
est égale au morphisme $H$ et le morphisme d'ensembles simpliciaux composé: 
$$
\xymatrix@C+15pt{
\Delta^2 \, \cong \, \star \times \Delta^2 \, = \, X_{\bullet,0} \times \Delta^2 \ar[r]^-{s_0^v \times \Delta^2} & X_{\bullet,1} \times \Delta^2 \ar[r]^-{\Phi} & \mathrm{N}(\G)
}
$$
est égale le morphisme constant à valeurs $\mathbb{1}$.

Alors le morphisme de déterminants $H'\circ H\colon\xymatrix@C-10pt{(D,T)\ar[r]&(D'',T'')}$ est le morphisme d'ensembles simpliciaux composé: 
$$
\xymatrix@C+12pt{X_{\bullet,1}\times \Delta^1 \ar[r]^-{X_{\bullet,1}\times \delta_1} & X_{\bullet,1}\times \Delta^2 \ar[r]^-{\Phi} & \mathrm{N}(\G) }\,.
$$

On définit un foncteur:
\begin{equation}\label{tututuma2}
\xymatrix@R=5pt@C+15pt{
{\bf MS}^{op} \, \times \, \text{$2$-${\bf Grp}$} \ar[r]  &  {\bf Grpd}\\
\big(X,\G\big) \;\, \ar@{}[r]|-{\longmapsto} & \underline{{\bf det}}_X(\G)}
\end{equation}
qui étend le foncteur \eqref{tututuma} comme suit: Si $f\colon\xymatrix@C-15pt{Y\ar[r]&X}$ est un morphisme de ${\bf MS}$ et $(\varphi,m^{\varphi})\colon\xymatrix@C-13pt{\G\ar[r]&\H}$ est un morphisme de $2$-groupes, on définit le foncteur :
$$
\xymatrix@C+34pt@R=1pt{
\underline{\bf det}_X(\G) \ar[r]^-{\underline{\bf det}_f(\varphi,m^{\varphi})} & \underline{\bf det}_Y(\H)\\
(D,T)  \ar[dddd]_-{H} &  (\overline{D},\overline{T}) \ar[dddd]^-{\overline{H}} \\  & \\  \ar@{}[r]|-{\longmapsto} & \\    & \\  
(D',T') & (\overline{D}',\overline{T}')}
$$
où $\overline{H}\colon\xymatrix@C-10pt{(\overline{D},\overline{T})\ar[r]&(\overline{D}',\overline{T}')}$ est le morphisme de déterminants défini par le composé de foncteurs:
$$
\vcenter{\xymatrix@C+13pt{
Y_{\bullet,1} \times \Delta^1 \ar[r]^-{f_{\bullet,1}}& X_{\bullet,1} \times \Delta^1 \ar[r]^-{H} & \mathrm{N}(\G)\ar[r]^-{\mathrm{N}(\varphi)} &  \mathrm{N}(\H) 
}}\,.
$$

De la Proposition \ref{representpri} et des Lemmes \ref{igrphom1MS00} et \ref{igrphom1MS} on en conclu:

\begin{corollaire}\label{fininininifer}
Les foncteurs $\mathrm{N}\Big(\underline{\bf det}_{\bullet_1}\big(\,\bullet_2\,\big)\Big)$ et $\underline{\mathrm{Hom}}^{(1)}_{\bf MS} \big( \, \bullet_1\, , \, \mathcal{N}_{\mathcal{S}}(\bullet_2) \, \big)$ de la catégorie produit ${\bf MS}^{op} \times \text{$2$-${\bf Grp}$}$ vers la catégorie des ensembles simpliciaux $\simp$ sont naturellement isomorphes.
\end{corollaire}

Rappelons que la catégorie homotopique des groupoïdes $h{\bf Grpd}$ est la catégorie dont les objets sont les groupoïdes et les morphismes sont les classes à isomorphisme naturel près des foncteurs entre eux: 
$$
\mathrm{Hom}_{h{\bf Grpd}}\big( \G , \H\big) \, = \, \pi_0\big(\mathcal{H}^{\mathcal{G}}\big)\,.
$$ 

Dans le Corollaire qui suit on note:
$$
\xymatrix@C+10pt{{\bf Grpd}\ar[r]^-{\pi}&h{\bf Grpd}}
$$ 
le foncteur canonique. Remarquons que $\underline{\mathrm{Hom}}_{\text{$2$-${\bf Grp}$}}$ est un enrichissement de la catégorie des $2$-groupes sur ${\bf Grpd}$, donc $\underline{\mathrm{Hom}}_{\text{$2$-${\bf Grp}$}}$ est aussi un enrichissement de $\text{$2$-${\bf Grp}$}$ sur la catégorie homotopique $h{\bf Grpd}$.  

\begin{corollaire}\label{repsteouni}
Si $X$ est un pré-monoïde de Segal le foncteur composé:
$$
\xymatrix@R=5pt@C+3pt{
\text{$2$-${\bf Grp}$} \ar[rr]^-{\underline{{\bf det}}_X(\,\bullet\,)} & &  {\bf Grpd}  \ar[r]^-{\pi}& h{\bf Grpd}
}
$$
est représentable par un $2$-groupe d'homotopie du pré-monoïde de Segal $X$.

En particulier si $(\varphi,m^{\varphi})\colon \xymatrix@C-10pt{\G\ar[r]& \H}$ est une $2$-équivalence faible de $2$-groupes le foncteur:
$$
\xymatrix@C+35pt{
 \underline{{\bf det}}_X(\G) \ar[r]^-{ \underline{{\bf det}}_X(\varphi,m^{\varphi})} &  \underline{{\bf det}}_X(\H)
}
$$
est une équivalence faible de groupoïdes et le foncteur induit:
$$
\xymatrix@R=5pt@C+40pt{
\text{$2$-$h{\bf Grp}$} \ar[r]^-{\pi_0\big(\underline{{\bf det}}_X(\,\bullet\,)\big)}   &  {\bf Ens}}
$$
est représentable aussi par un $2$-groupe d'homotopie du pré-monoïde de Segal $X$.

De façon analogue si on considère la catégorie ${\bf MS}$ munie de l'enrichissement $\underline{\mathrm{Hom}}^{(1)}_{\bf MS}$ le foncteur composé:
$$
\xymatrix@R=5pt@C+20pt{
{\bf MS}^{op} \ar[r]^-{\underline{{\bf det}}_{\,\bullet\,}(\G)}   &  {\bf Grpd} \ar[r]^-{\mathrm{N}} & \simp} 
$$
est représentable par le pré-monoïde de Segal $\mathcal{N}_{\mathcal{S}}(\G)$ pour tout $2$-groupe $\G$.

En particulier si $f\colon \xymatrix@C-10pt{X\ar[r]& Y}$ est une $2$-équivalence faible diagonale de pré-monoïdes de Segal le morphisme:
$$
\xymatrix@C+20pt{
\mathrm{N}\big(\underline{{\bf det}}_Y(\G)\big) \ar[r]^-{\mathrm{N}\big(\underline{{\bf det}}_f(\G)\big)} &  \mathrm{N}\big(\underline{{\bf det}}_X(\G)\big)
}
$$
est une $1$-équivalence faible d'ensembles simpliciaux et le foncteur induit:
$$
\xymatrix@R=5pt@C+30pt{
\mathrm{Ho}_2\big({\bf MS}\big)^{op} \ar[r]^-{h\mathrm{N}\big(\underline{{\bf det}}_{\,\bullet\,}(\G)\big)}   &  \mathrm{Ho}_1(\simp)}
$$
est représentable aussi par le pré-monoïde de Segal $\mathcal{N}_{\mathcal{S}}(\G)$.
\end{corollaire} 
\begin{proof}
Soit $X$ un pré-monoïde de Segal et $\Pi_2(X)$ un $2$-groupe d'homotopie de $X$. D'après la définition de $\Pi_2(X)$ il existe un isomorphisme $\xymatrix@C-10pt{\mathcal{N}_{\mathcal{S}}\big(\Pi_2(X)\big) \ar[r] & X}$ dans la catégorie $\mathrm{Ho}_2\big({\bf MS}\big)$

Vu que $({\bf MS},{\bf W}^{diag}_{2}, {\bf mono},{\bf fib}_{2}^{diag})$ munie de l'enrichissement $\underline{\mathrm{Hom}}^{(1)}_{\bf MS}$ est une catégorie de modèles simplicial d'ordre $1$ (voir la Proposition \ref{moduno}) il se suit que le foncteur $\underline{\mathrm{Hom}}^{(1)}_{\bf MS} \big( \bullet ,\mathcal{N}_{\mathcal{S}}(\G)\big)$ de la catégorie ${\bf MS}^{op}$ vers $\simp$ envoi les $2$-équivalences faibles diagonales entre des objets quelconques vers $1$-équivalences faibles d'ensembles simpliciaux. Donc il existe un isomorphisme dans $\mathrm{Ho}_1(\simp)$:
$$
\underline{\mathrm{Hom}}^{(1)}_{\bf MS} \Big(\mathcal{N}_{\mathcal{S}}\big(\Pi_2(X), \mathcal{N}_{\mathcal{S}}(\G)\big)\Big)  \; \cong \;
\underline{\mathrm{Hom}}^{(1)}_{\bf MS} \big(X,\mathcal{N}_{\mathcal{S}}(\G)\big)
$$
lequel est naturel en $\G$.

D'un autre côté d'après le Corollaire \ref{pleinffMS1} et le Lemme \ref{fininininifer} ils existent des isomorphismes naturels dans la catégorie $\simp$: 
$$
\mathrm{N}\Big(\underline{\mathrm{Hom}}_{\text{$2$-${\bf Grp}$}} \big(\Pi_2(X),\G\big)\Big) \; \cong \; 
\underline{\mathrm{Hom}}^{(1)}_{\bf MS} \Big(,\mathcal{N}_{\mathcal{S}}\big(\Pi_2(X)\big),\mathcal{N}_{\mathcal{S}}(\G)\Big)
$$
$$\text{et}$$
$$
\mathrm{N}\big(\underline{\bf det}_X(\G)\big)  \; \cong \; 
\underline{\mathrm{Hom}}^{(1)}_{\bf MS} \big(X,\mathcal{N}_{\mathcal{S}}(\G)\big)\,.
$$

Donc il existe un isomorphisme dans $h{\bf Grpd}$ naturel en $\G$:
$$
\underline{\bf det}_X(\G)\;\cong\; \underline{\mathrm{Hom}}_{\text{$2$-${\bf Grp}$}} \big(\Pi_2(X),\G\big)
$$
parce que le foncteur nerf $\mathrm{N}\colon\xymatrix@C-10pt{{\bf Grpd}\ar[r]&\simp}$ induit une équivalence de catégories (voir \ref{groupdnnerfe}):
$$
\xymatrix@C+10pt{h{\bf Grpd} \ar[r]^-{h\mathrm{N}} &\mathrm{Ho}_1(\simp)}
$$

Rappelons maintenant que si $(\varphi,m^{\varphi})\colon \xymatrix@C-10pt{\G\ar[r]& \H}$ est une $2$-équivalence faible de $2$-groupes:
$$
\xymatrix@C+18pt{\mathcal{N}_{\mathcal{S}}(\G)\ar[r]^-{\mathcal{N}_{\mathcal{S}}(\varphi,m^{\varphi})}& \mathcal{N}_{\mathcal{S}}(\H)}
$$ 
est une $2$-équivalence faible diagonale des objets fibrants de la catégorie de modèles simplicial $\big({\bf MS},{\bf W}^{diag}_{2}, {\bf mono},{\bf fib}_{2}^{diag},\underline{\mathrm{Hom}}^{(1)}_{\bf MS}\big)$. En particulier le morphisme induit:
$$
\xymatrix@C+20pt{
\underline{\mathrm{Hom}}^{(1)}_{\bf MS} \big(X,\mathcal{N}_{\mathcal{S}}(\G)\big) \ar[r]^-{\mathcal{N}_{\mathcal{S}}(\varphi,m^{\varphi})_*}& 
\underline{\mathrm{Hom}}^{(1)}_{\bf MS} \big(X,\mathcal{N}_{\mathcal{S}}(\H)\big)}
$$
est une $1$-équivalences faible de $1$-groupoïdes de Kan. 

Autrement dit:
$$
\xymatrix@C+25pt{
 \underline{{\bf det}}_X(\G) \ar[r]^-{ \underline{{\bf det}}_X(\varphi,m^{\varphi})} &  \underline{{\bf det}}_X(\H)
}
$$
est une équivalence faible de groupoïdes.

\end{proof}

On déduit des Corollaires \ref{equiv2grpMS2} et \ref{repsteouni}:

\begin{corollaire}\label{repsteouni3}
Si $X$ est un ensemble simplicial réduit, le foncteur induit (voir \ref{repsteouni2}):
$$
\xymatrix@R=5pt@C+40pt{
\text{$2$-$h{\bf Grp}$} \ar[r]^-{\pi_0\big(\underline{{\bf det}}_X(\,\bullet\,)\big)}   &  {\bf Ens}}
$$
est représentable par un $2$-groupe d'homotopie du pré-monoïde de Segal $X$.
\end{corollaire}

\bibliographystyle{amsalpha}
\bibliography{bib}
\end{document}